\documentclass[11pt]{article}
\usepackage{cmap}
\usepackage[T1]{fontenc}
\usepackage[utf8]{inputenc}
\usepackage{lmodern}
\usepackage{amsmath,amssymb,amsthm,mathtools,mathrsfs}
\usepackage{dsfont}
\usepackage{enumitem,array,booktabs,longtable,graphicx}
\usepackage{tikz}
\usetikzlibrary{arrows.meta,positioning,calc}
\usepackage[letterpaper,margin=1.15in,headheight=15pt]{geometry}
\usepackage{microtype}
\usepackage[section]{placeins}
\usepackage[font=small,labelfont=bf,labelsep=period]{caption}
\usepackage[dvipsnames]{xcolor}
\usepackage[backend=biber,style=alphabetic,sorting=nyt,maxbibnames=99,maxalphanames=4,natbib=true]{biblatex}
\appto{\bibsetup}{\emergencystretch=1.5em\hbadness=2000}

\usepackage{aliascnt}
\usepackage[colorlinks=true,linkcolor=MidnightBlue,citecolor=OliveGreen,urlcolor=BrickRed,bookmarksdepth=4]{hyperref}
\usepackage[nameinlink,noabbrev,capitalise]{cleveref}
\usepackage{fancyhdr}
\numberwithin{equation}{section}

\allowdisplaybreaks

\setlist[itemize]{leftmargin=2em,itemsep=0.25em,topsep=0.35em}
\setlist[enumerate]{leftmargin=2em,itemsep=0.25em,topsep=0.35em}
\newtheorem{theorem}{Theorem}[section]
\newaliascnt{proposition}{theorem}
\newtheorem{proposition}[proposition]{Proposition}
\aliascntresetthe{proposition}
\newaliascnt{lemma}{theorem}
\newtheorem{lemma}[lemma]{Lemma}
\aliascntresetthe{lemma}
\newaliascnt{corollary}{theorem}
\newtheorem{corollary}[corollary]{Corollary}
\aliascntresetthe{corollary}
\newaliascnt{conjecture}{theorem}

\aliascntresetthe{conjecture}
\newaliascnt{problem}{theorem}

\aliascntresetthe{problem}
\newaliascnt{hypothesis}{theorem}
\newtheorem{hypothesis}[hypothesis]{Hypothesis}
\aliascntresetthe{hypothesis}

\theoremstyle{definition}
\newaliascnt{definition}{theorem}
\newtheorem{definition}[definition]{Definition}
\aliascntresetthe{definition}
\newtheorem*{definition*}{Definition}
\newaliascnt{assumption}{theorem}

\aliascntresetthe{assumption}
\newaliascnt{example}{theorem}
\newtheorem{example}[example]{Example}
\aliascntresetthe{example}
\theoremstyle{remark}
\newaliascnt{remark}{theorem}
\newtheorem{remark}[remark]{Remark}
\aliascntresetthe{remark}
\crefname{theorem}{\mbox{Theorem}}{\mbox{Theorems}}
\Crefname{theorem}{\mbox{Theorem}}{\mbox{Theorems}}
\crefname{proposition}{\mbox{Proposition}}{\mbox{Propositions}}
\Crefname{proposition}{\mbox{Proposition}}{\mbox{Propositions}}
\crefname{lemma}{\mbox{Lemma}}{\mbox{Lemmas}}
\Crefname{lemma}{\mbox{Lemma}}{\mbox{Lemmas}}
\crefname{corollary}{\mbox{Corollary}}{\mbox{Corollaries}}
\Crefname{corollary}{\mbox{Corollary}}{\mbox{Corollaries}}
\crefname{conjecture}{\mbox{Conjecture}}{\mbox{Conjectures}}
\Crefname{conjecture}{\mbox{Conjecture}}{\mbox{Conjectures}}
\crefname{definition}{\mbox{Definition}}{\mbox{Definitions}}
\Crefname{definition}{\mbox{Definition}}{\mbox{Definitions}}
\crefname{assumption}{\mbox{Assumption}}{\mbox{Assumptions}}
\Crefname{assumption}{\mbox{Assumption}}{\mbox{Assumptions}}
\crefname{example}{\mbox{Example}}{\mbox{Examples}}
\Crefname{example}{\mbox{Example}}{\mbox{Examples}}
\crefname{remark}{\mbox{Remark}}{\mbox{Remarks}}
\Crefname{remark}{\mbox{Remark}}{\mbox{Remarks}}
\crefname{problem}{\mbox{Problem}}{\mbox{Problems}}
\Crefname{problem}{\mbox{Problem}}{\mbox{Problems}}
\crefname{hypothesis}{\mbox{Hypothesis}}{\mbox{Hypotheses}}
\Crefname{hypothesis}{\mbox{Hypothesis}}{\mbox{Hypotheses}}
\def\Xint#1{\mathchoice{\XXint\displaystyle\textstyle{#1}}{\XXint\textstyle\scriptstyle{#1}}%
{\XXint\scriptstyle\scriptscriptstyle{#1}}{\XXint\scriptscriptstyle\scriptscriptstyle{#1}}\!\int}
\def\XXint#1#2#3{{\setbox0=\hbox{$#1{#2#3}{\int}$}\vcenter{\hbox{$#2#3$}}\kern-.5\wd0}}
\def\fint{\Xint-}
\newcommand{\R}{\mathbb{R}}
\newcommand{\eps}{\varepsilon}
\newcommand{\cH}{\mathcal{H}}
\newcommand{\cK}{\mathcal{K}}
\newcommand{\cU}{\mathcal{U}}
\newcommand{\cO}{\mathcal{O}}

\newcommand{\cN}{\mathcal{N}}
\newcommand{\cE}{\mathcal{E}}

\newcommand{\Ad}{\operatorname{Ad}}

\newcommand{\rc}{\mathsf{c}}

\newcommand{\RCD}{\mathrm{RCD}}
\newcommand{\Fr}{\operatorname{Fr}}
\newcommand{\Isom}{\operatorname{Isom}}
\newcommand{\Lip}{\operatorname{Lip}}
\newcommand{\spec}{\operatorname{spec}}
\newcommand{\tr}{\operatorname{tr}}
\newcommand{\supp}{\operatorname{supp}}
\newcommand{\diam}{\operatorname{diam}}

\newcommand{\pr}{\mathrm{pr}}

\newcommand{\Hzero}{H^{1,2}_0}
\newcommand{\Hone}{H^{1,2}}
\newcommand{\T}{\mathsf{T}}
\newcommand{\Gr}{\operatorname{Gr}}
\newcommand{\Cov}{\operatorname{Cov}}
\newcommand{\cB}{\mathcal{B}}

\newcommand{\cD}{\mathcal{D}}
\newcommand{\av}{\operatorname{av}}

\newcommand{\ft}{\mathfrak{t}}
\newcommand{\cA}{\mathcal{A}}

\newcommand{\cR}{\mathcal{R}}

\newcommand{\cT}{\mathcal{T}}

\newcommand{\Rm}{\operatorname{Rm}}
\newcommand{\Ric}{\operatorname{Ric}}
\newcommand{\Scal}{\operatorname{Scal}}
\newcommand{\vol}{\operatorname{vol}}
\newcommand{\area}{\operatorname{area}}
\newcommand{\dis}{\operatorname{dis}}
\newcommand{\Int}{\operatorname{Int}}
\newcommand{\Cone}{\mathsf{C}}
\newcommand{\Z}{\mathbb{Z}}
\newcommand{\CP}{\mathbb{CP}}
\newcommand{\tB}{\widetilde{B}}

\newaliascnt{question}{theorem}

\aliascntresetthe{question}
\crefname{question}{\mbox{Question}}{\mbox{Questions}}
\Crefname{question}{\mbox{Question}}{\mbox{Questions}}
\newcommand{\cY}{\mathcal{Y}}

\newcommand{\sd}{\mathsf{d}}

\newcommand{\fg}{\mathfrak{g}}
\newcommand{\fk}{\mathfrak{k}}
\newcommand{\fm}{\mathfrak{m}}
\newcommand{\im}{\operatorname{im}}

\newcommand{\Sh}{\operatorname{Sh}}
\newcommand{\Bx}{\operatorname{Rect}}

\newcommand{\cW}{\mathcal{W}}

\newcommand{\sB}{\mathsf{B}}
\newcommand{\sV}{\mathsf{V}}

\newcommand{\sE}{\mathsf{E}}
\newcommand{\sQ}{\mathsf{Q}}

\newcommand{\cC}{\mathcal{C}}
\newcommand{\cP}{\mathcal{P}}

\newcommand{\sa}{\mathsf{a}}
\newcommand{\af}{\mathsf{f}}
\newcommand{\eGHto}{\xrightarrow{\mathrm{eGH}}}
\newcommand{\one}{\mathds{1}}

\crefname{section}{\mbox{Section}}{\mbox{Sections}}
\Crefname{section}{\mbox{Section}}{\mbox{Sections}}
\crefname{subsection}{\mbox{Section}}{\mbox{Sections}}
\Crefname{subsection}{\mbox{Section}}{\mbox{Sections}}
\crefname{subsubsection}{\mbox{Section}}{\mbox{Sections}}
\Crefname{subsubsection}{\mbox{Section}}{\mbox{Sections}}
\crefformat{enumi}{#2#1#3}
\Crefformat{enumi}{#2#1#3}
\crefrangeformat{enumi}{#3#1#4--#5#2#6}
\Crefrangeformat{enumi}{#3#1#4--#5#2#6}
\crefmultiformat{enumi}{#2#1#3}{ and~#2#1#3}{, #2#1#3}{ and~#2#1#3}
\Crefmultiformat{enumi}{#2#1#3}{ and~#2#1#3}{, #2#1#3}{ and~#2#1#3}

\usepackage{needspace}
\makeatletter
\pretocmd{\l@section}{\Needspace{3\baselineskip}}{}{}
\renewcommand*{\l@subsection}{\@dottedtocline{2}{1.5em}{3.2em}}
\makeatother
\title{The Singular Fibration of Collapsing RCD Spaces\\ with Reifenberg Local Covers}
\author{%
Shaosai Huang\thanks{Kspectra Research Inc., Toronto, ON, M2N 0G3, Canada.
Email: \texttt{arthur.foxie.huang@kspectra.ai}.}}
\date{October 20, 2026}
\hypersetup{pdftitle={The Singular Fibration of Collapsing RCD Spaces with Reifenberg Local Covers},pdfauthor={Shaosai Huang},pdflang={en-GB}}
\begin{document}
\maketitle
\begin{center}\begin{minipage}{0.86\textwidth}
\small\textbf{Abstract.} Let compact $\mathrm{RCD}(\kappa,N)$ spaces converge in the measured
Gromov--Hausdorff sense, and suppose that they have $(\delta,\rho)$-Reifenberg local covers
for one fixed $\rho>0$, with $|\kappa|\rho^2\le1$ and $\delta$ small. We prove that, along a
subsequence, the spaces are singular fibrations over the limit: their canonical frame spaces
converge to a topological manifold with a locally linear action of $O(n)$ whose quotient is the
limit; for large indices the frame spaces fibre over this manifold, by maps that commute with
$O(n)$, with fibres homeomorphic to nilmanifolds; and the induced maps to the limit are, near
every point, quotients of products by the isotropy group, with fibres homeomorphic to
infranilmanifolds. As applications, we extend a fibration theorem of Honda and Zhang to
non-smooth RCD spaces whose local covers have injectivity radius bounded below, and we show
that for aspherical spaces the isotropy groups are finite, the fibre groups inject into the
fundamental groups, and the limit is the quotient of a contractible manifold by a proper
cocompact action; under Rong's condition $A_s$ it is a closed aspherical manifold. Three
appendices treat smooth manifolds with $\mathrm{Ric}\ge-(n-1)$ and small local rewinding
Reifenberg error: a uniform short-time Ricci flow with $|\mathrm{Rm}|\le\alpha/t$, also for
complete non-compact data, where it is unique among complete flows with curvature at most $A/t$; a second
construction with affine nilmanifold bundles; and, conditional on Nifa's construction, examples
showing that rewinding volume control alone does not give continuous Gromov--Hausdorff
approximations with small error.
\end{minipage}\end{center}
\tableofcontents
\section{Introduction}\label{ch:intro}

We study collapsing sequences of $\RCD(\kappa,N)$ spaces whose local universal covers, at
one fixed radius, are uniformly close to Euclidean space at all smaller scales. Along a
subsequence such spaces are topological singular fibrations over their limit: canonical
metric frames resolve the limiting isotropy, harmonic coordinates construct an equivariant
bundle on the frame spaces, and fibre topology identifies its quotient. The fibres are
nilmanifolds on the frame spaces and infranilmanifolds on the spaces themselves. Two
applications follow: a fibration theorem under local bounded covering geometry, and the
structure of limits of aspherical spaces.

\subsection{The Problem and the Hypothesis}\label{sec:question}

Let compact $\RCD(\kappa,N)$ spaces converge in the measured Gromov--Hausdorff sense. We
ask when they fibre over the limit, and with which fibres. The model is the
bounded-curvature theory: for closed Riemannian manifolds collapsing with bounded sectional
curvature, the orthonormal frame bundles converge to a manifold with an $O(n)$-action and
fibre equivariantly over it with nilmanifold fibres, and taking the quotient gives singular
fibrations with infranil fibres \cite{Fuk88,Fuk89,CFG92}. An $\RCD$ space has no frame
bundle in this sense, and a lower Ricci bound alone does not give this structure: a product
with a shrinking round sphere already shows why the geometry of local covers matters. Rong
proposed adding the local rewinding Reifenberg property \cite{Ron22,HRW20}. We use the
following version for $\RCD$ spaces.

\begin{definition}[Reifenberg Local Covers]\label{m:fr:def:localcovers}\label{m:fr:def:oneradius}
Fix an integer $n\ge1$ and $\delta,\rho>0$. An $\RCD(\kappa,N)$ space $(X,d,m)$ has
\emph{$(\delta,\rho)$-Reifenberg local covers} if, for every $x\in X$ and
$0<t\le\rho/3$, the ball $B_t(\tilde x)$ in the universal cover of $B_\rho(x)$,
with its lifted length metric, admits a map to $B_t(0^n)$ that distorts distances
by at most $\delta t$ and has $\delta t$-dense image.
Here $\tilde x$ is a lift of $x$; the condition is independent of
the lift. Thus $\rho$ is the covering radius and $\rho/3$ the comparison scale.
\end{definition}

The covers and their lifted measures are constructed in
\cref{m:fr:sec:localcovers}. The condition concerns every centre and every smaller
scale, but only one covering radius. It assumes neither agreement of neighbouring
covers nor a uniform injectivity radius on the space itself. The reference measure need not
be Hausdorff measure, and $n$ may be smaller than $N$.

\begin{hypothesis}[The Collapsing Sequence]\label{hyp:sequence}
Let $n\ge1$ be an integer, $N\ge\max\{n,2\}$, $\kappa\in\R$ and $\delta,\rho>0$, with
$|\kappa|\rho^2\le1$. The compact $\RCD(\kappa,N)$ spaces $(X_i,d_i,m_i)$ have
$m_i(X_i)=1$ and $(\delta,\rho)$-Reifenberg local covers for the same $n,\delta,\rho$,
and converge in the measured Gromov--Hausdorff sense to a compact space $(X,d,m)$.
\end{hypothesis}

The curvature-scale bound is a hypothesis, not a scaling normalization:
$\kappa\rho^2$ is scale invariant. For $\kappa\ge0$ we may use the lower bound $0$.
The hypothesis does not require a dimension drop; it includes the non-collapsing
case. The bound $N\ge2$ is no restriction, since an $\RCD(\kappa,N')$ space with $N'\le N$ is an
$\RCD(\kappa,N)$ space; it gives $N>1$, which the Laplacian comparison in the proof of
\cref{m:fr:prop:completion} uses.

\emph{Smooth manifolds.} For Riemannian manifolds with $\Ric\ge-(n-1)$,
\crefrange{rf:app:flow}{rf:sec:example} follow a different route, through Ricci flow: small
local rewinding Reifenberg error gives a flow for a uniform time with $|\Rm|\le\alpha/t$, for
a prescribed $\alpha\in(0,1)$, on closed and on complete manifolds, and for collapsing closed
manifolds the bounded-curvature theory at a positive time gives singular fibrations with
affine structure group (\cref{rf:thm:main,rf:thm:noncompact,rf:thm:rong}). Conditional on
\cite{Nif26}, \cref{rf:thm:example} shows that non-collapsing of the local covers cannot
replace the small-error hypothesis. These three appendices contain, with corrections, the
previous arXiv version of this paper; the main text does not use them.

\subsection{The Singular Fibration Theorem}\label{sec:answer}\label{es:sec:results}

The main theorem resolves the isotropy of the limit by canonical frame spaces and describes
the spaces of the sequence as singular fibrations over the limit.

\begin{theorem}[Singular Fibration]\label{thm:main}
There are constants $\delta_1,c_1>0$, depending only on $n,N$, with the following
property. Assume \cref{hyp:sequence} with $\delta\le\delta_1$, and fix
$0<r\le c_1\rho$. After passing to a subsequence:
\begin{enumerate}[label=\textup{(\roman*)}]
\item \emph{Frame resolution.} The canonical frame spaces at this fixed scale,
$Y_i=\Fr_\rho(X_i)$, are compact principal $O(n)$-bundles over $X_i$.
They converge in a compact metrizable realization to a compact topological
manifold $Y$ with a locally linear $O(n)$-action and $Y/O(n)=X$.

\item \emph{Equivariant bundle.} There exist equivariant realizing maps
$F_i\colon Y_i\to Y$ (\cref{fl:def:realizing}), that is, continuous equivariant maps such that
$F_i(y_i)\to y$ whenever $y_i\to y\in Y$.
For any such maps and every $\eps>0$, all sufficiently large $i$ admit an
equivariant surjective fibre bundle $F_i^\natural\colon Y_i\to Y$ with closed
manifold fibres and $\sup_{Y_i}\sd(F_i^\natural,F_i)<\eps$, where $\sd$ is a fixed
compatible invariant metric on $Y$. The induced maps
$f_i^\natural,f_i\colon X_i\to X$ satisfy
$\sup_{X_i}d(f_i^\natural,f_i)<\eps$.

\item \emph{Local quotient form.} The bundles in \textup{(ii)} can be chosen such
that, for every large $i$, every $x\in X$ has a neighbourhood $O_x$ homeomorphic to
$B/I$, where $I$ is the isotropy group of a point of $Y$ over $x$ and $B$ is a ball
about $0$ in an orthogonal representation of $I$, with the following property. The
homeomorphism sends $x$ to $I0$, and there is a closed manifold $\mathsf F$ with a free
$I$-action such that the restriction of $f^\natural_i$ over $O_x$ is equivalent to
\[
\mathsf F\times_I B\longrightarrow B/I,\qquad [u,v]\longmapsto Iv,
\]
and the restriction of $F^\natural_i$ over the preimage of $O_x$ is equivariantly
equivalent to \mbox{$O(n)\times_I(\mathsf F\times B)\to O(n)\times_I B$}. The local data may
depend on $i$, and
\[
\dim B=n-\dim\mathsf F+\dim I,\qquad\dim(\mathsf F/I_v)=\dim\mathsf F-\dim I_v,
\]
where, for $v\in B$, $I_v$ is the isotropy group of $v$ in $I$ and $\mathsf F/I_v$ is the fibre of
$f^\natural_i$ at $Iv$.

\item \emph{Fibre identification.} Along a further common subsequence the maps
in \textup{(ii)} can be chosen so that every fibre of $F^\natural_i$, and every $\mathsf F$
in \textup{(iii)}, is homeomorphic to the nilmanifold of the fibre lattice of
a chart over the point. Every fibre of $f_i^\natural$ is homeomorphic to a
closed infranilmanifold.
\end{enumerate}
\end{theorem}

The frame spaces are constructed in \cref{ch:canonicalframes,ch:coveratlas}. The charts are those
of \cref{m:fr:def:chart}, the convergence of points in \textup{(ii)} is that of
\cref{m:fr:def:Yconv}, locally linear actions are defined in \cref{ls:def:locallylinear}, and
the fibre lattices of the charts are those of \cref{m:fr:ext:physical-lattice-realization},
taken along the common subsequence of \cref{fl:sec:malcev}. The constants may be chosen as
\[
\delta_1=\min\{\delta_*,\eps_{\mathrm{glue}}(n,N)/3\},\qquad c_1=c_*,
\]
where $\delta_*,c_*$ are the covering constants of \cref{m:fr:sec:coverbundle} and
$\eps_{\mathrm{glue}}$ is the final gluing threshold of \cref{es:sec:finitesetting}.
All depend only on $n,N$; neighbourhoods and eventual index bounds may also
depend on the sequence, frame scale and approximation tolerance.
The proof is completed in \cref{es:sec:proofsintro}.

\emph{Application: local bounded covering geometry.}
The completion of local covers also extends \cite[Theorem~1.7]{HZ26} to the local covering
setting of \cite[Remark~5.4]{HZ26}, without assuming that the spaces are smooth.

\begin{definition}[Local Bounded Covering Geometry]\label{hz:def:lbcg}
Let $\rho,\iota_0>0$. An $\RCD(\kappa,N)$ space has \emph{$(\rho,\iota_0)$-local bounded
covering geometry} if, for every $x\in X$, the universal cover of $B_\rho(x)$ with its lifted
length metric satisfies $\operatorname{Injrad}(q)\ge\iota_0$ for every
$q\in B_{\rho/2}(\tilde x)$. Here $\tilde x$ is a lift of $x$, and we use the metric
cut-locus definition of injectivity radius from \cite[Definition~2.24]{HZ26}, recalled in
\cref{m:fr:sec:completion}.
\end{definition}

This is $(\rho/2,\iota_0)$-locally bounded covering geometry in the notation of
\cite[Remark~5.4]{HZ26}. No Reifenberg or curvature-scale hypothesis is imposed here.

\begin{theorem}[Fibration under Local Bounded Covering Geometry]\label{hz:thm:fibration}
Let $\kappa\in\R$, $N\in[1,\infty)$ and $\rho,\iota_0>0$. Let $(X_i,d_i,m_i)$ be compact
$\RCD(\kappa,N)$ spaces with $m_i(X_i)=1$ and $(\rho,\iota_0)$-local bounded covering
geometry, converging in the measured Gromov--Hausdorff sense to a compact space $(X,d,m)$.
Then for all large $i$ there is a possibly singular fibration $X_i\to X$ with fibres
diffeomorphic to infranilmanifolds, induced by an $O(n)$-equivariant fibration from their
frame bundles; here $n$ is the essential dimension of $X_i$. If $X$ is smooth in the
metric-measure sense of \cite[Section~2.2]{HZ26}, the map is a smooth fibre bundle.
\end{theorem}
\begin{proof}
We may increase $N$ to $\max\{N,2\}$ and split the sequence according to its finitely many
possible essential dimensions. By \cref{hz:lem:scale}, the spaces have harmonic Riemannian structures
and a uniform positive covering regularity scale. These replace the two uses of global
bounded covering geometry in the proofs of \cite[Theorem~5.3 and Theorem~1.7]{HZ26}: the
local regularity of the metric and density, and the uniform scale bound in \cite[(5.6)]{HZ26}.
The smoothing construction and the subsequent frame-bundle argument therefore give the
assertion, with constants depending also on $\rho$.
\end{proof}

\emph{Application: limits of aspherical spaces.}
Rong asked when a collapsed limit of aspherical manifolds is an aspherical manifold, and
conjectured an answer for fundamental groups of the following type \cite[Conjecture~0.7]{Ron26}.

\begin{definition}[Type $A_s$]\label{as:def:As}
A finitely generated group $\Pi$ is of \emph{type $A_s$} if, for every normal virtually
nilpotent subgroup $\Lambda'$ of $\Pi$, every finite subgroup of $\Pi/\Lambda'$ has a normalizer
of finite index in $\Pi/\Lambda'$ \cite[Definition~0.4]{Ron26}.
\end{definition}

The fundamental group of a nilmanifold is of type $A_s$, and so is that of its product with a
closed hyperbolic manifold; that of the flat Klein bottle, which collapses to an interval, is
not. The singular fibration of \cref{thm:main} gives the topological part of Rong's conjecture
in the setting of \cref{hyp:sequence}, and more.

\begin{theorem}[Limits of Aspherical Spaces]\label{as:thm:aspherical}
In the setting of \cref{thm:main}, let $i$ be large enough for parts \textup{(ii)} to
\textup{(iv)}, and suppose that $X_i$ is aspherical. Then the following hold.
\begin{enumerate}[label=\textup{(\roman*)}]
\item Every isotropy group of the action of $O(n)$ on $Y$ is finite. The fundamental groups of
the fibres of $F^\natural_i$ and of $f^\natural_i$ inject into $\pi_1(X_i)$.
\item The space $X$ is homeomorphic to the quotient of a contractible topological manifold by a
properly discontinuous, cocompact action of a discrete group with finite stabilizers, and the
universal cover of $X_i$ is a fibre bundle over this manifold whose fibres are homeomorphic to
Euclidean spaces.
\item If $\pi_1(X_i)$ is of type $A_s$ \textup{(}\cref{as:def:As}\textup{)}, then $X$ is a closed aspherical
topological manifold.
\end{enumerate}
\end{theorem}

Part \textup{(iii)} is the topological conclusion of Rong's conjecture in this setting; parts
\textup{(i)} and \textup{(ii)} need no condition on the fundamental group. The proof is in
\cref{as:sec:aspherical}.

\subsection{Proof Strategy}\label{sec:chain}

The proof of \cref{thm:main} has four steps, described below with the sections that carry
them out; the last paragraph concerns the two applications.
\Cref{guide:sec:dependencies} gives linked dependency diagrams, the order of the
parameter choices, and a register distinguishing uniform constants from local scales.

\emph{Compatible frames (\cref{ch:canonicalframes,ch:frameconvergence,ch:coveratlas}).}
The frame operator, whose kernel is formed with the distance function alone, defines
canonical frames. Gauges trivialize the frame
bundle locally and control the convergence of frames, and closeness of balls to Euclidean balls at all small
scales makes the limit groups act freely on them. Agreement of neighbouring covers then gives a finite atlas and
the framed limit. \Cref{fl:sec:setting} records the resulting geometric data.

\emph{Coordinates, slices and realizing maps (\cref{ch:harmoniccomparison}).}
Injective harmonic charts admit comparison errors with small Lipschitz constants.
The canonical frame supplies an exact orthogonal isotropy representation; averaging
a harmonic comparison map makes it equivariant for this representation. Smoothing
in the group direction then gives linear slices as transverse level sets. The slices make the framed limit
a manifold with a locally linear $O(n)$-action, and an equivariant neighbourhood
retraction supplies realizing maps. Quantitative versions of the slices, with
equivariance errors quadratic in the distance and constants uniform under
recentring, prepare the construction on the approximating charts.

\emph{The bundle (\cref{ch:stages,ch:localbundle,ch:globalbundle}).}
The limit coordinates are transplanted to the approximating charts, where deck
transformations with small displacement receive isotropy labels (\cref{ch:stages}).
These labels are corrected to an exact cocycle, and twisted deck averages give an
equivariant coordinate, with limits given by Haar averages (\cref{ch:localbundle}). Convergent sections of the frame spaces realize the cocycle. Regular coordinates of such
sections give product charts, and their quantitative stability under coordinate changes
allows gluing over a finite system of tubes with bounded overlap; equivariant triviality
over balls then gives the local product form of the glued bundle.

\emph{The fibres (\cref{ch:lattices,ch:fibregroups}).}
\Cref{ch:lattices} constructs marked Malcev lattices, controls their period kernels and builds unwrapped coverings. \Cref{ch:fibregroups} fills spheres in the contractible nilpotent model and transfers the fillings through nested contractible sheets, identifying the intrinsic fibre groups and proving asphericity. Rigidity shows that the fibres of the bundles $Y_i\to Y$ are nilmanifolds. A principal-bundle homotopy argument gives asphericity of their isotropy quotients, and rigidity identifies these fibres of $X_i\to X$ as infranilmanifolds.

The separate application \cref{hz:thm:fibration} uses only the completion and its regularity
consequence \cref{hz:lem:scale}, both proved in \cref{m:fr:sec:completion}. The second
application, \cref{as:thm:aspherical}, uses \cref{thm:main} and the limit groups of the charts.
Its proof in \cref{as:sec:aspherical} lifts the bundle to a covering of $Y$ whose Borel
construction is an Eilenberg--MacLane space. A torus bundle that realizes the kernel of the
homomorphism from the fibre group to $\pi_1(X_i)$ has weakly contractible Borel construction,
so $O(n)$ acts freely on it, and its quotient is a contractible manifold. The fixed points of the
torus, Poincar\'e duality, and the triviality of isometries of $\RCD$ spaces that fix an open set
force this kernel to vanish. For type $A_s$, Smith theory then shows that the elements of the
deck group of the lifted quotient that have fixed points act trivially.

\subsection{Scope and Inputs}\label{sec:quoted}\label{sec:scope}

We delimit the conclusions, record the steps that hold more generally, and list the external
inputs. The conclusion of \cref{thm:main} is topological. It does not assert an affine structure group for
the global bundle. Finite isotropy actions on the nilmanifold fibres can be made
affine, but no such conclusion is claimed for positive-dimensional isotropy;
see \cref{es:rem:affine,es:rem:affinebundles}. The frame limit depends on the
subsequence, and the bundle maps involve choices (\cref{m:fr:ex:interleave}). The conclusion
of \cref{as:thm:aspherical} is also topological: it does not assert the quantitative Reifenberg
regularity of the limit metric requested in Rong's conjecture, and we do not prove
contractibility of the metric limit of the universal covers \cite[Conjecture~0.13]{Ron26}.

Several steps hold in more generality than the main theorem needs. The frame operator and the
frame bundle of \cref{m:fr:thm:A,m:fr:thm:bundle} are defined on every separable metric space
with a Borel measure that is finite on bounded sets and positive on balls, without any
curvature or dimension bound. The completion of local covers in \cref{m:fr:prop:completion}
holds on every $\RCD(\kappa,N)$ space with $N>1$, at every radius. The
linear slices of \cref{ls:thm:slices} hold for proper isometric actions of Lie groups that
preserve the measure, on invariant Reifenberg-flat open subsets of $\RCD$ spaces.

The principal external inputs of the main theorem are $\RCD$ compactness and harmonic calculus
\cite{GMS15,EKS15,Gig15,Jia14,AH18}; Reifenberg and harmonic-coordinate theory
\cite{CC97,GV25,HP23,BPS23}; equivariant convergence, covering, isometry-group and
Margulis theory \cite{FY92,MW19,Wan24b,Sos18,GSR19,Wan23}; nilprogressions and nilpotent
lattices \cite{Zam24,CG90,Dek17}; and slices, equivariant retractions, bundle
theory and topological rigidity
\cite{Pal61,Bre72,Mur83,Hus94,Las82,Sie72,KS77,FH83,FQ90,Luc10}. The aspherical application also
uses Borel constructions, equivariant torus bundles and Poincar\'e duality groups
\cite{Bor60,Hsi75,LMS83,May90,Kha18,Bro82}.

\section{Geometric and Analytic Preliminaries}\label{ch:preliminaries}

We fix the conventions, state the results on $\RCD$ spaces that we quote, and define
the forms of closeness to Euclidean balls used below, together with their
consequences. The last two subsections contain lemmas on rigid motions, harmonic
maps and group actions. The setting for the main $\RCD$ construction is \cref{hyp:sequence};
\cref{fl:sec:setting} records the atlas derived from it.

\subsection{Conventions}\label{m:book:ch:02}\label{m:fr:sec:prelim}\label{m:fr:sec:conv}

This subsection fixes notation and conventions.

The symbol $\Psi(\eps,\delta,\eta,\dots\,|\,a,b,\dots)$ denotes a non-negative function that
tends to $0$ as $\eps,\delta,\eta,\ldots\to0$, for fixed $a,b,\dots$, and $C(a,b,\dots)$ denotes
a positive constant depending only on $a,b,\dots$; the lists of small quantities and of
parameters are finite, and both symbols may change from line to line. Balls are open and
carry the restricted metric, and $B_s(0^n)$ is the ball in $\R^n$. A map
$f\colon A\to A'$ between metric spaces is an \emph{$\eps$-Gromov--Hausdorff
approximation}, or $\eps$-GH approximation, if $|d'(f(a),f(b))-d(a,b)|\le\eps$ for all
$a,b\in A$ and every point of $A'$ is within distance $\eps$ of $f(A)$. A \emph{metric measure space}
$(Z,d,m)$ is a separable metric space with a Borel measure that is finite on bounded sets and
positive on balls.
Averages $\fint$ are formed with the measure of the space in question, and $\omega_n$ is
the Lebesgue measure of $B_1(0^n)$. For a vector $v\in\R^n$, $|v|$ is the Euclidean norm.
For a matrix $M$, $\|M\|$ is the operator norm, $|M|$ the largest absolute value of an entry, and $|M|_{\mathrm{HS}}$ the Hilbert--Schmidt norm. We write $\one$ for the identity matrix, also when it is regarded as the identity element of $O(n)$; the identity element of every other group is $1$, and $\mathrm{id}$ is the identity map.

Isometries are surjective, $\Isom(Z)$ carries the compact-open topology, and actions of
topological groups are continuous. For a group
$G$ of isometries, $G_x$ is the isotropy group of $x$, and $Z/G$ carries the quotient
distance ${d(Gx,Gy)=\inf_{g\in G}d(x,gy)}$. We use pointed equivariant Gromov--Hausdorff
convergence ${(Z_i,z_i,\Gamma_i)\eGHto(Z,z,G)}$ in the sense of \cite{FY92}: there are
${\eps_i\to0}$ and maps ${\af_i\colon B_{1/\eps_i}(z_i)\to Z}$,
${\alpha_i\colon\Gamma_i(1/\eps_i)\to G(1/\eps_i)}$ and
${\beta_i\colon G(1/\eps_i)\to\Gamma_i(1/\eps_i)}$, where
$G(R):=\{g\in G:d(gz,z)<R\}$ and $\Gamma_i(R)$ is defined in the same way, such that
$\af_i(z_i)=z$, the map $\af_i$ distorts distances by at most $\eps_i$ and its image is
$\eps_i$-dense in $B_{1/\eps_i}(z)$, and
$d(\af_i(\gamma p),\alpha_i(\gamma)\af_i(p))\le\eps_i$ and $d(\af_i(\beta_i(g)p),g\af_i(p))\le\eps_i$ whenever the points involved lie in the balls of
radius $1/\eps_i$. Limits of closed subgroups exist after passing to subsequences, and
quotients converge: if closed normal subgroups $\Gamma_i'\lhd\Gamma_i$ converge to
$G'\lhd G$, then $(Z_i/\Gamma_i',\Gamma_i/\Gamma_i')\eGHto(Z/G',G/G')$
\cite{FY92}, \cite[Theorem~2.1 and Lemma~3.1]{Wan23}. Here $\Gamma_i/\Gamma_i'$ and $G/G'$ stand for
their images in the isometry groups of the quotient spaces, because their actions need not be
effective: if $Z=\R$, $G=\Isom(\R)$ and $G'$ is the group of translations, then $Z/G'$ is a
point and $G/G'\cong\mathbb Z/2$.

A sequence of elements $\gamma_i\in\Gamma_i$ \emph{converges in the eGH sense} to $g\in G$,
written $\gamma_i\eGHto g$, if the displacements
$d_i(\gamma_iz_i,z_i)$ are bounded and $d(\af_i(\gamma_iq),g\,\af_i(q))\to0$ uniformly for
$q\in B_R(z_i)$, for every fixed $R>0$. We use the same notion for sequences indexed by
an infinite set of indices. Convergence within a fixed group $G$ retains the unlabelled
arrow $g_j\to g$.

\subsection{RCD Spaces}\label{m:fr:sec:rcd}

We use the following facts about an $\RCD(\kappa,N)$ space $(Z,d,m)$ with $N<\infty$; nothing in
this subsection is new. The space is \emph{non-collapsed} if $N$ is an integer and
$m=\cH^N$ \cite{DPG18}. We write $\kappa$ for the curvature lower bound, to distinguish it from compact subgroups $K\subset G$. We write $\mathbb{D}$ for a distance that metrizes the pointed
measured Gromov--Hausdorff convergence of pointed $\RCD(-1,N)$ spaces whose measures are
normalized at the base point, either by $m(B_1(q))=1$ or by
$\int_{B_1(q)}(1-d(q,\cdot))\,dm=1$, and we call two such spaces \emph{$\eps$-close for
$\mathbb{D}$} if their distance is less than $\eps$. The class of these spaces is compact
by \ref{m:fr:R:geometry}, so two such distances are uniformly equivalent, and the
statements below do not depend on the choice.

\begin{enumerate}[label=\textup{(R\arabic*)},ref=\textup{(R\arabic*)},leftmargin=3em]
\item\label{m:fr:R:geometry} \emph{Geometry.} $Z$ is proper and geodesic. If $|\kappa|s^2\le1$,
the Bishop--Gromov inequality gives $m(B_{2t}(z))\le C_D\,m(B_t(z))$ for $t\le s$, with
$C_D=C_D(N)$, spheres are $m$-null, and a $(1,2)$-Poincar\'e inequality holds on balls of
radius at most $s$. The pointed $\RCD(-1,N)$ spaces whose measures are normalized at the
base point form a compact class for the pointed measured Gromov--Hausdorff convergence, and
limits of $\RCD(\kappa_j,N)$ spaces with $\kappa_j\to0$ are $\RCD(0,N)$ spaces \cite{GMS15}. An
$\RCD(0,N)$ space that contains a line is isomorphic, as a metric measure space, to the
product of the line with Lebesgue measure and an $\RCD(0,N-1)$ space, or a point
\cite{Gig13}.
\item\label{m:fr:R:isometries} \emph{Isometries.} $\Isom(Z)$ is a Lie group
\cite{Sos18,GSR19}; its closed subgroups act properly. An isometry that fixes a non-empty open set
pointwise is the identity. Indeed, $Z$ is non-branching \cite{Den25}; if $\varphi$ fixes
$B_s(x)$ pointwise and $\gamma$ is a geodesic from $x$ to $y$, then $\varphi\circ\gamma$
is a geodesic of the same length that agrees with $\gamma$ near $x$, hence equals
$\gamma$, and $\varphi(y)=y$.
\item\label{m:fr:R:calculus} \emph{Calculus.} For an open set $U\subset Z$, $\Hzero(U)$ is the
closure of the Lipschitz functions with compact support in $U$ in $\Hone(Z)$, and
$\widehat H^{1,2}_0(U)$ is the set of $f\in\Hone(Z)$ that vanish almost everywhere
outside $U$. A function $f\in\Hone_{\mathrm{loc}}(U)$ is in the domain of the local
Laplacian, with $\Delta f=g\in L^2_{\mathrm{loc}}(U)$, if
$\int\langle\nabla f,\nabla\psi\rangle\,dm=-\int g\psi\,dm$ for all Lipschitz
$\psi$ with compact support in $U$; it is \emph{harmonic} if $g=0$, and a map to $\R^j$
is harmonic if its components are. Linear combinations of harmonic functions are harmonic,
and if $\varphi$ is an isometry of $Z$ that preserves $m$ and $f$ is harmonic on $U$, then
$f\circ\varphi$ is harmonic on $\varphi^{-1}(U)$. For a locally Lipschitz map $f$ from an
open subset of $Z$ to $\R^n$ we write
$\Gr(f):=\bigl(\langle\nabla f_a,\nabla f_b\rangle\bigr)_{a,b=1}^n$ for its Gram matrix.
The Leibniz rule \cite[Theorem~4.29]{Gig15} gives, by induction on the degree, the chain rule
$\Delta(\Theta\circ f)=\sum_a\partial_a\Theta(f)\Delta f_a+\sum_{a,b}\partial_{ab}\Theta(f)\,\Gr(f)_{ab}$
for polynomials $\Theta$ and locally Lipschitz maps $f$ with locally bounded Laplacian;
approximation of $\Theta$ in $C^2$ on bounded sets and the stability of the Laplacian
\cite[Proposition~4.15]{Gig15} extend it to every $\Theta\in C^2$. For a harmonic Lipschitz map $f$ the right side is
$\sum_{a,b}\partial_{ab}\Theta(f)\,\Gr(f)_{ab}$, which is bounded on every set on which $f$
is bounded. A Lipschitz function $f$ satisfies $|\nabla f|=\operatorname{lip}f$ almost
everywhere \cite{Che99}; in particular $|\nabla d(x,\cdot)|=1$ almost everywhere.
\item\label{m:fr:R:estimates} \emph{Gradient estimate.} Let $|\kappa|s^2\le1$. A harmonic function $f$ on
$B_{2s}(z)$ has a locally Lipschitz representative, it satisfies
\begin{equation*}
\sup_{B_s(z)}|\nabla f|\le\frac{C_J}s\fint_{B_{2s}(z)}|f|\,dm,\qquad C_J=C_J(N),
\end{equation*}
by \cite[Theorem~1.1]{Jia14}, and it is Lipschitz on $B_{s/2}(z)$ with the constant
$\sup_{B_s(z)}|\nabla f|$.
\item\label{m:fr:R:splitting} \emph{Splitting maps.} For every $\delta>0$ there is
$\eps_B=\eps_B(N,\delta)>0$ with the following property. Let $q\in Z$ and $r>0$ with
$|\kappa|r^2\le\eps_B$, and suppose that $(Z,r^{-1}d,m,q)$, with the measure normalized at
$q$, is $\eps_B$-close for $\mathbb{D}$ to $\R^n$ with a multiple of Lebesgue measure.
Then there is a harmonic map $v\colon B_{5r}(q)\to\R^n$ whose components are Lipschitz
with a constant $C_N$ and which satisfies $\fint_{B_{5r}(q)}|\Gr(v)-\one|\,dm<\delta$
\cite[Definition~3.4, Corollary~3.10]{BPS23}.
\item\label{m:fr:R:almostisometry} \emph{Almost isometry.} There is a function
$\Psi_G(\eps|n,N,L)$ that tends to $0$ as $\eps\to0$, with the following property.
Suppose that $Z$ has the essential dimension $n$, let $x\in Z$ and $s>0$ with
$|\kappa|s^2\le\eps$, and let $f\colon B_{6ns}(x)\to\R^n$ be a harmonic map with
$|\nabla f_a|\le L$ almost everywhere and $\fint_{B_{6ns}(x)}|\Gr(f)-\one|\,dm\le\eps$. Then
$f$ changes the distances of the points of $B_s(x)$ by at most $\Psi_Gs$, and every
point of $B_s(f(x))$ has distance at most $\Psi_Gs$ from $f(B_s(x))$
\cite[Proposition~4.2]{HP23}.
\end{enumerate}

The \emph{essential dimension} is the integer $k$ such that the tangent cone at
$m$-almost every point is $\R^k$; it exists by the theorem of Bru\`e and Semola, as
recalled in \cite{HP23}. The reference measures in \ref{m:fr:R:splitting} and
\ref{m:fr:R:almostisometry} are arbitrary.

\subsection{Closeness to Euclidean Space}\label{m:fr:sec:closeness}\label{m:fr:sec:tools}
The Reifenberg conditions used below compare balls with balls of the same radius in a Euclidean
space of one fixed dimension. Fix an integer $n\ge1$, the dimension of this model space; in the
application it is the integer $n$ of \cref{m:fr:def:localcovers}. The first form asks for
closeness to Euclidean balls about all points of a set and at all scales up to a bound. It is
the form of the hypothesis on the local covers (\cref{m:fr:def:localcovers}) and of the
standing hypothesis on the regions of the harmonic charts (\cref{ch:harmoniccomparison}).

\begin{definition}[Close Balls, Reifenberg Sets and Points]\label{m:fr:def:reifenberg}
A ball $B_s(z)$ of a metric space is \emph{$\eps s$-close} to $B_s(0^n)$ if there is an
$\eps s$-GH approximation $B_s(z)\to B_s(0^n)$. Let $\eps,\lambda>0$. A subset $\Omega$ of a metric space $Z$ is
\emph{$(\eps,\lambda)$-Reifenberg} if the ambient ball $B_t(q)$ is $\eps t$-close to
$B_t(0^n)$ for every $q\in\Omega$ and $0<t\le\lambda$. A point $q\in Z$ is
\emph{$(\eps,\lambda)$-Reifenberg} if $\{q\}$ is. This is a metric condition;
any curvature normalization is stated separately.
\end{definition}

The other two forms concern a single ball together with the measure. The frames of
\cref{ch:canonicalframes} need an approximation of one ball that sends the centre to the
origin, and a lower bound for the measures of its smaller balls. The third form is the
one that remains available when the measure is not a Hausdorff measure, as on the
completed local covers of \cref{ch:coveratlas}.

\begin{definition}[Almost Euclidean Points]\label{m:fr:def:euclid}
Let $\eps>0$ and $r>0$. A point $x$ of a metric measure space $(Z,d,m)$ is
\emph{$\eps$-Euclidean at scale $r$} if there is a Borel $\eps r$-GH approximation
$\Phi\colon B_{4r}(x)\to B_{4r}(0^n)$ with $\Phi(x)=0$, and
$m(B_{r/8}(p))\ge32^{-n}\,m(B_{2r}(x))$ for every $p\in B_r(x)$.
\end{definition}

For the Lebesgue measure on $\R^n$ the ratio $m(B_{r/8}(p))/m(B_{2r}(x))$ is $16^{-n}$. The
smaller constant $32^{-n}$ leaves room for the errors in the application of
\cref{m:fr:lem:flat}. The third form compares the normalized measure of the ball, transported
by the approximation, with the normalized Lebesgue measure, and it adds a curvature
normalization at the scale of the ball.

\begin{definition}[Flat Points and Regions]\label{m:fr:def:flat}
Let $(Z,d,m)$ be an $\RCD(\kappa,N)$ space with $N<\infty$, and let $n\ge1$, $r>0$ and
$\eps>0$. A point $x\in Z$ is \emph{$\eps$-flat at scale $r$} if $|\kappa|r^2\le\eps$ and there
is a Borel $\eps r$-GH approximation $\Phi\colon B_{4r}(x)\to B_{4r}(0^n)$ with $\Phi(x)=0$
such that
\begin{equation*}
\Bigl|\fint_{B_{4r}(x)}g\circ\Phi\,dm-\fint_{B_{4r}(0^n)}g\,d\cH^n\Bigr|\le\eps
\end{equation*}
for every function $g$ on $B_{4r}(0^n)$ with $|g|\le1$ and $\Lip g\le1/r$, where the
second average is formed with the Lebesgue measure. An open subset $\Omega$ of $Z$ is \emph{$\eps$-flat
at scale $r$} if each of its points is. For $s>0$ we write $\Omega^{[s]}$ for the open
set of the points of $\Omega$ whose distance from $Z\setminus\Omega$ is larger than $s$,
and $\Omega^{[s]}:=Z$ if $\Omega=Z$.
\end{definition}

We now record the consequences of these definitions that are used later. Harmonic maps
that almost preserve distances take values in $\R^n$ rather than in a ball, and they are
compared with Euclidean balls through the following notion.

\begin{definition}[Centred Approximations]\label{m:fr:def:centred}
Let $B_s(z)$ be a ball of a metric space and $\eps>0$. A map $\Phi\colon B_s(z)\to\R^n$ is a
\emph{centred $\eps s$-GH approximation} of $B_s(z)$ if $\Phi(z)=0$,
\begin{equation*}
\bigl|\,|\Phi(p)-\Phi(q)|-d(p,q)\bigr|\le\eps s\qquad\text{for all }p,q\in B_s(z),
\end{equation*}
and every point of $B_s(0^n)$ has distance at most $\eps s$ from $\Phi(B_s(z))$.
\end{definition}

The approximations in \cref{m:fr:def:euclid,m:fr:def:flat} are centred approximations of
$B_{4r}(x)$ with values in $B_{4r}(0^n)$. The next lemma produces such
approximations from closeness.

\begin{lemma}[Centring]\label{m:fr:lem:centred}
Let $B_s(z)$ be a ball of a separable metric space and $0<\eps\le\tfrac1{100}$. Every
$\eps s$-GH approximation $f\colon B_s(z)\to B_s(0^n)$ satisfies $|f(z)|\le2\eps s$. If
$B_s(z)$ is $\eps s$-close to $B_s(0^n)$, there is a Borel $6\eps s$-GH approximation
$B_s(z)\to B_s(0^n)$ that takes $z$ to $0$.
\end{lemma}
\begin{proof}
Let $f$ be an $\eps s$-GH approximation with $f(z)\ne0$, and let $0<t<1$ and
$\sigma>\eps s$. The point $y:=-(1-t)s\,f(z)/|f(z)|$ lies in $B_s(0^n)$, so some $p$
satisfies $|f(p)-y|<\sigma$. Then
\begin{equation*}
s>d(p,z)\ge|f(p)-f(z)|-\eps s>|y-f(z)|-\sigma-\eps s=(1-t)s+|f(z)|-\sigma-\eps s ,
\end{equation*}
so $|f(z)|<t s+\sigma+\eps s$. Letting $t\to0$ and $\sigma\to\eps s$ gives
$|f(z)|\le2\eps s$.

Now let $f$ be an $\eps s$-GH approximation $B_s(z)\to B_s(0^n)$, let $\pi$ be the
nearest-point retraction of $\R^n$ onto the closed ball $\bar B_s(0^n)$, and put
$f_1:=\pi\circ(f-f(z))$. Since $|f(p)-f(z)|<s+2\eps s$ and $\pi$ moves such points by less
than $2\eps s$, the map $f_1$ takes $z$ to $0$ and changes distances by at most $5\eps s$.
For $y\in B_s(0^n)$ the point $(1-2\eps)y+f(z)$ lies in $B_s(0^n)$, so some $p$ satisfies
$|f(p)-(1-2\eps)y-f(z)|<\tfrac32\eps s$, and then $|f_1(p)-y|<4\eps s$, because $\pi$ is
$1$-Lipschitz and fixes $y$. Choose points $y_j$, $j\ge1$, that are dense in $B_s(0^n)$,
and points $p_j$ with $|f_1(p_j)-y_j|<4\eps s$. Since the space is separable, $B_s(z)$ is the
disjoint union of countably many Borel sets of diameter at most $\tfrac18\eps s$, among
them the sets $\{z\}$ and $\{p_j\}$. Choose a point in each set, the point itself in the
case of $\{z\}$ and $\{p_j\}$, and let $f_2$ be $(1-\tfrac18\eps)$ times the value of $f_1$
at the chosen point. Then $f_2$ is Borel, it takes values in $B_s(0^n)$ and $z$ to $0$, and
it changes distances by at most $5\eps s+\tfrac14\eps s+\tfrac14\eps s<6\eps s$, since the
values of $f_1$ lie in $\bar B_s(0^n)$. Every $y_j$ has distance less than
$4\eps s+\tfrac18\eps s$ from $f_2(p_j)$, so every point of $B_s(0^n)$ has distance at most
$5\eps s$ from the image of $f_2$.
\end{proof}

The Gromov--Hausdorff closeness of a ball implies closeness of the measures at smaller
scales, in every $\RCD$ space of finite dimension. The reason is that Euclidean space
carries only one $\RCD(0,N)$ measure up to a factor. We state this for balls, in the form
of \cref{m:fr:def:flat}, and for pointed spaces, in the form used by the harmonic-coordinate
estimates of \cite{HP23}.

\begin{lemma}[Measured Closeness]\label{ls:lem:measured}\label{m:fr:lem:flatness}
For every $\eps>0$ there are $\delta_H=\delta_H(n,N,\eps)>0$ and
$c_H=c_H(n,N,\eps)\in(0,1)$ with the following property. Let $(Z,d,m)$ be an
$\RCD(\kappa,N)$ space, let $q\in Z$ and $s>0$ with $|\kappa|s^2\le1$, and suppose that $B_s(q)$
is $\delta_Hs$-close to $B_s(0^n)$. Then for every $r$ with $\tfrac12c_Hs\le r\le c_Hs$:
\begin{enumerate}[label=\textup{(\roman*)}]
\item the point $q$ is $\eps$-flat at scale $r$;
\item the space $(Z,r^{-1}d,m,q)$, with the measure normalized at $q$ in either of the two
ways of \cref{m:fr:sec:rcd}, is $\eps$-close for $\mathbb{D}$ to $\R^n$ with the multiple of
Lebesgue measure that has the same normalization at $0$.
\end{enumerate}
In particular, put $\eps_V(n,N,\eps):=\min\{\delta_H,c_H^2\}$. If $|\kappa|r^2\le\eps_V$ and
$B_t(q)$ is $\eps_Vt$-close to $B_t(0^n)$ for all $t\in[r,r/\eps_V]$, then \textup{(i)} and
\textup{(ii)} hold at the scale $r$.
\end{lemma}
\begin{proof}
Suppose that, for some $\eps$, no choice of $\delta_H$ and $c_H$ has the property. Then for
every $j$ there are an
$\RCD(\kappa_j,N)$ space $(Z_j,d_j,m_j)$, a point $q_j$ and numbers $s_j$ and $r_j$ that satisfy
the hypotheses with $c_H=1/j$ and $\delta_H=1/j^2$, such that \textup{(i)} or \textup{(ii)}
fails at the scale $r_j$. We multiply the distances by $1/r_j$ and the measures by
constants that normalize them at $q_j$, in the way for which \textup{(ii)} fails if it does.
The new spaces are $\RCD(\kappa'_j,N)$ spaces with $|\kappa'_j|\le j^{-2}$, and the ball of radius
$R_j:=s_j/r_j\in[j,2j]$ about $q_j$ is $\tfrac2j$-close to the Euclidean ball of this
radius. An approximation that realizes this closeness moves $q_j$ by at most $\tfrac4j$
(\cref{m:fr:lem:centred}), so its restrictions show that for every $R>0$ the
closed balls of radius $R$ about $q_j$ converge to the closed ball of radius $R$ in $\R^n$
in the Gromov--Hausdorff sense. By \ref{m:fr:R:geometry} a subsequence converges in the
pointed measured Gromov--Hausdorff sense to an $\RCD(0,N)$ space
$(Z_\infty,d_\infty,m_\infty,q_\infty)$ with the same normalization, and its closed balls
about $q_\infty$ are the limits of those about $q_j$. So $\bar B_R(q_\infty)$ is isometric
to the closed ball of radius $R$ in $\R^n$ for every $R$. An isometry takes $q_\infty$ to
$0$, the only point with distance at most $R$ from all others. By the theorem of Arzel\`a
and Ascoli, a sequence of such isometries with $R\to\infty$ has a subsequence that
converges uniformly on compact sets, and the limit is an isometry of $Z_\infty$ onto
$\R^n$. Thus $m_\infty$ is a measure on $\R^n$ for which the Euclidean distance gives an
$\RCD(0,N)$ space. By the splitting theorem in \ref{m:fr:R:geometry}, applied to the lines
in the coordinate directions one after the other, $m_\infty$ is a multiple of Lebesgue
measure, and the normalization at $0$ determines the multiple. So the rescaled pointed
spaces converge to the model, and \textup{(ii)} holds at the scale $r_j$ for large $j$.

The convergence also gives Borel maps of the balls of radius $5$ about $q_j$ to $\R^n$ that
send $q_j$ to $0$, change distances by amounts that tend to $0$, and push the measures
$m_j$ forward to measures that converge weakly to $m_\infty$. Composing them with dilations
whose factors tend to $1$, we obtain Borel $\eps$-GH approximations $\Phi_j$ of the balls
$B_4(q_j)$ to $B_4(0^n)$ with $\Phi_j(q_j)=0$, for large $j$. Spheres are null sets for the
Lebesgue measure, and the functions $g$ on $B_4(0^n)$ with $|g|\le1$ and $\Lip g\le1$ form
a compact set for the uniform norm. So the averages of $g\circ\Phi_j$ over $B_4(q_j)$
converge to the average of $g$ over $B_4(0^n)$, uniformly in $g$. Hence $q_j$ is
$\eps$-flat at scale $1$ in the new space for large $j$, that is, at scale $r_j$ in the
original one. So both assertions hold for large $j$, which is a contradiction.

For the last assertion put $S:=r/c_H$. Since $\eps_V\le c_H^2\le c_H$, we have
$S\le r/\eps_V$, so $B_S(q)$ is $\eps_VS$-close, hence $\delta_HS$-close, to $B_S(0^n)$.
Moreover $|\kappa|S^2\le\eps_V/c_H^2\le1$ and $r=c_HS$, so \textup{(i)} and \textup{(ii)} hold
at the scale $r$.
\end{proof}

The intrinsic Reifenberg theorem of Cheeger and Colding \cite[Theorem~A.1.2]{CC97} makes the
regions of the first form topological manifolds. We use it in the local form of Gigli and
Violo \cite[Theorem~2.7]{GV25}, which parametrizes a ball of the first form by a Euclidean
ball, with moduli that depend only on the dimension. We state it for complete metric
spaces.

\begin{lemma}[Reifenberg Parametrizations]\label{m:fr:cor:localUIC}
There are numbers $\delta_U>0$, $C_U\ge1$, $\beta_U\in(0,1)$ and $\ell_U\in(0,\tfrac1{200}]$,
depending only on $n$, with the following property. Let $(Z,d)$ be a
complete metric space, $z\in Z$ and $R>0$, and suppose that the ball $B_R(z)$ is
$(\delta,R)$-Reifenberg, with $\delta\le\delta_U$. Then there is a map
$F\colon B_1(0^n)\to Z$ with $B_{R/2}(z)\subset F(B_1(0^n))\subset B_R(z)$ and
\begin{equation}\label{m:fr:eq:reifholder}
\frac{d(Fu,Fv)}{R}\le C_U\,|u-v|^{\beta_U},\qquad
|u-v|\le C_U\Bigl(\frac{d(Fu,Fv)}{R}\Bigr)^{\beta_U}\qquad(u,v\in B_1(0^n)).
\end{equation}
The restriction of $F$ to the open set $F^{-1}(B_{R/2}(z))$ is a homeomorphism onto
$B_{R/2}(z)$, and the point $v_0:=F^{-1}(z)$ has distance at least
$(2C_U)^{-1/\beta_U}$ from the boundary of $B_1(0^n)$. Moreover, there is an open set
$U\subset Z$ that is homeomorphic to $\R^n$, with $B_{\ell_UR}(z)\subset U\subset
B_{R/100}(z)$.
\end{lemma}
\begin{proof}
By scaling we may assume $R=1$. Let $\eps_R(n)$ and $M_R(n)$ be the constants of
\cite[Theorem~2.7]{GV25} for the dimension $n$, and let
$\delta_U:=\tfrac12\min\{\eps_R(n),(2M_R(n))^{-1}\}$. For $q\in B_1(z)$ and $0<t\le1$, the
Gromov--Hausdorff distance of $B_t(q)$ from $B_t(0^n)$ is less than $2\delta t$
\cite[Corollary~7.3.28]{BBI01}, hence less than $\eps_R(n)\,t$; for $q=z$ and $t=1$ it is
less than $(2M_R(n))^{-1}$. So the theorem applies to the ball $B_1(z)$ and gives a
bi-H\"older map $F\colon B_1(0^n)\to Z$ whose image contains $B_{1/2}(z)$. The map
constructed in \cite[Section~5.7]{GV25} takes values in $B_1(z)$, and it is the
composition of a dilation of $B_1(0^n)$, by a factor that is bounded below in terms of
$n$, with a map that satisfies the bounds (5.70) and (5.72) there. These bounds have the
exponents $1/(1+\alpha)$ and $1/(1-\alpha)$ for a number $\alpha=\alpha(n)\in(0,1)$, and
constants that depend only on $n$. Since $|u-v|<2$ for $u,v\in B_1(0^n)$, this gives
\eqref{m:fr:eq:reifholder} with $\beta_U:=1-\alpha$ and a constant $C_U=C_U(n)\ge1$.

A bi-H\"older map is a homeomorphism onto its image, so the restriction of $F$ to the
open set $F^{-1}(B_{1/2}(z))$ is a homeomorphism onto $B_{1/2}(z)$. Since $Z$ is
complete, $F$ and the bounds \eqref{m:fr:eq:reifholder} extend to the closed ball
$\bar B_1(0^n)$. If $u\in\partial B_1(0^n)$ and $F(u)\in B_{1/2}(z)$, then $F(u)=F(u')$
for some $u'\in B_1(0^n)$, and the second bound gives $u=u'$, which is impossible.
Hence $d(F(u),z)\ge\tfrac12$ for $u\in\partial B_1(0^n)$, and the first bound gives
$|u-v_0|\ge(2C_U)^{-1/\beta_U}$.

Let $\eta:=\tfrac12(100C_U)^{-1/\beta_U}$ and
$\ell_U:=\min\{\tfrac1{200},(\eta/(2C_U))^{1/\beta_U}\}$. Then $B_\eta(v_0)\subset
B_1(0^n)$, and $U:=F(B_\eta(v_0))$ lies in $B_{1/100}(z)$, by the first bound. So $U$ is
the image of an open ball under the homeomorphism of $F^{-1}(B_{1/2}(z))$ onto
$B_{1/2}(z)$, and it is open and homeomorphic to $\R^n$. If $y\in B_{\ell_U}(z)$, then
$y=F(u)$ for some $u\in B_1(0^n)$, and the second bound gives
$|u-v_0|<C_U\ell_U^{\beta_U}\le\eta/2$, so $y\in U$.
\end{proof}

In an $\RCD$ space the parametrizations also fix the essential dimension, so that the
integer $n$ of the hypotheses is the dimension of the space. The corollary is used for the
regions of the harmonic charts in \cref{ch:harmoniccomparison}.

\begin{corollary}[Reifenberg Regions]\label{ls:cor:region}
There is $\eps_D=\eps_D(n,N)>0$ with the following property. Let $(Z,d,m)$ be an
$\RCD(\kappa,N)$ space and $\Omega\subset Z$ a non-empty open set that is
$(\eps_D,\lambda)$-Reifenberg for some $\lambda>0$. Then $Z$ has the essential dimension $n$, and $\Omega$ is a
topological manifold of dimension $n$.
\end{corollary}
\begin{proof}
The spaces $\R^k$ with $0\le k\le N$ and $k\ne n$, with normalized Lebesgue measures, have
a positive distance $3\eps'$ for $\mathbb{D}$ from $\R^n$. Let $\eps_D\le\eps_V(n,N,\eps')$. By
\cref{ls:lem:measured} the rescalings of $Z$ at a point $q\in\Omega$ by the factors $1/s$, for
all small $s$, are $\eps'$-close to $\R^n$, so no tangent cone at $q$ is one of the spaces
$\R^k$ with $k\ne n$. Since $m(\Omega)>0$, the essential dimension is $n$.

For the second assertion let $\eps_D$ be also at most $\delta_U(n)$. Let $q_0\in\Omega$,
and choose $R\in(0,\lambda]$ with $B_R(q_0)\subset\Omega$. Then the ball $B_R(q_0)$ is
$(\eps_D,R)$-Reifenberg, and $Z$ is complete, because it is proper by
\ref{m:fr:R:geometry}. So \cref{m:fr:cor:localUIC} gives a homeomorphism of an open subset
of $B_1(0^n)$ onto $B_{R/2}(q_0)$. Thus every point of $\Omega$ has a neighbourhood that is
homeomorphic to an open subset of $\R^n$.
\end{proof}

The third form gives the second, for an arbitrary measure. The frame theorem of
\cref{m:book:ch:07} is stated for Euclidean points, while from \cref{m:fr:sec:gauges} on the
standing hypothesis is a flat region; the next lemma connects the two.

\begin{lemma}[Flat Points]\label{m:fr:lem:flat}
Let $(Z,d,m)$ be an $\RCD(\kappa,N)$ space with $N<\infty$, and let $x\in Z$ be $\eps$-flat at
scale $r$, with $\eps\le\eps(n)$. Then $x$ is $\eps$-Euclidean at scale $r$.
\end{lemma}
\begin{proof}
A map $\Phi$ as in \cref{m:fr:def:flat} is a Borel $\eps r$-GH approximation
$B_{4r}(x)\to B_{4r}(0^n)$ with $\Phi(x)=0$, as required in \cref{m:fr:def:euclid}, so it
remains to prove the thickness condition. Put $s:=\eps^{1/2}r$, and assume
$\eps\le2^{-8}$. If $g$ is a
function on $B_{4r}(0^n)$ with $|g|\le1$ and $\Lip g\le1/s$, then $\eps^{1/2}g$ is one of the
functions of \cref{m:fr:def:flat}, so
\begin{equation*}
\Bigl|\fint_{B_{4r}(x)}g\circ\Phi\,dm-\fint_{B_{4r}(0^n)}g\,d\cH^n\Bigr|\le\eps^{1/2}.
\end{equation*}
Let $p\in B_r(x)$ and $v:=\Phi(p)$; then $|v|<r+\eps r$. Consider the functions
\begin{equation*}
g_1(w):=\min\bigl\{1,s^{-1}\bigl(\tfrac18r-\eps r-|w-v|\bigr)_+\bigr\},\qquad
g_2(w):=\min\bigl\{1,s^{-1}\bigl(2r+\eps r+s-|w|\bigr)_+\bigr\}.
\end{equation*}
The function $g_1$ equals $1$ on the ball of radius $\tfrac18r-\eps r-s$ about $v$, which lies
in $B_{4r}(0^n)$. If $q\in B_{4r}(x)$ and $g_1(\Phi(q))>0$, then
$d(q,p)\le|\Phi(q)-v|+\eps r<\tfrac18r$; so $g_1\circ\Phi$ is at most the indicator function
of $B_{r/8}(p)$. The function $g_2$ vanishes outside the ball of radius $2r+\eps r+s$
about $0$, and $g_2\circ\Phi=1$ on $B_{2r}(x)$, because
$|\Phi(q)|\le d(q,x)+\eps r$. Hence
\begin{equation*}
\frac{m(B_{r/8}(p))}{m(B_{4r}(x))}\ge\Bigl(\frac{1-8\eps-8\eps^{1/2}}{32}\Bigr)^n-\eps^{1/2},
\qquad
\frac{m(B_{2r}(x))}{m(B_{4r}(x))}\le\Bigl(\frac{2+\eps+\eps^{1/2}}4\Bigr)^n+\eps^{1/2}.
\end{equation*}
As $\eps\to0$ the quotient of these bounds tends to $16^{-n}$, so
$m(B_{r/8}(p))\ge32^{-n}\,m(B_{2r}(x))$ when $\eps\le\eps(n)$.
\end{proof}

The three forms are related as follows. By the last assertion of \cref{m:fr:lem:flatness},
every point of an $(\eps_V,\lambda)$-Reifenberg subset of an $\RCD(\kappa,N)$ space, with
$\eps_V=\eps_V(n,N,\eps)$, is $\eps$-flat at every scale $r\le\eps_V\lambda$ with
$|\kappa|r^2\le\eps_V$, and by \cref{m:fr:lem:flat} an $\eps$-flat point with $\eps\le\eps(n)$ is
$\eps$-Euclidean at the same scale.

\subsection{Rigid Motions}\label{m:fr:sec:rigid}

Maps that almost preserve distances are close to rigid motions.

\begin{definition}[Rigid Motions and Their Parts]\label{m:fr:def:rigid}
A \emph{rigid motion} of $\R^n$ is a map \mbox{$E(v)=R_Ev+b_E$} with $R_E\in O(n)$ and $b_E\in\R^n$;
we call $R_E$ its \emph{rotation part} and $b_E$ its \emph{translation part}.
\end{definition}

Quantitative forms of the next lemma are
known: for an $\eps\,d(A)$-nearisometry of a bounded set $A\subset\R^n$ the distance to the
nearest isometry is at most $c(n)\sqrt\eps\,d(A)$, and for thick sets it is linear in $\eps$
\cite[2.2 and 3.3]{ATV01}, \cite[Section~4]{Vai02}. We need only the non-quantitative
statement for centred approximations, and we include its short compactness proof. The lemma
is used for the comparison of harmonic charts in \cref{ch:harmoniccomparison}.

\begin{lemma}[Rigid Motions]\label{ls:lem:compare}
For every $\eta\in(0,\tfrac16]$ there is $\eps_1=\eps_1(n,\eta)\in(0,\eta]$ with the following
property. Let $\Phi$ be a centred $\eps s$-GH approximation of a ball $B_s(w)$ of a metric
space, with $\eps\le\tfrac12\eps_1$.
\begin{enumerate}[label=\textup{(\roman*)}]
\item If $\beta\colon B_s(w)\to\R^n$ changes distances by at most $\eps s$, there is a rigid
motion $E$ with $|\beta(p)-E(\Phi(p))|\le2\eta s$ for all $p\in B_s(w)$.
\item For $\alpha\ge0$, if $E$ and $E'$ are rigid motions with $|E(\Phi(p))-E'(\Phi(p))|\le\alpha s$ for all
$p\in B_s(w)$, then $\|R_E-R_{E'}\|\le4\alpha$ and $|b_E-b_{E'}|\le\alpha s$.
\end{enumerate}
\end{lemma}
\begin{proof}
(i) By scaling we may assume that $s=1$. Suppose that the assertion fails for some $\eta$.
Then for every $j\ge1$ there are a centred $\tfrac1j$-GH approximation $\Phi_j$ of a ball
$B_1(w_j)$ and a map $\beta_j$ on $B_1(w_j)$ that changes distances by at most $\tfrac1j$,
such that no rigid motion $E$ satisfies $|\beta_j-E\circ\Phi_j|\le2\eta$. After a
translation of $\beta_j$ we may assume that $\beta_j(w_j)=0$; then $|\beta_j|\le2$. Let $Q$
be a countable dense subset of $B_1(0^n)$, and for $x\in Q$ choose $p^x_j\in B_1(w_j)$ with
$|\Phi_j(p^x_j)-x|\le\tfrac2j$. After passing to a subsequence, $\beta_j(p^x_j)$ converges to a
point $F(x)$ for every $x\in Q$, and
\begin{equation*}
|F(x)-F(x')|=\lim_j|\beta_j(p^x_j)-\beta_j(p^{x'}_j)|=\lim_jd(p^x_j,p^{x'}_j)=|x-x'| .
\end{equation*}
So $F$ extends to a map of the closed unit ball that preserves distances, and such a map
is the restriction of a rigid motion $E$. For large $j$ there is $p'_j\in B_1(w_j)$ with
$|\beta_j(p'_j)-E(\Phi_j(p'_j))|>2\eta$. Since $|\Phi_j(p'_j)|<1+\tfrac1j$, after passing to a
subsequence $\Phi_j(p'_j)$ converges to a point $x_0$ of the closed unit ball. For $x\in Q$,
\begin{equation*}
|\beta_j(p'_j)-\beta_j(p^x_j)|\le|\Phi_j(p'_j)-\Phi_j(p^x_j)|+\tfrac2j\longrightarrow|x_0-x|,
\end{equation*}
and $\beta_j(p^x_j)\to Ex$. Since $E$ preserves distances, the upper limit of
$|\beta_j(p'_j)-E(\Phi_j(p'_j))|$ is at most $2|x_0-x|$, which is smaller than $2\eta$ for a
suitable $x\in Q$. This is a contradiction.

(ii) Put $A:=R_E-R_{E'}$, so that $E(v)-E'(v)=Av+b_E-b_{E'}$. Since $\Phi(w)=0$, we get
$|b_E-b_{E'}|\le\alpha s$. Let $v$ be a unit vector, and choose points $p_\pm$ with
$|\Phi(p_\pm)\mp\tfrac12sv|\le2\eps s$. Then $|A(\Phi(p_+)-\Phi(p_-))|\le2\alpha s$, and
$\Phi(p_+)-\Phi(p_-)=sv+e$ with $|e|\le4\eps s\le\tfrac12s$, so
$s|Av|\le2\alpha s+\tfrac12s\|A\|$. Taking the supremum over $v$ gives $\|A\|\le4\alpha$.
\end{proof}

\subsection{Harmonic Maps and Group Actions}\label{ls:sec:prelim}

This subsection proves two elementary lemmas: averages of harmonic maps are harmonic,
and near a compact subgroup a Lie group is a product of a ball in a complement of its
Lie algebra with the subgroup. The first is used for the equivariant harmonic charts of
\cref{ls:sec:equivariant,ls:sec:smoothing}, the second for the slices of
\cref{ls:sec:slicesfromcharts} and for the component in the isotropy group of
\cref{es:def:component}.

\subsubsection{Integrals of Harmonic Maps}\label{ls:sec:integrals}

The equivariant harmonic charts of \cref{ls:prop:equivariant,ls:prop:smooth} are averages
of harmonic maps over compact groups. The next lemma says
that such averages are harmonic.

\begin{lemma}[Integrals]\label{ls:lem:integrals}
Let $U$ be an open subset of an $\RCD(\kappa,N)$ space, let $\Lambda$ be a compact metrizable
space with a finite Borel measure $\nu$, and let $F\colon\Lambda\times U\to\R$ be a bounded
continuous function such that $F(\ell,\cdot)$ is harmonic on $U$ and Lipschitz with a
constant $L$, for every $\ell\in\Lambda$. Then $f(p):=\int_\Lambda F(\ell,p)\,d\nu(\ell)$ is
harmonic on $U$ and Lipschitz with the constant $\nu(\Lambda)L$.
\end{lemma}
\begin{proof}
The Lipschitz bound is clear. Let $\psi$ be a Lipschitz function with compact support in
$U$, and let $U'$ be a bounded open set with $\operatorname{supp}\psi\subset U'$ and
$\bar U'\subset U$. The map $\ell\mapsto F(\ell,\cdot)$ is bounded in $H^{1,2}(U')$ and
continuous in $L^2(U')$ by dominated convergence, hence weakly continuous in $H^{1,2}(U')$;
its weak integral with respect to $\nu$ is $f$, and the continuous functional
$g\mapsto\int\langle\nabla g,\nabla\psi\rangle\,dm$ vanishes on every $F(\ell,\cdot)$, hence
on $f$.
\end{proof}
\subsubsection{Neighbourhoods of an Isotropy Group}\label{ls:sec:isotropy}

We record the structure of a Lie group near a compact subgroup and two consequences of
the properness of an action. Let $G$ be a Lie group with Lie algebra $\fg$, and $K$ a
compact subgroup with Lie algebra $\fk$. We fix a norm on $\fg$ that is invariant under
$\Ad(K)$, and we let $\fm$ be the orthogonal complement of $\fk$ for an invariant inner
product; so $\Ad(K)\fm=\fm$. For $r>0$ let $\fm(r)$ be the open ball of radius $r$ in
$\fm$. The product structure of the first part of the next lemma is the basis of the slices
of \cref{ls:sec:slicesfromcharts} and of the component in the isotropy group of \cref{es:def:component}.

\begin{lemma}[The Group near the Isotropy Group]\label{ls:lem:tube}
Let $G$ act properly and by isometries on a metric space $(Z,d)$, and let $K$ be the
isotropy group of a point $z$.
\begin{enumerate}[label=\textup{(\roman*)}]
\item There is $r_0>0$ such that $(\zeta,k)\mapsto\exp(\zeta)k$ is a homeomorphism of
$\fm(r_0)\times K$ onto an open neighbourhood of $K$ in $G$.
\item For every neighbourhood $U$ of $K$ in $G$ there is $\sigma>0$ such that
every $g\in G$ with $d(gz,z)<\sigma$ lies in $U$.
\item There is $c_0>0$ such that $d(\exp(\zeta)z,z)\ge c_0|\zeta|$ for all
$\zeta\in\fm(\tfrac12r_0)$.
\end{enumerate}
\end{lemma}
\begin{proof}
(i) The differential of the map at a point $(0,k)$ is an isomorphism, since
$\fg=\fm\oplus\fk$. So the map is a local diffeomorphism near $\{0\}\times K$. If it were
not injective on $\fm(r)\times K$ for any $r$, there would be pairs
$(\zeta_j,k_j)\ne(\zeta'_j,k'_j)$ with the same image and $\zeta_j,\zeta'_j\to0$; for a
subsequence $k_j\to k$, hence $k'_j\to k$, which contradicts the injectivity near
$(0,k)$. (ii) If $g_jz\to z$, a subsequence of $g_j$ converges, by properness, and its
limit fixes $z$. (iii) If $\fm=0$, the only element of $\fm(\tfrac12r_0)$ is $\zeta=0$, and the assertion holds
for every $c_0>0$. Let $\fm\ne0$, and let $c'>0$ be the minimum of $d(\exp(\zeta')z,z)$ over the
$\zeta'\in\fm$ with $\tfrac14r_0\le|\zeta'|\le\tfrac12r_0$, a non-empty compact set; it is positive because
$\exp(\zeta')\notin K$ for these $\zeta'$, by (i). For $\zeta\in\fm(\tfrac12r_0)$, $\zeta\ne0$,
there is an integer $j\ge1$ with $\tfrac14r_0\le j|\zeta|\le\tfrac12r_0$. The triangle
inequality and the invariance of the distance give
$d(g^jz,z)\le j\,d(gz,z)$ for $g:=\exp(\zeta)$, so
$d(\exp(\zeta)z,z)\ge c'/j\ge2c'|\zeta|/r_0$.
\end{proof}

\section{Canonical Metric Frames and Gauges}\label{ch:canonicalframes}

Canonical frames supply a choice-independent space on which local isometries act exactly. Gauges give its local trivializations, which converge along convergent sequences of spaces.
The frame spaces of the main theorem are formed from these frames in \cref{ch:coveratlas}, and
the convergence statements of this section underlie the limit constructions from
\cref{ch:frameconvergence} on.

\subsection{Canonical Frames}\label{m:book:ch:07}

A singular space need not have tangent frames. The frames defined here replace them: a frame
at a point $x$ is an orthonormal basis of a spectral subspace of dimension $n$ of an
integral operator on a ball of a fixed radius about $x$ whose kernel is formed with the
distance function alone (\cref{m:fr:def:operator,m:fr:def:frames}). So the space of frames is defined without
choices, and isometries that preserve the measure act on it. There are two main results: at almost Euclidean points frames exist and their frame
maps are almost linear (\cref{m:fr:thm:A}), and the frames at the frame-regular points form a
principal $O(n)$-bundle on which these isometries act (\cref{m:fr:thm:bundle}).

\subsubsection{The Frame Operator}\label{m:fr:sec:operator}

We define the frame operator and study it; the frames that it selects are defined in
\cref{m:fr:sec:model}. Throughout \cref{m:book:ch:07}, $(Z,d,m)$ is a metric measure space in the sense
of \cref{m:fr:sec:conv}. Finiteness of $m$ on bounded sets and positivity on balls make the
normalizing constant $V_x^{(r)}$ of \cref{m:fr:def:operator} finite and positive. Nothing else is assumed of $Z$ in
the definitions and in the two main results, \cref{m:fr:thm:A,m:fr:thm:bundle}: no
completeness, no curvature bound, and no hypothesis on the dimension of $Z$; the
convergence statements \cref{m:fr:lem:spectral,m:fr:prop:convergence} assume that the
limit space is proper. The generality is used in \cref{m:fr:sec:coverbundle}, where the
frames are taken on balls of the local covers, which need not be complete, and are
compared with the frames of the completed local covers, which are $\RCD$ spaces, through
\cref{m:fr:rem:local}. The dimension enters as the parameter $n$ of the frames and through
the hypothesis of \cref{m:fr:thm:A}, as explained below; from \cref{m:book:ch:08} on, the
standing hypothesis is an $\RCD(\kappa,N)$ space with an $\eps$-flat region at scale $r$.
We fix the integer $n\ge1$ of \cref{m:fr:sec:closeness} and a frame scale $r>0$,
and write $\psi(t):=\min\{1,(2-t)_+\}$ for $t\ge0$.

The integer $n$ is the dimension of the Euclidean model: the balls in
\cref{m:fr:def:reifenberg,m:fr:def:euclid,m:fr:def:flat} and in \cref{m:fr:def:localcovers}
are balls of $\R^n$. The frame operator is defined on every metric measure space, but the
frames are $n$-tuples: a point is frame-regular when exactly $n$ eigenvalues of its frame
operator lie above the threshold $\lambda_n$, with a gap in the spectrum below them
(\cref{m:fr:def:frames}). This is the only dimensional hypothesis of \cref{m:book:ch:07}. It
holds at the points which are $\eps_C$-Euclidean at scale $r$ in the sense of
\cref{m:fr:def:euclid} (\cref{m:fr:thm:A}), hence, for small parameters, at the flat points of
$\RCD(\kappa,N)$ spaces (\cref{m:fr:lem:flat}) and at small scales on their Reifenberg subsets
(\cref{m:fr:lem:flatness}); the set of frame-regular points may be empty. A point modelled on
$\R^k$ with $k\ne n$ is not frame-regular for the parameter $n$, because the frame operator of
the Lebesgue measure on $\R^k$ has exactly $k$ non-zero eigenvalues, all equal
(\cref{m:fr:lem:model}).

\begin{definition}[Frame Operator]\label{m:fr:def:operator}
Let $(Z,d,m)$ be a metric measure space and $x\in Z$. Put
\begin{equation*}
V_x^{(r)}:=\int_Z\psi\Bigl(\frac{d(x,q)}{r}\Bigr)dm(q),\qquad
\mu_x^{(r)}:=\frac1{V_x^{(r)}}\,\psi\Bigl(\frac{d(x,\cdot)}{r}\Bigr)m,
\end{equation*}
and
\begin{equation*}
\mathsf G_x^{(r)}(p,q):=\frac{d(x,p)^2+d(x,q)^2-d(p,q)^2}{2r^2}\qquad(p,q\in Z).
\end{equation*}
The \emph{frame operator} at $x$ is
\begin{equation}\label{m:fr:eq:K}
(\cK_x^{(r)}f)(p):=\int_Z\mathsf G_x^{(r)}(p,q)\,f(q)\,d\mu_x^{(r)}(q)
\qquad(f\in L^2(\mu_x^{(r)}),\ p\in Z).
\end{equation}
\end{definition}

The scale $r$ specifies the region sampled by the operator: the cut-off has full
weight on $B_r(x)$ and vanishes outside $B_{2r}(x)$. For $a,b>0$, replacing
$(d,m,r)$ by $(ad,bm,ar)$ multiplies $V_x^{(r)}$ by $b$ and leaves
$\mu_x^{(r)}$, $\mathsf G_x^{(r)}$ and $\cK_x^{(r)}$ unchanged under the identity
identification. At fixed $r$ we will write
$V_x$, $\mu_x$, $\mathsf G_x$ and $\cK_x$, omitting $(r)$ in statements and proofs except
in definitions and rescaling statements.

By the kernel bounds of \cref{m:fr:lem:kernel}, $\cK_x$ is a self-adjoint
Hilbert--Schmidt operator on $L^2(\mu_x)$, and \eqref{m:fr:eq:K} defines $\cK_xf$ as a
continuous function on all of $Z$.

Heuristically, the frame operator recovers the directions at $x$ from distances alone.
The kernel $\mathsf G_x$ is the polarization identity written with distances only. In $\R^n$ it
is the inner product, $\mathsf G_x(p,q)=\langle p-x,q-x\rangle/r^2$, so that
\[
(\cK_xf)(p)=\Bigl\langle\frac{p-x}{r},\,w_f\Bigr\rangle,\qquad
w_f:=\int\frac{q-x}{r}\,f(q)\,d\mu_x(q)\in\R^n .
\]
Hence, if $\mu_x$ is not concentrated on a hyperplane through $x$, the range of $\cK_x$ is
the $n$-dimensional space of the linear functions $p\mapsto\langle p-x,v\rangle/r$; the
operator vanishes on the orthogonal complement of this space in $L^2(\mu_x)$ and acts on it
as the second-moment matrix of $\mu_x$ (\cref{m:fr:lem:model}). Its eigenfunctions with
non-zero eigenvalues are thus coordinate functions centred at $x$, that is, directions at
$x$. For Lebesgue measure the second-moment matrix is a multiple of the identity, so every
orthonormal basis of $\R^n$ gives an eigenbasis and no direction is preferred. This is why
a frame will be an arbitrary basis of the eigenspace, normalized with $\cK_x$, and why the
frames at a point form one orbit of $O(n)$ (\cref{m:fr:def:frames}). On a general space
$\cK_x$ need not have finite rank, and the frames use the eigenfunctions above a spectral
gap instead. Finally, the cut-off $\psi(d(x,\cdot)/r)$ equals $1$ on $B_r(x)$ and decreases
linearly to $0$ on $B_{2r}(x)\setminus B_r(x)$ (\cref{m:fr:fig:weights}). Unlike the
indicator function of $B_r(x)$, it changes by at most $d(x,x')/r$ at every point when $x$
is replaced by $x'$. This makes the frame operators continuous in $x$ on every metric
measure space, even when spheres carry positive measure (\cref{m:fr:lem:cont}).

The kernel $\mathsf G_x$ is Schoenberg's Gram kernel with base point $x$ \cite{Sch35}, and $\cK_x$ is
the base-point form of classical multidimensional scaling \cite{YH38,Tor52} of the probability
measure $\mu_x$;
\cref{m:fr:lem:model} is the exactness of classical scaling on Euclidean configurations.
Classical scaling of metric measure spaces, with the kernel centred at the mean of the
measure, and its stability under convergence of metric measure spaces are studied in
\cite{ABK20,KST22,LM22}. Here the kernel is centred at the point $x$, the measure is
localized at $x$, the frames are taken in the eigenspaces above a spectral gap, and the
stability that we need, under GH approximations and weak convergence of
the localized measures, is \cref{m:fr:lem:spectral}.

\begin{center}
\begin{minipage}{\textwidth}
\captionsetup{type=figure}
\centering
\includegraphics[width=\textwidth]{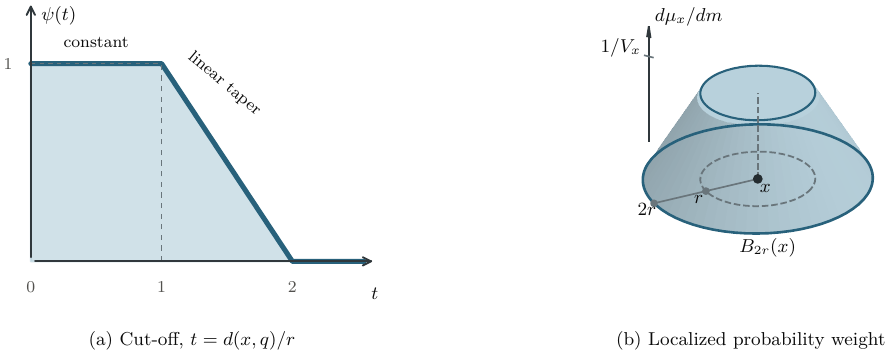}
\caption{The radial cut-off and the localized measure. The height in (b) is the density $d\mu_x/dm$, equal to $V_x^{-1}\psi(d(x,\cdot)/r)$: it is $1/V_x$ on $B_r(x)$,
decreases linearly in distance on the annulus, and is zero at and beyond $2r$.
The disk represents a geodesic ball schematically; neither Euclidean geometry
nor a density of $m$ with respect to planar area is assumed.}
\label{m:fr:fig:weights}
\end{minipage}
\end{center}

The kernel $\mathsf G_x$ is controlled by the triangle inequality alone. We record the bounds
that are used for the operator and for its dependence on $x$.

\begin{lemma}[Kernel Bounds]\label{m:fr:lem:kernel}
For all $x,x',p,p',q\in Z$,
\begin{align*}
|\mathsf G_x(p,q)|&\le\frac{d(x,p)\,d(x,q)}{r^2},\\
|\mathsf G_x(p,q)-\mathsf G_x(p',q)|&\le\frac{d(p,p')\bigl(d(x,p)+d(x,p')+d(x,q)\bigr)}{r^2},\\
|\mathsf G_x(p,q)-\mathsf G_{x'}(p,q)|&\le\frac{d(x,x')\bigl(d(x,p)+d(x',p)+d(x,q)+d(x',q)\bigr)}{2r^2}.
\end{align*}
In particular, $\cK_x$ is a self-adjoint Hilbert--Schmidt operator on $L^2(\mu_x)$ of norm
at most $4$, and for $f\in L^2(\mu_x)$ the function $\cK_xf$ of \eqref{m:fr:eq:K} satisfies
$|\cK_xf(p)|\le2r^{-1}d(x,p)\|f\|_{L^1(\mu_x)}$ and is $10r^{-1}\|f\|_{L^1(\mu_x)}$-Lipschitz
on $B_{4r}(x)$.
\end{lemma}
\begin{proof}
Put $a=d(x,p)$, $b=d(x,q)$ and $c=d(p,q)$. Since ${|a-b|\le c\le a+b}$, we have
${-2ab\le a^2+b^2-c^2\le2ab}$, which is the first bound. The other two follow from the
identity ${|s^2-t^2|=|s-t|(s+t)}$ and the triangle inequality. The measure $\mu_x$ is a probability
measure concentrated on $B_{2r}(x)$, where $|\mathsf G_x|<4$, and $\mathsf G_x$ is symmetric; this gives
the statement on $\cK_x$. The bounds on $\cK_xf$ follow from the first two inequalities
with $q\in B_{2r}(x)$.
\end{proof}
\subsubsection{The Euclidean Model and the Frames}\label{m:fr:sec:model}

We compute the frame operator of a measure on $\R^n$ and define the frames of a general space
by the picture that this gives. The next lemma is the precise
form of the heuristic of \cref{m:fr:sec:operator}, for an arbitrary measure on $\R^n$. The measure need not be the Lebesgue measure, which is why frames will be
normalized with $\cK_x$ and not with the $L^2$ inner product (\cref{m:fr:def:frames}).

\begin{lemma}[Euclidean Model]\label{m:fr:lem:model}
Let $\mu$ be a Borel probability measure on the closed ball $\bar B_{2r}(0^n)$ whose
second-moment matrix $S:=r^{-2}\int y\otimes y\,d\mu(y)$ is positive definite, and let
$\cK$ be the integral operator on $L^2(\mu)$ with kernel
$\mathsf G(y,y')=\langle y,y'\rangle/r^2$. For $v\in\R^n$ put $f_v:=\langle\cdot,v\rangle/r$.
\begin{enumerate}[label=\textup{(\roman*)}]
\item $\cK f_v=f_{Sv}$ and $\langle f_v,f_w\rangle_{L^2(\mu)}=\langle v,Sw\rangle$ for
all $v,w\in\R^n$, and the range of $\cK$ is ${\{f_v:v\in\R^n\}}$.
\item The non-zero eigenvalues of $\cK$, with multiplicities, are the eigenvalues of
$S$, and every eigenfunction with non-zero eigenvalue $\lambda$ is $f_v$ with
$Sv=\lambda v$.
\item If $v_1,\dots,v_n\in\R^n$ satisfy
$\langle f_{v_a},\cK f_{v_b}\rangle_{L^2(\mu)}=\delta_{ab}$, then the map
$y\mapsto r\,(\cK f_{v_1}(y),\dots,\cK f_{v_n}(y))$ is a linear isometry of $\R^n$.
\end{enumerate}
\end{lemma}
\begin{proof}
For $f\in L^2(\mu)$ we have $\cK f=f_w$ with $w=r^{-1}\int y'f(y')\,d\mu(y')$, and for
$f=f_v$ this vector is $Sv$. The formula for the inner product is the definition of $S$.
Since $\|f_v\|^2=\langle v,Sv\rangle>0$ for $v\ne0$, the map $v\mapsto f_v$ is injective,
and it intertwines $S$ with the restriction of $\cK$ to its range. This proves (i) and
(ii). In (iii), $r\,\cK f_{v_a}=\langle\cdot,Sv_a\rangle$, and
$\langle Sv_a,Sv_b\rangle=\langle f_{v_a},f_{Sv_b}\rangle=\delta_{ab}$, so the matrix with
rows $Sv_a$ is orthogonal.
\end{proof}

By \cref{m:fr:lem:model}, if a measure on $\R^n$ is not concentrated on a hyperplane
through $x$, the non-zero spectrum of its frame operator at $x$ consists of $n$ eigenvalues,
whose eigenfunctions are the coordinate functions, and all other eigenvalues vanish. On a general space we select these eigenvalues with the threshold $\lambda_n:=32^{-n-1}$.
At the almost Euclidean points of \cref{m:fr:def:euclid} the $n$ eigenvalues are close to
numbers that are at least $32^{-n}/16=2\lambda_n$, by the lower bound on the measures of
small balls there, and all other eigenvalues are close to $0$ (proof of \cref{m:fr:thm:A}). So the band from $\tfrac12\lambda_n$ to $\lambda_n$ is free of
spectrum at such points, and this band makes the spectral condition below an open
condition and the eigenspaces that it selects continuous in $x$ (\cref{m:fr:thm:bundle}).

\begin{definition}[Frames]\label{m:fr:def:frames}
Let $(Z,d,m)$ be a metric measure space and $x\in Z$. The point $x$ is
\emph{frame-regular at scale $r$} if
\begin{equation*}
\spec(\cK_x^{(r)})\subset\bigl(-\tfrac12\lambda_n,\tfrac12\lambda_n\bigr)\cup
(\lambda_n,\infty)
\end{equation*}
and the sum $E_x^{(r)} \le L^2(\mu_x^{(r)})$ of the eigenspaces with eigenvalues in $(\lambda_n,\infty)$ has
dimension $n$. A \emph{frame} at a frame-regular point $x$ is an $n$-tuple
$e=(e_1,\dots,e_n)$ of elements of $E_x^{(r)}$ with
$\langle e_a,\cK_x^{(r)}e_b\rangle_{L^2(\mu_x^{(r)})}=\delta_{ab}$, and its \emph{frame map} is
\begin{equation*}
u_e^{(r)}:=(r \cK_x^{(r)}e_1,\dots, r \cK_x^{(r)}e_n)\colon Z\longrightarrow\R^n.
\end{equation*}
We write $\cR^{(r)}\subset Z$ for the set of frame-regular points and $\Fr^{(r)}(Z)$ for the set
of pairs $(x,e)$, with the projection $\pr(x,e):=x$, the topology induced by the
injection $(x,e)\mapsto(x,u_e^{(r)})$ into $Z\times C(Z,\R^n)$, where $C(Z,\R^n)$ carries the
topology of uniform convergence on bounded sets, and the right action
$(e\cdot A)_b:=\sum_ae_aA_{ab}$ of $O(n)$.
\end{definition}

At fixed $r$ we likewise abbreviate $E_x^{(r)}$, $u_e^{(r)}$, $\cR^{(r)}$ and
$\Fr^{(r)}(Z)$ to $E_x$, $u_e$, $\cR$ and $\Fr(Z)$, omitting the superscript $(r)$. Under $(d,m,r)\mapsto(ad,bm,ar)$ the
spectral subspace $E_x^{(r)}$ is unchanged and $u_e^{(r)}$ is multiplied by $a$.
Thus the normalized proofs may take $r=1$ by dividing the metric and frame maps
by $r$. For a subset $\Omega\subset Z$ we write $\Fr(\Omega):=\pr^{-1}(\Omega)$ for the frames of $Z$
over the points of $\Omega$; they are formed with the distance and the measure of $Z$, not of
$\Omega$.

For Lebesgue measure on $\R^n$ the measure $\mu_x$ has a density proportional to
$\psi(|y-x|/r)$, so its second-moment matrix about $x$, as in \cref{m:fr:lem:model}, is
$\mathsf{s}_n$ times the identity, where
\begin{equation*}
\mathsf{s}_n:=\frac1n\left(\int_{\R^n}|y|^2\psi(|y|)\,dy\right)\left(\int_{\R^n}\psi(|y|)\,dy\right)^{-1}
\ge\frac1{(n+2)\,2^n},
\end{equation*}
the bound because $\psi=1$ on the unit ball and $\psi\le1$ on the ball of radius $2$. By
\cref{m:fr:lem:model} the $n$ non-zero eigenvalues of $\cK_x$ equal $\mathsf{s}_n$, which lies far
above the threshold $\lambda_n=32^{-n-1}$, so the range of $\cK_x$ is $E_x$. The
normalization $\langle e_a,\cK_xe_b\rangle=\delta_{ab}$ compensates for the second moments
of the measure, and the factor $r$ turns the values of $\cK_xe_a$ back into lengths; so the
frame maps are exactly the isometries $p\mapsto R(p-x)$ with $R\in O(n)$, by
\cref{m:fr:lem:model}(iii). The spectral definition retains this space of directions at
the almost Euclidean points of \cref{m:fr:def:euclid}.

\subsubsection{Spectral Stability and the Existence of Frames}\label{m:fr:sec:stability}

Our first result says that frames exist at almost Euclidean points and that their frame
maps are almost linear there.

\begin{theorem}[Frames at Almost Euclidean Points]\label{m:fr:thm:A}
Let $n\ge1$ and $r>0$. For every $\eta>0$ there is $\eps_C>0$, depending only on $n$
and $\eta$, such that every point $x$ of a metric measure space which is
$\eps_C$-Euclidean at scale $r$ in the sense of \cref{m:fr:def:euclid}, through a map
$\Phi$, is frame-regular, and for every frame $e$ at $x$ there is $R\in O(n)$ with
$\sup_{B_{3r}(x)}|u_e-R\circ\Phi|\le\eta r$.
\end{theorem}

\begin{proof}
Normalize $r=1$, write $\mu=\mu_x$ and $\cK=\cK_x$, and suppose that $x$ is
$\eps$-Euclidean through $\Phi$. On the same Hilbert space as $\cK$, compare it with
the finite-rank Euclidean Gram operator
\[
Mf:=\int\Phi f\,d\mu\in\R^n,\qquad L:=M^*M,
\qquad (Lf)(p)=\int\langle\Phi(p),\Phi(q)\rangle f(q)\,d\mu(q).
\]
Polarizing the distance errors gives
$|\mathsf G_x(p,q)-\langle\Phi(p),\Phi(q)\rangle|\le C\eps$ on $B_4(x)^2$.
Consequently $\|\cK-L\|\le C\eps$, and
$|(\cK-L)f(p)|\le C\eps\|f\|_{L^2(\mu)}$ on $B_3(x)$.

The matrix $S:=MM^*=\int\Phi\otimes\Phi\,d\mu$ satisfies $S\ge2\lambda_n\one$.
Indeed, for a unit vector $v$ choose $p$ with $|\Phi(p)-v/2|<2\eps$.
For small $\eps$, $p\in B_1(x)$ and $B_{1/8}(p)\subset B_1(x)$, so the thickness
condition in \cref{m:fr:def:euclid} gives $\mu(B_{1/8}(p))\ge32^{-n}$.
On this ball $\langle\Phi(q),v\rangle\ge1/2-1/8-3\eps\ge1/4$, proving
$\langle Sv,v\rangle\ge32^{-n}/16=2\lambda_n$.
The non-zero spectrum of $L$ is that of $S$. The min--max principle and
$\|\cK-L\|\le C\eps<\lambda_n/2$ therefore give exactly $n$ eigenvalues of $\cK$
above $\lambda_n$, with all remaining eigenvalues in
$(-\lambda_n/2,\lambda_n/2)$. Thus $x$ is frame-regular.

For a frame $e$, let $B$ be the matrix whose $a$-th row is $Me_a$, and put
$Uv:=\sum_a v_ae_a$ for $v\in\R^n$. The frame normalization gives
$U^*\cK U=\one$ and $\|U\|\le\lambda_n^{-1/2}$. Hence
\[
BB^{\T}=U^*LU=\one+U^*(L-\cK)U,\qquad
\|BB^{\T}-\one\|\le C\eps/\lambda_n.
\]
For small $\eps$, $R:=(BB^{\T})^{-1/2}B$ is orthogonal and
$\|B-R\|\le C(n)\eps$. The pointwise operator estimate gives
$|u_e-B\Phi|\le C(n)\eps$ on $B_3(x)$, and therefore
$|u_e-R\Phi|\le C(n)\eps$ there. Choose $\eps_C$ so that this bound is at most
$\eta$, and restore the scale $r$.
\end{proof}

The proof uses no curvature condition or completeness: the thickness assumption keeps
all $n$ Euclidean directions present, and the frame normalization makes their coordinate
matrix almost orthogonal. For convergence of frames to a possibly non-Euclidean limit,
we also need the following spectral stability lemma. It is an instance of the spectral
approximation theory for collectively compact sequences of operators, cf.\ \cite{Ans71};
we give the direct proof because the operators act on different spaces. The frame scale
$r$ is fixed throughout, and the indices from which the conclusions hold depend on it.

\begin{lemma}[Spectral Stability]\label{m:fr:lem:spectral}
For $j\in\mathbb N$ let $(Z_j,d_j,m_j)$ be metric measure spaces with points $x_j$, and
write $\mu_j$, $\mathsf G_j$, $\cK_j$ for the objects of \cref{m:fr:def:operator} at $x_j$. Let $(Y,d_Y)$
be a proper metric space with a point $y_0$, and let
$\Phi_j\colon B_{4r}(x_j)\to Y$ be Borel maps with $\Phi_j(x_j)=y_0$ and
$|d_Y(\Phi_j(p),\Phi_j(q))-d_j(p,q)|\le\eps_jr$, where $\eps_j\to0$. Suppose that the
measures $\nu_j:=(\Phi_j)_\#\mu_j$ converge weakly to a Borel probability measure
$\mu_\infty$ on $Y$, and let $\cK_\infty$ be the integral operator on $L^2(\mu_\infty)$
with kernel $\mathsf G_{y_0}$.
\begin{enumerate}[label=\textup{(\roman*)}]
\item For every continuous function $g\colon\R\to\R$ that vanishes on a neighbourhood
of $0$, the sum of $g$ over the eigenvalues of $\cK_j$, repeated according to
multiplicity, converges to the corresponding sum for $\cK_\infty$.
\item Let $f_j\in L^2(\mu_j)$ satisfy $\cK_jf_j=\lambda_jf_j$ and $\|f_j\|=1$, with
$\inf_j|\lambda_j|>0$. After passing to a subsequence, $\lambda_j\to\lambda_\infty$ and
there is $f_\infty\in L^2(\mu_\infty)$ with $\cK_\infty f_\infty=\lambda_\infty f_\infty$
and $\|f_\infty\|=1$ such that
\begin{equation}\label{m:fr:eq:eigconv}
\sup_{p\in B_{4r}(x_j)}\bigl|\cK_jf_j(p)-(\cK_\infty f_\infty)(\Phi_j(p))\bigr|
\longrightarrow0.
\end{equation}
If $f_j'$ is a second such sequence with limit $f_\infty'$, then
$\langle f_j,f_j'\rangle\to\langle f_\infty,f_\infty'\rangle$.
\end{enumerate}
\end{lemma}
\begin{proof}
Let $r=1$. The measures $\nu_j$ are concentrated on $\bar B_3(y_0)$, which is compact,
and $\mu_\infty$ is concentrated on $\bar B_2(y_0)$. For $p,q\in B_4(x_j)$ each of the
three squared distances in $\mathsf G_j(p,q)$ differs from the corresponding one in
$\mathsf G_{y_0}(\Phi_j(p),\Phi_j(q))$ by at most $17\eps_j$, so
\begin{equation}\label{m:fr:eq:kernelclose}
\bigl|\mathsf G_j(p,q)-\mathsf G_{y_0}(\Phi_j(p),\Phi_j(q))\bigr|\le26\,\eps_j
\qquad(p,q\in B_4(x_j)).
\end{equation}

(i) For $\ell\ge2$ the operator $\cK_j^\ell$ is of trace class, and
\begin{equation*}
\tr\bigl(\cK_j^\ell\bigr)=\int\!\!\cdots\!\!\int \mathsf G_j(p_1,p_2)\,\mathsf G_j(p_2,p_3)\cdots \mathsf G_j(p_\ell,p_1)\,
d\mu_j(p_1)\cdots d\mu_j(p_\ell),
\end{equation*}
and likewise for $\cK_\infty$. By \eqref{m:fr:eq:kernelclose} and the bound $|\mathsf G_j|\le4$ on the
support of $\mu_j$, the right-hand side differs by at most $26\,\ell\,5^{\ell-1}\eps_j$ from the
integral of
$\mathsf G_{y_0}(y_1,y_2)\cdots \mathsf G_{y_0}(y_\ell,y_1)$ against $\nu_j^{\otimes\ell}$. This integrand is
continuous, and $\nu_j^{\otimes\ell}\to\mu_\infty^{\otimes\ell}$ weakly on the compact set
$\bar B_3(y_0)^\ell$, so $\tr(\cK_j^\ell)\to\tr(\cK_\infty^\ell)$ for every $\ell\ge2$. All eigenvalues
lie in $[-4,4]$, and the sum of their squares is at most $16$, for $\cK_j$ and for
$\cK_\infty$. Given $g$ and $\delta>0$, the Weierstrass theorem applied to the
continuous function $g(t)/t^2$ gives a polynomial $P$ with $|g(t)-t^2P(t)|\le\delta t^2$
on $[-4,4]$. The sum of $g$ over the eigenvalues then differs from
$\tr(\cK^2P(\cK))$ by at most $16\delta$, for $\cK=\cK_j$ and for $\cK=\cK_\infty$, and
$\tr(\cK_j^2P(\cK_j))\to\tr(\cK_\infty^2P(\cK_\infty))$. Since $\delta$ is arbitrary,
(i) follows.

(ii) Since $|\lambda_j|\le4$, we may assume $\lambda_j\to\lambda_\infty\ne0$. The signed
measures $\sigma_j:=(\Phi_j)_\#(f_j\mu_j)$ on $\bar B_3(y_0)$ have total variation at most
$\|f_j\|_{L^1(\mu_j)}\le1$, so after passing to a subsequence they converge weakly to a
signed measure $\sigma$. Put $k(y):=\int \mathsf G_{y_0}(y,y')\,d\sigma(y')$. By
\eqref{m:fr:eq:kernelclose}, for $p\in B_4(x_j)$,
\begin{equation*}
\Bigl|\cK_jf_j(p)-\int \mathsf G_{y_0}(\Phi_j(p),y')\,d\sigma_j(y')\Bigr|\le26\,\eps_j .
\end{equation*}
By \cref{m:fr:lem:kernel} the functions $\mathsf G_{y_0}(y,\cdot)$ with $y\in\bar B_5(y_0)$ are
uniformly bounded and equicontinuous on $\bar B_3(y_0)$, so their integrals against
$\sigma_j$ converge to $k(y)$ uniformly in $y$. Hence
\begin{equation}\label{m:fr:eq:kconv}
\sup_{p\in B_4(x_j)}\bigl|\cK_jf_j(p)-k(\Phi_j(p))\bigr|\longrightarrow0 .
\end{equation}
Now $f_j=\lambda_j^{-1}\cK_jf_j$ almost everywhere with respect to $\mu_j$. For a bounded
continuous function $g$ on $Y$, \eqref{m:fr:eq:kconv} therefore gives
\begin{equation*}
\int g\,d\sigma_j=\lambda_j^{-1}\int g(\Phi_j(q))\,\cK_jf_j(q)\,d\mu_j(q)
=\lambda_j^{-1}\int g\,k\,d\nu_j+o(1)\longrightarrow\lambda_\infty^{-1}\int g\,k\,d\mu_\infty,
\end{equation*}
so $\sigma=f_\infty\mu_\infty$ with $f_\infty:=\lambda_\infty^{-1}k$. By the definition of
$k$ this says $\cK_\infty f_\infty=\lambda_\infty f_\infty$ and $k=\cK_\infty f_\infty$,
and \eqref{m:fr:eq:kconv} is \eqref{m:fr:eq:eigconv}. Finally,
\begin{equation*}
\langle f_j,f_j'\rangle=\frac{1}{\lambda_j\lambda_j'}\int\cK_jf_j\,\cK_jf_j'\,d\mu_j
=\frac{1}{\lambda_j\lambda_j'}\int k\,k'\,d\nu_j+o(1)\longrightarrow
\langle f_\infty,f_\infty'\rangle,
\end{equation*}
and the case $f_j'=f_j$ gives $\|f_\infty\|=1$.
\end{proof}

\subsubsection{The Frame Bundle}\label{m:fr:sec:bundle}

The frames at the frame-regular points form a bundle on which the measure-preserving
isometries act. We prove it here, and \cref{m:fr:sec:convergence} shows that frames converge when the
spaces do. The bundle structure descends to the local covers in \cref{m:fr:sec:coverbundle} and gives
the principal bundles of \cref{thm:main}, and the convergence statement is the basis of
\cref{m:fr:sec:frameconv}.

\begin{theorem}[The Frame Bundle]\label{m:fr:thm:bundle}
Let $n\ge1$ and $r>0$, and let $(Z,d,m)$ be a metric measure space. The set $\cR$ is
open, and the projection $\pr$ of $\Fr(Z)$ onto $\cR$ is a principal
$O(n)$-bundle. Every measure-preserving isometry $\varphi$ of $Z$ preserves $\cR$, the
map \mbox{$(x,e)\mapsto(\varphi(x),e\circ\varphi^{-1})$} is a bundle automorphism of
$\Fr(Z)$, and $u_{e\circ\varphi^{-1}}=u_e\circ\varphi^{-1}$.
\end{theorem}

The set $\cR$ may be empty. The frame
operators at different points act on different Hilbert spaces, so we first move them to a
fixed one.

Fix $x_0\in Z$ and put $D_0:=B_{r/2}(x_0)$, $B_0:=B_{3r}(x_0)$ and $v_0:=m(D_0)>0$. For
$x\in D_0$ the measure $\mu_x$ is concentrated on $B_{2r}(x)\subset B_0$ and
$V_x\ge v_0$. Put
\begin{equation*}
a_x:=\Bigl(\frac{1}{V_x}\,\psi\Bigl(\frac{d(x,\cdot)}{r}\Bigr)\Bigr)^{1/2},\qquad
J_x\colon L^2(\mu_x)\to L^2(B_0,m),\quad J_xf:=a_xf,\qquad
\widetilde\cK_x:=J_x\cK_xJ_x^* .
\end{equation*}
Then $J_x$ is a linear isometry and $\widetilde\cK_x$ is the integral operator on
$L^2(B_0,m)$ with kernel $a_x(p)\mathsf G_x(p,q)a_x(q)$. The following lemma is the continuity
statement behind the bundle structure.

\begin{lemma}[Continuity of the Frame Operators]\label{m:fr:lem:cont}
The operators $\widetilde\cK_x$ and $\cK_x$ have the same non-zero eigenvalues with the
same multiplicities, and $J_x$ maps eigenfunctions to eigenfunctions. There is a
function $\omega$ with $\omega(t)\to0$ as $t\to0$, depending on $r$, $v_0$ and $m(B_0)$,
such that for $x,x'\in D_0$ and $p\in Z$,
\begin{equation*}
\|\widetilde\cK_x-\widetilde\cK_{x'}\|\le\omega(d(x,x')),\qquad
\sup_{q\in B_0}\bigl|a_x(q)\mathsf G_x(p,q)-a_{x'}(q)\mathsf G_{x'}(p,q)\bigr|
\le\Bigl(1+\frac{d(x_0,p)}{r}\Bigr)\omega(d(x,x')) .
\end{equation*}
\end{lemma}
\begin{proof}
The first assertion holds because $J_x$ is an isometry onto its range and
$\widetilde\cK_x$ vanishes on the orthogonal complement of that range. For the second,
$\psi$ is $1$-Lipschitz, so $\psi(d(x,q)/r)$ changes by at most $d(x,x')/r$ when $x$ is
replaced by $x'$, and ${|V_x-V_{x'}|\le m(B_0)\,d(x,x')/r}$. Since ${|\sqrt s-\sqrt t|\le\sqrt{|s-t|}}$ and
$V_x\ge v_0$, the functions $a_x$ are bounded by $v_0^{-1/2}$, and $a_{x'}\to a_x$
uniformly on $B_0$ as $d(x,x')\to0$, at a rate depending only on $r$, $v_0$ and
$m(B_0)$. \Cref{m:fr:lem:kernel} bounds $\mathsf G_x(p,q)$ and $\mathsf G_x(p,q)-\mathsf G_{x'}(p,q)$ for
$q\in B_0$ by multiples of $1+d(x_0,p)/r$ and of $(1+d(x_0,p)/r)\,d(x,x')/r$. This gives
the second inequality, and the first follows because the operator norm is at most
$m(B_0)$ times the supremum of the kernel on $B_0\times B_0$.
\end{proof}

We now prove \cref{m:fr:thm:bundle}. The frame map determines the frame, because
$u_e=r\,\cK_xe$ and $\cK_x$ is injective on $E_x$; this is why $(x,e)\mapsto(x,u_e)$ is an
injection and defines a topology on $\Fr(Z)$.

\begin{proof}[Proof of \cref{m:fr:thm:bundle}]
\emph{Openness.} Let $x_0\in\cR$ and use the notation fixed before \cref{m:fr:lem:cont}. The set
$\Sigma:=\spec(\widetilde\cK_{x_0})$ is a compact subset of
$\cO:=(-\tfrac12\lambda_n,\tfrac12\lambda_n)\cup(\lambda_n,\infty)$; let
$\delta_{\mathrm{sp}}>0$ be its distance to $\R\setminus\cO$. For self-adjoint operators the spectrum
of $\widetilde\cK_x$ lies in the $\|\widetilde\cK_x-\widetilde\cK_{x_0}\|$-neighbourhood
of $\Sigma$. By \cref{m:fr:lem:cont} there is a neighbourhood $D_1\subset D_0$ of $x_0$ on
which this norm is less than $\delta_{\mathrm{sp}}/2$, so that $\spec(\widetilde\cK_x)\subset\cO$ and
the circle $\Gamma:=\{|\zeta-10|=10-\tfrac34\lambda_n\}$ stays at a positive
distance from $\spec(\widetilde\cK_x)$, uniformly in $x\in D_1$. The Riesz projection
\begin{equation*}
P_x:=\frac{1}{2\pi\mathrm i}\oint_\Gamma(\zeta-\widetilde\cK_x)^{-1}\,d\zeta
\end{equation*}
is the orthogonal projection onto $J_xE_x$, and it depends continuously on $x\in D_1$ in
the operator norm \cite[Chapter~IV, Section~3]{Kat95}. A norm-continuous family of projections
has locally constant rank, so $\dim E_x=n$ near $x_0$. Hence $\cR$ is open.

\emph{Local sections.} Fix a frame $e^0$ at $x_0$ and put
$\tilde e_a(x):=P_xJ_{x_0}e^0_a\in J_xE_x$. The matrix
$\gamma_{ab}(x):=\langle\tilde e_a(x),\widetilde\cK_x\tilde e_b(x)\rangle$ is continuous in
$x$ and equals the identity at $x_0$, so it is positive definite on a neighbourhood
$D_2\subset D_1$ of $x_0$. For $x\in D_2$ put
$\tilde e'(x):=\tilde e(x)\cdot\gamma(x)^{-1/2}$ and $e(x):=J_x^*\tilde e'(x)$. Since
$J_xJ_x^*$ is the identity on the range of $J_x$, we get
$\langle e_a(x),\cK_xe_b(x)\rangle=\langle\tilde e'_a(x),\widetilde\cK_x\tilde e'_b(x)\rangle
=\delta_{ab}$, so $e(x)$ is a frame at $x$, and
\begin{equation}\label{m:fr:eq:section}
u_{e(x)}(p)=r\int_{B_0}a_x(q)\,\mathsf G_x(p,q)\,\tilde e'(x)(q)\,dm(q)\qquad(p\in Z).
\end{equation}
The map $x\mapsto\tilde e'(x)\in L^2(B_0,m)^n$ is continuous, so \eqref{m:fr:eq:section} and the
second inequality of \cref{m:fr:lem:cont} show that $x\mapsto u_{e(x)}$ is continuous for
uniform convergence on bounded sets. Thus $x\mapsto(x,e(x))$ is a continuous section of
$\pr$ over $D_2$.

\emph{Local triviality.} The group $O(n)$ acts freely and transitively on the frames at
each point, because frames are the orthonormal bases of $E_x$ for the inner product
$\langle\cdot,\cK_x\cdot\rangle$, and $u_{e\cdot A}=A^{\T}u_e$. Hence
$(x,A)\mapsto(x,e(x)\cdot A)$ is a continuous bijection of $D_2\times O(n)$ onto
$\pr^{-1}(D_2)$. Its inverse sends $(x,e)$ to $(x,A)$ with
\begin{equation*}
A_{ab}=\langle e_a(x),\cK_xe_b\rangle_{L^2(\mu_x)}
=\frac1r\int_{B_0}a_x(q)\,\tilde e'_a(x)(q)\,u_{e,b}(q)\,dm(q),
\end{equation*}
which is continuous in $(x,e)$ for the topology of $\Fr(Z)$. So $\pr$ is a principal
$O(n)$-bundle.

\emph{Isometries.} Let $\varphi$ be a measure-preserving isometry of $Z$. Then
$V_{\varphi(x)}=V_x$, \mbox{$\varphi_\#\mu_x=\mu_{\varphi(x)}$} and
$\mathsf G_{\varphi(x)}(\varphi(p),\varphi(q))=\mathsf G_x(p,q)$, so the unitary map
$f\mapsto f\circ\varphi^{-1}$ from $L^2(\mu_x)$ to $L^2(\mu_{\varphi(x)})$ intertwines
$\cK_x$ with $\cK_{\varphi(x)}$, also as functions on $Z$:
$\cK_{\varphi(x)}(f\circ\varphi^{-1})=(\cK_xf)\circ\varphi^{-1}$. Hence
$\varphi(\cR)=\cR$, the tuple $e\circ\varphi^{-1}$ is a frame at $\varphi(x)$ whenever
$e$ is a frame at $x$, and \mbox{$u_{e\circ\varphi^{-1}}=u_e\circ\varphi^{-1}$}. The resulting map
of $\Fr(Z)$ is continuous, commutes with the action of $O(n)$, and covers $\varphi$.
\end{proof}

The isometry entered the proof only through its restriction to a ball about the point. We
record this local form, which is the one applied in \cref{ch:coveratlas}.

\begin{remark}[Local Naturality]\label{m:fr:rem:local}
The proof uses $\varphi$ only on $B_{2r}(x)$ for the operator and on the set where $u_e$
is evaluated. If $\varphi\colon B_{4r}(x)\to B_{4r}(x')$ is a measure-preserving isometry
of balls with $\varphi(x)=x'$, then $x$ is frame-regular if and only if $x'$ is, and
$u_{e\circ\varphi^{-1}}=u_e\circ\varphi^{-1}$ on $B_{4r}(x')$. In this form the
construction applies to deck transformations of incomplete local covers.
\end{remark}
\subsubsection{Convergence of Frames}\label{m:fr:sec:convergence}

Frames are stable under measured Gromov--Hausdorff convergence. We state this for one
sequence of base points; by the naturality just proved, it is compatible with isometric
actions. The statement is the basis of the convergence theorem of \cref{m:fr:sec:frameconv}.

\begin{proposition}[Convergence of Frames]\label{m:fr:prop:convergence}
Let $(Z_j,d_j,m_j)$ and $(Z,d,m)$ be metric measure spaces, with $Z$ proper, and let
$x_j\in Z_j$ and $x\in Z$. Let $\Phi_j\colon B_{4r}(x_j)\to Z$ be Borel maps with
$\Phi_j(x_j)=x$ and $|d(\Phi_j(p),\Phi_j(q))-d_j(p,q)|\le\eps_jr$, where $\eps_j\to0$, such
that $\int g\circ\Phi_j\,dm_j\to\int g\,dm$ for every continuous function $g$ with
compact support in $B_{4r}(x)$. If $x$ is frame-regular, then $x_j$ is frame-regular for
all large $j$, and every sequence of frames $e_j$ at $x_j$ has a subsequence along which
\begin{equation*}
\sup_{B_{4r}(x_j)}\bigl|u_{e_j}-u_e\circ\Phi_j\bigr|\longrightarrow0
\end{equation*}
for some frame $e$ at $x$.
\end{proposition}
\begin{proof}
We first show that $\nu_j:=(\Phi_j)_\#\mu_{x_j}$ converges weakly to $\mu_x$. Since $\psi$
is $1$-Lipschitz, $|\psi(d_j(x_j,p)/r)-\psi(d(x,\Phi_j(p))/r)|\le\eps_j$, and both terms
vanish outside $B_{3r}(x_j)$. The masses $m_j(B_{3r}(x_j))$ are bounded, because they are
at most $\int\chi\circ\Phi_j\,dm_j$ for a compactly supported continuous function
$\chi\ge0$ on $B_{4r}(x)$ that equals $1$ on $B_{7r/2}(x)$. Therefore, for every bounded
continuous function $g$ on $Z$,
\begin{equation*}
\int g(\Phi_j(p))\,\psi\Bigl(\frac{d_j(x_j,p)}{r}\Bigr)dm_j(p)
=\int\Bigl(g\,\psi\Bigl(\frac{d(x,\cdot)}{r}\Bigr)\Bigr)\circ\Phi_j\,dm_j+o(1)
\longrightarrow\int g\,\psi\Bigl(\frac{d(x,\cdot)}{r}\Bigr)dm .
\end{equation*}
With $g=1$ this gives $V_{x_j}\to V_x$, and the claim follows.

\Cref{m:fr:lem:spectral} now applies with $Y=Z$, $y_0=x$ and $\cK_\infty=\cK_x$.
In part (i), choose a non-negative continuous test equal to $1$ on
$[-4,-\lambda_n/2]\cup[\lambda_n/2,\lambda_n]$ and vanishing near $0$ and on
$\spec(\cK_x)$. Its spectral sum tends to $0$, excluding these intervals for large $j$.
A continuous test equal to $0$ on $[-4,\lambda_n/2]$ and to $1$ on $[\lambda_n,4]$
then counts the remaining large eigenvalues; the count tends to $n$ and hence equals
$n$ eventually. Thus $x_j$ is frame-regular.
For frames $e_j$ at $x_j$, choose orthonormal eigenfunctions $f_{j,k}$ spanning $E_{x_j}$
and write $e_{j,a}=\sum_k\beta^j_{ak}f_{j,k}$. The frame normalization bounds
$|\beta^j_{ak}|\le\lambda_n^{-1/2}$, so part (ii), after extraction, gives
$\beta^j\to\beta$, $\lambda_{j,k}\to\lambda_k$ and $f_{j,k}\to f_{\infty,k}$
in the sense of \eqref{m:fr:eq:eigconv}. The $f_{\infty,k}$ are orthonormal eigenfunctions of
$\cK_x$ with eigenvalues $\lambda_k>\lambda_n$, so they form a basis of $E_x$, the
tuple $e_a:=\sum_k\beta_{ak}f_{\infty,k}$ satisfies
$\langle e_a,\cK_xe_b\rangle=\sum_k\beta_{ak}\beta_{bk}\lambda_k=\delta_{ab}$, and
\eqref{m:fr:eq:eigconv} gives the uniform convergence of $u_{e_j}$ to $u_e\circ\Phi_j$.
\end{proof}
\subsection{Gauges and Trivializations of the Frame Bundle}\label{m:book:ch:08}

\label{m:fr:sec:gauges}

At a fixed scale $r$ the frame bundle is locally a product, and this subsection makes the product structure explicit. The spaces of this subsection are $\RCD(\kappa,N)$ spaces with $N<\infty$,
and the frames are considered over an open set $\Omega$ that is $\eps$-flat at scale $r$
(\cref{m:fr:def:flat}), so that all its points are frame-regular; this standing hypothesis
is stated before \cref{m:fr:def:gauge}. A map into $\R^n$ that is almost isometric at the scale $r$,
a gauge map, selects at every point of the core of its domain the frame that is best
aligned with it (\cref{m:fr:def:gauge}). The aligned frames form continuous sections
(\cref{m:fr:thm:gauges}). Frame maps are
gauge maps, so every frame at a point gives a trivialization over the intersection of
the flat region with the ball of radius $\tfrac12r$ about the point
(\cref{m:fr:sec:frameenergy}). \Cref{m:fr:sec:frameconv} shows that frames and aligned sections converge along a convergent sequence
of spaces, also together with convergent sequences of isometries
(\cref{m:fr:thm:frameconv}).

\subsubsection{Gauges of the Frame Bundle}\label{m:fr:sec:gaugemaps}
We now work at the scale $r$ of the frames. A map that is almost isometric at this scale
selects, at each point, the frame whose frame map is best aligned with it; we call such
a map a gauge map, because it trivializes the frame bundle over the points that lie
well inside its domain. The main result here, \cref{m:fr:thm:gauges}, says that the aligned
frames of a gauge map form a continuous section of the frame bundle, compatible with rigid
motions and with measure-preserving isometries. The alignment is the solution of an orthogonal
Procrustes problem for the covariance of a frame map with the gauge map
(\cref{m:fr:def:gauge}), and the estimates rest on the almost linearity of frame maps
(\cref{m:fr:thm:A,m:fr:lem:framegauge}).

For the rest of \cref{m:fr:sec:gauges}, let $(Z,d,m)$ be an $\RCD(\kappa,N)$ space with
$N<\infty$, and let $\Omega\subset Z$ be an open $\eps$-flat region at scale $r$
(\cref{m:fr:def:flat}). We take $\eps$ below a dimensional constant so that
\cref{m:fr:lem:flat,m:fr:thm:A} apply. Every point of $\Omega$ is frame-regular, so
$\Fr:=\pr^{-1}(\Omega)$ is a principal $O(n)$-bundle over $\Omega$. By the last assertion of \cref{m:fr:lem:flatness},
this setting includes the open Reifenberg subsets of $\RCD(\kappa,N)$ spaces at small
scales. For a point $x\in\Omega$ and a bounded Borel map $w$ into
$\R^n$ that is defined on $B_{2r}(x)$ we write
$\av_x^{(r)}(w):=\int w\,d\mu_x^{(r)}$ for its mean, again suppressing $(r)$ at fixed scale.
The alignment of a frame with a map is measured by a covariance.

\begin{definition}[Gauge Maps and Aligned Frames]\label{m:fr:def:gauge}
Let $\theta>0$ and let $U\subset Z$ be open.
\begin{enumerate}[label=\textup{(\alph*)}]
\item A Borel map $w\colon U\to\R^n$ is a \emph{gauge map with parameter $\theta$} if
$\bigl||w(p)-w(q)|-d(p,q)\bigr|\le\theta r$ for all $p,q\in U$ with $d(p,q)<6r$. The \emph{core} $U^\circ$ of $U$ is the open set of
all $x\in Z$ whose distance to $Z\setminus U$ is larger than $\tfrac52r$, with $U^\circ=Z$
if $U=Z$.
\item For $x\in Z$, a frame $e$ at $x$ and a bounded Borel map $w$ defined on
$B_{2r}(x)$, the \emph{covariance} is the matrix
\begin{equation*}
\Cov_x^{(r)}(e,w):=\frac1{r^2}\Bigl(\int_Zu_{e,a}^{(r)}\,
\bigl(w_b-\av_x^{(r)}(w_b)\bigr)\,d\mu_x^{(r)}\Bigr)_{a,b=1}^n .
\end{equation*}
The map $w$ is \emph{non-degenerate at $x$} if $\Cov_x^{(r)}(e,w)$ is invertible. The frame
\begin{equation*}
e_w^{(r)}(x):=e\cdot\operatorname{pol}\bigl(\Cov_x^{(r)}(e,w)\bigr),\qquad
\operatorname{pol}(C):=C\,(C^{\T}C)^{-1/2},
\end{equation*}
is then called the frame at $x$ \emph{aligned with $w$}.
\end{enumerate}
\end{definition}

At fixed $r$ we write $\Cov_x$ and $e_w$, again omitting $(r)$. Under
$(d,m,r,w)\mapsto(ad,bm,ar,aw)$ the covariance and aligned frame are unchanged.

Since $u_{e\cdot A}=A^{\T}u_e$, we have $\Cov_x(e\cdot A,w)=A^{\T}\Cov_x(e,w)$ for
$A\in O(n)$. Hence non-degeneracy does not depend on $e$, and $e_w(x)$ is the unique frame
at $x$ whose covariance with $w$ is symmetric and positive definite; in particular it
does not depend on $e$ either. Equivalently, $e_w(x)$ is the unique frame $e'$ at $x$ that
minimizes the distance in $L^2(\mu_x)$ between $w-\av_x(w)$ and $u_{e'}-\av_x(u_{e'})$: the
norm of the second map does not depend on $e'$, the inner product of the two maps is
$r^2$ times the trace of $\Cov_x(e',w)$, and the trace of $A^{\T}C$ over $A\in O(n)$ is
maximal exactly at the polar factor of an invertible matrix $C$, by the solution of the orthogonal
Procrustes problem \cite{Sch66}.

The next lemma compares gauge maps with frame maps. For a frame $e$ at $x$ we write $S_e$ for the matrix with the entries $\langle\cK_xe_a,\cK_xe_b\rangle_{L^2(\mu_x)}$, the second-moment matrix of $u_e/r$ with respect to $\mu_x$.

\begin{lemma}[Frame Maps and Gauge Maps]\label{m:fr:lem:framegauge}
For every $\eta\in(0,1)$ there are \mbox{$\eps_M=\eps_M(n,\eta)>0$} and \mbox{$\theta_M=\theta_M(n,\eta)>0$}
with the following property. Let $x$ be $\eps$-flat at scale $r$, with
\mbox{$\eps\le\eps_M$}, and let $e$ be a frame at $x$.
\begin{enumerate}[label=\textup{(\roman*)}]
\item The frame map $u_e$ satisfies $\bigl||u_e(p)-u_e(q)|-d(p,q)\bigr|\le\eta r$ for all
$p,q\in B_{3r}(x)$, so it is a gauge map with parameter $\eta$ on $B_{3r}(x)$, and
$|\av_x(u_e)|\le\eta r$.
\item The matrix $S_e$ is symmetric and $\lambda_n\le S_e\le4$.
\item If $w\colon B_{5r/2}(x)\to\R^n$ is a Borel map with
$\bigl||w(p)-w(q)|-d(p,q)\bigr|\le\theta_Mr$ for all points $p$ and $q$ of this ball, then
there are $A\in O(n)$ and $b\in\R^n$ with $|w-Au_e-b|\le\eta r$ on $B_{5r/2}(x)$.
\end{enumerate}
\end{lemma}
\begin{proof}
(i) Choose a Borel $\eps r$-GH approximation
\mbox{$\Phi\colon B_{4r}(x)\to B_{4r}(0^n)$} with $\Phi(x)=0$, using
\cref{m:fr:lem:flat}. By \cref{m:fr:thm:A},
if $\eps\le\eps_C(n,\tfrac14\eta)$, there is a matrix $R\in O(n)$ with
${|u_e-R\circ\Phi|\le\tfrac14\eta r}$ on $B_{3r}(x)$. Since $\eps\le\tfrac14\eta$, this
gives the first assertion. For the
second, $|\av_x(u_e)|\le\tfrac14\eta r+|\int\Phi\,d\mu_x|$, and the last term is at most
$\Psi(\eps|n)\,r$: along any sequence with $\eps\to0$, flatness and the continuous
cut-off $\psi$ imply that the push-forwards of $\mu_x$ under $\Phi$ converge to the
probability measure proportional to $\psi(|y|/r)\,dy$, whose mean is $0$.

(ii) Let $f_1,\dots,f_n$ be an orthonormal basis of $E_x$ consisting of eigenfunctions,
with eigenvalues $\lambda_k\in(\lambda_n,4]$, and write $e_a=\sum_k\beta_{ak}f_k$. With the
matrices $\beta=(\beta_{ak})$ and $\Lambda=\operatorname{diag}(\lambda_k)$, the frame
condition is $\beta\Lambda\beta^{\T}=\one$, so $\beta\Lambda^{1/2}$ is orthogonal, and
$S_e=\beta\Lambda^2\beta^{\T}$ is conjugate to $\Lambda$ by this orthogonal matrix.

(iii) Normalize $r=1$ and put $u=u_e$. Distances to $n+1$ well-positioned points
determine the alignment. By \cref{m:fr:thm:A} and the density of the GH approximation,
choose \mbox{$p_1,\dots,p_n\in B_{5/2}(x)$} so that the columns $u(p_a)$ of a matrix $U$
are within $\Psi(\eps|n)$ of the standard basis; also $u(x)=0$.
Put $b=w(x)$ and let $W$ have columns $w(p_a)-b$. Polarizing the distance
comparisons for $u$ and $w$, with $\delta:=\Psi(\eps|n)+\theta_M$, gives
\begin{align*}
\|W^{\T}W-U^{\T}U\|&\le C(n)\delta,\\
|W^{\T}(w(p)-b)-U^{\T}u(p)|&\le C(n)\delta
\qquad(p\in B_{5/2}(x)).
\end{align*}
For small $\delta$, both $U$ and $W$ have uniformly bounded inverses. Thus
$Q:=WU^{-1}$ satisfies $\|Q^{\T}Q-\one\|\le C(n)\delta$ and
$|w-b-Qu|\le C(n)\delta$ on the whole ball. Its polar part $A$ is orthogonal
and satisfies $\|A-Q\|\le C(n)\delta$. Choosing $\eps_M$ and $\theta_M$ small
enough and restoring $r$ proves (iii).
\end{proof}

We can now state the theorem on gauges: alignment with any gauge map produces a continuous
section of the frame bundle, which transforms under rigid motions of the gauge map and under
measure-preserving isometries as the frames do, and every frame is aligned with its own frame
map.

\begin{theorem}[Gauges]\label{m:fr:thm:gauges}
There are $\eps_G=\eps_G(n)>0$ and $\theta_G=\theta_G(n)>0$ with the following property. Let
$\Omega$ be $\eps$-flat at scale $r$ with $\eps\le\eps_G$, let $\theta\le\theta_G$, and let
$w\colon U\to\R^n$ be a gauge map with parameter $\theta$. Then $w$ is non-degenerate at every
$x\in U^\circ\cap\Omega$, and $x\mapsto(x,e_w(x))$ is a continuous section of $\Fr$ over
$U^\circ\cap\Omega$. For $R\in O(n)$, $b\in\R^n$ and a measure-preserving isometry $\varphi$
of $Z$, wherever the frames are defined,
\begin{equation*}
e_{Rw+b}(x)=e_w(x)\cdot R^{\T},\qquad
e_{w\circ\varphi^{-1}}(\varphi(x))=e_w(x)\circ\varphi^{-1}.
\end{equation*}
Every frame over $\Omega$ is aligned with its own frame map: $e_{u_e}(x)=e$ for every
frame $e$ at $x\in\Omega$.
\end{theorem}
\begin{proof}
Let $\eta\in(0,1)$ and suppose that $\eps\le\eps_M(n,\eta)$ and $\theta\le\theta_M(n,\eta)$, with
the constants of \cref{m:fr:lem:framegauge}. We show that the assertions hold when $\eta$ is
smaller than a dimensional constant; this proves the theorem.

\emph{Non-degeneracy.} Let $x\in U^\circ\cap\Omega$ and let $e$ be a frame at $x$.
The ball $B_{5r/2}(x)$ lies in $U$, so \cref{m:fr:lem:framegauge}(iii) gives $A\in O(n)$ and
$b\in\R^n$ with $w=Au_e+b+\xi$ and $|\xi|\le\eta r$ on $B_{5r/2}(x)$. The covariance is
linear in its second entry, and by the definitions
\begin{equation}\label{m:fr:eq:covself}
\Cov_x(e,u_e)=S_e-r^{-2}\,\av_x(u_e)\otimes \av_x(u_e) .
\end{equation}
This matrix is symmetric, and by \cref{m:fr:lem:framegauge} it is at least $\lambda_n-\eta^2$.
Since $|u_e|\le3r$ on $B_{5r/2}(x)$, we get
\begin{equation*}
\Cov_x(e,w)=\Cov_x(e,u_e)\,A^{\T}+\Cov_x(e,\xi)=(S_e+M)\,A^{\T},\qquad|M|\le C(n)\,\eta .
\end{equation*}
For small $\eta$ this matrix is invertible, so $w$ is non-degenerate at $x$.

\emph{Sections and identities.} Let $y\mapsto e^0(y)$ be a continuous local section of
$\Fr$ near $x$, which exists by \cref{m:fr:thm:bundle}. The frame maps $u_{e^0(y)}$ depend
continuously on $y$, uniformly on bounded sets, the densities of the measures $\mu_y$
with respect to $m$ depend continuously on $y$ in the supremum norm, and $w$ is
bounded on $B_{5r/2}(x)$. So $\Cov_y(e^0(y),w)$ is continuous in $y$ near $x$, and
${e_w(y)=e^0(y)\cdot\operatorname{pol}(\Cov_y(e^0(y),w))}$ is a continuous section. Next, ${\Cov_x(e,Rw+b)=\Cov_x(e,w)R^{\T}}$, so the covariance of $e_w(x)\cdot R^{\T}$ with
$Rw+b$
is the conjugate by $R$ of a positive definite symmetric matrix; this is the first
identity. An isometry $\varphi$ maps $\mu_x$ to $\mu_{\varphi(x)}$, and
$u_{e\circ\varphi^{-1}}=u_e\circ\varphi^{-1}$ by \cref{m:fr:thm:bundle}, so
$\Cov_{\varphi(x)}(e\circ\varphi^{-1},w\circ\varphi^{-1})=\Cov_x(e,w)$, which gives the
second identity. Finally \eqref{m:fr:eq:covself} is symmetric and positive definite, so $e_{u_e}(x)=e$.
\end{proof}

\subsubsection{Trivializations of the Frame Bundle}\label{m:fr:sec:frameenergy}
By \cref{m:fr:thm:gauges} the frame bundle is locally a product. We record the trivializations of $\Fr$ that this gives.

We keep the standing assumption and write $\Psi$ for functions $\Psi(\eps|n)$. We fix
$\eps_E=\eps_E(n)>0$ such that $\eps_E\le\eps_G$ and such that, for $\eps\le\eps_E$, the
frame maps are gauge maps with a parameter at most $\theta_G$, by \cref{m:fr:lem:framegauge};
in the rest of \cref{m:fr:sec:frameenergy} we assume $\eps\le\eps_E$. A gauge map $w\colon U\to\R^n$ with parameter at most $\theta_G$ gives the homeomorphism
\begin{equation*}
\ft_w\colon (U^\circ\cap\Omega)\times O(n)\longrightarrow\pr^{-1}(U^\circ\cap\Omega),\qquad\ft_w(y,A):=(y,e_w(y)\cdot A).
\end{equation*}

By \cref{m:fr:lem:framegauge} the frame maps $u_e$ of the frames $e$ at $x\in\Omega$ are gauge
maps on $B_{3r}(x)$ with a parameter $\Psi$, and the core of this ball contains
$B_{r/2}(x)$. Its trivialization is therefore defined over $B_{r/2}(x)\cap\Omega$.

\subsubsection{Convergence of Frames and Trivializations}\label{m:fr:sec:frameconv}
\Cref{m:fr:prop:convergence} gives the convergence of frames along subsequences.
Alignment gives the convergence of the sections of
\cref{m:fr:thm:gauges}, compatibly with the actions of isometries.

Fix $\eps_F=\eps_F(n)\le\eps_E$ so that frame maps at $\eps_F$-flat points have gauge
parameter at most $\tfrac12\theta_G$. Throughout this subsection and
\cref{m:fr:sec:funcconv}, $(Z_i,d_i,m_i)$ and $(Z,d,m)$ are $\RCD(\kappa,N)$ spaces,
and $(Z_i,z_i)\to(Z,z)$ with fixed Borel approximations $\af_i$ as in
\cref{m:fr:sec:conv}, along which
$\int g\circ\af_i\,dm_i\to\int g\,dm$ for every $g\in C_c(Z)$.
The open sets $\Omega_i\subset Z_i$ and $\Omega\subset Z$ are $\eps$-flat at scale $r$,
with $\eps\le\eps_F$, and every sequence $y_i\to y\in\Omega$ lies in $\Omega_i$ for
large $i$. We write $\Fr_i:=\pr^{-1}(\Omega_i)$, $\Fr:=\pr^{-1}(\Omega)$ and $u^i_e$ for
the frame maps on $Z_i$. Frames converge when their base points converge and their frame
maps converge uniformly, compared through the approximations.

\begin{definition}[Convergence of Frames]\label{m:fr:def:frameconv}
A sequence of frames $(y_i,e_i)\in\Fr_i$ \emph{converges} to $(y,e)\in\Fr$ if
\begin{equation}\label{m:fr:eq:frameconvdef}
\af_i(y_i)\to y\qquad\text{and}\qquad
\sup_{q\in B_{3r}(y_i)}\bigl|u^i_{e_i}(q)-u_e(\af_i(q))\bigr|\longrightarrow0 .
\end{equation}
\end{definition}

The limit is unique, because a frame is determined by the restriction of its frame map
to $B_{2r}(y)$. The next lemma shows that alignment is compatible with convergence.

\begin{lemma}[Convergence of Covariances]\label{m:fr:lem:covconv}
Let $(y_i,e_i)\in\Fr_i$ converge to $(y,e)\in\Fr$, let $w\colon B_{5r/2}(y)\to\R^n$ be
continuous and bounded, and let $w_i\colon B_{2r}(y_i)\to\R^n$ be Borel maps with
$\sup_{B_{2r}(y_i)}|w_i-w\circ \af_i|\to0$. Then $\Cov_{y_i}(e_i,w_i)\to\Cov_y(e,w)$.
\end{lemma}
\begin{proof}
Let $\nu_i$ be the image of the measure $\mu_{y_i}$ of $Z_i$ under $\af_i$. As in the proof
of \cref{m:fr:prop:convergence}, $\nu_i$ converges weakly to $\mu_y$, and for large $i$ the
measures $\nu_i$ are concentrated on a compact subset of $B_{5r/2}(y)$. By
\eqref{m:fr:eq:frameconvdef} and the hypothesis, the covariance $\Cov_{y_i}(e_i,w_i)$ differs
by a quantity that tends to $0$ from
\begin{equation*}
\frac1{r^2}\int u_e\otimes\Bigl(w-\int w\,d\nu_i\Bigr)\,d\nu_i ,
\end{equation*}
and this converges to $\Cov_y(e,w)$ because $u_e$ and $w$ are continuous.
\end{proof}

The convergence of covariances, together with the compactness of \cref{m:fr:prop:convergence},
gives the convergence theorem for the frame bundles: frames have convergent subsequences, every
limit frame is attained, aligned frames converge with their gauge maps, and the convergence is
compatible with convergent isometries. This is the form in which frames are used from
\cref{ch:frameconvergence} on.

\begin{theorem}[Convergence of Frame Bundles]\label{m:fr:thm:frameconv}
In the situation described above, the following holds.
\begin{enumerate}[label=\textup{(\roman*)}]
\item \emph{Subsequence compactness.} Every sequence of frames whose base points converge
to a point of $\Omega$ has a convergent subsequence.
\item \emph{Every frame is a limit.} If $(y,e)\in\Fr$ and $y_i\in Z_i$ with
$\af_i(y_i)\to y$, there are frames $e_i$ at $y_i$ for large $i$ such that $(y_i,e_i)$ converges to
$(y,e)$.
\item \emph{Aligned frames.} Let $w$ be a continuous gauge map on an open subset $U$ of $Z$
with a parameter $\theta<\theta_G$, let $y\in U^\circ\cap\Omega$ and $y_i\in Z_i$ with
$\af_i(y_i)\to y$, and let $w_i$ be Borel maps of $B_{5r/2}(y_i)$ to $\R^n$ such that the
functions ${|w_i-w\circ \af_i|}$, which are defined on this ball for large $i$, tend to $0$
uniformly on it. Then, for large $i$, the map $w_i$ is non-degenerate at $y_i$, and the frames
$(y_i,e_{w_i}(y_i))$ converge to $(y,e_w(y))$.
\item \emph{Isometries.} Let $\varphi_i$ be measure-preserving isometries of $Z_i$ such that
$d_i(\varphi_i(z_i),z_i)$ is bounded, and $\varphi$ a measure-preserving isometry of $Z$ such that
$d(\af_i\circ\varphi_i,\varphi\circ \af_i)\to0$ uniformly on $B_R(z_i)$ for every $R$. If $(y_i,e_i)$ converges to $(y,e)$, then
$(\varphi_i(y_i),e_i\circ\varphi_i^{-1})$ converges to $(\varphi(y),e\circ\varphi^{-1})$
whenever $\varphi(y)\in\Omega$.
\end{enumerate}
\end{theorem}
\begin{proof}
(i) Let $(y_i,e_i)$ be frames with $\af_i(y_i)\to y\in\Omega$. The assertion follows from
\cref{m:fr:prop:convergence}, applied to the maps $\af_i$ on $B_{4r}(y_i)$ with the value at
$y_i$ changed to $y$.

(ii) and (iii). Let $w$, $y$, $y_i$ and $w_i$ be as in (iii). For large $i$ the map $w_i$
changes distances on $B_{5r/2}(y_i)$ by less than $\theta_Gr$, so it is non-degenerate at
$y_i$, as in the proof of \cref{m:fr:thm:gauges}. Put $e_i:=e_{w_i}(y_i)$. Every subsequence
has a further subsequence along which $(y_i,e_i)$ converges to a frame $(y,e'')$. By
\cref{m:fr:lem:covconv} the covariances $\Cov_{y_i}(e_i,w_i)$, which are symmetric and
positive definite, converge to $\Cov_y(e'',w)$. This matrix is therefore symmetric and
positive semi-definite, and it is invertible because $w$ is non-degenerate at $y$. So $e''=e_w(y)$, and the whole sequence converges to $(y,e_w(y))$. This proves (iii). For
(ii) we apply (iii) to $w=u_e$, which is a continuous gauge map on
$B_{3r}(y)$ with $e_{u_e}(y)=e$, and to $w_i=u_e\circ \af_i$.

(iv) The frame map of $e_i\circ\varphi_i^{-1}$ is $u^i_{e_i}\circ\varphi_i^{-1}$. For
$q\in B_{3r}(\varphi_i(y_i))$ the point $\varphi_i^{-1}(q)$ lies in $B_{3r}(y_i)$, so by
\eqref{m:fr:eq:frameconvdef} the value $u^i_{e_i}(\varphi_i^{-1}(q))$ differs from
$u_e(\af_i(\varphi_i^{-1}(q)))$ by a quantity that tends to $0$ uniformly. The hypothesis,
applied at the point $\varphi_i^{-1}(q)$, shows that $\af_i(\varphi_i^{-1}(q))$ and
$\varphi^{-1}(\af_i(q))$ have distance tending to $0$, and $u_e$ is Lipschitz on $B_{4r}(y)$.
Since $u_e\circ\varphi^{-1}$ is the frame map of $e\circ\varphi^{-1}$, and
$\af_i(\varphi_i(y_i))\to\varphi(y)$, the assertion follows.
\end{proof}

\section{Convergence, Averaging and Free Frame Actions}\label{ch:frameconvergence}\label{m:book:ch:09}\label{m:fr:sec:mosco}

This section prepares the passage to the limit along a sequence of quotients
$Z_i/\Gamma_i$ that may collapse. We work over flat regions with arbitrary reference
measures, as in \cref{m:fr:sec:frameconv}, keeping support margins inside these regions.
Functions on the frame bundles of a convergent sequence are compared by continuous
convergence, and every continuous function with compact support over
$\Omega^{[r/2]}$ is the limit of a test sequence with
admissible supports (\cref{m:fr:lem:corrcharts}). Discrete groups of isometries
converge to a closed group, their normalized counting measures converge to a Haar
measure of the limit group (\cref{m:fr:lem:groupconv,m:fr:lem:counting}), and averages over
the groups preserve continuous convergence (\cref{m:fr:lem:groupaverages}), without a bound
on the quotient diameters. Finally, isometries act
freely on the frames over a region whose balls are close to Euclidean balls at all
scales up to eight times the frame scale (\cref{m:fr:prop:freeframes}). Throughout this section the spaces converge in the pointed measured Gromov--Hausdorff sense as in
\cref{m:fr:sec:frameconv}, and from \cref{m:fr:sec:quotientenergy} on groups of isometries
converge with them.

\subsection{Convergence of Functions on Frame Bundles}\label{m:fr:sec:funcconv}
We first say what it means that functions on the frame bundles $\Fr_i$ converge to a
function on $\Fr$. The definition uses only the convergence of frames
\eqref{m:fr:eq:frameconvdef}. In trivializations of the frame bundles it is continuous convergence of functions along a
convergent sequence of metric spaces.

Throughout this subsection we are in the situation of \cref{m:fr:sec:frameconv}.

It is convenient to realize the convergence in one space, as in the extrinsic approach
of \cite{GMS15}. Let $\cA$ be the disjoint
union of $Z$ and of all the spaces $Z_i$, with the distance $d_{\cA}$ that restricts to $d$ on $Z$
and to $d_i$ on $Z_i$, that is given by
\begin{equation*}
d_{\cA}(p,q):=\eps_i+\inf\bigl\{d_i(p,p')+d(\af_i(p'),q):p'\in B_{1/\eps_i}(z_i)\bigr\}
\qquad(p\in Z_i,\ q\in Z),
\end{equation*}
and that is the infimum of $d_{\cA}(p,q)+d_{\cA}(q,p'')$ over $q\in Z$ for points $p$ and $p''$ of
two different spaces $Z_i$ and $Z_j$. Since $\af_i$ changes distances by at most $\eps_i$,
this is a distance, $\cA$ is proper, and $d_{\cA}(p,\af_i(p))=\eps_i$ for
$p\in B_{1/\eps_i}(z_i)$. Thus points $q_i\in Z_i$ converge to $q\in Z$ in $\cA$ if and
only if $\af_i(q_i)\to q$, every point of $Z$ is such a limit, and every limit of points
$q_i\in Z_i$, with $i$ in an infinite set of indices, lies in $Z$. By the convergence of
the measures recalled in \cref{m:fr:sec:frameconv} and the uniform continuity of continuous
functions with compact support,
$\int_{Z_i}\varphi\,dm_i\to\int_Z\varphi\,dm$ for every such function $\varphi$ on
$\cA$. We use the same notions for the products $Z_i\times O(n)$ in $\cA\times O(n)$, and for frames, in
the sense of \eqref{m:fr:eq:frameconvdef}. Functions on the varying spaces are compared through
the convergence of points. The following definition is used for every convergence of points
that occurs in the paper: of points of $Z_i$ to points of $Z$, of frames
(\cref{m:fr:def:frameconv}), of group elements in \cref{m:fr:sec:quotientenergy}, and of points
of the frame spaces $Y_i$ (\cref{m:fr:def:Yconv}).

\begin{definition}[Continuous Convergence]\label{m:fr:def:contconv}
Let $\mathcal Q_i$, for $i\in\mathbb N$, and $\mathcal Q$ be sets with a notion of convergence $q_i\to q$ of points $q_i\in\mathcal Q_i$ to points $q\in\mathcal Q$. Functions $f_i\colon\mathcal Q_i\to\R$ \emph{converge continuously} to a function $f\colon\mathcal Q\to\R$ if
\begin{equation*}
f_i(q_i)\longrightarrow f(q)\qquad\text{whenever }q_i\to q .
\end{equation*}
\end{definition}

The functions are tested along all convergent sequences of points at once. In the realization
$\cA$, continuous convergence of functions $f_i$ on $Z_i$ to a continuous function $f$ on $Z$ says that the function equal to $f_i$ on $Z_i$ and to $f$ on $Z$ is continuous at every point of $Z$;
this is the picture behind the term. Every point of $Z$ and, by
\cref{m:fr:thm:frameconv}(ii), every frame in $\Fr$ is the limit of a sequence. So a sequence
that converges along an infinite set of indices can be completed to a convergent
sequence, continuous convergence can be used along subsequences, and the limit of a
continuously convergent sequence is continuous. All notions of this subsection are used also for subsequences of the spaces.

 The following margin observation keeps all bundle charts inside
the flat regions. Put $\Omega^{[0]}:=\Omega$ and $\Omega_i^{[0]}:=\Omega_i$.
For $s\ge0$, if $T\subset\Omega^{[s]}$ is compact,
there is $\eta>0$ such that, for large $i$, every point of $Z_i$ at $d_{\cA}$-distance at most
$\eta$ from $T$ in $\cA$ lies in $\Omega_i^{[s]}$. Choose $3\eta$ smaller than the
difference between $d(T,Z\setminus\Omega)$ and $s$. Otherwise points within $\eta$
of $T$ have points of $Z_i\setminus\Omega_i$ within $s+\eta$; a bounded subsequence
of the latter converges to a point of $\Omega$, contrary to the hypothesis on the
regions. This also covers $s=0$ by taking the point itself outside $\Omega_i$.

The supports of the functions that we compare are kept inside the regions by the following
condition.

\begin{definition}[Admissible Sets]\label{m:fr:def:admissible}
For $s>0$, closed sets $T_i\subset\Omega_i^{[s]}$ and a compact set
$T\subset\Omega^{[s]}$ are \emph{admissible with margin $s$} if $T_i\subset B_R(z_i)$
for one $R$ and every subsequential limit of points of $T_i$ belongs to $T$.
\end{definition}

Admissible sets eventually lie in every neighbourhood of $T$ in $\cA$. By the margin observation,
a sufficiently small closed neighbourhood of $T$ gives admissible enlargements of
these sets.

Trivializations reduce convergence of frames to convergence in such products.
Let \mbox{$(x_i,e_i)\to(x,e)$} with $x\in\Omega^{[r/2]}$. The closed ball of radius $r/2$
about $x$ lies in $\Omega$, so $B_{r/2}(x_i)\subset\Omega_i$ for large $i$ by the
margin observation. The frame maps $u^i_{e_i}$ and $u_e$
are gauge maps on $B_{3r}(x_i)$ and on $B_{3r}(x)$, and the cores of these balls contain
$B_{r/2}(x_i)$ and $B_{r/2}(x)$.

\begin{definition}[Corresponding Trivializations]\label{m:fr:def:corrtriv}
For frames $(x_i,e_i)\to(x,e)$ with $x\in\Omega^{[r/2]}$, the \emph{corresponding trivializations}
are the restrictions
\begin{equation*}
\ft_i\colon B_{r/2}(x_i)\times O(n)\to\pr^{-1}(B_{r/2}(x_i)),\qquad
\ft\colon B_{r/2}(x)\times O(n)\to\pr^{-1}(B_{r/2}(x))
\end{equation*}
of the trivializations of $u^i_{e_i}$ and of $u_e$ from \cref{m:fr:sec:frameenergy}.
\end{definition}

Test sequences are the functions along which we pass to the limit. We define them for any
convergence of points, as in \cref{m:fr:def:contconv}; on the frame bundles a support
condition keeps them inside the flat regions.

\begin{definition}[Test Sequences]\label{m:fr:def:testseq}
In the setting of \cref{m:fr:def:contconv}, with topologies on the sets $\mathcal Q_i$, a
\emph{test sequence} is a sequence of continuous functions $\Theta_i$ on $\mathcal Q_i$ that are bounded by one constant and converge continuously to a function $\Theta$ on $\mathcal Q$; the function $\Theta$ is the \emph{limit} of the test sequence.
\end{definition}

A test sequence on the frame bundles $\Fr_i$ has \emph{admissible supports} if its functions
vanish outside $\pr^{-1}(T_i)$ for sets $T_i$ and $T$ that are admissible with margin $r/2$
(\cref{m:fr:def:admissible}).
The limit of a test sequence on $\Fr_i$ is continuous, by what was said above. If the supports
are admissible, the limit $\Theta$ vanishes outside the compact set $\pr^{-1}(T)$, because
every frame over a point outside $T$ is the limit of frames over points outside the sets
$T_i$, as limits of points of these sets lie in $T$; so $\Theta$ has compact support in
$\pr^{-1}(\Omega^{[r/2]})$. Sums and products of test sequences are test sequences, with
admissible supports if the summands or factors have them.

\begin{lemma}[Corresponding Trivializations]\label{m:fr:lem:corrcharts}
In \textup{(i)} and \textup{(ii)}, $\ft_i$ and $\ft$ are corresponding trivializations over
$B_{r/2}(x_i)$ and $B_{r/2}(x)$.
\begin{enumerate}[label=\textup{(\roman*)}]
\item If points $(y_i,A_i)$ of $B_{r/2}(x_i)\times O(n)$ converge to a point $(y,A)$ of
$B_{r/2}(x)\times O(n)$, then the frames $\ft_i(y_i,A_i)$ converge to $\ft(y,A)$.
\item If frames $p_i\in\Fr_i$ converge to a frame $p$ with $\pr(p)\in B_{r/2}(x)$, then
$\pr(p_i)\in B_{r/2}(x_i)$ for large $i$, and $\ft_i^{-1}(p_i)$ converges to
$\ft^{-1}(p)$.
\item For every continuous function $\Theta$ with compact support in $\pr^{-1}(\Omega^{[r/2]})$
there is a test sequence with admissible supports and limit $\Theta$ whose functions have
compact support.
\end{enumerate}
\end{lemma}
\begin{proof}
(i) We apply \cref{m:fr:thm:frameconv}(iii) to the gauge map $w=u_e$ on $B_{3r}(x)$ and to
the maps $w_i=u^i_{e_i}$. For large $i$ the ball $B_{5r/2}(y_i)$ lies in $B_{3r}(x_i)$,
where $|w_i-w\circ \af_i|$ tends to $0$ uniformly by \eqref{m:fr:eq:frameconvdef}. Hence the
frames $\ft_i(y_i,\one)$, which are aligned with $w_i$, converge to $\ft(y,\one)$. Since
$u_{e'\cdot A}=A^{\T}u_{e'}$ for every frame $e'$, the frames
$\ft_i(y_i,A_i)=\ft_i(y_i,\one)\cdot A_i$ converge to $\ft(y,\one)\cdot A$.

(ii) The first assertion holds because $d_i(\pr(p_i),x_i)\to d(\pr(p),x)$. Write
$p_i=\ft_i(y_i,A_i)$ and $p=\ft(y,A)$. Then $y_i\to y$. If $A_i$ did not converge to
$A$, a subsequence would converge to some $A'\ne A$, and by (i) the frames $p_i$ would
converge along it to $\ft(y,A')\ne p$. This is impossible, because limits of frames are
unique.

(iii) Put $T^0:=\pr(\operatorname{supp}\Theta)\subset\Omega^{[r/2]}$. Choose
$x^1,\dots,x^m\in T^0$ such that the balls $B_{r/8}(x^k)$ cover $T^0$, and put
$\psi^k:=(\tfrac14r-d_{\cA}(\,\cdot\,,x^k))_+$ and
$\chi^k:=\psi^k/\max\{\tfrac18r,\sum_l\psi^l\}$, which are continuous functions on $\cA$.
Then $\chi^k$ vanishes outside the open ball of radius $\tfrac14r$ about $x^k$ in $\cA$,
and $\sum_k\chi^k=1$ on a neighbourhood of $T^0$ in $\cA$. By
\cref{m:fr:thm:frameconv}(ii) there are frames $(x^k_i,e^k_i)$ that converge to frames
$(x^k,e^k)$; let $\ft^k_i$ and $\ft^k$ be the corresponding trivializations. The function
$((\chi^k\circ\pr)\,\Theta)\circ\ft^k$, extended by $0$, is continuous on $Z\times O(n)$
and vanishes where $d(y,x^k)\ge\tfrac14r$. By the extension theorem of Tietze it is the
restriction of a bounded continuous function on $\cA\times O(n)$, and after
multiplication by a suitable function of the first variable we may assume that this
function, which we call $\vartheta^k$, vanishes at all points whose first component has
$d_{\cA}$-distance at least $\tfrac38r$ from $x^k$. Multiply it also by a continuous function of
the first variable that equals $1$ on $T^0$ and vanishes at $d_{\cA}$-distance at least $\eta$
from $T^0$, with $\eta$ small enough for the margin observation. Let $\Theta^k_i$ be the function on $\Fr_i$
that equals $\vartheta^k\circ(\ft^k_i)^{-1}$ on $\pr^{-1}(B_{r/2}(x^k_i))$ and $0$
elsewhere. For large $i$ it is continuous with compact support; for the remaining
indices we replace it by $0$. The closed $d_{\cA}$-neighbourhoods of radius $\eta$ of $T^0$ in $Z_i$
and $Z$ are admissible with margin $r/2$, so
$\Theta_i:=\sum_k\Theta^k_i$ has the required support. If frames $p_i$ converge to
$p$, then $\Theta^k_i(p_i)\to((\chi^k\circ\pr)\,\Theta)(p)$ for every $k$, by (ii) if
$\pr(p)\in B_{r/2}(x^k)$, and because both sides vanish for large $i$ otherwise. Since
$\sum_k\chi^k=1$ on a neighbourhood of $T^0$, summing over $k$ gives
$\Theta_i(p_i)\to\Theta(p)$. The functions $\Theta_i$ are bounded by
$\sum_k\sup|\vartheta^k|$, so they form a test sequence with admissible supports and limit $\Theta$.
\end{proof}

By \cref{m:fr:lem:corrcharts}(iii), every continuous function with compact support in
$\pr^{-1}(\Omega^{[r/2]})$ is the limit of a test sequence with admissible supports, and a test sequence that is
indexed by a subsequence can be completed to a test sequence, with admissible supports if it
has them.

\subsection{Convergence of the Groups and Averages}\label{m:fr:sec:quotientenergy}
We collect standard convergence and averaging tools for collapsing quotients
$X_i=Z_i/\Gamma_i$, with the uniform estimates needed below. Normalized counting
measures converge to Haar measure on the possibly positive-dimensional limit group
$G$, and cut-off functions allow averaging without choosing fundamental domains.

Throughout this subsection $(Z_i,z_i,\Gamma_i)\eGHto(Z,z,G)$ with
approximations $(\af_i,\alpha_i,\beta_i)$ as in \cref{m:fr:sec:conv}, where the spaces
$(Z_i,d_i,m_i)$ and $(Z,d,m)$ are $\RCD(\kappa,N)$, the maps $\af_i$ are Borel and the
reference measures converge along them. The groups $\Gamma_i$ are discrete groups of
measure-preserving isometries, which need not act freely, and $G$ is closed.
Where frames occur, we assume in addition invariant open sets $\Omega_i\subset Z_i$ and
$\Omega\subset Z$ as in \cref{m:fr:sec:frameconv}, flat at scale $r$.
The induced actions on the frame bundles are jointly continuous: isometries
converge uniformly on bounded sets, and frame maps are locally Lipschitz. The actions
are proper because the base actions are proper and the equivariant projections
have compact fibres. We write $\gamma_*$ and $g_*$ for the induced maps
$\varphi_*(x,e):=(\varphi(x),e\circ\varphi^{-1})$. No quotient-diameter bound is assumed.
Constants $C(R)$ depend only on $N$, $\kappa$ and $R$ and may change from line to line.
As in \cref{m:fr:sec:conv}, $\Gamma_i(R)$ and $G(R)$ are the
sets of elements that move the base point by less than $R$; the sets $\Gamma_i(R)$ are
finite, because discrete groups of isometries act properly.

The first lemma records consequences of the eGH convergence of \cite{FY92}; see also
\cite[Section~2.1]{Wan23}. The moving-element and frame formulations
are the ones used in \cref{m:fr:lem:limittransitions} and in the proof of \ref{fl:S:Yconv} in
\cref{m:fr:thm:geomain}.

\begin{lemma}[Convergence of Group Elements]\label{m:fr:lem:groupconv}
\begin{enumerate}[label=\textup{(\roman*)}]
\item Every sequence of elements $\gamma_i\in\Gamma_i$ with bounded displacements
$d_i(\gamma_iz_i,z_i)$ has a convergent subsequence.
\item If $\gamma_i\eGHto g$, then $\alpha_i(\gamma_i)\to g$ in $G$ and
$d_i(\gamma_iz_i,z_i)\to d(gz,z)$, and $g$ is the only limit of the sequence.
\item If $g^i\to g$ in $G$, then $\beta_i(g^i)\eGHto g$. In particular $\beta_i(g)\eGHto g$ for
every $g\in G$.
\item If $\gamma_i\eGHto g$ and $\gamma'_i\eGHto g'$, then $\gamma_i\gamma'_i\eGHto gg'$ and
$\gamma_i^{-1}\eGHto g^{-1}$.
\item The group $G$ preserves $m$. If $\gamma_i\eGHto g$ and points $p_i\in Z_i$ converge to $p\in Z$, then
$\gamma_ip_i\to gp$. If frames $p_i\in\Fr_i$ converge to $p\in\Fr$, then
$(\gamma_i)_*p_i\to g_*p$.
\end{enumerate}
\end{lemma}
\begin{proof}
On fixed bounded balls the approximation inequalities of \cref{m:fr:sec:conv}
apply for large $i$. They and the density of $\af_i$ give, for bounded-displacement
sequences, the equivalence $\gamma_i\eGHto g$ if and only if
$\alpha_i(\gamma_i)\to g$ uniformly on bounded sets. The images under $\alpha_i$ have
bounded displacement, so properness of the closed $G$-action gives (i).
The equivalence gives uniqueness in (ii), and evaluation at $z_i$ gives the
displacement assertion.

For (iii), the inequality for $\beta_i$ at $z_i$ bounds the displacements, and
\begin{equation*}
d\bigl(\af_i(\beta_i(g^i)q),g\,\af_i(q)\bigr)\le\eps_i+d\bigl(g^i\af_i(q),g\,\af_i(q)\bigr)\longrightarrow0
\end{equation*}
uniformly on bounded sets. Composition and inversion preserve this convergence:
for (iv), apply the defining estimates at $\gamma'_iq$ and $\gamma_i^{-1}q$,
respectively, which remain in bounded balls.

The point assertion in (v) follows from the definition and the isometry property.
For $g\in G$, the measure-preserving maps $\beta_i(g)$ converge to $g$; testing
against compactly supported continuous functions, whose pull-backs have uniformly
bounded support, gives $g_\#m=m$. The frame assertion now follows from
\cref{m:fr:thm:frameconv}(iv).
\end{proof}

Functions $\varphi_i$ on $\Gamma_i$ converge continuously to a function $\varphi$ on $G$ in the
sense of \cref{m:fr:def:contconv}, for the convergence $\gamma_i\eGHto g$ of group elements:
$\varphi_i(\gamma_i)\to\varphi(g)$ whenever $\gamma_i\eGHto g$. A sequence that
converges to $g$ along an infinite set of indices can be completed by the elements
$\beta_i(g)$ to a convergent sequence, by \cref{m:fr:lem:groupconv}(iii), so continuous
convergence can be used along subsequences. The limit $\varphi$ is continuous. Indeed, let
$g^k\to g$ in $G$. By \cref{m:fr:lem:groupconv}(iii) there are indices $i_1<i_2<\cdots$ with
$|\varphi_{i_k}(\beta_{i_k}(g^k))-\varphi(g^k)|<1/k$. Put $\gamma_i:=\beta_{i_k}(g^k)$ for
$i=i_k$ and $\gamma_i:=\beta_i(g)$ for the other indices. Then $\gamma_i\eGHto g$, again by
\cref{m:fr:lem:groupconv}(iii), so $\varphi_{i_k}(\beta_{i_k}(g^k))\to\varphi(g)$ and therefore
$\varphi(g^k)\to\varphi(g)$. If the functions $\varphi_i$ vanish outside $\Gamma_i(R)$,
then $\varphi$ vanishes outside the compact set of all $g$ with $d(gz,z)\le R$, by
\cref{m:fr:lem:groupconv}(ii) and~(iii). Continuity of each function alone does not give
continuous convergence.

\begin{example}[Identity Indicators]\label{m:fr:ex:identityindicator}
Let $\Gamma_i=i^{-1}\mathbb Z$ act on $\mathbb R$ by translations, so that
$(\mathbb R,0,\Gamma_i)\eGHto(\mathbb R,0,\mathbb R)$, and put
$f_i:=\mathbf1_{\{0\}}$ on $\Gamma_i$. Each $f_i$ is continuous because $\Gamma_i$ is
discrete, but both $0\eGHto0$ and $1/i\eGHto0$, whereas
\begin{equation*}
f_i(0)=1,\qquad f_i(1/i)=0.
\end{equation*}
Thus the functions $f_i$ have no continuously convergent limit, in particular not the
identity indicator on $G=\mathbb R$.
\end{example}

For $R>0$ let $\hat\chi_{R,i}:=\min\{1,(R+1-d_i(z_i,\cdot))_+\}$ on $Z_i$ and
$\hat\chi_R:=\min\{1,(R+1-d(z,\cdot))_+\}$ on $Z$. These functions are $1$-Lipschitz with
values in $[0,1]$, they equal $1$ on the balls of radius $R$ about the base points and
vanish outside the balls of radius $R+1$, and $\hat\chi_{R,i}$ converges continuously to
$\hat\chi_R$. We abbreviate $\hat\chi_i:=\hat\chi_{1,i}$ and $\hat\chi:=\hat\chi_1$,
and we count the groups with the weights
\begin{equation}\label{m:fr:eq:elldef}
\ell_i:=\sum_{\gamma\in\Gamma_i}\hat\chi_i(\gamma z_i).
\end{equation}
Thus $\#\Gamma_i(1)\le\ell_i\le\#\Gamma_i(2)$.

The next lemma is the normalized-counting version of the Haar-convergence
argument in \cite[Section~6.2, Propositions~6.6 and~6.8]{Cav23}.
Those propositions concern proper geodesic spaces and closed isometry groups,
without a CAT(0) or cocompactness assumption. Its multiplicity argument is used again in
\cref{th:lem:weights}(ii), and its normalization and invariance statements in
\cref{m:fr:prop:limitgroup,m:fr:ext:physical-good-model}.

\begin{lemma}[Normalization of the Counting Measures]\label{m:fr:lem:counting}
\begin{enumerate}[label=\textup{(\roman*)}]
\item For every $R\ge1$ and all $i$,
\begin{equation*}
\#\Gamma_i(R)\le C(R)\,\#\Gamma_i(1)\le C(R)\,\ell_i .
\end{equation*}
\item There is exactly one Haar measure $\nu_G$ on $G$ with
$\int_G\hat\chi(gz)\,d\nu_G(g)=1$. It is invariant under left and right translations and
under inversion, and
\begin{equation}\label{m:fr:eq:haarlimit}
\frac1{\ell_i}\sum_{\gamma\in\Gamma_i}\varphi_i(\gamma)\longrightarrow\int_G\varphi\,d\nu_G
\end{equation}
for all functions $\varphi_i$ on $\Gamma_i$ that are bounded by one constant, vanish
outside $\Gamma_i(R)$ for one $R$, and converge continuously to a function $\varphi$ on
$G$. In particular $G$ is unimodular.
\item For every $R>0$, all $i$ and all $q\in Z_i$, the number of elements
$\gamma\in\Gamma_i$ with $\gamma q\in B_R(z_i)$ is at most $C(R)\,\ell_i$. For every
$q\in Z$ the set of all $g\in G$ with $gq\in B_R(z)$ has $\nu_G$-measure at most $C(R)$.
\end{enumerate}
\end{lemma}
\begin{proof}
(i) The balls $B_{1/2}(\gamma z_i)$ with $\gamma\in\Gamma_i(R)$ lie in $B_{R+1/2}(z_i)$, and
each of them has the measure of $B_{1/2}(z_i)$. If two of them, for $\gamma$ and
$\gamma'$, intersect, then $\gamma^{-1}\gamma'\in\Gamma_i(1)$; so every point lies in the
balls of at most $\#\Gamma_i(1)$ elements $\gamma$. Hence
$\#\Gamma_i(R)\,m_i(B_{1/2}(z_i))\le\#\Gamma_i(1)\,m_i(B_{R+1/2}(z_i))$, and the first
inequality follows from the Bishop--Gromov inequality. The second one holds because
$\hat\chi_i(\gamma z_i)=1$ for $\gamma\in\Gamma_i(1)$.

(ii) Push the normalized counting measure on $\Gamma_i(1/\eps_i)$ to $G$ by
$\alpha_i$, obtaining $\nu_i$. The displacement error is at most $2\eps_i$, so
(i) gives $\nu_i(G(R))\le C(R+1)$ for large $i$. Properness therefore gives
subsequential vague limits $\nu$ on $G$.

For tests $\varphi_i\to\varphi$ as in (ii), compactness in
\cref{m:fr:lem:groupconv}(i)--(ii) gives
$\sup_{\Gamma_i(R+1)}|\varphi_i-\varphi\circ\alpha_i|\to0$:
any violating sequence has a convergent subsequence on which both terms tend to
the same value. All terms in the comparison vanish outside $\Gamma_i(R+1)$ for large $i$, so (i)
and vague convergence prove \eqref{m:fr:eq:haarlimit} along the subsequence.
The normalizing test $\hat\chi_i(\gamma z_i)$ gives
$\int\hat\chi(gz)\,d\nu(g)=1$, hence $\nu\ne0$.

For $\varphi\in C_c(G)$ pull it back by $\alpha_i$, extending by zero outside
$\Gamma_i(1/\eps_i)$. For each fixed $h\in G$, its left and right translates by
$\beta_i(h)$ and its composition with inversion have common bounded-displacement
support and converge continuously to the corresponding tests on $G$, by
\cref{m:fr:lem:groupconv}(ii)--(iv). Their counting sums are equal. Passing to
the limit proves left and right invariance and inversion invariance of $\nu$.
Thus $\nu$ is Haar; its normalization makes it unique, proving convergence of the
whole sequence and (ii).

(iii) If no element sends $q$ into $B_R(z_i)$ the count is zero. Otherwise choose one,
$\gamma_0$. Every other such $\gamma$ satisfies
$d_i(\gamma\gamma_0^{-1}z_i,z_i)<2R$, so the count is at most
$\#\Gamma_i(2R)\le C(R)\ell_i$ by (i), using radius $\max\{1,2R\}$ if needed.
The same argument on $Z$ puts the non-empty
set of contributing elements in $G(2R)g_0$. Its measure is at most $\nu_G(G(2R))$,
by right invariance. Applying \eqref{m:fr:eq:haarlimit} to
$\gamma\mapsto\hat\chi_{2R,i}(\gamma z_i)$ and using (i) bounds this by $C(R)$.
\end{proof}

Put
\begin{equation*}
\varrho_i(q):=\min_{\gamma\in\Gamma_i}d_i(\gamma q,z_i),\qquad
\varrho(q):=\min_{g\in G}d(gq,z),\qquad
Z_i(R):=\{\varrho_i<R\},\quad Z(R):=\{\varrho<R\}.
\end{equation*}
The minima exist by properness. We use a regional version of the cut-off functions
for proper actions described in \cite{Tu99survey}. The construction
below gives, beyond existence, the uniform support, Lipschitz and convergence estimates that
the proof of the last assertion of \ref{fl:S:Yconv} in \cref{m:fr:thm:geomain} uses together with
\cref{m:fr:lem:groupaverages}. Choose
$c_{\chi}=c_{\chi}(N,\kappa)>0$ with $c_{\chi}\le\#\Gamma_i(1)/\ell_i$, using
\cref{m:fr:lem:counting}(i). For $R>0$ put
\begin{equation*}
\bar\chi_{R,i}(q):=\frac1{\ell_i}\sum_{\gamma\in\Gamma_i}\hat\chi_{R+1,i}(\gamma q),
\qquad
\chi_{R,i}:=\frac{\hat\chi_{R+1,i}}{\max\{\bar\chi_{R,i},c_{\chi}\}},
\end{equation*}
and define $\bar\chi_R,\chi_R$ on $Z$ by replacing the normalized sum with the integral
against $\nu_G$.

\begin{lemma}[Cut-Off Functions for the Groups]\label{m:fr:lem:groupunity}
For each $R>0$, the following holds, with bounds depending only on $N,\kappa,R$.
\begin{enumerate}[label=\textup{(\roman*)}]
\item The invariant functions $\bar\chi_{R,i}$ are uniformly bounded and Lipschitz,
converge continuously to $\bar\chi_R$, and are at least $c_{\chi}$ on
$\{\varrho_i\le R\}$; the same assertions hold on the limit.
\item The functions $\chi_{R,i}$ and $\chi_R$ take values in $[0,1/c_{\chi}]$, are positive
on the balls of radius $R+2$ about the base points and vanish outside these balls.
They are uniformly Lipschitz, $\chi_{R,i}$ converges continuously to $\chi_R$, and
\begin{equation}\label{m:fr:eq:regionalunity}
\frac1{\ell_i}\sum_{\gamma\in\Gamma_i}\chi_{R,i}(\gamma q)=1\quad(q\in Z_i(R)),
\qquad
\int_G\chi_R(gq)\,d\nu_G(g)=1\quad(q\in Z(R)).
\end{equation}
The averages are at most $1$ everywhere, and equal $1$ also when $\varrho_i=R$,
respectively $\varrho=R$.
\end{enumerate}
\end{lemma}
\begin{proof}
The orbit-count bounds of \cref{m:fr:lem:counting}(iii) bound the number, or Haar measure,
of non-zero terms by $C(R)$ after normalization. Since each cut-off function is $1$-Lipschitz,
applying this bound at two points gives the uniform bound and Lipschitz estimate for
$\bar\chi_{R,i}$ and $\bar\chi_R$. Invariance follows from right invariance of the
counting and Haar measures. If $\varrho_i(q)\le R$, choose $\gamma_0$ with
$d_i(\gamma_0q,z_i)\le R$. For $\gamma\in\Gamma_i(1)$ we have
$\hat\chi_{R+1,i}(\gamma\gamma_0q)=1$, so $\bar\chi_{R,i}(q)\ge c_{\chi}$.
For $q_i\to q$, the functions $\gamma\mapsto\hat\chi_{R+1,i}(\gamma q_i)$ are uniformly
bounded, vanish outside $\Gamma_i(R')$ for one $R'$ and converge continuously to
$g\mapsto\hat\chi_{R+1}(gq)$. Thus \eqref{m:fr:eq:haarlimit} gives continuous
convergence. If $\varrho(q)\le R$, choose $q_i\to q$ and an element $g$ realizing
$\varrho(q)$. Then $\beta_i(g)q_i\to gq$; the terms indexed by
$\Gamma_i(1)\beta_i(g)$ in the average are at least $1-o(1)$, proving
$\bar\chi_R(q)\ge c_{\chi}$. Finally, the denominator defining $\chi_{R,i}$ is invariant,
so its group average is
$\bar\chi_{R,i}/\max\{\bar\chi_{R,i},c_{\chi}\}$, and likewise on the limit.
This proves (ii), including the bounds and convergence.
\end{proof}

Functions with bounded support on the frame bundles give functions on the quotient
bundles by averaging over the groups. For a Borel function $U$ on $\Fr$, respectively
on $\Fr_i$, that vanishes outside the preimage of a bounded set, we put
\begin{equation}\label{m:fr:eq:groupaverage}
\bar U(p):=\int_GU(g_*^{-1}p)\,d\nu_G(g),\qquad\text{respectively}\qquad
\bar U(p):=\frac1{\ell_i}\sum_{\gamma\in\Gamma_i}U(\gamma_*^{-1}p),
\end{equation}
at the points $p$ at which the integral of $|U|$ in place of $U$ is finite. The sums are
finite.

\begin{lemma}[Averages over the Groups]\label{m:fr:lem:groupaverages}
Let $R>0$, and let $U$ be a continuous function with compact support on $\Fr$
that vanishes outside $\pr^{-1}(B_R(z))$.
\begin{enumerate}[label=\textup{(\roman*)}]
\item The function $\bar U$ is defined everywhere and continuous, it satisfies
$\bar U\circ h_*=\bar U$ for every $h\in G$, and $\sup|\bar U|\le C(R)\sup|U|$.
\item The same holds on $\Fr_i$, with $z_i$ and $\Gamma_i$. If $\Theta_i$ is a test sequence with
admissible supports and limit $\Theta$, then the functions $\bar\Theta_i$ converge continuously to
$\bar\Theta$.
\end{enumerate}
\end{lemma}
\begin{proof}
(i) By \cref{m:fr:lem:counting}(iii) and inversion invariance, the contributing
elements in \eqref{m:fr:eq:groupaverage} have measure at most $C(R)$, giving the
bound. Left invariance gives $\bar U\circ h_*=\bar U$. As $p$ varies near a fixed
frame, all contributions lie in one compact subset of $G$; joint continuity of
the action then gives continuity of $\bar U$. The same proof applies to the
normalized sums on $\Gamma_i$.

(ii) Let $\Theta_i$ be a test sequence with admissible supports and limit $\Theta$, and let
frames $p_i\in\Fr_i$
converge to $p\in\Fr$. By \cref{m:fr:lem:groupconv}(iv) and~(v), the functions
$\gamma\mapsto\Theta_i(\gamma_*^{-1}p_i)$ on $\Gamma_i$ converge continuously to
$g\mapsto\Theta(g_*^{-1}p)$. They are bounded by one constant and vanish outside
$\Gamma_i(R')$ for one $R'$, so $\bar\Theta_i(p_i)\to\bar\Theta(p)$ by
\eqref{m:fr:eq:haarlimit}.
\end{proof}
\subsection{Flat Regions}\label{m:fr:sec:regions}
The preceding results were proved directly over flat regions with arbitrary reference
measures. We collect their scope for the local covers of \cref{m:book:ch:10}, then
specialize to the invariant regions defined by quotient distance.

\begin{proposition}[Localization]\label{m:fr:prop:localization}
Let $\eps_F=\eps_F(n)$ be the constant fixed in \cref{m:fr:sec:frameconv}.
\begin{enumerate}[label=\textup{(\roman*)}]
\item \emph{Frames and gauges.} Let $(Z,d,m)$ be $\RCD(\kappa,N)$ and let
$\Omega\subset Z$ be open and $\eps$-flat at scale $r$, with $\eps\le\eps_F$.
Then $\Fr(\Omega):=\pr^{-1}(\Omega)$ is a principal $O(n)$-bundle.
\Cref{m:fr:lem:framegauge,m:fr:thm:gauges} apply at the points of $\Omega$, with
gauge cores intersected with $\Omega$ and measure-preserving isometries.
\item \emph{Regions and test sequences.} Let $(Z_i,d_i,m_i,z_i)\to(Z,d,m,z)$
and open flat regions $\Omega_i,\Omega$ satisfy the measured convergence and
eventual-inclusion hypotheses of \cref{m:fr:sec:frameconv}.
Every continuous function with compact support in $\pr^{-1}(\Omega^{[r/2]})$
is the limit of a test sequence with admissible supports.
\item \emph{Convergence of frames.} In \textup{(ii)}, \cref{m:fr:thm:frameconv}
holds over $\Omega_i,\Omega$, and \cref{m:fr:lem:corrcharts}(i), (ii) holds for
corresponding trivializations at $x_i\to x$ with
$d(x,Z\setminus\Omega)\ge r$.
\end{enumerate}
\end{proposition}
\begin{proof}
Part (i) is \cref{m:fr:thm:A,m:fr:thm:bundle,m:fr:lem:flat} and the local gauge
construction of \cref{m:fr:sec:gauges}. Parts (ii) and (iii) are
\cref{m:fr:lem:corrcharts,m:fr:thm:frameconv}; the stated distance condition places
$x$ in $\Omega^{[r/2]}$, as required for corresponding trivializations.
\end{proof}

In the group setting of \cref{m:fr:sec:quotientenergy}, the quotient-distance functions
$\varrho_i$ converge continuously to $\varrho$. Indeed, if $q_i\to q$ and $g$ realizes
$\varrho(q)$, then $\beta_i(g)q_i\to gq$ by \cref{m:fr:lem:groupconv}, giving the upper
limit inequality. Elements $\gamma_i$ realizing $\varrho_i(q_i)$ have bounded
displacements, so subsequences converge to elements $h$ with
$d(hq,z)=\lim\varrho_i(q_i)$, giving the lower limit inequality.

Fix $R_\Omega\in(0,\infty]$ and suppose that $\Omega_i:=Z_i(R_\Omega)$ and $\Omega:=Z(R_\Omega)$ are $\eps$-flat at scale $r$, with
$\eps\le\eps_F$. Continuous convergence of $\varrho_i$ gives the eventual-inclusion
hypothesis on these regions. Thus the frame convergence, test sequences and averages
of the preceding subsections apply to $\Fr(\Omega_i):=\pr^{-1}(\Omega_i)$ and $\Fr(\Omega):=\pr^{-1}(\Omega)$. Their normalization is always
$\ell_i=\sum_\gamma\hat\chi_{1,i}(\gamma z_i)$, and the cut-off functions $\chi_{R,i},\chi_R$ of
\cref{m:fr:lem:groupunity} satisfy \eqref{m:fr:eq:regionalunity}. In particular the averages of test sequences
with admissible supports are continuous invariant functions, uniformly bounded, and converge continuously to
the average of the limit.

The actions on these frame bundles are proper but need not be free
(\cref{m:fr:ex:notfree}). For a frame $e$ at $q$, the relation
$e\circ g^{-1}=e\cdot A_e(g)$ defines a homomorphism $G_q\to O(n)$, which need not
be injective. Freeness requires the additional all-scales condition below.

\subsection{Free Actions on Frames}\label{m:fr:sec:freeframes}
For the flat regions of \cref{m:fr:sec:regions} the actions need not be free, as the
example below shows.
In this subsection we prove that they are free when, in addition, the balls of all
radii up to $8r$ about the points of the region are close to Euclidean balls. This is a
condition on distances only, so it passes to Gromov--Hausdorff limits, with a loss in
the constant; it holds for the limits of the local covers in \cref{m:fr:sec:covercollapse}.
The idea is the following. The isometries that fix a given frame form a group, and
each of its elements moves the points near the base point by a small amount; a group
of isometries of a Euclidean ball with this property is trivial.

\begin{example}[Frames That Are Fixed]\label{m:fr:ex:notfree}
Let $S$ be a round two-dimensional sphere of diameter at most $\eps r$, and let
$Z:=\R^n\times S$ with the product distance and the product measure, an
$\RCD(0,n+2)$ space. The projection to $\R^n$ shows that $Z$ is $C(n)\,\eps$-flat at scale
$r$, so for small $\eps$ its points are frame-regular. Let $G$ be the group of the
rotations of $S$ about an axis, acting on the second factor, and let $q$ be a point
whose second component lies on the axis. For a frame $e$
at $q$ and $g\in G$ we have $e\circ g^{-1}=e\cdot A(g)$ with a homomorphism
$A\colon G\to O(n)$, and $u_e\circ g^{-1}=A(g)^{\T}u_e$. By \cref{m:fr:thm:A} the frame map
$u_e$ is uniformly close on $B_{3r}(q)$ to the composition of the projection with a
linear isometry, and the projection is invariant under $G$. Testing at $n$ points whose
projected displacements from $q$ are $r$ times an orthonormal basis gives
$\|A(g)-\one\|\le C(n)\Psi(\eps|n)\le\tfrac12$
for all $g\in G$ when $\eps$ is small. A subgroup of $O(n)$ with this property is
trivial, because a non-trivial orthogonal matrix has a power that maps some unit vector
to a vector at distance at least $\sqrt2$. Hence every element of $G$ fixes every frame
at $q$: the action of $G$ on the frame bundle is not free, although only the identity
fixes a non-empty open subset of $Z$. The balls of $Z$ whose radii are smaller than the
diameter of $S$ are not close to balls of $\R^n$.
\end{example}

The next lemma excludes this behaviour when all small balls are close to Euclidean balls.
It is a statement about metric spaces: the absence of small subgroups in the isometry
group of $\R^n$, transported through the GH approximations. The same Euclidean fact,
reached by blowing up at points with Euclidean tangent cones, excludes small subgroups
of the isometry groups of $\RCD$ spaces in \cite{Sos18,GSR19}.

\begin{lemma}[Groups of Isometries with Small Displacements]\label{m:fr:lem:smallgroups}
For every $n\ge1$ there is $\delta_G=\delta_G(n)>0$ with the following property. Let
$(Z,d)$ be a metric space, $a\in Z$ and $t>0$, and suppose that $B_{8t}(a)$ is
$8\delta_Gt$-close to $B_{8t}(0^n)$. Let $H$ be a group of isometries of $Z$ such that
$d(gy,y)\le\tfrac14t$ for all $g\in H$ and all $y\in B_t(a)$. Then $d(gy,y)\le\tfrac18t$
for all $g\in H$ and all $y\in B_t(a)$.
\end{lemma}
\begin{proof}
We may assume that $t=1$. Suppose that the assertion fails for a sequence
$\delta_j\to0$ in place of $\delta_G$, with spaces $Z_j$, points $a_j$ and groups $H_j$.
Put $\eta_j:=8\delta_j$, let $\Phi_j\colon B_8(a_j)\to B_8(0^n)$ be $\eta_j$-GH
approximations, so that $|\Phi_j(a_j)|\le2\eta_j$ by \cref{m:fr:lem:centred}, and
choose maps $\Phi'_j\colon B_8(0^n)\to B_8(a_j)$ with $|\Phi_j(\Phi'_j(v))-v|\le\eta_j$. Then
$\Phi'_j$ changes distances by at most $3\eta_j$, and $d(\Phi'_j(\Phi_j(b)),b)\le2\eta_j$ for
all $b$. An element $g$ of $H_j$ moves $a_j$ by at most $\tfrac14$, so it
maps $B_s(a_j)$ into $B_{s+1/4}(a_j)$ for every $s$. For $g\in H_j$ we put
$\hat g:=\Phi_j\circ g\circ\Phi'_j$ on $B_4(0^n)$, which is possible for large $j$. Since
$g$ preserves distances, for large $j$ and all $g,g'\in H_j$ we have
\begin{gather*}
\bigl||\hat g(v)-\hat g(v')|-|v-v'|\bigr|\le4\eta_j\quad(v,v'\in B_4(0^n)),\qquad
|\hat g(v)-v|\le\tfrac14+2\eta_j\quad(|v|<1-4\eta_j),\\
\bigl|\widehat{gg'}(v)-\hat g(\hat g'(v))\bigr|\le3\eta_j\quad(v\in B_2(0^n)).
\end{gather*}
By assumption there are $g_j\in H_j$ and $y_j\in B_1(a_j)$ with
$d(g_jy_j,y_j)>\tfrac18$, and for $v_j:=\Phi_j(y_j)$ this gives
$|\hat g_j(v_j)-v_j|\ge\tfrac18-4\eta_j$. The maps $\hat g_j$ are bounded. By a diagonal
argument a subsequence converges on a countable dense subset of $B_4(0^n)$, and by the
first inequality the limit preserves distances and the convergence is uniform on
$B_4(0^n)$. So the limit is the restriction of an isometry $E$ of $\R^n$.
Passing to a further subsequence, $v_j$ converges to a point $v$ with $|v|\le1$ and
$|E(v)-v|\ge\tfrac18$; so $E$ is not the identity. For $v\in B_2(0^n)$ and $g\in H_j$ the
point $\hat g(v)$ lies in $B_4(0^n)$. So by the third inequality and induction on $k$,
the maps $\widehat{g_j^k}$ converge to $E^k$ uniformly on $B_2(0^n)$, for every $k\ge1$.
Since the powers of $g_j$ belong to $H_j$, the second inequality gives
\begin{equation*}
|E^k(v)-v|\le\tfrac14\qquad(|v|<1,\ k\ge1).
\end{equation*}

Write $E(v)=R_Ev+b_E$ and $E^k(v)=R_E^kv+b_k$. Applying the bound at $v$ and $-v$ gives
$2|(R_E^k-\one)v|\le1/2$, hence $\|R_E^k-\one\|\le1/4$ for every $k$.
A non-identity orthogonal matrix has an eigenvalue $-1$ or a non-trivial rotation
plane; some power then moves a unit vector by at least $\sqrt2$. Thus $R_E=\one$.
Now $b_k=kb_E$, and the bound at $v=0$ forces $b_E=0$. This contradicts $E\ne\mathrm{id}$.
\end{proof}

We now combine the lemma with the frames. The proof moves the estimate of the lemma
from the scale $r$ down to all smaller scales.

\begin{proposition}[Free Actions on Frames]\label{m:fr:prop:freeframes}
There are $\eps_P=\eps_P(n)>0$ and $\delta_P=\delta_P(n)>0$ with the following property.
Let $(Z,d,m)$ be an $\RCD(\kappa,N)$ space with $N<\infty$, let $G$ be a group of isometries
of $Z$ that preserve $m$, and let $\Omega\subset Z$ be an open set that is $\eps$-flat at
scale $r$, with $\eps\le\eps_P$, and $(\delta_P,8r)$-Reifenberg
(\cref{m:fr:def:reifenberg}). Let $q\in\Omega$, and let $e$ be a frame at
$q$. Then the identity is the only element $g$ of $G$ with $g(q)=q$ and
$e\circ g^{-1}=e$. In particular, if $\Omega$ is invariant under $G$, then $G$ acts
freely on $\Fr(\Omega)$, and the homomorphisms $A_e$ of the isotropy groups of the
points of $\Omega$ into $O(n)$ are injective.
\end{proposition}
\begin{proof}
Let $\delta_P:=\delta_G$, and let $\eps_P\le\min\{\eps_F,\tfrac1{16}\}$ be so small that \cref{m:fr:thm:A} applies with $\eta=\tfrac1{32}$ at $\eps_P$-flat points,
which are almost Euclidean by \cref{m:fr:lem:flat}. Then the frame map $u_e$ is within $\tfrac1{32}r$ of $R\circ\Phi$ on
$B_{3r}(q)$, for an orthogonal matrix $R$ and an $\eps r$-GH approximation $\Phi$, so it
changes the distances between the points of $B_{3r}(q)$ by at most $\tfrac18r$. Let $H$
be the set of all
$g\in G$ with $g(q)=q$ and $e\circ g^{-1}=e$. It is a group. For $g\in H$ we have
$u_e\circ g^{-1}=u_{e\circ g^{-1}}=u_e$ by \cref{m:fr:thm:bundle}, and $g$ maps $B_{3r}(q)$ onto
itself. So $|u_e(gy)-u_e(y)|=0$, and therefore
\begin{equation*}
d(gy,y)\le\tfrac18r\qquad(g\in H,\ y\in B_{3r}(q)).
\end{equation*}
Fix a point $a$ of $B_r(q)\cap\Omega$, and for $0<t\le r$ let $D(t)$ be the supremum of
$d(gy,y)$ over $g\in H$ and $y\in B_t(a)$. Then $D$ is non-decreasing, and
$D(r)\le\tfrac18r$, because $B_r(a)\subset B_{3r}(q)$. If $D(t)\le\tfrac18t$ for some
$t\le r$, then $D(\tfrac12t)\le\tfrac14\cdot\tfrac12t$, and \cref{m:fr:lem:smallgroups}, applied
to the ball $B_{t/2}(a)$ and the group $H$, gives $D(\tfrac12t)\le\tfrac18\cdot\tfrac12t$;
the ball $B_{4t}(a)$ is $4\delta_Pt$-close to $B_{4t}(0^n)$, as required, because $\Omega$
is $(\delta_P,8r)$-Reifenberg, $a\in\Omega$ and $4t\le8r$. By induction, $D(2^{-k}r)\le2^{-k}\,\tfrac18r$ for all $k\ge0$. Hence
$d(ga,a)=0$ for every $g\in H$. Thus the elements of $H$ fix the non-empty open set
$B_r(q)\cap\Omega$ pointwise, and they are the identity by \ref{m:fr:R:isometries}. The last
assertions follow from the frame action. Explicitly, for $g\in G_q$ the identity
$e\circ g^{-1}=e\cdot A_e(g)$ defines $A_e(g)$ uniquely. Composing two isometries and
using commutation with the right $O(n)$ action gives $A_e(gg')=A_e(g)A_e(g')$.
Continuity follows from that of the frame action, and its kernel is trivial by
what was just proved.
\end{proof}

The proposition applies to the regions of \cref{m:fr:sec:regions} when they are also
$(\delta_P,8r)$-Reifenberg.

\section{Compatible Local Covers and the Framed Limit}\label{ch:coveratlas}

We derive a compatible atlas and its framed limit from the local-cover hypothesis of \cref{m:fr:def:localcovers}. The result is \cref{m:fr:thm:geomain}, stated in \cref{fl:sec:setting}. Its assertions are the geometric data used throughout the bundle construction. \Cref{m:book:ch:10} constructs the atlas and its limit, \cref{m:book:ch:11} describes the limit groups, and \cref{fl:sec:tools} derives two tools for the later sections.

\subsection{The Atlas and the Framed Limit}\label{fl:sec:setting}\label{fl:sec:charts}

We state the theorem first. Its objects are constructed in \cref{m:book:ch:10}, and on a first
reading it may be taken as a summary of that construction. Each assertion is proved where its
objects are constructed, and the paragraph after the theorem says where.

Use the radii $R_0,R_1$ of \cref{m:fr:sec:localcovers}, the frame spaces $Y_i=\Fr_\rho(X_i)$ of
\cref{m:fr:sec:coverbundle}, and the finitely many centres and normalized charts of
\cref{m:fr:sec:covercollapse}. Suppress the centre superscript when working in one chart. With
the invariant distance $\varrho_i$ from the base orbit and the sets $Z_i(R)$ on which
$\varrho_i<R$, put
\[
\Omega_i=Z_i(R_0)\supset\Omega'_i=Z_i(R_0/2)
 \supset\Omega_i^\flat:=Z_i(R_0/4),
\]
and use the analogous notation on a limit chart. Since $\varrho_i$ is
$1$-Lipschitz, a ball of radius $R_0/4$ about a point of $\Omega_i^\flat$ lies in
$\Omega'_i$, and likewise in the limit. The frame scale obeys $128r\le R_1$.

\begin{theorem}[The Atlas and the Framed Limit]\label{m:fr:thm:geomain}
Assume \cref{hyp:sequence}, with $\delta\le\delta_*$ and $0<r\le c_*\rho$, where $\delta_*$, $c_*$ and the flatness threshold $\eps_\Omega$ are the constants of \cref{m:fr:sec:coverbundle}. After a common
subsequence the charts converge, with limit transitions between them
(\cref{m:fr:lem:limittransitions}), and the compact principal
$O(n)$-bundles $Y_i=\Fr_\rho(X_i)$ have a compact metrizable limit $Y$ (\cref{m:fr:def:Y}), with a continuous action of $O(n)$
and $Y/O(n)=X$. More precisely:
\begin{enumerate}[label=\textup{(S\arabic*)},ref=\textup{(S\arabic*)},leftmargin=3.2em]
\item\label{fl:S:stage} \emph{The charts.}
Each $Z_i$ is simply connected and $\RCD(-C_0\rho^{-2},N+1)$.
Its deck group $\Gamma_i$ acts freely by measure-preserving isometries, and
$\Omega_i$ is $(\delta,R_1)$-Reifenberg and $\eps_\Omega$-flat at scale $r$.

\item\label{fl:S:limit} \emph{Limit charts.}
$(Z_i,z_i,\Gamma_i)\eGHto(Z,z,G)$, with Borel approximations $\af_i$
along which the measures converge. Points $q_i\in Z_i$ converge to $q\in Z$ if
$\af_i(q_i)\to q$; every point of $Z$ is such a limit, bounded sequences have convergent
subsequences, and distances converge. The limit has the same RCD bound,
$\Omega$ is $(3\delta,R_1)$-Reifenberg and $\eps_\Omega$-flat at scale $r$,
and $\varrho_i\to\varrho$ continuously. The projections converge to a map $\pi$ from $\Omega$ to $X$, in the sense
that $\pi_i(q_i)\to\pi(q)$ whenever points $q_i\in\Omega_i$ converge to $q\in\Omega$. It satisfies $d(\pi(q),\pi(q'))=\min_{g\in G}\hat d(q,gq')$, and it maps
$Z(R)$ onto $B_R(a)$ for $R\le R_0$. The group $G$ is a closed Lie group
of measure-preserving isometries acting properly, with open identity component.

\item\label{fl:S:groups} \emph{Group convergence.}
The convergence of elements is that of \cref{m:fr:sec:quotientenergy}:
bounded-displacement sequences are precompact, products and inverses respect
limits, and \mbox{$\gamma_iq_i\to gq$} when $\gamma_i\eGHto g$ and $q_i\to q$.
Every $g\in G$ is a limit, and these assertions hold along subsequences
(\cref{m:fr:lem:groupconv}).

\item\label{fl:S:frames} \emph{Frame actions.}
Over $\Omega_i$ and $\Omega$, the canonical frames form principal $O(n)$-bundles.
The actions $\gamma(q,e)=(\gamma q,e\circ\gamma^{-1})$ and
$g(q,e)=(gq,e\circ g^{-1})$ commute with the right $O(n)$-action.
They are free and proper, and the $\Gamma_i$-action is properly discontinuous.
Write $\Fr_i=\Fr(\Omega'_i)$ and $\Fr=\Fr(\Omega')$.

\item\label{fl:S:frameconv} \emph{Frame convergence.}
Convergence is defined by \eqref{m:fr:eq:frameconvdef}. Frames over convergent
base points in $\Omega$ are precompact; every frame over $q$ is a limit of
frames over any prescribed $q_i\to q$. Limits are unique, and
$p_iA_i\to pA$ and $\gamma_ip_i\to gp$ for convergent frames, matrices and
group elements. Frames converging over $\Omega^\flat$ eventually lie over
$\Omega_i^\flat$.

\item\label{fl:S:product} \emph{Trivializations.}
For a frame $p=(q,e)$ over $\Omega_i^\flat$ or $\Omega^\flat$, the frame map $u_e$ is a gauge map on $B_{3r}(q)$ and gives the trivialization
\[
\ft_p\colon B_{r/2}(q)\times O(n)\longrightarrow\pr^{-1}(B_{r/2}(q)),
\qquad \ft_p(x,A)=(x,e_{u_e}(x)\cdot A).
\]
It is a homeomorphism onto this subset of $\Fr_i$, respectively of $\Fr$, with $\ft_p(q,\one)=p$ and
$\ft_p(x,AA')=\ft_p(x,A)\cdot A'$. Let frames $p_i=(q_i,e_i)$ over $\Omega^\flat_i$ converge to a frame
$p=(q,e)$ over $\Omega^\flat$. If points $(x_i,A_i)$ of $B_{r/2}(q_i)\times O(n)$ converge to a point
$(x,A)$ of $B_{r/2}(q)\times O(n)$, then $\ft_{p_i}(x_i,A_i)\to\ft_p(x,A)$. If frames $p'_i$ over
$\Omega_i$ converge to a frame $p'$ with $\pr(p')\in B_{r/2}(q)$, then $\pr(p'_i)\in B_{r/2}(q_i)$
for large $i$, and $\ft_{p_i}^{-1}(p'_i)\to\ft_p^{-1}(p')$.

\item\label{fl:S:charts} \emph{Chart projections.}
The equivariant covering $\tau_i\colon\Fr_i\to\pr^{-1}(B_{R_0/2}(a_i))$
has fibres the $\Gamma_i$-orbits. The equivariant map
$\tau\colon\Fr\to\pr^{-1}(B_{R_0/2}(a))$ is continuous, open and onto, with fibres
the $G$-orbits. Both satisfy $\pr\circ\tau_i=\pi_i\circ\pr$ and
$\pr\circ\tau=\pi\circ\pr$, and
\[
d_i(\pr(\tau_i(p)),a_i)=\varrho_i(\pr(p)),\qquad
d(\pr(\tau(p)),a)=\varrho(\pr(p)).
\]
The isotropy group in $O(n)$ of $\tau(q,e)$ is $A_e(G_q)$, where $G_q$ is the isotropy group of
$q$ in $G$ and $A_e$ is the injective homomorphism with $e\circ g^{-1}=e\cdot A_e(g)$. If $\iota$
is a limit transition from $a$ to $b$ based at $q$ (\cref{m:fr:def:limittransition}), then
$\tau^a(q',e')=\tau^b(\iota(q'),e'\circ\iota^{-1})$ for the frames $(q',e')$ with
$q'\in B_{8r}(q)\cap\Omega'^a$ and $\iota(q')\in\Omega'^b$.

\item\label{fl:S:Yconv} \emph{Framed quotient convergence.}
The convergence of \cref{m:fr:def:Yconv} has unique limits,
every sequence has a convergent subsequence, and every point of $Y$ is a
limit. It respects $O(n)$ and the projections to $X_i$ and $X$. If $y^k\to y$ in $Y$ and, for every $k$,
points $y^k_i\in Y_i$ converge to $y^k$ as $i\to\infty$, then $y^k_{i_k}\to y$ for suitable indices
$i_1<i_2<\cdots$.
If $y_i\to y$ over $B_{R_0/2}(a)$, every lift $p\in\tau^{-1}(y)$ is the
limit of lifts $p_i\in\tau_i^{-1}(y_i)$.
Every $g\in C(Y)$ is the limit of uniformly bounded continuous functions
$g_i$ on $Y_i$, in the sense $g_i(y_i)\to g(y)$ whenever $y_i\to y$.
\end{enumerate}
\end{theorem}

Assertion \ref{fl:S:stage} follows
from \cref{m:fr:prop:completion,m:fr:prop:coverscharts}, and the limit transitions are those of
\cref{m:fr:lem:limittransitions}. Assertion \ref{fl:S:limit} is proved in
\cref{m:fr:sec:covercollapse}, except the Lie-group assertion, which is \cref{m:fr:prop:limitgroup};
\ref{fl:S:groups} is \cref{m:fr:lem:groupconv}. The frame assertions \ref{fl:S:frames} and
\ref{fl:S:frameconv} follow from \cref{m:fr:prop:coverbundle} and the local frame theory of
\cref{m:fr:prop:localization}. The freeness of $G$ in \ref{fl:S:frames}, the assertions on $Y$,
and \ref{fl:S:product}, \ref{fl:S:charts} and \ref{fl:S:Yconv} are proved in
\cref{m:fr:sec:covercollapse}.

We retain the subsequence of the theorem throughout. Since the balls $B_{R_0/8}(a)$ cover
$X$, every point of $Y$ has a lift over some $\Omega^\flat$, and convergent lifts
eventually lie over $\Omega^\flat_i$; so the trivializations of \ref{fl:S:product} are
available at every point of $Y$.

\subsection{From Local Covers to a Compatible Atlas}\label{m:book:ch:10}

We construct the objects of \cref{m:fr:thm:geomain} and prove most of its assertions. We first complete each cover without changing its interior geometry, then prove exact agreement of neighbouring covers. Their frames descend to a principal
$O(n)$-bundle, and a finite-atlas compactness argument gives the framed limit.

\subsubsection{Lifted Metrics and the Covering Scales}\label{m:fr:sec:localcovers}

We construct the lifted metric and measure used in \cref{m:fr:def:localcovers}
and fix the scales for the atlas. For the atlas construction, the parameters
$n,N,\kappa,\rho$ are those of \cref{hyp:sequence}: the spaces are $\RCD(\kappa,N)$ spaces with
$N\ge\max\{n,2\}$ and $|\kappa|\rho^2\le1$, and $n$ is the dimension of the Euclidean balls in
\cref{m:fr:def:localcovers}. Every statement below about a space with $(\delta,\rho)$-Reifenberg
local covers assumes these bounds. The construction of the lifted metric and measure,
\cref{m:fr:prop:completion}, and its regularity consequence \cref{hz:lem:scale} instead use
only their stated hypotheses; they require neither Reifenberg local covers nor
\cref{hyp:sequence}.

Let $(X,d,m)$ be an $\RCD(\kappa,N)$ space, $x\in X$ and $s>0$. The ball $B_s(x)$ is connected
and locally path-connected. It is semi-locally simply connected: by the main theorem of
\cite{Wan24b}, for every $y$ and every $R>0$ with $B_R(y)\subset B_s(x)$ there is $t<R$ such that
every loop in $B_t(y)$ is contractible in $B_R(y)$. So
it has a simply connected covering space $\pi\colon\widetilde B_s(x)\to B_s(x)$, on which
the group of deck transformations acts freely, and transitively on the fibres of $\pi$. We
fix a point $\tilde x$ over $x$. The \emph{length} of a curve in $\widetilde B_s(x)$ is
the length of its image in $X$, and $\tilde d(q,q')$ is the infimum of the lengths of
the curves from $q$ to $q'$. Let $y\in B_s(x)$, let $t>0$ be such that $B_{3t}(y)$ lies in $B_s(x)$ and is evenly
covered, and let $U$ be a sheet over $B_{3t}(y)$. A curve that starts at a point of $U$
over $B_t(y)$ and leaves $U$ has length at least $2t$, and the geodesics of $X$ between
points of $B_t(y)$ are shorter than $2t$ and lie in $B_{2t}(y)$, so they lift to $U$.
Hence $\pi$ maps the part of $U$ over $B_t(y)$ isometrically onto $B_t(y)$. It follows
that $\tilde d$ is finite, because a curve can be replaced by a chain of lifts of
geodesics, that it is a distance which induces the topology of $\widetilde B_s(x)$,
that the map $\pi$ is $1$-Lipschitz, and that curves in $B_s(x)$ lift to curves of the same
length. Let $\tilde m$ be the Borel measure on $\widetilde B_s(x)$ whose restriction to
every open set on which $\pi$ is injective is mapped by $\pi$ to the restriction of $m$. The
deck transformations preserve $\tilde d$ and $\tilde m$. We call
$(\widetilde B_s(x),\tilde d,\tilde m)$, with the point $\tilde x$, the \emph{local cover
of radius $s$ at $x$}. It is determined by $X$, $x$ and $s$ up to a unique isomorphism
that commutes with the projections and preserves the base points. If $q$ lies over
$y$ and $B_t(y)\subset B_s(x)$, then $\pi$ maps $B_t(q)$ onto $B_t(y)$, because geodesics of
$X$ that start at $y$ and are shorter than $t$ lift.

In the terminology of \cref{m:fr:def:reifenberg}, the lifted centre in
\cref{m:fr:def:localcovers} is $(\delta,\rho/3)$-Reifenberg. For non-collapsed
$\RCD(\kappa,n)$ spaces this is the local rewinding condition of \cite{HHWZ26},
up to the convention for GH distance.
The bound $|\kappa|\rho^2\le1$ gives the uniform curvature bound in the completion and the
loop contraction below.
Let $L_B=L_B(N)$ be the uniform contraction
constant of \cref{m:fr:thm:localBC} below. We use the radii
\begin{equation*}
R_0:=\frac{L_B\rho}{200},\qquad
R_1:=\frac{R_0}{32},\qquad R_B:=\frac{L_B\rho}{10}.
\end{equation*}
The radius $R_B$ is that of the balls on which the covers of nearby points agree
(\cref{m:fr:lem:agreement}), $R_0$ is the radius of the charts (\cref{m:fr:def:chart}), and $R_1$
that of their transition maps and Reifenberg balls (\cref{m:fr:prop:coverscharts}); the factors
leave room for the margins in these proofs.

\subsubsection{Completion of a Local Cover}\label{m:fr:sec:completion}
A local cover is in general not complete, and the results of the preceding sections are
formulated for $\RCD$ spaces. We make it complete by a time change in the sense of \cite{HS22}: the
distance is multiplied by a factor that tends to infinity at the boundary of the ball
and equals $1$ on the ball of half the radius, and the measure by the square of this
factor. Such a change needs a bound from above, and none from below, for the Laplacian
of the logarithm of the factor, which the Laplacian comparison theorem provides for
functions of the distance from the centre. The dimension bound increases; we let it
increase by one.

\begin{proposition}[Completion of a Local Cover]\label{m:fr:prop:completion}
For every $N>1$ there is $C_0=C_0(N)\ge1$ with the following property. Let $(X,d,m)$ be an
$\RCD(\kappa,N)$ space, let $x\in X$ and $s>0$, and let
$(\widetilde B_s(x),\tilde d,\tilde m)$ be the local cover of radius $s$ at $x$, with
projection $\pi$ and group $\Gamma$ of deck transformations. Then there are a distance
$\hat d$ and a Borel measure $\hat m$ on $\widetilde B_s(x)$ with the following
properties.
\begin{enumerate}[label=\textup{(\roman*)}]
\item Put $\kappa_s:=-(C_0/5)(s^{-2}+|\kappa|)$. Then
$(\widetilde B_s(x),\hat d,\hat m)$ is an $\RCD(\kappa_s,N+1)$ space with the
topology of $\widetilde B_s(x)$, and $\Gamma$ is a discrete group of isometries of it
that acts freely and preserves $\hat m$. In particular, it is
$\RCD(-C_0s^{-2},N+1)$ whenever $|\kappa|s^2\le4$.
\item $\hat d\ge\tilde d$ and $\hat m\ge\tilde m$. The two measures agree on
$\pi^{-1}(B_{s/2}(x))$, and $\hat d(q,q')=\tilde d(q,q')$ whenever $\pi(q)\in B_{s/8}(x)$ and
$\tilde d(q,q')\le\tfrac14s$.
\item For all $q,q'\in \pi^{-1}(B_{s/8}(x))$,
\begin{equation*}
d\bigl(\pi(q),\pi(q')\bigr)=\min_{\gamma\in\Gamma}\hat d(q,\gamma q') .
\end{equation*}
\end{enumerate}
\end{proposition}
\begin{proof}
Put $B:=B_s(x)$ and $d_x:=d(x,\cdot)$.

\emph{Step 1: the factors.} Let $h\colon[0,1)\to[0,\infty)$ be the function that vanishes
on $[0,\tfrac12]$ and equals $-\log(1-(2t-1)^4)$ on $[\tfrac12,1)$. It is of class $C^2$,
non-decreasing and convex, its integral of $e^h$ over $[0,1)$ is infinite, and
\begin{equation}\label{m:fr:eq:hbounds}
e^{-h}\,h'\le8,\qquad e^{-2h}\,h''\le64 .
\end{equation}
Fix a smooth non-increasing function $\chi\colon\R\to(0,1]$ with $\chi=1$ on
$(-\infty,0]$ and $\chi(v)=e^{-3v}$ for $v\ge1$, let $T$ be its primitive with $T(v)=v$
for $v\le0$, and let $T_\infty$ be the limit of $T$ at infinity. For $j\in\mathbb{N}$ put
$h_j(t):=j+T(h(t)-j)$ for $t<1$ and $h_j(t):=j+T_\infty$ for $t\ge1$. Near
$t=1$ we have $h_j(t)=j+T_\infty-\tfrac13e^{3j}(1-(2t-1)^4)^3$ for $t<1$, so
$h_j$ is a bounded non-decreasing function of class $C^2$ on $[0,\infty)$. It
agrees with $h$ where $h\le j$, and $0\le h_j\le h$ on $[0,1)$. With $v:=h-j$ we
have $h_j'=\chi(v)\,h'$ and $h_j''\le\chi(v)\,h''$ on $[0,1)$, and
$e^{-2h_j}\,\chi(v)\le e^3e^{-2h}$: for $v\le0$ this is clear, and for $v\ge0$ it
holds because $h_j\ge j$ and $\chi(v)\le e^{3-2v}$. With \eqref{m:fr:eq:hbounds} we get
\begin{equation}\label{m:fr:eq:thetabounds}
e^{-2h_j}\,(h_j')^2\le64\,e^3,\qquad e^{-2h_j}\,h_j'\le8\,e^3,\qquad
e^{-2h_j}\,h_j''\le64\,e^3 .
\end{equation}

\emph{Step 2: the time change.} Put $w_j:=h_j(d_x/s)$ on $X$. This function is
bounded and Lipschitz, it vanishes on $B_{s/2}(x)$, and it is constant outside $B$. By
the Laplacian comparison theorem \cite[Corollary~5.15]{Gig15}, which applies because
$N>1$, the function $d_x$ has a
measure-valued Laplacian on $X\setminus\{x\}$ whose singular part is non-positive and
whose absolutely continuous part is at most $C(N)(s^{-1}+\sqrt{|\kappa|})$ on
$B\setminus B_{s/2}(x)$. Indeed, for $\kappa\ge0$ the bound is $(N-1)/d_x$, and for
$\kappa<0$ it is
\[
\sqrt{(N-1)|\kappa|}\,\coth\!\left(d_x\sqrt{|\kappa|/(N-1)}\right)
\le\sqrt{(N-1)|\kappa|}+\frac{2(N-1)}s,
\]
using $\coth t\le1+1/t$. Since $|\nabla d_x|=1$ almost everywhere and $w_j$ vanishes near $x$,
the chain rule \cite[Proposition~4.11]{Gig15} shows that $w_j$ has the measure-valued Laplacian
\begin{equation*}
\Delta w_j=\frac1s\,h_j'\Bigl(\frac{d_x}{s}\Bigr)\,\Delta d_x
+\frac1{s^2}\,h_j''\Bigl(\frac{d_x}{s}\Bigr)\,m
\end{equation*}
on $X$. Its singular part is non-positive, because $h_j'\ge0$, and its absolutely
continuous part is at most
$s^{-2}\bigl(C(N)(1+s\sqrt{|\kappa|})\,h_j'+h_j''\bigr)(d_x/s)$. Moreover
$|\nabla w_j|^2=s^{-2}\,h_j'(d_x/s)^2$. This measure vanishes outside the closure of $B$, which is compact. Let $\mu^+$ and $\mu^-$
be its positive and its negative part; $\mu^+$ has a bounded density $f^+$ with compact support.

The time-change theorem of \cite{HS22} used below applies to locally bounded functions that are
continuous almost everywhere, whose measure-valued Laplacian has a non-positive singular part,
and that lie in the domain of the measure-valued Laplacian in the sense of
\cite[Section~2.2]{HS22}: the Laplacian is a signed Borel measure $\mu$ that charges no set of
capacity zero, and $\int\tilde v\,d\mu=-\int\langle\nabla w_j,\nabla v\rangle\,dm$ for every
$v\in\Hone(X)$ whose quasi-continuous representative $\tilde v$ is $|\mu|$-integrable. The
function $w_j$ is bounded and Lipschitz, it differs from a constant by a function with compact
support, and the singular part of $\Delta w_j$ is non-positive. The identity holds for the
Lipschitz functions $v$ with compact support, by the definition of the Laplacian in
\cite{Gig15}, and it remains to extend it to all of $\Hone(X)$. The Cheeger energy of $X$ is a
strongly local Dirichlet form, and it is regular: the Lipschitz functions with compact support
are uniformly dense in the continuous functions with compact support, and they are dense in
$\Hone(X)$, by the density in energy \cite[Theorem~7.3 and Section~8.2]{AGS13} and because the
Cheeger energy is quadratic. So the potential theory of \cite[Chapter~2]{FOT11} applies, and it
suffices to show that $\mu^+$ and $\mu^-$ are measures of finite energy integral, that is,
\begin{equation*}
\int|v|\,d\mu^\pm\le C\,\|v\|_{\Hone(X)}\qquad\text{for all }v\in\Hone(X)\cap C_c(X)
\end{equation*}
with a constant $C$. For $\mu^+$ this holds with $C=\|f^+\|_{L^2}$. For $\mu^-$, let first
$u$ be a non-negative Lipschitz function with compact support. The definition of the
Laplacian gives
\begin{equation}\label{m:fr:eq:finiteenergy}
\int u\,d\mu^-=\int u\,f^+\,dm+\int\langle\nabla w_j,\nabla u\rangle\,dm
\le\|f^+\|_{L^2}\,\|u\|_{L^2}+\|\nabla w_j\|_{L^2}\,\|\nabla u\|_{L^2} ,
\end{equation}
and the norms of $f^+$ and of $|\nabla w_j|$ are finite, because these functions are bounded
and vanish outside the closure of $B$; in particular $\mu^-$ is a finite measure. For
$v\in\Hone(X)\cap C_c(X)$ we approximate $|v|$ by discrete convolutions \cite[Section~3]{KL02}:
for a maximal $\eps$-separated set $\{p_k\}$ and a partition of unity $\{\psi_k\}$ with
$\supp\psi_k\subset B_{2\eps}(p_k)$ and $\Lip\psi_k\le C/\eps$, where $C$ may depend on
the fixed local geometry, the functions
$v_\eps:=\sum_k\psi_k\fint_{B_{2\eps}(p_k)}|v|\,dm$ are non-negative and Lipschitz, their
supports lie in a fixed compact set for $\eps\le1$, they converge to $|v|$ uniformly, because
$|v|$ is uniformly continuous, and for small $\eps$ their norms in $\Hone(X)$ are at most a
constant times that of $|v|$, hence of $v$, by the doubling property and the Poincar\'e
inequality of \ref{m:fr:R:geometry}. Since $\mu^-$ is finite, \eqref{m:fr:eq:finiteenergy} for
$u=v_\eps$ gives the bound in the limit $\eps\to0$.

Measures of finite energy integral charge no set of capacity zero, and for them every
$v\in\Hone(X)$ has a quasi-continuous representative $\tilde v$ with
$\int|\tilde v|\,d\mu^\pm\le C'\,\|v\|_{\Hone(X)}$ \cite[Sections~2.1 and~2.2]{FOT11}. Both
sides of the identity $\int\tilde v\,d(\Delta w_j)=-\int\langle\nabla w_j,\nabla v\rangle\,dm$
are therefore continuous on $\Hone(X)$, and the identity extends from the dense subspace of the
Lipschitz functions with compact support to all $v\in\Hone(X)$. So $w_j$ belongs to the class of
functions that is admitted in \cite{HS22}. Since $e^{-2w_j}\kappa\ge-|\kappa|$ and
$s\sqrt{|\kappa|}\le1+|\kappa|s^2$, \eqref{m:fr:eq:thetabounds} allows us to choose
$C_0=C_0(N)$ large enough that
\begin{equation*}
e^{-2w_j}\Bigl(\kappa-(N-2)(N-1)\,|\nabla w_j|^2-\Delta_{\mathrm{ac}}w_j\Bigr)\ge\kappa_s
\end{equation*}
almost everywhere, where $\Delta_{\mathrm{ac}}w_j$ is the density of the absolutely
continuous part. By \cite[Theorem~2]{HS22}, applied with the dimension $N+1$,
the space $(X,d_j,m_j)$ is an $\RCD(\kappa_s,N+1)$ space, where $m_j:=e^{2w_j}m$ and $d_j(y,y')$ is
the infimum of the integrals $\int e^{w_j(c(t))}\,|\dot c(t)|\,dt$ over the absolutely
continuous curves $c$ from $y$ to $y'$ in $X$; this description of the distance is
\cite[Lemma~6]{HS22}.

\emph{Step 3: the limit.} On $B$ put $w:=h(d_x/s)$ and $m_w:=e^{2w}m$, and let
$d_w(y,y')$ be the infimum of the integrals $\int e^{w(c(t))}\,|\dot c(t)|\,dt$ over the
absolutely continuous curves from $y$ to $y'$ in $B$. Since $d_x$ is $1$-Lipschitz, a
curve in $X$ that starts at $x$ and reaches a point at distance $ts$ from $x$ has
$d_j$-length at least $s\int_0^te^{h_j}$. Let $L>0$, let $t_L\in(0,1)$ be the number
with $s\int_0^{t_L}e^h=3L$, and let $j\ge h(t_L)$, so that $h_j=h$ on $[0,t_L]$. Then
every point $y$ with $d_j(x,y)<3L$ lies in the closed ball $\bar B$ of radius $t_Ls$
about $x$, which is contained in $B$ and on which $w_j=w$. If $d_j(x,y)\le L$ and
$d_j(x,y')\le L$, then a curve from $y$ to $y'$ whose $d_j$-length is less than
$d_j(y,y')+L$ stays within $d_j$-distance less than $3L$ of $x$, hence in $\bar B$, and
its length for $w$ is its $d_j$-length. So $d_w(y,y')\le d_j(y,y')$, and the reverse
inequality holds on $B$ because $w_j\le w$ there. Hence the closed balls of radius $L$
about $x$ for $d_j$ and for $d_w$ coincide, and on them $d_j=d_w$ and $m_j=m_w$. These
balls are compact, because $(X,d_j)$ is proper, as $d\le d_j\le e^{j+T_\infty}d$. So $(B,d_w)$ is proper and complete, it
has the topology of $B$, and $(B,d_w,m_w,x)$ is the pointed measured Gromov--Hausdorff
limit of the spaces $(X,d_j,m_j,x)$. By the stability of the $\RCD$ condition
\cite{GMS15,EKS15} it is an $\RCD(\kappa_s,N+1)$ space.

\emph{Step 4: the cover.} By \cite[Theorem~1.1]{MW19} the universal cover of
$(B,d_w,m_w)$, with the lifted distance $\hat d$ and the lifted measure $\hat m$, is an
$\RCD(\kappa_s,N+1)$ space. As a covering space of $B$ it is $\widetilde B_s(x)$, and
the deck transformations preserve $\hat d$ and $\hat m$. They form a discrete group of
isometries, because the orbits are discrete and the action is free. The function $w$ is
continuous, so the length of a curve in $B$ for $d_w$ is its integral of $e^{w}$. Hence
$\hat d(q,q')$ is the infimum of the integrals of $e^{w\circ \pi}$ along the curves from
$q$ to $q'$, with respect to the length of curves defined in \cref{m:fr:sec:localcovers}, and
$\hat m=e^{2w\circ \pi}\,\tilde m$. Since $w\ge0$, and $w=0$ on $B_{s/2}(x)$, we get the
first two assertions of (ii). Let $\pi(q)\in B_{s/8}(x)$ and $\tilde d(q,q')\le\tfrac14s$. A
curve from $q$ to $q'$ of length less than $\tfrac38s$ is mapped by $\pi$ into
$B_{s/2}(x)$, so its integral of $e^{w\circ \pi}$ is its length. Hence
$\hat d(q,q')\le\tilde d(q,q')$, which proves (ii).

For (iii) let $q,q'\in \pi^{-1}(B_{s/8}(x))$. A geodesic of $X$ from $\pi(q)$ to $\pi(q')$ is
shorter than $\tfrac14s$ and lies in $B_{s/4}(x)$. Its lift that starts at $q$ ends at
$\gamma q'$ for some $\gamma\in\Gamma$, so $\tilde d(q,\gamma q')\le d(\pi(q),\pi(q'))$, and
by (ii) the same holds for $\hat d$. Conversely
$\hat d(q,\gamma q')\ge\tilde d(q,\gamma q')\ge d(\pi(q),\pi(q'))$ for every $\gamma$, because
$\pi$ is $1$-Lipschitz.
\end{proof}

We call $(\widetilde B_s(x),\hat d,\hat m)$ the \emph{completed local cover} of radius $s$
at $x$. By (ii), for a point $q$ over $B_{s/8}(x)$ the balls $B_t(q)$ with
$t\le\tfrac14s$ are the same for $\tilde d$ and for $\hat d$, and the two distances agree
on pairs of points of such a ball when $t\le\tfrac1{20}s$, because then the ball lies
over $B_{s/8+s/20}(x)$, where the argument of the proof of (ii) applies with curves that
are shorter than $\tfrac15s$.

The completion also makes local bounded covering geometry
(\cref{hz:def:lbcg}) accessible to the regularity theory of \cite{HZ26}. Recall the
metric definition \cite[Definition~2.24]{HZ26}:
\[
\operatorname{Cut}(q)=\{z:d(q,z)+d(z,w)>d(q,w)\text{ for every }w\ne z\}.
\]
Set $\operatorname{Injrad}(q):=\inf_{z\in\operatorname{Cut}(q)}d(q,z)$, with
$\inf\varnothing=+\infty$. These formulas apply also to the incomplete lifted
length metric. In particular, $z\notin\operatorname{Cut}(q)$ means that some $w\ne z$
satisfies $d(q,z)+d(z,w)=d(q,w)$.

We use the covering regularity scale $\mathfrak s^\alpha_\eta(y)$ of
\cite[Definition~5.1]{HZ26} with $k=0$: it is the supremum of the radii $r$ for which a
local diffeomorphism $B_{2r}(\hat y)\to B_{2r}(y)$, equipped with the lifted metric and
density, admits weighted-harmonic coordinates $u$ on $B_r(\hat y)$ with
$\|g^{ab}-\delta_{ab}\|_{C^0}+r^\alpha[g^{ab}]_{C^{0,\alpha}}<\eta$ for every $a,b$,
where $g^{ab}=\langle\nabla u_a,\nabla u_b\rangle$.

\begin{corollary}[Regularity under Local Bounded Covering Geometry]\label{hz:lem:scale}
Let $(X,d,m)$ be an $\RCD(\kappa,N)$ space with $N>1$, of essential dimension $n$, with
$(\rho,\iota_0)$-local bounded covering geometry.
\begin{enumerate}[label=\textup{(\roman*)}]
\item For every $p\in[1,\infty)$ and $\alpha\in(0,1)$, $X$ is a
$W^{1,p}_{\mathrm{loc}}\cap C^{0,\alpha}_{\mathrm{loc}}$-Riemannian manifold of dimension
$n$ in harmonic charts, and $m=e^{-f}\mathscr H^n$ with
$f\in W^{1,p}_{\mathrm{loc}}\cap C^{0,\alpha}_{\mathrm{loc}}$.
\item For every $\eta>0$ and $\alpha\in(0,1)$ there is $c_4>0$, depending only on
$N$, $|\kappa|\rho^2$, $\iota_0/\rho$, $\eta$ and $\alpha$, such that
$\mathfrak s^\alpha_\eta(y)\ge c_4\rho$ for every $y\in X$.
\end{enumerate}
\end{corollary}
\begin{proof}
Fix $y\in X$ and complete its radius-$\rho$ local cover by
\cref{m:fr:prop:completion}. Write $\pi$ for the projection and put
$U:=\pi^{-1}(B_{\rho/16}(y))$ and $\iota_1:=\min\{\iota_0,\rho/32\}$.

\emph{Injectivity radius.} Every $q\in U$ is a deck translate of a point of
$B_{\rho/16}(\tilde y)$, so $\operatorname{Injrad}_{\tilde d}(q)\ge\iota_0$.
Suppose $0<\hat d(q,z)<\iota_1$. By \cref{m:fr:prop:completion}(ii),
$\ell:=\hat d(q,z)=\tilde d(q,z)$, and there is $w\ne z$ with
$\tilde d(q,w)=\ell+\tilde d(z,w)$. Choose
$0<a<\min\{\tilde d(z,w),\rho/32\}$ and let
$S:=\{u:\tilde d(z,u)=a\}$. Since $z$ lies over $B_{3\rho/32}(y)$, the closed
$\tilde d$-ball of radius $a$ about $z$ is also its closed $\hat d$-ball and is compact.
Every curve from $z$ to $w$ crosses $S$, so the length-metric property and the triangle
inequality give
\[
\tilde d(z,w)=a+\min_{u\in S}\tilde d(u,w).
\]
For a minimizer $u$, the triangle inequalities imply $\tilde d(q,u)=\ell+a$.
The three distances between $q,z,u$ are less than $\rho/16$, and $q,z$ lie over
$B_{\rho/8}(y)$; hence \cref{m:fr:prop:completion}(ii) gives the same distances for
$\hat d$. Thus $z\notin\operatorname{Cut}_{\hat d}(q)$, proving
$\operatorname{Injrad}_{\hat d}(U)\ge\iota_1$.

\emph{Regularity and descent.} By \cref{m:fr:sec:localcovers} and
\cref{m:fr:prop:completion}(ii), $\pi$ on $U$ is a local metric-measure isometry.
The completed cover has essential dimension $n$, since on $U$ its tangent cones are
those of $X$. Apply \cite[Theorem~4.4 and Corollary~4.7]{HZ26} to this open subset of
the complete $\RCD(\kappa_\rho,N+1)$ space. The resulting harmonic charts, metric and
density descend by these local isometries. Their transition maps lift locally to
transition maps on the cover. This proves (i).

\emph{Uniform scale.} The ball $B_{\rho/32}(\tilde y;\hat d)$ lies in $U$. Dividing
distances by $\rho/32$ gives an
$\RCD(-C_0(1+|\kappa|\rho^2)/5120,N+1)$ space whose unit ball has injectivity radius
at least $\min\{32\iota_0/\rho,1\}$. Choose $p>(N+1)/(1-\alpha)$. The quantitative
part of \cite[Theorem~4.4]{HZ26} supplies harmonic charts with arbitrarily small
scale-normalized $L^p$ Hessian bounds and Gram coefficients close to the identity,
at a radius depending only on the listed parameters and the desired accuracy.
The gradient product rule bounds $|\nabla g^{ab}|$ by the Hessians of these
coordinates. Morrey's inequality, using the normalized measure and the local
Bishop--Gromov and Poincar\'e bounds, then gives an arbitrarily small
scale-invariant $C^{0,\alpha}$ bound for $g^{ab}-\delta_{ab}$ on a fixed smaller
ball. Choose its radius $r$ with $2r\le\rho/32$. The projection
$B_{2r}(\tilde y;\hat d)\to B_{2r}(y)$ is onto by
\cref{m:fr:sec:localcovers} and \cref{m:fr:prop:completion}(ii), and is a local
diffeomorphism for the structures in (i). It pulls back the metric and density to
those on the cover. These are precisely the data defining
$\mathfrak s^\alpha_\eta(y)$, and rescaling back proves (ii).
\end{proof}

We now return to the standing Reifenberg local-cover hypothesis and the atlas scales
fixed in \cref{m:fr:sec:localcovers}.

\subsubsection{Agreement of Nearby Covers and Flatness}\label{m:fr:sec:agreement}
The local covers of different balls are incomplete spaces whose domains depend on the
centre, and their outer balls need not agree. We prove that the balls of a fixed smaller
radius do agree: first, short loops of a cover are contracted within a ball of a fixed
larger radius of the same cover; then the covers are lifted to each other in both
directions. The contraction uses the curvature-dimension bound and the closeness of the
balls about the centres to Euclidean balls. The resulting isometries are determined by
their value at one point. They make the balls about all points near the centre of a cover
close to Euclidean balls, which therefore have parametrizations by Euclidean balls
(\cref{m:fr:cor:localUIC}), and they are the transition maps of the charts
(\cref{m:fr:prop:coverscharts}) and join the frames of the covers
(\cref{m:fr:sec:coverbundle}).

\begin{definition}[Loop Contraction]\label{m:fr:def:localBC}
Write $\widetilde B_\rho(x)$ for the local cover of radius $\rho$ at $x$ and $\pi_x$ for its projection. The covers have
\emph{loop contraction} with constants $0<L_B<b_B<1/8$ if the inclusion
\begin{equation}\label{m:fr:eq:localBC}
\pi_1(B_{L_B\rho}(\tilde x;\widetilde B_\rho(x)),\tilde x)
\longrightarrow\pi_1(B_{b_B\rho}(\tilde x;\widetilde B_\rho(x)),\tilde x)
\end{equation}
is zero for every $x$. The balls carry the lifted distances.
\end{definition}

Under the standing hypothesis the covers have this property, with constants that depend only
on $N$. The proof completes the cover, passes to the local cover of a smaller ball inside it,
and excludes deck transformations of small displacement there.

\begin{theorem}[Loop Contraction]\label{m:fr:thm:localBC}
There are $\delta_B=\delta_B(N)>0$ and numbers $L_B,b_B$ depending only on $N$,
with \mbox{$0<L_B<b_B<1/8$}, such that
every $\RCD(\kappa,N)$ space with $(\delta,\rho)$-Reifenberg local covers,
where $\delta\le\delta_B$ and $|\kappa|\rho^2\le1$, has loop contraction.
\end{theorem}
\begin{proof}
Rescale to $\rho=1$; the hypothesis $|\kappa|\rho^2\le1$ is not affected by this, so that $|\kappa|\le1$. Complete $\widetilde B_1(x)$ by \cref{m:fr:prop:completion}, obtaining an
$\RCD(-C,N+1)$ space $(Z,z,m_Z)$ with $C=C_0(N)$. The balls of radius at most $1/8$
about the centre carry the original distance and measure. Put
\[
 b_B:=\min\{1/128,1/(2\sqrt C)\}.
\]
Apply \cref{m:fr:prop:completion} to the local cover of $B_{b_B}(z;Z)$ and denote the
completed local cover by $(W,\omega)$, with projection $\pi_W\colon W\to B_{b_B}(z;Z)$.
Its topology is that of the simply connected local cover, which exists by the main theorem
of \cite{Wan24b}, as in \cref{m:fr:sec:localcovers}, and it is $\RCD(-C_0(N+1)\,b_B^{-2},N+2)$. Both completions fall under
the bounded-scale conclusion, since $Cb_B^2\le1/4$. The second conformal factor is at least
one and equals one over $B_{b_B/2}(z)$. Thus $\pi_W$ is $1$-Lipschitz for the distance
of $Z$, and curves over this smaller ball have unchanged lengths.

If there is a nonidentity deck transformation of $\pi_W$, its positive based
displacement has an attained minimum $d$ among all such transformations: the action
is free and properly discontinuous, and $W$ is proper. Suppose $d<h:=b_B/100$, and
join $\omega$ to $g\omega$ by a minimizing geodesic attaining this minimum.
Its projected loop $\gamma$ lies in $\bar B_{d/2}(z)$, where the conformal factor is one.
Its midpoint $v$ has $d_Z(z,v)=d/2$. Indeed, a shorter path to $v$ would split
$\gamma$ into two based loops of length less than $d$. They cannot both lift closed,
since their product represents $g\ne1$; one contradicts the definition of $d$. This is the
argument of Sormani's halfway lemma \cite{Sor00}.
All the paths used here stay in the unchanged core.

Apply \cite[Theorem~2.15]{MW19} to $(W,b_B^{-1}d_W)$. There are
$A=A(C_0(N+1),N+2)>0$ and $\alpha=\alpha(N+2)\in(0,1)$ such that
\[
 e_{\omega,g\omega}(w)\le A t^{1+\alpha}d
 \quad\text{if }d_W(w,\widetilde\gamma(1/2))<td,\quad 0<t<1/4.
\]
The endpoint distance after rescaling is $d/b_B<1$, as required by that theorem.
Choose $S\in(0,1/32)$ with $A(3S)^{1+\alpha}<S$.
At the scale $10d$ the distance of $Z$ is that of the local cover, since
$10d<b_B/10<1/8$, so by \cref{m:fr:def:localcovers} there is a $10\delta d$-GH
approximation $\Phi\colon B_{10d}(z)\to B_{10d}(0^n)$. Let $v_+\in B_{10d}(z)$ be a point
whose image lies within $10\delta d$ of the point of the ray from $\Phi(z)$ through
$\Phi(v)$ at distance $(\tfrac12+2S)d$ from $\Phi(z)$; it exists because the image of
$\Phi$ is $10\delta d$-dense. Then
\[
 d_Z(z,v_+)\ge(1/2+2S-c\delta)d,\qquad
 d_Z(v,v_+)\le(2S+c\delta)d ,
\]
with a universal constant $c$. Decrease $\delta_B$ so that $c\delta_B<S/2$.
A short minimizing curve from $v$ to $v_+$ stays in $B_{b_B/2}(z)$ and lifts from
$\widetilde\gamma(1/2)$ to a point $w_+$ with
$d_W(w_+,\widetilde\gamma(1/2))<3Sd$. Since both endpoints of the deck geodesic
project to $z$ and $\pi_W$ is $1$-Lipschitz,
\[
 e_{\omega,g\omega}(w_+)\ge2d_Z(z,v_+)-d
 \ge(4S-2c\delta)d>3Sd,
\]
contradicting the upper bound $A(3S)^{1+\alpha}d<Sd$.
Hence every nonidentity based deck displacement is at least $h$.

Put $L_B:=h/10=b_B/1000$. Short radial curves lift, so $\pi_W$ maps
$B_{L_B}(\omega)$ onto $B_{L_B}(z)$.
Two lifts in this ball with the same image differ by a deck transformation with
based displacement less than $2L_B<h$, so the restriction is injective.
It preserves pairwise distances as well. A minimizing curve between two points of
$B_{L_B}(z)$ stays in $B_{3L_B}(z)$ and lifts with unchanged length.
A different lift endpoint would yield a deck transformation with based displacement
less than $4L_B<h$, which is impossible. The reverse distance inequality follows
from the $1$-Lipschitz projection.
The inverse isometry is a section over $B_{L_B}(z)$. Every continuous loop there
lifts closed to the simply connected $W$ and contracts in $W$; projecting the
contraction makes it null-homotopic in $B_{b_B}(z)$. These balls are the corresponding
balls of $\widetilde B_1(x)$, with the original distance. Rescaling proves \eqref{m:fr:eq:localBC}.
\end{proof}

With loop contraction in hand, the covers of nearby centres can be lifted to each other.

\begin{lemma}[Agreement of Nearby Covers]\label{m:fr:lem:agreement}
Suppose the local covers have loop contraction with $L_B,b_B$. Let
$d(x,y)<L_B\rho/20$ and let $q\in \widetilde B_\rho(x)$ project to $y$.
There is exactly one continuous map
\[
 c_q\colon B_{R_B}(\tilde y;\widetilde B_\rho(y))\to \widetilde B_\rho(x),
 \qquad \pi_x c_q=\pi_y,\qquad c_q(\tilde y)=q .
\]
It is a measure-preserving isometry onto $B_{R_B}(q;\widetilde B_\rho(x))$.
If the hypothesis of \cref{m:fr:def:localcovers} holds, $B_t(q)$ is $\delta t$-close to
$B_t(0^n)$ for $0<t\le R_B$. Two compositions of such maps that lift the same
projection and have the same value at one point agree on every connected set on which
both are defined.
\end{lemma}
\begin{proof}
First take $q$ to be the endpoint of the lift from $\tilde x$ of a short minimizing
geodesic from $x$ to $y$, so $d(\tilde x,q)<L_B\rho/20$.
Since $b_B+L_B/20<1$, the projection of $B_{b_B\rho}(\tilde y)$ lies in $B_\rho(x)$,
and the projection of $B_{b_B\rho}(\tilde x)$ lies in $B_\rho(y)$.
By loop contraction, every loop in either ball of radius $L_B\rho$ is null-homotopic in
the ball of radius $b_B\rho$ of its cover, which projects into the ball of radius $\rho$
about the other centre.
The covering lifting criterion, with change of base point to $q$ in the second
ball, gives
\[
 \lambda\colon B_{L_B\rho}(\tilde y)\to \widetilde B_\rho(x),\quad \lambda(\tilde y)=q,\qquad
 \lambda'\colon B_{L_B\rho}(\tilde x)\to \widetilde B_\rho(y),\quad \lambda'(q)=\tilde y,
\]
commuting with the projections to $X$. They are local measure-preserving
isometries and preserve lengths of curves in their domains.

Put $a=L_B/10$, so $R_B=a\rho$.
For $0<s\le3a\rho$, a curve of length less than $s$ that starts at the centre stays in the
domain of $\lambda$, respectively of $\lambda'$. Consequently
$\lambda(B_s(\tilde y))\subset B_s(q)$ and $\lambda'(B_s(q))\subset B_s(\tilde y)$.
The target-side balls lie in the domain of $\lambda'$, because the sum $d(\tilde x,q)+3a\rho$ is less than \mbox{$7L_B\rho/20<L_B\rho$}.
The compositions $\lambda'\lambda$ and $\lambda\lambda'$ are lifts of the projections to $X$ that agree with
the identity at one point. Uniqueness of lifts on these connected balls
makes the compositions the identity.
Thus $\lambda$ and $\lambda'$ are inverse bijections on the radius-$a\rho$ balls.
Almost minimizing curves between two points of either ball stay in the corresponding
radius-$(3a\rho+\eta)$ ball for small $\eta>0$. Applying $\lambda$ and $\lambda'$ to these curves
and letting $\eta\downarrow0$ proves equality of pairwise distances.
The local measure identities give measure preservation.
For any other lift $q$, postcompose with a deck transformation sending the short-path
lift to $q$. The uniqueness of lifts with a prescribed value shows that the result does not depend
on this choice.
Finally, two compositions of such maps that lift the same projection and agree at one
point agree on their connected common domain, by the same uniqueness.
\end{proof}

We collect what we have obtained in the object with which the rest of the paper works.

\begin{definition}[The Chart at a Point]\label{m:fr:def:chart}
For $x\in X$ let $(Z_x,\hat d,\hat m)$ be the completed
local cover of radius $\rho$ at $x$ of \cref{m:fr:prop:completion}, with projection
$\pi_x\colon Z_x\to B_\rho(x)$, base point $z_x$ over $x$ and group $\Gamma_x$ of deck
transformations, and put
\begin{equation*}
\Omega_x:=\pi_x^{-1}\bigl(B_{R_0}(x)\bigr),\qquad\Omega'_x:=\pi_x^{-1}\bigl(B_{R_0/2}(x)\bigr).
\end{equation*}
The \emph{chart} of $X$ at $x$ is $(Z_x,z_x,\Gamma_x)$, with $\pi_x$ and these two sets.
\end{definition}

The larger set is the region on which the chart is flat, and the smaller one is
the part that we use to describe $X$; the margin between them contains the balls of
radius $R_1$ about the points of $\Omega'_x$.

\begin{proposition}[Charts]\label{m:fr:prop:coverscharts}
For every $\eps>0$ there are numbers ${\delta_C=\delta_C(n,N,\eps)>0}$ and
${c_C=c_C(n,N,\eps)\in(0,\tfrac1{3200}]}$ with the following property. Let $(X,d,m)$ be an
$\RCD(\kappa,N)$ space with $|\kappa|\rho^2\le1$ that has $(\delta,\rho)$-Reifenberg local
covers, with $\delta\le\delta_C$.
\begin{enumerate}[label=\textup{(\roman*)}]
\item \emph{Quotients.} For every $x\in X$, the space $Z_x$ is an
$\RCD(-C_0\rho^{-2},N+1)$ space, and $\Gamma_x$ is a discrete group of isometries of
$Z_x$ that acts freely and preserves $\hat m$. The set $\Omega_x$ consists of the points
$q$ with $\min_\gamma\hat d(\gamma q,z_x)<R_0$, we have
$d(\pi_x(q),\pi_x(q'))=\min_\gamma\hat d(q,\gamma q')$ for $q,q'\in\Omega_x$, with the minima
over $\Gamma_x$, and $\pi_x$ maps the restriction of $\hat m$ to an open subset of
$\Omega_x$ on which it is injective to the restriction of $m$ to its image.
\item \emph{Flatness.} For every $x\in X$, the region $\Omega_x$ of $Z_x$ is
$(\delta,R_1)$-Reifenberg, and it is an $\eps$-flat region of $Z_x$ at
every scale $r\le c_C\rho$.
\item \emph{Transitions.} Let $x,x'\in X$, $q\in\Omega_x$ and $q'\in\Omega_{x'}$ with
$\pi_x(q)=\pi_{x'}(q')$. Then there is exactly one continuous map
$\iota\colon B_{R_1}(q)\to Z_{x'}$ with $\pi_{x'}\circ\iota=\pi_x$ and $\iota(q)=q'$. It
maps $B_{R_1}(q)$ onto $B_{R_1}(q')$, preserves the distances $\hat d$ between all
pairs of points, and maps $\hat m$ to $\hat m$. If $x=x'$, it is the restriction of an
element of $\Gamma_x$.
\end{enumerate}
\end{proposition}
\begin{proof}
Let
\[
\delta_C:=\min\{\delta_B(N),\delta_U(n),
                   \tfrac14\delta_H(n,N+1,\eps)\}.
\]
The factor $1/4$ leaves room for the limit balls in the proof of \ref{fl:S:limit}. With
$\delta_U$ in the minimum, \cref{m:fr:cor:localUIC} applies to the balls of radius at most
$R_1$ in the regions $\Omega_x$, by~(ii).
Part~(i) is \cref{m:fr:prop:completion} with \mbox{$s=\rho$}, since
$R_0<\rho/8$ and $|\kappa|\rho^2\le1$.
For $q\in\Omega_x$ put $y=\pi_x(q)$. Then \mbox{$d(x,y)<R_0<L_B\rho/20$}.
By \cref{m:fr:thm:localBC,m:fr:lem:agreement}, the ball of radius $R_B$ about $q$ in the
local cover is isometric to the ball of the same radius about the centre of $\widetilde B_\rho(y)$,
even when $q$ is far from $z_x$. Hence it is $\delta t$-close to $B_t(0^n)$ at every
smaller radius $t$.
The completion preserves its distances and measures, because
$R_0+R_B<\rho/8$ and $2R_B<\rho/4$.
This proves the first assertion of (ii), in fact for $t\le R_B$.

Put $s_C:=\min\{R_1,\rho\,C_0^{-1/2}\}$, let
$c_H=c_H(n,N+1,\eps)$, and choose
\[
c_C:=\min\{R_0/(4096\rho),\,c_Hs_C/(64\rho)\}.
\]
If $r\le c_C\rho$, apply \cref{m:fr:lem:flatness} with $s=r/c_H$ in $Z_x$.
Then $s\le s_C/64$, the required ball lies in the compared region, and
$C_0\rho^{-2}s^2\le1/4096$. The lemma gives $\eps$-flatness at scale $r$.
Notice also that $128r\le R_1$.

For (iii) put $y=\pi_x(q)=\pi_{x'}(q')$. Both centres $x$ and $x'$ have distance less than
$R_0$ from $y$, so \cref{m:fr:lem:agreement} gives the two isometries of the ball of
radius $R_B$ of $\widetilde B_\rho(y)$ onto the corresponding balls about $q$ and $q'$.
Their composition $c_{q'}c_q^{-1}$, restricted to $B_{R_1}(q)$, is the required
$\iota$. It is unique because a lift with a prescribed value at one point is unique on
a connected set. The inverse and the compositions of three such maps agree on their
connected common domains for the same reason. If $x=x'$, the deck transformation sending $q$ to $q'$
is such a lift.
\end{proof}
\subsubsection{The Frame Bundle of the Local Covers}\label{m:fr:sec:coverbundle}
The frames of the local covers at the points of $X$ form a bundle over $X$ itself. Over
a ball of radius $\tfrac12R_0$ it is the quotient of the frame bundle of a chart by the
group of the chart. We first record the invariance under isometries of balls.

Let $\eps_\Omega$ denote the minimum of the constants
$\eps_F(n)$ of \cref{m:fr:prop:localization}, $\eps_P(n)$ of
\cref{m:fr:prop:freeframes}, and the gauge thresholds
$\eps_G(n)$ and $\eps_M(n,\theta_G(n))$ of
\cref{m:fr:thm:gauges,m:fr:lem:framegauge}; the gauge thresholds are included so that the
frame maps on the limits of the charts are gauge maps, as the proof of \cref{ls:thm:limits}
requires. Put
\begin{equation*}
\delta_*:=\min\bigl\{\delta_C(n,N,\eps_\Omega),\,\tfrac13\delta_P(n)\bigr\},\qquad
c_*:=c_C(n,N,\eps_\Omega),
\end{equation*}
with the constants of \cref{m:fr:prop:coverscharts,m:fr:prop:freeframes}. The spaces have
$(\delta,\rho)$-Reifenberg local covers with
$\delta\le\delta_*$, and the scale $r$ of the frames satisfies $0<r\le c_*\rho$; thus $128r\le R_1$. By \cref{m:fr:prop:coverscharts}, the regions $\Omega_x$ of the charts are
then $\eps_\Omega$-flat at scale $r$, so \cref{m:fr:prop:localization} applies to them, because $\eps_\Omega\le\eps_F$.

\begin{lemma}[Isometries of Balls]\label{m:fr:lem:localnatural}
Let $(Z,d,m)$ and $(Z',d',m')$ be $\RCD$ spaces of finite dimension, let
$\Omega\subset Z$ and $\Omega'\subset Z'$ be open sets that are $\eps_\Omega$-flat at scale $r$,
and let $B_{16r}(q)\subset\Omega$ and $B_{16r}(q')\subset\Omega'$. Let $\iota$ be a bijective
map of $B_{16r}(q)$ onto $B_{16r}(q')$ that preserves distances, with $\iota(q)=q'$, and
maps $m$ to $c\,m'$ for a constant $c>0$. Then the map
$\iota_*(y,e):=(\iota(y),e\circ\iota^{-1})$ is a homeomorphism of
$\pr^{-1}(B_{12r}(q))$ onto $\pr^{-1}(B_{12r}(q'))$. It commutes with the action of $O(n)$.
\end{lemma}
\begin{proof}
For $y\in B_{12r}(q)$ the map $\iota$ restricts to an isometry of $B_{4r}(y)$ onto
$B_{4r}(\iota(y))$ that multiplies the measures by a constant. The frame operator is
formed with the normalized measure, so \cref{m:fr:rem:local} applies: $e\circ\iota^{-1}$ is a frame at $\iota(y)$, and its frame map is $u_e\circ\iota^{-1}$ on $B_{4r}(\iota(y))$. A sequence of frames converges in a
frame bundle if and only if the base points converge and the frame maps converge
uniformly on the balls of radius $3r$ about them: the fibres over a compact set form a
compact set, and a frame is determined by the restriction of its frame map to the ball
of radius $2r$ about its base point, as noted after \eqref{m:fr:eq:frameconvdef}. So $\iota_*$ and its inverse are continuous. It commutes with $O(n)$, because $(e\cdot A)\circ\iota^{-1}=(e\circ\iota^{-1})\cdot A$.
\end{proof}

We now define the bundle. Let $y\in X$, and let $M_y$ be the ball $B_{R_1}(\tilde y)$
of the local cover of radius $\rho$ at $y$, with the distance $\tilde d$ and the
measure $\tilde m$. It is a metric measure space with a base point, and it is
determined by $y$ up to a unique isomorphism over $X$.

\begin{definition}[Frames of the Local Covers]\label{m:fr:def:coverbundle}
A \emph{frame of the local covers at $y\in X$} is a frame at scale $r$, for the
threshold $\lambda_n$ and in the sense of \cref{m:fr:def:frames}, of the metric measure space
$M_y$ at its base point $\tilde y$. We
write $\Fr_\rho^{(r)}(X)$, with covering radius $\rho$ and frame scale $r$, for the set of pairs $(y,e)$ of a point $y\in X$ and a frame
$e$ of the local covers at $y$, with the projection $\pr(y,e):=y$ and the right action
$(y,e)\cdot A:=(y,e\cdot A)$ of $O(n)$. For $x\in X$ the \emph{chart map} at $x$ is
\begin{equation*}
\tau_x\colon\Fr^{(r)}(\Omega'_x)\longrightarrow\Fr_\rho^{(r)}(X),\qquad
\tau_x(q,e):=\bigl(\pi_x(q),\,e\circ c_q\bigr),
\end{equation*}
where $c_q\colon B_{R_1}(\tilde y)\to B_{R_1}(q)$ is the restriction of the map of \cref{m:fr:lem:agreement}
for $y=\pi_x(q)$.
\end{definition}

At fixed $r$ we write $\Fr_\rho(X)$ for $\Fr_\rho^{(r)}(X)$, following the scale
convention after \cref{m:fr:def:frames}; the subscript $\rho$ still records the covering radius.

By \cref{m:fr:lem:agreement} and the remark after \cref{m:fr:prop:completion}, $c_q$ is an
isometry of $M_y$ onto the ball $B_{R_1}(q)$ of the chart that preserves the measures,
so by \cref{m:fr:rem:local} the base point of $M_y$ is frame-regular and $e\circ c_q$ is a
frame at it. The frames use the balls of radius $2r$ only, and the radius $R_1$ of
$M_y$ plays no role.

\begin{proposition}[The Bundle]\label{m:fr:prop:coverbundle}
Let $(X,d,m)$ be a compact $\RCD(\kappa,N)$ space with $|\kappa|\rho^2\le1$ and
$(\delta,\rho)$-Reifenberg local covers, where $\delta\le\delta_*$, and let
\mbox{$0<r\le c_*\rho$}.
\begin{enumerate}[label=\textup{(\roman*)}]
\item \emph{Charts.} For every $x\in X$ the map $\tau_x$ commutes with the action of
$O(n)$, its image is $\pr^{-1}(B_{R_0/2}(x))$, and $\tau_x(q,e)=\tau_x(q',e')$ if and only
if $(q',e')=\gamma_*(q,e)$ for some $\gamma\in\Gamma_x$. If $x'\in X$ and
$\tau_x(q,e)=\tau_{x'}(q',e')$, then $e'=e\circ\iota^{-1}$ for the map $\iota$ of
\cref{m:fr:prop:coverscharts}\textup{(iii)}, and $\tau_x=\tau_{x'}\circ\iota_*$ at the frames
whose base points lie in $B_{8r}(q)\cap\Omega'_x$ and are mapped into $\Omega'_{x'}$ by
$\iota$.
\item \emph{Topology.} There is exactly one topology on $\Fr_\rho(X)$ for which all
maps $\tau_x$ are continuous and open. With it $\Fr_\rho(X)$ is a compact metrizable
space, $\pr$ is a principal $O(n)$-bundle, and $\tau_x$ induces a homeomorphism of
$\Gamma_x\backslash\Fr(\Omega'_x)$ onto $\pr^{-1}(B_{R_0/2}(x))$.
\end{enumerate}
\end{proposition}
\begin{proof}
(i) The first assertion is clear, since $(e\cdot A)\circ c_q=(e\circ c_q)\cdot A$. Let
$y\in B_{R_0/2}(x)$, let $e'$ be a frame of the local covers at $y$, and let $q$ be a
point of $Z_x$ over $y$. Then $q\in\Omega'_x$, the tuple $e:=e'\circ c_q^{-1}$ is a frame
at $q$ by \cref{m:fr:rem:local}, and $\tau_x(q,e)=(y,e')$. So the image of $\tau_x$ is
$\pr^{-1}(B_{R_0/2}(x))$. If $q$ and $q'$ lie over $y$ and $\gamma\in\Gamma_x$ maps $q$ to
$q'$, then $c_{q'}=\gamma\circ c_q$ by the uniqueness in \cref{m:fr:lem:agreement}, so
$e\circ c_q=e'\circ c_{q'}$ if and only if $e'=e\circ\gamma^{-1}$. In the same way, if
$q'\in\Omega'_{x'}$ lies over $y$ and $\iota$ is the map of
\cref{m:fr:prop:coverscharts}(iii), then $c_{q'}=\iota\circ c_q$, which gives
$e'=e\circ\iota^{-1}$. For $q_2\in B_{8r}(q)$ the maps $c_{\iota(q_2)}$ and
$\iota\circ c_{q_2}$ are lifts of the projection of the local cover of radius $\rho$ at
$\pi_x(q_2)$ under $\pi_{x'}$ that agree at the base point, so they agree on the ball of
radius $\tfrac12R_1$ about it. Hence
$\tau_{x'}(\iota(q_2),e_2\circ\iota^{-1})=\tau_x(q_2,e_2)$ for every frame $e_2$ at $q_2$,
if $q_2\in\Omega'_x$ and $\iota(q_2)\in\Omega'_{x'}$.

(ii) We call a subset of $\Fr_\rho(X)$ open if its preimages under all maps $\tau_x$
are open. Then the maps $\tau_x$ are continuous. They are open: if
$V\subset\Fr(\Omega'_x)$ is open, then by (i) the preimage of $\tau_x(V)$ under
$\tau_{x'}$ is a union of images of open subsets of the open set $\Gamma_xV$ under maps
$\iota_*$, which are homeomorphisms by \cref{m:fr:lem:localnatural}; the balls
$B_{16r}(q)$ with $q\in\Omega'_x$ lie in $\Omega_x$. A topology with these
two properties is the one just defined, because the images of the maps $\tau_x$ cover
$\Fr_\rho(X)$. By (i), $\tau_x$ induces a continuous open bijection of
$\Gamma_x\backslash\Fr(\Omega'_x)$ onto the open set $\pr^{-1}(B_{R_0/2}(x))$. The group
$\Gamma_x$ acts freely and properly discontinuously on $\Omega'_x$, and
$\Fr(\Omega'_x)$ is a principal $O(n)$-bundle over $\Omega'_x$ by \cref{m:fr:thm:bundle}, on
which $\Gamma_x$ acts by automorphisms. So the quotient is a principal $O(n)$-bundle over
$\Gamma_x\backslash\Omega'_x$, which $\pi_x$ maps homeomorphically onto $B_{R_0/2}(x)$ by
\cref{m:fr:prop:coverscharts}(i). Hence $\pr$ is a principal $O(n)$-bundle over the compact
metrizable space $X$, and its total space is compact and metrizable.
\end{proof}
\subsubsection{Collapsing Sequences}\label{m:fr:sec:covercollapse}
We now let the spaces vary. Along a convergent sequence of spaces with Reifenberg local
covers the charts converge, after passing to a subsequence, to spaces with groups that
need not be discrete. The limits of the
charts are glued by the limits of the transition maps to a compact space $Y$ over the
limit space, on which $O(n)$ acts.
The space $Y$ depends on the subsequence (\cref{m:fr:ex:interleave}). No collapse is
assumed.

Assume \cref{hyp:sequence}, with $\delta\le\delta_*$ and $0<r\le c_*\rho$.
Choose Borel $\eps_i$-GH approximations $\phi_i\colon X_i\to X$, with
$\eps_i\to0$, such that the images of the measures $m_i$ converge weakly to $m$. Then
$(X,d,m)$ is an $\RCD(\kappa,N)$ space \cite{GMS15,EKS15}. It need not have Reifenberg local
covers. Points $x_i\in X_i$ \emph{converge} to $x\in X$ if $\phi_i(x_i)\to x$. We write $Y_i:=\Fr_\rho(X_i)$ for the frame spaces of \cref{m:fr:sec:coverbundle}. We fix finitely many points of $X$, which we call \emph{centres}, such that
the balls of radius $\tfrac18R_0$ about them cover $X$, and for every centre $a$ points
$a_i\in X_i$ that converge to $a$. We denote centres by $a$, $b$ and $c$, and a sum over a centre runs over all centres. For the chart of $X_i$ at $a_i$ we
write $Z^a_i$, $z^a_i$, $\Gamma^a_i$, $\pi^a_i$, $\Omega^a_i$, $\Omega'^a_i$ and $\tau^a_i$, and we
normalize its measure: $\hat m^a_i$ is the measure of the chart divided by the measure
of the ball of radius $R_0$ about $z^a_i$, which we denote by $1/\rc^a_i$. As in \cref{m:fr:sec:quotientenergy}, $\varrho(q)$ denotes the minimal distance of the points of the orbit of $q$ from the base point, and $Z(R)$ is the set where $\varrho<R$. 

We first pass to a subsequence along which, for every centre $a$, the charts at the points
$a_i$ converge as in \ref{fl:S:limit}: $(Z^a_i,z^a_i,\Gamma^a_i)\eGHto(Z^a,z^a,G^a)$. We write $\hat m^a$ for the
limit measure, $\Omega^a:=Z^a(R_0)$ for the region and $\pi^a\colon\Omega^a\to X$ for the limit
projection. The measures converge as in \cref{m:fr:prop:localization}\textup{(ii)},
and the charts are in the situation of \cref{m:fr:sec:regions}, with $R_\Omega=R_0$.

\begin{proof}[Proof of \ref{fl:S:limit}, except the Lie-group assertion]
\emph{The limits of the charts.} By \cref{m:fr:prop:coverscharts}(i) the charts are $\RCD(-C_0\rho^{-2},N+1)$ spaces with
discrete groups of isometries that preserve the measures. By the compactness theorem
for such spaces \cite{GMS15}, a subsequence converges in the pointed measured
Gromov--Hausdorff sense. The approximations can be chosen Borel by replacing each
with a map constant on the cells of a countable Borel partition of vanishing mesh;
this changes the approximation error by a quantity tending to zero. By \cite{FY92}
a further subsequence of the groups
converges for these approximations. The limit is an $\RCD$ space with the same constants
by the stability theorem, and limits of isometries that preserve the measures $\hat m^a_i$
preserve $\hat m^a$. As in \cref{m:fr:sec:regions}, the functions $\varrho$ of the charts
converge continuously to that of the limit. So a point $q\in\Omega^a$ is the limit of
points $q_i\in\Omega^a_i$. By \cref{m:fr:prop:coverscharts}(ii) the balls $B_t(q_i)$ with
$t\le R_1$ are $\delta t$-close to $B_t(0^n)$, and they converge to $B_t(q)$ in the
Gromov--Hausdorff sense. The composition of an approximation of $B_t(q)$ by $B_t(q_i)$
whose error tends to $0$ with a $\delta t$-GH approximation of $B_t(q_i)$ changes
distances by at most $\delta t+o(1)$, and its image is dense up to $2\delta t+o(1)$. Hence
$B_t(q)$ is $3\delta t$-close to $B_t(0^n)$, and since
$3\delta\le\delta_H(n,N+1,\eps_\Omega)$ the second step of the proof of
\cref{m:fr:prop:coverscharts}(ii) shows that $q$ is $\eps_\Omega$-flat at scale $r$.

\emph{The projections.} The maps $\pi^a_i$ are $1$-Lipschitz on $\Omega^a_i$. By the theorem of Arzel\`a and
Ascoli, applied on a countable dense subset of $\Omega^a$, we may pass to a subsequence
for which the limit $\pi^a(q)$ exists for all $q\in\Omega^a$ and all sequences
$q_i\to q$. By \cref{m:fr:prop:coverscharts}(i),
$d_i(\pi^a_i(q_i),\pi^a_i(q'_i))$ is the minimum of $\hat d(q_i,\gamma q'_i)$ over
$\gamma\in\Gamma^a_i$, and these minima converge to the minimum of $\hat d(q,gq')$ over
$g\in G^a$, by the argument for the functions $\varrho$ in \cref{m:fr:sec:regions}. This
gives the formula. Since $\pi^a_i(z^a_i)=a_i$, we get $\pi^a(z^a)=a$ and
$d(\pi^a(q),a)=\varrho(q)$, so $\pi^a$ maps $Z^a(R)$ into $B_R(a)$. If $y\in B_R(a)$ and
$y_i\in X_i$ converge to $y$, there are points $q_i$ over $y_i$ with
$\hat d(q_i,z^a_i)=d_i(y_i,a_i)$, a subsequence of them converges to a point
$q\in Z^a(R)$, and $\pi^a(q)=y$.
\end{proof}

For the rest of \cref{m:book:ch:10} we keep this subsequence, and we put
${\Omega'^a:=Z^a(\tfrac12R_0)}$. We write $\chi_R$ for
the cut-off functions of \cref{m:fr:lem:groupunity} for the sequence of the charts at the
points $a_i$ and its limit.

The charts and their limits carry the frame trivializations of \ref{fl:S:product}. They are the
coordinates in which small deck transformations are compared with lattice elements and boxes
in \cref{ch:lattices}. In the proof the centre is suppressed, as in the theorem.

\begin{proof}[Proof of \ref{fl:S:product}]
By \cref{m:fr:sec:frameenergy} the frame map $u_e$ is a gauge map on $B_{3r}(q)$. The ball
$B_{r/2}(q)$ lies in $\Omega'_i$, respectively in $\Omega'$, because
$\varrho_i(q)<\tfrac14R_0$ and $r<\tfrac14R_0$. The trivialization $\ft_w$ defined by a gauge map $w$ is the
homeomorphism $(x,A)\mapsto(x,e_w(x)\cdot A)$ over the core of its domain
(\cref{m:fr:sec:frameenergy}), and the core of $B_{3r}(q)$ contains $B_{r/2}(q)$. Every frame is
aligned with its own frame map, $e_{u_e}(q)=e$ (\cref{m:fr:thm:gauges}), which
gives $\ft_p(q,\one)=p$. The identity $\ft_p(x,AA')=\ft_p(x,A)\cdot A'$ holds by the definition.

For the convergence, the trivializations $\ft_{p_i}$ and $\ft_p$ are corresponding
trivializations in the sense of \cref{m:fr:def:corrtriv}, and the two convergence assertions are
parts~(i) and~(ii) of \cref{m:fr:lem:corrcharts}. They hold for frames over the flat regions by
\cref{m:fr:prop:localization}\textup{(iii)}, because the distance of $q$ from $Z\setminus\Omega$
is larger than~$r$.
\end{proof}

The limits of the charts have to be glued, by limits of the transitions between the charts
of $X_i$.

\begin{definition}[Limit Transitions]\label{m:fr:def:limittransition}
For centres $a$ and $b$ let $\cT^{ba}_i$ be the set of the maps $\iota$ of
\cref{m:fr:prop:coverscharts}(iii) for the charts of $X_i$ at $a_i$ and at $b_i$ and for points
$q\in\Omega'^a_i$ and $q'\in\Omega'^b_i$ over the same point of $X_i$; we say that $\iota$
is \emph{based} at $q$. A sequence of maps $\iota_i\in\cT^{ba}_i$ based at points $q_i$,
indexed by an infinite set of indices, \emph{converges} to a map
$\iota\colon B_{R_1}(q)\to Z^b$ if $q_i\to q$ and $\iota_i(q'_i)\to\iota(q')$ whenever
points $q'_i\in B_{R_1}(q_i)$ converge to $q'\in B_{R_1}(q)$. A \emph{limit transition}
from $a$ to $b$, based at $q\in\Omega'^a$, is the limit $\iota$ of such a sequence with
$\iota(q)\in\Omega'^b$.
\end{definition}

The next lemma says that limit transitions exist, that they are
unique up to the groups, and that they can be inverted and composed.

\begin{lemma}[Limits of Transitions]\label{m:fr:lem:limittransitions}
After passing to a further subsequence, the following holds for all centres $a$, $b$
and $c$.
\begin{enumerate}[label=\textup{(\roman*)}]
\item Let $\iota_i\in\cT^{ba}_i$ be maps whose base points $q_i$ converge to a point
$q\in Z^a$ and for which the points $\iota_i(q_i)$ stay at bounded distance from the
base points $z^b_i$. Then a subsequence converges. Its limit $\iota$ preserves
distances, maps $B_{R_1}(q)$ onto $B_{R_1}(\iota(q))$, satisfies
$\pi^b\circ\iota=\pi^a$ at the points of $B_{R_1}(q)\cap\Omega^a$ that it maps into
$\Omega^b$, and maps $\hat m^a$ to $\rc^{ba}\,\hat m^b$ for a constant $\rc^{ba}>0$
that depends only on $a$ and $b$. If frames $(q_i,e_i)$ converge to
$(q,e)$, then the frames $(\iota_i(q_i),e_i\circ\iota_i^{-1})$ converge to
$(\iota(q),e\circ\iota^{-1})$ along the subsequence.
\item If $q\in\Omega'^a$ and $\pi^a(q)\in B_{R_0/2}(b)$, there is a limit transition from
$a$ to $b$ based at $q$, as a limit along a subsequence of any given infinite set of
indices. If $\iota$ and $\iota'$ are two limit transitions from $a$ to $b$ based at
$q$, then $\iota'=g\circ\iota$ for some $g\in G^b$. If $\iota$ is one, $g\in G^b$ and
$g'\in G^a$, then $g\circ\iota\circ g'^{-1}$ is a limit transition based at $g'q$.
\item The limit transitions from $a$ to $a$ are the restrictions of the elements of
$G^a$ to the balls of radius $R_1$ about the points of $\Omega'^a$. The inverse of a
limit transition from $a$ to $b$ is a limit transition from $b$ to $a$. If $\iota$ is a
limit transition from $a$ to $b$ based at $q$ and $\iota'$ is one from $b$ to $c$ based
at $\iota(q)$, then $\iota'\circ\iota$ agrees on $B_{R_1/2}(q)$ with a limit
transition from $a$ to $c$ based at $q$.
\end{enumerate}
\end{lemma}
\begin{proof}
We first pass to a subsequence along which the quotients $\rc^a_i/\rc^b_i$ converge in
$[0,\infty]$ for all $a$ and $b$.

(i) The maps $\iota_i$ preserve distances, between balls of radius $R_1$ in proper
spaces that converge. By the theorem of Arzel\`a and Ascoli a subsequence converges,
together with the inverse maps. The limit preserves distances and is onto. If
$q'\in B_{R_1}(q)\cap\Omega^a$ and $\iota(q')\in\Omega^b$, and $q'_i\to q'$, then
$q'_i\in\Omega^a_i$ and $\iota_i(q'_i)\in\Omega^b_i$ for large $i$, and
$\pi^b_i(\iota_i(q'_i))=\pi^a_i(q'_i)$; so $\pi^b(\iota(q'))=\pi^a(q')$ by \ref{fl:S:limit}. The map $\iota_i$ sends $\hat m^a_i$ to $\rc^a_i/\rc^b_i$ times
$\hat m^b_i$. The measures of the balls $B_{R_1}(q_i)$ and of their images converge to
the measures of $B_{R_1}(q)$ and of its image, which are positive and finite, because
spheres are null sets. So the limit $\rc^{ba}$ of $\rc^a_i/\rc^b_i$ is positive and
finite, and $\iota$ sends $\hat m^a$ to $\rc^{ba}\,\hat m^b$. For the frames we argue
as in the proof of \cref{m:fr:thm:frameconv}(iv): by \cref{m:fr:rem:local} the frame map of
$e_i\circ\iota_i^{-1}$ is $u_{e_i}\circ\iota_i^{-1}$ on the ball of radius $4r$ about
$\iota_i(q_i)$, the maps $\iota_i^{-1}$ converge uniformly to $\iota^{-1}$, and
$u_e\circ\iota^{-1}$ is the frame map of $e\circ\iota^{-1}$ on the ball of radius $4r$
about $\iota(q)$.

(ii) Let $q_i\in Z^a_i$ converge to $q$. For large $i$ we have $q_i\in\Omega'^a_i$, and
the points $y_i:=\pi^a_i(q_i)$, which converge to $\pi^a(q)$, lie in $B_{R_0/2}(b_i)$. By
\cref{m:fr:prop:coverscharts}(i) there are points $s_i\in\Omega'^b_i$ over $y_i$ with
$\hat d(s_i,z^b_i)=d_i(y_i,b_i)$, and maps $\iota_i\in\cT^{ba}_i$ based at $q_i$ with
$\iota_i(q_i)=s_i$. By (i) a subsequence converges to a map $\iota$, and
$\hat d(\iota(q),z^b)=d(\pi^a(q),b)<\tfrac12R_0$. So $\iota$ is a limit transition.

For every pair $(a,b)$ we fix a countable dense subset of the set of the points
$q\in\Omega'^a$ with $\pi^a(q)\in B_{R_0/2}(b)$ and, for each of its points, a sequence
$\iota_i$ as above. By (i) and a diagonal argument we may pass to a subsequence along
which all these sequences converge; we call their limits \emph{reference transitions}.
Let now $\iota$ be a limit transition from $a$ to $b$ based at $q$, the limit of maps
$\iota_i$ based at $q_i$ along an infinite set of indices. Choose a reference transition
$\iota^*$, the limit of maps $\iota^*_i$ based at $q^*_i$, whose base point $q^*$ has
distance less than $\tfrac14R_1$ from $q$. For large $i$ the point $q_i$ lies in
$B_{R_1/2}(q^*_i)$, and by the uniqueness in \cref{m:fr:prop:coverscharts}(iii),
$\iota_i=\gamma_i\circ\iota^*_i$ on $B_{R_1/2}(q_i)$ for the element
$\gamma_i\in\Gamma^b_i$ that maps $\iota^*_i(q_i)$ to $\iota_i(q_i)$. The displacements of
$\gamma_i$ at the base points are bounded, so by \cref{m:fr:lem:groupconv} a subsequence
converges to some $g\in G^b$, and $\iota=g\circ\iota^*$ on $B_{R_1/2}(q)$. Hence two
limit transitions $\iota$ and $\iota'$ based at $q$ satisfy $\iota'=g\circ\iota$ on
$B_{R_1/2}(q)$ for some $g\in G^b$. Both sides preserve distances on $B_{R_1}(q)$,
so they map the geodesics that start at $q$ to geodesics of $Z^b$, and these agree
near their initial points. Since $Z^b$ is non-branching by \ref{m:fr:R:isometries}, the two
maps agree on $B_{R_1}(q)$. For the last assertion of (ii) we note that, by
\cref{m:fr:lem:groupconv}(iii)--(v), the maps
$\beta_i(g)\circ\iota_i\circ\beta_i(g')^{-1}$, which belong to $\cT^{ba}_i$, converge to
$g\circ\iota\circ g'^{-1}$.

(iii) The elements of $\cT^{aa}_i$ are restrictions of elements of $\Gamma^a_i$, by
\cref{m:fr:prop:coverscharts}(iii). If they converge, the displacements of these elements at
the base points are bounded, a subsequence of them converges to some $g\in G^a$ by
\cref{m:fr:lem:groupconv}, and the limit is the restriction of $g$; conversely the
restrictions of $\beta_i(g)$ converge to that of $g$. If $\iota_i\to\iota$, then
$\iota_i^{-1}\to\iota^{-1}$, which gives the second assertion. For the third one let
$\iota$ be the limit of maps $\iota_i$ along an infinite set of indices. By (ii) there
is a limit transition $\iota''$ from $b$ to $c$ based at $\iota(q)$ that is the limit
of maps $\iota''_i\in\cT^{cb}_i$ along a subsequence of this set, and
$\iota'=g\circ\iota''$ for some $g\in G^c$. For large $i$ the map
$\iota''_i\circ\iota_i$ is defined on $B_{R_1/2}(q_i)$, and it agrees there with the
element of $\cT^{ca}_i$ based at $q_i$ that has the same value at $q_i$; this element
exists because the value lies over a point that converges to
$\pi^a(q)=\pi^c(\iota'(\iota(q)))$, which lies in $B_{R_0/2}(c)$. A subsequence
of these elements converges to a limit transition $\iota'''$ from $a$ to $c$ with
$\iota'''=\iota''\circ\iota$ on $B_{R_1/2}(q)$, and $g\circ\iota'''$ is a limit
transition by (ii) that agrees with $\iota'\circ\iota$ on this ball.
\end{proof}

We pass to the subsequence of this lemma. The limit of the frame spaces is obtained by gluing
the frame bundles of the limit charts along the limit transitions; by
\cref{m:fr:lem:limittransitions}(iii) the relation in the following definition is an
equivalence relation.

\begin{definition}[The Limit of the Frame Spaces]\label{m:fr:def:Y}
On the disjoint union of the frame bundles $\Fr(\Omega'^a)$ of all centres $a$, a frame $(q,e)$
over $\Omega'^a$ and a frame $(q',e')$ over $\Omega'^b$ are \emph{equivalent} if there is a limit
transition $\iota$ from $a$ to $b$ based at $q$ with $\iota(q)=q'$ and $e'=e\circ\iota^{-1}$. The
\emph{limit of the frame spaces} $Y_i$ is the set $Y$ of equivalence classes with the quotient
topology, and
\begin{equation*}
\tau^a\colon\Fr(\Omega'^a)\longrightarrow Y
\end{equation*}
are the restrictions of the quotient map.
\end{definition}

Two frames over $\Omega'^a$ are equivalent if and only if they lie in one orbit of $G^a$.
Points of the spaces $Y_i$ converge to a point of $Y$ when they come from convergent frames
of one chart.

\begin{definition}[Convergence in the Frame Spaces]\label{m:fr:def:Yconv}
A sequence of points $y_i\in Y_i$ \emph{converges} to $y\in Y$ if there are a centre $a$ and
frames $(q_i,e_i)\in\Fr(\Omega'^a_i)$ that converge to a frame $(q,e)\in\Fr(\Omega'^a)$, in the
sense of \cref{m:fr:def:frameconv} for the charts, such that $y_i=\tau^a_i(q_i,e_i)$ for
large $i$ and $y=\tau^a(q,e)$.
\end{definition}

Convergent frames can be chosen in the orbits of the deck groups as soon as they can be chosen
along subsequences.

\begin{lemma}[Selection in Orbits]\label{m:fr:lem:orbitselect}
Let $a$ be a centre, $p=(q,e)$ a frame over $\Omega'^a$, $\mathcal I$ an infinite set of
indices, and $p'_i\in\Fr(\Omega'^a_i)$ for $i\in\mathcal I$. If every infinite subset of $\mathcal I$
contains an infinite subset along which suitable frames in the orbits $\Gamma^a_ip'_i$ converge
to $p$, then there are frames $p_i\in\Gamma^a_ip'_i$ that converge to $p$ along $\mathcal I$.
\end{lemma}
\begin{proof}
For a frame $p''=(q'',e'')$ over $\Omega'^a_i$ put
\begin{equation*}
D_i(p''):=d_{\cA}(q'',q)+\sup_{B_{3r}(q'')}\bigl|u^i_{e''}-u_e\circ \af_i\bigr|,
\end{equation*}
where $d_{\cA}$ is the distance of \cref{m:fr:sec:funcconv} on the disjoint union of $Z^a$ and
the spaces $Z^a_i$, for which $q''_i\to q$ if and only if $d_{\cA}(q''_i,q)\to0$, and the
supremum is taken to be $+\infty$ when $B_{3r}(q'')$ does not lie in the domain of $\af_i$. By
\eqref{m:fr:eq:frameconvdef}, frames $p''_i$ converge to $p$ if and only if $D_i(p''_i)\to0$. By
hypothesis every infinite subset of $\mathcal I$ contains an infinite subset along which the
infima of $D_i$ over the orbits $\Gamma^a_ip'_i$ tend to $0$, so these infima tend to $0$ along
$\mathcal I$. Frames $p_i\in\Gamma^a_ip'_i$ whose values of $D_i$ lie within $1/i$ of the infima
converge to $p$.
\end{proof}

Put $\pr(\tau^a(q,e)):=\pi^a(q)$, and let $O(n)$ act on $Y$ through its right action on the
frames. We now prove the assertions of \cref{m:fr:thm:geomain} on $Y$: those of the preamble, those
on $\tau$ in \ref{fl:S:charts}, \ref{fl:S:Yconv} except its last sentence, and the freeness of $G$ in
\ref{fl:S:frames}.

\begin{proof}[Proof of the assertions on $Y$]
\emph{The space $Y$ and the maps $\tau$.} Let $V\subset\Fr(\Omega'^a)$ be open, and let $(q',e')$ be a frame over $\Omega'^b$
that is equivalent to a frame $(q,e)\in V$, through a limit transition $\iota$, the
limit of maps $\iota_i$ based at points $q_i$. Let $q_2\in B_{8r}(q)\cap\Omega'^a$ with
$\iota(q_2)\in\Omega'^b$, and let $q_{2,i}\to q_2$. For large $i$ the element of
$\cT^{ba}_i$ based at $q_{2,i}$ with the value $\iota_i(q_{2,i})$ agrees with $\iota_i$
on $B_{R_1/2}(q_{2,i})$, and by \cref{m:fr:lem:limittransitions}(i) a subsequence of
these elements converges to a limit transition based at $q_2$ that agrees with $\iota$
on $B_{R_1/2}(q_2)$. So the frames $(\iota(q_2),e_2\circ\iota^{-1})$ are equivalent to
the frames $(q_2,e_2)$, which is the assertion on limit transitions in \ref{fl:S:charts}. By \cref{m:fr:lem:localnatural}
the map
$\iota_*$ is a homeomorphism near $(q,e)$. Hence the set of the frames over $\Omega'^b$
that are equivalent to frames in $V$ is open, and the maps $\tau^a$ are open. They are
continuous by the definition of the topology. Since the frames over $\Omega'^a$ that are
equivalent are those in one orbit of $G^a$, the map $\tau^a$ induces a homeomorphism of
$G^a\backslash\Fr(\Omega'^a)$ onto an open subset of $Y$. The action of $G^a$ on
$\Fr(\Omega'^a)$ is proper, because $G^a$ is a closed group of isometries and the
fibres of $\pr$ are compact; so this quotient is a Hausdorff space.

By \cref{m:fr:lem:limittransitions}(i), $\pr$ is well defined and continuous.
It is onto by \ref{fl:S:limit}. Let $(q,e)$ and $(q',e')$ be frames over
$\Omega'^a$ and $\Omega'^b$ with $\pi^a(q)=\pi^b(q')$. By \cref{m:fr:lem:limittransitions}(ii)
there is a limit transition $\iota$ from $a$ to $b$ based at $q$, and
$\pi^b(\iota(q))=\pi^b(q')$. By the formula for $\pi^b$ in \ref{fl:S:limit} the points
$\iota(q)$ and $q'$ lie in one orbit of $G^b$, so after composing $\iota$ with an
element of $G^b$ we have $\iota(q)=q'$. Then $e\circ\iota^{-1}$ and $e'$ are frames at
$q'$, which differ by an element of $O(n)$. So the two classes lie in one orbit of
$O(n)$: the fibres of $\pr$ are the orbits, and the image of $\tau^a$ is
$\pr^{-1}(B_{R_0/2}(a))$. Two points of $Y$ in different fibres are separated by $\pr$,
and two points of one fibre lie in the image of one map $\tau^a$, which is open and a
Hausdorff space. So $Y$ is a Hausdorff space. The frames over the closed ball of radius
$\tfrac18R_0$ about $z^a$ form a compact set, and by \ref{fl:S:limit} its image
under $\tau^a$ is the preimage under $\pr$ of the closed ball of radius $\tfrac18R_0$
about $a$; these images cover $Y$. So $Y$ is compact, and it has a
countable base. The action of $O(n)$ is continuous because the maps $\tau^a$ are open
and commute with it.

\emph{Convergence.} Let $y_i\to y$ by means of the chart at $a$ and of frames $(q_i,e_i)\to(q,e)$, and
let ${\pr(y)\in B_{R_0/2}(b)}$. As in the proof of \cref{m:fr:lem:limittransitions}(ii) there
are maps $\iota_i\in\cT^{ba}_i$ based at $q_i$, for large $i$, and then
$y_i=\tau^b_i(s_i,e'_i)$ with $s_i:=\iota_i(q_i)$ and $e'_i:=e_i\circ\iota_i^{-1}$, by
\cref{m:fr:prop:coverbundle}(i). By \cref{m:fr:lem:limittransitions}(i), every subsequence has a
subsequence along which the frames $(s_i,e'_i)$ converge to
$(\iota(q),e\circ\iota^{-1})$ for a limit transition $\iota$, that is, to a frame in the
orbit $(\tau^b)^{-1}(y)$ of $G^b$. Fix $(s,e')$ in this orbit. By \cref{m:fr:lem:groupconv}(iii) and~(v), every subsequence admits a
further one along which suitable frames in the orbits $\Gamma^b_i(s_i,e'_i)$ converge to $(s,e')$. So
by \cref{m:fr:lem:orbitselect} there are $\gamma_i\in\Gamma^b_i$ such that $(\gamma_i)_*(s_i,e'_i)$
converges to $(s,e')$, which is the assertion on the charts at $b_i$ and $b$. For the diagonal
sequences, let $y=\tau^a(q,e)$. Since $\tau^a$ is open, $y^k=\tau^a(q^k,e^k)$ for large
$k$, with frames that converge to $(q,e)$, and by what we have shown
$y^k_i=\tau^a_i(q^k_i,e^k_i)$ with frames that converge to $(q^k,e^k)$ as $i\to\infty$.
Choose $i_k$ such that the base points of these frames for $i=i_k$ have $d_{\cA}$-distance at
most $1/k$ from $q^k$ in the space of \cref{m:fr:sec:funcconv}, and such that the supremum in
\eqref{m:fr:eq:frameconvdef} formed with $(q^k,e^k)$ is at most $1/k$. The frame maps of the
frames $e^k$ converge to that of $e$ uniformly on bounded sets, so the frames
$(q^k_{i_k},e^k_{i_k})$ converge to $(q,e)$.

If $y_i\to y$ and $y_i\to y'$, then $\pr(y)=\pr(y')$ by \ref{fl:S:limit}, and
by what we have shown both limits are obtained in one chart, from frames
$(q_i,e_i)\to(q,e)$ and $(q'_i,e'_i)\to(q',e')$ with
$(q'_i,e'_i)=(\gamma_i)_*(q_i,e_i)$ for elements $\gamma_i\in\Gamma^a_i$, by
\cref{m:fr:prop:coverbundle}(i). Their displacements at the base points are bounded, a
subsequence converges to some $g\in G^a$, and $(q',e')=g_*(q,e)$ by
\cref{m:fr:lem:groupconv}(v). So $y=y'$. Given points $y_i\in Y_i$, a subsequence of the
points $\pr(y_i)$ converges to a point of the ball $B_{R_0/8}(a)$ of some centre $a$.
Along it $y_i=\tau^a_i(q_i,e_i)$ with $\hat d(q_i,z^a_i)<\tfrac14R_0$ for large $i$, and
by \cref{m:fr:prop:convergence} a subsequence of these frames converges to a frame over
$\Omega'^a$. Every frame over $\Omega'^a$ is a limit of frames over the sets
$\Omega'^a_i$ by \cref{m:fr:thm:frameconv}(ii), which holds for the regions by
\cref{m:fr:prop:localization}(iii). If $A_i\to A$ in $O(n)$, then $u^i_{e_i\cdot A_i}=A_i^{\T}u^i_{e_i}$ converges to
$A^{\T}u_e=u_{e\cdot A}$ whenever $(q_i,e_i)\to(q,e)$. Passing to classes gives the
statement on $O(n)$. Convergence of the base points follows from \ref{fl:S:limit}.

\emph{Isotropy.} By \ref{fl:S:limit} the region $\Omega^a$ is $\eps_\Omega$-flat at scale $r$,
with $\eps_\Omega\le\eps_P$, and $(3\delta,R_1)$-Reifenberg, hence $(\delta_P,8r)$-Reifenberg,
because $3\delta\le\delta_P$ and $8r\le R_1$. The group
$G^a$ preserves $\hat m^a$ and the set $\Omega^a$. So \cref{m:fr:prop:freeframes} shows that
$G^a$ acts freely on $\Fr(\Omega^a)$ and that the homomorphisms $A_e$ are injective. A
matrix $A$ fixes $\tau^a(q,e)$ if and only if $(q,e\cdot A)=g_*(q,e)$ for some $g\in G^a$,
that is, if and only if $g\in G^a_q$ and $A=A_e(g)$.

\end{proof}

For a Borel function $F$ on $Y$
and a centre $a$ the \emph{lift} $F^a$ is the function on $\Fr(\Omega^a)$ that equals
$F\circ\tau^a$ over $\Omega'^a$ and vanishes elsewhere. Functions on
the spaces $Y_i$ are compared with functions on $Y$ as in \cref{m:fr:sec:funcconv}: a
\emph{test sequence} on the spaces $Y_i$ is one in the sense of \cref{m:fr:def:testseq}, for the
convergence of \cref{m:fr:def:Yconv}, so that $\Theta_i(y_i)\to\Theta(y)$ whenever $y_i\to y$. No support condition is needed, because the
spaces $Y_i$ and $Y$ are compact. The limit of a test sequence is continuous, by the argument of \cref{m:fr:sec:funcconv} and
the assertion on diagonal sequences in \ref{fl:S:Yconv}. Test sequences are obtained from the charts by means of a partition of unity on the
base. For a centre $a$ put $\psi^a:=(\tfrac14R_0-d(a,\cdot))_+$ and
$\chi^a:=\psi^a/\sum_b\psi^b$. Since the balls of radius $\tfrac18R_0$ about the centres
cover $X$, the sum of the functions $\psi^b$ is at least $\tfrac18R_0$; so the functions
$\chi^a$ are continuous, they vanish outside $B_{R_0/4}(a)$, and their sum is $1$.

\begin{proof}[Proof of the last assertion of \ref{fl:S:Yconv}]
We show that every continuous function on $Y$ is the limit of a test sequence. Let $\Theta$ be a continuous function on $Y$, and let $a$ be a centre. We write
$\varrho_i$ and $\varrho$ for the functions of \cref{m:fr:sec:quotientenergy} on $Z^a_i$ and on
$Z^a$. Let $H$ be the lift of $(\chi^a\circ\pr)\,\Theta$ to the chart at $a$. Since
$d(\pi^a(q),a)=\varrho(q)$ by \ref{fl:S:limit} and $\chi^a$ vanishes outside
$B_{R_0/4}(a)$, the function $H$ vanishes over the points with
$\varrho\ge\tfrac14R_0$. So it is continuous on $\Fr(\Omega^a)$, and it is invariant
under $G^a$, because the fibres of $\tau^a$ are the orbits of $G^a$. Put
$U:=(\chi_{R_0/4}\circ\pr)\,H$. This function is continuous, and it vanishes outside the
preimage of the set of the points $q$ of a closed ball about $z^a$ with
$\varrho(q)\le\tfrac14R_0$. This set is compact, and its distance from the complement of
$\Omega^a$ is at least $\tfrac34R_0$, because $\varrho$ is $1$-Lipschitz. If
$\varrho(\pr(p))<\tfrac14R_0$, then by the invariance of $H$, by
\eqref{m:fr:eq:regionalunity} and by the invariance of the Haar measure $\nu_{G^a}$ of $G^a$ under inversion,
\begin{equation*}
\bar U(p)=H(p)\int_{G^a}\chi_{R_0/4}\bigl(g^{-1}\pr(p)\bigr)\,d\nu_{G^a}(g)=H(p),
\end{equation*}
and at the other frames both $\bar U$ and $H$ vanish. So $\bar U=H$.

By \cref{m:fr:prop:localization}(ii) the function $U$ is the limit of a test sequence $U_i$ with
admissible supports on the frame bundles of the regions $\Omega^a_i$, and by \cref{m:fr:lem:groupaverages} for the
regions the averages $\bar U_i$ are continuous, invariant under $\Gamma^a_i$ and bounded by
one constant, and they converge continuously to $H$. Let $\zeta$ be a continuous function on $[0,\infty)$ with values in $[0,1]$ that equals $1$ on $[0,\tfrac14R_0]$ and vanishes on
$[\tfrac38R_0,\infty)$, and put $V_i:=(\zeta\circ\varrho_i\circ\pr)\,\bar U_i$. The functions
$\varrho_i$ converge continuously to $\varrho$, and $H$ vanishes where
$\varrho\ge\tfrac14R_0$; so the functions $V_i$ converge continuously to $H$ as well. They
are continuous and invariant under $\Gamma^a_i$, and they vanish over the points with
$\varrho_i\ge\tfrac38R_0$. By \cref{m:fr:prop:coverbundle}(i) and~(ii) the restriction of
$V_i$ to $\Fr(\Omega'^a_i)$ is the composition with $\tau^a_i$ of a continuous function
on $\pr^{-1}(B_{R_0/2}(a_i))$, and this function vanishes outside the preimage of the
closed ball of radius $\tfrac38R_0$ about $a_i$, because the distance of
$\pi^a_i(q)$ from $a_i$ is $\varrho_i(q)$ by \cref{m:fr:prop:coverscharts}(i). Extended by zero,
it is a continuous function $\Theta^a_i$ on $Y_i$, and these functions are bounded by one
constant.

Let $y_i\to y$. We claim that $\Theta^a_i(y_i)$ converges to
$\chi^a(\pr(y))\,\Theta(y)$. If $\pr(y)\in B_{R_0/2}(a)$, then by \ref{fl:S:Yconv} there are frames $p_i\in\Fr(\Omega'^a_i)$ with
$\tau^a_i(p_i)=y_i$ for large $i$ that converge to a frame $p$ with $\tau^a(p)=y$, and
$\Theta^a_i(y_i)=V_i(p_i)\to H(p)$, which is the assertion. Otherwise
$d(\pr(y),a)\ge\tfrac12R_0$. Since $\pr(y_i)\to\pr(y)$, the distance of $\pr(y_i)$ from
$a_i$ is larger than $\tfrac38R_0$ for large $i$, so $\Theta^a_i(y_i)=0$, and
$\chi^a(\pr(y))=0$. The sum of the functions $\chi^a$ over the centres is $1$, so the
functions $\Theta_i:=\sum_a\Theta^a_i$ form a test sequence with limit $\Theta$.
\end{proof}

By \ref{fl:S:charts} the isotropy groups of $O(n)$ on $Y$ are those of the limiting chart groups
$G^a$.

The passage to a subsequence in the proof of \ref{fl:S:limit} and in
\cref{m:fr:lem:limittransitions} cannot be
avoided, even when the spaces converge. The reason is that only the convergence of the
spaces $X_i$ is assumed, and not, as in \cref{m:fr:sec:quotientenergy}, the equivariant convergence of
covers.

\begin{example}[The Limit Depends on the Subsequence]\label{m:fr:ex:interleave}
Let $n=N=2$ and $\kappa=0$, let $\delta_i\to0$, and let $X_i$ be the flat torus
$\R^2/\delta_i\mathbb{Z}^2$ for even $i$, and for odd $i$ the flat Klein bottle that is the
quotient of $\R^2$ by the group generated by $(x,y)\mapsto(x+\delta_i,-y)$ and
$(x,y)\mapsto(x,y+\delta_i)$, both with the normalized area. In a local cover of a
complete flat manifold the balls about the base point whose radii are at most a third
of the radius of the cover are Euclidean balls, so these spaces have
$(\delta,\rho)$-Reifenberg local covers for every $\delta,\rho>0$, and they converge to a point. The frames of the local
covers are the orthonormal frames, and $Y_i$ is the quotient of $\R^2\times O(2)$ by the
group. For even $i$ it is $X_i\times O(2)$, which has two components. For odd $i$ it is
connected, because the first generator reverses the orientation.
Along the even indices the groups of the charts converge to the group of the
translations and $Y=O(2)$; along the odd indices they converge to the group generated
by the translations and a reflection, and $Y$ is the quotient of $O(2)$ by the left
multiplication with a reflection, a circle.
\end{example}

\subsection{The Limit Groups}\label{m:book:ch:11}\label{m:fr:sec:fields}\label{m:fr:sec:limitgroup}

The groups of deck transformations of the charts converge to closed groups of isometries
of the limits of the charts. This subsection records the structure of these limit groups
(\cref{m:fr:prop:limitgroup}): they are Lie groups that act properly, their identity
components are nilpotent by a result of J.~Wang, their isotropy groups are compact, and
near every point the quotient by the limit group is the quotient by an open subgroup that
is generated by small elements.

Throughout this subsection we are in the group setting of \cref{m:fr:sec:quotientenergy},
without flat regions: \cref{m:fr:lem:groupconv}(i)--(v) holds for points, the functions
$\varrho_i$ converge continuously to $\varrho$ (\cref{m:fr:sec:regions}), and
\cref{m:fr:lem:counting} holds with the numbers $\ell_i$ of \eqref{m:fr:eq:elldef} and the
Haar measure $\nu_G$ for which $\int_G\hat\chi_1(gz)\,d\nu_G(g)=1$, with constants $C(R)$
that depend on $N$, $\kappa$ and $R$. We write $G_0$ for the identity component of $G$, and $\cA$ for the space of
\cref{m:fr:sec:funcconv} formed with $Z$ and the spaces $Z_i$.

\begin{proposition}[The Limit Group]\label{m:fr:prop:limitgroup}
\begin{enumerate}[label=\textup{(\roman*)}]
\item \emph{Nilpotency.} $G$ is a unimodular Lie group that acts properly on $Z$ and
preserves $m$, and $G_0$ is nilpotent.
\item \emph{The central torus.} $G_0$ has a largest compact subgroup $\mathbb T$. It is a torus
in the centre of $G_0$ and a normal subgroup of $G$, the quotient $G_0/\mathbb T$ is a simply
connected nilpotent Lie group, and $G_0/\mathbb T$ acts freely on $Z/\mathbb T$.
\item \emph{Isotropy.} For every $q\in Z$ the isotropy group $G_q$ is compact,
$G_q\cap G_0$ is a subgroup of $\mathbb T$ that fixes every point of the orbit $G_0q$, and the
quotient of $G_q$ by $G_q\cap G_0$ is a finite group that embeds into $G/G_0$. The orbit
$G_0q$ is homeomorphic, by the orbit map, to the connected nilpotent Lie group
$G_0/(G_q\cap G_0)$.
\item \emph{The quotient near a point.} For every $q\in Z$ the set $G_0G_q$ is an open
subgroup of $G$, and there is $\eps_q>0$ with the following properties for
$0<\eps\le\eps_q$. The group $G_0G_q$ is generated by the elements $g\in G$ with
$d(gq,q)<\eps$. If $g\in G$ and $p,p'\in B_{\eps/4}(q)$ with $gp=p'$, then $g\in G_0G_q$.
The projection of $Z/(G_0G_q)$ onto $Z/G$ maps the ball of radius $\tfrac14\eps$ about
the class of $q$ isometrically onto the ball of radius $\tfrac14\eps$ about the class of
$q$.
\end{enumerate}
\end{proposition}
\begin{proof}
(i) By \ref{m:fr:R:isometries}, $\Isom(Z)$ is a Lie group whose closed subgroups act properly,
so $G$ is a Lie group that acts properly. It preserves $m$ by \cref{m:fr:lem:groupconv}\textup{(v)}, and it
is unimodular by \cref{m:fr:lem:counting}(ii). The spaces $Z_i$ are proper and geodesic by
\ref{m:fr:R:geometry}, the groups $\Gamma_i$ are discrete subgroups of $\Isom(Z_i)$, and
their limit $G$ is a Lie group. Hence $G_0$ is nilpotent by \cite[Lemma~5.2]{Wan23}.

(ii) By \cite[Corollary~2.5]{Wan23}, every compact subgroup of the connected nilpotent
Lie group $G_0$ lies in its centre, there is a compact subgroup $\mathbb T$ that contains all
others, $\mathbb T$ is a torus, and $G_0/\mathbb T$ is diffeomorphic to a Euclidean space. The
conjugation by an element of $G$ is an automorphism of $G_0$, and it maps $\mathbb T$ to a
compact subgroup, which lies in $\mathbb T$; applied to the inverse element this shows that $\mathbb T$
is normal in $G$. Let $g\in G_0$ and $q\in Z$ with $gq\in \mathbb Tq$. Then $t^{-1}g$ fixes $q$
for some $t\in \mathbb T$, so it lies in the compact subgroup $G_q\cap G_0$ of $G_0$, hence in
$\mathbb T$, and $g\in \mathbb T$. So $G_0/\mathbb T$ acts freely on $Z/\mathbb T$.

(iii) The group $G_q$ is compact because the action is proper. Its subgroup
$G_q\cap G_0$ is a compact subgroup of $G_0$, so by (ii) it lies in $\mathbb T$ and in the
centre of $G_0$. For $h\in G_q\cap G_0$ and $g\in G_0$ this gives $h(gq)=g(hq)=gq$.
Since $G_0$ is open in $G$, the group $G_q\cap G_0$ is open in the compact group $G_q$;
so it has finite index, and the quotient embeds into $G/G_0$. A central subgroup is
normal, so $G_0/(G_q\cap G_0)$ is a connected nilpotent Lie group. The orbit map from it
to $G_0q$ is continuous and bijective, and it is proper because $G_0$ acts properly; so
it is a homeomorphism.

(iv) The set $G_0G_q$ is a subgroup because $G_0$ is normal, and it is open because it
is a union of cosets of the open subgroup $G_0$. For $\eps>0$ let $W_\eps$ be the open
set of all $g\in G$ with $d(gq,q)<\eps$. The group generated by $W_\eps$ is open, hence
closed, so it contains $G_0$, and it contains $G_q$. The sets of all $g$ with
$d(gq,q)\le\eps$ are compact, because the action is proper, and they decrease to $G_q$
as $\eps\to0$. Since $G_0G_q$ is an open set that contains $G_q$, there is $\eps_q>0$
such that $W_\eps\subset G_0G_q$ for $\eps\le\eps_q$. For these $\eps$ the group
generated by $W_\eps$ is $G_0G_q$. If $p,p'\in B_{\eps/4}(q)$ and $gp=p'$, then
\begin{equation*}
d(gq,q)\le d(gq,gp)+d(p',q)<\tfrac12\eps ,
\end{equation*}
so $g\in W_\eps\subset G_0G_q$. For the last assertion let $p,p'\in B_{\eps/4}(q)$. The
distance of their classes in $Z/G$ is $d(p,gp')$ for some $g\in G$, it is at most
$d(p,p')<\tfrac12\eps$, and
\begin{equation*}
d(gq,q)\le d(gq,gp')+d(gp',p)+d(p,q)<\eps .
\end{equation*}
So $g\in G_0G_q$, and the distance of the classes in $Z/(G_0G_q)$ is the same. Every
point of the ball of radius $\tfrac14\eps$ about the class of $q$, in either quotient,
is the class of a point of $B_{\eps/4}(q)$, because the minimum in the definition of the
quotient distance is attained. So the projection maps the one ball onto the other and
preserves distances.
\end{proof}

Part~(i) gives the Lie-group assertion of \ref{fl:S:limit}, and with it the proof of
\cref{m:fr:thm:geomain} is complete.

\subsection{Derived Tools}\label{fl:sec:tools}

We derive two tools for the later sections from \cref{m:fr:thm:geomain}: a common compact
realization of the convergence and Reifenberg neighbourhoods in the charts.

\subsubsection{A Common Realization of the Convergence}\label{fl:sec:realization}

The convergence of \cref{fl:S:Yconv}, that of \cref{m:fr:def:Yconv}, is defined through
frames of charts. To speak of
neighbourhoods of a point of $Y$ that contain parts of the spaces $Y_i$, we realize it as
the convergence in one compact space. We first note that it is determined by
subsequences.

\begin{lemma}[Subsequences]\label{fl:lem:subseq}
Let $\mathcal I$ be an infinite set of indices, let $y_i\in Y_i$ for $i\in \mathcal I$, and let $y\in Y$. If
every infinite subset of $\mathcal I$ contains an infinite subset along which $y_i$ converges to
$y$, then $y_i$ converges to $y$ along $\mathcal I$.
\end{lemma}
\begin{proof}
By \cref{fl:S:Yconv} the points $\pr(y_i)$ converge to $\pr(y)$ along such subsets, hence
along $\mathcal I$, because the convergence in $X$ is that of the points $\phi_i(\pr(y_i))$ in a
metric space, for the approximations $\phi_i\colon X_i\to X$. Choose a centre $a$ with
$\pr(y)\in B_{R_0/2}(a)$. For large $i\in \mathcal I$ the point $\pr(y_i)$ lies in $B_{R_0/2}(a_i)$,
so $y_i=\tau_i(p'_i)$ for a frame $p'_i\in\Fr_i$ by \cref{fl:S:charts}. Fix
$p\in\tau^{-1}(y)$. By hypothesis and \cref{fl:S:Yconv}, every infinite subset of $\mathcal I$
contains an infinite subset along which suitable frames in $\tau_i^{-1}(y_i)=\Gamma_ip'_i$ converge
to $p$. By \cref{m:fr:lem:orbitselect} there are frames $p_i\in\tau_i^{-1}(y_i)$ that converge to
$p$ along $\mathcal I$, and so $y_i\to y$ along $\mathcal I$.
\end{proof}

The realization is obtained from countably many test sequences, which separate the points of
$Y$ and define a compact metric on the disjoint union whose convergent sequences are those of
\ref{fl:S:Yconv}.

\begin{proposition}[Realization]\label{fl:prop:realization}
There is a metric $\sd$ on the disjoint union $\cY$ of $Y$ and of the spaces $Y_i$ with the
following properties.
\begin{enumerate}[label=\textup{(\roman*)}]
\item The space $(\cY,\sd)$ is compact, $\sd$ induces the given topologies of $Y$ and of
every $Y_i$, and every $Y_i$ is open and closed in $\cY$.
\item Let $I$ be an infinite set of indices, $y_i\in Y_i$ for $i\in I$ and $y\in Y$. Then
$\sd(y_i,y)\to0$ along $I$ if and only if $y_i$ converges to $y$ along $I$ in the sense of
\cref{fl:S:Yconv}.
\end{enumerate}
The topology of $\cY$ is determined by these two properties.
\end{proposition}
\begin{proof}
The space $C(Y)$ is separable. Let $g^1,g^2,\dots$ be a dense sequence in its unit ball, let
$(g^k_i)_i$ be test sequences with the limits $g^k$, which exist by \cref{fl:S:Yconv}, and for every $i$
let $h^1_i,h^2_i,\dots$ be continuous functions on $Y_i$ with values in $[0,1]$ that separate
its points. Define functions on $\cY$ by $H^0:=1/i$ on $Y_i$ and $H^0:=0$ on $Y$, by
$G^k:=g^k_i$ on $Y_i$ and $G^k:=g^k$ on $Y$, and by $H^k:=h^k_i/i$ on $Y_i$ and $H^k:=0$ on
$Y$, and put
\begin{equation*}
\sd(u,v):=|H^0(u)-H^0(v)|+\sum_{k\ge1}2^{-k}\min\bigl\{1,|G^k(u)-G^k(v)|\bigr\}
+\sum_{k\ge1}2^{-k}\,|H^k(u)-H^k(v)|.
\end{equation*}
This is a metric, because $H^0$ separates the different spaces, the functions $h^k_i$
separate the points of $Y_i$, and the functions $g^k$ those of $Y$. Its restrictions to $Y$
and to the spaces $Y_i$ are continuous, as uniformly convergent series of continuous
functions, so they induce the given topologies of these compact spaces; and $Y_i$ is open,
since the ball of radius $1/i-1/(i+1)$ about a point of $Y_i$ lies in $Y_i$.

For (ii) let $y_i\to y$ along $I$ in the sense of \cref{fl:S:Yconv}. Then $H^0(y_i)\to0$,
$H^k(y_i)\le1/i$ and $G^k(y_i)=g^k_i(y_i)\to g^k(y)$ for every $k$, so $\sd(y_i,y)\to0$ by
dominated convergence. Conversely, let $\sd(y_i,y)\to0$ along $I$. By \cref{fl:S:Yconv}
every infinite subset of $I$ contains an infinite subset along which $y_i$ converges to
some $y'\in Y$; then $\sd(y_i,y')\to0$ along it, so $y'=y$, and $y_i\to y$ along $I$ by
\cref{fl:lem:subseq}. A sequence in $\cY$ has infinitely many terms in $Y$ or in one of the
compact spaces $Y_i$, or a subsequence in spaces $Y_{i_k}$ with $i_k\to\infty$, which has a
convergent subsequence by \cref{fl:S:Yconv} and (ii); so $\cY$ is compact, and the sets $Y_i$
are closed. Finally, (i) and (ii) determine which sequences converge, hence the topology of
the metrizable space $\cY$. This is the usual realization of a convergent sequence of
compact spaces in one compact metric space, cf.\ \cite[Chapter~7]{BBI01}.
\end{proof}

By (i), the union of $Y$ and of the spaces $Y_i$ with $i\ge i_1$ is open in $\cY$ for
every $i_1$. An open neighbourhood of a point of $Y$ may therefore be replaced
by its intersection with such a union. The action of the
orthogonal group extends to the realization.

\begin{lemma}[Invariant Metrics]\label{fl:lem:invariant}
The actions of $O(n)$ on $Y$ and on the spaces $Y_i$ form a continuous action on $\cY$,
and there is a metric on $\cY$ that induces its topology and is invariant under this
action.
\end{lemma}
\begin{proof}
At the points of the open sets $Y_i$ the continuity is that of the action on $Y_i$. A
sequence in $\cY$ that converges to a point of $Y$ consists, up to finitely many terms in
each $Y_i$, of points of $Y$ and of points of spaces $Y_i$ with $i\to\infty$, and for
such sequences the continuity is stated in \cref{fl:S:charts,fl:S:Yconv} and in
\cref{fl:prop:realization}. For a metric $\sd_0$ that induces the topology, the maximum of $\sd_0(uA,vA)$ over
$A\in O(n)$ is an invariant metric, and it induces the topology because $O(n)$ is
compact.
\end{proof}

\subsubsection{Reifenberg Neighbourhoods}\label{fl:sec:reifenberg}

By \cref{m:fr:cor:localUIC}, the Reifenberg condition provides open sets homeomorphic to
$\R^n$ between comparable balls. We record this for the charts of the sequence. These sets
provide the contractible sets for the lifted homotopies of \cref{ch:fibregroups}.

\begin{lemma}[Reifenberg Sets]\label{fl:lem:reifenberg}
There is $\ell_0=\ell_0(n)\in(0,\tfrac12]$ with the following property. For every
$q\in\Omega'_i$ and every $0<s\le\tfrac14r$ there is an open set $U\subset Z_i$ that is
homeomorphic to $\R^n$, with
\begin{equation*}
B_{\ell_0s}(q)\subset U\subset B_s(q).
\end{equation*}
\end{lemma}
\begin{proof}
Put $R:=100s$ and $\ell_0:=100\,\ell_U$, with the constant $\ell_U(n)\in(0,\tfrac1{200}]$ of
\cref{m:fr:cor:localUIC}. Since $\varrho_i(q)<\tfrac12R_0$ and $R\le25r<\tfrac12R_0$, the
ball $B_R(q)$ lies in $\Omega_i$, and $R\le25r\le R_1$. So by \cref{fl:S:stage} the ball
$B_R(q)$ is $(\delta,R)$-Reifenberg, where $\delta\le\delta_*\le\delta_U(n)$,
and $Z_i$ is complete. \Cref{m:fr:cor:localUIC} gives an open set $U\subset Z_i$ that is
homeomorphic to $\R^n$, with $B_{\ell_UR}(q)\subset U\subset B_{R/100}(q)$, that is,
$B_{\ell_0s}(q)\subset U\subset B_s(q)$.
\end{proof}

\section{Harmonic Coordinates, Linear Slices and Realizing Maps}\label{ch:harmoniccomparison}\label{ch:slices}

We construct and control the limit coordinates used in
\cref{ch:stages,ch:localbundle,ch:globalbundle}. The local hypothesis is
\ref{ls:R}, stated in \cref{ls:sec:harmonic}: a Reifenberg region in an
$\RCD(\kappa,N)$ space. Harmonic charts are injective there, with comparison
errors that have small Lipschitz constants relative to the charts. This control
preserves injectivity under later corrections; uniform smallness alone would not.

\Cref{ls:sec:harmonic,th:sec:lipschitz} establish these analytic estimates, and
\cref{ls:sec:equivariant} uses canonical frames to make the charts isotropy-equivariant. The slice theorem,
stated in \cref{ls:sec:intro}, follows from differentiability along the group
direction (\cref{ls:sec:smoothing}) and the slice criterion
(\cref{ls:sec:slicesfromcharts}). \Cref{ls:sec:limits,fl:sec:anr} apply
the slices to the framed limit of \cref{ch:coveratlas}, giving its manifold
structure and equivariant realizing maps. Finally, quantitative slice estimates
(\cref{th:sec:limitchart}) and recentring (\cref{th:sec:along}) prepare the construction on the charts $Z_i$ in \cref{ch:stages}. Group hypotheses and
accuracy-dependent constants are specified at each step.

\subsection{Harmonic Charts Are Injective}\label{ls:sec:harmonic}

To prove injectivity with a quantitative H\"older bound for the inverse, we compare
gradients with those of a harmonic comparison map and apply
the local embedding theorem of Honda and Peng \cite[Proposition~4.7]{HP23}.
We use the following standing hypothesis, for a number $\eps>0$.

\begin{enumerate}[label=\textup{(R$_\eps$)},ref=\textup{(R$_\eps$)},leftmargin=3.2em]
\item\label{ls:R}\label{th:R} $(Z,d,m)$ is an $\RCD(\kappa,N)$ space, $\lambda>0$, and $\Omega\subset Z$ is an open
set that is $(\eps,\lambda)$-Reifenberg (\cref{m:fr:def:reifenberg}).
\end{enumerate}

\begin{definition}[Harmonic Charts]\label{ls:def:hchart}
Assume \ref{ls:R}, and let $z\in\Omega$, $s>0$ and $\theta>0$. A \emph{harmonic chart with
parameter $\theta$ at $(z,s)$} is a harmonic map $u\colon B_{8s}(z)\to\R^n$ with $u(z)=0$ for
which there is a map $\Phi\colon B_{8s}(z)\to\R^n$ that changes distances by at most
$\theta s$ and satisfies $|u-\Phi|\le\theta s$ on $B_{8s}(z)$. For $c>0$, the scale $s$ is
\emph{$c$-admissible at $z$} if $B_{9s}(z)\subset\Omega$, $s\le c\lambda$ and $|\kappa|s^2\le c$.
\end{definition}

A harmonic chart with parameter $\theta$ changes distances by at most $3\theta s$. Its
restriction to $B_{8s'}(z)$, for $s'<s$, is a harmonic chart at $(z,s')$ with the
parameter $\theta s/s'$, and $s'$ is $c$-admissible if $s$ is. Every sufficiently small
scale is $c$-admissible at a given point of $\Omega$.

The first step is the existence of a harmonic chart with a small parameter, which we take
from the quoted results.

\begin{lemma}[Comparison Maps]\label{ls:lem:existence}
For every $\delta\in(0,1)$ there are $\eps_X=\eps_X(n,N,\delta)>0$ and
$c_X=c_X(n,N,\delta)>0$ with the following property. Assume \ref{ls:R} with $\eps\le\eps_X$,
and let $s$ be $c_X$-admissible at $z\in\Omega$. Then $Z$ has the essential dimension
$n$, the set $\Omega$ is a topological manifold of dimension $n$, and there is a harmonic
map $v\colon B_{48ns}(z)\to\R^n$ with $v(z)=0$ whose components are Lipschitz with the
constant $C_N$ of \ref{m:fr:R:splitting}, such that
\begin{equation*}
\fint_{B_{48ns}(z)}|\Gr(v)-\one|\,dm<\delta
\end{equation*}
and such that the restriction of $v$ to $B_{8s}(z)$ is a centred $8\Psi_G(\delta|n,N,C_N)\,s$-GH approximation of $B_{8s}(z)$. In particular this restriction is a harmonic chart with parameter
$8\Psi_G(\delta|n,N,C_N)$ at $(z,s)$.
\end{lemma}
\begin{proof}
Let $\eps_B:=\eps_B(N,\delta)$ be the constant of \ref{m:fr:R:splitting}, and let
$\eps_X$ be the smaller of the numbers $\eps_V(n,N,\eps_B)$ and $\eps_D(n,N)$ of
\cref{ls:lem:measured,ls:cor:region}. Let $c_X>0$ be so small that $10\,n\,c_X\le\eps_X$ and
that $|\kappa|s^2\le c_X$ implies $|\kappa|(48ns)^2\le\min\{\eps_X,\eps_B,\delta\}$. The first two
assertions hold by \cref{ls:cor:region}. Put $s':=\tfrac{48}5ns$. Then $s'\le\eps_X\lambda$, so
by \cref{ls:lem:measured} the space $(Z,(s')^{-1}d,m,z)$, with the measure normalized at $z$, is $\eps_B$-close for $\mathbb D$ to $\R^n$ with a multiple of Lebesgue measure, and
\ref{m:fr:R:splitting} gives a harmonic map $v$ on the ball of radius $5s'=48ns$ about $z$,
with the Lipschitz bound and the bound for the Gram matrix; we subtract $v(z)$. By \ref{m:fr:R:almostisometry}, with
$8s$ in place of $s$ and $L=C_N$, the map $v$ changes the distances of the points of
$B_{8s}(z)$ by at most $\Psi_G\cdot8s$, and its image is dense up to this amount in
$B_{8s}(0^n)$.
\end{proof}

The H\"older exponent in the next proposition is free, and the constants depend on it.

\begin{proposition}[Harmonic Charts Are Injective]\label{ls:prop:injective}
For every $\eta\in(0,\tfrac14]$ there are positive numbers $\eps_A$, $\theta_A$ and $c_A$
that depend on $n$, $N$ and $\eta$, and there are $C_1=C_1(n,N)$ and $c_2=c_2(N)>0$, with
the following property. Assume \ref{ls:R} with $\eps\le\eps_A$, let $s$ be
$c_A$-admissible at $z\in\Omega$, and let $u$ be a harmonic chart with parameter
$\theta\le\theta_A$ at $(z,s)$.
\begin{enumerate}[label=\textup{(\roman*)}]
\item We have $|\nabla u_a|\le C_1$ almost everywhere on $B_{4s}(z)$ for all $a$, and $u$
is Lipschitz on $B_{2s}(z)$ with the constant $\sqrt n\,C_1$.
\item We have
\begin{equation*}
|u(p)-u(q)|\ge c_2\,\Bigl(\frac{d(p,q)}s\Bigr)^{\eta}\,d(p,q)\qquad
\text{for all }p,q\in B_{s/(100n)}(z),
\end{equation*}
and $u$ maps $B_{s/(100n)}(z)$ homeomorphically onto an open subset of $\R^n$.
\end{enumerate}
\end{proposition}
\begin{proof}
\emph{Constants.} Let $C_D$, $C_J$ and $C_N$ be the constants of
\ref{m:fr:R:geometry}, \ref{m:fr:R:estimates} and \ref{m:fr:R:splitting}, and put
\begin{gather*}
C_1:=\sqrt n\,C_N+4C_J,\qquad C_2:=8\sqrt n\,C_NC_J+16C_J^2,\\
c_2:=(4\sqrt{C_D})^{-1}.
\end{gather*}
Let
$\Psi_H(\eps|N)$ be the error in \cite[Proposition~4.7]{HP23}, with the
base-point normalization displayed in \cite[Corollary~4.8]{HP23}.
Choose $\eps_H\in(0,1)$ so small that
$\psi:=\Psi_H(\eps_H|N)\le\min\{\eta,\tfrac12\}$, and then choose
$\tilde\eta\in(0,\tfrac16]$ and $\delta\in(0,1)$ with
\begin{equation*}
C_2\tilde\eta+n\,C_D^{\,j}\delta<\eps_H\qquad\text{and}\qquad
8\,\Psi_G(\delta|n,N,C_N)\le\theta_A:=\eps_1(n,\tilde\eta),
\end{equation*}
where $j$ is the integer with $49n\le2^j<98n$ and $\eps_1$ is the function of
\cref{ls:lem:compare}. Let
$\eps_A:=\min\{\eps_X(n,N,\delta),\eps_V(n,N,\eps_H)\}$, with the constants of
\cref{ls:lem:existence,ls:lem:measured}. Finally let $c_A\in(0,c_X(n,N,\delta)]$ be so small that $10\,n\,c_A\le\eps_A$, and that
$|\kappa|(49ns)^2\le\min\{\eps_H,\eps_A\}$ whenever $|\kappa|s^2\le c_A$. The doubling property and the estimates of
\ref{m:fr:R:estimates} then hold on the balls used below. By \cref{ls:lem:existence} the
essential dimension is $n$ and $\Omega$ is a topological manifold of dimension $n$.

\emph{Step 1: a comparison map.} Let $v$ be the map of \cref{ls:lem:existence} for the
number $\delta$. Its restriction to $B_{8s}(z)$ is a centred $\eps'\cdot8s$-GH approximation of $B_{8s}(z)$ with $\eps':=\Psi_G(\delta|n,N,C_N)$, and ${\eps'\le\tfrac18\theta_A\le\tfrac12\eps_1(n,\tilde\eta)}$.

\emph{Step 2: gradients.} The map $u$ changes the distances of the points of $B_{8s}(z)$
by at most $3\theta s\le\tfrac12\eps_1(n,\tilde\eta)\cdot8s$. By
\cref{ls:lem:compare}\textup{(i)} there is a rigid motion $E$ such that the harmonic map
$w:=u-E\circ v$ satisfies $|w|\le16\tilde\eta s$ on $B_{8s}(z)$. By
\ref{m:fr:R:estimates}, $|\nabla w_a|\le4C_J\tilde\eta$ on $B_{4s}(z)$. The map $\tilde v:=R_Ev$
has $|\nabla\tilde v_a|\le\sqrt n\,C_N$ and $|\Gr(\tilde v)-\one|\le n\,|\Gr(v)-\one|$, and $u$
differs from $\tilde v+w$ by a constant. So $|\nabla u_a|\le C_1$ and
$|\Gr(u)-\Gr(\tilde v)|\le C_2\tilde\eta$ on $B_{4s}(z)$, which gives (i). Doubling gives
\begin{equation}\label{ls:eq:gramstart}
\fint_{B_{4s}(z)}|\Gr(u)-\one|\,dm
\le C_2\tilde\eta+n\,\frac{m(B_{48ns}(z))}{m(B_{4s}(z))}\,\delta
\le C_2\tilde\eta+n\,C_D^{\,j}\delta<\eps_H .
\end{equation}

\emph{Step 3: the H\"older estimate.} Rescale to
$(Z,\bar d,\bar m):=(Z,s^{-1}d,m(B_s(z))^{-1}m)$ and $F:=u/s$ on $B_4^{\bar d}(z)$.
This space is $\RCD(\kappa s^2,N)$, hence $\RCD(-\eps_H,N)$. For every
$y\in B_{3s}(z)$ and $0<t\le1$, put $s'=ts$ and
$\eps'_V:=\eps_V(n,N,\eps_H)$. We have $y\in\Omega$, $s'/\eps'_V\le\lambda$ and $|\kappa|(s')^2\le\eps'_V$. Thus \ref{ls:R} and \cref{ls:lem:measured} give
\begin{equation*}
\mathbb D\bigl((Z,(s')^{-1}d,m(B_{s'}(y))^{-1}m,y),
(\R^n,|\cdot|,\omega_n^{-1}\cH^n,0)\bigr)<\eps_H ,
\end{equation*}
with $\mathbb D$ the distance of \cite{HP23}, as allowed in \cref{m:fr:sec:rcd}. The
gradients and the Gram matrix of $F$ for $\bar d$ are those of $u$ for $d$,
$\bar\Delta F_a=s\Delta u_a=0$, and $F$ is bounded on $B_4^{\bar d}(z)$, since $u(z)=0$ and
$u$ changes distances by at most $3\theta s$. For every smooth $\Theta$ on $\R^n$, the local
chain rule of \ref{m:fr:R:calculus}, the gradient bound in (i) and the boundedness of $F$
give
\begin{equation*}
\bar\Delta(\Theta\circ F)
=\sum_{a,b}\partial_{ab}\Theta(F)\,\Gr(u)_{ab}\in L^\infty(B_{4s}(z),m).
\end{equation*}
So $F$ is regular in the sense of \cite[Definition~3.1]{HP23}, and by
\eqref{ls:eq:gramstart} and the closeness above the hypotheses of
\cite[Proposition~4.7]{HP23}, in the normalization of \cite[Corollary~4.8]{HP23}, hold on
$B_4^{\bar d}(z)$. That proposition yields
\begin{equation*}
|u(p)-u(q)|\ge(1-\psi)\Bigl(\frac{d(p,q)}s\Bigr)^\psi d(p,q)
\qquad(p,q\in B_s(z)).
\end{equation*}
For distinct $p,q\in B_{s/(100n)}(z)$ we have $d(p,q)/s<1$, so $\psi\le\eta$ and
$1-\psi\ge\tfrac12\ge c_2$ give (ii). The case $p=q$ is immediate.
The bound makes $u$ injective on this ball, an open subset of the topological
$n$-manifold $\Omega$. Invariance of domain gives the asserted homeomorphism
onto an open subset of $\R^n$.
\end{proof}

\subsection{Harmonic Functions in Harmonic Charts}\label{th:sec:lipschitz}

We show that harmonic functions are Lipschitz in harmonic charts. This gives the
comparison of two charts, or a chart and its composition with a measure-preserving
isometry, by a rigid motion and an error with small relative Lipschitz constant.

\subsubsection{Notation and Constants}\label{th:sec:constants}

We fix the notation and the constants used in the Lipschitz comparisons of this and the later
sections. We use hypothesis \ref{ls:R} and \cref{ls:def:hchart}. Rigid motions and their rotation parts are as in \cref{m:fr:def:rigid}. For maps $\mathsf e$ and $u$ on a
set $A$ and a number $L\ge0$ we say that \emph{$\mathsf e$ is $L$-Lipschitz with respect
to $u$ on $A$} if
\begin{equation*}
|\mathsf e(p)-\mathsf e(q)|\le L\,|u(p)-u(q)|\qquad\text{for all }p,q\in A .
\end{equation*}

Let $C_D=C_D(N)\ge2$ be the doubling constant of \textup{\ref{m:fr:R:geometry}}, and put
\begin{equation*}
\eta_0=\eta_0(n):=\min\bigl\{\tfrac14,\tfrac12\eps_1(n,\tfrac1{100})\bigr\},\qquad
\eta:=\min\bigl\{\eta_0(n),\,(4+16\log_2C_D)^{-1}\bigr\}.
\end{equation*}
Here $\eps_1$ is the comparison threshold of \cref{ls:lem:compare}.
The assertions of \cref{ls:prop:injective} hold for every exponent in $(0,\tfrac14]$,
with constants that depend on the exponent. We write $\eps_A$, $\theta_A$ and $c_A$ for the
smaller of the constants of that proposition for the two exponents $\eta_0(n)$ and
$\eta$; they depend on $n$ and $N$, and so do the constants $C_1$ and $c_2$ there. So if
\ref{th:R} holds with $\eps\le\eps_A$, the scale $s$ is $c_A$-admissible at $z$, and $u$ is a
harmonic chart with parameter at most $\theta_A$ at $(z,s)$, then $u$ maps
$B_{s/(100n)}(z)$ homeomorphically onto an open subset of $\R^n$, and
\begin{equation}\label{th:eq:holder}
c_2\,\Bigl(\frac{d(p,q)}s\Bigr)^{\eta}\,d(p,q)\le|u(p)-u(q)|\le\sqrt n\,C_1\,d(p,q)
\qquad\text{for all }p,q\in B_{s/(100n)}(z).
\end{equation}
The estimates below allow any fixed $\eps_*\in(0,\tfrac1{200}]$.
For the application to \cref{thm:main} we choose once and for all
\begin{equation}\label{es:eq:kappasharp}
M_n:=539^n,\qquad
\eps_*:=\frac1{240\,((n+1)M_n+1)}.
\end{equation}
Here $(n+1)M_n$ is the overlap bound proved in \cref{es:lem:system};
the choice leaves a total Lipschitz-error budget of $1/4$ for gluing.
With this choice the subsequent universal accuracy thresholds depend only on $n,N$;
the orbit-dependent radii remain local (\cref{guide:sec:choices}).
We put
\begin{equation*}
\tilde\eta_\flat:=\frac{\eps_*\theta_A}{256},\qquad
\theta_\flat:=\min\bigl\{\tfrac14\theta_A,\tfrac1{12}\eps_1(n,\tilde\eta_\flat)\bigr\},
\end{equation*}
with the function $\eps_1$ of \cref{ls:lem:compare}. With the function $\Psi_G$ and the constant $C_N$ of \textup{\ref{m:fr:R:almostisometry}},
\textup{\ref{m:fr:R:splitting}} we choose a number $\delta_\flat\in(0,1)$ such that
\begin{equation*}
\Psi_G(\delta_\flat|n,N,C_N)\le\min\bigl\{\tfrac12\eps_1(n,\tilde\eta_\flat),\tfrac1{24}\theta_\flat\bigr\},
\end{equation*}
and we let
\begin{equation*}
\eps_\flat:=\min\{\eps_A,\eps_X(n,N,\delta_\flat)\},\qquad
c_\flat:=\min\{c_A,c_X(n,N,\delta_\flat)\},
\end{equation*}
with the constants of \cref{ls:lem:existence}.

\begin{definition}[The Comparison Map]\label{th:def:comparison}
Let \ref{th:R} hold with $\eps\le\eps_\flat$, and let $s$ be $c_\flat$-admissible at $z$. A
\emph{comparison map} at $(z,s)$ is a harmonic map $v\colon B_{48ns}(z)\to\R^n$ with $v(z)=0$ and
with $C_N$-Lipschitz components whose restriction to $B_{8s}(z)$ is a centred
$4\eps_1(n,\tilde\eta_\flat)\,s$-GH approximation of $B_{8s}(z)$ and a harmonic chart with parameter
$\tfrac13\theta_\flat$ at $(z,s)$.
\end{definition}

Comparison maps exist by \cref{ls:lem:existence}.

\subsubsection{The Lipschitz Estimate}\label{th:sec:lipproof}

The estimate comes from the injectivity of harmonic charts, applied to a chart that is
perturbed by the given function. It turns the uniform comparison of \cref{ls:lem:compare} into
the Lipschitz comparison that \cref{ch:stages,ch:globalbundle} require.

\begin{lemma}[Harmonic Functions Are Lipschitz in Harmonic Charts]\label{th:lem:lipschitz}
Assume \ref{th:R} with $\eps\le\eps_A$, let $s$ be $c_A$-admissible at $z\in\Omega$, and let
$u$ be a harmonic chart with parameter at most $\tfrac12\theta_A$ at $(z,s)$. Let $j\ge1$
and $\alpha>0$, and let $\mathsf e\colon B_{8s}(z)\to\R^j$ be a harmonic map with
$|\mathsf e-\mathsf e(z)|\le\alpha s$ on $B_{8s}(z)$. Then $\mathsf e$ is $(2\alpha/\theta_A)$-Lipschitz with
respect to $u$ on $B_{s/(100n)}(z)$.
\end{lemma}
\begin{proof}
Let $p,q\in B_{s/(100n)}(z)$ with $\mathsf e(p)\ne \mathsf e(q)$. Then $p\ne q$, so $u(p)\ne u(q)$ by
\eqref{th:eq:holder}. Let $\nu$ and $\omega$ be the unit vectors in the directions of
$\mathsf e(p)-\mathsf e(q)$ and of $u(p)-u(q)$, and put
\begin{equation*}
f:=\langle \mathsf e-\mathsf e(z),\nu\rangle,\qquad\sigma:=\frac{|u(p)-u(q)|}{|\mathsf e(p)-\mathsf e(q)|},\qquad
\hat u:=u-\sigma f\,\omega .
\end{equation*}
The function $f$ is harmonic with $|f|\le\alpha s$, and $f(p)-f(q)=|\mathsf e(p)-\mathsf e(q)|$. The map
$\hat u$ is harmonic on $B_{8s}(z)$, it vanishes at $z$, and
\begin{equation*}
\hat u(p)-\hat u(q)=|u(p)-u(q)|\,\omega-\sigma\,|\mathsf e(p)-\mathsf e(q)|\,\omega=0 .
\end{equation*}
Let $\Phi$ be a map on $B_{8s}(z)$ that changes distances by at most $\tfrac12\theta_As$
with $|u-\Phi|\le\tfrac12\theta_As$. Then $|\hat u-\Phi|\le(\tfrac12\theta_A+\sigma\alpha)s$. If
$\sigma\alpha\le\tfrac12\theta_A$, the map $\hat u$ is a harmonic chart with parameter at most
$\theta_A$ at $(z,s)$, hence injective on $B_{s/(100n)}(z)$, which contradicts
$\hat u(p)=\hat u(q)$. So $\sigma\alpha>\tfrac12\theta_A$, that is,
$|\mathsf e(p)-\mathsf e(q)|<(2\alpha/\theta_A)\,|u(p)-u(q)|$.
\end{proof}

\subsubsection{Two Harmonic Charts, and Isometries}\label{th:sec:twocharts}

We compare a harmonic chart with a second harmonic map that almost preserves distances.
The two differ by a rigid motion up to a uniformly small error, by \cref{ls:lem:compare}, and by \cref{th:lem:lipschitz} the error has a small Lipschitz constant with respect to the harmonic chart.

\begin{lemma}[Two Harmonic Charts]\label{th:lem:twocharts}
Assume \ref{th:R} with $\eps\le\eps_\flat$, and let $s$ be $c_\flat$-admissible at
$z\in\Omega$. Let $u$ be a harmonic chart with parameter at most $2\theta_\flat$ at $(z,s)$,
and let $u'\colon B_{8s}(z)\to\R^n$ be a harmonic map that changes distances by at most
$12\theta_\flat s$. Then there is a rigid motion $E$ such that $\mathsf e:=u'-E\circ u$ satisfies
$|\mathsf e|\le32\,\tilde\eta_\flat s$ on $B_{8s}(z)$, and $\mathsf e$ is $\tfrac12\eps_*$-Lipschitz with
respect to $u$ on $B_{s/(100n)}(z)$.
\end{lemma}
\begin{proof}
Let $v$ be a comparison map at $(z,s)$. Its restriction to $B_{8s}(z)$ is a centred $\eps'\cdot8s$-GH approximation of $B_{8s}(z)$ for $\eps':=\tfrac12\eps_1(n,\tilde\eta_\flat)$. The maps $u$ and
$u'$ change the distances of the points of $B_{8s}(z)$ by at most $12\theta_\flat s$, which
is at most $\eps'\cdot8s$. By \cref{ls:lem:compare}\textup{(i)}, applied twice at $(z,8s)$,
there are rigid motions $E_1$ and $E_2$ with $|u-E_1\circ v|\le16\,\tilde\eta_\flat s$ and
$|u'-E_2\circ v|\le16\,\tilde\eta_\flat s$ on $B_{8s}(z)$. Put $E:=E_2\circ E_1^{-1}$. Since
$E$ is an isometry of $\R^n$,
\begin{equation*}
|u'-E\circ u|\le|u'-E_2\circ v|+|E_1\circ v-u|\le32\,\tilde\eta_\flat s .
\end{equation*}
The map $\mathsf e=u'-R_Eu-b_E$ is harmonic, and $|\mathsf e-\mathsf e(z)|\le64\,\tilde\eta_\flat s$. Since
$2\theta_\flat\le\tfrac12\theta_A$, \cref{th:lem:lipschitz} applies with $\alpha=64\,\tilde\eta_\flat$,
and $2\alpha/\theta_A=\tfrac12\eps_*$.
\end{proof}

Applied to the composition of the chart with an isometry, which is again harmonic and almost
preserves distances, the lemma gives the form of the isometries in a harmonic chart.

\begin{corollary}[Isometries in a Harmonic Chart]\label{th:cor:isometries}
Assume \ref{th:R} with $\eps\le\eps_\flat$, and let $2s$ be $c_\flat$-admissible at
$z\in\Omega$. Let $u$ be a harmonic chart with parameter at most $\theta_\flat$ at $(z,2s)$,
and let $\varphi$ be an isometry of $Z$ that preserves $m$, with $d(\varphi(z),z)\le8s$.
Then there is a rigid motion $E_\varphi$ such that $\mathsf e_\varphi:=u\circ\varphi-E_\varphi\circ u$,
which is defined on $B_{8s}(z)$, satisfies $|\mathsf e_\varphi|\le32\,\tilde\eta_\flat s$ there and is
$\tfrac12\eps_*$-Lipschitz with respect to $u$ on $B_{s/(100n)}(z)$.
\end{corollary}
\begin{proof}
The isometry maps $B_{8s}(z)$ into $B_{16s}(z)$, so $u\circ\varphi$ is defined on
$B_{8s}(z)$. It is harmonic by \textup{\ref{m:fr:R:calculus}}, and it changes distances by at
most $3\theta_\flat\cdot2s$, as $u$ does. The restriction of $u$ to $B_{8s}(z)$ is a harmonic
chart with parameter $2\theta_\flat$ at $(z,s)$, and $s$ is $c_\flat$-admissible. So
\cref{th:lem:twocharts} applies to $u$ and $u':=u\circ\varphi$.
\end{proof}

\subsection{Equivariant Harmonic Charts}\label{ls:sec:equivariant}

At a point fixed by a compact group, the canonical frame space already carries an exact
orthogonal representation. We align a harmonic comparison map with a frame map and average
the harmonic map against this representation. Thus exact equivariance comes from the frames,
while harmonicity comes from the comparison map.

\begin{proposition}[Equivariant Harmonic Charts]\label{ls:prop:equivariant}
For every $\theta\in(0,1)$ there are positive numbers $\eps_K$ and $c_K$ that depend on
$n$, $N$ and $\theta$, with the following property. Assume \ref{ls:R} with $\eps\le\eps_K$, let
$z\in\Omega$, let $K$ be a compact group that acts continuously on $Z$ by isometries that
preserve $m$ and fix $z$, and let $s$ be $c_K$-admissible at $z$. Then there are a continuous
homomorphism $\mathsf R\colon K\to O(n)$ and a harmonic chart $u$ with parameter $\theta$ at
$(z,s)$ such that $u(kp)=\mathsf R(k)\,u(p)$ for all $k\in K$ and $p\in B_{8s}(z)$.
\end{proposition}
\begin{proof}
\emph{Reduction.} Since the action is continuous, it defines a continuous homomorphism of $K$
onto a subgroup $\bar K$ of the isometry group of $Z$ with the compact-open topology. This
group is metrizable and acts continuously on $Z$, because $Z$ is proper. Hence $\bar K$ is a
compact metrizable group that acts continuously on $Z$ by isometries that preserve $m$ and fix
$z$. If $\bar{\mathsf R}$ and $u$ have the asserted properties for $\bar K$, then $u$ and the
composition of $\bar{\mathsf R}$ with the homomorphism $K\to\bar K$ have them for $K$. So we may
assume that $K$ is metrizable, as \cref{ls:lem:integrals} requires.

\emph{The auxiliary frame.} Put $a:=\theta/128$ and
$\hat\eps:=\min\{\tfrac12\eps_1(n,a),\theta/8\}$. Choose $\eps_f>0$ so small that
\cref{m:fr:lem:flat,m:fr:thm:A} and \cref{m:fr:lem:framegauge}\textup{(i)} give frames and frame-map
distance error at most $\hat\eps\widehat r$ at an $\eps_f$-flat point of scale $\widehat r$.
Here we use the auxiliary scale $\widehat r:=4s$, independently of the fixed frame scale
in \cref{ch:coveratlas}. Put $\eps'_V:=\eps_V(n,N,\eps_f)$ from \cref{ls:lem:measured},
choose $\delta\in(0,1)$ with $\Psi_G(\delta|n,N,C_N)\le\hat\eps$, and take
\[
\eps_K:=\min\{\eps'_V,\eps_X(n,N,\delta)\},\qquad
c_K:=\min\{\eps'_V/16,c_X(n,N,\delta)\}.
\]
Then $\widehat r/\eps'_V\le\lambda$ and $|\kappa|\widehat r^2\le\eps'_V$, so
\cref{ls:lem:measured} makes $z$ $\eps_f$-flat at scale $\widehat r$. Choose a frame $e$ there
and write $u_e:=u_e^{(\widehat r)}$.

\emph{The exact representation.} Since $K$ fixes $z$ and preserves $m$, it preserves
$\mu_z^{(\widehat r)}$ and commutes with $\cK_z^{(\widehat r)}$. By
\cref{m:fr:thm:bundle}, the identities
\[
e\circ k^{-1}=e\cdot\mathsf R(k),\qquad
u_e(kp)=\mathsf R(k)u_e(p)
\]
define a homomorphism $\mathsf R\colon K\to O(n)$. It is continuous: if $k_j\to k$, the maps
$u_e\circ k_j^{-1}$ converge uniformly on bounded sets to $u_e\circ k^{-1}$, which is
convergence in the frame topology.

\emph{The harmonic average.} By \cref{ls:lem:existence} there is a centred harmonic
$\hat\eps\cdot8s$-GH approximation $v$ on $B:=B_{8s}(z)$, with $C_N$-Lipschitz components.
On $B\subset B_{3\widehat r}(z)$ the map $u_e$ changes distances by at most
$\hat\eps\widehat r\le\hat\eps\cdot8s$. Thus \cref{ls:lem:compare}\textup{(i)}, with $v$
as the approximation, gives a rigid motion $E$ with $|u_e-E\circ v|\le16as$ on $B$. Frame maps vanish at their base point, because $\mathsf G_z(z,\cdot)=0$, and $v(z)=0$; so $|b_E|\le16as$, and $|u_e-w|\le32as$ for the centred harmonic map $w:=R_Ev$. Define
\[
u(p):=\int_K\mathsf R(k)^{-1}w(kp)\,dk\qquad(p\in B),
\]
with Haar measure of mass one. The ball $B$ is $K$-invariant, and the integrands are
bounded continuous harmonic maps with uniformly Lipschitz components. Hence $u$ is
harmonic by \cref{ls:lem:integrals}; Haar invariance gives exact equivariance, and $u(z)=0$.
Since $u_e$ is equivariant, $|u-u_e|\le32as$, whence $|u-w|\le64as=\theta s/2$.
The map $w$ changes distances by at most $8\hat\eps s\le\theta s$, so it is the comparison
map required to make $u$ a harmonic chart with parameter $\theta$.
\end{proof}

\subsection{Linear Slices}\label{ls:sec:intro}\label{ls:sec:results}

We use equivariant harmonic charts to construct slices for proper isometric Lie group
actions. Actions are continuous by convention. For a point $z$ with isotropy group $K$, the group $K$ acts on
the quotient $\fg/\fk$ of the Lie algebras by the adjoint representation. For a space
$S$ with an action of $K$, the twisted product $G\times_KS$ is the quotient of
$G\times S$ by the action $k\cdot(g,q):=(gk^{-1},kq)$, and $[g,q]$ denotes the class of
$(g,q)$; the group $G$ acts on it by left multiplication on the first factor.

\begin{theorem}[Linear Slices]\label{ls:thm:slices}
For all integers $n\ge1$ and all $N<\infty$ there is $\eps_L=\eps_L(n,N)>0$ with the
following property. Let $(Z,d,m)$ be an $\RCD(\kappa,N)$ space, let $G$ be a Lie group that
acts on $Z$ properly and by isometries that preserve $m$, and let $\Omega\subset Z$ be an
open invariant set that is $(\eps_L,\lambda)$-Reifenberg for some $\lambda>0$
(\cref{m:fr:def:reifenberg}). Let $z\in\Omega$, and let $K$ be the isotropy group of $z$. Then there are a continuous
homomorphism $\mathsf R\colon K\to O(n)$, a decomposition $\R^n=T\oplus V$ into orthogonal
subspaces that are invariant under $\mathsf R(K)$, an open neighbourhood $U$ of $z$ in
$\Omega$ that is invariant under $K$, a map $h\colon U\to\R^n$ with $h(z)=0$, and a
number $r_1>0$ with the following properties.
\begin{enumerate}[label=\textup{(\roman*)}]
\item The map $h$ is a homeomorphism onto an open subset of $\R^n$ that contains the
closed ball of radius $r_1$ about $0$, and $h(kp)=\mathsf R(k)\,h(p)$ for all $k\in K$ and
$p\in U$.
\item The representation of $K$ on $T$ is isomorphic to its adjoint representation on
$\fg/\fk$. In particular $\dim V=n-\dim G+\dim K$.
\item Put $S:=h^{-1}(V\cap B_{r_1}(0^n))$. The map $[g,q]\mapsto gq$ is a homeomorphism of
$G\times_KS$ onto an open neighbourhood of the orbit $Gz$ in $Z$.
\end{enumerate}
\end{theorem}

The harmonic chart $h$
is Lipschitz and bi-H\"older, by \cref{ls:prop:injective,ls:prop:smooth}. We name the
notions that the theorem provides. Tubes and slices are those of Bredon \cite[Chapter~II]{Bre72},
and locally linear actions are his locally smooth actions \cite[Chapter~IV, Section~1]{Bre72}.

\begin{definition}[Tubes and Slices]\label{ls:def:slice}
Let a Lie group $G$ act on a space $Z$. A \emph{tube} in $Z$ is a $G$-equivariant homeomorphism
$\Xi$ of a twisted product $G\times_IB$, for a closed subgroup $I$ of $G$ and a space $B$ with an
action of $I$, onto an open subset of $Z$; it is a tube about the orbits in its image, and it is
\emph{linear} if $B$ is a ball about $0$ in an orthogonal representation of $I$. Now let the
action be proper, and let $K$ be the isotropy group of a point $z$. A \emph{slice} at $z$ is a
$K$-invariant subset $S$ of $Z$ that contains $z$, with the subspace topology, such that
$[g,q]\mapsto gq$ is a homeomorphism of $G\times_KS$ onto an open neighbourhood of the orbit $Gz$,
that is, a tube about $Gz$. The slice is \emph{linear} if there is a $K$-equivariant
homeomorphism of $S$ onto a ball about $0$ in an orthogonal representation of $K$; the tube is
then linear after this identification.
\end{definition}

Thus $S$ in \cref{ls:thm:slices}(iii) is a linear slice: $h$ maps it
homeomorphically and $K$-equivariantly onto the ball $V\cap B_{r_1}(0^n)$,
where $K$ acts through $\mathsf R$. The corresponding linear tube has domain
$G\times_K(V\cap B_{r_1}(0^n))$.

\begin{definition}[Locally Linear Actions]\label{ls:def:locallylinear}
An action of a Lie group by homeomorphisms is \emph{locally linear} if every point has a
linear slice.
\end{definition}

The action of $G$ on $\Omega$ in \cref{ls:thm:slices} is therefore locally linear.

\subsection{Differentiability along the Group Direction}\label{ls:sec:smoothing}

Averaging the chart of \cref{ls:prop:equivariant} against a smooth density on $G$
makes it differentiable along the group direction, as in G{\aa}rding's construction of smooth
vectors \cite{Gar47}. We also prove that its orbit derivative is injective transverse
to the isotropy algebra. In this subsection and the next we assume, for $\eps>0$:

\begin{enumerate}[label=\textup{(G$_\eps$)},ref=\textup{(G$_\eps$)},leftmargin=3.2em]
\item\label{ls:G}\label{th:G} Hypothesis \ref{ls:R} holds, $G$ is a Lie group that acts on $Z$ properly and
by isometries that preserve $m$, the set $\Omega$ is invariant under $G$, and $z$ is a
point of $\Omega$ with the isotropy group $K$.
\end{enumerate}

The group $K$ is compact, because the action is proper, and it fixes $z$; so
\cref{ls:prop:equivariant} applies to it. For $0<\theta\le\theta_A$ we fix
the constants
\begin{equation*}
\eps_L(\theta):=\min\{\eps_A,\eps_K(n,N,\tfrac1{16}\theta)\},\qquad
c_L(\theta):=\min\{c_A,c_K(n,N,\tfrac1{16}\theta)\},
\end{equation*}
where $\eps_A$, $\theta_A$ and $c_A$ retain the common choices of
\cref{th:sec:constants}, which in particular apply for $\eta=\eta_0(n)$,
and $\eps_K$, $c_K$ are those of \cref{ls:prop:equivariant}; they depend on
$n$, $N$ and $\theta$. We write $\eps_L:=\eps_L(\theta_A)$ and $c_L:=c_L(\theta_A)$ for
the constants of the slice theorem; these depend on $n$ and $N$ only. We use the norm on $\fg$ and the complement $\fm$ of $\fk$ of
\cref{ls:sec:isotropy}, and the numbers $r_0$ and $c_0$ of \cref{ls:lem:tube}.

\begin{proposition}[Harmonic Charts That Are Differentiable along the Group]\label{ls:prop:smooth}\label{th:lem:chart}
Let $0<\theta\le\theta_A$. Assume \ref{ls:G} with $\eps\le\eps_L(\theta)$, let $s$ be
$c_L(\theta)$-admissible at $z$, and let $\mathsf R$ be the homomorphism of
\cref{ls:prop:equivariant} for the parameter $\theta/16$. Then there are
a harmonic map $h\colon B_{3s/2}(z)\to\R^n$ with $h(z)=0$, numbers $r_2>0$ and
$0<C_3<\infty$, and linear maps $L_p\colon\fg\to\R^n$ for $p\in B_s(z)$, with the following
properties.
\begin{enumerate}[label=\textup{(\roman*)}]
\item The restriction of $h$ to $B_s(z)$ is a harmonic chart with parameter $\theta$ at
$(z,\tfrac18s)$, and \mbox{$h(kp)=\mathsf R(k)\,h(p)$} for all $k\in K$ and $p\in B_{3s/2}(z)$.
\item For all $p\in B_s(z)$ and all $\zeta\in\fg$ with $|\zeta|<r_2$ the point $\exp(\zeta)p$ lies
in $B_{3s/2}(z)$, and
\begin{equation*}
\bigl|h(\exp(\zeta)p)-h(p)-L_p\zeta\bigr|\le C_3\,|\zeta|^2 .
\end{equation*}
\item For all $p,p'\in B_s(z)$ and $\zeta\in\fg$ we have
$|L_p\zeta-L_{p'}\zeta|\le C_3\,d(p,p')\,|\zeta|$.
\item We have $L_z(\Ad(k)\zeta)=\mathsf R(k)L_z\zeta$ for $k\in K$ and $\zeta\in\fg$; the map $L_z$
vanishes on $\fk$, and it is injective on $\fm$.
\end{enumerate}
\end{proposition}
\begin{proof}
\emph{The density.} Let $u$ be the harmonic chart of \cref{ls:prop:equivariant}
with parameter $\theta/16$ at $(z,s)$. By \cref{ls:prop:injective}\textup{(i)} it is Lipschitz
on $B_{2s}(z)$ with the constant $C_4:=\sqrt n\,C_1$. Put
$\varepsilon_{\mathrm{av}}:=\min\{\tfrac12,\theta/(32C_4)\}$. The function
$D(g):=\max\{d(gp,p):p\in\bar B_{2s}(z)\}$ is continuous on $G$, since the action is
continuous and the closed ball is compact, and $D(g^{-1})=D(g)=D(kgk^{-1})$ for $k\in K$.
So
\begin{equation*}
\mathcal V:=\{g\in G:D(g)<\tfrac12\varepsilon_{\mathrm{av}} s\}
\end{equation*}
is an open neighbourhood of the identity that is symmetric and invariant under conjugation
by $K$, and its closure is compact because the action is proper. For $g,g'\in\mathcal V$
and $p\in B_{3s/2}(z)$ the point $g'p$ lies in $B_{2s}(z)$, so $d(gg'p,p)<\varepsilon_{\mathrm{av}} s$. Let
$d_rg$ be a right Haar measure on $G$. Conjugation by an element of $K$ maps it to a
multiple of itself, and the factor is a continuous homomorphism of the compact group $K$
into the positive numbers, hence equal to one. We choose a smooth function $\varphi\ge0$
on $G$ with compact support in $\mathcal V$, with $\int\varphi\,d_rg=1$ and with
$\varphi(kgk^{-1})=\varphi(g)$ for $k\in K$; the average over $K$ of the conjugates of a
function without the last property has all of them.

\emph{The harmonic chart.} For $p\in B_{3s/2}(z)$ let
\begin{equation*}
h_0(p):=\int_G\varphi(g)\,u(gp)\,d_rg ,\qquad h(p):=h_0(p)-h_0(z).
\end{equation*}
The functions $(g,p)\mapsto u_a(gp)$ are continuous and bounded on the product of the
support of $\varphi$ with $B_{3s/2}(z)$, and for fixed $g$ they are harmonic, by
\ref{m:fr:R:calculus}, and Lipschitz with the constant $C_4$. So $h_0$ is harmonic by
\cref{ls:lem:integrals}, and $|h_0-u|\le\tfrac12C_4\varepsilon_{\mathrm{av}} s$. For $k\in K$ we have
$u(gkp)=\mathsf R(k)\,u(k^{-1}gk\,p)$, and the substitution $g\mapsto kgk^{-1}$ gives
$h_0(kp)=\mathsf R(k)h_0(p)$. In particular $\mathsf R(k)$ fixes $h_0(z)$, so $h$ has the same
property, and $|h_0(z)|\le\tfrac12C_4\varepsilon_{\mathrm{av}} s$ because $u(z)=0$. Let $\Phi$ be a map on
$B_{8s}(z)$ that changes distances by at most $\theta s/16$ with $|u-\Phi|\le\theta s/16$. On
$B_s(z)$ we get $|h-\Phi|\le(C_4\varepsilon_{\mathrm{av}}+\theta/16)s\le\tfrac18\theta s$, and
$\theta s/16\le\tfrac18\theta s$. This is (i).

\emph{Differentiation.} Choose $r_2>0$ such that $\exp(\zeta)\in\mathcal V$ for $|\zeta|\le r_2$.
Let $p\in B_s(z)$ and $|\zeta|<r_2$. Then $\exp(\zeta)p\in B_{3s/2}(z)$, and the invariance of
$d_rg$ under right translations gives
\begin{equation*}
h_0(\exp(\zeta)p)=\int_G\varphi(g)\,u(g\exp(\zeta)p)\,d_rg
=\int_G\varphi\bigl(g\exp(-\zeta)\bigr)\,u(gp)\,d_rg .
\end{equation*}
The function $(\zeta,g)\mapsto\varphi(g\exp(-\zeta))$ is smooth, and for $|\zeta|\le r_2$ it
vanishes unless $g$ lies in the compact closure $Q$ of the set of the products of two
elements of $\mathcal V$; for $g\in Q$ and $p\in B_s(z)$ the point $gp$ lies in
$B_{2s}(z)$. Let $\varphi_1(\zeta,g)$ be the derivative of $t\mapsto\varphi(g\exp(-t\zeta))$
at $t=0$, which is linear in $\zeta$. By the formula of Taylor there is a constant
$C_\varphi$ with
\begin{equation*}
\bigl|\varphi(g\exp(-\zeta))-\varphi(g)-\varphi_1(\zeta,g)\bigr|\le C_\varphi|\zeta|^2,\qquad
|\varphi_1(\zeta,g)|\le C_\varphi|\zeta|\qquad(g\in G,\ |\zeta|\le r_2),
\end{equation*}
and both functions vanish for $g\notin Q$. We define
\begin{equation*}
L_p\zeta:=\int_G\varphi_1(\zeta,g)\,u(gp)\,d_rg\qquad(p\in B_s(z),\ \zeta\in\fg).
\end{equation*}
Since $|u|\le2C_4s$ on $B_{2s}(z)$ and $u$ is Lipschitz there with the constant $C_4$,
assertions (ii) and (iii) hold with
$C_3:=1+C_\varphi\,C_4\max\{2s,1\}\int_Q d_rg$.

\emph{The derivative at $z$.} Let $k\in K$ and $\zeta\in\fg$, and let $t>0$ be small. Since
$k^{-1}$ fixes $z$, the point $\exp(t\Ad(k)\zeta)z=k\exp(t\zeta)k^{-1}z$ is the image of
$\exp(t\zeta)z$ under $k$, so
\begin{equation*}
h\bigl(\exp(t\Ad(k)\zeta)z\bigr)=\mathsf R(k)\,h\bigl(\exp(t\zeta)z\bigr).
\end{equation*}
By (ii) at $p=z$ the two sides are
$tL_z(\Ad(k)\zeta)$ and $t\mathsf R(k)L_z\zeta$ up to errors of the order $t^2$, which gives the
first assertion of (iv). For $\zeta\in\fk$ we have $h(\exp(t\zeta)z)=0$, so $L_z\zeta=0$.

Suppose that $L_z\zeta=0$ for some $\zeta\in\fm$ with $\zeta\ne0$, and let $p_t:=\exp(t\zeta)z$ for
small $t>0$. By (ii), $|h(p_t)|\le C_3t^2|\zeta|^2$. By \cref{ls:lem:tube}\textup{(iii)},
$d(p_t,z)\ge c_0t|\zeta|$, and $p_t\to z$ as $t\to0$. By (i) and
\cref{ls:prop:injective}\textup{(ii)} for the scale $\tfrac18s$, for small $t$,
\begin{equation*}
|h(p_t)|=|h(p_t)-h(z)|\ge c_2\Bigl(\frac{8\,d(p_t,z)}s\Bigr)^{\eta_0}d(p_t,z)
\ge c_2\Bigl(\frac{8c_0|\zeta|}s\Bigr)^{\eta_0}c_0|\zeta|\;t^{1+\eta_0}.
\end{equation*}
Since $\eta_0<1$, the two bounds contradict each other as $t\to0$. So $L_z$ is injective
on $\fm$.
\end{proof}

\begin{remark}[Smaller Parameters]\label{ls:rem:smallpar}
The parameter $\theta$ in \cref{ls:prop:smooth} is free, with the corresponding
thresholds $\eps_L(\theta)$ and $c_L(\theta)$. The slice theorem uses $\theta=\theta_A$; the construction on the charts $Z_i$ uses the smaller parameters fixed in
\cref{th:sec:slicechart,th:sec:along}.
\end{remark}

\subsection{Slices from Equivariant Coordinates}\label{ls:sec:slicesfromcharts}\label{ls:sec:proof}

The orbit derivative is transverse to the level set defining a slice. A single expansion
gives uniqueness of the slice intersection, existence by Brouwer's theorem, and the
quadratic normal error used below. The uniqueness argument uses only the orbit components
of $h$, so it survives correction of the normal components in
\cref{ch:localbundle,ch:globalbundle}.

Fix the data of \cref{ls:prop:smooth}, and put $T:=L_z(\fm)$ and $V:=T^\perp$.
These subspaces are invariant under $\mathsf R(K)$, by part~\textup{(iv)} of that
proposition. Let $P_T$ and $P_V$ be the orthogonal projections onto $T$ and $V$,
and choose $c_z>0$
such that $|L_z\zeta|\ge c_z|\zeta|$ for $\zeta\in\fm$. Parts~\textup{(ii),(iii)} of
\cref{ls:prop:smooth} give, for $p\in B_s(z)$ and $|\zeta|<r_2$,
\begin{equation}\label{ls:eq:orbitexpansion}
\begin{aligned}
h(\exp(\zeta)p)&=h(p)+L_z\zeta+\mathsf e_z(p,\zeta),\\
|\mathsf e_z(p,\zeta)|&\le C_3\bigl(|\zeta|+d(p,z)\bigr)|\zeta|.
\end{aligned}
\end{equation}
Thus the orbit derivative is uniformly close to $L_z$ near $z$, and at $z$ its
normal component vanishes. Choose
\begin{equation*}
0<r_3\le\min\Bigl\{\tfrac12r_0,\tfrac12r_2,\frac{c_z}{4C_3}\Bigr\},\qquad
0<\rho_s\le\min\Bigl\{\tfrac12\sigma_s,s,\frac{c_z}{4C_3}\Bigr\},
\end{equation*}
where $\sigma_s>0$ is chosen so that $d(gz,z)<\sigma_s$ implies
$g\in\exp(\fm(r_3))K$, as permitted by \cref{ls:lem:tube}\textup{(ii)}.

\begin{lemma}[Slices from Equivariant Coordinates]\label{ls:lem:slicesfromcharts}\label{th:lem:slices}
Let $U\subset B_{\rho_s}(z)$ be an open neighbourhood of $z$ that is invariant
under $K$, and let $\check h\colon U\to\R^n$ be injective and continuous, with
$\check h(z)=0$ and
\begin{equation*}
\check h(kp)=\mathsf R(k)\check h(p),\qquad P_T\check h(p)=P_Th(p)
\qquad(k\in K,\ p\in U).
\end{equation*}
Then $\check h$ is a homeomorphism onto an open subset of $\R^n$. For every $r_1>0$
such that $\check h(U)$ contains the closed ball of radius $r_1$ about $0$, put
$S:=\check h^{-1}(V\cap B_{r_1}(0^n))$. The set $S$ is invariant under $K$, and
$\check h$ restricts to an equivariant homeomorphism of $S$ onto
$V\cap B_{r_1}(0^n)$. If $q,q'\in S$ and $gq=q'$ for $g\in G$, then $g\in K$.
Moreover
\begin{equation*}
\Upsilon\colon G\times_KS\longrightarrow Z,\qquad [g,q]\longmapsto gq,
\end{equation*}
is an equivariant homeomorphism onto an open neighbourhood of $Gz$.
\end{lemma}
\begin{proof}
The assumptions of \cref{ls:prop:smooth} include the thresholds for
\cref{ls:prop:injective}; its proof, through \cref{ls:lem:existence}, gives that
$\Omega$ is an $n$-dimensional topological manifold. Since $U\subset\Omega$,
invariance of domain makes $\check h$ a homeomorphism onto an open subset of $\R^n$.
Equivariance makes its image invariant under $\mathsf R(K)$ and gives the stated
properties of $S$.

Suppose $gq=q'$ for $q,q'\in S$. Then
$d(gz,z)\le d(gz,gq)+d(q',z)<2\rho_s\le\sigma_s$, so
$g=\exp(\zeta)k$ for $\zeta\in\fm(r_3)$ and $k\in K$. The point $p:=kq$
lies in $S$, and $\exp(\zeta)p=q'$. By \eqref{ls:eq:orbitexpansion},
\begin{equation*}
|\mathsf e_z(p,\zeta)|\le C_3(r_3+\rho_s)|\zeta|
\le\tfrac12c_z|\zeta|\le\tfrac12|L_z\zeta|.
\end{equation*}
Since $\check h(q'),\check h(p)\in V$ and $P_T\check h=P_Th$, the projection of
$h(q')-h(p)$ onto $T$ is zero. As $L_z\zeta\in T$, the displayed estimate
implies $|L_z\zeta|\le\tfrac12|L_z\zeta|$. Thus $\zeta=0$ and $g=k\in K$.

The map $\Upsilon$ is well defined, continuous and equivariant. If $gq=g'q'$,
then $g'^{-1}g\in K$ by the preceding argument, and
$[g,q]=[g',q']$, so $\Upsilon$ is injective. The principal bundle $G\to G/K$
makes $G\times_KS$ a fibre bundle over the manifold $G/K$, with fibre an open
ball of $V$. Its dimension is
$\dim G-\dim K+\dim V=\dim\fm+\dim V=n$, because
$L_z|_{\fm}\colon\fm\to T$ is an isomorphism. Its image lies in the invariant
manifold $\Omega$, so a second application of invariance of domain proves that
$\Upsilon$ is a homeomorphism onto an open subset of $Z$. Since $z\in S$, that
subset contains $Gz$.
\end{proof}

The slice theorem now follows by applying \cref{ls:lem:slicesfromcharts}
to the chart of \cref{ls:prop:smooth}.
\begin{proof}[Proof of \cref{ls:thm:slices}]
Let $\eps_L$ and $c_L$ be the constants of \cref{ls:sec:smoothing}. The hypothesis of the
theorem is \ref{ls:G} with $\eps=\eps_L$. Every sufficiently small scale is $c_L$-admissible
at $z$; we fix one, $s$, and let $\mathsf R$, $h$, $r_2$, $C_3$ and $L_p$ be as in
\cref{ls:prop:smooth} with $\theta=\theta_A$.

\emph{The decomposition.} Let $T=L_z(\fm)$, $V=T^\perp$ and $c_z$ be as in
\cref{ls:sec:slicesfromcharts}; since $L_z$ vanishes on $\fk$, also $T=L_z(\fg)$. By
\cref{ls:prop:smooth}\textup{(iv)}, $L_z$ restricts to a linear isomorphism of $\fm$ onto $T$
that commutes with the actions of $K$ by $\Ad$ and by $\mathsf R$, and the representation of $K$
on $\fm$ is isomorphic to that on $\fg/\fk$. This is (ii).

\emph{The harmonic chart and the slice.} Choose $r_3$, $\sigma_s$ and $\rho_s$ as in
\cref{ls:sec:slicesfromcharts}, decreasing $\rho_s$ so that ${\rho_s\le s/(800n)}$, and put $U:=B_{\rho_s}(z)$.
The set $U$ is invariant under $K$. By \cref{ls:prop:smooth}\textup{(i)} and
\cref{ls:prop:injective}\textup{(ii)} for the scale $\tfrac18s$, the map $h$ restricts to
a homeomorphism of $U$ onto an open subset of $\R^n$, which contains $0=h(z)$ and is
invariant under $\mathsf R(K)$. We choose $r_1>0$ such that it contains the closed ball of
radius $r_1$ about $0$. This is (i), for the restriction of $h$ to $U$. The set
$S:=h^{-1}(V\cap B_{r_1}(0^n))$ is therefore the slice furnished by
\cref{ls:lem:slicesfromcharts} with $\check h=h$. That lemma gives (iii).
\end{proof}

\begin{remark}[Smaller Slices]\label{ls:rem:smaller}
The number $r_1$ in \cref{ls:thm:slices} can be replaced by any smaller positive number.
Assertions (i) and (ii) are not affected, and the proof of (iii) uses only that $S$ is
the preimage in $U$ of a ball of $V$ about $0$.
\end{remark}

\subsection{The Limit of the Frame Spaces}\label{ls:sec:limits}

The slice theorem now applies to the limit charts of \cref{m:fr:thm:geomain}, with
dimension bound $N+1$ and Reifenberg error $3\delta$. An isotropy-equivariant section
of the frame bundle over a slice identifies the local model of the framed limit.

\begin{theorem}[The Limit of the Frame Spaces]\label{ls:thm:limits}
In the setting of \cref{m:fr:thm:geomain}, suppose that $3\delta\le\eps_L(n,N+1)$. Then every point $y\in Y$, with isotropy group
$I\subset O(n)$, has an open neighbourhood
that is invariant under $O(n)$ and equivariantly homeomorphic to $O(n)\times_IB$, where
$B$ is an open ball about $0$ in a Euclidean space on which $I$ acts through a continuous
homomorphism into the orthogonal group. In particular $Y$ is a topological manifold, and
the action of $O(n)$ on $Y$ is locally linear (\cref{ls:def:locallylinear}), that is, locally
smooth in the sense of \cite[Chapter~IV]{Bre72}.
\end{theorem}

We fix a centre $a$ and use the facts \ref{fl:S:stage}--\ref{fl:S:Yconv} of \cref{fl:sec:setting} with their notation: the limit chart $(Z,\hat d,\hat m)$ is an $\RCD$ space with the dimension
parameter $N+1$, the limit group $G$ is a Lie group that
acts properly by isometries that preserve $\hat m$, the regions $\Omega\supset\Omega'$ are
open and invariant, and $\Omega$ is $(3\delta,R_1)$-Reifenberg. So if
$3\delta\le\eps_L(n,N+1)$, the hypotheses of \cref{ls:thm:slices} hold for every point
$z\in\Omega$. The balls $B_{R_0/8}(a)$ of the centres cover $X$
(\cref{m:fr:thm:geomain}).

We regard $Y$ and the frame spaces as left $O(n)$-spaces, through
$A\cdot y:=y\cdot A^{-1}$.

\begin{proof}[Proof of \cref{ls:thm:limits}]
Let $y\in Y$, and choose a centre $a$ with $\pr(y)\in B_{R_0/8}(a)$. By \ref{fl:S:charts} there is a point $(z,e)$ of $\Fr=\pr^{-1}(\Omega')$ with $\tau^a(z,e)=y$. Let $K$ be the isotropy
group of $z$.

\emph{The isotropy group.} For $k\in K$ let $A_k\in O(n)$ be the matrix with
$e\circ k^{-1}=e\cdot A_k$. By \ref{fl:S:charts}, $k\mapsto A_k$ is an injective homomorphism of
$K$ onto the isotropy group $I$ of $y$. It is continuous, because $G$ acts continuously on the
frames, so it is an isomorphism, since $K$ is compact.

\emph{An equivariant section.} The region $\Omega$ is $\eps_\Omega$-flat at scale $r$ by
\ref{fl:S:limit}, and $\eps_\Omega$ lies below the gauge thresholds $\eps_G(n)$ and
$\eps_M(n,\theta_G(n))$ by its choice in \cref{m:fr:sec:coverbundle}; so
\cref{m:fr:lem:framegauge,m:fr:thm:gauges} hold on it, by
\cref{m:fr:prop:localization}\textup{(i)}. Let $w$ be the restriction of the frame map $u_e$
to $B_{3r}(z)$. By \cref{m:fr:lem:framegauge}\textup{(i)} it is a gauge map with parameter
$\theta_G(n)$ whose core contains $B_{r/2}(z)$. By \cref{m:fr:thm:bundle} and \ref{fl:S:frames},
$u_e\circ k^{-1}=u_{e\circ k^{-1}}=u_{e\cdot A_k}=A_k^{\T}u_e$ for $k\in K$, so
$w\circ k^{-1}=A_k^{-1}w$. Let $\sigma(q):=e_w(q)$ be the frame aligned with $w$ at
$q\in B_{r/2}(z)$; by \cref{m:fr:thm:gauges} it is a continuous section, and $\sigma(z)=e$.
By \cref{m:fr:thm:gauges} we have $e_{Rw}=e_w\cdot R^{\T}$ for $R\in O(n)$ and
$e_{w\circ k^{-1}}(kq)=e_w(q)\circ k^{-1}$. With $R=A_k^{-1}$ this gives
\begin{equation}\label{ls:eq:section}
\sigma(q)\circ k^{-1}=\sigma(kq)\cdot A_k\qquad\text{for all }k\in K\text{ and }
q\in B_{r/2}(z).
\end{equation}

\emph{The frame space over a tube.} Let $\mathsf R$, $V$, $h$, $r_1$ and $S$ be as in
\cref{ls:thm:slices} for the point $z$; by \cref{ls:rem:smaller} we may assume that
$S\subset B_{r/2}(z)\cap\Omega'$. The group $K$ acts on $S\times O(n)$ by the rule
$k\cdot(q,A):=(kq,A_kA)$. The map
\begin{equation*}
G\times_K\bigl(S\times O(n)\bigr)\longrightarrow\pr^{-1}(GS),\qquad
[g,(q,A)]\longmapsto g\cdot\bigl(q,\sigma(q)\cdot A\bigr),
\end{equation*}
is well defined by \eqref{ls:eq:section}. It is a map of principal $O(n)$-bundles over the
homeomorphism $G\times_KS\to GS$ of \cref{ls:thm:slices}\textup{(iii)}, hence a
homeomorphism, and it commutes with the actions of $G$. By \ref{fl:S:charts} the map $\tau^a$ is open and induces a homeomorphism of $\Fr/G$ onto its image, and $\pr^{-1}(GS)$ is an
open invariant subset of $\Fr$. So $\tau^a$ induces a homeomorphism of the quotient of
$S\times O(n)$ by $K$ onto an open subset of $Y$ that is invariant under $O(n)$; it
commutes with the right multiplication on the factor $O(n)$, and it takes the class of
$(z,\one)$ to $y$.

\emph{The model.} Let $B:=V\cap B_{r_1}(0^n)$, and let $I$ act on $B$ through the
homomorphism $A_k\mapsto\mathsf R(k)$. The map $h$ identifies the quotient of $S\times O(n)$
by $K$ with the quotient of $B\times O(n)$ by the action $A_k\cdot(v,A):=(\mathsf R(k)v,A_kA)$ of $I$,
and $[v,A]\mapsto[A^{-1},v]$ identifies the latter with $O(n)\times_IB$, compatibly with
the left actions of $O(n)$. So $y$ has an invariant open neighbourhood that is
equivariantly homeomorphic to $O(n)\times_IB$. This space is the quotient of the manifold
$O(n)\times B$ by a free smooth action of the compact Lie group $I$, hence a manifold,
and the action of $O(n)$ on $Y$ is locally linear by \cref{ls:def:locallylinear}.
\end{proof}

\subsection{Retraction Properties and Realizing Maps}\label{fl:sec:anr}

The locally linear manifold $Y$ admits an equivariant neighbourhood retraction.
Restricting it to $Y_i$ produces equivariant realizing maps, with no further subsequence.

\begin{definition}[Realizing Maps]\label{fl:def:realizing}
Continuous maps $F_i\colon Y_i\to Y$, defined for all large $i$, are
\emph{realizing maps} if $F_i(y_i)\to y$ whenever points $y_i\in Y_i$ converge to
$y\in Y$ (\cref{m:fr:def:Yconv}). They are \emph{equivariant realizing maps} if, in addition, every $F_i$
commutes with the $O(n)$-actions.
\end{definition}

The induced quotient maps are the maps $f_i$ of \cref{thm:main}.

\begin{theorem}[Retraction Property and Equivariant Realizing Maps]\label{fl:thm:anr}
In the geometric setting of \cref{fl:sec:setting}, assume the hypotheses of
\cref{ls:thm:limits}.
\begin{enumerate}[label=\textup{(\roman*)}]
\item The space $Y$ is an $O(n)$-absolute neighbourhood retract; in particular,
it is an absolute neighbourhood retract.
\item There are realizing maps that are equivariant for $O(n)$.
\item Equivariant realizing maps $F_i$ induce continuous maps $f_i\colon X_i\to X$ with
$\pr\circ F_i=f_i\circ\pr$, and $f_i(x_i)\to x$ whenever points $x_i\in X_i$ converge to
$x\in X$.
\end{enumerate}
\end{theorem}

Part~(i) is what makes (ii) possible. The convergence $Y_i\to Y$ is defined through approximations that are not necessarily continuous, and in the common realization $\cY$ of
\cref{fl:prop:realization} the nearest-point maps onto $Y$ need not be continuous either. A
continuous map from a neighbourhood of $Y$ onto $Y$ that fixes the points of $Y$, a neighbourhood
retraction, does the job: its domain contains $Y_i$ for large $i$, and its restrictions to the
spaces $Y_i$ are realizing maps. This is the construction of maps between close spaces in
\cite{Pet90}. The retraction has to commute with $O(n)$, so that the maps descend to the
quotients in~(iii).

Spaces with such retractions are the absolute neighbourhood retracts (ANR). A metrizable space is an
\emph{ANR} if, wherever it lies as a closed subset of a metrizable space, it is a retract of a
neighbourhood. The model is a compact submanifold of a Euclidean space, onto which the
nearest-point projection retracts a tubular neighbourhood; the ANR property replaces this
projection where there is no smooth structure. Manifolds and finite polyhedra are ANRs, while a
space with wild local topology, such as the topologist's sine curve, is not. For a compact Lie
group $G$, a metrizable $G$-space is a \emph{$G$-ANR} if every realization of it as a closed
invariant subspace of a metrizable $G$-space admits an equivariant retraction from an invariant
neighbourhood, which may be taken open. This asks more: near a point the retraction has to
respect its isotropy group. The local linearity of the action on $Y$, which the slice theorem
provides, is what gives it.

\begin{proof}[Proof of \cref{fl:thm:anr}]
(i) By \cref{ls:thm:limits}, $Y$ is a compact metrizable manifold
with a locally linear $O(n)$-action. This is local smoothness in the terminology of
\cite[p.~501]{Mur83}: each point has a slice (\cref{ls:def:slice}) equivariantly homeomorphic
to a Euclidean representation of its isotropy group. The slices in
\cref{ls:thm:limits} are open balls in such representations;
radial expansion identifies them equivariantly with the full representations.
Thus \cite[Theorem~8.8]{Mur83} applies and makes $Y$ an $O(n)$-ANR.
It is also an ordinary ANR by \cite[Corollary~8.2\textup{(2)}]{Mur83}.

(ii) By \cref{fl:prop:realization,fl:lem:invariant}, $\cY$ is a compact metrizable
$O(n)$-space with an invariant metric, in which $Y$ is closed and invariant. By (i) there
is a continuous equivariant retraction $F\colon U\to Y$ of an invariant open neighbourhood
$U$ of $Y$ in $\cY$. For large $i$ the space $Y_i$ lies in $U$, since otherwise points $y_k\in Y_{i_k}\setminus U$ with $i_k\to\infty$ would have a subsequence converging to a
point of $Y$. The maps $F_i:=F|_{Y_i}$ are equivariant, and they are realizing maps: if $y_i\to y\in Y$, then $F_i(y_i)=F(y_i)\to F(y)=y$.

(iii) The fibres of $\pr$ on $Y_i$ and on $Y$ are the orbits, and $\pr$ is an open map on
$Y_i$; so an equivariant continuous map $F_i$ induces a continuous map $f_i$ of the
quotients. Let $x_i\in X_i$ converge to $x\in X$, and choose points $y_i\in Y_i$ over
them. By \cref{fl:S:Yconv} every infinite set of indices contains an infinite subset along
which $y_i$ converges to a point $y'$ of $Y$. Then $\pr(y_i)\to\pr(y')$ along it, so
$\pr(y')=x$, and $f_i(x_i)=\pr(F_i(y_i))$ converges to $\pr(y')=x$ along it. Hence
$f_i(x_i)\to x$.
\end{proof}

\subsection{Quantitative Slice Estimates}\label{th:sec:limitchart}

Here we specialize the harmonic chart of \cref{ls:prop:smooth} and the slice criterion of
\cref{ls:lem:slicesfromcharts} to the accuracy needed on the charts $Z_i$. We prove that the normal components are equivariant
up to an error quadratic in the distance and that nearby orbits meet the slice.
We keep hypothesis \ref{ls:G} and the notation of \cref{ls:sec:smoothing}.

\subsubsection{The Harmonic Chart and Its Decomposition}\label{th:sec:slicechart}

The harmonic chart is that of \cref{ls:prop:smooth}, with the parameter $\theta_\flat$ of
\cref{th:sec:constants} in place of the parameter used there. Let
\begin{equation*}
\eps_S:=\min\{\eps_\flat,\eps_L(\theta_\flat)\},\qquad c_S:=\min\{c_\flat,c_L(\theta_\flat)\},
\end{equation*}
with the constants $\eps_L(\theta)$ and $c_L(\theta)$ of \cref{ls:sec:smoothing}; they depend on
$n$, $N$ and the fixed number $\eps_*$.

Assume \ref{th:G} with $\eps\le\eps_S$, and choose
$s_0>0$ so that $8s_0$ is $c_S$-admissible at $z$. Apply
\cref{ls:prop:smooth} with $s=8s_0$ and $\theta=\theta_\flat$. Throughout this subsection and the constructions of \cref{ch:stages,ch:localbundle}, references to that proposition use this convention:
$h$ is defined on $B_{12s_0}(z)$, its restriction to $B_{8s_0}(z)$ has harmonic
chart parameter $\theta_\flat$ at $(z,s_0)$, and the maps $L_p$ and the estimates
\textup{(ii),(iii)} apply for $p\in B_{8s_0}(z)$. 
For this chart use $T$, $V$, their orthogonal projections and $c_z$ as in
\cref{ls:sec:slicesfromcharts}, and put
\begin{equation*}
V_0:=\{v\in V:\mathsf R(k)v=v\text{ for all }k\in K\},\qquad V':=V\cap V_0^\perp.
\end{equation*}
Let $P_0$ and $P'$ be the orthogonal projections onto $V_0$ and $V'$.
By \cref{ls:prop:smooth}\textup{(iv)}, all four subspaces are invariant under
$\mathsf R(K)$; their projections commute with $\mathsf R(K)$, and
$P_0\,\mathsf R(k)=P_0$. Put $k_0:=\dim V_0$ and $m':=\dim V'$. The representation of $K$ on $V$ is the slice representation at $z$
(\cref{ls:thm:slices}); the subspace $V_0$ is the direction of the stratum through $z$
within the slice, and $V'$ is normal to it in $V$.

\subsubsection{Almost Equivariance}\label{th:sec:almostequiv}

The harmonic chart commutes with $K$. The next lemma says that its components in $V$ commute with
all elements of $G$ that move a point near $z$ to a point near $z$, up to the square of
the distance from $z$. The reason is that the derivative of the harmonic chart along the group at
$z$ takes values in $T$.

\begin{lemma}[Almost Equivariance]\label{th:lem:almostequiv}
In the situation of \cref{th:lem:chart} there are numbers $\rho_2$ and $C_5$, with
$0<\rho_2\le4s_0$, that have the following property. Let $g\in G$ and $p\in B_{\rho_2}(z)$ with
$gp\in B_{\rho_2}(z)$, and let $D$ be the larger of the numbers $d(p,z)$ and $d(gp,z)$.
Then there are a unique $\zeta\in\fm(\tfrac12r_0)$ and a unique $k\in K$ with
$g=\exp(\zeta)k$, and
\begin{equation*}
|\zeta|\le\frac{2D}{c_0},\qquad
\bigl|h(gp)-\mathsf R(k)\,h(p)-L_z\zeta\bigr|\le C_5D^2,\qquad
\bigl|P_Vh(gp)-\mathsf R(k)\,P_Vh(p)\bigr|\le C_5D^2 .
\end{equation*}
\end{lemma}
\begin{proof}
Let $r_*:=\min\{\tfrac12r_0,r_2\}$. By \cref{ls:lem:tube}\textup{(ii)} there is
$\sigma_z>0$ such that every $g\in G$ with $d(gz,z)<\sigma_z$ lies in
$\exp(\fm(r_*))K$. Put $\rho_2:=\min\{\tfrac12\sigma_z,4s_0\}$ and
$C_5:=C_3(4c_0^{-2}+2c_0^{-1})$. Let $g$ and $p$ be as in the statement. Then
$d(gz,z)\le d(gz,gp)+d(gp,z)\le2D<\sigma_z$, so $g=\exp(\zeta)k$ with $\zeta\in\fm(r_*)$ and
$k\in K$, and the pair is unique by \cref{ls:lem:tube}\textup{(i)}. Since
$\exp(\zeta)z=gz$, \cref{ls:lem:tube}\textup{(iii)} gives $|\zeta|\le2D/c_0$. Apply
\eqref{ls:eq:orbitexpansion}, with $s=8s_0$, to $p':=kp\in B_{8s_0}(z)$. Since
$d(p',z)=d(p,z)$ and $\exp(\zeta)p'=gp$,
\begin{equation*}
|\mathsf e_z(p',\zeta)|\le C_3(2D/c_0+D)(2D/c_0)=C_5D^2,
\end{equation*}
and $h(p')=\mathsf R(k)h(p)$. The last inequality follows, because $L_z\zeta$ lies in $T$ and
$P_V$ commutes with $\mathsf R(k)$.
\end{proof}

\subsubsection{Slices from Corrected Coordinates}\label{th:sec:slicesfromcharts}

We use \cref{ls:lem:slicesfromcharts} for corrected coordinates, with the following
choice of its radii. Let
\begin{equation*}
r_3:=\min\Bigl\{\tfrac12r_0,\ \tfrac12r_2,\ \frac{c_z}{4C_3}\Bigr\},
\end{equation*}
made smaller, if necessary, so that $d(\exp(\zeta)p,p)\le s_0/(200n)$ for all
$\zeta\in\overline{\fm(r_3)}$ and $p\in B_{8s_0}(z)$; this is possible because the action is
continuous and closed balls are compact. Let $\sigma_s>0$ be such that every $g\in G$ with
$d(gz,z)<\sigma_s$ lies in $\exp(\fm(r_3))K$, by \cref{ls:lem:tube}\textup{(ii)}, and put
\begin{equation*}
\rho_s:=\min\Bigl\{\tfrac12\sigma_s,\ \frac{c_z}{4C_3},\ \rho_2\Bigr\}.
\end{equation*}

Since $\rho_s\le\rho_2\le4s_0<8s_0=s$, these choices satisfy every radius
condition of \cref{ls:sec:slicesfromcharts}. Thus the criterion applies to any
injective continuous $K$-equivariant $\check h$ on a $K$-invariant open
$U\subset B_{\rho_s}(z)$ with $\check h(z)=0$ and $P_T\check h=P_Th$.

The following lemma moves a point near $z$ along its orbit until the
components of the harmonic chart in $T$ vanish.

\begin{lemma}[Orbits Meet the Slice]\label{th:lem:meets}
In the situation of \cref{th:lem:chart} let $p\in B_{\rho_s}(z)$, and suppose that
$t:=|P_Th(p)|$ satisfies $2t\le c_zr_3$. Then there is $\zeta\in\fm$ with
$|\zeta|\le2t/c_z$ such that the point $q:=\exp(\zeta)p$ satisfies
\begin{equation*}
P_Th(q)=0,\qquad|P_Vh(q)-P_Vh(p)|\le t,\qquad d(q,p)\le\frac{s_0}{200\,n} .
\end{equation*}
\end{lemma}
\begin{proof}
Let $\Lambda\colon T\to\fm$ be the inverse of the restriction of $L_z$ to $\fm$; its norm
is at most $1/c_z$. Put $R:=2t/c_z\le r_3$, and for $\zeta\in\fm$ with $|\zeta|\le R$ let
\begin{equation*}
\Phi(\zeta):=-\Lambda\bigl(P_Th(p)\bigr)-\Lambda\bigl(P_T\mathsf e_z(p,\zeta)\bigr).
\end{equation*}
By \eqref{ls:eq:orbitexpansion}, with $s=8s_0$,
$|\mathsf e_z(p,\zeta)|\le C_3(R+\rho_s)|\zeta|\le\tfrac12c_z|\zeta|$, so
$|\Phi(\zeta)|\le t/c_z+\tfrac12R=R$. The map $\Phi$ is continuous, so by the fixed point
theorem of Brouwer there is $\zeta$ with $\Phi(\zeta)=\zeta$ and $|\zeta|\le R$. Applying
$L_z$ gives $L_z\zeta=-P_Th(p)-P_T\mathsf e_z(p,\zeta)$, that is, $P_Th(\exp(\zeta)p)=0$. Moreover
$P_Vh(q)-P_Vh(p)=P_V\mathsf e_z(p,\zeta)$ has norm at most $\tfrac12c_z|\zeta|\le t$, and the bound for
$d(q,p)$ holds by the choice of $r_3$.
\end{proof}

\subsection{Uniform Constants under Recentring}\label{th:sec:along}\label{th:sec:uniform}

The preceding estimates, and the constants of the construction on the charts $Z_i$ in \cref{ch:stages} that are built on them, are centred at $z$. One harmonic chart serves
all nearby points with the same isotropy group, after a translation, so these constants
can be chosen uniformly on compact sets of such points. This supplies the local
uniformity used in \cref{ch:globalbundle} (\cref{es:lem:uniform}).
The harmonic chart is taken at twice the scale and with half the parameter, as
\cref{ls:rem:smallpar} allows; let
\begin{equation*}
\eps'_S:=\min\{\eps_S,\eps_L(\tfrac12\theta_\flat)\},\qquad c'_S:=\min\{c_S,c_L(\tfrac12\theta_\flat)\},
\end{equation*}
with the numbers $\eps_S$ and $c_S$ of \cref{th:sec:slicechart} and the constants of that remark; so
$\eps'_S\le\eps_S$ and $c'_S\le c_S$.

\begin{lemma}[Recentring]\label{th:lem:recentre}
Assume \ref{th:G} with $\eps\le\eps'_S$, and let $s_0>0$ be such that $16s_0$ is
$c'_S$-admissible at $z$. Then there are a homomorphism $\mathsf R$, a harmonic map
$h\colon B_{24s_0}(z)\to\R^n$, numbers $r_2$ and $C_3$ and linear maps $L_p$ with the
properties of \cref{ls:prop:smooth} for $s=16s_0$ and
$\theta=\tfrac12\theta_\flat$. For these data there are numbers
$\rho_z\in(0,s_0]$ and $r'_2\in(0,r_2]$ with the following property. Let
$z'\in B_{\rho_z}(z)$ be a point whose isotropy group is $K$. Then all assertions of
\cref{ls:prop:smooth} hold for the point $z'$, $s=8s_0$ and
$\theta=\theta_\flat$, with the same $\mathsf R$, $C_3$
and $L_p$, with $r'_2$ in place of $r_2$, and with the map $h-h(z')$. Moreover
$|L_{z'}\zeta|\ge\tfrac12c_z|\zeta|$ for $\zeta\in\fm$, and the numbers $r_0$, $c_0$,
$\sigma_z$ and $\sigma_s$ that enter \cref{th:sec:almostequiv,th:sec:slicesfromcharts} can be
chosen independently of $z'$.
\end{lemma}
\begin{proof}
The data exist by \cref{ls:prop:smooth} for $s=16s_0$ and
$\theta=\tfrac12\theta_\flat$. Let $\rho_z\le\min\{s_0,c_z/(2C_3)\}$, and let
$z'$ be as in the statement. The ball $B_{12s_0}(z')$ lies in $B_{24s_0}(z)$, and
$B_{8s_0}(z')$ lies in $B_{16s_0}(z)$. Let $\Phi$ be a map on $B_{16s_0}(z)$ that changes
distances by at most $\tfrac12\theta_\flat\cdot2s_0$ with $|h-\Phi|\le\theta_\flat s_0$. Then
$\Phi-h(z')$ shows that the restriction of $h-h(z')$ to $B_{8s_0}(z')$ is a harmonic
chart with parameter $\theta_\flat$ at $(z',s_0)$. Since $z'$ is fixed by $K$, we have
$\mathsf R(k)h(z')=h(kz')=h(z')$, so $h-h(z')$ commutes with $K$. This is \textup{(i)}. The
assertions \textup{(ii)} and \textup{(iii)} for the points of $B_{8s_0}(z')$ are those
for $h$, once $r'_2$ is so small that $\hat d(\exp(\zeta)p,p)<4s_0$ for $|\zeta|<r'_2$ and
$p\in B_{16s_0}(z)$, which is possible by compactness. By \textup{(ii)} the map
$L_{z'}$ is the derivative of $\zeta\mapsto h(\exp(\zeta)z')$ at $0$. For $\zeta\in\fk$ this
curve is constant, so $L_{z'}$ vanishes on $\fk$; and
$h(\exp(\Ad(k)\zeta)z')=h(k\exp(\zeta)z')=\mathsf R(k)h(\exp(\zeta)z')$ gives
$L_{z'}\circ\Ad(k)=\mathsf R(k)L_{z'}$. By \textup{(iii)},
$|L_{z'}\zeta|\ge|L_z\zeta|-C_3\rho_z|\zeta|\ge\tfrac12c_z|\zeta|$ for $\zeta\in\fm$.
This is \textup{(iv)}.

The number $r_0$ of \cref{ls:lem:tube}\textup{(i)} depends on $G$, $K$ and $\fm$ only.
If $\fm=0$, the number $c_0$ is arbitrary. Otherwise it is obtained in the proof of
\cref{ls:lem:tube}\textup{(iii)} from the
minimum of $\hat d(\exp(\zeta)z,z)$ over the $\zeta\in\fm$ with
$\tfrac14r_0\le|\zeta|\le\tfrac12r_0$; this minimum changes by at most $2\rho_z$ when
$z$ is replaced by $z'$, so half of it serves for all $z'$ if $\rho_z$ is small. If
every $g$ with $\hat d(gz,z)<\sigma$ lies in a neighbourhood $U$ of $K$ and
$\rho_z\le\tfrac14\sigma$, then every $g$ with $\hat d(gz',z')<\tfrac12\sigma$ lies in
$U$, because $\hat d(gz,z)\le\hat d(gz',z')+2\rho_z$.
\end{proof}

\section{Transplanted Coordinates and Deck Labels}\label{ch:stages}

We transfer the controlled limit coordinates of \cref{ch:harmoniccomparison}
to the approximating charts $Z_i$. \Cref{th:sec:stages} fixes the standing setting and transplants the harmonic chart. \Cref{es:sec:lambda} defines the
component in the isotropy group of the elements of $G$ near it, \cref{es:sec:labels} uses
it to label the deck transformations with small displacement, and
\cref{th:sec:stageequiv} compares the actions of the labels in the transplanted
coordinates. The scale and the cut-off weights of \cref{th:sec:averages} then prepare
the exact cocycle and the corrected coordinates of \cref{ch:localbundle}.

\subsection{Coordinates on the Charts of the Sequence}\label{th:sec:stages}

We transplant the harmonic chart of \cref{th:sec:limitchart} through convergent
harmonic charts on $Z_i$. In the transplanted chart, deck transformations that move $z_i$
by at most $8s_0$ act by orthogonal maps up to an error with a small Lipschitz constant
(\cref{th:lem:transplant}).

\subsubsection{The Standing Setting}\label{th:sec:standing}\label{es:sec:setting}\label{es:sec:standing}

We use the common geometric setting and subsequence of \cref{fl:sec:setting},
including the cover, group and frame data of
\textup{\ref{fl:S:stage}--\ref{fl:S:Yconv}}, and all convergence assertions are those of
\textup{\ref{fl:S:limit}}, \textup{\ref{fl:S:groups}} and \textup{\ref{fl:S:Yconv}}. Fix a centre $a$ and
suppress its superscript. The approximations $\af_i$ are those of
\textup{\ref{fl:S:limit}}; they are defined on balls $B_{1/\eps_i}(z_i)$ with $\eps_i\to0$. The orbit projection $\pi\colon\Omega'\to X$ and the maps $\pi_i$, $\tau_i$,
$\tau$ satisfy the projection identities in \textup{\ref{fl:S:charts}}. In particular \ref{th:R}
holds with the scale $\lambda=R_1$ and with errors $\delta$ on the spaces $Z_i$ and $3\delta$ on the limit. All constants of
\cref{th:sec:lipschitz,th:sec:limitchart} are taken for the dimension bound $N+1$, and we assume
\begin{equation}\label{th:eq:threshold}
3\delta\le\eps_T(n,N):=\min\{\eps_S(n,N+1),\,\eps_L(n,N+1)\}.
\end{equation}
The dependence of $\eps_T$ and $\eps_S$ on the already fixed $\eps_*$
is suppressed; all local thresholds use that same value. The constant $\eps_L$ is that of
\cref{ls:thm:slices}, so \eqref{th:eq:threshold} gives the hypothesis $3\delta\le\eps_L(n,N+1)$ of
\cref{ls:thm:limits}. Thus \cref{ls:thm:limits} applies, and there are equivariant realizing maps
$F_i\colon Y_i\to Y$, with induced maps $f_i\colon X_i\to X$ (\cref{fl:thm:anr}). We fix
such maps and an invariant metric $\sd$ on $Y$, and we regard $Y$ and the frame spaces as
left $O(n)$-spaces, as in \cref{ls:sec:limits}. All realizing maps used in the local
construction and the global gluing are equivariant; the term \emph{realizing} itself has
the generic meaning of \cref{fl:def:realizing}.

We fix a point $y\in Y$ with $\pr(y)\in B_{R_0/4}(a)$ and a frame $(z,e)\in\Fr$ with
$\tau(z,e)=y$, and we put $x:=\pr(y)=\pi(z)$; then $z\in\Omega^\flat$. Let $K$ be the isotropy group of $z$ in $G$ and $I$
that of $y$ in $O(n)$. For $k\in K$ let $A_k\in O(n)$ be the matrix with
$e\circ k^{-1}=e\cdot A_k$; the map $k\mapsto A_k$ is an isomorphism of $K$ onto $I$
(proof of \cref{ls:thm:limits}). The setting \ref{th:G} holds at $z$ with
$\eps=3\delta$. We fix $s_0>0$ so small that $8s_0$ is $c_S$-admissible at $z$ and
$18s_0\le\tfrac18R_0$. For this scale fix the harmonic chart, decomposition,
orbit derivative and radii of \cref{th:sec:limitchart}, with the conventions
$s=8s_0$ and $\theta=\theta_\flat$ stated there. We write $B_V(r)$ for the open ball of radius $r$ about $0$ in $V$; note that $\dim V=n-\dim\fm=n-\dim G+\dim K$.
We choose points $z_i\in Z_i$ that converge to $z$; for large $i$ the ball
$B_{18s_0}(z_i)$ lies in $\Omega'_i$, and the scale $2s_0$ is $c_\flat$-admissible at $z_i$.
We write
\begin{equation*}
s_1:=\frac{s_0}{100\,n},\qquad B^i:=B_{s_1}(z_i),\qquad B:=B_{s_1}(z),\qquad
C_6:=2\sqrt n\,C_1,\qquad c_3:=\tfrac14c_2 .
\end{equation*}
As in \cref{m:fr:sec:conv}, but about the points $z_i$, put
\begin{equation*}
\Gamma_i(t):=\{\gamma\in\Gamma_i:\hat d(\gamma z_i,z_i)<t\}\qquad(t>0).
\end{equation*}
This set is invariant under inversion, and
$\Gamma_i(t)\Gamma_i(t')\subset\Gamma_i(t+t')$.

\subsubsection{Harmonic Charts That Converge}\label{th:sec:convcharts}

We choose harmonic charts on the spaces $Z_i$, pass to a limit, and compare the limit with the
harmonic chart $h$. The comparison is a map with a small Lipschitz constant relative to the
chart, and its extension to $\R^n$ transplants $h$ to the charts $Z_i$ with the same H\"older
control.

\begin{lemma}[Harmonic Charts on $Z_i$]\label{th:lem:stagecharts}
After passing to a subsequence there are, for all large $i$, harmonic maps
$v_i\colon B_{16s_0}(z_i)\to\R^n$ with $C_N\sqrt n$-Lipschitz components and a map
$v\colon B_{16s_0}(z)\to\R^n$ with the following properties.
\begin{enumerate}[label=\textup{(\roman*)}]
\item The map $v_i$ is a harmonic chart with parameter $\tfrac13\theta_\flat$ at
$(z_i,2s_0)$.
\item If points $p_i\in B_{16s_0}(z_i)$ converge to $p\in B_{16s_0}(z)$, then
$v_i(p_i)\to v(p)$.
\item The map $v$ is a harmonic chart with parameter $\theta_\flat$ at $(z,2s_0)$.
\item The map $\mathsf e:=h-v$ is $\tfrac12\eps_*$-Lipschitz with respect to $v$ on $B$, and
$\mathsf e(z)=0$.
\end{enumerate}
\end{lemma}
\begin{proof}
Let $v^0_i$ be a comparison map at $(z_i,2s_0)$ in $Z_i$ (\cref{th:def:comparison}); its
components are $C_N$-Lipschitz, and it is a harmonic chart with parameter
$\tfrac13\theta_\flat$ at $(z_i,2s_0)$. The maps $v^0_i$ vanish at $z_i$ and are uniformly
Lipschitz, so a subsequence converges in the sense of (ii) to a Lipschitz map $v^0$ on
$B_{16s_0}(z)$, by the theorem of Arzel\`a and Ascoli for converging spaces. The
components of $v^0$ are harmonic: they are limits of harmonic functions that are
uniformly Lipschitz, so the convergence is strong in $L^2$ on balls, and
\cite[Corollary~4.5]{AH18} applies. The maps
$v^0_i$ change the distances of the points of $B_{16s_0}(z_i)$ by at most $2\theta_\flat s_0$,
and distances converge, so $v^0$ changes distances by at most $2\theta_\flat s_0$. Hence
$v^0$ is a harmonic chart with parameter $\theta_\flat$ at $(z,2s_0)$, with $\Phi:=v^0$ in
the definition, and its restriction to $B_{8s_0}(z)$ is a harmonic chart with parameter
$2\theta_\flat$ at $(z,s_0)$. The map $h$ changes the distances of the points of
$B_{8s_0}(z)$ by at most $3\theta_\flat s_0$. By \cref{th:lem:twocharts} at $(z,s_0)$ there is a
rigid motion $E$ such that $h-E\circ v^0$ is $\tfrac12\eps_*$-Lipschitz with respect to
$v^0$ on $B$. Put $v_i:=R_Ev^0_i$ and $v:=R_Ev^0$. Then (i), (ii) and (iii) hold, since
$R_E$ is orthogonal, and $h-v$ differs from $h-E\circ v^0$ by the constant $b_E$, so it
is $\tfrac12\eps_*$-Lipschitz with respect to $v$; it vanishes at $z$ because $h$ and
$v$ do.
\end{proof}

We fix the subsequence and the maps of \cref{th:lem:stagecharts}. By \eqref{th:eq:holder} the
map $v$ is injective on $B$. So $\mathsf e\circ v^{-1}$ is a map of the subset $v(B)$ of $\R^n$ to
$\R^n$ that is $\tfrac12\eps_*$-Lipschitz, and by the theorem of Kirszbraun
\cite[2.10.43]{Fed69} it extends to a map $\hat{\mathsf e}\colon\R^n\to\R^n$ with the same Lipschitz
constant.

\begin{definition}[The Transplanted Chart]\label{th:def:transplant}
Put $\hat h:=\mathrm{id}+\hat{\mathsf e}\colon\R^n\to\R^n$. The \emph{transplanted chart} on $Z_i$ is
\begin{equation*}
h_i:=\hat h\circ v_i\colon B_{16s_0}(z_i)\to\R^n .
\end{equation*}
\end{definition}

Then $(1-\tfrac12\eps_*)|\nu-\nu'|\le|\hat h(\nu)-\hat h(\nu')|\le(1+\tfrac12\eps_*)|\nu-\nu'|$
for all $\nu,\nu'$, and $\hat h\circ v=h$ on $B$. The maps $h_i$ are the transplants of the
harmonic chart $h$ to the spaces $Z_i$; they are not harmonic.

\begin{lemma}[The Transplanted Chart]\label{th:lem:transplant}
For all large $i$ the following holds.
\begin{enumerate}[label=\textup{(\roman*)}]
\item The map $h_i$ is injective on $B^i$ and maps it onto an open subset of $\R^n$, we
have $h_i(z_i)=0$, and
\begin{equation}\label{th:eq:stageholder}
c_3\,\Bigl(\frac{\hat d(p,q)}{s_0}\Bigr)^{\eta}\,\hat d(p,q)\le|h_i(p)-h_i(q)|\le C_6\,\hat d(p,q)
\qquad\text{for all }p,q\in B^i .
\end{equation}
The same holds for $h$ on $B$.
\item For $0<t<s_1$ the image of $B_t(z_i)$ under $h_i$ contains the open ball of radius
$c_3(t/s_0)^{\eta}t$ about $0$, and so does the image of $B_t(z)$ under $h$.
\item If points $p_i\in B_{16s_0}(z_i)$ converge to $p\in B$, then $h_i(p_i)\to h(p)$.
\item For every $\gamma\in\Gamma_i$ with $\hat d(\gamma z_i,z_i)\le8s_0$ there is
$R_\gamma\in O(n)$ such that $h_i\circ\gamma-R_\gamma h_i$ is $3\eps_*$-Lipschitz with
respect to $h_i$ on $B^i$.
\end{enumerate}
\end{lemma}
\begin{proof}
(i) The ball $B^i$ lies in the ball of radius $2s_0/(100n)$ about $z_i$, on which
\eqref{th:eq:holder} holds for $v_i$ at the scale $2s_0$. Since $\hat h$ changes distances
by factors between $1-\tfrac12\eps_*$ and $1+\tfrac12\eps_*$, this gives
\eqref{th:eq:stageholder}, and $h_i(z_i)=\hat h(0)=\hat{\mathsf e}(v(z))=\mathsf e(z)=0$. The map $v_i$ is a
homeomorphism of $B^i$ onto an open set, and $\hat h$ is a homeomorphism of $\R^n$, by
the invariance of domain and because it is proper. For $h$ we use \eqref{th:eq:holder} at
the scale $s_0$.

(ii) Let $U_t:=h_i(B_t(z_i))$, an open set that contains $0$, and suppose that a point
$\nu$ with $|\nu|<c_3(t/s_0)^{\eta}t$ does not lie in $U_t$. On the segment from $0$ to $\nu$
there is a first point $\nu_1$ outside $U_t$. It is the limit of points $h_i(p_j)$ with
$p_j\in B_t(z_i)$, and a subsequence of the points $p_j$ converges to a point $p$ of the
closed ball, which is a compact subset of $B^i$. Then $h_i(p)=\nu_1$, so $p$ does not lie
in $B_t(z_i)$, that is, $\hat d(p,z_i)=t$, and \eqref{th:eq:stageholder} gives
$|\nu_1|\ge c_3(t/s_0)^{\eta}t>|\nu|$, a contradiction. The proof for $h$ is the same.

(iii) We have $v_i(p_i)\to v(p)$, and $\hat h(v(p))=h(p)$.

(iv) The scale $2s_0$ is $c_\flat$-admissible at $z_i$, and $\gamma$ is an isometry that
preserves $\hat m_i$. By \cref{th:cor:isometries} there is a rigid motion $E_\gamma$ such that
$\mathsf e_\gamma:=v_i\circ\gamma-E_\gamma\circ v_i$ is $\tfrac12\eps_*$-Lipschitz with respect to
$v_i$ on $B^i$. Let $R_\gamma$ be the rotation part of $E_\gamma$. On $B^i$,
\begin{equation*}
h_i\circ\gamma-R_\gamma h_i=b_{E_\gamma}+\mathsf e_\gamma+\hat{\mathsf e}\circ(E_\gamma\circ v_i+\mathsf e_\gamma)
-R_\gamma\,\hat{\mathsf e}\circ v_i .
\end{equation*}
The three maps on the right that are not constant are Lipschitz with respect to $v_i$
with the constants $\tfrac12\eps_*$, $\tfrac12\eps_*(1+\tfrac12\eps_*)$ and
$\tfrac12\eps_*$, and $|v_i(p)-v_i(q)|\le(1-\tfrac12\eps_*)^{-1}|h_i(p)-h_i(q)|$. Since
$\eps_*\le\tfrac1{200}$, the sum of the constants, divided by $1-\tfrac12\eps_*$, is less
than $3\eps_*$.
\end{proof}

\subsection{The Component in the Isotropy Group}\label{es:sec:lambda}

An element of $G$ near $K$ is the product of the exponential of a small vector of $\fm$
and of an element of $K$, and we call the second factor its component in $K$
(\cref{es:def:component}). This subsection fixes the notation for this component and for its failure to be
multiplicative. By \cref{m:fr:prop:limitgroup}\textup{(ii),(iii)}, the identity
component $K_0$ of $K$ lies in the central torus of $G_0$ and is therefore a torus.
So $\fk$ is abelian, the exponential map $\exp\colon\fk\to K_0$ is a
surjective homomorphism with a discrete kernel $\Lambda_K$, and
$k\exp(\mathsf b)k^{-1}=\exp(\Ad(k)\mathsf b)$ for $k\in K$ and $\mathsf b\in\fk$. We use on
$\fk$ the restriction of the norm of \cref{ls:sec:isotropy}, which is invariant under $\Ad(K)$,
write $\fk(r)$ for its open ball of radius $r$, and put
\begin{equation*}
r_K:=\tfrac14\min\{|\mathsf b|:\mathsf b\in\Lambda_K,\ \mathsf b\ne0\},
\end{equation*}
and $r_K:=1$ if $\Lambda_K=0$. Then $\exp$ is injective on $\fk(2r_K)$ and maps this ball
homeomorphically onto an open neighbourhood of $1$ in $K$. Since $\mathsf R\circ\exp$ is a
homomorphism of $\fk$ to $O(n)$, there is a constant $C_{\mathsf R}$ with
\begin{equation}\label{es:eq:crho}
\|\mathsf R(\exp\mathsf b)-\one\|\le C_{\mathsf R}\,|\mathsf b|\qquad\text{for all }\mathsf b\in\fk .
\end{equation}

The map $(\zeta,k)\mapsto\exp(\zeta)k$ of $\fm\times K$ to $G$ is smooth, its differential
is invertible at the points $(0,k)$, it commutes with the right multiplication by $K$,
and it is injective on $\fm(r_0)\times K$ (\cref{ls:lem:tube}\textup{(i)}). Since $K$ is
compact, there is $r_{\mathsf k}\in(0,\tfrac12r_0]$ such that it is a diffeomorphism of
$\fm(2r_{\mathsf k})\times K$ onto an open neighbourhood $U_K$ of $K$ in $G$.

\begin{definition}[The Component in the Isotropy Group]\label{es:def:component}
Every $g\in U_K$ has a unique decomposition $g=\exp(\zeta)k$ with $\zeta\in\fm(2r_{\mathsf k})$
and $k\in K$. The factor $k$ is the \emph{component of $g$ in $K$}, written $\mathsf k(g)$; this
defines the map $\mathsf k\colon U_K\to K$.
\end{definition}

This map is smooth, it is the identity on $K$, and
\begin{equation}\label{es:eq:lambdaK}
\mathsf k(gk)=\mathsf k(g)\,k,\qquad\mathsf k(kg)=k\,\mathsf k(g)\qquad\text{for }g\in U_K,\ k\in K,
\end{equation}
the second because $k\exp(\zeta)=\exp(\Ad(k)\zeta)k$ and $\Ad(k)$ preserves the norm of
$\fm$. For small $\zeta,\zeta'\in\fm$ the element $\exp(\zeta)\exp(\zeta')$ lies in
$U_K$, and its component in $K$ depends smoothly on $(\zeta,\zeta')$ and equals $1$ if
$\zeta=0$ or $\zeta'=0$. Hence there are numbers $r_c\in(0,r_{\mathsf k}]$ and $C_c<\infty$
and a smooth map $\varsigma$ of $\fm(2r_c)\times\fm(2r_c)$ to $\fk(r_K)$ with
\begin{equation}\label{es:eq:varsigma}
\mathsf k\bigl(\exp(\zeta)\exp(\zeta')\bigr)=\exp\bigl(\varsigma(\zeta,\zeta')\bigr),\qquad
|\varsigma(\zeta,\zeta')|\le C_c\,|\zeta|\,|\zeta'|
\end{equation}
for $|\zeta|,|\zeta'|\le r_c$; the bound holds because $\varsigma$ vanishes on the two
axes, so that $\varsigma(\zeta,\zeta')$ is the integral of the mixed second derivative
of $\varsigma$ over the rectangle with the corners $0$ and $(\zeta,\zeta')$. Let
$g=\exp(\zeta)k$ and $g'=\exp(\zeta')k'$ with $k,k'\in K$ and $|\zeta|,|\zeta'|\le r_c$.
Then $gg'$ is the product of $\exp(\zeta)\exp(\Ad(k)\zeta')$ and $kk'$, and therefore
\begin{equation}\label{es:eq:lambdaproduct}
gg'\in U_K\qquad\text{and}\qquad
\mathsf k(gg')=\exp\bigl(\varsigma(\zeta,\Ad(k)\zeta')\bigr)\,\mathsf k(g)\,\mathsf k(g') .
\end{equation}
So $\mathsf k$ is a homomorphism up to an element of the torus $K_0$ whose size is
quadratic in the distance from $K$; when $K$ is finite, $\fk=0$ and $\mathsf k$ is
multiplicative on the elements with $|\zeta|,|\zeta'|\le r_c$.

By \cref{ls:lem:tube}\textup{(ii)} and~\textup{(iii)} there is $\sigma_{\mathsf k}>0$ such
that every $g\in G$ with $\hat d(gz,z)<\sigma_{\mathsf k}$ is of the form $g=\exp(\zeta)k$
with $k\in K$ and $\zeta\in\fm(r_c)$, and then
\begin{equation}\label{es:eq:zetabound}
|\zeta|\le\frac{\hat d(gz,z)}{c_0},
\end{equation}
because $\exp(\zeta)z=gz$. We put
\begin{equation*}
t_0:=\tfrac13\sigma_{\mathsf k},\qquad\sE_0(t):=\frac{2C_c\,t^2}{c_0^2}\qquad(t>0).
\end{equation*}
The labels of the next subsection are defined for deck transformations that move $z_i$ by less
than $t_0$; the factor $\tfrac13$ keeps the limit elements that occur there, which move $z$ by
at most $2t_0$, inside $U_K$. For elements that move $z$ by less than $t\le t_0$,
\eqref{es:eq:zetabound} and \eqref{es:eq:varsigma} bound the vector $\varsigma(\zeta,\Ad(k)\zeta')$
in \eqref{es:eq:lambdaproduct} by $C_ct^2/c_0^2<\sE_0(t)$, and the labels inherit this bound
(\cref{es:lem:defect}).

\subsection{Labels}\label{es:sec:labels}\label{es:sec:labeldef}

We label each small deck transformation by the component in $K$ of a nearby limit
element. Convergence of these labels gives almost equivariance
(\cref{th:sec:stageequiv}) and the controlled multiplicative defect of
\cref{es:sec:defect}.

\begin{lemma}[Labels]\label{es:lem:labels}
For all large $i$ there are maps $\gamma\mapsto k_\gamma$ of $\Gamma_i(t_0)$ to $K$, with
$k_1=1$, that have the following property. If elements $\gamma_i\in\Gamma_i(t_0)$, for
$i$ in an infinite set of indices, converge to an element $g\in G$, then $g\in U_K$ and
$k_{\gamma_i}\to\mathsf k(g)$.
\end{lemma}
\begin{proof}
Use the group approximations $\alpha_i$ of \cref{m:fr:sec:conv}, normalized by
$\alpha_i(1)=1$, and put \mbox{$k_\gamma:=\mathsf k(\alpha_i(\gamma))$}.
These values are defined on $\Gamma_i(t_0)$ for all large $i$: otherwise a
bounded-displacement sequence $\gamma_i$ with $\alpha_i(\gamma_i)\notin U_K$ has a
subsequence converging to $g$ with \mbox{$\hat d(gz,z)\le t_0<\sigma_{\mathsf k}$}.
Then $g\in U_K$ and $\alpha_i(\gamma_i)\to g$ by
\cref{m:fr:lem:groupconv}\textup{(ii)}, a contradiction. That lemma applies also at the points $z_i$ and $z$, whose distances from the base points of the chart are bounded independently of $i$. Finally $k_1=1$, and if $\gamma_i\eGHto g$, the same lemma and continuity of
$\mathsf k$ give $k_{\gamma_i}\to\mathsf k(g)$.
\end{proof}

We fix such maps and call $k_\gamma$ the \emph{deck label} of $\gamma$, or its label. They are not
canonical; only the property of \cref{es:lem:labels} is used. For
$\gamma\in\Gamma_i(t_0)$ we write
\begin{equation*}
\Ad_\gamma:=\Ad(k_\gamma)|_{\fk},
\end{equation*}
an isometry of $\fk$. Since $\fk$ is abelian, $\Ad(k)$ acts trivially on $\fk$ for
$k\in K_0$, so $\Ad_\gamma$ depends only on the class of $k_\gamma$ in the finite group
$K/K_0$.

\subsection{Almost Equivariance on the Charts \texorpdfstring{$Z_i$}{Z\_i}}\label{th:sec:stageequiv}

The components of the transplanted chart in $V$ inherit, for the labels of
\cref{es:lem:labels}, the almost equivariance of \cref{th:lem:almostequiv}: first uniformly
and then, by a difference quotient, with a small Lipschitz constant. Put
\begin{equation*}
\rho_3:=\tfrac12\min\{\rho_2,\rho_s,s_1\},\qquad
\sE_1(t):=\frac{16\,C_5\,(2s_0)^{\eta}}{c_3}\;t^{1-\eta}\qquad(t>0).
\end{equation*}

\begin{lemma}[Labels in the Transplanted Chart]\label{th:lem:stageequiv}\label{es:lem:labelchart}
Let $0<t\le\min\{\rho_3,\tfrac12t_0\}$. For all large $i$ and every
$\gamma\in\Gamma_i$ for which $\gamma B_t(z_i)$ meets $B_t(z_i)$, the label
$k_\gamma$ of \cref{es:lem:labels} is defined and has the following properties.
\begin{enumerate}[label=\textup{(\roman*)}]
\item We have $|P_Vh_i(\gamma p)-\mathsf R(k_\gamma)P_Vh_i(p)|\le2C_5t^2$ for all
$p\in B_t(z_i)$ with $\gamma p\in B_t(z_i)$.
\item If $\hat d(\gamma z_i,z_i)\le\tfrac12t$, then the map
$P_Vh_i\circ\gamma-\mathsf R(k_\gamma)P_Vh_i$ is $(6\eps_*+\sE_1(t))$-Lipschitz with respect to
$h_i$ on $B^i$.
\end{enumerate}
\end{lemma}
\begin{proof}
(i) If $\gamma B_t(z_i)$ meets $B_t(z_i)$, then
$\hat d(\gamma z_i,z_i)<2t\le t_0$, so its label is defined. Suppose the bound
fails for infinitely many $i$, for elements $\gamma_i$ and points
$p_i\in B_t(z_i)$ with $\gamma_ip_i\in B_t(z_i)$ and
\begin{equation*}
|P_Vh_i(\gamma_ip_i)-\mathsf R(k_{\gamma_i})\,P_Vh_i(p_i)|>2C_5t^2 .
\end{equation*}
Their displacements are bounded by $2t$. After taking a subsequence,
$\gamma_i\eGHto g\in G$ and $p_i\to p$, with $\hat d(p,z)\le t$ and
$\hat d(gp,z)\le t$. By \cref{es:lem:labels}, $g\in U_K$ and
$k_{\gamma_i}\to k:=\mathsf k(g)$. Both $p$ and $gp$ lie in $B$ and in
$B_{\rho_2}(z)$, since $t\le\rho_3\le\tfrac12\min\{s_1,\rho_2\}$.
The decomposition $g=\exp(\zeta)k$ defining $\mathsf k(g)$ has
$\zeta\in\fm(r_c)\subset\fm(r_0/2)$, so it is the unique decomposition of
\cref{th:lem:almostequiv}, and its factor in $K$ is this $k$.
By \cref{th:lem:transplant}\textup{(iii)}, the left-hand side converges to
$|P_Vh(gp)-\mathsf R(k)P_Vh(p)|\le C_5t^2$, a contradiction.

(ii) Let $R_\gamma$ be as in \cref{th:lem:transplant}\textup{(iv)}, and put
$\ell_t:=c_3(t/(2s_0))^{\eta}\,\tfrac12t$. By \cref{th:lem:transplant}\textup{(ii)} the image
of $B_{t/2}(z_i)$ under $h_i$ contains the open ball of radius $\ell_t$ about $0$. Let
$\nu$ be a vector of length $\tfrac12\ell_t$, and let $q\in B_{t/2}(z_i)$ be the point with
$h_i(q)=\nu$. The points $z_i$, $q$, $\gamma z_i$ and $\gamma q$ lie in $B_t(z_i)$. We
write the difference of $P_Vh_i\circ\gamma$ at $q$ and at $z_i$ in two ways, with
\cref{th:lem:transplant}\textup{(iv)} and with (i):
\begin{equation*}
\bigl|\bigl(P_VR_\gamma-\mathsf R(k_\gamma)P_V\bigr)\nu\bigr|\le3\eps_*|\nu|+4C_5t^2 .
\end{equation*}
Hence the operator norm of $P_VR_\gamma-\mathsf R(k_\gamma)P_V$ is at most
$3\eps_*+8C_5t^2/\ell_t=3\eps_*+\sE_1(t)$. For $p,q\in B^i$ the difference of the
values of $P_Vh_i\circ\gamma-\mathsf R(k_\gamma)P_Vh_i$ at $p$ and at $q$ is the sum of
$(P_VR_\gamma-\mathsf R(k_\gamma)P_V)(h_i(p)-h_i(q))$ and of the difference of the values of
$P_V(h_i\circ\gamma-R_\gamma h_i)$, so its norm is at most
$(6\eps_*+\sE_1(t))\,|h_i(p)-h_i(q)|$.
\end{proof}

\subsection{Averages over the Deck Groups}\label{th:sec:averages}

A function on $Z_i$ descends to $X_i$ when it is invariant under $\Gamma_i$.
We construct normalized cut-off weights for averaging over $\Gamma_i$ and bound
their total variation in the transplanted chart. These estimates are used in \cref{ch:localbundle} to correct the labels and to construct the equivariant coordinate.

\subsubsection{The Scale}\label{th:sec:scale}\label{es:sec:scale}

All constructions take place in a coordinate ball of radius $\ell$, and their errors are
functions of $\ell$ that tend to zero. Put $\beta:=1/(1+\eta)$, and for $\ell>0$ let
\begin{gather*}
d_+(\ell):=s_0\Bigl(\frac{3\ell}{c_3s_0}\Bigr)^{\beta},\qquad t_\ell:=4\,d_+(\ell),\qquad
\sE_2(\ell):=2C_5t_\ell^2,\\
C_\#(\ell):=C_D^4\Bigl(\frac{6C_6d_+(\ell)}{\ell}\Bigr)^{2\log_2C_D},\qquad
\sE_3(\ell):=\sE_1(t_\ell)+\frac{5\,C_\#(\ell)\,\sE_2(\ell)}{\ell},
\end{gather*}
where $C_D=C_D(N+1)$ is the doubling constant of \cref{th:sec:constants} and $\sE_1$ is the
function of \cref{th:sec:stageequiv}. With the numbers of \cref{es:sec:lambda} we also put
\begin{equation*}
\sE_4(\ell):=C_{\mathsf R}\,\sE_0(t_\ell)\Bigl(1+10\,C_\#(\ell)+\frac{5\,C_\#(\ell)\,\sE_2(\ell)}{\ell}\Bigr).
\end{equation*}
Here $\sE_0(t)=2C_ct^2/c_0^2$ is the function of \cref{es:sec:lambda}: it bounds the failure of
the component in $K$ to be multiplicative on elements that move $z$ by less than $t$
\eqref{es:eq:lambdaproduct}, and hence the defects of the labels (\cref{es:lem:defect}); at
$t=t_\ell$ it bounds all the defects that occur at the scale $\ell$.
The functions $\sE_2$, $\sE_3$ and $\sE_4$ are the errors of the equivariant coordinate of
\cref{es:sec:twisted}, which averages the maps $\mathsf R(k_{\gamma^{-1}})^{-1}P_Vh_i\circ\gamma^{-1}$
after a correction of the labels by a cocycle. These maps differ from $P_Vh_i$ by at most
$\sE_2(\ell)$, their weighted average differs from $P_Vh_i$ by a map that is
$(6\eps_*+\sE_3(\ell))$-Lipschitz with respect to $h_i$, and the correction adds at most
$\sE_4(\ell)$ to this constant (\cref{th:lem:weights,es:prop:coordinate}).

\begin{lemma}[The Errors Are Small]\label{th:lem:errors}
As $\ell\to0$ the numbers $\sE_2(\ell)/\ell$, $C_\#(\ell)\sE_2(\ell)/\ell$, $\sE_3(\ell)$ and $\sE_4(\ell)$ tend to $0$.
\end{lemma}
\begin{proof}
Up to constant factors, $\sE_2(\ell)/\ell$ is $\ell^{2\beta-1}$, and
$C_\#(\ell)\sE_2(\ell)/\ell$ is the power of $\ell$ with the exponent
\begin{equation*}
2\beta-1+2(\beta-1)\log_2C_D=\beta\,(1-\eta-2\eta\log_2C_D),
\end{equation*}
since $2\beta-1=\beta(1-\eta)$ and $\beta-1=-\beta\eta$ for $\beta=1/(1+\eta)$. Both exponents are positive, the second because $\eta\,(4+16\log_2C_D)\le1$. Finally $\sE_1(t_\ell)$ is a
constant times $t_\ell^{1-\eta}$. Up to a constant factor $\sE_0(t_\ell)$ is $\sE_2(\ell)$, so $\sE_4(\ell)$
tends to $0$ as well.
\end{proof}

We fix $\ell_*>0$ such that the following holds for all $0<\ell\le\ell_*$: we have $t_\ell\le\rho_3$ and $30\,d_+(\ell)\le s_0$, the doubling inequality
$\hat m_i(B_{2t})\le C_D\,\hat m_i(B_t)$ holds for the balls of radius $t\le30\,d_+(\ell)$
in $Z_i$, the ball of radius $d_+(\ell)$ about $z$ lies in $B_{r/2}(z)$, where $r$ is the
frame scale, and
\begin{equation}\label{th:eq:ellsmall}
3\ell<c_3\Bigl(\frac{s_1}{2s_0}\Bigr)^{\eta}\frac{s_1}2,\qquad\ell\le\tfrac12c_zr_3,\qquad
\sE_3(\ell)\le\eps_*,\qquad 10 \sE_2(\ell)\le r_\ell:=\frac{\ell}{800},
\end{equation}
and, for the cocycle construction of \cref{ch:localbundle},
\begin{equation}\label{es:eq:ellone}
\begin{gathered}
2t_\ell\le t_0,\qquad\ell\le c_z\,r_c,\qquad
\sE_0(t_\ell)\le\min\Bigl\{\tfrac18r_K,\ \frac1{8000\,C_{\mathsf R}}\Bigr\},\\
\sE_3(\ell)+\sE_4(\ell)\le\eps_* .
\end{gathered}
\end{equation}
The number $r_\ell$ is the unit in which the sizes of the tube at $y$ (\cref{es:def:tube}) are measured.
For $0<t\le3$ we put
\begin{equation*}
\cC^i_t:=\{p\in B^i:|h_i(p)|<t\ell\},\qquad \cC_t:=\{p\in B:|h(p)|<t\ell\}.
\end{equation*}
By \cref{th:lem:transplant}\textup{(i)} and~\textup{(ii)} and the first inequality of
\eqref{th:eq:ellsmall}, the map $h_i$ is a homeomorphism of $\cC^i_t$ onto the open ball of
radius $t\ell$ about $0$ in $\R^n$, and $h$ is one of $\cC_t$ onto the same ball. By
\eqref{th:eq:stageholder},
\begin{equation}\label{th:eq:metricsize}
\cC^i_3\subset B_{d_+(\ell)}(z_i),\qquad \cC_3\subset B_{d_+(\ell)}(z),
\end{equation}
and these balls have radius at most $\tfrac18s_1$ and at most $\tfrac14\rho_s$. If points
$p_i\in B^i$ converge to a point $p\in \cC_t$, then $p_i\in \cC^i_t$ for large $i$, by
\cref{th:lem:transplant}\textup{(iii)}. From now on $\ell\in(0,\ell_*]$ is fixed, and
``large $i$'' may depend on it.

\subsubsection{Weights}\label{th:sec:weights}

We record the normalization and weighted-sum estimate used throughout the construction. The
weights $\varpi_\gamma$ constructed afterwards, from a cut-off function in the transplanted chart,
are those with which the labels are corrected in \cref{es:sec:cocycle} and the equivariant
coordinate is formed in \cref{es:sec:twisted}. The weighted-sum estimate
\eqref{th:eq:weightedaverage} uses neither the group action nor values of the maps
$\mathsf e_\gamma$ other than those at the two points, so it applies to locally defined maps.

\begin{lemma}[Normalized Weights and Averages]\label{th:lem:normalizeweights}
Let a group $\Gamma$ act on a space, let $U$ be open, and let
$w_\gamma\colon U\to[0,\infty)$ be continuous functions whose supports are locally
finite and whose sum $\bar w$ is positive. Then $\varpi_\gamma:=w_\gamma/\bar w$ are continuous,
their supports are locally finite, their sum is $1$, and
\begin{equation}\label{th:eq:normalizeweights}
\sum_\gamma|\varpi_\gamma(p)-\varpi_\gamma(q)|
\le\frac2{\bar w(p)}\sum_\gamma|w_\gamma(p)-w_\gamma(q)|\qquad(p,q\in U).
\end{equation}
If $w_\gamma(\gamma'p)=w_{\gamma'^{-1}\gamma}(p)$ whenever $p,\gamma'p\in U$, then $\varpi_\gamma(\gamma'p)=\varpi_{\gamma'^{-1}\gamma}(p)$ there as well.
For any such normalized family of weights and maps $\mathsf e_\gamma$ to a normed vector
space, put $\bar{\mathsf e}:=\sum_\gamma\varpi_\gamma \mathsf e_\gamma$.
For a metric $d$ on $U$ and $p,q\in U$, suppose the maps with a non-zero weight at either point satisfy
$|\mathsf e_\gamma(q)-\mathsf e_\gamma(p)|\le Ld(p,q)$ and $|\mathsf e_\gamma(p)-\nu|\le M$ for some vector $\nu$.
Then
\begin{equation}\label{th:eq:weightedaverage}
|\bar{\mathsf e}(q)-\bar{\mathsf e}(p)|\le Ld(p,q)+M\sum_\gamma|\varpi_\gamma(q)-\varpi_\gamma(p)|.
\end{equation}
\end{lemma}
\begin{proof}
Put $\Delta:=\sum_\gamma|w_\gamma(p)-w_\gamma(q)|$. The left-hand side of \eqref{th:eq:normalizeweights} is at most
\begin{equation*}
\frac{\Delta}{\bar w(p)}+\frac{|\bar w(q)-\bar w(p)|}{\bar w(p)}\le\frac{2\Delta}{\bar w(p)}.
\end{equation*}
Continuity follows from local finiteness and positivity, and the covariance follows
by reindexing the sum for $\bar w$. For the last assertion use
$\sum_\gamma(\varpi_\gamma(q)-\varpi_\gamma(p))=0$:
\[
\bar{\mathsf e}(q)-\bar{\mathsf e}(p)=\sum_\gamma\varpi_\gamma(q)(\mathsf e_\gamma(q)-\mathsf e_\gamma(p))
+\sum_\gamma(\varpi_\gamma(q)-\varpi_\gamma(p))(\mathsf e_\gamma(p)-\nu).
\]
\end{proof}

Let $\chi_0\colon[0,\infty)\to[0,1]$ be
the function that equals $1$ on $[0,2]$, equals $3-t$ on $[2,3]$ and vanishes on
$[3,\infty)$, and let $\chi\colon Z_i\to[0,1]$ be the function that equals
$\chi_0(|h_i|/\ell)$ on $B^i$ and vanishes elsewhere. It is continuous, its support lies
in $B_{d_+(\ell)}(z_i)$, and it equals $1$ on $\cC^i_2$. For $p\in Z_i$ the sum
$\bar\chi(p):=\sum_{\gamma\in\Gamma_i}\chi(\gamma^{-1}p)$ is locally finite, by proper
discontinuity and the compact support of $\chi$, and it is positive on the open invariant set $\Gamma_i\,\cC^i_3$. For $p\in\Gamma_i\,\cC^i_3$ and $\gamma,\gamma'\in\Gamma_i$ let
\begin{equation}\label{es:eq:weights}
\varpi_\gamma(p):=\frac{\chi(\gamma^{-1}p)}{\bar\chi(p)},\qquad\sum_\gamma\varpi_\gamma(p)=1,\qquad
\varpi_\gamma(\gamma'p)=\varpi_{\gamma'^{-1}\gamma}(p);
\end{equation}
the weights are continuous, and the two identities hold by \cref{th:lem:normalizeweights}.
If a point $p$ and its image $\gamma p$ lie in $\cC^i_3$, then $\gamma\in\Gamma_i(\tfrac12t_\ell)$,
because both points lie in $B_{d_+(\ell)}(z_i)$; in particular this holds if $p\in \cC^i_3$ and
$\varpi_{\gamma^{-1}}(p)>0$.
For a continuous map $u$ of the closure of $\cC^i_3$ to a Euclidean space let
\begin{equation*}
\hat u(p):=\sum_\gamma\varpi_\gamma(p)\,u(\gamma^{-1}p)\qquad(p\in\Gamma_i\,\cC^i_3),
\end{equation*}
the sum over the elements with $\chi(\gamma^{-1}p)>0$. The map $\hat u$ is continuous and invariant under $\Gamma_i$.

\begin{lemma}[Weights]\label{th:lem:weights}
For all large $i$ the following holds. Let $p,q\in \cC^i_1$, and let $\gamma\in\Gamma_i$ be
an element for which $\chi(\gamma^{-1}p)$ or $\chi(\gamma^{-1}q)$ is positive.
\begin{enumerate}[label=\textup{(\roman*)}]
\item The points $\gamma^{-1}p$ and $\gamma^{-1}q$ lie in $B_{t_\ell}(z_i)$, and
$\hat d(\gamma^{-1}z_i,z_i)\le\tfrac12t_\ell$. With the label $k:=k_{\gamma^{-1}}$ of \cref{es:lem:labels}, the map
\begin{equation}\label{es:eq:epsgamma}
\tilde\eps_\gamma:=\mathsf R(k)^{-1}\,P_Vh_i\circ\gamma^{-1}-P_Vh_i
\end{equation}
satisfies
\begin{equation}\label{es:eq:epsbounds}
\max\{|\tilde\eps_\gamma(p)|,|\tilde\eps_\gamma(q)|\}\le\sE_2(\ell),\qquad
|\tilde\eps_\gamma(p)-\tilde\eps_\gamma(q)|\le\bigl(6\eps_*+\sE_1(t_\ell)\bigr)\,|h_i(p)-h_i(q)| .
\end{equation}
\item We have
\begin{equation}\label{es:eq:thetalip}
\sum_\gamma|\varpi_\gamma(p)-\varpi_\gamma(q)|\le\frac{5\,C_\#(\ell)}{\ell}\,|h_i(p)-h_i(q)| .
\end{equation}
\end{enumerate}
\end{lemma}
\begin{proof}
(i) Suppose that $\chi(\gamma^{-1}p)>0$; the other case is symmetric. By
\eqref{th:eq:metricsize} the points $p$, $q$ and $\gamma^{-1}p$ lie in $B_{d_+(\ell)}(z_i)$.
So $\hat d(\gamma^{-1}z_i,z_i)\le\hat d(\gamma^{-1}z_i,\gamma^{-1}p)+\hat d(\gamma^{-1}p,z_i)
<2d_+(\ell)=\tfrac12t_\ell$, and $\hat d(\gamma^{-1}q,z_i)<3d_+(\ell)<t_\ell$. Thus $p$ and
$q$ lie in $B_{t_\ell}(z_i)$ together with their images under $\gamma^{-1}$, and
\cref{th:lem:stageequiv} applies to $\gamma^{-1}$ with
$t=t_\ell\le\min\{\rho_3,t_0/2\}$. Part (i) of that
lemma gives the first bound, since $\mathsf R(k)$ is orthogonal and $2C_5t_\ell^2=\sE_2(\ell)$,
and part (ii) gives the second.

(ii) For $p'\in Z_i$ and $t>0$ let $N_{p'}(t)$ be the number of the elements
$\gamma\in\Gamma_i$ with $\hat d(\gamma p',p')<t$. For $p'\in B^i$ and
$0<t\le R$ with $R+t/2\le30d_+(\ell)$, a direct overlap estimate gives
\begin{equation}\label{th:eq:countdoubling}
N_{p'}(R)\le
\frac{\hat m_i(B_{R+t/2}(p'))}{\hat m_i(B_{t/2}(p'))}\,N_{p'}(t)
\le C_D^{\lceil\log_2(1+2R/t)\rceil}\,N_{p'}(t).
\end{equation}
This is the packing argument of the proof of \cref{m:fr:lem:counting}\textup{(i)}: the sets
$\gamma B_{t/2}(p')$ with $\hat d(\gamma p',p')<R$ lie in
$B_{R+t/2}(p')$ and have multiplicity at most $N_{p'}(t)$: two that contain the
same point have $\gamma'^{-1}\gamma$ moving $p'$ by less than $t$.
Integrating their indicators proves the first inequality, since the action preserves
$\hat m_i$. The second uses doubling at radii less than $R+t/2$.

Let $p,q\in \cC^i_1$. If $\hat d(\gamma^{-1}p,p)<\ell/C_6$, then $\gamma^{-1}p$ lies in
$B^i$ and $|h_i(\gamma^{-1}p)|<2\ell$ by \eqref{th:eq:stageholder}, so $\chi(\gamma^{-1}p)=1$;
hence $\bar\chi(p)\ge N_p(\ell/C_6)$. If $\chi(\gamma^{-1}p)>0$, then
$\hat d(\gamma^{-1}p,p)<2d_+(\ell)$, and if $\chi(\gamma^{-1}q)>0$, then
$\hat d(\gamma^{-1}q,q)<2d_+(\ell)$ and so $\hat d(\gamma^{-1}p,p)<6d_+(\ell)$. So at most
$2N_p(6d_+(\ell))$ elements $\gamma$ have $\chi(\gamma^{-1}p)>0$ or
$\chi(\gamma^{-1}q)>0$. By \eqref{th:eq:stageholder},
$Q:=6C_6d_+(\ell)/\ell\ge1$. Apply
\eqref{th:eq:countdoubling} with $R=6d_+(\ell)$ and $t=\ell/C_6$.
Here $R+t/2\le9d_+(\ell)$ and
$\lceil\log_2(1+2Q)\rceil\le4+2\log_2Q$, so
\begin{equation*}
N_p(6d_+(\ell))\le C_\#(\ell)\,N_p(\ell/C_6)
\le C_\#(\ell)\,\bar\chi(p) .
\end{equation*}
For such $\gamma$ the function $\chi_0$ is Lipschitz with the constant $1$, and by (i)
and \cref{th:lem:transplant}\textup{(iv)}
\begin{equation*}
|\chi(\gamma^{-1}p)-\chi(\gamma^{-1}q)|\le\ell^{-1}\,|h_i(\gamma^{-1}p)-h_i(\gamma^{-1}q)|
\le(1+3\eps_*)\,\ell^{-1}\,|h_i(p)-h_i(q)| .
\end{equation*}
Consequently \eqref{th:eq:normalizeweights} gives
\begin{equation*}
\sum_\gamma|\varpi_\gamma(p)-\varpi_\gamma(q)|
\le4(1+3\eps_*)\,C_\#(\ell)\,\ell^{-1}\,|h_i(p)-h_i(q)|
\le5C_\#(\ell)\,\ell^{-1}\,|h_i(p)-h_i(q)|,
\end{equation*}
because $\eps_*\le1/200$.
\end{proof}

\section{Cocycles and Admissible Sections}\label{ch:localbundle}

The labels give a cocycle with values in the isotropy group, which depends on the point
(\cref{es:sec:defect,es:sec:cocycle}), and an equivariant coordinate on a slice in the sense of \cref{ls:def:slice}
(\cref{es:sec:coordinate}). With them we fix the domain of the coordinate and describe all
nearby points by linear slices of the representation (\cref{es:sec:model}). Sections of the
frame spaces with exactly this cocycle, the admissible sections, exist and are unique up to
small corrections (\cref{es:sec:frames}). A triviality theorem for maps with product charts
over a ball of a representation completes the local tools (\cref{es:sec:product});
\cref{ch:globalbundle} glues with them.

Throughout this section we are in the standing setting of \cref{th:sec:standing}, with the
objects fixed there and in \cref{es:sec:lambda,es:sec:labels,th:sec:averages}: one centre $a$,
the point $y\in Y$ with its frame $(z,e)$ and the isotropy groups $K\cong I$, the harmonic
chart $h$ at the scale $s_0$ and the transplanted charts $h_i$ on the charts $Z_i$
(\cref{th:def:transplant}), the component $\mathsf k$ (\cref{es:def:component}), the numbers $t_0$ and $r_K$ and the
functions $\sE_j$, the labels $k_\gamma$, the fixed scale $\ell\in(0,\ell_*]$ with the sets
$\cC^i_t$ and the weights $\varpi_\gamma$. All statements refer to the subsequence of
\cref{th:lem:stagecharts}; the averages below converge without further extraction.

\subsection{The Defect}\label{es:sec:defect}

The labels of a product differ from the product of the labels by an element of the
torus. We bound it and record the identity that it satisfies, which is that of a
cocycle of degree two with twisted coefficients.

\begin{lemma}[The Defect]\label{es:lem:defect}
Let $0<t\le\tfrac12t_0$ with $\sE_0(t)\le r_K$. For all large $i$ and all
$\gamma,\gamma'\in\Gamma_i(t)$ there is a unique $\omega_i(\gamma,\gamma')\in\fk$ with
\begin{equation*}
k_{\gamma\gamma'}=\exp\bigl(\omega_i(\gamma,\gamma')\bigr)\,k_\gamma\,k_{\gamma'}\qquad\text{and}
\qquad|\omega_i(\gamma,\gamma')|\le\sE_0(t).
\end{equation*}
\end{lemma}
\begin{proof}
The product lies in $\Gamma_i(2t)\subset\Gamma_i(t_0)$, so its label is defined, and the
uniqueness holds because $\exp$ is injective on $\fk(2r_K)$. Suppose that the existence
fails for infinitely many $i$, for elements $\gamma_i$ and $\gamma'_i$. After passing to
a subsequence they converge to elements $g$ and $g'$ of $G$ that move $z$ by at most
$t$, and $\gamma_i\gamma'_i\eGHto gg'$. Since $t<\sigma_{\mathsf k}$, we have $g=\exp(\zeta)k$ and
$g'=\exp(\zeta')k'$ with $k,k'\in K$ and $\zeta,\zeta'\in\fm(r_c)$, and
$|\zeta|,|\zeta'|\le t/c_0$ by \eqref{es:eq:zetabound}. By
\cref{es:lem:labels} and \eqref{es:eq:lambdaproduct} the elements
$k_{\gamma_i\gamma'_i}(k_{\gamma_i}k_{\gamma'_i})^{-1}$ converge to
\begin{equation*}
\mathsf k(gg')\,\bigl(\mathsf k(g)\mathsf k(g')\bigr)^{-1}=\exp\bigl(\varsigma(\zeta,\Ad(k)\zeta')\bigr),
\end{equation*}
and $|\varsigma(\zeta,\Ad(k)\zeta')|\le C_ct^2/c_0^2<\sE_0(t)$ by \eqref{es:eq:varsigma}.
The image of $\fk(\sE_0(t))$ under $\exp$ is an open subset of $K$ that contains this
limit, so it contains $k_{\gamma_i\gamma'_i}(k_{\gamma_i}k_{\gamma'_i})^{-1}$ for large $i$
of the subsequence, a contradiction.
\end{proof}

The defect is a two-cochain with values in $\fk$, defined only on pairs of small deck
transformations. Computing the label of a triple product in two ways gives its cocycle
identity, in which the labels act on $\fk$ through the adjoint action.

\begin{lemma}[The Identity of the Defect]\label{es:lem:cocycleid}
Let $0<t\le\tfrac14t_0$ with $\sE_0(2t)\le\tfrac12r_K$. For all large $i$ and all
$\gamma,\gamma',\gamma''\in\Gamma_i(t)$ we have $\Ad_{\gamma\gamma'}=\Ad_\gamma\Ad_{\gamma'}$,
$\Ad_1=\mathrm{id}$, $\omega_i(\gamma,1)=\omega_i(1,\gamma)=0$ and
\begin{equation}\label{es:eq:twococycle}
\omega_i(\gamma\gamma',\gamma'')+\omega_i(\gamma,\gamma')
=\omega_i(\gamma,\gamma'\gamma'')+\Ad_\gamma\,\omega_i(\gamma',\gamma''),
\end{equation}
where the defects are those of \cref{es:lem:defect} for the number $2t$.
\end{lemma}
\begin{proof}
All products of two of the three elements lie in $\Gamma_i(2t)$, and
$\gamma\gamma'\gamma''\in\Gamma_i(3t)\subset\Gamma_i(t_0)$, so the labels and the defects in
the statement are defined, and the defects have norm at most $\sE_0(2t)$. Since
$\exp(\omega_i(\gamma,\gamma'))$ lies in $K_0$, which acts trivially on $\fk$, the first
identity holds, and $\Ad_1=\mathrm{id}$ because $k_1=1$. From $k_1=1$ and the
uniqueness in \cref{es:lem:defect} we get $\omega_i(\gamma,1)=\omega_i(1,\gamma)=0$. We
compute the label of $\gamma\gamma'\gamma''$ in two ways. On the one hand it is
$\exp(\omega_i(\gamma\gamma',\gamma''))\,k_{\gamma\gamma'}k_{\gamma''}
=\exp(\omega_i(\gamma\gamma',\gamma'')+\omega_i(\gamma,\gamma'))\,k_\gamma k_{\gamma'}k_{\gamma''}$.
On the other hand it is
\begin{equation*}
\exp\bigl(\omega_i(\gamma,\gamma'\gamma'')\bigr)\,k_\gamma\,
\exp\bigl(\omega_i(\gamma',\gamma'')\bigr)\,k_{\gamma'}k_{\gamma''}
=\exp\bigl(\omega_i(\gamma,\gamma'\gamma'')+\Ad_\gamma\,\omega_i(\gamma',\gamma'')\bigr)\,
k_\gamma k_{\gamma'}k_{\gamma''} .
\end{equation*}
The two exponents have norm at most $2\sE_0(2t)\le r_K$, and $\exp$ is injective on
$\fk(2r_K)$.
\end{proof}

When $K$ is finite, $\fk=0$, all defects vanish, and \cref{es:lem:defect} says that the
labels are multiplicative on $\Gamma_i(t)$.

\subsection{The Cocycle}\label{es:sec:cocycle}\label{es:sec:correction}

We correct the labels by averaging their defects with the weights of
\cref{th:sec:weights}. The identity \eqref{es:eq:twococycle} makes the correction
exact; compare the degree-two averaging of \cite[Proposition~1]{Cra03} and
\cite{GKR74actions,CFG92}. For $\gamma\in\Gamma_i(\tfrac12t_\ell)$ and
$p\in \cC^i_3$ let
\begin{equation}\label{es:eq:adef}
\sa_i(\gamma,p):=\sum_\mu\varpi_\mu(p)\,\omega_i(\gamma,\mu)\in\fk,\qquad
\varkappa_i(\gamma,p):=\exp\bigl(\sa_i(\gamma,p)\bigr)\,k_\gamma\in K,
\end{equation}
the sum over the elements $\mu\in\Gamma_i$ with $\varpi_\mu(p)>0$. For these elements
$\mu^{-1}p$ lies in $\cC^i_3$, so $\mu\in\Gamma_i(\tfrac12t_\ell)$ and
$\gamma\mu\in\Gamma_i(t_\ell)$. The defect in \eqref{es:eq:adef} is
that of \cref{es:lem:defect} for $t=t_\ell$; by \eqref{es:eq:ellone} this lemma and
\cref{es:lem:cocycleid}, for $t=\tfrac12t_\ell$, apply for large $i$. We recall that
$\gamma\in\Gamma_i(\tfrac12t_\ell)$ whenever $p$ and $\gamma p$ lie in $\cC^i_3$.

\begin{proposition}[The Cocycle]\label{es:prop:cocycle}
For all large $i$ the following holds.
\begin{enumerate}[label=\textup{(\roman*)}]
\item For every $\gamma\in\Gamma_i(\tfrac12t_\ell)$ the map $\sa_i(\gamma,\cdot)$ is continuous
on $\cC^i_3$, we have $|\sa_i(\gamma,p)|\le\sE_0(t_\ell)$, and
\begin{equation*}
|\sa_i(\gamma,p)-\sa_i(\gamma,q)|\le\frac{5\,\sE_0(t_\ell)\,C_\#(\ell)}{\ell}\,|h_i(p)-h_i(q)|
\qquad\text{for all }p,q\in \cC^i_1 .
\end{equation*}
\item We have $\varkappa_i(1,p)=1$ for $p\in \cC^i_3$. If $\gamma,\gamma'\in\Gamma_i$ and the points
$p$, $\gamma'p$ and $\gamma\gamma'p$ lie in $\cC^i_3$, then
\begin{equation}\label{es:eq:cocycle}
\varkappa_i(\gamma\gamma',p)=\varkappa_i(\gamma,\gamma'p)\,\varkappa_i(\gamma',p).
\end{equation}
\item For every $\gamma\in\Gamma_i(\tfrac12t_\ell)$ and $p\in \cC^i_3$ the restriction of
$\Ad(\varkappa_i(\gamma,p))$ to $\fk$ is $\Ad_\gamma$.
\end{enumerate}
\end{proposition}
\begin{proof}
(i) Continuity follows from local finiteness. Positivity and normalization give
$|\sa_i|\le\sE_0(t_\ell)$. Apply \eqref{th:eq:weightedaverage} to the constant vectors $\mathsf e_\mu=\omega_i(\gamma,\mu)$, with $L=0$, $\nu=0$ and $M=\sE_0(t_\ell)$, and use
\eqref{es:eq:thetalip} for the Lipschitz estimate.

(ii) Since $\omega_i(1,\mu)=0$ and $k_1=1$, we have $\varkappa_i(1,p)=1$. Let $p$,
$\gamma'p$ and $\gamma\gamma'p$ lie in $\cC^i_3$; then $\gamma$, $\gamma'$ and $\gamma\gamma'$ lie in
$\Gamma_i(\tfrac12t_\ell)$. If $\varpi_\mu(p)>0$, then $\mu^{-1}p\in \cC^i_3$,
so $\mu$ and $\gamma'\mu$ also lie in $\Gamma_i(\tfrac12t_\ell)$:
the latter carries $\mu^{-1}p$ to $\gamma'p$. Reindexing the first sum by
\eqref{es:eq:weights} and applying \eqref{es:eq:twococycle} to
$(\gamma,\gamma',\mu)$ therefore gives
\begin{align*}
\sa_i(\gamma,\gamma'p)+\Ad_\gamma\sa_i(\gamma',p)
&=\sum_\mu\varpi_\mu(p)
 \bigl[\omega_i(\gamma,\gamma'\mu)+\Ad_\gamma\omega_i(\gamma',\mu)\bigr]\\
&=\sum_\mu\varpi_\mu(p)
 \bigl[\omega_i(\gamma\gamma',\mu)+\omega_i(\gamma,\gamma')\bigr].
\end{align*}
Since the weights sum to one, this is
\begin{equation*}
\sa_i(\gamma,\gamma'p)+\Ad_\gamma\,\sa_i(\gamma',p)=\sa_i(\gamma\gamma',p)+\omega_i(\gamma,\gamma').
\end{equation*}
Since $\fk$ is abelian and $k_\gamma\exp(\mathsf b)=\exp(\Ad_\gamma\mathsf b)\,k_\gamma$,
\begin{equation*}
\varkappa_i(\gamma,\gamma'p)\,\varkappa_i(\gamma',p)
=\exp\bigl(\sa_i(\gamma,\gamma'p)+\Ad_\gamma\,\sa_i(\gamma',p)\bigr)\,k_\gamma k_{\gamma'}
=\exp\bigl(\sa_i(\gamma\gamma',p)\bigr)\,\exp\bigl(\omega_i(\gamma,\gamma')\bigr)\,k_\gamma k_{\gamma'},
\end{equation*}
and the right-hand side is $\varkappa_i(\gamma\gamma',p)$ by the definition of the defect.

(iii) The element $\exp(\sa_i(\gamma,p))$ lies in $K_0$, which acts trivially on $\fk$.
\end{proof}

So $\varkappa_i$ is a cocycle on the part of the action of $\Gamma_i$ that is seen in the
coordinate ball $\cC^i_3$, with values in $K$; it differs from the labels by elements of
the torus of size $\sE_0(t_\ell)$, and it varies slowly: by (i) and \eqref{es:eq:ellone} a
change of the point by the size $\ell$ of the ball, measured in the coordinates,
changes $\sa_i$ by at most $\eps_*/(2C_{\mathsf R})$. By \eqref{es:eq:cocycle} for $\gamma^{-1}$ and
$\gamma$,
\begin{equation*}
\varkappa_i(\gamma^{-1},\gamma p)=\varkappa_i(\gamma,p)^{-1}\qquad\text{if }p,\gamma p\in \cC^i_3 .
\end{equation*}
When $K$ is finite, $\sa_i=0$ and $\varkappa_i(\gamma,p)=k_\gamma$.

\subsection{The Equivariant Coordinate}\label{es:sec:coordinate}

Twisting each summand by the inverse cocycle gives an exactly equivariant coordinate.
Its Lipschitz error is controlled by \eqref{th:eq:weightedaverage}; its limit commutes
with $K$ and parametrizes a slice at $z$.

\subsubsection{The Twisted Average}\label{es:sec:twisted}

We define the coordinate and state its properties. It is the object to which the admissible
sections of \cref{es:sec:frames} and the gluing of \cref{ch:globalbundle} refer. For $p\in \cC^i_1$ let
\begin{equation}\label{es:eq:xidef}
\tilde\xi^{\,0}_i(p):=\sum_\gamma\varpi_\gamma(p)\,\mathsf R\bigl(\varkappa_i(\gamma^{-1},p)\bigr)^{-1}
P_Vh_i(\gamma^{-1}p),\qquad
\tilde\xi_i:=\tilde\xi^{\,0}_i-P_0\,\tilde\xi^{\,0}_i(z_i),
\end{equation}
the sum over the elements $\gamma$ with $\varpi_\gamma(p)>0$; for these $\gamma^{-1}p$ lies
in $\cC^i_3$, so the cocycle is defined at $(\gamma^{-1},p)$.

\begin{proposition}[The Equivariant Coordinate]\label{es:prop:coordinate}
For all large $i$ the map $\tilde\xi_i\colon \cC^i_1\to V$ is continuous and has the
following properties.
\begin{enumerate}[label=\textup{(X\arabic*)},ref=\textup{(X\arabic*)},leftmargin=3.2em]
\item\label{es:X:equiv} $\tilde\xi_i(\gamma p)=\mathsf R\bigl(\varkappa_i(\gamma,p)\bigr)\,\tilde\xi_i(p)$
whenever $\gamma\in\Gamma_i$ and both $p$ and $\gamma p$ lie in $\cC^i_1$.
\item\label{es:X:lip} On $\cC^i_1$ we have $|\tilde\xi_i-P_Vh_i|\le\tfrac3{10}r_\ell$ and
$|P'\tilde\xi_i-P'h_i|\le\tfrac15r_\ell$, and $\tilde\xi_i-P_Vh_i$ is $7\eps_*$-Lipschitz with
respect to $h_i$.
\item\label{es:X:tang} The map $P_0\,\tilde\xi_i$ is invariant:
$P_0\,\tilde\xi_i(\gamma p)=P_0\,\tilde\xi_i(p)$ in the situation of \ref{es:X:equiv}, and \mbox{$P_0\,\tilde\xi_i(z_i)=0$}.
\end{enumerate}
\end{proposition}
\begin{proof}
The map is continuous, because the weights and the cocycle are continuous and the sum
is locally finite.

\ref{es:X:equiv} Let $p$ and $\gamma'p$ lie in $\cC^i_1$. In the sum for
$\tilde\xi^{\,0}_i(\gamma'p)$ we use the index $\mu:=\gamma'^{-1}\gamma$; then
$\varpi_\gamma(\gamma'p)=\varpi_\mu(p)$ by \eqref{es:eq:weights} and
$\gamma^{-1}\gamma'=\mu^{-1}$. The points $p$, $\gamma'p$ and $\mu^{-1}p=\gamma^{-1}(\gamma'p)$
lie in $\cC^i_3$, so by \eqref{es:eq:cocycle}
\begin{equation*}
\varkappa_i(\mu^{-1},p)=\varkappa_i(\gamma^{-1},\gamma'p)\,\varkappa_i(\gamma',p),\qquad\text{that is,}\qquad
\varkappa_i(\gamma^{-1},\gamma'p)^{-1}=\varkappa_i(\gamma',p)\,\varkappa_i(\mu^{-1},p)^{-1}.
\end{equation*}
Hence
\begin{equation*}
\tilde\xi^{\,0}_i(\gamma'p)=\mathsf R\bigl(\varkappa_i(\gamma',p)\bigr)\sum_\mu\varpi_\mu(p)\,
\mathsf R\bigl(\varkappa_i(\mu^{-1},p)\bigr)^{-1}P_Vh_i(\mu^{-1}p)
=\mathsf R\bigl(\varkappa_i(\gamma',p)\bigr)\,\tilde\xi^{\,0}_i(p).
\end{equation*}
The vector $P_0\,\tilde\xi^{\,0}_i(z_i)$ lies in $V_0$ and is fixed by $\mathsf R(K)$, so the same
holds for $\tilde\xi_i$.

\ref{es:X:tang} Since $P_0$ commutes with $\mathsf R(K)$ and $P_0\,\mathsf R(k)=P_0$, the map
$P_0\,\tilde\xi^{\,0}_i$ is the average of $P_0h_i$ with the weights $\varpi_\gamma$.
Its invariance follows from \ref{es:X:equiv}, and the normalization in
\eqref{es:eq:xidef} gives $P_0\,\tilde\xi_i(z_i)=0$.

\ref{es:X:lip} Let $\gamma$ be an element with $\varpi_\gamma(p)>0$ or $\varpi_\gamma(q)>0$ for
points $p,q\in \cC^i_1$; then $\gamma\in\Gamma_i(\tfrac12t_\ell)$. We write $k:=k_{\gamma^{-1}}$
and
\begin{equation*}
\mathsf b_\gamma:=\Ad(k)^{-1}\,\sa_i(\gamma^{-1},\cdot)\colon \cC^i_1\to\fk,\qquad
U_\gamma:=\mathsf R\bigl(\exp(-\mathsf b_\gamma)\bigr)\colon \cC^i_1\to O(n).
\end{equation*}
By the definitions \eqref{es:eq:adef} and \eqref{es:eq:epsgamma}, and since
$k^{-1}\exp(-\mathsf b)\,k=\exp(-\Ad(k)^{-1}\mathsf b)$,
\begin{equation*}
\mathsf R\bigl(\varkappa_i(\gamma^{-1},\cdot)\bigr)^{-1}P_Vh_i\circ\gamma^{-1}
=\mathsf R(k)^{-1}\mathsf R\bigl(\exp(-\sa_i(\gamma^{-1},\cdot))\bigr)\,\mathsf R(k)\bigl(P_Vh_i+\tilde\eps_\gamma\bigr)
=U_\gamma\bigl(P_Vh_i+\tilde\eps_\gamma\bigr).
\end{equation*}
Put $A:=C_{\mathsf R}\sE_0(t_\ell)$ and $L_\varpi:=5C_\#(\ell)/\ell$.
By \cref{es:prop:cocycle}\textup{(i)} and \eqref{es:eq:crho},
$\|U_\gamma-\one\|\le A$ and
$\|U_\gamma(p)-U_\gamma(q)\|\le A\,L_\varpi|h_i(p)-h_i(q)|$.
The latter follows by factoring $U_\gamma(q)$ and applying \eqref{es:eq:crho}
to $\mathsf b_\gamma(q)-\mathsf b_\gamma(p)$, since $\fk$ is abelian. Thus
\begin{equation}\label{es:eq:xisplit}
\tilde\xi^{\,0}_i-P_Vh_i=\sum_\gamma\varpi_\gamma \mathsf e_\gamma,\qquad
\mathsf e_\gamma:=(U_\gamma-\one)P_Vh_i+U_\gamma\tilde\eps_\gamma .
\end{equation}
For every index active at $p$ or $q$, \eqref{es:eq:epsbounds} and $|P_Vh_i|<\ell$ give
\[
|\mathsf e_\gamma|\le A\ell+\sE_2(\ell),\qquad
|\mathsf e_\gamma(p)-\mathsf e_\gamma(q)|
\le\bigl(6\eps_*+\sE_1(t_\ell)+A+A\,L_\varpi(\ell+\sE_2(\ell))\bigr)|h_i(p)-h_i(q)|.
\]
Apply \eqref{th:eq:weightedaverage} with $\nu=0$ and use \eqref{es:eq:thetalip}.
The Lipschitz constant of \eqref{es:eq:xisplit} is at most
\[
6\eps_*+\sE_1(t_\ell)+L_\varpi\sE_2(\ell)+A(1+2L_\varpi\ell+L_\varpi\sE_2(\ell))
=6\eps_*+\sE_3(\ell)+\sE_4(\ell)\le7\eps_*.
\]
Also $A\ell,\sE_2(\ell)\le r_\ell/10$ by \eqref{th:eq:ellsmall} and \eqref{es:eq:ellone}, so
$|\tilde\xi^{\,0}_i-P_Vh_i|\le r_\ell/5$. Since $U_\gamma$ fixes $V_0$,
\mbox{$P_0\tilde\xi^{\,0}_i(z_i)=\sum_\gamma\varpi_\gamma(z_i)P_0\tilde\eps_\gamma(z_i)$} has norm
at most $\sE_2(\ell)\le r_\ell/10$. Subtracting this constant gives \ref{es:X:lip}, including
the sharper bound on the $V'$-component.
\end{proof}

We put
\begin{equation*}
\check h_i:=P_Th_i+\tilde\xi_i\colon \cC^i_1\to\R^n .
\end{equation*}
By \ref{es:X:lip},
\begin{equation}\label{es:eq:psibilip}
(1-7\eps_*)\,|h_i(p)-h_i(q)|\le|\check h_i(p)-\check h_i(q)|\le(1+7\eps_*)\,|h_i(p)-h_i(q)|
\qquad(p,q\in \cC^i_1).
\end{equation}
So $\check h_i$ is injective and continuous on $\cC^i_1$, an open subset of a topological
manifold of dimension $n$ (\cref{ls:lem:existence}); by the invariance of domain it is a
homeomorphism of $\cC^i_1$ onto an open subset of $\R^n$, and $\tilde\xi_i=P_V\check h_i$ is a coordinate projection for $\check h_i$. In the coordinate $h_i$ the map $\check h_i$ is a map
$\nu\mapsto\nu+\mathsf e(\nu)$ of the open ball of radius $\ell$ with $\mathsf e$ Lipschitz with the
constant $7\eps_*$ and $|\mathsf e(0)|\le\tfrac3{10}r_\ell$; as in
\cref{th:lem:transplant}\textup{(ii)} its image contains the open ball of radius
$(1-7\eps_*)\ell$ about $\mathsf e(0)$, hence the closed ball of radius $\tfrac34\ell$ about $0$.

\subsubsection{The Limit Coordinate}\label{es:sec:limitcoordinate}

Normalized deck averages converge to Haar averages by \cref{m:fr:lem:counting}\textup{(ii)}.
This identifies both the coordinate and the cocycle limits. Write $\chi_i$ for the cut-off
$\chi$ of \cref{th:sec:weights}, and let $\chi_\infty$ equal $\chi_0(|h|/\ell)$ on $B$ and
zero elsewhere. As for $\Gamma_i(t)$ in \cref{th:sec:standing}, put
\[
G(t):=\{g\in G:\hat d(gz,z)<t\}\quad(t>0),\qquad
\bar\chi_\infty(p):=\int_G\chi_\infty(g'^{-1}p)\,d\nu_G(g')\quad(p\in\cC_3).
\]
For $g,g'\in G(\tfrac12t_\ell)$ let
\[
\omega_\infty(g,g'):=\log\!\left[\mathsf k(gg')
 \bigl(\mathsf k(g)\mathsf k(g')\bigr)^{-1}\right]\in\fk.
\]
The logarithm is the inverse of $\exp$ on $\fk(2r_K)$: it exists by
\eqref{es:eq:lambdaproduct} and has norm at most $\sE_0(t_\ell)$.
If $p\in\cC_3$ and $\chi_\infty(g'^{-1}p)>0$, then $g'\in G(\tfrac12t_\ell)$ by
\eqref{th:eq:metricsize}. Since $\chi_\infty(p)>0$, the denominator $\bar\chi_\infty(p)$ is
positive. Thus, on $G(\tfrac12t_\ell)\times\cC_3$, we may define
\begin{equation}\label{es:eq:haarcocycle}
\begin{aligned}
\sa_\infty(g,p)&:=
 \frac{\int_G\chi_\infty(g'^{-1}p)\,\omega_\infty(g,g')\,d\nu_G(g')}{\bar\chi_\infty(p)},\\
\varkappa_\infty(g,p)&:=\exp(\sa_\infty(g,p))\,\mathsf k(g).
\end{aligned}
\end{equation}
Weighted integrands are set to zero where the cut-off vanishes. They are continuous
with locally common compact support in $G$, so these maps are continuous and
$|\sa_\infty|\le\sE_0(t_\ell)$. For $k\in K$, \eqref{es:eq:lambdaK} gives
$\omega_\infty(k,g')=0$, hence $\varkappa_\infty(k,p)=k$.

\begin{lemma}[The Limit Coordinate]\label{es:lem:limit}
There is a map $\xi\colon \cC_1\to V$ such that
$\tilde\xi_i(p_i)\to\xi(p)$ whenever points $p_i\in \cC^i_1$ converge to $p\in \cC_1$. It has
the following properties.
\begin{enumerate}[label=\textup{(\roman*)}]
\item On $\cC_1$ we have $|\xi-P_Vh|\le\tfrac3{10}r_\ell$, the map $\xi-P_Vh$ is
$7\eps_*$-Lipschitz with respect to $h$, and $\xi(z)=0$.
\item We have $\xi(kp)=\mathsf R(k)\,\xi(p)$ for all $k\in K$ and $p\in \cC_1$.
\item If $g\in G$ and both $p$ and $gp$ lie in $\cC_1$, then $\xi(gp)\in\mathsf R(K)\,\xi(p)$.
\end{enumerate}
\end{lemma}
\begin{proof}
The cutoffs $\chi_i$ converge continuously to $\chi_\infty$, by
\cref{th:lem:transplant}\textup{(iii)} and their common interior support.
If $\gamma_i\eGHto g\in G(\tfrac12t_\ell)$ and $\mu_i\eGHto g'\in G(\tfrac12t_\ell)$, convergence of
products and labels, followed by the common logarithm branch, gives
$\omega_i(\gamma_i,\mu_i)\to\omega_\infty(g,g')$.
For $p_i\to p\in\cC_3$, the tests
$\mu\mapsto\chi_i(\mu^{-1}p_i)$ and
$\mu\mapsto\chi_i(\mu^{-1}p_i)\omega_i(\gamma_i,\mu)$ therefore converge continuously
to their Haar-integral counterparts. They have common bounded-displacement support;
at its boundary the cut-off tends to zero and the defects remain bounded.
Applying \cref{m:fr:lem:counting}\textup{(ii)} to numerator and denominator of
\eqref{es:eq:adef}, with the counting normalization cancelled, proves
\begin{equation}\label{es:eq:haarconvergence}
\sa_i(\gamma_i,p_i)\longrightarrow\sa_\infty(g,p),\qquad
\varkappa_i(\gamma_i,p_i)\longrightarrow\varkappa_\infty(g,p).
\end{equation}
No further subsequence is needed.

Apply the same argument, componentwise, to the bounded weighted summands in
\eqref{es:eq:xidef}. It gives the limit
\[
\begin{aligned}
\xi^0(p)&:=
\frac{\int_G\chi_\infty(g'^{-1}p)\,
 \mathsf R(\varkappa_\infty(g'^{-1},p))^{-1}P_Vh(g'^{-1}p)\,d\nu_G(g')}
 {\bar\chi_\infty(p)},\\
\xi(p)&:=\xi^0(p)-P_0\xi^0(z).
\end{aligned}
\]
Every contributing $g'$ and $g'^{-1}$ lies in $G(\tfrac12t_\ell)$, so the formula is defined.
The estimates of (i) pass to the limit from \ref{es:X:lip}.
For $p,gp\in\cC_1$, choose $\gamma_i\eGHto g$ and $p_i\to p$. Then $g\in G(\tfrac12t_\ell)$,
and for large $i$ the points $p_i,\gamma_ip_i$ lie in $\cC^i_1$. Passing to the limit in
\begin{equation}\label{es:eq:stagerelation}
\tilde\xi_i(\gamma_ip_i)=\mathsf R\bigl(\varkappa_i(\gamma_i,p_i)\bigr)\,\tilde\xi_i(p_i)
\end{equation}
gives $\xi(gp)=\mathsf R(\varkappa_\infty(g,p))\xi(p)$, proving (iii).
For $g=k\in K$, the set $\cC_1$ is $K$-invariant and
$\varkappa_\infty(k,p)=k$, proving (ii). Finally $\xi(z)$ is fixed by $K$, hence lies
in $V_0$, and its $P_0$-component vanishes by the normalization; so $\xi(z)=0$.
\end{proof}

By \cref{es:lem:limit}\textup{(i)} the map
\begin{equation*}
\check h:=P_Th+\xi\colon \cC_1\to\R^n
\end{equation*}
is injective and continuous, with $\check h(z)=0$ and $P_T\check h=P_Th$, and by (ii) it commutes
with the actions of $K$. As for $\check h_i$, its image contains the closed ball of radius
$\tfrac12\ell$ about $0$. The set $\cC_1$ is an open neighbourhood of $z$ in $B_{\rho_s}(z)$
that is invariant under $K$. So \cref{th:lem:slices} applies to $U:=\cC_1$, the map $\check h$
and $r_1:=\tfrac12\ell$, and gives the first part of the following lemma.

\begin{lemma}[The Slice and Its Coordinate]\label{es:lem:slice}
Let
\begin{equation*}
S:=\bigl\{q\in \cC_1:P_Th(q)=0,\ |\xi(q)|<\tfrac12\ell\bigr\},\qquad O:=\pi(S).
\end{equation*}
\begin{enumerate}[label=\textup{(\roman*)}]
\item The set $S$ is invariant under $K$, the map $\xi$ restricts to a homeomorphism of
$S$ onto $B_V(\tfrac12\ell)$ that commutes with the actions of $K$, and
$[g,q]\mapsto gq$ is a homeomorphism of $G\times_KS$ onto an open subset of $Z$ that
commutes with the actions of $G$. If $g\in G$ maps a point of $S$ to a point of $S$,
then $g\in K$.
\item The set $O$ is an open neighbourhood of $x$ in $X$, and there is a homeomorphism
\begin{equation*}
\bar\xi\colon O\longrightarrow B_V(\tfrac12\ell)/\mathsf R(K)\qquad\text{with}\qquad
\bar\xi\bigl(\pi(q)\bigr)=\mathsf R(K)\,\xi(q)\quad\text{for }q\in S .
\end{equation*}
\item Let $p\in \cC_{1/2}$ with $|\xi(p)|\le9r_\ell$. Then there is $\zeta\in\fm$ with
$|\zeta|\le2|P_Th(p)|/c_z$ such that $q:=\exp(\zeta)p$ lies in $S$, and
$\pi(p)\in O$ with $\bar\xi(\pi(p))=\mathsf R(K)\,\xi(p)$.
\end{enumerate}
\end{lemma}
\begin{proof}
(i) is \cref{th:lem:slices}.

(ii) The set $S$ lies in $B_{d_+(\ell)}(z)$, hence in $\Omega'$. By (i) two points of $S$
have the same image under $\pi$ if and only if they lie in one orbit of $K$, that is, if
and only if their images under $\xi$ lie in one orbit of $\mathsf R(K)$. So $\bar\xi$ is well
defined and bijective. It is a homeomorphism, and $O$ is open. Indeed, by (i) the orbit space $S/K$ of
$G\times_KS$ is homeomorphic to the open subset $GS/G$ of $\Omega/G$, and $\xi$ induces a
homeomorphism of $S/K$ onto $B_V(\tfrac12\ell)/\mathsf R(K)$. By \ref{fl:S:limit} the map $\pi$ induces an isometry of $\Omega/G$,
with the quotient metric, onto $B_{R_0}(a)$, and the quotient metric induces the quotient
topology because the action is proper. So $O$, the image of $GS/G$, is open in $X$, and
$\bar\xi$ is the composition of these homeomorphisms.

(iii) The number $t:=|P_Th(p)|$ is at most $\tfrac12\ell$, so $2t\le c_zr_3$, and $p$
lies in $B_{\rho_s}(z)$. By \cref{th:lem:meets} there is $\zeta\in\fm$ with
$|\zeta|\le2t/c_z$ such that $q:=\exp(\zeta)p$ satisfies $P_Th(q)=0$,
$|P_Vh(q)-P_Vh(p)|\le t$ and $\hat d(q,p)\le\tfrac12s_1$. So $q\in B$, and by
\cref{es:lem:limit}\textup{(i)}
\begin{equation*}
|h(q)|=|P_Vh(q)|\le|\xi(p)|+\tfrac3{10}r_\ell+\tfrac12\ell<\ell ,
\end{equation*}
that is, $q\in \cC_1$. By \cref{es:lem:limit}\textup{(iii)} we have $\xi(q)\in\mathsf R(K)\xi(p)$,
so $|\xi(q)|\le9r_\ell<\tfrac12\ell$ and $q\in S$. Hence $\pi(p)=\pi(q)\in O$ and
$\bar\xi(\pi(p))=\mathsf R(K)\xi(q)=\mathsf R(K)\xi(p)$.
\end{proof}

The set $S$ is a slice at $z$ and $\xi$ is its coordinate; it commutes with $K$, and by
\cref{es:lem:limit}\textup{(iii)} it describes the map $\pi$ on the whole coordinate ball:
two points of $\cC_1$ in one orbit of $G$ have images under $\xi$ in one orbit of
$\mathsf R(K)$. The proof of
\cref{ls:thm:limits}, applied to $S$ and $\xi$, gives the continuous
section $\sigma$ of the frame space over $B_{r/2}(z)$ that is aligned with the frame map of
$e$, with
\begin{equation}\label{es:eq:section}
\sigma(z)=e,\qquad\sigma(q)\circ k^{-1}=\sigma(kq)\cdot A_k\qquad\text{for }k\in K,
\end{equation}
and a homeomorphism
\begin{equation*}
\Xi\colon O(n)\times_IB_V(\tfrac12\ell)\longrightarrow\pr^{-1}(O)\subset Y,\qquad
\Xi([\one,v])=\tau\bigl(q,\sigma(q)\bigr)\quad\text{for }q\in S\text{ with }\xi(q)=v,
\end{equation*}
that commutes with the actions of $O(n)$; here $I$ acts on $V$ by
$A_k\cdot v:=\mathsf R(k)v$, and $\Xi([\one,0])=y$.

\begin{definition}[The Tube]\label{es:def:tube}
The \emph{tube} at $y$ is this linear tube $\Xi$ of $O(n)\times_IB_V(\tfrac12\ell)$ onto the open
neighbourhood $\pr^{-1}(O)$ of the orbit of $y$ (\cref{ls:def:slice}). Its \emph{data} are the centre $a$, the frame
$(z,e)$ over $y$ with the isotropy groups $K\cong I$, the scale $\ell$, the harmonic chart $h$, the
slice $S$, the coordinate $\xi$ and the section $\sigma$.
\end{definition}

We have
\begin{equation}\label{es:eq:prtube}
\pr\bigl(\Xi([A,v])\bigr)=\bar\xi^{-1}\bigl(\mathsf R(K)v\bigr)\qquad\text{for }A\in O(n),\ v\in B_V(\tfrac12\ell).
\end{equation}

\subsection{The Domain and Nearby Points}\label{es:sec:model}

In this subsection we fix the domain of the equivariant coordinate in the charts $Z_i$
(\cref{es:lem:domain}), and we show that the slice at $z$ contains slices at all of its
points, whatever their isotropy groups, and that their translates by $G$ give slices at all points
of $GS$ (\cref{es:thm:nearbyprecise}).

\subsubsection{The Domain}\label{es:sec:domain}

We fix the domain of the coordinate; all later constructions in the tube take place on it. Let
\begin{equation*}
r_V:=8r_\ell,\qquad
D^i:=\bigl\{p\in \cC^i_{1/2}:|P_Th_i(p)|<\tfrac14\ell,\ |\tilde\xi_i(p)|<r_V\bigr\},\qquad
\cW_i:=\pi_i(D^i),
\end{equation*}
and let $B_T(r)$ be the open ball of radius $r$ about $0$ in $T$. The set $\cW_i$ is open
in $X_i$, because $\pi_i$ is an open map.

\begin{lemma}[The Domain]\label{es:lem:domain}
For all large $i$ the following holds.
\begin{enumerate}[label=\textup{(\roman*)}]
\item The map $\check h_i$ maps $D^i$ homeomorphically onto
$B_T(\tfrac14\ell)\times B_V(r_V)$.
\item If $p\in D^i$ and $\gamma\in\Gamma_i$ with $\gamma p\in \cC^i_1$, then
$|\tilde\xi_i(\gamma p)|=|\tilde\xi_i(p)|$.
\item For every $p\in D^i$ there is $\gamma\in\Gamma_i$ with $\gamma p\in D^i\cap \cC^i_{1/16}$.
\end{enumerate}
\end{lemma}
\begin{proof}
(i) Let $\nu\in B_T(\tfrac14\ell)\times B_V(r_V)$; then $|\nu|<208\,r_\ell$. By
\cref{es:sec:twisted} there is $p\in \cC^i_1$ with $\check h_i(p)=\nu$, and by
\eqref{es:eq:psibilip} and \ref{es:X:lip}
\begin{equation*}
(1-7\eps_*)\,|h_i(p)|\le|\check h_i(p)-\check h_i(z_i)|\le|\nu|+\tfrac3{10}r_\ell,
\end{equation*}
so $|h_i(p)|<\tfrac12\ell$, since $\eps_*\le\tfrac1{200}$. Thus $D^i$ is the preimage of the
product under the homeomorphism $\check h_i$.

(ii) This is \ref{es:X:equiv}, because $\mathsf R$ is orthogonal.

(iii) Suppose that the assertion fails for infinitely many $i$, at points $p_i\in D^i$.
After passing to a subsequence, $p_i\to p$. The estimate in (i) gives
$|h(p)|\le(208+\tfrac3{10})r_\ell /(1-7\eps_*)<\tfrac12\ell$, and
$|\xi(p)|\le8r_\ell$. By \cref{es:lem:slice}\textup{(iii)} there is $g\in G$ such that
$q:=gp\in S$ and $|\xi(q)|=|\xi(p)|\le8r_\ell$. Since $P_Th(q)=0$,
\[
|h(q)|\le|\xi(q)|+\tfrac3{10}r_\ell\le\tfrac{83}{10}r_\ell<\tfrac1{16}\ell.
\]
Choose $\gamma_i\in\Gamma_i$ with $\gamma_i\eGHto g$. Then $\gamma_ip_i\to q$, so
$\gamma_ip_i\in \cC^i_{1/16}$ for large $i$. By (ii), \mbox{$|\tilde\xi_i(\gamma_ip_i)|=|\tilde\xi_i(p_i)|<r_V$}, and consequently
$\gamma_ip_i\in D^i$, a contradiction.
\end{proof}

\subsubsection{Nearby Points}\label{es:sec:nearby}

The slice at a point of $S$ is the part of $S$ over a linear slice of the representation
at its image under $\xi$; so one tube describes all strata near $x$ (\cref{es:lem:strata}). For $v\in V$ let $K_v$
be the group of the elements $k\in K$ with $\mathsf R(k)v=v$, and let $N_v\subset V$ be the
orthogonal complement of the tangent space of the
orbit $\mathsf R(K)v$ at $v$; it is invariant under $\mathsf R(K_v)$, and by the theorem on
linear tubes \cite[Chapter~VI, Section~2]{Bre72} there is $r_\perp>0$ such that
\begin{equation}\label{es:eq:lineartube}
K\times_{K_v}\bigl(v+B_{N_v}(r_\perp)\bigr)\longrightarrow V,\qquad[k,v+w]\longmapsto\mathsf R(k)(v+w),
\end{equation}
is a homeomorphism onto an open neighbourhood of the orbit that commutes with $K$; here
$B_{N_v}(r_\perp)$ is the open ball of radius $r_\perp$ about $0$ in $N_v$.

\begin{theorem}[Nearby Points]\label{es:thm:nearbyprecise}
In the setting of \cref{es:sec:setting}, let $v$
be a point of $B_V(\tfrac12\ell)$, let $r_\perp>0$ be such that \eqref{es:eq:lineartube} is a
homeomorphism onto an open subset of $B_V(\tfrac12\ell)$, and put
\begin{equation*}
\Sigma_v:=v+B_{N_v}(r_\perp),\qquad
O^v:=\bar\xi^{-1}\bigl(\mathsf R(K)\Sigma_v/\mathsf R(K)\bigr).
\end{equation*}
Then the set $O^v$ is an open neighbourhood of the point $x':=\bar\xi^{-1}(\mathsf R(K)v)$
in $X$, and $\bar\xi$ induces a homeomorphism of $O^v$ onto $\Sigma_v/\mathsf R(K_v)$. The set
$S_v:=\{q\in S:\xi(q)\in\Sigma_v\}$ is a slice for the action of $G$ at the point
$q_v\in S$ with $\xi(q_v)=v$, the isotropy group of $q_v$ is $K_v$, and $\xi-v$ maps
$S_v$ homeomorphically onto $B_{N_v}(r_\perp)$, compatibly with $K_v$.
\end{theorem}
\begin{proof}
The image $\mathsf R(K)\Sigma_v$ of \eqref{es:eq:lineartube} is open and invariant, so $O^v$
is open, and the quotient of $K\times_{K_v}\Sigma_v$ by $K$ is $\Sigma_v/\mathsf R(K_v)$. By
\cref{es:lem:slice}\textup{(i)} the coordinate $\xi$ identifies the action of $K$ on $S$ with
that of $\mathsf R(K)$ on $B_V(\tfrac12\ell)$. So the isotropy group of $q_v$ in $K$ is
$K_v$, and it is the isotropy group in $G$ by the last assertion of
\cref{es:lem:slice}\textup{(i)}. Moreover $KS_v$ is an open subset of $S$ that is
homeomorphic to $K\times_{K_v}S_v$, so $GS_v$ is homeomorphic to
$G\times_K(K\times_{K_v}S_v)=G\times_{K_v}S_v$ and open in $Z$; that is, $S_v$ is a slice
at $q_v$.
\end{proof}

The isotropy group of $\Xi([\one,v])$ is $I_v=\{A_k:k\in K_v\}$, because $\Xi$ is injective and
commutes with $O(n)$.

\subsection{Admissible Sections of the Frame Spaces}\label{es:sec:frames}

We construct sections of the reduction of a realizing map over the tube
(\cref{es:sec:reduction}) whose cocycle is exactly $\varkappa_i$, the admissible sections
(\cref{es:def:admissible,es:lem:admissible}). Convergent frames first give sections whose cocycles are close
to $\varkappa_i$; an averaged coboundary makes them exactly compatible. The coordinate
$\tilde\xi_i$ itself is independent of these sections, and the coordinate of an admissible
section is close to it; \cref{es:lem:modall} makes the two equal near the centre of the
tube.

\subsubsection{The Reduction of the Realizing Map}\label{es:sec:reduction}

A map into the image of a tube that commutes with the group is determined by its part over the
slice of the tube. This is the standard correspondence between spaces over a homogeneous space
and spaces with an action of the isotropy group \cite[Chapter~II, Section~3]{Bre72}; we record
it for the maps of the frame spaces into the tubes of their limit, where it is used here and in
\cref{ch:globalbundle}. All maps are continuous.

\begin{proposition}[The Local Form]\label{ef:prop:local}
Let $F\colon Y_i\to Y$ be a map that commutes with $O(n)$, let $I\subset O(n)$ be a closed
subgroup and $B$ a space with an action of $I$, and let $\Xi\colon O(n)\times_IB\to Y$ be a tube
with image $U$ (\cref{ls:def:slice}). Put $W:=F^{-1}(\Xi([\one,B]))$ and let $\varphi\colon W\to B$ be
the map with $\Xi([\one,\varphi(w)])=F(w)$. Then $I$ acts freely on $W$, the map $\varphi$ commutes
with $I$, and $[A,w]\mapsto A\cdot w$ is an $O(n)$-homeomorphism of $O(n)\times_IW$ onto
$F^{-1}(U)$ under which $F$ corresponds to $[A,w]\mapsto\Xi([A,\varphi(w)])$.
\end{proposition}
\begin{proof}
Let $\pi_I\colon F^{-1}(U)\to O(n)/I$ be the composition of $F$ with $\Xi^{-1}$ and the
projection $[A,b]\mapsto AI$. It commutes with $O(n)$, and $W=\pi_I^{-1}(I)$, because the
classes $[\one,b]$ are the points of $O(n)\times_IB$ over the coset $I$. By the correspondence just
recalled \cite[Chapter~II, Section~3]{Bre72}, $W$ is invariant under $I$ and
$[A,w]\mapsto A\cdot w$ is an $O(n)$-homeomorphism of $O(n)\times_IW$ onto $F^{-1}(U)$. The group
$I$ acts freely on $W$ because $O(n)$ acts freely on $Y_i$. The map $\varphi$ is continuous
because $b\mapsto\Xi([\one,b])$ is a homeomorphism onto its image. For $A\in O(n)$ we have
$F(A\cdot w)=A\cdot\Xi([\one,\varphi(w)])=\Xi([A,\varphi(w)])$, and for $A\in I$ the class $[A,b]$
equals $[\one,A\cdot b]$; so $\varphi$ commutes with $I$, and $F$ corresponds to
$[A,w]\mapsto\Xi([A,\varphi(w)])$.
\end{proof}

We apply this to equivariant realizing maps over the tube and fix the notation for
frames. Let $F=F^{(i)}$ be equivariant realizing maps, and let
\begin{equation*}
W^F_i:=F^{-1}\bigl(\Xi([\one,B_V(\tfrac12\ell)])\bigr),\qquad
\varphi^F\colon W^F_i\to B_V(\tfrac12\ell)\quad\text{with}\quad\Xi([\one,\varphi^F(w)])=F(w).
\end{equation*}
By \cref{ef:prop:local} the group $I$ acts freely on $W^F_i$, the map $\varphi^F$
commutes with $I$, and \mbox{$[A,w]\mapsto A\cdot w$} is a homeomorphism of
$O(n)\times_IW^F_i$ onto $F^{-1}(\pr^{-1}(O))$ under which $F$ corresponds to
$[A,w]\mapsto\Xi([A,\varphi^F(w)])$. If two points of $W^F_i$ lie over the same point of
$X_i$, they differ by an element of $O(n)$ that maps a point of $\Xi([\one,B_V])$ to such a
point, hence by an element of $I$. The set $\Xi([\one,B_V(\tfrac12\ell)])$ is closed in
$\pr^{-1}(O)$, so $W^F_i$ is closed in the open set $F^{-1}(\pr^{-1}(O))$.

The frames, the actions of $O(n)$, $\Gamma_i$ and $G$ on them and the chart maps $\tau_i$ and
$\tau$ are those of \ref{fl:S:frames} and \ref{fl:S:charts}. Since $\Gamma_i$ acts freely on $Z_i$,
the map $\tau_i$ is injective on the frames over one point. Over a point $p$ of $\Omega'$ the
map $\tau$ identifies the frames $(p,e)$ and $(p,e\circ k^{-1})$ for the elements $k$ of the
isotropy group $G_p$, and no other frames over $p$; so it is injective on the frames over
$p$ exactly when $G_p$ is trivial, and in general its fibres over $p$ are the orbits of $G_p$.
We choose frames $(z_i,e_i)$ that converge to $(z,e)$ and use their frame trivializations
$\ft_i$ on $B_{r/2}(z_i)\times O(n)$ and $\ft$ on $B_{r/2}(z)\times O(n)$ of \ref{fl:S:product}; the balls $B_{d_+(\ell)}(z_i)$ and $B_{d_+(\ell)}(z)$ lie in their domains, by the
choice of $\ell_*$ in \cref{th:sec:scale}.

\subsubsection{Smaller Tubes and Other Realizing Maps}\label{es:sec:realizing}

The gluing shrinks the tube and changes the realizing map. We introduce the nested domains and
the cut-off function used for this, and we show that the preimage of a smaller tube under any
realizing map lies over the corresponding domain. We work on one fixed tube under
\eqref{th:eq:threshold}. For $0<t\le8$ put
\begin{gather*}
D^i(t):=\{p\in D^i:|\tilde\xi_i(p)|<t r_\ell\},\qquad
\cW_i(t):=\pi_i(D^i(t)),\\
D^i_c:=\{p\in D^i:|P_Th_i(p)|\le\ell/8,\ |\tilde\xi_i(p)|\le6r_\ell\},
\end{gather*}
and for $0<t\le10$ put $O(t):=\bar\xi^{-1}(B_V(t r_\ell)/\mathsf R(K))$. This is an open subset
of $X$, because $10r_\ell<\tfrac12\ell$; the comparison of tubes in \cref{ch:globalbundle} uses the
larger range.
Every point of the closure of $\cW_i(6)$ in $\cW_i$ has a lift in $D^i_c$,
by \cref{es:lem:domain}\textup{(iii)}.

Let $\chi_3\colon[0,\infty)\to[0,1]$ equal $1$ on $[0,2]$, vanish on
$[4,\infty)$ and be linear in between, and put
\[
\chi^\xi_i:=\chi_3(|\tilde\xi_i|/r_\ell)\colon D^i\to[0,1].
\]

\begin{lemma}[Smaller Tubes]\label{es:lem:smallerreduction}
Let $F=F^{(i)}$ be realizing maps. For fixed $0<t<t'<8$, for all large $i$,
\[
F^{-1}\bigl(\Xi(O(n)\times_I B_V(t r_\ell))\bigr)
\subset \pr^{-1}(\cW_i(t')).
\]
\end{lemma}
\begin{proof}
Otherwise choose points $u_i$ in the set on the left whose projections lie
outside $\cW_i(t')$. Write $F(u_i)=\Xi([A_i,v_i])$, where $|v_i|<t r_\ell$.
A subsequence has $A_i\to A$, $v_i\to v$ and $u_i\to u$; since $F$ is
realizing, $u=\Xi([A,v])$. Let $q\in S$ satisfy $\xi(q)=v$.
By the frame convergence in \textup{\ref{fl:S:Yconv}}, the points $u_i$
have lifts to frames over $q_i\to q$ in this chart. Thus
$\pi_i(q_i)=\pr(u_i)$, $P_Th_i(q_i)\to0$ and
$|\tilde\xi_i(q_i)|\to|v|\le t r_\ell<t'r_\ell$. For large $i$ we have
$q_i\in D^i(t')$, a contradiction.
\end{proof}

\subsubsection{Orbital Coordinates}\label{es:sec:limitsection}

The slice and the exponentials of $\fm$ give product coordinates on an
open set containing the limits of the domains $D^i$. They are the coordinates in which the
limit section of \cref{es:sec:sigmahat} is written.

\begin{lemma}[Orbital Coordinates]\label{es:lem:limitsection}
Let $U_S$ be the set of the points $\exp(\zeta)q$ with $\zeta\in\fm(r_c)$ and $q\in S$.
\begin{enumerate}[label=\textup{(\roman*)}]
\item The map $(\zeta,q)\mapsto\exp(\zeta)q$ is a homeomorphism of $\fm(r_c)\times S$ onto
$U_S$, and $U_S$ is open in $Z$.
\item Every point $p\in \cC_{1/2}$ with $|P_Th(p)|\le\tfrac14\ell$ and $|\xi(p)|\le r_V$ lies
in $U_S$, and $p=\exp(\zeta)q$ with $q\in S$ and $|\zeta|\le\ell/(2c_z)$.
\end{enumerate}
\end{lemma}
\begin{proof}
(i) If $\exp(\zeta)q=\exp(\zeta')q'$, then $\exp(\zeta')^{-1}\exp(\zeta)$ maps the point
$q$ of $S$ to the point $q'$ of $S$, so it is an element $k$ of $K$, by
\cref{es:lem:slice}\textup{(i)}. Then $\exp(\zeta)=\exp(\zeta')k$, and $\zeta=\zeta'$, $k=1$
and $q=q'$, because the map of \cref{es:sec:lambda} is injective on $\fm(2r_{\mathsf k})\times K$.
So the map is injective. It is continuous, and its domain and the open subset
$\Omega$ of $Z$ are topological manifolds of dimension $\dim\fm+\dim V=n$; by the
invariance of domain it is a homeomorphism onto an open set.

(ii) By \cref{es:lem:slice}\textup{(iii)} there is $\zeta\in\fm$ with
$|\zeta|\le2|P_Th(p)|/c_z\le\ell/(2c_z)$ and $\exp(\zeta)p\in S$, and
$\ell/(2c_z)\le r_c$ by \eqref{es:eq:ellone}. So $p=\exp(-\zeta)\exp(\zeta)p$ lies in $U_S$.
\end{proof}

\subsubsection{The Limit Cocycle}\label{es:sec:limitcocycle}

The Haar limit \eqref{es:eq:haarcocycle} satisfies the exact cocycle identity on the
arrows whose endpoints lie in a smaller coordinate ball. Let
\begin{equation*}
\mathcal G:=\bigl\{(g,p)\in G\times \cC_{3/4}:gp\in \cC_{3/4}\bigr\}.
\end{equation*}
By \eqref{th:eq:metricsize}, $\mathcal G\subset G(\tfrac12t_\ell)\times\cC_3$, so the convergence
\eqref{es:eq:haarconvergence} applies on this domain.

\begin{lemma}[The Limit Cocycle]\label{es:lem:limitcocycle}
The maps $\sa_\infty$ and $\varkappa_\infty$ of \eqref{es:eq:haarcocycle}, restricted to
$\mathcal G$, are continuous, satisfy $|\sa_\infty|\le\sE_0(t_\ell)$ and the convergence
\eqref{es:eq:haarconvergence}, and have the following properties.
\begin{enumerate}[label=\textup{(\roman*)}]
\item If $(g',p)$ and $(g,g'p)$ lie in $\mathcal G$, then
$\varkappa_\infty(gg',p)=\varkappa_\infty(g,g'p)\,\varkappa_\infty(g',p)$.
\item We have $\varkappa_\infty(k,p)=k$ for $k\in K$ and $p\in \cC_{3/4}$.
\item We have $\xi(gp)=\mathsf R(\varkappa_\infty(g,p))\,\xi(p)$ for $(g,p)\in\mathcal G$.
\end{enumerate}
\end{lemma}
\begin{proof}
The convergence was proved in \cref{es:lem:limit}. For (i) choose
$\gamma_i\eGHto g$, $\gamma'_i\eGHto g'$ and $p_i\to p$.
The three points $p_i$, $\gamma'_ip_i$ and $\gamma_i\gamma'_ip_i$ eventually lie in
$\cC^i_{3/4}$, so \eqref{es:eq:cocycle} applies. Its limit, using
$\gamma_i\gamma'_i\eGHto gg'$, is (i). Part (ii) follows from
$\omega_\infty(k,g')=0$, and (iii) follows by passing to the limit in
\eqref{es:eq:stagerelation}.
\end{proof}

\subsubsection{A Cocycle-Compatible Limit Section}\label{es:sec:sigmahat}

We use the limit cocycle of \cref{es:lem:limitcocycle} to extend the
section $\sigma$ of \eqref{es:eq:section} from the slice to the domain of the orbital
coordinates, so that $G$ acts on the extension through $\varkappa_\infty$.
For $p\in U_S$ let \mbox{$\zeta(p)\in\fm(r_c)$} and $\pi_S(p)\in S$ be the elements with $p=\exp(\zeta(p))\,\pi_S(p)$; they depend continuously on
$p$ (\cref{es:lem:limitsection}\textup{(i)}). Let $\hat U_S:=U_S\cap \cC_{1/2}$. The slice $S$
lies in $\cC_{3/4}$, because $|h|\le|\xi|+\tfrac3{10}r_\ell$ on it; so for $p\in \hat U_S$ the
pair $(\exp(\zeta(p)),\pi_S(p))$ lies in $\mathcal G$, and we put
\begin{equation*}
\begin{split}
\varkappa_S(p)&:=\varkappa_\infty\bigl(\exp(\zeta(p)),\pi_S(p)\bigr)\in K,\\
\hat\sigma(p)&:=\exp(\zeta(p))\cdot(\pi_S(p),\sigma(\pi_S(p)))\cdot A_{\varkappa_S(p)}^{-1}
\qquad(p\in \hat U_S).
\end{split}
\end{equation*}

\begin{lemma}[The Limit Section with the Coordinate $\xi$]\label{es:lem:sigmahat}
The map $\hat\sigma$ is a continuous section of the frame space over $\hat U_S$, and the
following holds.
\begin{enumerate}[label=\textup{(\roman*)}]
\item For $p\in \hat U_S$ we have $\tau(\hat\sigma(p))=\Xi([\one,\xi(p)])$.
\item If $g\in G$ and $p$ and $gp$ lie in $\hat U_S$, then
$g\cdot\hat\sigma(p)=\hat\sigma(gp)\cdot A_{\varkappa_\infty(g,p)}$.
\item The map $v\mapsto\Xi([\one,v])$ is a homeomorphism of $B_V(\tfrac12\ell)$ onto a closed
subset of $\pr^{-1}(O)$.
\end{enumerate}
\end{lemma}
\begin{proof}
The maps $\zeta$, $\pi_S$ and $\varkappa_\infty$ are continuous, so $\varkappa_S$ and $\hat\sigma$
are.

(i) By \cref{es:lem:limitcocycle}\textup{(iii)}, $\xi(p)=\mathsf R(\varkappa_S(p))\,\xi(\pi_S(p))$. The
map $\tau$ is $G$-invariant and commutes with the right action of $O(n)$,
and $A\cdot y=y\cdot A^{-1}$. Hence
\begin{equation*}
\tau\bigl(\hat\sigma(p)\bigr)=A_{\varkappa_S(p)}\cdot\Xi\bigl([\one,\xi(\pi_S(p))]\bigr)
=\Xi\bigl([A_{\varkappa_S(p)},\xi(\pi_S(p))]\bigr)=\Xi\bigl([\one,\mathsf R(\varkappa_S(p))\,\xi(\pi_S(p))]\bigr).
\end{equation*}

(ii) The element $k':=\exp(-\zeta(gp))g\exp(\zeta(p))$ maps $\pi_S(p)\in S$
to $\pi_S(gp)\in S$, so $k'\in K$ by \cref{es:lem:slice}\textup{(i)}.
The identity \eqref{es:eq:section} for $\sigma$ and the commuting left and right actions
on frames give
\begin{equation*}
g\cdot\hat\sigma(p)
=\exp(\zeta(gp))\cdot(\pi_S(gp),\sigma(\pi_S(gp)))\cdot A_{k'}A_{\varkappa_S(p)}^{-1}
=\hat\sigma(gp)\cdot A_{\varkappa_S(gp)\,k'\,\varkappa_S(p)^{-1}} .
\end{equation*}
We write $q:=\pi_S(p)$. By \cref{es:lem:limitcocycle}\textup{(i)} and~\textup{(ii)}, applied to
points of $\cC_{3/4}$,
\begin{equation*}
\varkappa_S(gp)\,k'=\varkappa_\infty\bigl(\exp(\zeta(gp)),k'q\bigr)\,\varkappa_\infty(k',q)
=\varkappa_\infty\bigl(g\exp(\zeta(p)),q\bigr)=\varkappa_\infty(g,p)\,\varkappa_S(p),
\end{equation*}
which gives (ii).

(iii) The map $v\mapsto[\one,v]$ is injective, because $A_k=\one$ only for $k=1$; it is
continuous and closed, because the quotient map of $O(n)\times B_V(\tfrac12\ell)$ by the
compact group $I$ is closed; and $\Xi$ is a homeomorphism onto $\pr^{-1}(O)$.
\end{proof}

\subsubsection{Correcting Sections to the Cocycle}\label{es:sec:sections}

We choose sections that converge to the limit section of \cref{es:lem:sigmahat}
(\cref{es:lem:sections}) and correct them by a coboundary, an average as in
\cref{es:sec:cocycle}, so that the deck transformations act on them through the cocycle
$\varkappa_i$ (\cref{es:lem:averagedcoboundary,es:lem:sectioncorrection}). The correction
preserves the convergence.

\begin{lemma}[Sections That Converge]\label{es:lem:sections}
Let $F=F^{(i)}$ be equivariant realizing maps. For all large $i$ there is a continuous map
$s^0_i\colon D^i\to W^F_i$ with
$\pr\circ s^0_i=\pi_i$ and with the following property. Let $\tilde e_i(p)$ be the
frame over $p$ that $\tau_i$ maps to $s^0_i(p)$. If points
$p_i\in D^i$ converge to a point $p$, then $p\in \hat U_S$, and the frames $\tilde e_i(p_i)$
converge to $\hat\sigma(p)$.
\end{lemma}
\begin{proof}
\emph{The limits.} Let points $p_i\in D^i$, for $i$ in an infinite set of indices,
converge to a point $p$. By the estimate in the proof of \cref{es:lem:domain}\textup{(i)}
we have $|h_i(p_i)|\le216\,r_\ell$, so $|h(p)|\le216\,r_\ell$ and $p\in \cC_{1/2}$; moreover
$|P_Th(p)|\le\tfrac14\ell$ and $|\xi(p)|\le r_V$. Hence $p\in U_S$, and $p=\exp(\zeta)q$
with $q\in S$ and $|\zeta|\le\ell/(2c_z)$, by \cref{es:lem:limitsection}\textup{(ii)}. Let
$\cD$ be the set of all such limits. We define a continuous map
$\mathsf A\colon \hat U_S\to O(n)$ by $\ft(p,\mathsf A(p))=\hat\sigma(p)$;
its domain contains $\cD$.

\emph{Frames in the reduction.} We claim: if points $p_i\in D^i$, for $i$ in an
infinite set of indices, converge to $p\in\cD$, then for large $i$ there are matrices
$A'_i$ with $\tau_i(\ft_i(p_i,A'_i))\in W^F_i$ and $A'_i\to\mathsf A(p)$. Indeed, by
\textup{\ref{fl:S:frameconv}}, \textup{\ref{fl:S:Yconv}} there are frames $p'_i$ over $p_i$ that converge
to $\hat\sigma(p)$, the points $y_i:=\tau_i(p'_i)$ converge to
$\tau(\hat\sigma(p))=\Xi([\one,\xi(p)])$, and so do the points $F(y_i)$, because the maps
$F$ are realizing maps. The image of $\Xi$ is open, and
$O(n)\times B_V(\tfrac12\ell)\to O(n)\times_IB_V(\tfrac12\ell)$ is a principal bundle; so
for large $i$ we have $F(y_i)=\Xi([A_i,v_i])$ with matrices $A_i\to\one$ and vectors
$v_i\in B_V(\tfrac12\ell)$. The point $A_i^{-1}\cdot y_i=\tau_i(p'_i\cdot A_i)$ lies in
$W^F_i$, the frames $p'_i\cdot A_i$ converge to $\hat\sigma(p)=\ft(p,\mathsf A(p))$, and
their coordinates $(p_i,A'_i)$ in the frame trivializations $\ft_i$ converge to $(p,\mathsf A(p))$.

\emph{The reduction in the frame trivialization.} For $p\in D^i$ let $\mathfrak B_i(p)$ be the
set of the matrices $A$ with $\tau_i(\ft_i(p,A))\in W^F_i$. By the claim, applied to
convergent subsequences, it is not empty for any $p\in D^i$ when $i$ is large; in
particular $\cW_i\subset\pr(W^F_i)$. By \cref{es:sec:reduction} the set $\mathfrak B_i(p)$ is
then one coset $A\,I$, and the set of the pairs $(p,A)$ with $A\in\mathfrak B_i(p)$ is
closed in $D^i\times O(n)$; so $p\mapsto\mathfrak B_i(p)$ is a continuous map of $D^i$ to
$O(n)/I$. By the claim again, if $p_i\in D^i$ converge to $p\in\cD$, then the cosets
$\mathfrak B_i(p_i)$ converge to $\mathsf A(p)\,I$.

\emph{The section.} Fix a bi-invariant metric on $O(n)$. There is $r_I>0$ such that
every $A\in O(n)$ whose distance from a coset $C=A'I$ is less than $r_I$ has a unique
nearest point $\mathrm{n}(C,A)$ in $C$, and $\mathrm{n}$ is continuous; this holds for
the coset $I$, a compact submanifold, and then for all cosets by the invariance of the
metric. Let $\psi_i:=h^{-1}\circ h_i\colon \cC^i_1\to \cC_1$; it is continuous, and
$\psi_i(p_i)\to p$ if $p_i\to p$. For $p\in D^i$ the point $\psi_i(p)$ has the same
coordinates as $p$, so $|h(\psi_i(p))|\le216r_\ell<\ell/2$,
$|P_Th(\psi_i(p))|<\tfrac14\ell$ and
$|\xi(\psi_i(p))|\le|P_Vh_i(p)|+\tfrac3{10}r_\ell<9r_\ell$, and $\psi_i(p)\in \hat U_S$ by
\cref{es:lem:slice}\textup{(iii)} and the proof of \cref{es:lem:limitsection}\textup{(ii)}. The
distance of $\mathsf A(\psi_i(p))$ from $\mathfrak B_i(p)$ tends to $0$ uniformly for
$p\in D^i$: otherwise there are points $p_i\in D^i$ at which it stays bounded from
below; a subsequence converges to a point $p\in\cD$, the points $\psi_i(p_i)$ converge
to $p$, and $\mathsf A(\psi_i(p_i))\to\mathsf A(p)$, which contradicts the convergence
of the cosets. So for large $i$
\begin{equation*}
\mathsf A_i(p):=\mathrm{n}\bigl(\mathfrak B_i(p),\mathsf A(\psi_i(p))\bigr),\qquad
\tilde e_i(p):=\ft_i\bigl(p,\mathsf A_i(p)\bigr),\qquad s^0_i(p):=\tau_i\bigl(\tilde e_i(p)\bigr)
\end{equation*}
are defined and continuous on $D^i$, the point $s^0_i(p)$ lies in $W^F_i$ over
$\pi_i(p)$, and $\mathsf A_i(p_i)\to\mathsf A(p)$ if $p_i\to p$; then
$\tilde e_i(p_i)\to\ft(p,\mathsf A(p))=\hat\sigma(p)$.
\end{proof}

Let $F$ be equivariant realizing maps, and let $s$ be a continuous map of $D^i$ to $W^F_i$
with $\pr\circ s=\pi_i$. For $\gamma\in\Gamma_i$ and $p\in D^i$ with $\gamma p\in D^i$ the
points $s(\gamma p)$ and $s(p)$ of $W^F_i$ lie over the same point of $X_i$, so there is a unique $c^s(\gamma,p)\in I$ with
\begin{equation}\label{es:eq:sectioncocycle}
s(\gamma p)=c^s(\gamma,p)\cdot s(p),
\end{equation}
by \cref{es:sec:reduction}. It is continuous in $p$, and
$c^s(\gamma\gamma',p)=c^s(\gamma,\gamma'p)\,c^s(\gamma',p)$ when $p$, $\gamma'p$ and
$\gamma\gamma'p$ lie in $D^i$. For the frames $\tilde e(p)$ over $p$ with
$\tau_i(\tilde e(p))=s(p)$ the relation \eqref{es:eq:sectioncocycle} says that
\begin{equation}\label{es:eq:framecocycle}
\gamma\cdot\tilde e(p)=\tilde e(\gamma p)\cdot c^s(\gamma,p),
\end{equation}
because $\tau_i$ is injective on the frames over $\gamma p$, the frames
$\gamma\cdot\tilde e(p)$ and $\tilde e(p)$ have the same image under $\tau_i$, and
$A\cdot y=y\cdot A^{-1}$.

\begin{lemma}[Averaged Coboundaries]\label{es:lem:averagedcoboundary}
Let a discrete group $\Gamma$ act on a space, let $D$ be an open set, and let $\mathcal P$ be
the set of the pairs $(\gamma,p)\in\Gamma\times D$ with $\gamma p\in D$. Let $\mathbb E$ be a
finite-dimensional Euclidean vector space, let $\alpha_\gamma$, for the elements $\gamma$
that occur in $\mathcal P$, be orthogonal maps of $\mathbb E$, and let $b\colon\mathcal P\to \mathbb E$ be a
bounded map that is continuous in $p$ for fixed $\gamma$. Suppose that
$\alpha_{\gamma\gamma'}=\alpha_\gamma\alpha_{\gamma'}$ and
\begin{equation}\label{es:eq:bcocycle}
b(\gamma\gamma',p)=b(\gamma,\gamma'p)+\alpha_\gamma b(\gamma',p)
\end{equation}
whenever $(\gamma',p)$ and $(\gamma,\gamma'p)$ lie in $\mathcal P$. Suppose that continuous
non-negative weights $\varpi_\nu$ on $D$ with locally finite supports sum to one, vanish at
the points $p$ with $\nu^{-1}p\notin D$, and satisfy $\varpi_\nu(\gamma p)=\varpi_{\gamma^{-1}\nu}(p)$
for $(\gamma,p)\in\mathcal P$, and that the functions $p\mapsto\varpi_\nu(p)\,b(\nu,\nu^{-1}p)$,
extended by zero, are continuous on $D$.
Then the continuous map
\[
g(p):=\sum_\nu\varpi_\nu(p)b(\nu,\nu^{-1}p)
\]
satisfies
\begin{equation*}
b(\gamma,p)=g(\gamma p)-\alpha_\gamma g(p)\quad\text{for }(\gamma,p)\in\mathcal P,\qquad
\|g\|_\infty\le\|b\|_\infty .
\end{equation*}
\end{lemma}
\begin{proof}
This is the degree-one averaging identity; compare \cite[Proposition~1]{Cra03}.
Continuity follows from local finiteness and the zero-extension hypothesis, and
the norm bound follows from positivity of the weights.
For $p,\gamma p\in D$, reindex by $\nu=\gamma\mu$. The cocycle law on the
three points $\mu^{-1}p,p,\gamma p$ gives
\[
b(\gamma\mu,\mu^{-1}p)
 =b(\gamma,p)+\alpha_\gamma b(\mu,\mu^{-1}p).
\]
Multiplication by $\varpi_\mu(p)$ and summation proves the identity.
\end{proof}

Applied to the exponents of the cocycle of a section that converges to the limit section, the
averaged coboundary corrects the section to one with the exact cocycle.

\begin{lemma}[Correction of a Section]\label{es:lem:sectioncorrection}
Let $s^0_i\colon D^i\to W^F_i$ be a continuous section over $\pi_i$, where
$F$ is an equivariant realizing map. Suppose that its cocycle is
\begin{equation}\label{es:eq:bdef}
c^{s^0_i}(\gamma,p)
 =A_{\exp(\mathsf b_i(\gamma,p))\,\varkappa_i(\gamma,p)},\qquad
|\mathsf b_i(\gamma,p)|<\tfrac12r_K .
\end{equation}
Then there is a continuous map $\mathsf g_i\colon D^i\to\fk$ with
$\|\mathsf g_i\|_\infty\le\sup|\mathsf b_i|$ such that
$s_i:=A_{\exp(-\mathsf g_i)}\cdot s^0_i$ satisfies
\begin{equation*}
s_i(\gamma p)=A_{\varkappa_i(\gamma,p)}\cdot s_i(p)
\qquad(p,\gamma p\in D^i).
\end{equation*}
In particular, if $\sup|\mathsf b_i|\to0$, the correction preserves convergence
of the associated frames.
\end{lemma}
\begin{proof}
The logarithm in \eqref{es:eq:bdef} is unique and continuous, since
$k\mapsto A_k$ is injective and $\exp$ is injective on $\fk(2r_K)$.
The two multiplicative cocycle laws, the abelianness of $\fk$, and
$\Ad(\varkappa_i(\gamma,p))|_\fk=\Ad_\gamma$ give
\[
\exp\bigl(\mathsf b_i(\gamma\gamma',p)\bigr)
=\exp\bigl(\mathsf b_i(\gamma,\gamma'p)
       +\Ad_\gamma\mathsf b_i(\gamma',p)\bigr).
\]
Both exponents have norm less than $r_K$, so this is the additive law
\eqref{es:eq:bcocycle}. The maps $\Ad_\gamma$ satisfy $\Ad_{\gamma\gamma'}=\Ad_\gamma\Ad_{\gamma'}$ for these
pairs, by \cref{es:prop:cocycle}\textup{(ii),(iii)}.

We construct the weights. Choose a continuous function
$\chi_1\colon[0,\infty)\to[0,1]$ equal to one on $[0,1/8]$ and zero on
$[3/16,\infty)$, and put $\chi'(p):=\chi_1(|P_Th_i(p)|/\ell)$ on $D^i$.
For $p\in D^i$ set
\[
\bar\chi'(p):=\sum_{\nu:\,\nu^{-1}p\in D^i}\chi'(\nu^{-1}p),\qquad
\varpi'_\nu(p):=
 \begin{cases}\chi'(\nu^{-1}p)/\bar\chi'(p),&\nu^{-1}p\in D^i,\\
 0,&\text{otherwise}.
 \end{cases}
\]
By \cref{es:lem:domain}\textup{(iii)} an orbit representative lies in
$D^i\cap \cC^i_{1/16}$, where $\chi'=1$, so $\bar\chi'\ge1$.
The set $D^i$ lies in a compact subset of $\cC^i_{1/2}$, by the bound in the proof of
\cref{es:lem:domain}\textup{(i)}; so the sums are locally finite, by proper
discontinuity. For continuity, suppose $p_j\to p$ in $D^i$ and $\nu^{-1}p_j\in D^i$ with
$\chi'(\nu^{-1}p_j)>0$. The points $\nu^{-1}p_j$ converge to $\nu^{-1}p$, which lies in
this compact subset of $\cC^i_{1/2}$ and satisfies
\[
|P_Th_i(\nu^{-1}p)|\le\tfrac3{16}\ell,\qquad
|\tilde\xi_i(\nu^{-1}p)|=|\tilde\xi_i(p)|<r_V,
\]
hence lies in $D^i$. Thus the functions $\chi'(\nu^{-1}\cdot)$, extended by zero, are
continuous on $D^i$, and so are their products with the bounded functions
$p\mapsto\mathsf b_i(\nu,\nu^{-1}p)$.
By \cref{th:lem:normalizeweights}, the normalized weights are continuous,
sum to one and satisfy $\varpi'_\nu(\gamma p)=\varpi'_{\gamma^{-1}\nu}(p)$.
Apply \cref{es:lem:averagedcoboundary} with $\alpha_\gamma:=\Ad_\gamma$, $b=\mathsf b_i$ and these weights.
Finally
\[
s_i(\gamma p)
=A_{\exp(-\mathsf g_i(\gamma p)+\mathsf b_i(\gamma,p)
                  +\Ad_\gamma\mathsf g_i(p))\,\varkappa_i(\gamma,p)}
 \cdot s_i(p)
=A_{\varkappa_i(\gamma,p)}\cdot s_i(p).
\]
If $\sup|\mathsf b_i|\to0$, then $\exp(-\mathsf g_i)\to1$ uniformly, proving
the last assertion.
\end{proof}

\subsubsection{Admissible Sections}\label{es:sec:admissible}

We keep the tube and the fixed harmonic-chart subsequence, and describe the
sections of the reduction of an equivariant realizing map that are compatible with the cocycle $\varkappa_i$ and the limit section. Let $F=F^{(i)}$
be equivariant realizing maps, with the sets $W^F_i$ and the maps $\varphi^F$ of
\cref{es:sec:reduction}.

\begin{definition}[Admissible Sections and Their Coordinates]\label{es:def:admissible}
An \emph{admissible section} of $F$ in the tube
consists of continuous maps $s_i\colon D^i\to W^F_i$ with $\pr\circ s_i=\pi_i$, for all
large $i$, with the following two properties.
\begin{enumerate}[label=\textup{(A\arabic*)},ref=\textup{(A\arabic*)},leftmargin=3.2em]
\item\label{es:A:cocycle} $s_i(\gamma p)=A_{\varkappa_i(\gamma,p)}\cdot s_i(p)$ whenever
$\gamma\in\Gamma_i$ and $p,\gamma p\in D^i$.
\item\label{es:A:frames} If points $p_i\in D^i$ converge to a point $p$, then
$p\in \hat U_S$, and the frames $\tilde e_i(p_i)$ over $p_i$ that $\tau_i$ maps to
$s_i(p_i)$ converge to $\hat\sigma(p)$.
\end{enumerate}
Its \emph{coordinate} is the map $w_i:=\varphi^F\circ s_i\colon D^i\to B_V(\tfrac12\ell)$.
\end{definition}

Admissible sections exist and are unique up to small equivariant corrections, and their
coordinates are equivariant and uniformly close to $\tilde\xi_i$.

\begin{lemma}[Admissible Sections]\label{es:lem:admissible}\label{es:lem:rightsection}
Let $F=F^{(i)}$ be equivariant realizing maps.
\begin{enumerate}[label=\textup{(\roman*)}]
\item There is an admissible section of $F$ in the tube.
\item If $s_i$ and $s'_i$ are admissible sections of $F$, then for large $i$ there are
continuous maps $u_i\colon D^i\to\fk(r_K)$ with
\begin{equation}\label{es:eq:twoadmissible}
s'_i=A_{\exp(u_i)}\cdot s_i,\qquad u_i(\gamma p)=\Ad_\gamma\,u_i(p)\quad
\text{for }p,\gamma p\in D^i,\qquad\sup_{D^i}|u_i|\to0 .
\end{equation}
Conversely, if $s_i$ is admissible and $u_i$ are continuous maps with the last two
properties, then $A_{\exp(u_i)}\cdot s_i$ is admissible.
\item The coordinate of an admissible section satisfies
$w_i(\gamma p)=\mathsf R(\varkappa_i(\gamma,p))\,w_i(p)$ for $p,\gamma p\in D^i$, and the
supremum of $|w_i-\tilde\xi_i|$ over $D^i$ tends to $0$.
\end{enumerate}
\end{lemma}
\begin{proof}
(i) By \cref{es:lem:sections} there are continuous maps $s^0_i\colon D^i\to W^F_i$ over
$\pi_i$ whose frames $\tilde e^0_i$ satisfy \ref{es:A:frames}.

Let $c_i:=c^{s^0_i}$ be their cocycle, as in \eqref{es:eq:sectioncocycle}. We claim that for
large $i$ and all $\gamma$ and $p$ with $p,\gamma p\in D^i$ there is a unique
$\mathsf b_i(\gamma,p)\in\fk(\tfrac12r_K)$ with
\begin{equation*}
c_i(\gamma,p)=A_k\quad\text{for}\quad k=\exp\bigl(\mathsf b_i(\gamma,p)\bigr)\,\varkappa_i(\gamma,p),
\qquad\text{and}\qquad\sup|\mathsf b_i|\to0 .
\end{equation*}
Otherwise there are $\eps>0$ and, for infinitely many $i$, elements $\gamma_i$ and
points $p_i$ for which no such vector of norm less than $\eps$ exists. After passing to
a subsequence, $\gamma_i\eGHto g$, $p_i\to p$ and $\gamma_ip_i\to gp$, with
$p,gp\in \hat U_S$. By \textup{\ref{fl:S:frameconv}} and \cref{es:lem:sigmahat}\textup{(ii)} the
frames $\gamma_i\cdot\tilde e^0_i(p_i)$ converge to
$g\cdot\hat\sigma(p)=\hat\sigma(gp)\cdot A_{\varkappa_\infty(g,p)}$, and the frames
$\tilde e^0_i(\gamma_ip_i)$ converge to $\hat\sigma(gp)$. In the frame
trivializations, \eqref{es:eq:framecocycle} expresses the second coordinate of
$\gamma_i\cdot\tilde e^0_i(p_i)$ as that of $\tilde e^0_i(\gamma_ip_i)$ multiplied
on the right by $c_i(\gamma_i,p_i)$. Taking the quotient of these two convergent
matrices gives $c_i(\gamma_i,p_i)\to A_{\varkappa_\infty(g,p)}$. Since
$\varkappa_i(\gamma_i,p_i)\to\varkappa_\infty(g,p)$ by \cref{es:lem:limitcocycle}, the elements
$k_i$ with $A_{k_i}=c_i(\gamma_i,p_i)$ satisfy $k_i\varkappa_i(\gamma_i,p_i)^{-1}\to1$, a
contradiction.

Apply \cref{es:lem:sectioncorrection} to $s^0_i$ and $\mathsf b_i$. It gives
$\|\mathsf g_i\|_\infty\le\sup|\mathsf b_i|\to0$ and the section \mbox{$s_i:=A_{\exp(-\mathsf g_i)}\cdot s^0_i$} with the property \ref{es:A:cocycle}. Its frames
are those of $s^0_i$, multiplied by matrices that tend to $\one$ uniformly; so \ref{es:A:frames}
holds.

(ii) Two points of $W^F_i$ over one point of $X_i$ differ by an element of $I$, so
$s'_i(p)=A_{m_i(p)}\cdot s_i(p)$ with a continuous map $m_i\colon D^i\to K$. If
$p_i\to p$, the frames of both sections at $p_i$ converge to $\hat\sigma(p)$, and in
the frame trivializations of \cref{es:sec:reduction} their second coordinates differ by the
factor $A_{m_i(p_i)}^{-1}$; so $A_{m_i(p_i)}\to\one$. Hence $m_i\to1$ uniformly, and for large $i$ there is a continuous map $u_i$ to $\fk(r_K)$
with $m_i=\exp(u_i)$ and $\sup|u_i|\to0$. By \ref{es:A:cocycle} for both sections,
$m_i(\gamma p)\,\varkappa_i(\gamma,p)=\varkappa_i(\gamma,p)\,m_i(p)$, that is,
$\exp(u_i(\gamma p))=\exp(\Ad_\gamma u_i(p))$, by
\cref{es:prop:cocycle}\textup{(iii)}; since $\exp$ is injective on $\fk(2r_K)$, this gives
the second relation of \eqref{es:eq:twoadmissible}. The converse follows from the same
computations.

(iii) The map $\varphi^F$ commutes with $I$, which gives the first assertion. Suppose
that there are a positive number $\eps$ and, for infinitely many $i$, points
$p_i\in D^i$ at which the vectors $w_i$ and $\tilde\xi_i$ have distance at least $\eps$; after passing to a
subsequence, $p_i\to p\in \hat U_S$. By \ref{es:A:frames}, \textup{\ref{fl:S:Yconv}} and
\cref{es:lem:sigmahat}\textup{(i)} the points $s_i(p_i)$ converge to $\Xi([\one,\xi(p)])$.
Since the maps $F$ are realizing, the points $F(s_i(p_i))=\Xi([\one,w_i(p_i)])$ converge to
the same point, and $w_i(p_i)\to\xi(p)$ by \cref{es:lem:sigmahat}\textup{(iii)}. But
$\tilde\xi_i(p_i)\to\xi(p)$ by \cref{es:lem:limit}, a contradiction.
\end{proof}

\subsection{Equivariant Triviality}\label{es:sec:product}

In the theory with bounded curvature the local model is $\mathsf F\times_KV$: the bundle over the
ball of the slice representation is a product, compatibly with the isotropy group. In
this subsection we prove that a proper map over a ball of a representation is such a
product if it has product charts, in which the group is a factor and the map is a
projection twisted by the group (\cref{es:prop:trivial}); the gluing of
\cref{ch:globalbundle} applies this to the glued maps over each tube. The proof is
topological, and it proceeds by induction on the dimension of the part of $V$ on which
$K$ acts non-trivially: a blow-up of the ball along the fixed subspace reduces the
statement to slices, on which a smaller group acts.

\subsubsection{Product Charts}\label{es:sec:productcharts}

We fix the class of maps for which the induction works. Let $K$ be a compact Lie group
and $B$ a smooth manifold with a smooth left action of $K$. Let $Q$ be a second
countable Hausdorff space with a free right action of $K$, and let
$\varphi\colon Q\to B$ be a continuous map with $\varphi(q\cdot k)=k^{-1}\varphi(q)$.

\begin{definition}[Product Charts]\label{es:def:productchart}
A \emph{product chart} of $\varphi$ is a homeomorphism $c$ of $U\times K$ onto an open
subset of $Q$, where $U$ is a topological manifold, such that
\begin{equation}\label{es:eq:productchart}
c(u,kk')=c(u,k)\cdot k',\qquad\varphi\bigl(c(u,k)\bigr)=k^{-1}\beta(u)
\end{equation}
for a continuous map $\beta\colon U\to B$ with the following property: every point of
$U$ has an open neighbourhood that is homeomorphic to a product $U'\times U''$ of
topological manifolds in such a way that $\beta$ is the projection to $U''$, followed by
a homeomorphism of $U''$ onto an open subset of $B$. We say that $\varphi$ \emph{has
product charts} if the images of its product charts cover $Q$.
\end{definition}

If $\varphi$ has product charts, then $Q$ is a topological manifold, $\varphi$ is a topological
submersion, and $Q\to Q/K$ is a principal bundle.

Two operations preserve product charts.

\begin{lemma}[Pull-Backs and Slices]\label{es:lem:pullslice}
Let $\varphi\colon Q\to B$ have product charts.
\begin{enumerate}[label=\textup{(\roman*)},beginpenalty=10000]
\item Let $\tilde B$ be a smooth manifold with a smooth action of $K$, and let
$\mathsf c\colon\tilde B\to B$ be a continuous map that commutes with $K$. Then the
projection of
\begin{equation*}
\mathsf c^*Q:=\bigl\{(q,b)\in Q\times\tilde B:\varphi(q)=\mathsf c(b)\bigr\},\qquad
(q,b)\cdot k:=(q\cdot k,k^{-1}b),
\end{equation*}
to $\tilde B$ has product charts.
\item Let $b\in B$, let $H$ be its isotropy group, and let $\Sigma$ be a slice at $b$
(\cref{ls:def:slice}) that is a smooth submanifold. Then the restriction of $\varphi$ to
$Q_\Sigma:=\varphi^{-1}(\Sigma)$, with the action of $H$, has product charts, and
$[q,k]\mapsto q\cdot k$ is a homeomorphism of $Q_\Sigma\times_HK$ onto
$\varphi^{-1}(K\Sigma)$.
\end{enumerate}
\end{lemma}
\begin{proof}
(i) Let $c$ be a product chart with the map $\beta$, and let
$\hat U:=\{(u,b)\in U\times\tilde B:\beta(u)=\mathsf c(b)\}$. Where $\beta$ is the
projection of $U'\times U''$ to an open subset $U''$ of $B$, the set $\hat U$ is
$U'\times\mathsf c^{-1}(U'')$; so $\hat U$ is a topological manifold, and its projection
$\hat\beta$ to $\tilde B$ has the property required of $\beta$. The map
$\hat c((u,b),k):=(c(u,k),k^{-1}b)$ has values in $\mathsf c^*Q$, by
\eqref{es:eq:productchart}; it is a homeomorphism onto the preimage of the image of $c$,
with the inverse $(q,b')\mapsto((u,kb'),k)$ for $q=c(u,k)$, and it satisfies
\eqref{es:eq:productchart} for $\hat\beta$.

(ii) Let $\mathcal O\subset K/H$ be an open neighbourhood of the coset $H$ on which
$K\to K/H$ has a continuous section $\sigma_H$ with $\sigma_H(H)=1$, and let $\pi_H\colon K\Sigma\to K/H$ be the projection induced by the tube $K\times_H\Sigma\to K\Sigma$. The map
$(o,x)\mapsto\sigma_H(o)x$ is a homeomorphism of $\mathcal O\times\Sigma$ onto the
open set $\pi_H^{-1}(\mathcal O)$. Let $c$ be a product chart with the map $\beta$, and
let $U_{\mathcal O}$ be the open set of the points $u\in U$ with
$\beta(u)\in\pi_H^{-1}(\mathcal O)$. For $u\in U_{\mathcal O}$ and $k\in K$ the point
$c(u,k)$ lies in $Q_\Sigma$ if and only if $\beta(u)\in k\Sigma$, that is, if and only
if $k\in\sigma_H(\pi_H(\beta(u)))\,H$. So
\begin{equation*}
c_\Sigma(u,h):=c\bigl(u,\sigma_H(\pi_H(\beta(u)))\,h\bigr)\qquad(u\in U_{\mathcal O},\ h\in H)
\end{equation*}
is a homeomorphism of $U_{\mathcal O}\times H$ onto an open subset of $Q_\Sigma$, and
\eqref{es:eq:productchart} holds for $c_\Sigma$, the group $H$ and the map
$\beta_\Sigma(u):=\sigma_H(\pi_H(\beta(u)))^{-1}\beta(u)$ to $\Sigma$. Where $\beta$ is
the projection of $U'\times U''$ to $U''\subset\pi_H^{-1}(\mathcal O)$, and $U''$
corresponds to a product $\mathcal O_1\times\Sigma_1$ of open subsets of $\mathcal O$ and
$\Sigma$, the map $\beta_\Sigma$ is the projection of $(U'\times\mathcal O_1)\times\Sigma_1$ to $\Sigma_1$. Every point of $Q_\Sigma$ lies in the image of such a chart, after a
translation of $c$ by an element of $K$. The last assertion holds because
$\varphi^{-1}(K\Sigma)\to K\Sigma\to K/H$ is a map that commutes with $K$ and $Q_\Sigma$ is
the preimage of the coset $H$ \cite[Chapter~II, Section~3]{Bre72}.
\end{proof}

\subsubsection{Homotopies of the Base}\label{es:sec:cht}

The next lemma is the covering homotopy theorem in the form that we use. Without the group
action it is \cite[Chapter~4, Theorems~9.6 and~9.8]{Hus94}, stated there for numerable
principal bundles; the proof uses only local product structures, and it applies when the
sets and functions are chosen invariant. We indicate the steps.

\begin{lemma}[Products along an Interval]\label{es:lem:cht}
Let $M$ be a smooth manifold with a smooth action of $K$, let $J\subset\R$ be an
interval that contains $0$, and let $E$ be a Hausdorff space with a free right action
of $K$ and a continuous map $\varphi\colon E\to M\times J$ with
$\varphi(e\cdot k)=k^{-1}\varphi(e)$, where $K$ acts trivially on $J$. Assume that every
point $(x,t)$ has the following property: there are a slice $\Sigma$ at $x$ for the
isotropy group $H$ of $x$, an interval $J_t\subset J$ that is a neighbourhood of $t$ in
$J$, a space $\mathsf F'$ with a right action of $H$, and a homeomorphism of
$\varphi^{-1}(\Sigma\times J_t)$ onto $\mathsf F'\times\Sigma\times J_t$ over $\Sigma\times J_t$
that commutes with $H$, for the action $(f,x',t')\cdot h:=(f\cdot h,h^{-1}x',t')$. Then
there is a homeomorphism of $E$ onto $E_0\times J$ over $M\times J$ that commutes with
$K$ and is the identity over $M\times\{0\}$, where $E_0:=\varphi^{-1}(M\times\{0\})$.
\end{lemma}
\begin{proof}
Over $K\Sigma\times J_t$ the space $E$ is $(\mathsf F'\times\Sigma\times J_t)\times_HK$, by the
argument for the last assertion of \cref{es:lem:pullslice}\textup{(ii)}. So $E$ is a
$K$-bundle over $M\times J$ that is locally trivial, in the equivariant sense, over the
tubes $K\Sigma\times J_t$, and these tubes form a numerable cover, because $M/K$ is
paracompact. The assertion is the covering homotopy theorem for numerable bundles over the
product of a space with an interval \cite[Chapter~4, Theorem~9.8]{Hus94}, in its equivariant
form, which is proved in the same way with invariant sets and invariant functions: Lashof
\cite{Las82}, and Bierstone \cite{Bie73} for smooth bundles. In the terms of
\cite{Hus94}, the product structures over the sets $U_j\times[a,b]$ of a locally finite
invariant cover define maps $R_j$ of $\varphi^{-1}(M\times[a,b])$ that commute with $K$, are
the identity outside the support of an invariant function $u_j$, where $\max_ju_j=1$, and
move the fibre over $(x,t)$ homeomorphically onto the fibre over
$(x,\min\{t,b-u_j(x)(b-a)\})$; the composition of the maps $R_j$ moves every fibre over
$(x,t)$ onto the fibre over $(x,a)$, and the products over the intervals of a decomposition
of $J$ are composed, starting from $0$.
\end{proof}

\subsubsection{Triviality over a Ball}\label{es:sec:trivialball}

The proof is the induction announced above. If the group acts trivially on the ball, the
induced map of the orbit space is a proper submersion over a contractible base, hence a
product. Otherwise a blow-up along the fixed subspace and the covering homotopy theorem of
\cref{es:lem:cht} reduce the statement to slices, on which the part of $V$ with a
non-trivial action has smaller dimension.

\begin{proposition}[Equivariant Triviality]\label{es:prop:trivial}
Let $K$ be a compact Lie group, let $V$ be a Euclidean space with an orthogonal
representation of $K$, let $B\subset V$ be an open ball about $0$, with the action of
$K$, and let $\varphi\colon Q\to B$ be a proper map with product charts. Then
$\mathsf F:=\varphi^{-1}(0)$ is a closed manifold on which $K$ acts freely, and there is a
homeomorphism
\begin{equation*}
\Phi\colon Q\longrightarrow \mathsf F\times B,\qquad\Phi(q\cdot k)=\Phi(q)\cdot k\quad\text{for}\quad
(f,v)\cdot k:=(f\cdot k,k^{-1}v),
\end{equation*}
whose second component is $\varphi$ and with $\Phi(f)=(f,0)$ for $f\in \mathsf F$.
\end{proposition}
\begin{proof}
The set $\mathsf F$ is compact and invariant, and in a product chart it is
$\beta^{-1}(0)\times K$; so it is a closed manifold. Let $V_0$ be the subspace of the
vectors that are fixed by $K$ and $V'$ its orthogonal complement. We use induction on
$m:=\dim V'$, for all groups and representations at once.

\emph{The case $m=0$.} Here $K$ acts trivially on $B$, and $\varphi$ is constant on the
orbits. The space $Q/K$ is a Hausdorff space with a countable base, and by
\eqref{es:eq:productchart} it is locally homeomorphic to the manifolds $U$, in such a way
that the map induced by $\varphi$ corresponds to $\beta$. So $Q/K$ is a topological
manifold, and the induced map to $B$ is a proper topological submersion with compact
manifold fibres, hence a fibre bundle \cite[Corollary~6.14]{Sie72}; manifolds satisfy the
hypotheses of that corollary on the fibres by \cite[Example~1.3 and Theorem~2.3]{Sie72}. Since $B$ is contractible, there is a homeomorphism of
$Q/K$ onto $(\mathsf F/K)\times B$ over $B$. The principal bundle $Q\to(\mathsf F/K)\times B$ is isomorphic to the product of its restriction $\mathsf F\to\mathsf F/K$ to $(\mathsf F/K)\times\{0\}$ with $B$, by an isomorphism over $(\mathsf F/K)\times B$ that is the identity over $(\mathsf F/K)\times\{0\}$
\cite[Chapter~4, Section~9]{Hus94}. This is the map $\Phi$.

\emph{The step.} Let $m\ge1$, and assume the assertion for all smaller values. Let $R$ be the radius of $B$. The ball $B$ is homeomorphic to the product of the balls $B_0$
and $B'$ of radius $R$ in $V_0$ and in $V'$, by a map that commutes with $K$ and is the
identity on $V_0$ and on $V'$; by \cref{es:lem:pullslice}\textup{(i)} we may assume that
$B=B_0\times B'$. Let $\mathbb S(V')$ be the unit sphere of $V'$, and let
\begin{equation*}
M:=B_0\times\mathbb S(V'),\qquad J:=(-R,R),\qquad
\mathsf c\colon M\times J\to B,\quad\mathsf c(v_0,u,s):=(v_0,su).
\end{equation*}
The map $\mathsf c$ is smooth and commutes with $K$, so $E:=\mathsf c^*Q$ and its
projection to $M\times J$ have product charts, by \cref{es:lem:pullslice}\textup{(i)}, and
the projection is proper. We verify the hypothesis of \cref{es:lem:cht}. Let
$(v_0,u,t)\in M\times J$, and let $H$ be the isotropy group of $u$. By the slice theorem
for the smooth action of $K$ on the sphere \cite[Chapter~VI, Section~2]{Bre72} there is
a slice at $u$ in $\mathbb S(V')$ that is, by a diffeomorphism that commutes with $H$, an
open ball about $0$ in the subspace $W_u$ of $V'$ that is orthogonal to $u$ and to the
orbit of $u$. Its product with a small ball about $v_0$ in $B_0$ is a slice $\Sigma$ at
$(v_0,u)$ in $M$, and $\Sigma\times J_t$, for a small open interval $J_t$ about $t$, is
a slice at $(v_0,u,t)$ in $M\times J$. It is an open ball about $0$ in the
representation $V_0\oplus W_u\oplus\R$ of $H$, up to a homeomorphism that commutes with
$H$ and to which \cref{es:lem:pullslice}\textup{(i)} applies, and the part of this space on which $H$ acts non-trivially lies in $W_u$, so its
dimension is at most $m-1$. By \cref{es:lem:pullslice} the restriction of $E$ to
$\Sigma\times J_t$ has product charts for the group $H$, and it is proper. By the
assumption of the induction it is homeomorphic to $\mathsf F'\times\Sigma\times J_t$ over
$\Sigma\times J_t$, compatibly with $H$, where $\mathsf F'$ is the fibre over $(v_0,u,t)$.

So \cref{es:lem:cht} gives a homeomorphism of $E$ onto $E_0\times J$ over $M\times J$ that
commutes with $K$ and is the identity over $M\times\{0\}$. The space $E_0$ is the set
of the pairs $(q,(v_0,u,0))$ with $\varphi(q)=(v_0,0)$, that is, the product of
$Q_0:=\varphi^{-1}(B_0\times\{0\})$ and $\mathbb S(V')$ over $B_0$. The restriction of $\varphi$ to
$Q_0$ has product charts, by \cref{es:lem:pullslice}\textup{(i)} for the inclusion of
$B_0$, and $K$ acts trivially on $B_0$; by the case $m=0$ there is a homeomorphism of
$Q_0$ onto $\mathsf F\times B_0$ over $B_0$ that commutes with $K$ and is the identity on $\mathsf F$.
Let $\mathfrak r_0\colon Q_0\to \mathsf F$ be its first component. Together we obtain a continuous map
\begin{equation*}
\tilde\Phi\colon E\longrightarrow \mathsf F,\qquad\tilde\Phi\bigl((q,b)\cdot k\bigr)=\tilde\Phi(q,b)\cdot k,
\qquad\tilde\Phi\bigl(q,(v_0,u,0)\bigr)=\mathfrak r_0(q),
\end{equation*}
that maps every fibre of $E\to M\times J$ homeomorphically onto $\mathsf F$. For $q\in Q$ with
$\varphi(q)=(v_0,v')$ let
\begin{equation*}
\mathfrak r(q):=\tilde\Phi\bigl(q,(v_0,v'/|v'|,|v'|)\bigr)\quad\text{if }v'\ne0,\qquad
\mathfrak r(q):=\mathfrak r_0(q)\quad\text{if }v'=0.
\end{equation*}
This map commutes with $K$, and it maps every fibre of $\varphi$ homeomorphically onto
$\mathsf F$. It is continuous: this is clear where $v'\ne0$; and if points $q_j$ converge to a
point $q$ with $v'=0$, a subsequence of the unit vectors $v'_j/|v'_j|$, for the indices
with $v'_j\ne0$, converges to some $u$, the points $(q_j,(v_{0,j},v'_j/|v'_j|,|v'_j|))$
converge to $(q,(v_0,u,0))$ in $E$, and the values of $\tilde\Phi$ converge to
$\mathfrak r_0(q)$. Hence $\Phi:=(\mathfrak r,\varphi)$ is a continuous bijection of $Q$ onto
$\mathsf F\times B$ that commutes with $K$. It is proper, because $\varphi$ is, and the spaces are
locally compact Hausdorff spaces; so it is a homeomorphism, and $\Phi(f)=(f,0)$ on $\mathsf F$
because $\mathfrak r_0$ is the identity there.
\end{proof}

\section{Uniform Tubes and Global Gluing}\label{ch:globalbundle}

We assemble the gluing tools and apply them to the cocycles and admissible sections of
\cref{ch:localbundle} (\cref{es:prop:cocycle,es:def:admissible}). First, with the single-tube data of \cref{es:sec:standing},
we modify a realizing map inside one tube (\cref{es:def:tube}) and give a regularity criterion
for product charts in the sense of \cref{es:def:productchart}
(\cref{es:sec:localextension,es:sec:regularcoordinates}).
We then choose a finite system of tubes with dimension-controlled overlap
(\cref{es:sec:gluing}) and modify the map one tube at a time, transferring its
regular coordinates across overlaps (\cref{es:sec:gluingall}).
The result is \cref{es:thm:gluingall}.

\subsection{Modification in a Tube}\label{es:sec:localextension}

We prepare the gluing by modifying a realizing map inside one tube;
\cref{es:sec:regularcoordinates} then gives a criterion for the resulting map to have
product charts. The modification makes
its admissible coordinate equal to $\tilde\xi_i$ near the centre. In this subsection and
the next we keep the setting of \cref{th:sec:standing}, the tube $\Xi$ of
\cref{es:sec:limitcoordinate} with the sets $D^i(t)$, $\cW_i(t)$ and the function $\chi^\xi_i$
of \cref{es:sec:realizing}, the admissible sections of \cref{es:def:admissible} and the
harmonic-chart subsequence of \cref{th:lem:stagecharts}.

\begin{lemma}[The Modification in a Tube]\label{es:lem:modall}
Let $F=F^{(i)}$ be equivariant realizing maps with an admissible section $s_i$ and its coordinate
$w_i$, and let $\chi^\xi_i=\chi_3(|\tilde\xi_i|/r_\ell)$ be the function of
\cref{es:sec:realizing}. For all large $i$ there is exactly one map $F'\colon Y_i\to Y$ that
commutes with $O(n)$, agrees with $F$ outside $\pr^{-1}(\cW_i)$, and satisfies
\begin{equation*}
F'\bigl(s_i(p)\bigr)=\Xi\bigl([\one,w'_i(p)]\bigr),\qquad
w'_i:=(1-\chi^\xi_i)\,w_i+\chi^\xi_i\,\tilde\xi_i,\qquad\text{for }p\in D^i .
\end{equation*}
It is continuous, it agrees with $F$ outside $\pr^{-1}(\cW_i(4))$, the supremum of
$\sd(F',F)$ over $Y_i$ tends to $0$, and the maps $F'$ are realizing. The maps $s_i$
form an admissible section of $F'$ with the coordinate $w'_i$, which equals
$\tilde\xi_i$ on $D^i(2)$.
\end{lemma}
\begin{proof}
By \cref{es:lem:domain}\textup{(ii)}, \ref{es:X:equiv} and
\cref{es:lem:admissible}\textup{(iii)}, the map $w'_i$ satisfies
\begin{equation*}
w'_i(\gamma p)=\mathsf R(\varkappa_i(\gamma,p))\,w'_i(p)\quad\text{when }p,\gamma p\in D^i,\qquad\text{and}\qquad|w'_i|<9r_\ell<\ell/2
\end{equation*}
for large $i$. By \ref{es:A:cocycle} and freeness of the
$O(n)$-action, $A\cdot s_i(p)=A'\cdot s_i(p')$ if and only if
$p'=\gamma p$ and $A=A'A_{\varkappa_i(\gamma,p)}$ for some $\gamma$. Thus
\[
F'(A\cdot s_i(p)):=\Xi([A,w'_i(p)])
\]
is well defined and uniquely determined on $\pr^{-1}(\cW_i)$. It is continuous:
locally, $\pi_i$ is a homeomorphism and $s_i$ trivializes the principal
$O(n)$-bundle, so $(A,p)\mapsto A\cdot s_i(p)$ is an open map.

Where $\chi^\xi_i=0$ this definition agrees with $F$. The closure of
$\cW_i(4)$ is compact in $\cW_i$: by \cref{es:lem:domain}\textup{(iii)}, a
sequence in $\cW_i(4)$ has lifts in $\cC^i_{1/16}$ with $|\tilde\xi_i|<4r_\ell$;
their limits lie in the compact set
$\overline \cC^i_{1/16}\cap\{|\tilde\xi_i|\le4r_\ell\}\subset D^i$.
Therefore extending by $F$ outside $\pr^{-1}(\cW_i)$ is
continuous. Since
$|w'_i-w_i|\le|\tilde\xi_i-w_i|\to0$ uniformly, uniform continuity of $\Xi$ on the image of \mbox{$O(n)\times\overline B_V(9r_\ell)$} gives
$\sup\sd(F',F)\to0$. In particular $F'$ is realizing. Its definition puts $s_i$
in $W^{F'}_i$ with coordinate $w'_i$; the conditions \ref{es:A:cocycle} and
\ref{es:A:frames} concern the section alone, so it remains admissible.
Finally $\chi^\xi_i=1$ on $D^i(2)$.
\end{proof}

\subsection{Regular Coordinates and Submersions}\label{es:sec:regularcoordinates}

The gluing needs a criterion under which a realizing map is a topological submersion
over a tube; it is a Lipschitz condition on the coordinate of an admissible section.

\begin{definition}[Regular Sections]\label{es:def:regular}
An admissible section $s_i$ of $F$ is \emph{$\mathsf L$-regular at a point $p\in D^i$} if its
coordinate $w_i$ has the property that $w_i-P_Vh_i$ is $\mathsf L$-Lipschitz with respect to
$h_i$ on a neighbourhood of $p$.
\end{definition}

The next lemma turns this quantitative condition into local product charts; the
coordinate is replaced inside a tube by \cref{es:lem:modall}.

\begin{lemma}[Regular Sections Give Submersions]\label{es:lem:regularall}
Let $F$ be equivariant realizing maps, let $i$ be large, let $U\subset D^i$ be an open set with
$\gamma U\cap U=\emptyset$ for all $\gamma\ne1$ whose image under $h_i$ is convex, and let
$s\colon U\to W^F_i$ be a continuous map with $\pr\circ s=\pi_i$ such that
$\varphi^F\circ s-P_Vh_i$ is locally $\mathsf L$-Lipschitz with respect to $h_i$ on $U$, for
some $\mathsf L\le\tfrac12$. Then $c(p,k):=A_k^{-1}\cdot s(p)$ is a homeomorphism of
$U\times K$ onto $W^F_i\cap\pr^{-1}(\pi_i(U))$ and a product chart of $\varphi^F$ in the
sense of \cref{es:sec:productcharts}, for the right action $w\cdot k:=A_k^{-1}\cdot w$ and
$\beta:=\varphi^F\circ s$; and $F$ is a topological submersion on
$\pr^{-1}(\pi_i(U))$.
\end{lemma}
\begin{proof}
Put $w:=\varphi^F\circ s$. Convexity of $h_i(U)$ makes $w-P_Vh_i$
globally $\mathsf L$-Lipschitz with respect to $h_i$. Consequently
\[
\bigl|(P_Th_i+w)(p)-(P_Th_i+w)(q)\bigr|
\ge(1-\mathsf L)|h_i(p)-h_i(q)|.
\]
The map $P_Th_i+w$ is therefore injective and, by invariance of domain, a
homeomorphism onto an open subset of $T\oplus V$. In these coordinates $w$ is
projection to $V$, so it is a topological submersion.

Since $\pi_i|_U$ is a homeomorphism and $s$ is a section of the principal $K$-bundle
$W^F_i$ over $\pi_i(U)$, the map $c$ is a homeomorphism. It satisfies \mbox{$c(p,kk')=c(p,k)\cdot k'$} and \mbox{$\varphi^F(c(p,k))=\mathsf R(k)^{-1}w(p)$}, so it is the required product chart.
The associated $O(n)$-space is locally a product over $O(n)/I$; in these product charts
$F$ is the product of the identity on $O(n)/I$ with $\varphi^F$.
Hence $F$ is a topological submersion.
\end{proof}

\subsection{Uniform Tubes for Gluing}\label{es:sec:gluing}

The gluing uses one finite system of tubes with a dimension-controlled overlap
bound. We establish the common data and comparisons for arbitrary isotropy groups.

\subsubsection{Setting and Tubes}\label{es:sec:finitesetting}

We fix the setting for the gluing: all centres at once, the threshold, the common scale and
the tubes. With the constants fixed in \eqref{es:eq:kappasharp}, the overlap bound
$(n+1)M_n$ of \cref{es:lem:system}\textup{(ii)} allows a regularity loss of
$60\eps_*$ per tube while keeping the total below $1/4$.
We use all centres of \cref{m:fr:sec:covercollapse} at once, with their charts, along the subsequence of \cref{m:fr:thm:geomain}, along which all
charts converge and their transitions have limit transitions; the objects of \cref{es:sec:standing} for the
centre $a$ carry the superscript $a$. We assume \eqref{th:eq:threshold} with the constants
that belong to this value of $\eps_*$, where the number $\eps_S$ in the definition of $\eps_T$ is replaced by the number $\eps'_S\le\eps_S$ of \cref{th:sec:uniform}; the resulting threshold is denoted by $\eps_{\mathrm{glue}}(n,N)$, and it
depends on $n$ and $N$ only, with $\eps_{\mathrm{glue}}\le\eps_T$. So in place of \eqref{th:eq:threshold}
we assume
\begin{equation}\label{es:eq:thresholdG}
3\delta\le\eps_{\mathrm{glue}}(n,N);
\end{equation}
this is the hypothesis of \cref{es:thm:gluingall} and follows from
$\delta\le\delta_1(n,N)$ in \cref{thm:main}. We fix a number $s_0>0$ such that $16s_0$ is
$c'_S$-admissible at every point of the regions $\Omega^{\flat,a}$ of all centres, with the
number $c'_S$ of \cref{th:sec:uniform}, and such that $36s_0\le\tfrac18R_0$ and
$16s_0\le\tfrac12R_1$; this is possible because the regions $\Omega^a$ contain the
balls of radius $\tfrac34R_0$ about the points of $\Omega^{\flat,a}$ and the curvature
bound and $R_1$ are fixed. Since $c'_S\le c_S$, the scale $8s_0$ is then $c_S$-admissible at these points, as
\cref{th:sec:standing} requires. By \cref{th:lem:recentre} the harmonic chart of \cref{th:lem:chart} then exists at the scale $s_0$ at every point of these regions.

A \emph{tube} at a point $x\in X$ is a tube of \cref{es:def:tube}, with its data, at a point $y$
over $x$, for a centre $a$ with $x\in B_{R_0/4}(a)$, built from the data of \cref{es:sec:standing}
at the scale $s_0$ and a number $\ell\in(0,\ell_*]$. We
mark the objects of a tube $\alpha$ with the index $\alpha$: the charts $Z^\alpha_i$ and $Z^\alpha$ with their groups, the
group $K_\alpha$, the spaces $T_\alpha$ and $V_\alpha$, the transplanted charts $h^\alpha_i$ and the harmonic charts $h^\alpha$, the coordinates $\tilde\xi^\alpha_i$ and $\xi^\alpha$, the labels $k^\alpha_\gamma$,
the slice $S^\alpha$, the tube $\Xi_\alpha$, the domains $D^{i,\alpha}$ and so on; the unit
$r_\ell=\ell/800$ depends on the tube. When a tube at $x$ is fixed, we use $x$ itself as its index
and speak of the tube $x$. For finitely many tubes we pass to a common
harmonic-chart subsequence as in \cref{th:lem:stagecharts}; the coordinate and
cocycle limits then follow from the Haar formulas of \cref{es:lem:limit,es:lem:limitcocycle}, without further extraction.
The smaller sets $D^{i,\alpha}(t)$, $\cW^\alpha_i(t)$, $O_\alpha(t)$ and
$D^{i,\alpha}_c$ have the definitions of \cref{es:sec:realizing}, with the tube index added.

As fixed in \cref{th:sec:standing}, the realizing maps below are equivariant
(\cref{fl:def:realizing}).

\subsubsection{Transfer between Tubes}\label{es:sec:transfer}

We compare the harmonic charts of two tubes along the chart transitions;
\cref{es:lem:transferred} uses the comparison to carry regularity (\cref{es:def:regular}) from
one tube to another. Let $\alpha$ and $\beta$ be
tubes, which may coincide. We use the transitions of \cref{m:fr:prop:coverscharts}\textup{(iii)},
which we call \emph{chart transitions}: for points $p\in\Omega'^{\alpha}_i$ and $p'\in\Omega'^{\beta}_i$
over the same point of $X_i$ there is exactly one chart transition $\iota_i$ with
$\iota_i(p)=p'$, a map of $B_{R_1}(p)$ onto $B_{R_1}(p')$ that preserves distances and
satisfies $\pi^\beta_i\circ\iota_i=\pi^\alpha_i$. It acts on the frames by
$\iota_{i*}(q,e):=(\iota_i(q),e\circ\iota_i^{-1})$, and $\tau^\beta_i\circ\iota_{i*}=\tau^\alpha_i$ at the
frames over a neighbourhood of $p$ (\cref{m:fr:prop:coverbundle}\textup{(i)}). If points $p_i$
converge and the points $\iota_i(p_i)$ stay bounded, a subsequence of the maps $\iota_i$
converges to a limit transition $\iota$ in the sense of \cref{m:fr:def:limittransition}, which
satisfies $\tau^\beta\circ\iota_*=\tau^\alpha$ near the limit of the points $p_i$; and if frames
over points $q_i\to q$ converge to a frame $\mathsf p$, their images under $\iota_{i*}$
converge to $\iota_*(\mathsf p)$ (\cref{m:fr:lem:limittransitions} and \ref{fl:S:charts}). For
tubes with the same centre the transitions are restrictions of elements of the groups.
The harmonic charts of two tubes are compared along a transition as follows.

\begin{lemma}[Two Harmonic Charts along a Transition]\label{es:lem:twochartstransition}
Let $\alpha$ and $\beta$ be tubes.
\begin{enumerate}[label=\textup{(\roman*)}]
\item Let $\iota$ be a limit transition from $\alpha$ to $\beta$ that maps a point of
$\cC^\alpha_3$ into $\cC^\beta_3$. Then there is $R\in O(n)$ such that
$h^\beta\circ\iota-R\,h^\alpha$ is $\tfrac12\eps_*$-Lipschitz with respect to $h^\alpha$ on
$\cC^\alpha_3$.
\item Let $\iota_i$ be chart transitions from $\alpha$ to $\beta$, for $i$ in an infinite
set of indices, that map points of $\cC^{i,\alpha}_3$ into $\cC^{i,\beta}_3$. Then for all
large $i$ there are matrices $R_i\in O(n)$ such that $h^\beta_i\circ\iota_i-R_i\,h^\alpha_i$
is $3\eps_*$-Lipschitz with respect to $h^\alpha_i$ on $\cC^{i,\alpha}_3$. If the maps
$\iota_i$ converge to a limit transition $\iota$ and $R_i\to R$, then
$h^\beta\circ\iota-R\,h^\alpha$ is $3\eps_*$-Lipschitz with respect to $h^\alpha$ on
$\cC^\alpha_3$.
\end{enumerate}
\end{lemma}
\begin{proof}
The two tubes have the same scale $s_0$, and $16s_0\le\tfrac12R_1$. Harmonic functions
remain harmonic under composition with transitions, because a transition preserves
distances and maps the measure to a constant multiple of the measure. The coordinate
balls $\cC_3$ of the tubes lie in the balls of radius $d_+(\ell)$ about the points
$z_\alpha$, $z_\beta$, $z^\alpha_i$ and $z^\beta_i$, by \eqref{th:eq:metricsize}.

(i) The domain of $\iota$ contains $B_{4s_0}(z_\alpha)$, and $\iota$ maps this ball into
$B_{8s_0}(z_\beta)$. The restriction of $h^\alpha$ is a harmonic chart with parameter
$2\theta_\flat$ at $(z_\alpha,\tfrac12s_0)$, and $h^\beta\circ\iota-h^\beta(\iota(z_\alpha))$ is a
harmonic map on $B_{4s_0}(z_\alpha)$ that changes distances by at most $3\theta_\flat s_0$. By
\cref{th:lem:twocharts} at $(z_\alpha,\tfrac12s_0)$ the two differ by a rigid motion up to a
map that is $\tfrac12\eps_*$-Lipschitz with respect to $h^\alpha$ on $B_{s_1/2}(z_\alpha)$,
which contains $\cC^\alpha_3$.

(ii) Let $v^\alpha_i$ and $v^\beta_i$ be the harmonic charts of \cref{th:lem:stagecharts} for
the two tubes. The map $\iota_i$ sends $B_{8s_0}(z^\alpha_i)$ into the ball of radius
$16s_0$ about $z^\beta_i$, on which $v^\beta_i$ is defined and changes distances by at most
$2\theta_\flat s_0$. So $u':=v^\beta_i\circ\iota_i-v^\beta_i(\iota_i(z^\alpha_i))$ is a harmonic
map on $B_{8s_0}(z^\alpha_i)$ that changes distances by at most $12\theta_\flat s_0$, and the
restriction of $v^\alpha_i$ is a harmonic chart with parameter at most $2\theta_\flat$ at
$(z^\alpha_i,s_0)$. By \cref{th:lem:twocharts} there is a rigid motion $E_i$ such that
$v^\beta_i\circ\iota_i-E_i\circ v^\alpha_i$ is $\tfrac12\eps_*$-Lipschitz with respect to
$v^\alpha_i$ on $B_{s_1}(z^\alpha_i)$. Let $R_i$ be its rotation part. We have
$h^\alpha_i=\hat h^\alpha\circ v^\alpha_i$ and $h^\beta_i=\hat h^\beta\circ v^\beta_i$ with maps
$\hat h^\alpha$ and $\hat h^\beta$ that differ from the identity by maps with the Lipschitz
constant $\tfrac12\eps_*$, and the computation in the proof of
\cref{th:lem:transplant}\textup{(iv)} shows that $h^\beta_i\circ\iota_i-R_i\,h^\alpha_i$ is
$3\eps_*$-Lipschitz with respect to $h^\alpha_i$ on $B_{s_1}(z^\alpha_i)$, at the points
that $\iota_i$ maps into $B_{s_1}(z^\beta_i)$; these include the points of $\cC^{i,\alpha}_3$.
If points $p_i\in \cC^{i,\alpha}_3$ converge to $p\in \cC^\alpha_3$, then
$\iota_i(p_i)\to\iota(p)$ and $h^\beta_i(\iota_i(p_i))\to h^\beta(\iota(p))$, by
\cref{th:lem:transplant}\textup{(iii)}; so the last assertion is the limit of the first.
\end{proof}

For converging chart transitions we refer to the conclusion of (ii) as
\begin{equation}\label{es:eq:stagetwocharts}
h^\beta_i\circ\iota_i-R_i h^\alpha_i
\quad\text{is $3\eps_*$-Lipschitz with respect to }h^\alpha_i
\text{ on }\cC^{i,\alpha}_3.
\end{equation}

\subsubsection{A System of Tubes}\label{es:sec:system}

We choose finitely many tubes $x$ whose sets $O_x(1)$ cover $X$ and of which only a number bounded
in terms of $n$ meet at any point. The tubes at points with different isotropy groups
have different sizes, and the bound is obtained stratum by stratum. For $x\in X$ let
$k_0(x)$ be the dimension of the space $V_0$ of a tube at $x$.

\begin{lemma}[Strata]\label{es:lem:strata}
Let $\alpha$ be a tube, let $v\in B_{V_\alpha}(\tfrac12\ell)$, let $K_v$ be the group of
the elements $k\in K_\alpha$ with $\mathsf R(k)v=v$, and let $N_v\subset V_\alpha$ be the
orthogonal complement of the tangent space at $v$ of the orbit $\mathsf R(K_\alpha)v$; when
$K_\alpha$ is finite, $N_v=V_\alpha$. Then $k_0(\bar\xi_\alpha^{-1}(\mathsf R(K_\alpha)v))$ is the
dimension of the space of the vectors of $N_v$ that are fixed by $\mathsf R(K_v)$. So
$k_0\ge k_0(x_\alpha)$ on $O_\alpha$, with equality exactly at the points
$\bar\xi_\alpha^{-1}(v)$ with $v\in V_{0,\alpha}$; the number $k_0(x)$ does not depend on the
tube at $x$; and the sets $X^{(\le k)}$ of the points with $k_0(x)\le k$ are closed.
\end{lemma}
\begin{proof}
We drop the index. The isotropy group of a point $\Xi([A,v'])$ of $Y$ is the conjugate
by $A$ of the group $I_{v'}$ of the matrices $A_k$ with $k\in K_{v'}$. Let $r_\perp>0$ be such that \eqref{es:eq:lineartube} is a homeomorphism onto an open subset of
$B_V(\tfrac12\ell)$. By \cref{es:thm:nearbyprecise}, the image of \mbox{$v+B_{N_v}(r_\perp)$} under $\bar\xi^{-1}$ is an open
neighbourhood of $x':=\bar\xi^{-1}(\mathsf R(K)v)$ that is homeomorphic to
$B_{N_v}(r_\perp)/\mathsf R(K_v)$, and the points of $Y$ over the image of $v+w$ have isotropy
groups conjugate to $I_{v+w}$, where $K_{v+w}\subset K_v$ is the group of the elements
of $K_v$ that fix $w$. A closed subgroup of a compact Lie group that is conjugate to
the group, in a larger group, is the whole group, because it has the same dimension
and the same number of components. So the points near $x'$ over which the isotropy groups are conjugate to those over $x'$ are the
images of the points $v+w$ with $w$ fixed by $\mathsf R(K_v)$, and they form a set that is
homeomorphic to a ball in the space of the vectors of $N_v$ fixed by $\mathsf R(K_v)$; this is the
description of the orbit types near a point by the slice representation
\cite[Chapter~2]{DK00}. For $v=0$ this shows that $k_0(x)$ is the
dimension at $x$ of the set of the points with the same type of isotropy groups, which
does not depend on the tube; for general $v$ it is the first assertion. The tangent
space of the orbit at $v$ is orthogonal to $V_0$ and to $v$, because the matrices
$\mathsf R(k)$ are orthogonal; so the space of the vectors of $N_v$ fixed by $\mathsf R(K_v)$
contains $V_0$ and $v$, and it equals $V_0$ if and only if $v\in V_0$. If points $x_m$
with $k_0(x_m)\le k$ converge to $x$, they lie in the set $O$ of a tube at $x$ for large
$m$, so $k_0(x)\le k_0(x_m)\le k$.
\end{proof}

On a compact set on which $k_0$ is constant the tubes can be given one scale, by the recentring of \cref{th:lem:recentre}: one harmonic chart serves all nearby
points with the same isotropy group. The bounds \eqref{th:eq:stageholder} are the same for all tubes, because
the scale $s_0$ is common and $c_3$, $C_6$ and $\eta$ are universal.

\begin{lemma}[Uniform Data]\label{es:lem:uniform}
Let $\mathcal S\subset X$ be a compact set on which $k_0$ is constant. Then there are a number
$\ell_{\mathcal S}>0$ and, for every $x\in \mathcal S$, a tube at $x$ with $\ell=\ell_{\mathcal S}$.
\end{lemma}
\begin{proof}
Let $x\in \mathcal S$. We choose a centre $a$ with $x\in B_{R_0/8}(a)$ and a point
$z_x\in\Omega^{\flat,a}$ over $x$, with the data and the number $\rho_{z_x}$ of
\cref{th:lem:recentre} for the scale $s_0$; that lemma applies by the choice of $s_0$. Let
$K$ be the isotropy group of $z_x$, and let $S$ be a slice at $z_x$ as in
\cref{th:lem:slices}. The isotropy group of a point of $S$ is a subgroup of $K$, and it
is $K$ exactly at the points that are fixed by $K$. If $x'\in \mathcal S$ is close to $x$, then
$k_0(x')=k_0(x)$, so by \cref{es:lem:strata} the isotropy groups of the points of $Y$ over $x'$
are conjugate to those of the points over $x$; as in the proof of that lemma, the point
of $S$ over $x'$ then has the isotropy group $K$. Hence there is $\rho_x>0$ such that
every point of $\mathcal S\cap B_{\rho_x}(x)$ is the image of a point $z'$ with the isotropy
group $K$ and $\hat d(z',z_x)<\min\{\rho_{z_x},\tfrac1{16}R_0\}$; then
$z'\in\Omega^{\flat,a}$.

Finitely many of the balls $B_{\rho_x}(x)$ cover $\mathcal S$; let $x_1,\dots,x_m$ be their
centres. To a point of $\mathcal S\cap B_{\rho_{x_j}}(x_j)$ and the corresponding point $z'$ of
the chart of $x_j$ we give the data of \cref{th:lem:recentre} at $z'$: the harmonic chart $h-h(z')$, the homomorphism $\mathsf R$, the complement $\fm$ with its norm, and the numbers
$r_0$, $c_0$ and $C_3$ of the point $z_{x_j}$. By that lemma all assertions of
\cref{th:lem:chart} hold at $z'$ for the scale $s_0$, the number $c_{z'}$ can be taken
equal to $\tfrac12c_{z_{x_j}}$, and the numbers that enter \cref{th:lem:almostequiv,th:lem:slices} can be chosen independently of $z'$. The constants of
\cref{th:lem:stagecharts} to \cref{th:lem:stageequiv} and the functions of the scale in \cref{th:lem:errors} are
determined by these data, by the scale $s_0$ and by $n$ and $N$; so they are the same for the points that belong to one $x_j$. The numbers $C_{\mathsf R}$, $r_K$, $r_c$ and $C_c$ of
\cref{es:sec:lambda} depend on $K$, $\mathsf R$, $\fm$ and the norms only. The number
$\sigma_{\mathsf k}$ can be chosen independently of $z'$: if every $g$ with
$\hat d(gz_{x_j},z_{x_j})<\sigma$ lies in a given neighbourhood of $K$ and
$\rho_{z_{x_j}}\le\tfrac14\sigma$, then every $g$ with $\hat d(gz',z')<\tfrac12\sigma$
lies in it, because $\hat d(gz_{x_j},z_{x_j})\le\hat d(gz',z')+2\rho_{z_{x_j}}$. So the
conditions on $\ell$ of \cref{th:sec:scale}, including \eqref{th:eq:ellsmall} and
\eqref{es:eq:ellone}, hold for all these points with one number, which we call $\ell_{\mathcal S}$.
\end{proof}

For a tube $x$ and $t>0$ let $\hat O_x(t)$ be the image under $\pi^x$ of the set of
the points $q\in \cC^x_{1/2}$ with $P^x_Th^x(q)=0$ and $|h^x(q)|<t\,r_\ell$. These sets do not involve the sequence. We use three facts. Let $x$ and $x'$ be tubes.
\begin{enumerate}[label=\textup{(\alph*)}]
\item For $0<t\le9$ we have $\hat O_x(t)\subset O_x(\tfrac{26}{25}t)$ and
$O_x(t)\subset\hat O_x(\tfrac{26}{25}t)$, for every limit $\xi^x$ as in \cref{es:lem:limit}.
Indeed, $\xi^x$ and $h^x$ vanish at $z_x$ and $\xi^x-P^x_Vh^x$ is $7\eps_*$-Lipschitz
with respect to $h^x$; so $(1-7\eps_*)|h^x|\le|\xi^x|\le(1+7\eps_*)|h^x|$ at the points
with $P^x_Th^x=0$, and the points of $S^x$ are such points.
\item If $p\in \cC^x_{1/2}$ and $|h^x(p)|<t\,r_\ell$ with $t\le50$, then
$\pi^x(p)\in\hat O_x(2t)$. Indeed, by \cref{th:lem:meets} there is a point $q$ in the
orbit of $p$ with $P^x_Th^x(q)=0$ at which $|h^x|$ is at most the sum of
$|P^x_Vh^x(p)|$ and $|P^x_Th^x(p)|$, hence less than $2t\,r_\ell$, as in the proof of
\cref{es:lem:slice}\textup{(iii)}.
\item If $\iota$ is a limit transition from $x$ to $x'$ that maps a point of $\cC^x_3$
into $\cC^{x'}_3$, then $\pi^{x'}\circ\iota=\pi^x$ and
\begin{equation*}
|h^{x'}(\iota(p))-h^{x'}(\iota(q))|\le(1+\eps_*)\,|h^x(p)-h^x(q)|\qquad\text{for }p,q\in \cC^x_3
\end{equation*}
(\cref{es:lem:twochartstransition}\textup{(i)}); such a transition exists and can be
chosen with $\iota(p)=p'$ for given points $p\in \cC^x_3$ and $p'\in \cC^{x'}_3$ over the same
point of $X$ (\cref{m:fr:lem:limittransitions}\textup{(ii)}).
\end{enumerate}

\begin{lemma}[A System of Tubes]\label{es:lem:system}
There are finitely many tubes $x_1,\dots,x_J$ with the following properties, along
every subsequence for which \cref{th:lem:stagecharts} holds for all of them.
\begin{enumerate}[label=\textup{(\roman*)}]
\item The sets $O_j(1)$ cover $X$.
\item Every point of $X$ lies in at most $(n+1)M_n$ of the sets $O_j(8)$.
\end{enumerate}
\end{lemma}
\begin{proof}
For a tube $x$ let $U_x$ be the interior of $\hat O_x(\tfrac{19}{20})$. By (a), for one
limit $\xi^x$, the open set $O_x(\tfrac{21}{25})$ contains $\hat O_x(\tfrac45)$ and lies
in $\hat O_x(\tfrac{19}{20})$, and the latter lies in $O_x(1)$; so
\begin{equation}\label{es:eq:Ux}
\hat O_x(\tfrac45)\subset U_x\subset O_x(1)\qquad\text{for every limit }\xi^x .
\end{equation}
For tubes $x$ and $x'$ with the same number $\ell$ we write $x'\prec x$ if
$x'\in\hat O_x(\tfrac13)$.

We define finite sets $\mathfrak F_0,\mathfrak F_1,\dots,\mathfrak F_n$ of tubes such that
$X^{(\le k)}$ lies in the union of the sets $U_x$ for
$x\in\mathfrak F_0\cup\dots\cup\mathfrak F_k$. Let $k\ge0$, and suppose that the sets with
smaller indices are defined. Let $\mathcal S_k$ be the set of the points of $X^{(\le k)}$ outside
the sets $U_x$ for the tubes chosen so far. It is compact, by \cref{es:lem:strata}, and
$k_0(x)=k$ on it. If it is empty, let $\mathfrak F_k:=\emptyset$. Otherwise let
$\ell_k:=\ell_{\mathcal S_k}$ be the number of \cref{es:lem:uniform}, with its tubes at the points of
$\mathcal S_k$, and let $r_\ell:=\ell_k/800$. If $x,x'\in \mathcal S_k$ and $d(x,x')<r_\ell/(8C_6)$, there is a
point $p$ over $x'$ in the chart of $x$ with $\hat d(p,z_x)=d(x,x')$ by \ref{fl:S:limit}, so $|h^x(p)|<\tfrac18r_\ell$ by \eqref{th:eq:stageholder} in the
limit, and $x'\prec x$ by (b). Hence a set of points of $\mathcal S_k$ no two of which are related
is finite; let $\mathfrak F_k$ be a maximal one. Let $x\in \mathcal S_k$. Then $x\in U_x$ if
$x\in\mathfrak F_k$, by \eqref{es:eq:Ux}; otherwise there is $x_j\in\mathfrak F_k$ with
$x\prec x_j$ or $x_j\prec x$. In the first case $x\in U_j$, by \eqref{es:eq:Ux}. In the
second case let $q\in \cC^x_{1/2}$ be a point over $x_j$ with $|h^x(q)|<\tfrac13r_\ell$, and
let $\iota$ be a limit transition from $x$ to $x_j$ with $\iota(q)=z_j$. By (c),
$|h^j(\iota(z_x))|<\tfrac{17}{50}r_\ell$, and $\iota(z_x)$ lies over $x$; so
$x\in\hat O_j(\tfrac{17}{25})\subset U_j$, by (b) and \eqref{es:eq:Ux}. Thus $\mathcal S_k$, and with
it $X^{(\le k)}$, lies in the union of the sets $U_x$. Since $k_0(x)\le n$ for all $x$, the
tubes of all sets $\mathfrak F_k$ together satisfy (i), by \eqref{es:eq:Ux}.

(ii) Fix $k$ and a point $x_*\in X$, and let $j_0$ be a tube of $\mathfrak F_k$ with
$x_*\in O_{j_0}(8)$. Let $j\in\mathfrak F_k$ be another tube with $x_*\in O_j(8)$. By (a)
there are points $p_*$ and $p'_*$ over $x_*$ in the coordinate balls of $j$ and of $j_0$
at which $|h^j|$ and $|h^{j_0}|$ are less than $\tfrac{208}{25}r_\ell$. Let $\iota$ be a limit
transition from $j$ to $j_0$ with $\iota(p_*)=p'_*$. By (c),
$|h^{j_0}(\iota(z_j))|<\tfrac{208}{25}r_\ell+(1+\eps_*)\tfrac{208}{25}r_\ell<\tfrac{167}{10}r_\ell$, and by (b)
there is a point $q_j\in \cC^{j_0}_{1/2}$ over $x_j$ with $P^{j_0}_Th^{j_0}(q_j)=0$ and
$|h^{j_0}(q_j)|<\tfrac{167}5r_\ell$. It lies in $S^{j_0}$, by the estimate in (a), and since $k_0(x_j)=k_0(x_{j_0})$, the
vector $\xi^{j_0}(q_j)$ lies in $V_{0,j_0}$ by \cref{es:lem:strata}. Hence $q_j$ is fixed by
$K_{j_0}$, and $u_j:=h^{j_0}(q_j)$ is a vector of $V_{0,j_0}$, because $h^{j_0}$ commutes
with the actions of $K_{j_0}$. For $j=j_0$ we put $q_{j_0}:=z_{j_0}$ and $u_{j_0}:=0$. Let
$j\ne j'$ be two tubes of this kind, and suppose that $|u_j-u_{j'}|<\tfrac18r_\ell$. Let
$\iota'$ be a limit transition from $j_0$ to $j'$ with $\iota'(q_{j'})=z_{j'}$. By (c),
$|h^{j'}(\iota'(q_j))|<\tfrac{13}{100}r_\ell$, and this point lies over $x_j$; by (b),
$x_j\in\hat O_{j'}(\tfrac13)$, that is, $x_j\prec x_{j'}$, which contradicts the choice of
$\mathfrak F_k$. So the vectors $u_j$ lie in the ball of radius $\tfrac{167}5r_\ell$ of a space of
dimension at most $n$ and have pairwise distances at least $\tfrac18r_\ell$; their number
is at most $(16\cdot\tfrac{167}5+1)^n\le M_n$. Summing over $k=0,\dots,n$ gives (ii).
\end{proof}

\subsection{Gluing along the Limit}\label{es:sec:gluingall}

We modify an equivariant map one tube at a time, carrying regularity (\cref{es:def:regular})
between tubes.
Admissible sections differ by small torus elements whose exponents can oscillate
(\cref{es:lem:admissible}). We therefore average the exponents of transferred
sections. Unlike the twisted linear average of \cref{es:sec:twisted}, this non-linear
average requires a remainder estimate that remains valid near fixed directions;
\cref{es:lem:averages} provides it without a Lipschitz bound on the exponents.

\subsubsection{Statement}\label{es:sec:allstatement}

Here we state the gluing theorem, the main result of the construction; its proof
occupies the rest of the section. We use the setting and notation of \cref{es:sec:finitesetting}: the threshold \eqref{es:eq:thresholdG}, the constants \eqref{es:eq:kappasharp}, the
scale $s_0$, tubes, the sets
$D^{i,\alpha}(t)$, $\cW^\alpha_i(t)$ and $O_\alpha(t)$, and realizing maps.

\begin{theorem}[Gluing]\label{es:thm:gluingall}
In the setting of \cref{es:sec:standing}, for the constant $\eps_*$ of
\eqref{es:eq:kappasharp} and all centres, assume \eqref{es:eq:thresholdG}. Then every infinite
set of indices contains an infinite subset along which, for all large $i$, there are maps
$F^\natural_i\colon Y_i\to Y$ with the following properties.
\begin{enumerate}[label=\textup{(\roman*)}]
\item The maps $F^\natural_i$ are continuous and surjective, they commute with $O(n)$,
and the supremum of $\sd(F^\natural_i,F_i)$ over $Y_i$ tends to $0$.
\item Every $F^\natural_i$ is a proper topological submersion and a fibre bundle whose
fibres are closed manifolds; over the image of the chart of a centre their dimension is
the dimension of the limit group of this chart.
\item For every $x\in X$ there are a point $y$ over $x$ with the isotropy group $I$, an
orthogonal representation of $I$ on a Euclidean space $V$, an open ball $B$ about $0$ in $V$, and a linear tube $\Xi$ of $O(n)\times_IB$ onto an open
neighbourhood of the orbit of $y$ in $Y$ (\cref{ls:def:slice}), with $\Xi([\one,0])=y$, such that
the following holds for all large $i$, with an index bound that does not depend on $x$:
there are a closed manifold $\mathsf F_i$ with a
free action of $I$ and a homeomorphism of $O(n)\times_I(\mathsf F_i\times B)$ onto the
preimage of the image of $\Xi$ under $F^\natural_i$ that commutes with $O(n)$ and under
which $F^\natural_i$ corresponds to the projection to $O(n)\times_IB$, followed by $\Xi$.
\item Let $x_\alpha$ be a tube of the system of \cref{es:lem:system}, with its group $K$, its
representation $\mathsf R$ on $V$, its unit $r_\ell$ and the maps $\Xi$ and $\bar\xi$ of
\cref{es:sec:limitcoordinate}. For $v\in B_V(r_\ell)$ let $K_v$ be the group of the elements of
$K$ that fix $v$, and $N_v$ the orthogonal complement in $V$ of the tangent space at $v$ of
the orbit $\mathsf R(K)v$; so $\dim N_v=\dim V-\dim K+\dim K_v$. For all large $i$ there is a
closed manifold $\mathsf F_i$ with a free action of $K$, homeomorphic to the fibre of
$F^\natural_i$ over $\Xi([\one,0])$, such that the following holds for every $v\in B_V(r_\ell)$ and
every $r_\perp>0$ for which \eqref{es:eq:lineartube} is a homeomorphism onto an open subset of
$B_V(5r_\ell)$: there is a homeomorphism of $(f^\natural_i)^{-1}(O^v)$, where
$f^\natural_i\colon X_i\to X$ is the map induced by $F^\natural_i$ and
$O^v:=\bar\xi^{-1}\bigl(\mathsf R(K)(v+B_{N_v}(r_\perp))/\mathsf R(K)\bigr)$, onto
$\mathsf F_i\times_{K_v}B_{N_v}(r_\perp)$ under which $f^\natural_i$ corresponds to the map
$[u,w]\mapsto\bar\xi^{-1}(\mathsf R(K)(v+w))$. In particular the fibre of $f^\natural_i$ over
$\bar\xi^{-1}(\mathsf R(K)v)$ is homeomorphic to $\mathsf F_i/K_v$, and it is a closed manifold of
dimension $\dim\mathsf F_i-\dim K_v$.
\end{enumerate}
\end{theorem}

The proof is given in \cref{es:sec:allproof}. The maps
$f^\natural_i\colon X_i\to X$ that the maps $F^\natural_i$ induce satisfy
$\sup d(f^\natural_i,f_i)\to0$, and over the neighbourhood $\pr(\Xi(O(n)\times_IB))$ of
$x$, which is homeomorphic to $B/I$, the map $f^\natural_i$ is the map
$\mathsf F_i\times_IB\to B/I$. This is the local model of \cite{Fuk89,CFG92}, apart from the
identification of $\mathsf F_i$ with a nilmanifold, which \cref{es:sec:fibres} carries out
(\cref{es:cor:fibresnil}).

\subsubsection{Transferred Sections}\label{es:sec:transferred}

We carry a section and the regularity of its coordinate from one tube to another,
using the limit section $\hat\sigma$ and one fixed transition for each pair of tubes.
The resulting continuous section has convergent frames, but the deck transformations
need not act on it through the cocycle of the second tube.

Let $\alpha$ and $\beta$ be tubes, and let $\iota_i$ be chart transitions from $\alpha$
to $\beta$, for all large $i$, that converge to a limit transition $\iota$ and map a
point of $\cC^{i,\alpha}_3$ into $\cC^{i,\beta}_3$ (\cref{es:sec:transfer}). By
\cref{es:lem:twochartstransition}\textup{(ii)} there are matrices $R_i\in O(n)$ with
\eqref{es:eq:stagetwocharts}; after passing to a subsequence they converge to a matrix
$R$, and $h^\beta\circ\iota-R\,h^\alpha$ is $3\eps_*$-Lipschitz with respect to
$h^\alpha$ on $\cC^\alpha_3$. Let
\begin{equation*}
D^{i,\alpha}_\beta:=\bigl\{q\in D^{i,\alpha}:\iota_i(q)\in D^{i,\beta}\bigr\},
\end{equation*}
an open subset of $D^{i,\alpha}$. Let $F=F^{(i)}$ be realizing maps, and let $s^\alpha_i$
be an admissible section of $F$ in the tube $\alpha$, with the coordinate $w_i$ and the
frames $\tilde e_i$. We recall the coordinates $\check h^\alpha=P^\alpha_Th^\alpha+\xi^\alpha$ on $\cC^\alpha_1$, whose image contains the closed ball of radius $\tfrac12\ell_\alpha$. For
large $i$ we have $|w_i|<9r_{\ell_\alpha}$ on $D^{i,\alpha}$
(\cref{es:lem:admissible}\textup{(iii)}), and we put
\begin{equation}\label{es:eq:lift}
\psi_i:=(\check h^\alpha)^{-1}\circ\bigl(P^\alpha_Th^\alpha_i+w_i\bigr)\colon D^{i,\alpha}\to \cC^\alpha_1,
\qquad\text{so that}\qquad
P^\alpha_Th^\alpha\circ\psi_i=P^\alpha_Th^\alpha_i,\quad\xi^\alpha\circ\psi_i=w_i .
\end{equation}

\Needspace{5\baselineskip}
\begin{lemma}[Transferred Sections]\label{es:lem:transferred}
For all large $i$ the following holds.
\begin{enumerate}[label=\textup{(\roman*)}]
\item For $q\in D^{i,\alpha}$ the point $\psi_i(q)$ lies in $\hat U_{S^\alpha}$, and for
$q\in D^{i,\alpha}_\beta$ the point $\iota(\psi_i(q))$ lies in $\hat U_{S^\beta}$. So there is a
continuous map $\mathsf A_i\colon D^{i,\alpha}_\beta\to O(n)$ with
\begin{equation}\label{es:eq:Adef}
\iota_*\bigl(\hat\sigma^\alpha(\psi_i(q))\bigr)\cdot\mathsf A_i(q)
=\hat\sigma^\beta\bigl(\iota(\psi_i(q))\bigr)\qquad(q\in D^{i,\alpha}_\beta).
\end{equation}
\item The map $s^\iota_i(q):=\mathsf A_i(q)^{-1}\cdot s^\alpha_i(q)$ on $D^{i,\alpha}_\beta$ has values in
$W^{F,\beta}_i$, the point $s^\iota_i(q)$ lies over $\pi^\beta_i(\iota_i(q))$, and
\begin{equation}\label{es:eq:transferredcoordinate}
\varphi^{F,\beta}\bigl(s^\iota_i(q)\bigr)=\xi^\beta\bigl(\iota(\psi_i(q))\bigr).
\end{equation}
\item Let $s^\beta_i$ be an admissible section of $F$ in the tube $\beta$. Then there
are continuous maps $\mu_i\colon D^{i,\alpha}_\beta\to\fk_\beta(r_{K_\beta})$ with
\begin{equation*}
s^\iota_i(q)=A_{\exp(\mu_i(q))}\cdot s^\beta_i\bigl(\iota_i(q)\bigr)\quad(q\in D^{i,\alpha}_\beta),\qquad\sup_{D^{i,\alpha}_\beta}|\mu_i|\to0 .
\end{equation*}
\item Let $0\le\mathsf L\le\tfrac14$, and let $q\in D^{i,\alpha}_\beta$ be a point at which $s^\alpha_i$ is
$\mathsf L$-regular. Then the map
$w''_i:=\varphi^{F,\beta}\circ s^\iota_i\circ\iota_i^{-1}$ on the open set $\iota_i(D^{i,\alpha}_\beta)$ has
the property that $w''_i-P^\beta_Vh^\beta_i$ is $(\mathsf L+40\eps_*)$-Lipschitz with
respect to $h^\beta_i$ on a neighbourhood of $\iota_i(q)$.
\end{enumerate}
\end{lemma}
\begin{proof}
(i) For $q\in D^{i,\alpha}$ the vector $\check h^\alpha(\psi_i(q))$ has a component in
$T_\alpha$ of norm less than $\tfrac14\ell_\alpha$ and a component in $V_\alpha$ of norm
less than $9r_{\ell_\alpha}$. As in the proof of \cref{es:lem:domain}\textup{(i)} this gives
$\psi_i(q)\in \cC^\alpha_{1/2}$, and $\psi_i(q)\in U_{S^\alpha}$ by
\cref{es:lem:slice}\textup{(iii)} and the proof of \cref{es:lem:limitsection}\textup{(ii)}. If
points $q_i\in D^{i,\alpha}$ converge to $q_\infty$, then $\psi_i(q_i)\to q_\infty$,
because $\check h^\alpha(\psi_i(q_i))=P^\alpha_Th^\alpha_i(q_i)+w_i(q_i)$ converges to
$\check h^\alpha(q_\infty)$, by \cref{es:lem:admissible}\textup{(iii)} and \cref{es:lem:limit}.
Suppose that $\iota(\psi_i(q_i))\notin \hat U_{S^\beta}$ for infinitely many $i$ and points
$q_i\in D^{i,\alpha}_\beta$. A subsequence converges to a point $q_\infty$, and
$\iota(\psi_i(q_i))\to\iota(q_\infty)$, which is the limit of the points
$\iota_i(q_i)\in D^{i,\beta}$ and so lies in the open set $\hat U_{S^\beta}$, by
\ref{es:A:frames} for any admissible section in the tube $\beta$; a contradiction. Both
sides of \eqref{es:eq:Adef} are frames over the point $\iota(\psi_i(q))$, and they depend
continuously on $q$.

(ii) Since $F$ commutes with $O(n)$ and $\tau^\beta\circ\iota_*=\tau^\alpha$, we get from
\eqref{es:eq:lift}, \cref{es:lem:sigmahat}\textup{(i)} for both tubes and \eqref{es:eq:Adef}
\begin{align*}
F\bigl(s^\iota_i(q)\bigr)&=\mathsf A_i(q)^{-1}\cdot\Xi_\alpha\bigl([\one,w_i(q)]\bigr)
=\tau^\alpha\bigl(\hat\sigma^\alpha(\psi_i(q))\cdot\mathsf A_i(q)\bigr)\\
&=\tau^\beta\bigl(\hat\sigma^\beta(\iota(\psi_i(q)))\bigr)
=\Xi_\beta\bigl([\one,\xi^\beta(\iota(\psi_i(q)))]\bigr).
\end{align*}
The point $s^\iota_i(q)$ lies over $\pi^\alpha_i(q)=\pi^\beta_i(\iota_i(q))$.

(iii) We have
$s^\iota_i(q)=\tau^\alpha_i(\tilde e_i(q)\cdot\mathsf A_i(q))=\tau^\beta_i(\iota_{i*}(\tilde e_i(q))\cdot\mathsf A_i(q))$.
Let $q_i\in D^{i,\alpha}_\beta$ converge to $q_\infty$. The frames $\tilde e_i(q_i)$ converge to
$\hat\sigma^\alpha(q_\infty)$ by \ref{es:A:frames}, their images under $\iota_{i*}$
converge to $\iota_*(\hat\sigma^\alpha(q_\infty))$ (\cref{es:sec:transfer}), and the
matrices $\mathsf A_i(q_i)$ converge to the matrix $\mathsf A_\infty$ that maps the
frame $\iota_*(\hat\sigma^\alpha(q_\infty))$ to the frame $\hat\sigma^\beta(\iota(q_\infty))$,
because $\psi_i(q_i)\to q_\infty$. So the frames over $\iota_i(q_i)$ that $\tau^\beta_i$
maps to $s^\iota_i(q_i)$ converge to $\hat\sigma^\beta(\iota(q_\infty))$, as the frames of
$s^\beta_i$ at $\iota_i(q_i)$ do. The assertion follows as in the proof of
\cref{es:lem:admissible}\textup{(ii)}.

(iv) Let $U\subset D^{i,\alpha}_\beta$ be an open neighbourhood of $q$ on which
$\mathsf e_1:=w_i-P^\alpha_Vh^\alpha_i$ is $\mathsf L$-Lipschitz with respect to $h^\alpha_i$. By
\eqref{es:eq:transferredcoordinate}, $w''_i\circ\iota_i=\xi^\beta\circ\iota\circ\psi_i$ on
$U$. Write $\mathsf e_2=\xi^\alpha-P^\alpha_Vh^\alpha$,
$\mathsf e_3=h^\beta\circ\iota-Rh^\alpha$ and $\mathsf e_4=\xi^\beta-P^\beta_Vh^\beta$.
By \cref{es:lem:limit}, $\mathsf e_2$ and $\mathsf e_4$ are $7\eps_*$-Lipschitz with respect to $h^\alpha$ and $h^\beta$, and $\mathsf e_3$ is $3\eps_*$-Lipschitz with respect to $h^\alpha$.
The identity
\[
h^\alpha\circ\psi_i=h^\alpha_i+\mathsf e_1-\mathsf e_2\circ\psi_i
\]
gives
$\operatorname{Lip}_{h^\alpha_i}(h^\alpha\circ\psi_i)
\le(1+\mathsf L)/(1-7\eps_*)\le13/10$ and
$\operatorname{Lip}_{h^\alpha_i}(h^\alpha\circ\psi_i-h^\alpha_i)
\le\mathsf L+10\eps_*$.
Using $h^\beta\circ\iota=Rh^\alpha+\mathsf e_3$ and
$\xi^\beta=P^\beta_Vh^\beta+\mathsf e_4$, we obtain
\[
\operatorname{Lip}_{h^\alpha_i}
 \bigl(\xi^\beta\circ\iota\circ\psi_i-P^\beta_VR h^\alpha_i\bigr)
 \le\mathsf L+24\eps_* .
\]
For large $i$, $\|RR_i^{-1}-\one\|\le\eps_*$. The comparison of the harmonic charts in \eqref{es:eq:stagetwocharts} now gives
\[
\operatorname{Lip}_{h^\beta_i}
 \bigl(w''_i-P^\beta_Vh^\beta_i\bigr)
 \le\frac{\mathsf L+24\eps_*+\eps_*+3\eps_*}{1-3\eps_*}
 \le\mathsf L+40\eps_* ,
\]
as claimed.
\end{proof}

\subsubsection{Averages of Regular Coordinates}\label{es:sec:averagesreg}

We fix a tube and drop its index. The coordinate of an averaged exponent is a
convex combination up to a remainder. Differences of regular coordinates control
the variation of each exponent multiplied by the component on which it acts;
this weighted control is sufficient even when that component vanishes.

Since $K_0$ is a torus, $V$ is the orthogonal sum of the space of the vectors that are
fixed by $\mathsf R(K_0)$ and of at most $n$ subspaces $V_\lambda$, each with a complex
structure and a linear form $\lambda$ on $\fk$, such that
$\mathsf R(\exp\mathsf b)v=e^{\mathrm i\lambda(\mathsf b)}v$ for $v\in V_\lambda$ and
$\mathsf b\in\fk$. By \eqref{es:eq:crho}, $|\lambda(\mathsf b)|\le C_{\mathsf R}|\mathsf b|$. We use
that $|e^{\mathrm i\theta}-e^{\mathrm i\theta'}|\ge\tfrac2\pi|\theta-\theta'|$ for $|\theta-\theta'|\le\pi$, and we say
that a map is \emph{locally $\mathsf L$-Lipschitz} on a set if every point of the set has a
neighbourhood on which it is $\mathsf L$-Lipschitz with respect to $h_i$.

\begin{lemma}[Averages of Regular Coordinates]\label{es:lem:averages}
Let $U\subset \cC^i_1$ be open, let $w\colon U\to V$ be a continuous map with $|w|\le\ell$,
and let $\eps>0$ with $C_{\mathsf R}\,\eps\le\tfrac18$. For $m$ in a set of indices let
$c_m\colon U\to[0,1]$ be continuous functions whose supports form a locally finite
family and whose sum is $1$, and let $v_m\colon U_m\to\fk$ be continuous maps with
$|v_m|\le\eps$ on open sets $U_m$ that contain the supports of the functions $c_m$.
Suppose that there are numbers $L\ge0$ and $\mathsf L\in[0,1]$ such that
\begin{equation*}
\sum_m|c_m(p)-c_m(q)|\le\frac L\ell\,|h_i(p)-h_i(q)|\qquad\text{for all }p,q\in U,
\end{equation*}
and such that for every $m$ the map $\mathsf R(\exp v_m)\,w-P_Vh_i$ is locally
$\mathsf L$-Lipschitz on $U_m$. Then $\bar v:=\sum_mc_m\,v_m$ is continuous on $U$ with
$|\bar v|\le\eps$, and
\begin{equation*}
\mathsf R(\exp\bar v)\,w-P_Vh_i\quad\text{is locally }
\bigl(\mathsf L+20\,n\,C_{\mathsf R}\,(1+L)\,\eps\bigr)\text{-Lipschitz on }U .
\end{equation*}
\end{lemma}
\begin{proof}
Let $p_0\in U$. There is an open neighbourhood $U_0$ of $p_0$ that meets the supports of
only finitely many functions $c_m$, all of which contain $p_0$; let $\mathcal I_0$ be
the set of their indices. We choose $U_0$ inside the sets $U_m$ for $m\in\mathcal I_0$
and so small that the maps $w_m:=\mathsf R(\exp v_m)\,w$, $m\in\mathcal I_0$, are
$\mathsf L$-Lipschitz on $U_0$ up to $P_Vh_i$. On $U_0$ the sum for $\bar v$ runs over
$\mathcal I_0$; so $\bar v$ is continuous. Fix $m_0\in\mathcal I_0$, points $p,q\in U_0$,
and write $[f]:=f(q)-f(p)$ and $\Delta:=h_i(q)-h_i(p)$. All sums below run over
$\mathcal I_0$.

\emph{The convex combination.} Let $w^c:=\sum_mc_m\,w_m$.
Apply \eqref{th:eq:weightedaverage} to $\mathsf e_m=w_m-P_Vh_i$, centred at $\mathsf e_{m_0}(p)$.
Since $|w_m-w_{m_0}|\le2C_{\mathsf R}\eps\ell$, it gives
$|[w^c-P_Vh_i]|\le(\mathsf L+2C_{\mathsf R} L\eps)|\Delta|$.

\emph{The remainder.} The components of $\mathsf R(\exp\bar v)\,w$, of $w^c$ and of
$w_{m_0}$ in the space of the vectors fixed by $\mathsf R(K_0)$ agree. Fix a subspace
$V_\lambda$, let $\omega:=(w_{m_0})_\lambda$ be the component of $w_{m_0}$ in it, and let
\begin{equation*}
\theta_m:=\lambda(v_m-v_{m_0}),\qquad\bar\theta:=\sum_mc_m\,\theta_m,\qquad
\mathsf r:=e^{\mathrm i\bar\theta}-\sum_mc_m\,e^{\mathrm i\theta_m} .
\end{equation*}
Then $|\theta_m|\le2C_{\mathsf R}\eps\le\tfrac14$, the components of $\mathsf R(\exp\bar v)\,w$ and of
$w^c$ in $V_\lambda$ are $e^{\mathrm i\bar\theta}\omega$ and $\sum_mc_me^{\mathrm i\theta_m}\omega$,
and their difference is $\mathsf r\,\omega$. We estimate $[\mathsf r\,\omega]=\mathsf r(q)[\omega]+[\mathsf r]\,\omega(p)$.
First, $|[\omega]|\le(1+\mathsf L)|\Delta|$, and
\begin{equation*}
|\mathsf r|=\Bigl|\sum_mc_m\bigl(e^{\mathrm i(\theta_m-\bar\theta)}-1-\mathrm i(\theta_m-\bar\theta)\bigr)\Bigr|
\le\tfrac12\max_m|\theta_m-\bar\theta|^2\le8\,C_{\mathsf R}^2\eps^2 .
\end{equation*}
Second, for $m\in\mathcal I_0$ the vector $[e^{\mathrm i\theta_m}]\,\omega(p)$ is the
component in $V_\lambda$ of
\begin{equation*}
[w_m]-\mathsf R\bigl(\exp(v_m-v_{m_0})(q)\bigr)\,[w_{m_0}],
\end{equation*}
and this vector differs from
$(\one-\mathsf R(\exp(v_m-v_{m_0})(q)))P_V\Delta$ by at most $2\mathsf L|\Delta|$; so its norm is at
most $(2\mathsf L+2C_{\mathsf R}\eps)|\Delta|$, and
\begin{equation}\label{es:eq:weightedbound}
|[\theta_m]|\,|\omega(p)|\le\tfrac\pi2\,(2\mathsf L+2C_{\mathsf R}\eps)\,|\Delta|\le4\,|\Delta| .
\end{equation}
Third, view $\mathsf r$ as a function of the phases and of the weights in the
probability simplex. Put $a:=2C_{\mathsf R}\eps$. Its phase derivatives satisfy
$|\partial_{\theta_m}\mathsf r|\le2a c_m$. For a tangent vector $\dot c$ to the simplex,
the identity $\sum_m\dot c_m=0$ gives
\begin{equation*}
D_c\mathsf r[\dot c]
=\sum_m\dot c_m\bigl(\mathrm i\theta_m e^{\mathrm i\bar\theta}-e^{\mathrm i\theta_m}+1\bigr),
\qquad |D_c\mathsf r[\dot c]|\le\tfrac32a^2\sum_m|\dot c_m|.
\end{equation*}
Subtracting the linear terms bounds each coefficient by
\mbox{$|\theta_m||\bar\theta|+\tfrac12|\theta_m|^2\le\tfrac32a^2$}.
Integrate first along the phase segment with weights $c_m(q)$, and then along the
weight segment with phases $\theta_m(p)$. By \eqref{es:eq:weightedbound} and the
weight-variation hypothesis, the resulting bound is
\[
|[\mathsf r]|\,|\omega(p)|\le(16C_{\mathsf R}\eps+6C_{\mathsf R}^2\eps^2L)|\Delta|.
\]
This quadratic cancellation controls the remainder even where $\omega$ vanishes.
Together, and since
$C_{\mathsf R}\eps\le\tfrac18$,
\begin{equation*}
|[\mathsf r\,\omega]|\le\bigl(16\,C_{\mathsf R}^2\eps^2+16\,C_{\mathsf R}\eps+6\,C_{\mathsf R}^2\eps^2L\bigr)|\Delta|
\le C_{\mathsf R}\eps\,(18+2L)\,|\Delta| .
\end{equation*}
Summing over the at most $n$ subspaces $V_\lambda$ and adding the bound for $w^c$ we obtain
the assertion, since $2C_{\mathsf R} L\eps+n\,C_{\mathsf R}\eps\,(18+2L)\le20\,n\,C_{\mathsf R}(1+L)\,\eps$.
\end{proof}

The bound does not depend on the number of the sections that are averaged; this is
needed because they will be indexed by deck transformations. The second lemma moves a
regular coordinate by a deck transformation.

\begin{lemma}[Translates of Regular Coordinates]\label{es:lem:translate}
For all large $i$ the following holds. Let $U\subset \cC^i_1$ be open, let $\nu\in\Gamma_i$
with $\nu U\subset \cC^i_1$, and let $w''\colon U\to V$ be a continuous map with
$|w''|\le9r_\ell$ such that $w''-P_Vh_i$ is locally $\mathsf L$-Lipschitz on $U$, for some
$\mathsf L\le\tfrac12$. Then the map
\begin{equation*}
w_\nu\colon\nu U\to V,\qquad w_\nu(p'):=\mathsf R\bigl(\varkappa_i(\nu,\nu^{-1}p')\bigr)\,w''(\nu^{-1}p'),
\end{equation*}
has the property that $w_\nu-P_Vh_i$ is locally $(\mathsf L+10\eps_*)$-Lipschitz on $\nu U$.
\end{lemma}
\begin{proof}
We have $\nu\in\Gamma_i(\tfrac12t_\ell)$. Let $q,\bar q\in U$ be close to each other, and
put $\Delta:=h_i(\bar q)-h_i(q)$ and $\Delta':=h_i(\nu\bar q)-h_i(\nu q)$. By
\cref{th:lem:transplant}\textup{(iv)}, $|\Delta|\le(1-3\eps_*)^{-1}|\Delta'|$, and by
\cref{es:lem:labelchart}\textup{(ii)} for $t=t_\ell$ the vectors $P_V\Delta'$ and
$\mathsf R(k_\nu)P_V\Delta$ differ by at most $(6\eps_*+\sE_1(t_\ell))|\Delta|$. Now
\begin{equation*}
w_\nu(\nu\bar q)-w_\nu(\nu q)=\mathsf R\bigl(\varkappa_i(\nu,\bar q)\bigr)\bigl(w''(\bar q)-w''(q)\bigr)
+\bigl[\mathsf R\bigl(\varkappa_i(\nu,\bar q)\bigr)-\mathsf R\bigl(\varkappa_i(\nu,q)\bigr)\bigr]\,w''(q).
\end{equation*}
In the first term $w''(\bar q)-w''(q)$ differs from $P_V\Delta$ by at most
$\mathsf L|\Delta|$, and $\mathsf R(\varkappa_i(\nu,\bar q))$ differs from $\mathsf R(k_\nu)$ by at
most $C_{\mathsf R}\sE_0(t_\ell)$ in norm, by \eqref{es:eq:crho} and
\cref{es:prop:cocycle}\textup{(i)}. By the same proposition the matrix in the second term
has norm at most $5\,C_{\mathsf R}\sE_0(t_\ell)C_\#(\ell)\,\ell^{-1}|\Delta|$, and
$|w''(q)|\le\tfrac9{800}\ell$. With \eqref{es:eq:ellone} we get
\begin{equation*}
\bigl|w_\nu(\nu\bar q)-w_\nu(\nu q)-P_V\Delta'\bigr|
\le\bigl(\mathsf L+6\eps_*+\sE_3(\ell)+\sE_4(\ell)\bigr)|\Delta|
\le\frac{\mathsf L+7\eps_*}{1-3\eps_*}\,|\Delta'|\le(\mathsf L+10\eps_*)\,|\Delta'|,
\end{equation*}
since $\mathsf L\le\tfrac12$ and $\eps_*\le\tfrac1{200}$.
\end{proof}

\subsubsection{Weights for a Pair of Tubes}\label{es:sec:pairweights}

The transferred section of \cref{es:lem:transferred} is defined on the set $\iota_i(D^{i,\alpha}_\beta)$,
which is not invariant under the deck transformations of the second tube. To average
its translates we need weights on the second tube that are supported in the translates
of this set. Let $\alpha$ and $\beta$ be tubes such that $\cW^\alpha_i\cap\cW^\beta_i$ is
not empty for all large $i$. We choose points of these intersections, lifts
$q_{0,i}\in D^{i,\alpha}\cap \cC^{i,\alpha}_{1/16}$ and
$p_{0,i}\in D^{i,\beta}\cap \cC^{i,\beta}_{1/16}$ of them (\cref{es:lem:domain}\textup{(iii)}),
and the chart transitions $\iota_i$ from $\alpha$ to $\beta$ with
$\iota_i(q_{0,i})=p_{0,i}$. After passing to a subsequence they converge, and we are in
the situation of \cref{es:sec:transferred}. The maps $\iota_i^{-1}$ are chart transitions from $\beta$ to $\alpha$, so \eqref{es:eq:stagetwocharts} holds for them as well. Since $\eps_*\le\tfrac1{303}$ by
\eqref{es:eq:kappasharp}, this gives
\begin{equation}\label{es:eq:bothways}
\tfrac{99}{100}\,|\Delta|\le\bigl|h^\beta_i(\iota_i(q))-h^\beta_i(\iota_i(\bar q))\bigr|
\le\tfrac{101}{100}\,|\Delta|,\qquad\Delta:=h^\alpha_i(q)-h^\alpha_i(\bar q),
\end{equation}
for $q,\bar q\in \cC^{i,\alpha}_3$, and for $q,\bar q\in\iota_i^{-1}(\cC^{i,\beta}_3)$. Let
$\chi_2\colon[0,\infty)\to[0,1]$ be the function that equals $1$ on $[0,\tfrac15]$,
vanishes on $[\tfrac9{40},\infty)$ and is linear in between, and let
\begin{equation*}
D^{(\alpha)}:=\bigl\{p'\in D^{i,\beta}:\pi^\beta_i(p')\in\cW^\alpha_i(\tfrac{15}2)\bigr\},\qquad
\chi''(p'):=\chi_2\Bigl(\frac{|P^\beta_Th^\beta_i(p')|}{\ell_\beta}\Bigr)\,
\chi_2\Bigl(\frac{|P^\alpha_Th^\alpha_i(\iota_i^{-1}(p'))|}{\ell_\alpha}\Bigr)
\end{equation*}
for $p'\in\iota_i(D^{i,\alpha}_\beta)$, and $\chi''(p'):=0$ for the other points of $D^{i,\beta}$. For
$p'\in D^{(\alpha)}$ and $\nu\in\Gamma^\beta_i$ we put
\begin{equation*}
\varpi''_\nu(p'):=\frac{\chi''(\nu^{-1}p')}{\bar\chi''(p')}\quad\text{if }\nu^{-1}p'\in D^{i,\beta},
\qquad\varpi''_\nu(p'):=0\quad\text{otherwise},
\end{equation*}
where $\bar\chi''(p')$ is the sum of the numbers $\chi''(\nu^{-1}p')$ over the elements
$\nu$ with $\nu^{-1}p'\in D^{i,\beta}$.

\begin{lemma}[Weights for a Pair of Tubes]\label{es:lem:pairweights}
For all large $i$ the following holds.
\begin{enumerate}[label=\textup{(\roman*)}]
\item The restriction of $\chi''$ to $D^{(\alpha)}$ is continuous, and the closure in
$D^{(\alpha)}$ of the set of its points with $\chi''>0$ lies in $\iota_i(D^{i,\alpha}_\beta)$.
\item For every $p'\in D^{(\alpha)}$ we have $\bar\chi''(p')\ge1$. The functions
$\varpi''_\nu$ are continuous on $D^{(\alpha)}$, near a point only finitely many of them
do not vanish, their sum is $1$, and $\varpi''_\nu(\gamma p')=\varpi''_{\gamma^{-1}\nu}(p')$
if $p',\gamma p'\in D^{(\alpha)}$.
\item There is a number $L''$ that does not depend on $i$ such that every point of
$D^{(\alpha)}$ has a neighbourhood on which
\begin{equation*}
\sum_\nu|\varpi''_\nu(p')-\varpi''_\nu(\bar p')|\le\frac{L''}{\ell_\beta}\,|h^\beta_i(p')-h^\beta_i(\bar p')| .
\end{equation*}
\end{enumerate}
\end{lemma}
\begin{proof}
(i) The function is continuous on the open set $\iota_i(D^{i,\alpha}_\beta)$. Let $p'_m\in\iota_i(D^{i,\alpha}_\beta)$
converge to a point $p'\in D^{(\alpha)}$ outside $\iota_i(D^{i,\alpha}_\beta)$. The points
$q_m:=\iota_i^{-1}(p'_m)$ of $D^{i,\alpha}$ converge to $q:=\iota_i^{-1}(p')$, which does
not lie in $D^{i,\alpha}$. By \cref{es:lem:domain}\textup{(i)}, either
$|P^\alpha_Th^\alpha_i(q)|=\tfrac14\ell_\alpha$ or $|\tilde\xi^\alpha_i(q)|=r_V$. The point $q$
lies over a point of $\cW^\alpha_i(\tfrac{15}2)$, so $|\tilde\xi^\alpha_i(q)|<\tfrac{15}2r_\ell$
by \cref{es:lem:domain}\textup{(ii)}. Hence $|P^\alpha_Th^\alpha_i(q_m)|>\tfrac9{40}\ell_\alpha$
for large $m$, and $\chi''(p'_m)=0$.

(ii) Let $p'\in D^{(\alpha)}$ lie over $b$. Suppose that $\ell_\alpha\le\ell_\beta$. Let
$q\in D^{i,\alpha}\cap \cC^{i,\alpha}_{1/16}$ be a lift of $b$, and put $p'':=\iota_i(q)$. By
\eqref{es:eq:bothways}, $|h^\beta_i(p'')-h^\beta_i(p_{0,i})|\le\tfrac{101}{800}\ell_\alpha$, so
$|h^\beta_i(p'')|<\tfrac{19}{100}\ell_\beta$. The point $p''$ is a lift of $b$ in
$\cC^{i,\beta}_1$; so $p''=\nu^{-1}p'$ for some $\nu\in\Gamma^\beta_i$, we have
$|\tilde\xi^\beta_i(p'')|=|\tilde\xi^\beta_i(p')|<r_V$ by \cref{es:lem:domain}\textup{(ii)},
and $p''\in D^{i,\beta}$. Hence $q\in D^{i,\alpha}_\beta$ and $\chi''(p'')=1$. If
$\ell_\alpha>\ell_\beta$, let $p''\in D^{i,\beta}\cap \cC^{i,\beta}_{1/16}$ be a lift of $b$
and $q:=\iota_i^{-1}(p'')$; in the same way $|h^\alpha_i(q)|<\tfrac{19}{100}\ell_\alpha$, the
point $q$ lies in $D^{i,\alpha}$ because it lies over $\cW^\alpha_i(\tfrac{15}2)$, and
$\chi''(p'')=1$. So $\bar\chi''\ge1$. The sums are locally finite, by proper discontinuity, because
$D^{i,\beta}$ lies in a compact subset of $\cC^{i,\beta}_{1/2}$, as in the proof of
\cref{es:lem:sectioncorrection}. The argument there for the extension by zero, together
with (i), shows that the functions $\chi''(\nu^{-1}\cdot)$, extended by zero, are
continuous on $D^{(\alpha)}$; \cref{th:lem:normalizeweights} gives the normalization and
covariance of the weights.

(iii) Let $C_6$ be the universal constant of \eqref{th:eq:stageholder},
let $\ell_-:=\min\{\ell_\alpha,\ell_\beta\}$ and $c'':=\ell_-/(200\,C_6)$, and let $p''$ be
the point of the proof of (ii). If $\gamma\in\Gamma^\beta_i$ and
$\hat d(\gamma p'',p'')<c''$, then $h^\beta_i(\gamma p'')$ and
$h^\alpha_i(\iota_i^{-1}(\gamma p''))$ differ from their values for $\gamma=1$ by at most
$\tfrac1{150}\ell_-$, and the argument of (ii) gives $\chi''(\gamma p'')=1$. So
$\bar\chi''(p')$ is at least the number $N_{p'}(c'')$ of the elements $\gamma$ with
$\hat d(\gamma p',p')<c''$, which is the same for $p'$ and for $p''$. Let $\bar p'$ be a
point of $D^{(\alpha)}$ with $\hat d(p',\bar p')<d_+(\ell_\beta)$. If $\chi''(\nu^{-1}p')$
or $\chi''(\nu^{-1}\bar p')$ is positive, then $\hat d(\nu^{-1}p',p')<6d_+(\ell_\beta)$,
because $D^{i,\beta}$ lies in $B_{d_+(\ell_\beta)}(z^\beta_i)$. Apply
\eqref{th:eq:countdoubling} in the chart of the tube $\beta$ with $t=c''$ and $R=6d_+(\ell_\beta)$; here $t\le R$ and $R+t/2\le9d_+(\ell_\beta)$.
The number of these elements is at most $C''\,\bar\chi''(p')$, with a constant $C''$
that depends on $c''$ and $d_+(\ell_\beta)$ but not on $i$.
Near a point of $D^{(\alpha)}$ each function
$\chi''\circ\nu^{-1}$ either vanishes identically or is the product of two functions
that are Lipschitz with respect to $h^\beta_i$ with the constants $41/\ell_\beta$ and
$42/\ell_\alpha$, by (i), \eqref{es:eq:bothways} and \cref{th:lem:transplant}\textup{(iv)}.
Thus \eqref{th:eq:normalizeweights} gives
\begin{equation*}
\sum_\nu|\varpi''_\nu(p')-\varpi''_\nu(\bar p')|
\le\frac{2\cdot83\,C''}{\ell_-}\,\bigl|h^\beta_i(p')-h^\beta_i(\bar p')\bigr|
\end{equation*}
for $\bar p'$ near $p'$. This is (iii) with $L'':=166\,C''\ell_\beta/\ell_-$.
\end{proof}

\subsubsection{Proof of the Gluing Theorem}\label{es:sec:allproof}

We glue the maps over the tubes of a system. At each step we carry, for every tube,
an admissible section of the current map
with a regular coordinate; the section of a tube changes at every step that concerns
it, by an average as in \cref{es:lem:averages}.

Let an infinite set of indices be given, and let $x_1,\dots,x_J$ be the tubes of
\cref{es:lem:system}. We pass to a common harmonic-chart subsequence as above,
and then to one along which for every
ordered pair of different tubes either $\cW^\alpha_i\cap\cW^\beta_i$ is empty for all
large $i$, or it is not empty for all large $i$ and the transitions of
\cref{es:sec:pairweights} converge, together with the matrices $R_i$ of
\cref{es:sec:transferred}; in the second case we say that the tubes \emph{meet}.
Put
\begin{equation*}
\mathsf L(m):=60\,\eps_*\,m,\qquad\text{so that}\qquad\mathsf L\bigl((n+1)M_n\bigr)+60\,\eps_*=1/4
\end{equation*}
by \eqref{es:eq:kappasharp}. For a tube $x_j$, with its unit $r_\ell$, let
\begin{equation*}
\Sigma^j_i:=\bigl\{b\in\cW^j_i:|\tilde\xi^j_i(b)|\le4r_\ell\bigr\},\qquad
\hat\Sigma^j_i:=\bigl\{b\in\cW^j_i:|\tilde\xi^j_i(b)|\le7r_\ell\bigr\},\qquad
\cD^j_i:=\bigcup_{l\le j}\cW^l_i(2),
\end{equation*}
where $|\tilde\xi^j_i|$ is regarded as a function on $\cW^j_i$. The first two sets are
compact. Let $\chi_4\colon[0,\infty)\to[0,1]$ be the function that equals $1$ on $[0,5]$,
vanishes on $[7,\infty)$ and is linear in between, and let $\chi^j\colon X_i\to[0,1]$ be
the continuous function that equals $\chi_4(|\tilde\xi^j_i|/r_\ell)$ on $\cW^j_i$ and $0$
elsewhere; its support lies in $\hat\Sigma^j_i$. For $b\in X_i$ let $m_j(b)$ be the
number of the indices $l\le j$ with $b\in\hat\Sigma^l_i$. For large $i$,
\begin{equation}\label{es:eq:coverall}
\cD^J_i=X_i\qquad\text{and}\qquad m_J(b)\le(n+1)M_n\quad\text{for all }b\in X_i.
\end{equation}
Indeed, if $b_i$ lay outside every $\cW^j_i(2)$, choose frames $u_i$ over $b_i$
converging to $u\in Y$. By \cref{es:lem:system}\textup{(i)}, $\pr(u)$ lies in
some $O_j(1)$. In this tube, represent $u$ by a frame over a point $q\in S^j$
with $|\xi^j(q)|<r_\ell$. Frame convergence provides lifts $q_i$ of $b_i$ converging
to $q$; then $q_i\in D^{i,j}(2)$ for large $i$, a contradiction.
For the multiplicity bound, otherwise a subsequence has points $b_i$ in
$\hat\Sigma^j_i$ for a fixed set of more than $(n+1)M_n$ indices $j$. In each tube,
\cref{es:lem:domain}\textup{(iii)} gives lifts in $\cC^{i,j}_{1/16}$; their
limits project to the common limit of $b_i$ and satisfy $|\xi^j|\le7r_\ell$.
By \cref{es:lem:slice}\textup{(iii)}, this limit belongs to every corresponding
$O_j(8)$, contradicting \cref{es:lem:system}\textup{(ii)}.

We define realizing maps $F^0,\dots,F^J$ and prove, by induction on $j$, the following
claim.
\begin{enumerate}[label=\textup{(R$_j$)},ref=\textup{(R$_j$)},leftmargin=3.4em]
\item\label{es:Rj} For large $i$ and every tube $x_\alpha$ of the system there is an
admissible section $s^{\alpha,j}_i$ of $F^j$ in the tube $x_\alpha$ that is
$\mathsf L(m_j(b))$-regular at every point of $D^{i,\alpha}$ that lies over a point $b$ of
$\cD^j_i$.
\end{enumerate}
Let $F^0:=F_i$. The set $\cD^0_i$ is empty, so (R$_0$) holds by
\cref{es:lem:admissible}\textup{(i)}. Let $j\ge1$, and assume (R$_{j-1}$). Let $F^j$ be the
map $F'$ of \cref{es:lem:modall} for the maps $F^{j-1}$, the tube $x_j$ and the section
$s^{j,j-1}_i$. The maps $F^j$ are realizing, $\sup\sd(F^j,F^{j-1})\to0$, and
$F^j=F^{j-1}$ outside $\pr^{-1}(\Sigma^j_i)$. We construct the sections of \ref{es:Rj}.

\emph{The tube $x_j$.} Let $s^{j,j}_i:=s^{j,j-1}_i$; by \cref{es:lem:modall} it is an
admissible section of $F^j$ with the coordinate $w'=(1-\chi^\xi_i)w+\chi^\xi_i\tilde\xi^j_i$,
where $w$ is the coordinate of $s^{j,j-1}_i$ for $F^{j-1}$. Let $p\in D^{i,j}$ lie over
$b\in\cD^j_i$. If $|\tilde\xi^j_i(p)|<2r_\ell$, then $w'=\tilde\xi^j_i$ near $p$, and the
section is $7\eps_*$-regular at $p$ by \ref{es:X:lip}; here $m_j(b)\ge1$. Otherwise
$b\notin\cW^j_i(2)$, so $b\in\cD^{j-1}_i$, the number $\mathsf L:=\mathsf L(m_{j-1}(b))$ is at
least $60\eps_*$, and $w-P_Vh_i$ is $\mathsf L$-Lipschitz near $p$. If
$|\tilde\xi^j_i(p)|>4r_\ell$, then $w'=w$ near $p$. If
$2r_\ell\le|\tilde\xi^j_i(p)|\le4r_\ell$, then $b\in\hat\Sigma^j_i$ and $m_j(b)=m_{j-1}(b)+1$,
and \eqref{th:eq:weightedaverage}, applied to $w-P_Vh_i$ and
$\tilde\xi^j_i-P_Vh_i$ and centred at the former, bounds the local Lipschitz constant by
\[
\mathsf L+2\operatorname{Lip}_{h_i}(\chi^\xi_i)\sup|\tilde\xi^j_i-w|
\le\mathsf L+(1+7\eps_*)\eps_*<\mathsf L+2\eps_*.
\]
Here $\operatorname{Lip}_{h_i}(\chi^\xi_i)\le(1+7\eps_*)/(2r_\ell)$,
$\sup|\tilde\xi^j_i-w|\le\eps_*r_\ell$ for large $i$ by
\cref{es:lem:admissible}\textup{(iii)}, and $\mathsf L\ge60\eps_*$. In all cases
\begin{equation}\label{es:eq:tubej}
s^{j,j}_i\text{ is }\mathsf L'(b)\text{-regular at }p,\qquad
\mathsf L'(b):=\begin{cases}7\eps_*&\text{if }b\in\cW^j_i(2),\\
\mathsf L(m_{j-1}(b))+2\eps_*&\text{if }b\in\cD^{j-1}_i\setminus\cW^j_i(2),\end{cases}
\end{equation}
and $\mathsf L'(b)\le\mathsf L(m_j(b))$ if $b\in\hat\Sigma^j_i$, while the section is
$\mathsf L(m_{j-1}(b))$-regular at $p$ if $b\notin\hat\Sigma^j_i$. This is \ref{es:Rj} for the
tube $x_j$.

\emph{A tube that does not meet $x_j$.} Let $\alpha\ne j$ be such that $\cW^\alpha_i$ and
$\cW^j_i$ are disjoint for large $i$. Then $F^j=F^{j-1}$ on $\pr^{-1}(\cW^\alpha_i)$, the
section $s^{\alpha,j}_i:=s^{\alpha,j-1}_i$ is an admissible section of $F^j$, and for the
points $b$ of $\cW^\alpha_i$ we have $m_j(b)=m_{j-1}(b)$, and $b\in\cD^j_i$ only if
$b\in\cD^{j-1}_i$.

\emph{A tube that meets $x_j$.} Let $\alpha\ne j$ be a tube that meets $x_j$, with the
transitions $\iota_i$ from $x_j$ to $x_\alpha$ of \cref{es:sec:pairweights}. Let $s_{\mathrm{ref}}$ be an
admissible section of $F^j$ in the tube $x_\alpha$, with the coordinate $w_0$
(\cref{es:lem:admissible}\textup{(i)}). We compare two kinds of sections with it.

The first is the section $s^{\alpha,j-1}_i$ of (R$_{j-1}$) on the open set $D^{\mathrm o}$
of the points of $D^{i,\alpha}$ that do not lie over $\Sigma^j_i$. Over these points
$F^j=F^{j-1}$, so it has values in $W^{F^j,\alpha}_i$ there, and it has the properties
\ref{es:A:cocycle} and \ref{es:A:frames} at the points of $D^{\mathrm o}$. By the proof of
\cref{es:lem:admissible}\textup{(ii)} there are continuous maps
$\mu^{\mathrm o}\colon D^{\mathrm o}\to\fk_\alpha$ with
\begin{equation*}
s^{\alpha,j-1}_i=A_{\exp(\mu^{\mathrm o})}\cdot s_{\mathrm{ref}}\quad\text{on }D^{\mathrm o},\qquad
\mu^{\mathrm o}(\gamma p')=\Ad_\gamma\,\mu^{\mathrm o}(p')\quad\text{for }p',\gamma p'\in D^{\mathrm o},
\qquad\sup|\mu^{\mathrm o}|\to0 .
\end{equation*}

The second kind are the translates of the transferred section. Let $D^{i,j}_\alpha\subset D^{i,j}$,
$s^\iota_i$, $\mu_i$ and $w''_i$ be the objects of \cref{es:lem:transferred} for the maps $F^j$,
the tubes $x_j$ and $x_\alpha$, the section $s^{j,j}_i$ and the section $s_{\mathrm{ref}}$. For
$\nu\in\Gamma^\alpha_i$ let $U_\nu$ be the open set of the points $p'\in D^{i,\alpha}$ with
$\nu^{-1}p'\in\iota_i(D^{i,j}_\alpha)$, and let
\begin{equation*}
\mu_\nu(p'):=\Ad_\nu\,\mu_i\bigl(\iota_i^{-1}(\nu^{-1}p')\bigr)\qquad(p'\in U_\nu).
\end{equation*}
By \ref{es:A:cocycle} for $s_{\mathrm{ref}}$ and \cref{es:prop:cocycle}\textup{(iii)},
\begin{equation*}
A_{\exp(\mu_\nu(p'))}\cdot s_{\mathrm{ref}}(p')=A_{\varkappa^\alpha_i(\nu,\nu^{-1}p')}\cdot s^\iota_i\bigl(\iota_i^{-1}(\nu^{-1}p')\bigr)
\qquad(p'\in U_\nu);
\end{equation*}
so the left-hand side is a continuous map of $U_\nu$ to $W^{F^j,\alpha}_i$ over
$\pi^\alpha_i$ whose coordinate $\mathsf R(\exp\mu_\nu)\,w_0$ is the map $w_\nu$ of
\cref{es:lem:translate} for $w''=w''_i$. The norms of the maps $\mu_\nu$ are at most $\sup|\mu_i|$, which tends to $0$, and
these maps satisfy $\mu_{\gamma\nu}(\gamma p')=\Ad_\gamma\,\mu_\nu(p')$.

Let $D^{(j)}$ and $\varpi''_\nu$ be the set and the weights of \cref{es:sec:pairweights} for
the pair of $x_j$ and $x_\alpha$, and write $\chi^j$ also for the function
$\chi^j\circ\pi^\alpha_i$ on $D^{i,\alpha}$. Its support lies in $D^{(j)}$, and the set
where it is less than $1$ lies in $D^{\mathrm o}$. We put
\begin{equation*}
u:=(1-\chi^j)\,\mu^{\mathrm o}+\chi^j\sum_\nu\varpi''_\nu\,\mu_\nu\colon D^{i,\alpha}\to\fk_\alpha,
\qquad s^{\alpha,j}_i:=A_{\exp(u)}\cdot s_{\mathrm{ref}} ,
\end{equation*}
where the first term is omitted where $\chi^j=1$ and the second where $\chi^j=0$. By
\cref{es:lem:pairweights}\textup{(i)} and~\textup{(ii)} the map $u$ is well defined and
continuous, $u(\gamma p')=\Ad_\gamma u(p')$ for $p',\gamma p'\in D^{i,\alpha}$, and
$\sup|u|\to0$. By \cref{es:lem:admissible}\textup{(ii)}, $s^{\alpha,j}_i$ is an admissible
section of $F^j$.

Let $p'\in D^{i,\alpha}$ lie over $b\in\cD^j_i$. If $b\notin\hat\Sigma^j_i$, then
$\chi^j=0$ near $p'$, so $s^{\alpha,j}_i=s^{\alpha,j-1}_i$ near $p'$; moreover
$b\in\cD^{j-1}_i$ and $m_j(b)=m_{j-1}(b)$, and the section is $\mathsf L(m_j(b))$-regular
at $p'$ by (R$_{j-1}$). Let $b\in\hat\Sigma^j_i$; then $m_j(b)=m_{j-1}(b)+1$. We apply
\cref{es:lem:averages} on a small neighbourhood $U$ of $p'$ to the map $w_0$, the weight
$1-\chi^j$ with the exponent $\mu^{\mathrm o}$, and the weights $\chi^j\varpi''_\nu$
with the exponents $\mu_\nu$. By \cref{es:lem:pairweights}\textup{(iii)}, \ref{es:X:lip},
\eqref{es:eq:bothways} and \cref{th:lem:transplant}\textup{(iv)} the weights satisfy the
hypothesis of that lemma with a number $L$ that does not depend on $i$, and their
supports lie in $D^{\mathrm o}$ and in the sets $U_\nu$. If $U$ is small, the weights
that do not vanish on $U$ have $p'$ in their supports, and we check the regularity of
the corresponding sections at $p'$.

If $p'$ lies in the support of $1-\chi^j$, then $b\notin\cW^j_i(2)$, so
$b\in\cD^{j-1}_i$, and the section $s^{\alpha,j-1}_i$, whose coordinate is
$\mathsf R(\exp\mu^{\mathrm o})\,w_0$, is $\mathsf L(m_{j-1}(b))$-regular at $p'$. If $p'$ lies
in the support of $\chi^j\varpi''_\nu$, then $q:=\iota_i^{-1}(\nu^{-1}p')$ is a point of
$D^{i,j}_\alpha$ over $b$, and $s^{j,j}_i$ is $\mathsf L'(b)$-regular at $q$ by \eqref{es:eq:tubej}.
Since $\mathsf L'(b)\le\tfrac14$, \cref{es:lem:transferred}\textup{(iv)} shows that
$w''_i-P_Vh^\alpha_i$ is $(\mathsf L'(b)+40\eps_*)$-Lipschitz near $\nu^{-1}p'$, and by
\cref{es:lem:translate} the map $w_\nu-P_Vh^\alpha_i$ is
$(\mathsf L'(b)+50\eps_*)$-Lipschitz near $p'$; here $|w''_i|\le9r_\ell$ for large $i$ by
\cref{es:lem:admissible}\textup{(iii)}. So all these maps are $\mathsf L_1$-Lipschitz near
$p'$ up to $P_Vh^\alpha_i$, with
\begin{equation*}
\mathsf L_1:=\begin{cases}57\,\eps_*&\text{if }b\in\cW^j_i(2)\setminus\cD^{j-1}_i,\\
\mathsf L(m_{j-1}(b))+52\,\eps_*&\text{if }b\in\cD^{j-1}_i .\end{cases}
\end{equation*}
By \cref{es:lem:averages}, with $\eps$ the larger of $\sup|\mu^{\mathrm o}|$ and
$\sup|\mu_i|$, the coordinate $\mathsf R(\exp u)\,w_0$ of $s^{\alpha,j}_i$ is
$(\mathsf L_1+\eps_*)$-regular at $p'$ for large $i$, and
$\mathsf L_1+\eps_*\le\mathsf L(m_{j-1}(b)+1)=\mathsf L(m_j(b))$. This proves \ref{es:Rj}, and the
induction is complete.

\emph{Conclusion.} Let $F^\natural_i:=F^J$. These maps are continuous, they commute with
$O(n)$, and $\sup\sd(F^\natural_i,F_i)\to0$. By \eqref{es:eq:coverall} and (R$_J$), for every
tube $x_\alpha$ the section $s^{\alpha,J}_i$ is $\tfrac14$-regular at every point of
$D^{i,\alpha}$. Every point $u\in Y_i$ lies over a point of some $\cW^\alpha_i$; by
\cref{es:lem:regularall}, applied to $s^{\alpha,J}_i$ on a small neighbourhood of a lift,
$F^\natural_i$ is a topological submersion near $u$. It is proper because
$Y_i$ is compact, so its image is both open and closed in $Y$. Realizing convergence
implies that the image meets every component of $Y$ for large $i$; there are only
finitely many components since $Y$ is a compact manifold. Thus it is surjective.
Its fibres are closed manifolds of dimension
$\dim T+\dim K$, as seen in the local product charts; this is the dimension of
the limit group of the chart. A proper topological submersion whose fibres are compact
manifolds is a fibre bundle \cite[Corollary~6.14]{Sie72}, see also \cite[Essay~II]{KS77};
manifolds satisfy the hypotheses of that corollary on the fibres by
\cite[Example~1.3 and Theorem~2.3]{Sie72}. This proves \textup{(i)} and \textup{(ii)}.

For \textup{(iii)} and \textup{(iv)} let $x_\alpha$ be a tube of the system. We drop the index
$\alpha$ and write $s$ and $w$ for the section $s^{\alpha,J}_i$ and its coordinate. Let
\begin{equation*}
W'_i:=(F^\natural_i)^{-1}\bigl(\Xi([\one,B_V(5r_\ell)])\bigr),\qquad\varphi'\colon W'_i\to B_V(5r_\ell)\quad
\text{with}\quad\Xi([\one,\varphi'(u)])=F^\natural_i(u).
\end{equation*}
By \cref{es:lem:smallerreduction}, with $t=5$ and $t'=6$, $W'_i$ lies over
$\cW_i(6)$ for large $i$. Hence every point of $W'_i$ is
$A_k^{-1}\cdot s(p)$ for some $k\in K$ and $p\in D^i$. The group $K$ acts freely
on it from the right by $u\cdot k:=A_k^{-1}\cdot u$, and
$\varphi'(u\cdot k)=\mathsf R(k)^{-1}\varphi'(u)$. The map $\varphi'$ is proper:
for a compact $C\subset B_V(5r_\ell)$, its preimage is the closed subset
$(F^\natural_i)^{-1}(\Xi([\one,C]))$ of compact $Y_i$.
It has product charts by \cref{es:lem:regularall}. By \cref{es:prop:trivial} there are a closed manifold
$\mathsf F_i$ with a free action of $K$ and a homeomorphism of $W'_i$ onto
$\mathsf F_i\times B_V(5r_\ell)$ that commutes with $K$ and has the second component
$\varphi'$; in particular $\mathsf F_i$ is homeomorphic to $\varphi'^{-1}(0)$, the fibre of
$F^\natural_i$ over $\Xi([\one,0])$. Let $v\in B_V(r_\ell)$, let $K_v$, $N_v$ and $r_\perp$ be as in
\textup{(iv)}, and let $I_v:=\{A_k:k\in K_v\}$, the isotropy group of $\Xi([\one,v])$. The map
$[A,w]\mapsto\Xi([A,v+w])$ is a linear tube of $O(n)\times_{I_v}B_{N_v}(r_\perp)$ onto an open
invariant subset of $Y$, because
$O(n)\times_I(K\times_{K_v}(v+B_{N_v}(r_\perp)))$ is $O(n)\times_{I_v}(v+B_{N_v}(r_\perp))$. By
\cref{ef:prop:local} the preimage of its image under $F^\natural_i$ is the image of
$O(n)\times_{I_v}\varphi'^{-1}(v+B_{N_v}(r_\perp))$ under $[A,u]\mapsto A\cdot u$, and
$\varphi'^{-1}(v+B_{N_v}(r_\perp))$ corresponds to $\mathsf F_i\times(v+B_{N_v}(r_\perp))$, compatibly
with $K_v$.

\emph{Part \textup{(iii)}.} Let $x\in X$. By \cref{es:lem:system}\textup{(i)} there are a tube
$x_\alpha$ of the system and a vector $v_x\in V_\alpha$ with $|v_x|<r_\ell$ such that $x$ is the
point $\bar\xi_\alpha^{-1}(\mathsf R(K_\alpha)v_x)$. The point $y:=\Xi([\one,v_x])$ lies over $x$, its
isotropy group is $I_{v_x}$, and $\Xi_x([A,w]):=\Xi([A,v_x+w])$ defines the linear tube of
the preceding paragraph for $v=v_x$, about the orbit of $y$. These data
do not depend on $i$; they are the data of \textup{(iii)}, with $N_{v_x}$ as the space $V$ and
the action of $I_{v_x}$ through $\mathsf R$. The preceding paragraph gives \textup{(iii)}, for the
manifold $\mathsf F_i$ with the action of $K_{v_x}$; its index bound is the one for the
finitely many tubes of the system, so it does not depend on $x$.

\emph{Part \textup{(iv)}.} The quotient by $O(n)$ of the preimage of the image of
$[A,w]\mapsto\Xi([A,v+w])$ under $F^\natural_i$ is $(f^\natural_i)^{-1}(O^v)$, by
$\pr\circ F^\natural_i=f^\natural_i\circ\pr$ and \eqref{es:eq:prtube}, and the quotient of
$O(n)\times_{I_v}(\mathsf F_i\times(v+B_{N_v}(r_\perp)))$ by $O(n)$ is
$\mathsf F_i\times_{K_v}B_{N_v}(r_\perp)$; under these identifications $f^\natural_i$ is the map
$[u,w]\mapsto\bar\xi^{-1}(\mathsf R(K)(v+w))$. The space $N_v$ is the orthogonal complement of the
tangent space of the orbit $\mathsf R(K)v$, whose dimension is $\dim K-\dim K_v$. The fibre over
$\bar\xi^{-1}(\mathsf R(K)v)$ is $\mathsf F_i/K_v$. It is a closed manifold of dimension
$\dim\mathsf F_i-\dim K_v$: in a product chart $(q,k)\mapsto A_k^{-1}\cdot s(q)$ of $\varphi'$
on $U\times K$, as in \cref{es:lem:regularall}, the points over this fibre are those with
$w(q)\in\mathsf R(K)v$. They form a submanifold of $U$ of dimension $\dim T+\dim K-\dim K_v$ in
the coordinates $P_Th_i+w$, and $\pi_i$ maps it homeomorphically onto an open subset of the
fibre. \qed

\section{Fibre Lattices and Unwrapped Coverings}\label{ch:lattices}

The bundles of \cref{es:thm:gluingall} exist, but their fibres have not yet been identified. We construct their algebraic models: marked Malcev lattices and period kernels (\cref{m:book:ch:12}), boxes that compare small deck transformations with small lattice elements (\cref{fl:sec:boxes}), and unwrapped coverings with lattice monodromy (\cref{fl:sec:covering}). \Cref{ch:fibregroups} uses these coverings to identify the fibres and their isotropy quotients.

Throughout this section we are in the setting of \cref{m:fr:sec:covercollapse}:
\cref{hyp:sequence} with $\delta\le\delta_*$ and $0<r\le c_*\rho$, the centres and charts fixed
there, and the subsequence of \cref{m:fr:thm:geomain}. For a chart with limit group $G$ we
write $d=\dim G_0$ for the dimension of the identity component of $G$; it is the Hirsch length
of the fibre lattice and the dimension of its nilmanifold (\cref{fl:lem:nilmanifold}).
Realizing maps are those of \cref{fl:def:realizing}.

\subsection{Nilpotent Presentations and Periods}\label{m:book:ch:12}\label{m:fr:sec:nilprogressions}

The small deck transformations will be marked by a torsion-free nilpotent lattice
$\Lambda_i$, the fibre lattice, whose marking kernel consists of the periods
(\cref{m:fr:ext:physical-lattice-realization}). We construct the lattice using nilprogressions and approximate
groups, prove convergence of the normalized period kernels, and describe the nilmanifold of the
lattice (\cref{fl:sec:nilmanifold}).

Fix a non-principal ultrafilter $\mathfrak u$ on $\mathbb N$: a collection of infinite subsets, closed under finite intersections and passage to supersets, containing exactly one of $A$ and $\mathbb N\setminus A$ for every $A\subseteq\mathbb N$. Its members are called \emph{$\mathfrak u$-large}. A statement holds \emph{for $\mathfrak u$-almost every $i$} if its set of valid indices belongs to $\mathfrak u$. The notation $x_i\to x$ as $i\to\mathfrak u$ means that $\{i:x_i\in U\}\in\mathfrak u$ for every neighbourhood $U$ of $x$: this is convergence along the ultrafilter. Every cofinite set is $\mathfrak u$-large, so ordinary convergence implies convergence along $\mathfrak u$ to the same limit.

We use the definitions of an ultra approximate group with no small subgroups (NSS) and of a non-degenerate ultra nilprogression in \cite[Definitions~2.58 and~2.61]{Zam24}; the definition of a good Lie model is recalled below. For a progression
\mbox{$P=P(u_1,\ldots,u_d;N_1,\ldots,N_d)$} in $C$-normal form, its ordered grid is
\[
 G(\varepsilon P)=
 \{u_1^{m_1}\cdots u_d^{m_d}:|m_j|\le\varepsilon N_j/C\}.
\]
The sharp sign below denotes the embedding of the grid into the Malcev lattice of \cite[Theorem~2.53]{Zam24}, which preserves products in the sense recalled there.

\subsubsection{Nilprogressions in the Deck Groups}
\label{m:fr:ext:physical-nilprogressions}

We show that the small deck transformations form approximate groups without small subgroups,
with the limit group as a good Lie model, and we derive from this the fibre lattice with its
marking and its polynomial coordinates. Fix one chart $(Z_i,z_i,\Gamma_i)\eGHto(Z,z,G)$ of \ref{fl:S:limit}, with approximations $(\af_i,\alpha_i,\beta_i)$ as in \cref{m:fr:sec:conv}, and suppress its centre superscript. The deck action on the completed local cover $Z_i$ is free by \cref{m:fr:prop:completion}; this freeness is what excludes small subgroups in \cref{m:fr:ext:deck-NSS}. The threshold $\delta_*$ of \cref{m:fr:sec:coverbundle} satisfies ${\delta_*\le\delta_C\le\delta_U(n)}$, by its definition and the proof of \cref{m:fr:prop:coverscharts}, and the ball $B_{R_1}(z_i)$ lies in the region $\Omega_i$; so \cref{m:fr:cor:localUIC} with
$R=R_1$ gives maps $\Phi_i\colon B_1(0^n)\to B_{R_1}(z_i)$ with the bounds
\eqref{m:fr:eq:reifholder}, whose constants do not depend on $i$, and with images that
contain $B_{R_1/2}(z_i)$. Write
\[
 D_i(\gamma)=\sup_{q\in\overline B_{R_1/2}(z_i)}d_i(\gamma q,q),\qquad
 D(g)=\sup_{q\in\overline B_{R_1/2}(z)}d(gq,q),\qquad
 S_i(t)=\{\gamma\in\Gamma_i:D_i(\gamma)<t\}.
\]
Distances in this paragraph are the completed distances, which agree with the original
distances on $B_{R_1}(z_i)$ by \cref{m:fr:prop:completion}.

\begin{lemma}[Uniform Exclusion of Small Deck Subgroups]
\label{m:fr:ext:deck-NSS}
There is $\theta_{\rm top}>0$, independent of $i$, such that $S_i(\theta_{\rm top})$ contains no
non-trivial subgroup for all sufficiently large $i$.
\end{lemma}
\begin{proof}
For every isometry, the triangle inequality gives
\[
 D_i(\gamma\gamma')\le D_i(\gamma)+D_i(\gamma'),\qquad
 D_i(\gamma^{-1})=D_i(\gamma).
\]
These statements do not require the ball to be invariant.  A subgroup
$H\subset S_i(\theta_{\rm top})$ is finite: its elements have bounded displacement at
$z_i$, and a discrete isometry group acts properly on the proper space
$Z_i$.

Let $v_i:=\Phi_i^{-1}(z_i)$. By \eqref{m:fr:eq:reifholder} and the bound for the distance of
$v_i$ from the boundary of $B_1(0^n)$ in \cref{m:fr:cor:localUIC}, there is $s>0$ with
$\overline B_{2s}(v_i)\subset B_1(0^n)$ and ${\Phi_i(\overline B_{2s}(v_i))\subset
B_{R_1/8}(z_i)}$. Choose $\theta_{\rm top}<R_1/8$ so small that
$C_U(2\theta_{\rm top}/R_1)^{\beta_U}<s$; then two points of $\Phi_i(B_1(0^n))$ at distance
less than $2\theta_{\rm top}$ have preimages at distance less than $s$. Neither $s$ nor
$\theta_{\rm top}$ depends on $i$.
The open saturation
$H \Phi_i(B_{3s/2}(v_i))$ lies in $B_{R_1/4}(z_i)$ and is $H$-invariant.
Its pull-back is an invariant open Euclidean domain containing
$\overline B_s(v_i)$. Every $H$-orbit there has parameter diameter
less than $s$: if $y=hx$, with $x$ in the original image, then
$d_i(ky,y)=d_i(h^{-1}khx,x)<\theta_{\rm top}$ for $k\in H$.

If $H$ is non-trivial, Cauchy's theorem gives an element of prime order.
Apply the local form of Newman's theorem \cite[Theorem~2.7.4]{Tao14} to its action in
this coordinate domain; it applies to homeomorphisms, so the parametrization need not be
differentiable. After translating and rescaling the parameter ball, that theorem
says that a prime-period homeomorphism whose orbits all have diameter
less than the ball radius fixes its centre.  This contradicts the free
deck action on $Z_i$.  Hence $H$ is trivial.
\end{proof}

We now recall the definition of a good Lie model. For a subset $S$ of a group and an integer
$m\ge1$, write $S^m:=\{s_1\cdots s_m:s_j\in S\}$ for its $m$-fold product set.
If $1\in S$, this also includes all products of fewer than $m$ elements.

\Needspace{5\baselineskip}
\begin{definition}[Good Lie Models]\label{m:fr:def:goodmodel}
Let $S=\prod_{\mathfrak u}S_i$ be an ultra approximate group. An ultraproduct
$\prod_{\mathfrak u}A_i$ of subsets $A_i$ of the groups that contain the sets $S_i$ is called
\emph{internal}. A \emph{good Lie model} for $S$ \cite[Definition~2.60]{Zam24} is a connected
local Lie group $G'$ and a local-group morphism $\sigma\colon S^8\to G'$ preserving identity,
inverses and products whenever the factors and their product lie in $S^8$, with these
properties:
\begin{enumerate}[label=\textup{(M\arabic*)},ref=\textup{(M\arabic*)},leftmargin=3.2em]
\item\label{m:fr:model:compact} The image $\sigma(S)$ has compact closure in $G'$.
\item\label{m:fr:model:neighbourhood} There is an open identity neighbourhood $U_0\subset G'$
with $U_0\subset\sigma(S)$ and $\sigma^{-1}(U_0)\subset S$.
\item\label{m:fr:model:internal} If $F\subset U\subset U_0$, with $F$ compact and $U$ open, then
some internal set $S'\subset S$ satisfies \mbox{$\sigma^{-1}(F)\subset S'\subset\sigma^{-1}(U)$}.
\end{enumerate}
All preimages are taken in $S^8$.
\end{definition}

Together with the counting estimates of \cref{m:fr:lem:counting}, \cref{m:fr:ext:deck-NSS} makes
the limit group a good Lie model for the sets of small deck transformations.

\begin{proposition}[A Good Model without Small Subgroups]
\label{m:fr:ext:physical-good-model}
For a sufficiently small fixed $\theta>0$, the finite sets $S_i:=S_i(\theta)$ are
$C$-approximate groups for one $C$, contain no non-trivial subgroups, and
their ultraproduct $S$ has the connected nilpotent Lie group $G_0$ as a
good Lie model
\[
 \sigma\colon S^8\longrightarrow G_0,\qquad
 \sigma((\gamma_i)_{\mathfrak u})=\lim_{i\to\mathfrak u}\alpha_i(\gamma_i).
\]
\end{proposition}
\begin{proof}
The convergence of group elements and equivariance on the ball imply
$D_i(\gamma_i)\to D(g)$ whenever $\gamma_i\eGHto g$.  Indeed, distance
distortion and the density of the approximations compare the two suprema;
the displacement functions have Lipschitz constant at most two.  The
same argument gives convergence uniformly for bounded sets of convergent
group elements.  The function $D$ is continuous and proper on $G$,
since $D(g)\ge d(gz,z)$.  It vanishes only at the identity: an isometry
fixing an open set is the identity, as recorded in
\ref{m:fr:R:isometries}.

Choose $\theta<\theta_{\rm top}$ so small that $\{D\le8\theta\}\subset G_0$ and
$\{D\le4\theta\}$ lies in a Lie-group neighbourhood containing no non-trivial
subgroup.  The sets $S_i$ are finite and symmetric.  To check the
approximate-group bound we use Ruzsa's covering argument \cite[Lemma~2.14]{TV06}: choose a
maximal family of disjoint translates
$x_jS_i(\theta/4)$ with $x_j\in S_i(2\theta)$.  They lie in $S_i(9\theta/4)$, and
maximality gives
\[
 S_i(2\theta)\subset\bigcup_jx_jS_i(\theta/2)
             \subset\bigcup_jx_jS_i(\theta).
\]
The number of translates is at most
$\#S_i(9\theta/4)/\#S_i(\theta/4)$.  Its numerator divided by $\ell_i$ is bounded
by the based-displacement counting estimate of \cref{m:fr:lem:counting}.
Its denominator divided by $\ell_i$ is bounded below by a positive
constant: apply the convergence \eqref{m:fr:eq:haarlimit} of the normalized counting measures
(\cref{m:fr:lem:counting}\textup{(ii)}) to a non-zero non-negative compactly supported continuous function on $G$
supported in $\{D<\theta/8\}$.  Equivariance places the support of its pull-back to $\Gamma_i$ in $S_i(\theta/4)$ for large $i$.  Its positive Haar integral gives the lower
bound.  Thus $S_i^2\subset S_i(2\theta)$ is covered by $C$ translates of $S_i$.
\cref{m:fr:ext:deck-NSS} shows that there are no small subgroups.

Products and inverses commute with the ultralimit by
\cref{m:fr:lem:groupconv}(iv), so $\sigma$ is a morphism of local groups.
All its arguments have displacement at most $8\theta$ in the limit, and hence
its image is in $G_0$. We verify \ref{m:fr:model:compact}--\ref{m:fr:model:internal}.
The image $\sigma(S)$ is
precompact, because $\{D\le \theta\}$ is compact.  A sufficiently small open
identity neighbourhood $U_0\subset\{D<\theta/2\}$ is contained in
$\sigma(S)$, using the approximants $\beta_i(g)$; moreover
$\sigma^{-1}(U_0)\subset S$, by convergence of $D_i$.  Finally, if
$F\subset U\subset U_0$, with $F$ compact and $U$ open, choose an open
$U'$ with $F\subset U'\subset\overline{U'}\subset U$.  The finite sets
\[
 S'_i=\{\gamma\in S_i:\alpha_i(\gamma)\in U'\}
\]
have ultraproduct $S'$ satisfying
$\sigma^{-1}(F)\subset S'\subset\sigma^{-1}(U)$. This verifies \ref{m:fr:model:internal}.
\end{proof}

With a good Lie model in hand, the structure theorem of \cite{Zam24} provides nilprogressions
in the deck groups. The Malcev lattice attached to such a progression maps onto the group
generated by it; this map is the marking, and the lattice is the fibre lattice.

\begin{theorem}[The Fibre Lattice and Its Marking]
\label{m:fr:ext:physical-lattice-realization}\label{fl:N:lattice}\label{fl:N:marking}
After passing to a subsequence there are nilprogressions
\[
 P_i=P(u_{1,i},\ldots,u_{d,i};N_{1,i},\ldots,N_{d,i})\subset\Gamma_i
\]
in one $C_P$-normal form, with $\min_jN_{j,i}\to\infty$, and torsion-free Malcev lattices $\Lambda_i$
of Hirsch length $d$, the \emph{fibre lattices}, with surjective homomorphisms, the \emph{markings},
\[
 \vartheta_i\colon\Lambda_i\longrightarrow
 H_i:=\langle u_{1,i},\ldots,u_{d,i}\rangle\subset\Gamma_i.
\]
The elements of $\ker\vartheta_i$ are the \emph{periods} of the marking. For sufficiently small fixed $\varepsilon>0$, the marking is injective on the grid $G(\varepsilon P_i)^\sharp$; every $\gamma\in\Gamma_i$ with $D_i(\gamma)$ sufficiently small belongs to $G(\varepsilon P_i)$, and so does every member of a sequence $\gamma_i\in\Gamma_i$ that converges to the identity, for large $i$.  For one fixed $b>0$ and all large $i$,
\[
 H_i=\langle S_i(b)\rangle,
 \qquad \ker\vartheta_i\cap G(\varepsilon P_i)^\sharp=\{1\}.
\]
\end{theorem}

Thus the group generated by the deck transformations of small displacement is a nilpotent
quotient of a lattice of Hirsch length $d$. The marking need not be injective; its kernel, the
periods, is the subject of \cref{m:fr:ext:period-kernels}.

\begin{proof}
Apply \cite[Theorem~7.1]{Zam24} to \cref{m:fr:ext:physical-good-model}.
It gives a non-degenerate ultra nilprogression $P$ of rank $d$ in $S^4$, an ultraproduct of
nilprogressions $P_i\subset\Gamma_i$. We need that the preimage of a neighbourhood of the identity lies in the grid. For fixed
small $\varepsilon$, the final paragraph of the proof of that theorem shows
$\sigma^{-1}(1)\subset G(\varepsilon P)$
\cite[proof of Theorem~7.1, p.~791]{Zam24}.
We claim that a fixed displacement neighbourhood is captured: $S_i(b)\subset G(\varepsilon P_i)$ for some $b>0$, depending on $\varepsilon$, and $\mathfrak u$-almost every $i$. Otherwise all the nested sets
\[
 E_k:=\{i:S_i(\min\{\theta/2,1/k\})\not\subset G(\varepsilon P_i)\}
 \qquad(k\ge1)
\]
would be $\mathfrak u$-large. For $i\in E_1$, let $k_i$ be the largest $k\le i$ with
$i\in E_k$, and choose a witness
$\gamma_i\in S_i(\min\{\theta/2,1/k_i\})\setminus G(\varepsilon P_i)$;
put $k_i=1$ and $\gamma_i=1$ otherwise.
Then $k_i\to\infty$ along $\mathfrak u$, so $D_i(\gamma_i)\to0$ and
$\sigma((\gamma_i)_{\mathfrak u})=1$. Yet $\gamma_i\notin G(\varepsilon P_i)$ on the
$\mathfrak u$-large set $E_1$, contradicting the kernel containment. Thus, for some fixed
$0<b<\theta$, depending on $\varepsilon$,
$S_i(b)\subset G(\varepsilon P_i)$ on a $\mathfrak u$-large set. All the witnesses lie in
$S_i$, so the argument stays inside the domain of $\sigma$.
The strong basis in \cite[Remark~7.3]{Zam24} makes every thickness
tend to infinity along $\mathfrak u$. Indeed, a $\mathfrak u$-bounded
$N_{j,i}$ would make $\lfloor tN_{j,i}/C_P\rfloor=0$ for every
sufficiently small standard $t$, forcing $\exp(tv_j)=1$ on an
interval and hence $v_j=0$, contrary to the basis property.
The goodness of the progression and the product-preserving embedding of the grid are \cite[Theorem~2.53]{Zam24}. The defining commutator relations of the abstract lattice hold for the elements $u_{j,i}$ of $\Gamma_i$; assigning the generators to these elements therefore defines $\vartheta_i$. Its image is $H_i$, and the embedding of the grid is a right inverse of $\vartheta_i$ on the grid, which gives the stated intersection with the kernel.

For every fixed integer $m$, non-degeneracy puts $u_j^m$ in $P\subset S^4$.
Consequently every power of $\sigma(u_j)$ lies in $\{D\le4\theta\}$, a neighbourhood of the identity in the Lie group $G_0$ that contains no non-trivial subgroup.  Thus $\sigma(u_j)=1$.  It follows that
$u_{j,i}\in S_i(b)$ for $\mathfrak u$-almost every $i$.  This gives both
inclusions in $H_i=\langle S_i(b)\rangle$.

The inclusions just proved, increasing thickness, and all the convergence statements below
hold along an ordinary further subsequence, by diagonal selection from
$\mathfrak u$-large sets. Include the capture inclusions for a countable sequence
$\varepsilon\downarrow0$; monotonicity of the grids then gives capture for every fixed
small $\varepsilon$, with its own $b$.
If $d=0$, the construction has trivial progression, lattice and group $H_i$, with the evident conventions for empty coordinates.
\end{proof}

We write $u_1,\dots,u_d$ for the generators of $\Lambda_i$ that the marking maps to
$u_{1,i},\dots,u_{d,i}$. The defining relators of $\Lambda_i$ are
\begin{equation}\label{fl:eq:relators}
[u_j,u_k]\,\Bigl(\prod_{l>k}u_l^{c_{jkl}}\Bigr)^{-1}\qquad(1\le j<k\le d),
\end{equation}
with integers $c_{jkl}$ depending on $i$ and satisfying
$|c_{jkl}|\le C_P\,N_{l,i}/(N_{j,i}N_{k,i})$
\cite[Definition~2.52 and Theorem~2.53]{Zam24}.

The lattice is now described in coordinates. The normal forms of its elements, normalized by
the lengths of the progression, embed it into $\R^d$ with a polynomial group law, and these
laws converge to the law of the universal covering of $G_0$.

\begin{proposition}[Polynomial Coordinates]
\label{m:fr:ext:physical-polynomial-coordinates}\label{fl:N:coord}
For the lattices $\Lambda_i$ above, after passing to a further subsequence, the following
holds. Every element of $\Lambda_i$ has a unique normal form $u_1^{m_1}\cdots u_d^{m_d}$ with
$m\in\Z^d$. There are normalizations $L_{j,i}\to\infty$ of the coordinates, with
$L_{j,i}/N_{j,i}$ converging to a positive number that does not depend on $j$, and
polynomial multiplication laws $Q_i$ on $\mathbb R^d$, of degree bounded in terms of $d$,
such that $u_1^{m_1}\cdots u_d^{m_d}\mapsto L_i^{-1}m$, with
$L_i:=\operatorname{diag}(L_{1,i},\dots,L_{d,i})$, is an injective homomorphism of $\Lambda_i$
into $(\R^d,Q_i)$ and $Q_i\to Q$ coefficientwise. The limiting law $Q$ is that of the
universal cover $\pi\colon\widetilde G_0\to G_0$ in coordinates
$\varphi\colon\R^d\to\widetilde G_0$ of the second kind for a strong Malcev basis. The laws
are globally triangular, their polynomial inverses converge coefficientwise, and Lebesgue
measure is left and right Haar measure for all these laws. Moreover there are a coordinate rectangle
$\mathsf X\subset\R^d$ about $0$ and a neighbourhood $\mathcal V$ of the identity in $G_0$ with the
following properties.
\begin{enumerate}[label=\textup{(\roman*)}]
\item If elements of the grid $G(\varepsilon P_i)^\sharp$ with exponent vectors $m_i$ satisfy
$L_i^{-1}m_i\to x\in\mathsf X$, then their markings converge to $\pi(\varphi(x))$.
\item If elements $\gamma_i\in\Gamma_i$ converge to $g\in \mathcal V$, then for large $i$ the element
$\gamma_i$ is the marking of an element of the grid $G(\varepsilon P_i)^\sharp$ whose exponent
vector $m_i$ satisfies $L_i^{-1}m_i\to x$, where $x$ is the point of $\mathsf X$ with
$\pi(\varphi(x))=g$.
\end{enumerate}
Both convergences are uniform on compact subsets of $\mathsf X$ and of $\mathcal V$, and $\mathsf X$ can
be chosen such that the products and inverses of its points lie in a larger coordinate
rectangle.
\end{proposition}
\begin{proof}
\cite[Remark~7.3]{Zam24} gives a strong Malcev basis $v_1,\ldots,v_d$ with
\[
 \exp(tv_j)=\sigma(u_j^{\lfloor tN_j/C_P\rfloor})
 \quad(0\le t\le1).
\]
Fix
$\varepsilon$ below the grid-embedding threshold so small that
$\pi$ is injective on the coordinate rectangle
corresponding to $G(2\varepsilon P)$, and let $b$ be its capture constant from
\cref{m:fr:ext:physical-lattice-realization}. Then choose $\delta_{\rm c}\ll\varepsilon$
so that $G(\delta_{\rm c} P)^2\subset G(\varepsilon P)$, and use the
normalizations of \cite[Section~8]{Zam24},
\[
 L_{j,i}=N_0\Big\lfloor\frac{\delta_{\rm c} N_{j,i}}{C_PN_0}\Big\rfloor,
 \qquad
 \varphi(x)=\exp(\delta_{\rm c} x_1v_1)\cdots\exp(\delta_{\rm c} x_dv_d).
\]
Here $N_0$ is a fixed sufficiently large interpolation denominator, and
$L_{j,i}/N_{j,i}\to\delta_{\rm c}/C_P$.
The map $\varphi\colon\mathbb R^d\to\widetilde G_0$ is a diffeomorphism.
If $\widetilde Q_i$ denotes the Malcev polynomial of the lattice, put
\[
 Q_i(x,y)=L_i^{-1}\widetilde Q_i(L_ix,L_iy),\qquad
 Q(x,y)=\varphi^{-1}(\varphi(x)\varphi(y)).
\]
First, the strong-basis formula holds for exponents that vary with $i$.
If $h_i=u_{j,i}^{e_i}$ with $e_i=o(N_{j,i})$, every fixed power
$h_i^m$ belongs to $P_i\subset S_i^4$ for all large $i$.
Every limit $h$ of a subsequence therefore has all its integer powers in $\{D\le4\theta\}$, a neighbourhood of the identity in $G_0$ without non-trivial subgroups, so $h=1$. Consequently $u_{j,i}^{o(N_{j,i})}\eGHto1$ in the convergence of group elements of \cref{m:fr:lem:groupconv}.  For exponents with
$C_Pm_{j,i}/N_{j,i}\to t$, subtract
$\lfloor tN_{j,i}/C_P\rfloor$; the difference is $o(N_{j,i})$.
The strong-basis formula and convergence of products give
$u_{j,i}^{m_{j,i}}\eGHto\exp(tv_j)$.  Negative $t$ follow by inversion.
Initially this is an ultrafilter statement for every fixed $t$.
Compactness of a fixed normalized exponent interval and the same
$o(N_{j,i})$ argument upgrade it to a uniform ultrafilter statement:
a violating choice has an ultralimit $t$ and reduces to that fixed $t$.
Selecting indices where the uniform error tends to zero makes it an
ordinary uniform statement.

It follows that the marking of every element of the grid whose normalized coordinates converge to a point $x$ of the rectangle converges to $\pi(\varphi(x))$.  On a fixed rational interpolation
grid for $(x,y)$, with denominator dividing $N_0$, the vectors $L_ix$
and $L_iy$ are integral.  The product of the corresponding lattice elements belongs to the larger grid. Its exponents, which are unique, have bounded normalized coordinates, and the limit of its marking is $\pi(\varphi(x)\varphi(y))$.  The local injectivity chosen above forces
every coordinate cluster point to be $Q(x,y)$.  Hence
$Q_i(x,y)\to Q(x,y)$ at every interpolation point.  The common degree
bound in \cite[Theorem~2.53]{Zam24} and finite-dimensional polynomial interpolation
now give coefficient convergence, as in \cite[Proposition~8.2]{Zam24}.

The relators \eqref{fl:eq:relators} give the global triangularity needed here: they show
that $[u_j,u_k]$ lies in $\langle u_{k+1},\ldots,u_d\rangle$ for $j<k$, which makes
$(u_1,\ldots,u_d)$ a Malcev basis, and in the corresponding coordinates of the second kind the group law
$\widetilde Q_i$ of the Malcev completion is a polynomial of triangular form, its inverse is
polynomial, and Lebesgue measure is a bi-invariant Haar measure \cite[Chapter~1]{CG90}.
Diagonal scaling preserves these properties, so
\[
 (Q_i(x,y))_j=x_j+y_j+
 P_{i,j}(x_1,\ldots,x_{j-1},y_1,\ldots,y_{j-1})
\]
for polynomials $P_{i,j}$, and the inverse $J_i$ of $Q_i$ is given by the recursion
\[
 (J_i(x))_j=-x_j-P_{i,j}\bigl(x_{<j},(J_i(x))_{<j}\bigr).
\]
So the inverse laws have degrees
bounded in terms of $d$ and the common product degree, their coefficients converge by finite
polynomial composition, and coefficient convergence gives convergence of every derivative on
compact coordinate sets.

Choose an open coordinate rectangle $\mathsf X$ about $0$ whose closure is inside the
coordinate rectangles defining the normalized grids for all large $i$. The convergence just proved
gives \textup{(i)}. For \textup{(ii)}, choose an identity neighbourhood
\[
 \mathcal V\Subset\{D<b/2\}\cap\pi(\varphi(\mathsf X)).
\]
If $\gamma_i\in\Gamma_i$ converges to $g\in\mathcal V$, then $D_i(\gamma_i)\to D(g)<b/2$,
so capture gives $\gamma_i\in G(\varepsilon P_i)$ for large $i$. Its unique grid exponents satisfy
$|m_{j,i}|\le\varepsilon N_{j,i}/C_P$.  The normalized vectors
$L_i^{-1}m_i$ are bounded.  Every convergent subsequence has limit $x$
in the larger injectivity rectangle, with $g=\pi(\varphi(x))$. Hence $x$ is the unique
preimage of $g$ in $\mathsf X$, and the full sequence converges to it. This proves
\textup{(ii)}.

Choose nested compact coordinate rectangles whose limiting products and inverses lie
in the interiors of larger ones.  Convergence of $Q_i$ and of their inverses makes these inclusions valid for large $i$. Uniformity of the convergence in both directions on these rectangles follows from the same compactness argument, and can be included in the selection of the common subsequence.
\end{proof}

\subsubsection{Normalized Coordinates and Convergence of the Period Kernels}
\label{m:fr:ext:period-kernels}

We compare the lattices with the universal covering of $G_0$ in the normalized coordinates.
Lattice elements with convergent coordinates have convergent markings, and convergent deck
transformations lift to lattice elements with convergent coordinates; together these show
that the kernels of the markings converge to the kernel of the covering homomorphism. Write $\mathcal N_i=(\mathbb R^d,Q_i)$ and
$\mathcal N=(\mathbb R^d,Q)$, with inversions $J_i$ and $J$, and identify $\mathcal N$ with
$\widetilde G_0$ by $\varphi$.  Let
\[
 \mathbf{c}_i\colon\Lambda_i\longrightarrow\mathcal N_i,\qquad
 \mathbf{c}_i(u_1^{m_1}\cdots u_d^{m_d})=L_i^{-1}m,
 \qquad
 \pi\colon\mathcal N\longrightarrow G_0
\]
be the normalized coordinates of the lattice and the covering homomorphism. Define the normalized period groups by
\[
 \mathsf K_i:=\mathbf{c}_i(\ker\vartheta_i),\qquad \mathsf K:=\ker\pi.
\]
As the kernel of a covering homomorphism of connected Lie groups, $\mathsf K$ is discrete and
central in $\mathcal N$.

\begin{lemma}[Convergence of Lattice Elements with Bounded Coordinates]
\label{m:fr:ext:global-normalized-words}\label{fl:N:limit}
If $\lambda_i\in\Lambda_i$ and the coordinates
\mbox{$\mathbf{c}_i(\lambda_i)$} converge to a point $x\in\mathbb R^d$, then $\vartheta_i(\lambda_i)\eGHto\pi(x)$ in the convergence of group elements of \cref{m:fr:lem:groupconv}. The convergence is uniform on each fixed compact set of normalized coordinates and each fixed bounded ball of $Z_i$.
\end{lemma}
\begin{proof}
Write $\lambda_i=u_1^{m_{1,i}}\cdots u_d^{m_{d,i}}$.  Bounded normalized
coordinates give $|m_{j,i}|\le CL_{j,i}$ for a fixed $C$.  Choose an
integer $N_*$, independent of $i$, so large that the segments $[-C/N_*,C/N_*]\,e_j$ lie in the interior of the rectangle $\mathsf X$ of \cref{m:fr:ext:physical-polynomial-coordinates}.  Split each $m_{j,i}$ into $N_*$
integer summands whose differences are at most one.  Each factor
$u_j^{m_{j,i}^{(h)}}$ is in the grid and has normalized coordinate tending to $(x_j/N_*)e_j$. Its marking therefore converges to $\exp(\delta_{\rm c}(x_j/N_*)v_j)$.  The product of these factors in $\Gamma_i$ is $\vartheta_i(\lambda_i)$.  Repeated use of
\cref{m:fr:lem:groupconv}(iv) gives its limit
\[
 \prod_{j=1}^d\exp(\delta_{\rm c} x_jv_j)=\pi(x).
\]
All factors have bounded displacement at the base point, so the finite product converges in the same sense. If uniformity on a compact set of coordinates and a bounded ball failed, choose
violating lattice elements, pass to a convergent coordinate subsequence,
and apply the assertion just proved, which includes uniformity on each bounded ball by the definition of the convergence of group elements.  This is a
contradiction.
\end{proof}

Conversely, convergent deck transformations lift. Near the identity this is
\cref{m:fr:ext:physical-polynomial-coordinates}\textup{(ii)}; a translation by a lattice element
with convergent coordinates extends it to all of $G_0$.

\begin{lemma}[Lifts of Convergent Deck Transformations]
\label{m:fr:npcq:lifts}
If elements $\gamma_i\in\Gamma_i$ converge to $g\in G_0$, then, for every $x\in\mathcal N$ with $\pi(x)=g$, there are
$\lambda_i\in\Lambda_i$ with
\[
 \vartheta_i(\lambda_i)=\gamma_i,
 \qquad \mathbf{c}_i(\lambda_i)\longrightarrow x.
\]
The conclusion holds for all sufficiently large $i$ and along arbitrary
subsequences. In particular, $\gamma_i$ eventually belongs to $H_i$.
\end{lemma}
\begin{proof}
Let $\lambda_i^0\in\Lambda_i$ be the element whose exponent vector is obtained by rounding
$L_ix$, so that $\mathbf{c}_i(\lambda_i^0)\to x$. \Cref{m:fr:ext:global-normalized-words} gives $\vartheta_i(\lambda_i^0)\eGHto g$. Thus the elements
\[
 \gamma'_i:=\vartheta_i(\lambda_i^0)^{-1}\gamma_i
\]
converge to the identity. \Cref{m:fr:ext:physical-polynomial-coordinates}\textup{(ii)}, applied to them and their limit $1$, supplies lifts
$\lambda_i'$ with $\vartheta_i(\lambda_i')=\gamma'_i$ and
$\mathbf{c}_i(\lambda_i')\to0$. Set
$\lambda_i=\lambda_i^0\lambda_i'$. Its marking is $\gamma_i$, and the convergence of the multiplication laws gives
$\mathbf{c}_i(\lambda_i)\to Q(x,0)=x$. The same inputs hold along
subsequences, so the assertion does too.
\end{proof}

The periods are the lattice elements with trivial marking. Applied to them, the two lemmas
identify the limits of their normalized coordinates with the kernel of the covering
homomorphism, in the following sense.

\begin{theorem}[Convergence of the Period Kernels]
\label{m:fr:ext:period-convergence}\label{fl:N:periods}
The groups $\mathsf K_i$ converge to $\mathsf K$ in the following precise senses.
\begin{enumerate}[label=\textup{(\roman*)}]
\item If $k_i\in \mathsf K_i$ and $k_i\to x$ in $\R^d$, then $x\in \mathsf K$.
\item For every $x\in \mathsf K$ there are $k_i\in \mathsf K_i$ with $k_i\to x$.
\item There is a fixed coordinate neighbourhood $\sV$ of $0$ with $\mathsf K_i\cap\sV=\{0\}$ for all large $i$.  Every $x\in \mathsf K$ has a
neighbourhood containing exactly one point of $\mathsf K_i$ for all large $i$;
this point converges to $x$.
\end{enumerate}
\end{theorem}
\begin{proof}
For (i), write $k_i=\mathbf{c}_i(\lambda_i)$ with
$\vartheta_i(\lambda_i)=1$.  \Cref{m:fr:ext:global-normalized-words} gives
$\pi(x)=1$, so $x\in \mathsf K$.

For (ii), let $x\in\mathsf K$, so that $\pi(x)=1$. \Cref{m:fr:npcq:lifts}, applied to the
trivial elements $\gamma_i=1$, gives $\lambda_i\in\ker\vartheta_i$ with
$\mathbf{c}_i(\lambda_i)\to x$; put $k_i:=\mathbf{c}_i(\lambda_i)$.

For (iii), choose a fixed small coordinate rectangle $\sV$ contained in the normalized coordinate domain of the grid for all large $i$. A lattice point of $\mathsf K_i$ in this rectangle is trivial because the marking is injective on the grid.  Fix
$x\in \mathsf K$ and choose a relatively compact neighbourhood $U_x$ so small
that the following quotient set lies in a compact subset of $\sV$:
\[
 \{Q(y,J(z)):y,z\in\overline U_x\}\subset\sV.
\]
Uniform polynomial and inverse convergence make the analogous inclusion
valid for $Q_i,J_i$.  Two points of
$\mathsf K_i\cap U_x$ would have quotient in $\mathsf K_i\cap\sV$ and must therefore
coincide.  Recovery in (ii) gives one point for all large $i$.  Choose
$U_x$ to contain no other point of the discrete set $\mathsf K$; (i), compactness,
and uniqueness then make that point converge to $x$.
\end{proof}

\subsubsection{The Nilmanifold of the Lattice}\label{fl:sec:nilmanifold}

The normalized coordinate space carries the group law $Q_i$, and the lattice acts on it by
left translation. The quotient is the model nilmanifold of the fibres, and the next lemma
records its properties. In the group $\cN_i=(\R^d,Q_i)$ of \cref{m:fr:ext:period-kernels} we
write $xy=Q_i(x,y)$ and $x^{-1}=J_i(x)$. The lattice acts by $\lambda x=\mathbf{c}_i(\lambda)x$,
and $\lambda\mathsf C:=\mathbf{c}_i(\lambda)\mathsf C$ for sets $\mathsf C\subset\R^d$.

\begin{lemma}[The Nilmanifold]\label{fl:lem:nilmanifold}
\begin{enumerate}[label=\textup{(\roman*)}]
\item Every $x\in\cN_i$ has the form $x=\mathbf{c}_i(\lambda)b$ for
$\lambda\in\Lambda_i$ and $b\in\prod_j[0,1/L_{j,i})$.
\item The action of $\Lambda_i$ on $\cN_i$ is free, properly discontinuous and
cocompact. The quotient $\operatorname{Nil}_i=\Lambda_i\backslash\cN_i$ is a compact
smooth nilmanifold of dimension $d$, with universal covering $\cN_i\cong\R^d$ and
fundamental group $\Lambda_i$.
\end{enumerate}
\end{lemma}
\begin{proof}
(i) The triangular law in \cref{fl:N:coord} determines the integers $m_j$ and the
numbers $b_j\in[0,1/L_{j,i})$ successively from
$x=\mathbf{c}_i(u_1^{m_1}\cdots u_d^{m_d})b$: the $j$-th coordinate is
$m_j/L_{j,i}+b_j$ plus a quantity determined by earlier coordinates.

(ii) Translations act freely, and properly discontinuously: if a translate $\lambda\mathsf C$ of a compact set $\mathsf C$ meets $\mathsf C$, then $\mathbf{c}_i(\lambda)\in\mathsf C\mathsf C^{-1}$, which contains only finitely many points of the discrete set $\mathbf{c}_i(\Lambda_i)=\prod_jL_{j,i}^{-1}\Z$. Part~(i) gives a compact fundamental set, so $\Lambda_i$ is a lattice in $\cN_i$, and the quotient is a compact nilmanifold with the universal covering $\cN_i\cong\R^d$ and the fundamental group $\Lambda_i$; this is the theorem of Malcev, see \cite[Chapter~II]{Rag72} and \cite[Section~5.1]{CG90}.
\end{proof}

\subsection{Boxes and Small Deck Transformations}\label{fl:sec:boxes}\label{fl:sec:malcev}\label{m:m:sec:relations}\label{m:fr:sec:malcevalgebras}

Boxes in the frame trivializations compare small deck transformations with lattice elements
of small coordinates, and frames with neighbourhoods in $\cY$. These comparisons are used at
all centres together, so we first pass to one subsequence for all charts.

Apply \cref{m:book:ch:12} successively at the finitely many centres to obtain one common subsequence, relabelled by $i$. Restrict the data of \cref{fl:sec:setting}, the realizing maps and the compact realization $\cY$ to these indices. The assertions \ref{fl:S:stage}--\ref{fl:S:Yconv}, equivariance, local product identities and estimates persist. Harmonic-chart subsequences for \cref{th:lem:stagecharts} are taken within these indices; the Haar limits in \cref{es:lem:limit,es:lem:limitcocycle} require no further extraction, and \cref{es:thm:gluingall} applies to every infinite subset. All remaining lattice and fibre assertions use this subsequence.

At each centre fix the lattice, marking and normalized laws above and the period-gap neighbourhood $\sV$ of \cref{fl:N:periods}\textup{(iii)}. Suppress the centre superscript. Since $\mathbf c_i(u_j)=e_j/L_{j,i}$ and the inverse laws converge, $\mathbf c_i(u_j^{\pm1})\to0$, and $\vartheta_i(u_j)\eGHto1$ by \cref{fl:N:limit}.

\subsubsection{Boxes}\label{fl:sec:boxdef}

Boxes are the images of products of balls under the frame trivializations. For $0<\eta\le1$ let
\begin{equation*}
O_\eta:=\{A\in O(n):|A-\one|_{\mathrm{HS}}<\eta\}.
\end{equation*}
This is an open neighbourhood of $\one$. It is contractible, so the sets $N_p(U,\eta)$ of
\cref{fl:def:boxes} are contractible whenever $U$ is; \cref{fl:sec:higher} uses this. Indeed,
a matrix $A\in O_\eta$ has no eigenvalue $-1$, so
it has determinant $1$ and a principal logarithm $\log A$, a skew-symmetric matrix that
depends continuously on $A$ and whose eigenvalues $\mathrm i\omega_k$ satisfy
$|\omega_k|<\pi$. Since
\begin{equation*}
|\exp(t\log A)-\one|_{\mathrm{HS}}^2=\sum_k4\sin^2\bigl(\tfrac12t\omega_k\bigr)
\end{equation*}
is non-decreasing in $t\in[0,1]$, the map $(A,t)\mapsto\exp((1-t)\log A)$ contracts $O_\eta$ to $\{\one\}$ within
$O_\eta$.

\begin{definition}[Boxes]\label{fl:def:boxes}
A \emph{level} is a pair $(s,\eta)$ with $0<s\le\tfrac12r$ and $0<\eta\le1$; we write
$(s,\eta)<(s',\eta')$ if $s<s'$ and $\eta<\eta'$, and $(s,\eta)\le(s',\eta')$ if $s\le s'$
and $\eta\le\eta'$. Let $p=(q,e)$ be a frame over $\Omega^\flat_i$ or over $\Omega^\flat$.
For an open set $U\subset B_{r/2}(q)$ and $0<\eta\le1$ put
\begin{equation*}
N_p(U,\eta):=\ft_p(U\times O_\eta),\qquad N_p(s,\eta):=N_p(B_s(q),\eta).
\end{equation*}
The set $N_p(s,\eta)$ is the \emph{box} of level $(s,\eta)$ about $p$.
\end{definition}

A box is an open neighbourhood of $p$ in $\Fr_i$, respectively in $\Fr$. It is
path-connected, because balls of the geodesic spaces $Z_i$ and $Z$ are path-connected. The next lemma relates
the boxes about convergent frames.

\begin{lemma}[Boxes and Convergence]\label{fl:lem:boxconv}
Let frames $p_i$ over $\Omega^\flat_i$ converge to a frame $p$ over $\Omega^\flat$, and let
$(s,\eta)<(s',\eta')$ be levels.
\begin{enumerate}[label=\textup{(\roman*)}]
\item If frames $x_i$ over $\Omega_i$ converge to a frame $x\in N_p(s,\eta)$, then
$x_i\in N_{p_i}(s,\eta)$ for large $i$.
\item Every sequence of frames $x_i\in N_{p_i}(s,\eta)$ has a subsequence that converges to
a frame in $N_p(s',\eta')$.
\item If levels $(s_k,\eta_k)$ tend to $(0,0)$, indices $i_k$ tend to infinity and
$x_k\in N_{p_{i_k}}(s_k,\eta_k)$, then $x_k$ converges to $p$.
\end{enumerate}
\end{lemma}
\begin{proof}
Write $p_i=(q_i,e_i)$ and $p=(q,e)$. (i) Since $\pr(x)\in B_s(q)$,
\ref{fl:S:product} gives $\pr(x_i)\in B_{r/2}(q_i)$ for large $i$ and
$\ft_{p_i}^{-1}(x_i)=(b_i,A_i)\to(b,A)=\ft_p^{-1}(x)$. Here $\hat d(b,q)<s$ and
$|A-\one|_{\mathrm{HS}}<\eta$, and distances converge, so $\hat d(b_i,q_i)<s$ and $|A_i-\one|_{\mathrm{HS}}<\eta$ for
large $i$.

(ii) Write $x_i=\ft_{p_i}(b_i,A_i)$. The points $b_i$ are bounded, so along a subsequence
$b_i\to b\in Z$ and $A_i\to A$, with $\hat d(b,q)\le s<s'$ and $|A-\one|_{\mathrm{HS}}\le\eta<\eta'$. By \ref{fl:S:product} the frames $x_i$ converge to $\ft_p(b,A)\in N_p(s',\eta')$.

(iii) Here the coordinates $(b_k,A_k)$ of $x_k$ satisfy $\hat d(b_k,q_{i_k})<s_k$ and
$|A_k-\one|_{\mathrm{HS}}<\eta_k$, so they converge to $(q,\one)$, and
$x_k\to\ft_p(q,\one)=p$ by \ref{fl:S:product}.
\end{proof}

\subsubsection{Small Deck Transformations and Small Lattice Elements}\label{fl:sec:capture}

The next two lemmas compare the deck transformations that move a box only slightly with
the lattice elements that have small coordinates. The first one is where the freeness of
the limit action on frames is used.

\begin{lemma}[Lifting Small Deck Transformations]\label{fl:lem:capture}
Let frames $p_i$ over $\Omega^\flat_i$ converge to a frame $p$ over $\Omega^\flat$, and let
$\sB$ be a neighbourhood of $0$ in $\R^d$. There are a level $(s,\eta)$ and $i_0$ such that
for all $i\ge i_0$ and $\gamma\in\Gamma_i$ the following holds: if
\begin{equation*}
\gamma\,N_{p_i}(s,\eta)\cap N_{p_i}(s,\eta)\ne\emptyset,
\end{equation*}
then $\gamma=\vartheta_i(\lambda)$ for some $\lambda\in\Lambda_i$ with $\mathbf{c}_i(\lambda)\in\sB$.
\end{lemma}
\begin{proof}
Suppose that this fails. Then there are levels $(s_k,\eta_k)\to(0,0)$, indices
$i_k\to\infty$, elements $\gamma_k\in\Gamma_{i_k}$ that are not of the form
$\vartheta_{i_k}(\lambda)$ with $\mathbf{c}_{i_k}(\lambda)\in\sB$, and frames $x_k$ and
$x'_k=\gamma_kx_k$ in $N_{p_{i_k}}(s_k,\eta_k)$. By \cref{fl:lem:boxconv}\textup{(iii)} both
$x_k$ and $x'_k$ converge to $p$. The displacements of $\gamma_k$ at the base points are
bounded, because $\gamma_k$ maps the point $\pr(x_k)$, which stays at bounded distance from
the base point, to $\pr(x'_k)$. By \cref{fl:S:groups} a subsequence of $\gamma_k$ converges to
some $g\in G$, and by \cref{fl:S:frameconv} the frames $\gamma_kx_k$ converge to $gp$ along
it. Limits are unique, so $gp=p$, and $g$ is the identity because $G$ acts freely on the frames by \cref{fl:S:frames}. Thus $\gamma_k\eGHto1\in G_0$ along the subsequence, and
\cref{m:fr:npcq:lifts} with $x=0$ gives elements $\lambda_k\in\Lambda_{i_k}$ with
$\vartheta_{i_k}(\lambda_k)=\gamma_k$ and $\mathbf{c}_{i_k}(\lambda_k)\to0$. For large $k$ this
contradicts the choice of $\gamma_k$.
\end{proof}

Conversely, lattice elements with small coordinates map a box into any larger box.

\begin{lemma}[Small Lattice Elements]\label{fl:lem:small}
Let frames $p_i$ over $\Omega^\flat_i$ converge to a frame $p$ over $\Omega^\flat$, and let
$(s,\eta)<(s',\eta')$ be levels. There are a neighbourhood $\sB$ of $0$ in $\R^d$ and $i_0$
such that
\begin{equation*}
\vartheta_i(\lambda)\,N_{p_i}(s,\eta)\subset N_{p_i}(s',\eta')
\end{equation*}
for all $i\ge i_0$ and all $\lambda\in\Lambda_i$ with $\mathbf{c}_i(\lambda)\in\sB$.
\end{lemma}
\begin{proof}
Suppose that this fails, and choose a level $(s'',\eta'')$ between the two given ones.
Then there are indices $i_k\to\infty$, elements $\lambda_k\in\Lambda_{i_k}$ with
$\mathbf{c}_{i_k}(\lambda_k)\to0$ and frames $x_k\in N_{p_{i_k}}(s,\eta)$ such that
$\vartheta_{i_k}(\lambda_k)x_k$ does not lie in $N_{p_{i_k}}(s',\eta')$. By
\cref{fl:lem:boxconv}\textup{(ii)} a subsequence of the frames $x_k$ converges to a frame
$x\in N_p(s'',\eta'')$. By \cref{fl:N:limit} the elements $\vartheta_{i_k}(\lambda_k)$ converge to
$\pi(0)=1$, so $\vartheta_{i_k}(\lambda_k)x_k\to x$ by \cref{fl:S:frameconv}, and
\cref{fl:lem:boxconv}\textup{(i)} shows that these frames lie in $N_{p_{i_k}}(s'',\eta'')$ for
large $k$. This is a contradiction.
\end{proof}

\subsubsection{Boxes and Neighbourhoods}\label{fl:sec:nbhd}

The images of the boxes in the spaces $Y_i$ form, for large $i$, a base of neighbourhoods
of a point of $Y$ in the realization $\cY$, in the following sense. This is how statements
about boxes are transferred to the neighbourhoods of \cref{fl:sec:results}.

\begin{lemma}[Boxes and Neighbourhoods]\label{fl:lem:nbhd}
Let $y\in Y$ with $\pr(y)\in B_{R_0/4}(a)$, let $p\in\tau^{-1}(y)$, and let frames $p_i$
over $\Omega^\flat_i$ converge to $p$.
\begin{enumerate}[label=\textup{(\roman*)}]
\item For every neighbourhood $\cU$ of $y$ in $\cY$ there are a level $(s,\eta)$ and $i_0$
with $\tau_i(N_{p_i}(s,\eta))\subset\cU$ for all $i\ge i_0$.
\item For every level $(s,\eta)$ there are an open neighbourhood $\cU$ of $y$ in $\cY$ and
an index $i_0$ such that $\cU\cap Y_i$ is contained in $\tau_i(N_{p_i}(s,\eta))$ for all
$i\ge i_0$.
\end{enumerate}
\end{lemma}
\begin{proof}
(i) Otherwise there are levels $(s_k,\eta_k)\to(0,0)$, indices $i_k\to\infty$ and frames
$x_k$ in the boxes $N_{p_{i_k}}(s_k,\eta_k)$ with $\tau_{i_k}(x_k)\notin\cU$. By
\cref{fl:lem:boxconv}\textup{(iii)}, $x_k\to p$, so $\tau_{i_k}(x_k)\to y$ in the sense of
\cref{fl:S:Yconv} and hence in $\cY$ by \cref{fl:prop:realization}; this is a contradiction.

(ii) Let $\cU_k$ be the open ball of radius $1/k$ about $y$ in $(\cY,\sd)$. If the
assertion fails for all these balls, there are indices $i_k\to\infty$ and points
$w_k\in\cU_k\cap Y_{i_k}$ outside $\tau_{i_k}(N_{p_{i_k}}(s,\eta))$. Then $w_k\to y$ in
$\cY$, hence in the sense of \cref{fl:S:Yconv}, and by \cref{fl:S:Yconv} there are frames
$x_k\in\tau_{i_k}^{-1}(w_k)$ that converge to $p$. By \cref{fl:lem:boxconv}\textup{(i)} they
lie in $N_{p_{i_k}}(s,\eta)$ for large $k$, a contradiction.
\end{proof}

\subsection{The Unwrapped Covering}\label{fl:sec:covering}

We glue lattice-labelled copies of a small open set of frames. Copies whose labels differ by a period lie over the same frames, but the covering keeps them apart. The lattice acts by deck transformations, simply transitively on each fibre, and supplies the monodromy used in \cref{fl:thm:fibre}.

\subsubsection{Faithful Coordinate Neighbourhoods}\label{fl:sec:lattice}\label{fl:sec:coordnbhd}

We fix three nested coordinate neighbourhoods of $0$ on which the marking is injective; they
define the sheets of the unwrapped coverings (\cref{fl:def:unwrapped}).

\begin{lemma}[Faithful Coordinate Neighbourhoods]\label{fl:lem:coordnbhd}\label{fl:lem:faithful}
There are bounded open neighbourhoods $\sB_0,\sB_1,\sB_2$ of $0$ with
$\overline{\sB_0}\subset\sB_1$, $\overline{\sB_1}\subset\sB_2\subset\sV$, such that, for all large $i$ and also for the limit laws $(Q,J)$:
\begin{enumerate}[label=\textup{(\alph*)}]
\item $Q_i(J_i(x),y)\in\sV$ for $x,y\in\sB_2$;
\item $Q_i(x,y),J_i(x)\in\sB_2$ for $x,y\in\sB_1$;
\item $\vartheta_i$ is injective on $\mathbf c_i^{-1}(\sB_2)$, $\mathsf K\cap\sV=\{0\}$, and $\pi$ is injective on $\sB_2$.
\end{enumerate}
Moreover $\mathbf c_i(u_j^{\pm1})\in\sB_0$ for all $j$ and all large $i$.
\end{lemma}
\begin{proof}
Choose $\sB_2$ and then $\sB_1$ so that the limit inclusions in (a),(b) hold with the closures of the domains and with images in the stated open targets, and take any $\sB_0$ with closure in $\sB_1$. Compactness and uniform convergence of $Q_i,J_i$ give the same inclusions for large $i$. The last assertion follows from $\mathbf c_i(u_j^{\pm1})\to0$.

If two lattice elements with coordinates in $\sB_2$ have the same marking, their quotient has coordinates in $\sV$ by (a) and lies in the period kernel, hence is trivial by \cref{fl:N:periods}\textup{(iii)}. If $x\in\mathsf K\cap\sV$, recovery in \cref{fl:N:periods}\textup{(ii)} gives points of $\mathsf K_i$ converging to $x$; they lie in $\sV$ eventually and are zero, so $x=0$. Applying (a) to two points with the same $\pi$-image proves the remaining injectivity.
\end{proof}

Fix these neighbourhoods at every centre and restrict henceforth to the indices for which the lemma holds.

\subsubsection{Unwrapping Sets and the Construction}\label{fl:sec:construction}

The following small-overlap condition makes lattice labels consistent on intersecting sheets.

\begin{definition}[Unwrapping Sets]\label{fl:def:small}
An open subset $V$ of $\Fr_i$ is an \emph{unwrapping set} if every $\gamma\in\Gamma_i$ with
$\gamma V\cap V\ne\emptyset$ is of the form $\gamma=\vartheta_i(\lambda)$ for some
$\lambda\in\Lambda_i$ with $\mathbf{c}_i(\lambda)\in\sB_1$.
\end{definition}

By \cref{fl:lem:faithful}\textup{(c)} the element $\lambda$ is the only lift of $\gamma$ with coordinates
in $\sB_2$. Open subsets of unwrapping sets are unwrapping sets. By
\cref{fl:lem:capture}, boxes of a sufficiently small level about convergent frames are
unwrapping sets for large $i$. The relation in the following definition is an equivalence
relation by \cref{fl:prop:covering}\textup{(i)}.

\begin{definition}[Unwrapped Coverings]\label{fl:def:unwrapped}
For an unwrapping set $V$ let
\begin{equation*}
\widetilde{\cE}_i(V):=\bigl\{(\lambda,x)\in\Lambda_i\times\Fr_i:x\in\vartheta_i(\lambda)V\bigr\},
\end{equation*}
with the topology induced by the product of the discrete group $\Lambda_i$ and $\Fr_i$, and
let
\begin{equation*}
(\lambda,x)\sim(\lambda',x')\quad:\Longleftrightarrow\quad x=x'\ \text{ and }\
\mathbf{c}_i(\lambda^{-1}\lambda')\in\sB_2 .
\end{equation*}
The \emph{unwrapped covering} of $V$ is the set $\cE_i(V)$ of equivalence classes with the
quotient topology, with the class $[\lambda,x]$ of $(\lambda,x)$ and the map
\begin{equation*}
\pi_V\colon \cE_i(V)\longrightarrow\Fr_i,\qquad\pi_V[\lambda,x]:=x.
\end{equation*}
\end{definition}

The image of $\pi_V$ is the union $H_iV$ of the sets $hV$ with $h\in H_i$. The group $\Lambda_i$ acts on $\widetilde{\cE}_i(V)$ by
$\nu(\lambda,x):=(\nu\lambda,\vartheta_i(\nu)x)$, and this action preserves the relation.

\begin{proposition}[The Unwrapped Covering]\label{fl:prop:covering}
Let $V\subset\Fr_i$ be an unwrapping set.
\begin{enumerate}[label=\textup{(\roman*)},leftmargin=2.5em]
\item The relation $\sim$ is an equivalence relation, and two elements $(\lambda,x)$ and
$(\lambda',x)$ of $\widetilde{\cE}_i(V)$ are equivalent if and only if
$\mathbf{c}_i(\lambda^{-1}\lambda')\in\sB_1$.
\item For every $\lambda\in\Lambda_i$ the map
$\sigma_\lambda\colon\vartheta_i(\lambda)V\to \cE_i(V)$, $x\mapsto[\lambda,x]$, is a
homeomorphism onto an open subset $\Sh_V(\lambda)$, the \emph{sheet} of $\lambda$, with
inverse $\pi_V$. Two sheets are disjoint, or $\mathbf{c}_i(\lambda^{-1}\lambda')\in\sB_1$ for
their indices. The space $\cE_i(V)$ is a Hausdorff space.
\item The map $\pi_V$ is a covering map of $\cE_i(V)$ onto $H_iV$: for every $\lambda$, the
preimage of $\vartheta_i(\lambda)V$ is the disjoint union of the sheets $\Sh_V(\lambda k)$ with
$k\in\ker\vartheta_i$.
\item The formula $\nu[\lambda,x]:=[\nu\lambda,\vartheta_i(\nu)x]$ defines a free action of
$\Lambda_i$ on $\cE_i(V)$ by homeomorphisms, with $\nu\,\Sh_V(\lambda)=\Sh_V(\nu\lambda)$
and $\pi_V(\nu\hat x)=\vartheta_i(\nu)\pi_V(\hat x)$ for $\hat x\in \cE_i(V)$.
\item The map $\tau_i\circ\pi_V$ is a covering map of $\cE_i(V)$ onto the open subset
$\tau_i(V)$ of $Y_i$. Its fibres are the orbits of $\Lambda_i$, so $\Lambda_i$ acts by deck
transformations, simply transitively on every fibre; if $\cE_i(V)$ is connected, $\Lambda_i$
is the group of all deck transformations.
\item If $V'\subset V$ is open, then $[\lambda,x]\mapsto[\lambda,x]$ is a homeomorphism of
$\cE_i(V')$ onto the open subset $(\tau_i\circ\pi_V)^{-1}(\tau_i(V'))$ of $\cE_i(V)$. It
commutes with the actions of $\Lambda_i$ and maps $\Sh_{V'}(\lambda)$ into
$\Sh_V(\lambda)$.
\item If $V$ is path-connected and contains, for some frame $x\in V$, the frames
$\vartheta_i(u_j)x$ with $1\le j\le d$, then $\cE_i(V)$ is path-connected.
\end{enumerate}
\end{proposition}
\begin{proof}
(i) Let $(\lambda,x)$ and $(\lambda',x)$ lie in $\widetilde{\cE}_i(V)$. Then $\vartheta_i(\lambda)^{-1}x$
lies in $V$ and in $\vartheta_i(\lambda^{-1}\lambda')V$, so
$\vartheta_i(\lambda^{-1}\lambda')=\vartheta_i(\mu)$ for some $\mu$ with $\mathbf{c}_i(\mu)\in\sB_1$,
because $V$ is an unwrapping set. If $\mathbf{c}_i(\lambda^{-1}\lambda')\in\sB_2$, then
$\lambda^{-1}\lambda'=\mu$ by \cref{fl:lem:faithful}\textup{(c)}. This proves the second assertion. The
relation is reflexive. It is symmetric and transitive by the second assertion and
\cref{fl:lem:coordnbhd}\textup{(b)}: if $\mathbf{c}_i(\lambda^{-1}\lambda')$ and
$\mathbf{c}_i(\lambda'^{-1}\lambda'')$ lie in $\sB_1$, then
$\mathbf{c}_i(\lambda'^{-1}\lambda)=J_i(\mathbf{c}_i(\lambda^{-1}\lambda'))$ and
$\mathbf{c}_i(\lambda^{-1}\lambda'')=Q_i(\mathbf{c}_i(\lambda^{-1}\lambda'),\mathbf{c}_i(\lambda'^{-1}\lambda''))$
lie in $\sB_2$.

(ii) The map $\sigma_\lambda$ is continuous and injective. For an open set
$O\subset\vartheta_i(\lambda)V$ the preimage of $\sigma_\lambda(O)$ in $\widetilde{\cE}_i(V)$ is the set of
all $(\lambda',x)$ with $x\in O\cap\vartheta_i(\lambda')V$ and
$\mathbf{c}_i(\lambda^{-1}\lambda')\in\sB_2$, which is open. So $\sigma_\lambda$ is an open
map, and its inverse is the restriction of $\pi_V$; in particular $\pi_V$ is continuous.
If two sheets have a common point $[\lambda,x]=[\lambda',x]$, then
$\mathbf{c}_i(\lambda^{-1}\lambda')\in\sB_1$ by (i). Two points of $\cE_i(V)$ with different
images under $\pi_V$ are separated by preimages of disjoint open subsets of $\Fr_i$. Two
different points $[\lambda,x]$ and $[\lambda',x]$ lie in the sheets of $\lambda$ and
$\lambda'$, and these are disjoint: a common point $[\lambda,x'']=[\lambda',x'']$ would
give $\mathbf{c}_i(\lambda^{-1}\lambda')\in\sB_2$ and hence $[\lambda,x]=[\lambda',x]$.

(iii) Let $[\lambda',x]\in \cE_i(V)$ with $x\in\vartheta_i(\lambda)V$. As in (i) there is $\mu$
with $\mathbf{c}_i(\mu)\in\sB_1$ and $\vartheta_i(\lambda^{-1}\lambda')=\vartheta_i(\mu)$. Then
$k:=\lambda^{-1}\lambda'\mu^{-1}$ lies in $\ker\vartheta_i$, the pair $(\lambda k,x)$ lies in
$\widetilde{\cE}_i(V)$ because $\vartheta_i(\lambda k)=\vartheta_i(\lambda)$, and $(\lambda k)^{-1}\lambda'=\mu$;
so $[\lambda',x]=[\lambda k,x]$ lies in $\Sh_V(\lambda k)$. For $k\ne k'$ in $\ker\vartheta_i$
the sheets of $\lambda k$ and $\lambda k'$ are disjoint, since otherwise $k^{-1}k'$ would
be an element of $\ker\vartheta_i$ with coordinates in $\sB_2\subset\sV$, which is trivial by
\cref{fl:N:periods}\textup{(iii)}. Each of these sheets is mapped homeomorphically onto
$\vartheta_i(\lambda k)V=\vartheta_i(\lambda)V$ by $\pi_V$, and the open sets $\vartheta_i(\lambda)V$
cover $H_iV$.

(iv) The map is well defined and continuous, with inverse given by $\nu^{-1}$, and it
maps $\sigma_\lambda(x)$ to $\sigma_{\nu\lambda}(\vartheta_i(\nu)x)$. If $\nu[\lambda,x]=[\lambda,x]$,
then $\vartheta_i(\nu)x=x$, so $\vartheta_i(\nu)=1$ because $\Gamma_i$ acts freely on $\Fr_i$, and
$\lambda^{-1}\nu\lambda$ is an element of $\ker\vartheta_i$ with coordinates in $\sB_2$; hence
$\nu=1$ by \cref{fl:N:periods}\textup{(iii)}.

(v) The set $H_iV$ is open and invariant under $H_i$. If $\gamma\in\Gamma_i$ and
$\gamma H_iV$ meets $H_iV$, then $\gamma hV$ meets $h'V$ for some $h,h'\in H_i$, so
$h'^{-1}\gamma h$ lies in $H_i$ because $V$ is an unwrapping set, and $\gamma\in H_i$. Let
$x_0\in V$. Since $\Gamma_i$ acts freely and properly discontinuously, there is an open
neighbourhood $O\subset V$ of $x_0$ with $\gamma O\cap O=\emptyset$ for all $\gamma\ne1$.
By \cref{fl:S:charts}, $\tau_i$ maps $O$ homeomorphically onto the open set $W:=\tau_i(O)$,
and by what we have shown the part of $\tau_i^{-1}(W)$ in $H_iV$ is the disjoint union of
the sets $hO$ with $h\in H_i$. By (iii), applied to the subsets $\vartheta_i(\lambda)O$ of
$\vartheta_i(\lambda)V$, the preimage of $hO$ under $\pi_V$, for $h=\vartheta_i(\lambda)$, is the
disjoint union of the sets $\sigma_{\lambda k}(\vartheta_i(\lambda)O)$ with $k\in\ker\vartheta_i$.
When $\lambda$ runs through representatives of the cosets of $\ker\vartheta_i$ and $k$ through
$\ker\vartheta_i$, the product $\lambda k$ runs once through $\Lambda_i$. Hence
\begin{equation*}
(\tau_i\circ\pi_V)^{-1}(W)=\bigsqcup_{\lambda\in\Lambda_i}\sigma_\lambda(\vartheta_i(\lambda)O),
\end{equation*}
each of these open sets is mapped homeomorphically onto $W$, and $\nu\in\Lambda_i$ maps the
set with index $\lambda$ onto the set with index $\nu\lambda$. The sets $W$ cover
$\tau_i(V)$, which is open because $\tau_i$ is an open map.

(vi) The set $\widetilde{\cE}_i(V')$ is a subset of $\widetilde{\cE}_i(V)$ on which the two relations agree, so
the map is well defined, injective and continuous, and it commutes with the actions. It
maps $\sigma^{V'}_\lambda(O)$ to $\sigma^V_\lambda(O)$ for open sets
$O\subset\vartheta_i(\lambda)V'$, so it is open. Let $[\lambda,x]\in \cE_i(V)$ with
$\tau_i(x)\in\tau_i(V')$, that is, $x\in\gamma V'$ for some $\gamma\in\Gamma_i$. Then
$\vartheta_i(\lambda)^{-1}\gamma V$ meets $V$, so $\vartheta_i(\lambda)^{-1}\gamma=\vartheta_i(\mu)$ with
$\mathbf{c}_i(\mu)\in\sB_1$, and $x\in\vartheta_i(\lambda\mu)V'$. Hence $(\lambda\mu,x)$ lies in
$\widetilde{\cE}_i(V')$, and it is equivalent to $(\lambda,x)$ in $\widetilde{\cE}_i(V)$. This proves that the
image is the preimage of $\tau_i(V')$.

(vii) Every sheet is homeomorphic to $V$ and hence path-connected. For every $j$ the
frame $\vartheta_i(u_j)x$ lies in $V$ and in $\vartheta_i(u_j)V$, and $\mathbf{c}_i(u_j)\in\sB_2$ by
\cref{fl:lem:coordnbhd}; so $[1,\vartheta_i(u_j)x]=[u_j,\vartheta_i(u_j)x]$ is a common point of the
sheets of $1$ and of $u_j$. Applying $\nu\in\Lambda_i$, we see that the sheets of $\nu$
and of $\nu u_j$ meet. Since the elements $u_j$ generate $\Lambda_i$, every sheet is
joined to the sheet of $1$ by a chain of sheets in which consecutive ones meet.
\end{proof}

Lifting loops to the unwrapped covering assigns to each loop an element of the fibre
lattice $\Lambda_i$.

\begin{proposition}[Monodromy]\label{fl:prop:monodromy}\label{fl:sec:monodromy}
Let $V$ be an unwrapping set and $x_0\in V$. For a loop $c$ in $\tau_i(V)$ at $\tau_i(x_0)$, let
$\hat c$ be its lift to $\cE_i(V)$ under $\tau_i\circ\pi_V$ that starts at $\hat x_0:=[1,x_0]$, and define
$\mathrm{mon}_V^{x_0}([c])\in\Lambda_i$ by $\hat c(1)=\mathrm{mon}_V^{x_0}([c])\,\hat x_0$.
\begin{enumerate}[label=\textup{(\roman*)}]
\item $\mathrm{mon}_V^{x_0}\colon\pi_1(\tau_i(V),\tau_i(x_0))\to\Lambda_i$ is a homomorphism, with path products traversed from left to right.
\item For open $V'\subset V$ containing $x_0$, $\mathrm{mon}_{V'}^{x_0}$ is the composition of $\mathrm{mon}_V^{x_0}$
with the map induced by the inclusion of $\tau_i(V')$ in $\tau_i(V)$.
\item If $\cE_i(V)$ is path-connected, $\mathrm{mon}_V^{x_0}$ is onto.
\end{enumerate}
\end{proposition}
\begin{proof}
The covering and restriction properties of \cref{fl:prop:covering}\textup{(v),(vi)} give (i),(ii) by path lifting \cite[Section~1.3]{Hat02}: the lift of a product is the first lift followed by the appropriate deck translate of the second. For (iii), project a path from $\hat x_0$ to $\nu\hat x_0$ for each $\nu\in\Lambda_i$.
\end{proof}

\section{Fibre Topology and Isotropy Descent}\label{ch:fibregroups}\label{ch:fibres}\label{es:sec:fibres}

The unwrapped coverings now identify the fibres. Controlled homotopies (\cref{fl:sec:higher}) determine the intrinsic fundamental groups and the higher homotopy of small pieces (\cref{fl:sec:results}). For the bundles $F^\natural_i$ of \cref{es:thm:gluingall} this shows that the fibres are nilmanifolds (\cref{es:sec:fibreframes}); the fibres of the induced maps $f^\natural_i\colon X_i\to X$ are their quotients by the isotropy groups, hence infranilmanifolds (\cref{es:sec:fibrestages}). \Cref{es:sec:proofsintro} proves \cref{thm:main}.

Throughout we retain the setting and common subsequence of \cref{ch:lattices}. For the local homotopy arguments fix one centre, suppress its superscript, and choose reference frames $p_i\to p$ over $\Omega_i^\flat$ and $\Omega^\flat$.

\subsection{Controlled Homotopies}\label{fl:sec:higher}\label{fl:sec:capturelevel}

A lifted sphere may visit arbitrarily many sheets. We map its sheet labels into the contractible nilpotent group $\cN_i$, fill there, and transfer the disk back through a tower of contractible sheets while fixing its boundary. By \cref{fl:lem:capture}, choose a \emph{lifting level} $(s_c,\eta_c)$ such that, for large $i$,
\[
\gamma N_{p_i}(s_c,\eta_c)\cap N_{p_i}(s_c,\eta_c)\ne\emptyset
\quad\Longrightarrow\quad
\gamma=\vartheta_i(\lambda),\qquad\mathbf c_i(\lambda)\in\sB_0.
\]
Every open subset of this box is therefore an unwrapping set. The group law of $\cN_i$ and the action of $\Lambda_i$ on it are written as in \cref{fl:sec:nilmanifold}.

\subsubsection{Complexes of Small Lattice Elements}\label{fl:sec:latticecomplex}

The combinatorial counterpart of the cover of an unwrapped covering by its sheets is the
following complex. For $\mathbf b\in(0,\infty)^d$ write
\[
\Bx(\mathbf b):=\prod_{j=1}^d(-\mathbf b_j,\mathbf b_j),\qquad
A_i(\mathbf b):=\mathbf{c}_i^{-1}(\Bx(\mathbf b)).
\]
The finite set $A_i(\mathbf b)$ need not be closed under inversion. Let $\cP_i(\mathbf b)$ be the simplicial complex
whose vertices are the elements of $\Lambda_i$ and whose simplices are the finite
non-empty sets $\sigma\subset\Lambda_i$ with
\begin{equation*}
\lambda^{-1}\lambda'\in A_i(\mathbf b)\qquad\text{for all }\lambda\ne\lambda'\text{ in }\sigma .
\end{equation*}
The group $\Lambda_i$ acts on it by left multiplication. The action is free on the set of
simplices, because an element that maps a simplex to itself has a power that fixes a
vertex, and $\Lambda_i$ is torsion-free. The complex is locally finite, because
$A_i(\mathbf b)$ is finite. We write $\cP_i(\mathbf b)^{(q)}$ for its $q$-skeleton, and
$\mathbf b\le\mathbf b'$ for the inequality in every component; then $\cP_i(\mathbf b)$ is a subcomplex
of $\cP_i(\mathbf b')$. The next lemma maps the skeleta of these
complexes to $\cN_i$ with control. Since the translations of $\cN_i$ are not affine maps of
$\R^d$, it extends maps by cones to the identity only on simplices with the vertex $1$ and
transports these extensions by the group.

\begin{lemma}[Controlled Maps]\label{fl:lem:controlled}
Let $q\ge0$, and let $\mathsf C$ be a neighbourhood of $0$ in $\R^d$. There are
$\mathbf b_0\in(0,\infty)^d$ and $i_0$ such that for every $i\ge i_0$ there is an equivariant
map
\begin{equation*}
a\colon|\cP_i(\mathbf b_0)^{(q)}|\longrightarrow\cN_i
\end{equation*}
with $a(\lambda)=\mathbf{c}_i(\lambda)$ for every vertex $\lambda$ and
$a(\sigma)\subset\lambda\mathsf C$ for every simplex $\sigma$ of dimension at most $q$ and
every vertex $\lambda$ of $\sigma$.
\end{lemma}
\begin{proof}
The map
$(x,x',t,c)\mapsto Q(J(x),t\cdot Q(x',c))$ is continuous and maps $(0,0,t,0)$ to $0$ for
all $t\in[0,1]$. Hence there are compact neighbourhoods
$\mathsf C_0\subset\mathsf C_1\subset\dots\subset\mathsf C_q\subset\mathsf C$ of $0$ and
$\mathbf b_0$ such that
\begin{equation}\label{fl:eq:controlled}
J(x)\,\bigl(t\cdot(x'c)\bigr)\in\mathsf C_k\qquad
\text{for }x,x'\in\overline{\Bx(\mathbf b_0)},\ t\in[0,1],\ c\in\mathsf C_{k-1},\ 1\le k\le q,
\end{equation}
for the group law of $\cN$, with the images in the interior of $\mathsf C_k$; one chooses
$\mathsf C_q$ first and then descends. By the uniform convergence in \cref{fl:N:coord},
\eqref{fl:eq:controlled} holds for the group law of $\cN_i$ for all large $i$.

We define $a$ on the skeleta, with $a(\sigma)\subset\lambda\mathsf C_k$ for every
$k$-simplex $\sigma$ and every vertex $\lambda$ of $\sigma$. On the vertices
$a(\lambda):=\mathbf{c}_i(\lambda)$, which lies in $\lambda\mathsf C_0$. Let $1\le k\le q$, and
let $a$ be defined on the $(k-1)$-skeleton. Every orbit of $k$-simplices contains a
simplex $\sigma_0$ with the vertex $1$; choose one in every orbit. Its vertices $\mu$
satisfy $\mathbf{c}_i(\mu)\in\Bx(\mathbf b_0)$, and every face of $\sigma_0$ contains one of
them, so $a(\partial\sigma_0)$ lies in the union of the sets $\mu\mathsf C_{k-1}$. Extend
$a$ to $\sigma_0$ as the cone with the value $0$ at the barycentre:
$a((1-t)\,b_{\sigma_0}+t\,\xi):=t\cdot a(\xi)$ for $\xi\in\partial\sigma_0$. By
\eqref{fl:eq:controlled}, $\mathbf{c}_i(\mu)^{-1}a(\sigma_0)\subset\mathsf C_k$ for every vertex
$\mu$ of $\sigma_0$. Define $a$ on the simplices $\nu\sigma_0$ by the equivariance. This
is consistent with the values on the $(k-1)$-skeleton, and
$a(\nu\sigma_0)=\nu a(\sigma_0)\subset\nu\mu\mathsf C_k$ for the vertices $\nu\mu$ of
$\nu\sigma_0$.
\end{proof}

\subsubsection{Towers and Carriers}\label{fl:sec:towers}

A tower provides, for every lattice element, nested contractible sheets. A small change of the
element is absorbed by passing to the next larger sheet. The carrier lemma then fills a cell once its boundary lies in that sheet, preserving the maps already
defined on its faces. Equivariance makes this construction uniform over all lattice elements.
One tower level is used per increase in cell dimension, not per cell or sheet, so
the tower length is independent of how many sheets a lifted sphere visits.

\begin{definition}[Towers]\label{fl:def:tower}
Let $M\ge1$, and let $\mathbf b'\in(0,\infty)^d$ with $\Bx(\mathbf b')\subset\sB_1$. A \emph{tower}
of length $M$ for $\mathbf b'$ consists, for all large $i$, of open sets of frames
\begin{equation*}
V^{[0]}_i\subset V^{[1]}_i\subset\dots\subset V^{[M]}_i\subset N_{p_i}(s_c,\eta_c)
\end{equation*}
with the following properties.
\begin{enumerate}[label=\textup{(T\arabic*)},ref=\textup{(T\arabic*)},leftmargin=3.2em]
\item\label{fl:T:contractible} Every set $V^{[j]}_i$ is contractible, and $V^{[0]}_i$
contains a box $N_{p_i}(s,\eta)$ whose level $(s,\eta)$ does not depend on $i$.
\item\label{fl:T:move} $\vartheta_i(\mu)V^{[j]}_i\subset V^{[j+1]}_i$ for all $\mu\in A_i(\mathbf b')$
and $0\le j<M$.
\end{enumerate}
\end{definition}

For a tower we write $\Sh^{[j]}(\lambda)$ for the sheet of $\lambda$ in $\cE_i(V^{[j]}_i)$,
regarded as a subset of $\cE_i(V^{[M]}_i)$ by \cref{fl:prop:covering}\textup{(vi)}. These
sheets are contractible, $\nu\Sh^{[j]}(\lambda)=\Sh^{[j]}(\nu\lambda)$, and
\begin{equation}\label{fl:eq:towersheets}
\Sh^{[j]}(\lambda)\subset\Sh^{[j+1]}(\lambda),\qquad
\Sh^{[j]}(\lambda\mu)\subset\Sh^{[j+1]}(\lambda)\quad\text{for }\mu\in A_i(\mathbf b').
\end{equation}
Indeed, for a point $[\lambda\mu,x]$ of $\Sh^{[j]}(\lambda\mu)$ the frame $x$ lies in
$\vartheta_i(\lambda\mu)V^{[j]}_i\subset\vartheta_i(\lambda)V^{[j+1]}_i$ by \cref{fl:T:move}, and
$(\lambda,x)$ is equivalent to $(\lambda\mu,x)$ because $\mathbf{c}_i(\mu)\in\sB_2$. Towers
exist below every level:

\begin{lemma}[Existence of Towers]\label{fl:lem:towers}
For every $M\ge1$ and every level $(s^\circ,\eta^\circ)$ there are $\mathbf b'$ and a tower of
length $M$ for $\mathbf b'$ with $V^{[M]}_i\subset N_{p_i}(s^\circ,\eta^\circ)$.
\end{lemma}
\begin{proof}
Let $\ell_0=\ell_0(n)$ be the constant of \cref{fl:lem:reifenberg}. Choose levels
$(s^{[0]},\eta^{[0]}),\dots,(s^{[M]},\eta^{[M]})$ with $(s^{[M]},\eta^{[M]})\le(s^\circ,\eta^\circ)$,
$(s^{[M]},\eta^{[M]})\le(s_c,\eta_c)$ and $s^{[M]}\le\tfrac14r$, the radius bound of that lemma, such that
\begin{equation*}
(s^{[j]},\eta^{[j]})<(\ell_0s^{[j+1]},\eta^{[j+1]})\qquad(0\le j<M).
\end{equation*}
By \cref{fl:lem:reifenberg} there are open sets $U^{[j]}_i\subset Z_i$ homeomorphic to $\R^n$ with
\[
 B_{\ell_0s^{[j]}}(\pr(p_i))\subset U^{[j]}_i\subset B_{s^{[j]}}(\pr(p_i)).
\]
Put $V^{[j]}_i:=N_{p_i}(U^{[j]}_i,\eta^{[j]})$. These sets are contractible, because $U^{[j]}_i$ and
$O_{\eta^{[j]}}$ are (\cref{fl:sec:boxdef}), and
\begin{equation*}
N_{p_i}(\ell_0s^{[j]},\eta^{[j]})\subset V^{[j]}_i\subset N_{p_i}(s^{[j]},\eta^{[j]})\subset
N_{p_i}(\ell_0s^{[j+1]},\eta^{[j+1]})\subset V^{[j+1]}_i .
\end{equation*}
By \cref{fl:lem:small}, applied to the $M$ pairs of levels, there is a neighbourhood $\sB$ of
$0$ such that for large $i$ and $\mathbf{c}_i(\mu)\in\sB$ the transformation $\vartheta_i(\mu)$
maps $N_{p_i}(s^{[j]},\eta^{[j]})$ into $N_{p_i}(\ell_0s^{[j+1]},\eta^{[j+1]})$ for all $j<M$. Every
$\mathbf b'$ with $\Bx(\mathbf b')\subset\sB\cap\sB_1$ has the required property; the box in
\ref{fl:T:contractible} is $N_{p_i}(\ell_0s^{[0]},\eta^{[0]})$.
\end{proof}

The maps are constructed cell by cell on complexes whose cells carry lattice elements that
change little between a cell and its faces.

\begin{definition}[Labellings]\label{fl:def:labelling}
A \emph{labelling} for $\mathbf b'$ of a CW complex $\mathcal L$ with a cellular action of $\Lambda_i$ assigns to every
cell $c$ a lattice element $\lambda_c\in\Lambda_i$, its \emph{cell label}, with $\lambda_{\nu c}=\nu\lambda_c$, such that
\begin{equation}\label{fl:eq:labelling}
\lambda_c^{-1}\lambda_{c'}\in A_i(\mathbf b')\cup\{1\}
\end{equation}
whenever $c'$ is a cell of lower dimension that meets the closure of $c$.
\end{definition}

The following lemma is the device by which maps are constructed; it is the equivariant form
of the extension of maps along contractible carriers.

\begin{lemma}[Carriers]\label{fl:lem:carriers}
Fix a tower of length $M$ for $\mathbf b'$ and a large index $i$. Let $\mathcal L$ be a CW complex
with a cellular action of $\Lambda_i$ that is free on the set of cells, with a labelling
$c\mapsto\lambda_c$, and let $\mathcal L_0$ be an invariant subcomplex. Let $m_0\ge0$ be an integer with
$m_0+\dim c\le M$ for all cells $c$ of $\mathcal L$, and let
$\varphi_0\colon \mathcal L_0\to \cE_i(V^{[M]}_i)$ be an equivariant map with
\begin{equation}\label{fl:eq:carrier}
\varphi_0(\bar c)\subset\Sh^{[m_0+\dim c]}(\lambda_c)
\end{equation}
for every cell $c$ of $\mathcal L_0$. Then $\varphi_0$ extends to an equivariant map
$\varphi\colon \mathcal L\to \cE_i(V^{[M]}_i)$ that satisfies \eqref{fl:eq:carrier} for every cell $c$
of $\mathcal L$.
\end{lemma}
\begin{proof}
We extend the map over the cells that are not in $\mathcal L_0$, by induction on their dimension,
and we treat one cell $c$ in every orbit; the map on the other cells of the orbit is then
defined by the equivariance, which is consistent because the action is free on the cells,
and it satisfies \eqref{fl:eq:carrier} because $\nu\Sh^{[j]}(\lambda)=\Sh^{[j]}(\nu\lambda)$.
For a cell $c$ of dimension $0$ put $\varphi(c):=\lambda_c\hat p_i$, with the point
$\hat p_i:=[1,p_i]$ of $\Sh^{[0]}(1)$; here $p_i\in V^{[0]}_i$ by \ref{fl:T:contractible}. This
is a point of $\Sh^{[0]}(\lambda_c)$. Let $c$ be a cell of dimension $k\ge1$ with a characteristic map
$\chi\colon D^k\to \mathcal L$. The map $\varphi\circ\chi$ is defined on the boundary of $D^k$, and
its image lies in the union of the sets $\varphi(\bar c')$ over the finitely many cells
$c'$ of lower dimension that meet $\bar c$. By \eqref{fl:eq:carrier} for these cells, by
\eqref{fl:eq:labelling} and by \eqref{fl:eq:towersheets},
\begin{equation*}
\varphi(\bar c')\subset\Sh^{[m_0+k-1]}(\lambda_{c'})\subset\Sh^{[m_0+k]}(\lambda_c).
\end{equation*}
The sheet on the right is contractible, so $\varphi\circ\chi$ extends from the boundary
to a map of $D^k$ into it, and this defines $\varphi$ on $\bar c$. The resulting map is
continuous by the weak topology of the CW complex.
\end{proof}

\subsubsection{Liftable Spheres}\label{fl:sec:liftablespheres}
Maps of spheres into the image of a small unwrapping set that lift to its unwrapped covering
are null-homotopic in the image of a larger contractible unwrapping set. This is the filling
step behind both theorems of \cref{fl:sec:intro}.

\begin{lemma}[Liftable Spheres]\label{fl:lem:liftablespheres}
For every $q\ge1$ and level $(s^\circ,\eta^\circ)$, there are a fixed level $(s',\eta')$ and, for all large $i$, a contractible unwrapping set $V_i^{\rm out}$ with
\[
V'_i:=N_{p_i}(s',\eta')\subset V_i^{\rm out}\subset N_{p_i}(s^\circ,\eta^\circ)
\]
such that every map $f\colon\mathbb S^j\to\tau_i(V'_i)$, $0\le j\le q$, that lifts to $\cE_i(V'_i)$ is null-homotopic in $\tau_i(V_i^{\rm out})$.
Thus any two points of $\tau_i(V'_i)$ are joined there, every sphere of dimension $2\le j\le q$ contracts there, and so does every loop with trivial monodromy.
\end{lemma}
\begin{proof}
Choose a tower of length $M=q+2$ inside $N_{p_i}(s^\circ,\eta^\circ)$ by
\cref{fl:lem:towers}, and put $V_i^{\rm out}=V_i^{[M]}$. Choose relatively compact open
neighbourhoods $\sB,\mathsf D$ of $0$ so small that, for the limiting group law,
\[
 \overline{\sB}\,\overline{\sB}^{\,-1}\cup
 \overline{\sB}\,\overline{\mathsf D}^{\,-1}\subset\Bx(\mathbf b').
\]
These inclusions hold also in $\cN_i$ for large $i$, by \cref{fl:N:coord}.
Moreover $\prod_j[0,1/L_{j,i})\subset\sB$ for large $i$, so the translates
$\lambda\sB$, $\lambda\in\Lambda_i$, cover $\cN_i$ by \cref{fl:lem:nilmanifold}\textup{(i)}.
With $\mathsf C=\mathsf D$, \cref{fl:lem:controlled} supplies $\mathbf b_0$ and a map on $|\cP_i(\mathbf b_0)^{(q)}|$; for the componentwise minimum $\mathbf b:=\min\{\mathbf b_0,\mathbf b'\}$ its restriction is an equivariant map $a\colon|\cP_i(\mathbf b)^{(q)}|\to\cN_i$.
By \cref{fl:lem:capture}, shrink the fixed box in \ref{fl:T:contractible} to $V'_i$
so that every $\gamma$ with $\gamma V'_i\cap V'_i\ne\emptyset$ has a lift in
$A_i(\mathbf b)$. These choices and the index bound are independent of the sphere.

Let $\hat f\colon\mathbb S^j\to \cE_i(V'_i)$ be a lift. Triangulate $\mathbb S^j$ and subdivide so that
the image of the closed star of every vertex $v$ lies in a sheet
$\Sh_{V'_i}(g(v))$, for a lattice element $g(v)$. This is possible by the
Lebesgue number lemma for the inverse images of the sheets. If $v,v'$ are vertices
of one simplex, their sheets meet along the image of that simplex. Thus
$\mathbf{c}_i(g(v)^{-1}g(v'))\in\sB_1$ by \cref{fl:prop:covering}\textup{(ii)}.
The choice of $V'_i$ gives $\mu\in A_i(\mathbf b)$ with
$\vartheta_i(g(v)^{-1}g(v'))=\vartheta_i(\mu)$, and the injectivity on $\sB_2$ (\cref{fl:lem:faithful}\textup{(c)}) gives $g(v)^{-1}g(v')=\mu$. Consequently the vertex
map extends to a simplicial map $g\colon\mathbb S^j\to \cP_i(\mathbf b)^{(q)}$.

For each simplex $\sigma$ choose a vertex and write its sheet label as $\lambda_\sigma$.
Then $\hat f(\sigma)\subset\Sh^{[0]}(\lambda_\sigma)$ and
$ag(\sigma)\subset\lambda_\sigma\mathsf D$.
Since $\cN_i\cong\R^d$, extend $ag$ to a map $H\colon D^{j+1}\to\cN_i$.
Cone the boundary triangulation and subdivide the resulting disk until each closed
simplex $c$ has $H(c)\subset\mu_c\sB$ for some $\mu_c\in\Lambda_i$.
This is possible by the Lebesgue number lemma; the induced boundary triangulation
refines the original one. If $c'$ is a face of $c$ and $y\in c'$, then
\[
 \mathbf c_i(\mu_c^{-1}\mu_{c'})
 =\bigl(\mathbf c_i(\mu_c)^{-1}H(y)\bigr)
  \bigl(\mathbf c_i(\mu_{c'})^{-1}H(y)\bigr)^{-1}
 \in\sB\sB^{-1}\subset\Bx(\mathbf b').
\]
Thus the cell labels satisfy \eqref{fl:eq:labelling}. On a boundary simplex
$c\subset\sigma$, the same calculation, using $H=ag$ on the boundary, gives
$\mathbf c_i(\mu_c^{-1}\lambda_\sigma)\in\sB\mathsf D^{-1}\subset\Bx(\mathbf b')$.
Consequently \eqref{fl:eq:towersheets} gives the exact boundary control
\[
 \hat f(c)\subset\Sh^{[0]}(\lambda_\sigma)\subset\Sh^{[1]}(\mu_c).
\]
Apply \cref{fl:lem:carriers} to $\Lambda_i\times D^{j+1}$ and its boundary, with the
free action on the first factor, cell labels $(\nu,c)\mapsto\nu\mu_c$, and the prescribed
boundary map $(\nu,x)\mapsto\nu\hat f(x)$. The carrier condition holds with $m_0=1$,
and $1+\dim D^{j+1}=j+2\le M$. Restricting the extension to
$\{1\}\times D^{j+1}$ fills $\hat f$ in $\cE_i(V_i^{\rm out})$.
Projection proves the assertion. The filling and the subdivisions may depend on $f$;
the number and size of their sheet labels do not enter the index bound.

A map from $\mathbb S^0$ always lifts, and its null-homotopy joins its two values.
Spheres of dimension at least $2$ lift because they are simply connected.
A loop lifts to a closed loop precisely when its monodromy is trivial, by the definition of
the monodromy and \cref{fl:prop:covering}\textup{(v)}; by \cref{fl:prop:monodromy}\textup{(ii)} the
condition is unchanged under the nested restrictions and enlargements of unwrapping
sets used here.
\end{proof}

\subsection{Intrinsic Fibre Topology}\label{fl:sec:intro}\label{fl:sec:results}\label{fl:sec:chartfree}

We identify the fundamental groups of small pieces of $Y_i$ with the fibre lattices and show that spheres of dimension $0$ and $2\le j\le q$ in small pieces are null-homotopic in larger ones; both statements pass to preimages of realizing maps. The results concern nested neighbourhoods fixed before $i\to\infty$; they identify image groups, not the fundamental group of every individual small preimage. In this subsection, statements carry the centre superscript $a$, and proofs suppress it.

\begin{definition}[Pieces]\label{fl:def:piece}
A \emph{piece} of $Y_i$ is a set $\cU\cap Y_i$ for an open subset $\cU$ of the common realization
$\cY$. For open sets $\cU\subset\cU'$ of $\cY$ and $y_i\in\cU\cap Y_i$, put
\begin{equation*}
\Pi_i(\cU,\cU'):=\im\bigl(\pi_1(\cU\cap Y_i,y_i)\longrightarrow\pi_1(\cU'\cap Y_i,y_i)\bigr),
\end{equation*}
the image of the map induced by the inclusion of the smaller piece in the larger one.
\end{definition}

First, monodromy identifies the image of the fundamental group of a small piece in a
larger one with the fibre lattice.

\begin{theorem}[Fibre Groups]\label{fl:thm:fibre}
Let $y\in Y$, and let $a$ be a centre with $\pr(y)\in B_{R_0/4}(a)$. There is an open
neighbourhood $\cU^0$ of $y$ in $\cY$ such that for all points $y_i\in Y_i$ that converge
to $y$ there are, for all large $i$, homomorphisms
\begin{equation*}
\mathrm{mon}_i\colon\pi_1(\cU^0\cap Y_i,y_i)\longrightarrow\Lambda^a_i
\end{equation*}
with the following property. For every open neighbourhood $\cU'\subset\cU^0$ of $y$ there is
an open neighbourhood $\cU_1\subset\cU'$ of $y$, which does not depend on the points $y_i$,
such that for every open neighbourhood $\cU\subset\cU_1$ of $y$ and all large $i$, the
composition of $\mathrm{mon}_i$ with the map induced by the inclusion of $\cU'\cap Y_i$ in
$\cU^0\cap Y_i$ restricts to an isomorphism of $\Pi_i(\cU,\cU')$ onto $\Lambda^a_i$.
\end{theorem}
\begin{proof}
Fix a frame $p\in\tau^{-1}(y)$ and reference frames $p_i\to p$
over $\Omega^\flat_i$. Let $(s_c,\eta_c)$ be their lifting level and put
$V_i:=N_{p_i}(s_c,\eta_c)$. By \cref{fl:lem:nbhd}\textup{(ii)}
there is an open neighbourhood $\cU^0$ of $y$ in $\cY$ with
$\cU^0\cap Y_i\subset\tau_i(V_i)$ for large $i$.

Given $y_i\to y$, choose frames $x_i\in\tau_i^{-1}(y_i)$ with $x_i\to p$,
by \cref{fl:S:Yconv}. They lie in every fixed box about $p_i$
for large $i$, by \cref{fl:lem:boxconv}\textup{(i)}. Define $\mathrm{mon}_i$ by
restricting the monodromy $\mathrm{mon}^{x_i}_{V_i}$ to $\cU^0\cap Y_i$.

Let $\cU'\subset\cU^0$ be an open neighbourhood of $y$. By
\cref{fl:lem:nbhd}\textup{(i)}, choose a box about $p_i$ whose image
lies in $\cU'\cap Y_i$ for large $i$. Apply \cref{fl:lem:liftablespheres} with $q=1$ inside this box, obtaining $V'_i\subset V_i^{\rm out}$ with $\tau_i(V_i^{\rm out})\subset\cU'\cap Y_i$.
By \cref{fl:lem:nbhd}\textup{(ii)}, choose an open neighbourhood
$\cU_1\subset\cU'$ of $y$ with $\cU_1\cap Y_i\subset\tau_i(V'_i)$.
These choices do not depend on $y_i$.

For any open neighbourhood $\cU\subset\cU_1$ of $y$,
\cref{fl:lem:nbhd}\textup{(i)} provides a smaller box $V''_i\subset V'_i$
about $p_i$ such that, for large $i$,
\begin{equation*}
\tau_i(V''_i)\subset\cU\cap Y_i\subset\tau_i(V'_i)
\subset\tau_i(V_i^{\rm out})\subset\cU'\cap Y_i
\subset\cU^0\cap Y_i\subset\tau_i(V_i).
\end{equation*}
The path-connected box $V''_i$ contains $x_i$ and
$\vartheta_i(u_j)p_i$, for all $j$, since these frames converge to $p$.
Thus $\cE_i(V''_i)$ is path-connected by \cref{fl:prop:covering}\textup{(vii)},
and its monodromy is onto $\Lambda_i$ by \cref{fl:prop:monodromy}\textup{(iii)}.
Monodromy is compatible with restriction, so $\mathrm{mon}_i$ maps
$\Pi_i(\cU,\cU')$ onto $\Lambda_i$.

If a loop in $\cU\cap Y_i$ has trivial $\mathrm{mon}_i$, its lift to $\cE_i(V'_i)$
is closed. By \cref{fl:lem:liftablespheres} it is null-homotopic in
$\tau_i(V_i^{\rm out})\subset\cU'\cap Y_i$. Hence $\mathrm{mon}_i$ is injective on
$\Pi_i(\cU,\cU')$.
\end{proof}

The groups $\Pi_i(\cU,\cU')$ are intrinsic; the isomorphisms $\mathrm{mon}_i$ use a chart and reference frames. Since $\mathrm{mon}_i$ is fixed before shrinking the neighbourhoods, the resulting isomorphisms commute with inclusion whenever both are defined.

To pass from pieces to preimages of realizing maps (\cref{fl:def:realizing}), we compare
realizing maps with the identity in the common realization $\cY$: their uniform displacement
tends to zero, and the preimages of small sets lie between pieces.

\begin{lemma}[Realizing Maps]\label{fl:lem:realizing}
Let $\sd$ be a metric on $\cY$ that induces its topology, and let $F_i\colon Y_i\to Y$ be
continuous maps, defined for all large $i$.
\begin{enumerate}[label=\textup{(\roman*)}]
\item The maps $F_i$ are realizing maps if and only if $\sup_{u\in Y_i}\sd(u,F_i(u))\to0$.
\end{enumerate}
Let $(F_i)$ be realizing maps.
\begin{enumerate}[label=\textup{(\roman*)},resume]
\item If $W\subset Y$ is open and its closure lies in an open subset $\cU$ of $\cY$, then
$F_i^{-1}(W)\subset\cU$ for large $i$.
\item If $\cU\subset\cY$ and $W\subset Y$ are open and $W$ contains the points of $Y$ in
the closure of $\cU$, then $\cU\cap Y_i\subset F_i^{-1}(W)$ for large $i$.
\end{enumerate}
\end{lemma}
\begin{proof}
(i) Uniform displacement tending to zero implies the realizing property by the triangle inequality. Conversely, if the suprema stay above $c>0$ along a subsequence, choose corresponding points $u_i$. Compactness gives a further subsequence converging to $u\in Y$; by \cref{fl:S:Yconv,fl:prop:realization} it extends to a convergent sequence, so the realizing property gives $F_i(u_i)\to u$, a contradiction.

(ii) Otherwise there are indices $i_k\to\infty$ and points $u_k\in Y_{i_k}$ outside $\cU$ with
$F_{i_k}(u_k)\in W$. A subsequence converges to a point $u\in Y$, and $F_{i_k}(u_k)\to u$ by (i),
so $u\in\overline W\subset\cU$. Since $\cU$ is open, $u_k\in\cU$ for large $k$, a contradiction.

(iii) Otherwise there are indices $i_k\to\infty$ and points $u_k\in\cU\cap Y_{i_k}$ with
$F_{i_k}(u_k)\notin W$. A subsequence converges to a point $u\in\overline{\cU}\cap Y\subset W$,
and $F_{i_k}(u_k)\to u$ by (i). The points $F_{i_k}(u_k)$ lie in the closed set $Y\setminus W$,
so $u\notin W$, a contradiction.
\end{proof}

In this way \cref{fl:thm:fibre} passes to preimages of realizing maps.

\begin{proposition}[Fundamental Groups of Realizing Preimages]\label{fl:prop:preimagegroups}
For every $y\in Y$ and every centre $a$ with
$\pr(y)\in B_{R_0/4}(a)$ there is an
open neighbourhood $W^0$ of $y$ in $Y$ such that for every open neighbourhood
$W'\subset W^0$ of $y$ there is an open neighbourhood $W^1\subset W'$ of $y$ with the
following property: for all realizing maps $(F_i)$, every open neighbourhood
$W\subset W^1$ of $y$, all points $y_i\in Y_i$ that converge to $y$ and all large $i$,
the image of $\pi_1(F_i^{-1}(W),y_i)$ in $\pi_1(F_i^{-1}(W'),y_i)$ is isomorphic to
$\Lambda^a_i$, by the isomorphism of \cref{fl:thm:fibre}.
\end{proposition}
\begin{proof}
Let $\cU^0$ be the neighbourhood of \cref{fl:thm:fibre}, and let $W^0$ be an open
neighbourhood of $y$ in $Y$ whose closure lies in $\cU^0$. Let $W'\subset W^0$ be an open
neighbourhood of $y$. Choose an open neighbourhood $\cU'\subset\cU^0$ of $y$ in $\cY$
such that the points of $Y$ in its closure lie in $W'$, let $\cU_1\subset\cU'$ be the
neighbourhood that \cref{fl:thm:fibre} provides for $\cU'$, and let $W^1\subset W'$ be an
open neighbourhood of $y$ in $Y$ whose closure lies in $\cU_1$. These sets do not depend
on the maps. Let $(F_i)$ be realizing maps, let $W\subset W^1$ be an open neighbourhood
of $y$, and choose an open neighbourhood $\cU\subset\cU_1$ of $y$ in $\cY$ such that the
points of $Y$ in its closure lie in $W$. By \cref{fl:lem:realizing}, for large $i$,
\begin{equation*}
\cU\cap Y_i\subset F_i^{-1}(W)\subset\cU_1\cap Y_i\subset\cU'\cap Y_i\subset F_i^{-1}(W')
\subset\cU^0\cap Y_i .
\end{equation*}
Let $y_i\to y$; then $y_i\in\cU$ for large $i$. The homomorphism $\mathrm{mon}_i$ is defined on
the fundamental group of $F_i^{-1}(W')$ through the last inclusion. By \cref{fl:thm:fibre}
for the neighbourhoods $\cU\subset\cU'$, the classes of the loops in $\cU\cap Y_i$ are
mapped onto $\Lambda_i$; so $\mathrm{mon}_i$ maps the image of $\pi_1(F_i^{-1}(W),y_i)$ in
$\pi_1(F_i^{-1}(W'),y_i)$ onto $\Lambda_i$. By \cref{fl:thm:fibre} for the neighbourhoods
$\cU_1\subset\cU'$, a loop in $\cU_1\cap Y_i$ at $y_i$ whose class is mapped to the identity
is null-homotopic in $\cU'\cap Y_i$; so $\mathrm{mon}_i$ is injective on this image.
\end{proof}

The identification of the fibres in \cref{es:sec:fibreframes} also uses the following control.

\begin{theorem}[Connectedness and Higher Homotopy of the Pieces]\label{fl:thm:higher}
In the situation of \cref{fl:thm:fibre}, let $q\ge2$. For every open neighbourhood
$\cU'\subset\cU^0$ of $y$ there is an open neighbourhood $\cU_2\subset\cU'$ of $y$
such that, for every open neighbourhood $\cU\subset\cU_2$ of $y$ and all large $i$:
\begin{enumerate}[label=\textup{(\roman*)}]
\item Any two points of $\cU\cap Y_i$ are joined by a path in $\cU'\cap Y_i$.
\item Every continuous map $\mathbb S^j\to\cU\cap Y_i$, for $2\le j\le q$, is
null-homotopic in $\cU'\cap Y_i$.
\end{enumerate}
The same assertions hold for realizing maps:
for every open neighbourhood $W'\subset W^0$ of $y$ in $Y$, with $W^0$ as in
\cref{fl:prop:preimagegroups}, there is an open neighbourhood $W^2\subset W'$ of $y$
such that the assertions hold with $F_i^{-1}(W)$ and $F_i^{-1}(W')$ in place of
$\cU\cap Y_i$ and $\cU'\cap Y_i$, for all realizing maps $(F_i)$, every open
neighbourhood $W\subset W^2$ of $y$, and all large $i$.
\end{theorem}
\begin{proof}
Use the reference frames $p_i\to p$ and the box $V_i$ of their lifting level from the proof of \cref{fl:thm:fibre}. Given $\cU'\subset\cU^0$, apply \cref{fl:lem:nbhd}\textup{(i)} and \cref{fl:lem:liftablespheres} to obtain $V'_i\subset V_i^{\rm out}$ with $\tau_i(V_i^{\rm out})\subset\cU'\cap Y_i$.
By \cref{fl:lem:nbhd}\textup{(ii)}, choose
$\cU_2\subset\cU'$ with $\cU_2\cap Y_i\subset\tau_i(V'_i)$.
Any map from $\mathbb S^0$, or from $\mathbb S^j$ with $2\le j\le q$, into
$\cU\cap Y_i\subset\tau_i(V'_i)$ lifts to $\cE_i(V'_i)$. The conclusion
follows from \cref{fl:lem:liftablespheres}.

For realizing maps, use the neighbourhood construction in the proof of
\cref{fl:prop:preimagegroups}, with $\cU_2,W^2$ in place of $\cU_1,W^1$.
The choices again precede the maps. For all realizing maps $(F_i)$ and every
open neighbourhood $W\subset W^2$ of $y$, \cref{fl:lem:realizing} gives, as in that proof, for all large $i$,
\[
F_i^{-1}(W)\subset\cU_2\cap Y_i,\qquad
\cU'\cap Y_i\subset F_i^{-1}(W').
\]
Apply the preceding sphere argument to maps into the first set; its paths and
null-homotopies lie in $\cU'\cap Y_i$, hence in the second set.
\end{proof}

The neighbourhoods depend on $q$; for bundle fibres it will suffice to take $q\ge\dim Y_i$.

\subsection{Nilmanifold Fibres}\label{es:sec:rigidity}\label{es:sec:fibreframes}

For the rest of this section impose also the gluing threshold \eqref{es:eq:thresholdG}, as in \cref{es:thm:gluingall}. Then $Y$ is a topological manifold by \cref{ls:thm:limits}. We first identify fibres of arbitrary realizing maps that are bundles locally, and then apply this to $F_i^\natural$.

An \emph{infranilmanifold} is the quotient of a
simply connected nilpotent Lie group $\cN'$ by a torsion-free discrete cocompact
subgroup of $\cN'\rtimes C$, for a compact group $C$ of automorphisms of $\cN'$; it is a
closed smooth manifold with a contractible universal covering space, and it is finitely
covered by a nilmanifold \cite[Theorem~3.4]{Dek17}. We use the following theorem.
\begin{enumerate}[label=\textup{(Rig)},ref=\textup{(Rig)},leftmargin=3.4em]
\item\label{es:Rig} A closed topological manifold that is homotopy equivalent to an
infranilmanifold is homeomorphic to it.
\end{enumerate}
In dimensions at least $5$ this is the theorem of Farrell and Hsiang on closed
aspherical manifolds with virtually nilpotent fundamental groups
\cite[Theorem~5.1]{FH83}. In dimension $4$ it is a theorem of Freedman and Quinn for
aspherical manifolds whose fundamental groups are poly-(finite or cyclic)
\cite[Section~11.5]{FQ90}. In dimension $3$ it follows from results of Waldhausen and
Turaev and from the geometrization theorem, and in dimensions $1$ and $2$ from the
classification of closed manifolds \cite[Remark~4.5]{Luc10}. We also use that
topological manifolds have the homotopy type of CW complexes \cite{Mil59} and that the
homology of a manifold vanishes above its dimension \cite[Section~3.3]{Hat02}.

\begin{theorem}[Fibres of Realizing Maps]\label{es:thm:fibrenil}
In the setting of \cref{es:sec:standing}, let $F'_i\colon Y_i\to Y$ be realizing maps,
for $i$ in an infinite subset of the common subsequence of \cref{fl:sec:malcev}, let $y\in Y$, let $a$ be a centre with
$\pr(y)\in B_{R_0/4}(a)$, and let $W_*$ be an open neighbourhood of $y$ in $Y$ such that
for all large $i$ the restriction of $F'_i$ to the preimage of $W_*$ is a fibre bundle
over $W_*$ whose fibre $\mathsf F_i$ is a closed manifold. Then for all large $i$ the
manifold $\mathsf F_i$ is connected, its universal covering space is contractible, its
fundamental group is isomorphic to $\Lambda^a_i$, and it is homeomorphic to
$\operatorname{Nil}^a_i$.
\end{theorem}
\begin{proof}
Let $q:=\max\{2,\dim Y_i\}=\max\{2,n+\tfrac12n(n-1)\}$; the dimension of
$\mathsf F_i$ is at most $q$. Let
$W^0$ be the neighbourhood of $y$ of \cref{fl:prop:preimagegroups}, let
$W'\subset W^0\cap W_*$ be an open neighbourhood of $y$ homeomorphic to a
Euclidean space, let $W^1\subset W'$ be the neighbourhood of
\cref{fl:prop:preimagegroups}, let $W^2\subset W'$ be the neighbourhood of
\cref{fl:thm:higher} for the number $q$ and $W'$, and let
$W\subset W^1\cap W^2$ be an open neighbourhood of $y$ homeomorphic to a Euclidean
space. Put
\begin{equation*}
U_i:=F'^{-1}_i(W)\subset U'_i:=F'^{-1}_i(W').
\end{equation*}
A fibre bundle over a contractible paracompact space is trivial: its associated principal
bundle is numerable, hence trivial \cite[Chapter~4, Section~9 and Corollary~10.3]{Hus94}. So there is a homeomorphism of $U'_i$ onto
$\mathsf F_i\times W'$ over $W'$, which maps $U_i$ onto $\mathsf F_i\times W$. Hence the
inclusion of $U_i$ in $U'_i$ is a homotopy equivalence, and both spaces are homotopy
equivalent to $\mathsf F_i$. We choose points $y_i\in Y_i$ that converge to $y$; they
lie in $U_i$ for large $i$.

By \cref{fl:thm:higher}, any two points of $U_i$ can be joined in $U'_i$, and the
images of $\pi_j(U_i,y_i)$ in $\pi_j(U'_i,y_i)$ vanish for $2\le j\le q$. Since the
inclusion is a homotopy equivalence, $U'_i$ and $\mathsf F_i$ are connected and
$\pi_j(\mathsf F_i)=0$ in these degrees. By \cref{fl:prop:preimagegroups,fl:thm:fibre}, the image
of $\pi_1(U_i,y_i)$ in $\pi_1(U'_i,y_i)$ is isomorphic to $\Lambda^a_i$.
The inclusion induces an isomorphism, so $\pi_1(\mathsf F_i)\cong\Lambda^a_i$.

The universal covering space of $\mathsf F_i$ is a simply connected manifold of
dimension at most $q$ whose homotopy groups vanish in the degrees up to $q$. By the
theorem of Hurewicz its homology vanishes in the positive degrees up to $q$, and it
vanishes above the dimension; so all its homotopy groups vanish, and it is
contractible, by the theorem of Whitehead. Thus $\mathsf F_i$ and
$\operatorname{Nil}^a_i$ are connected spaces of the homotopy type of CW complexes with
contractible universal covering spaces and isomorphic fundamental groups; so they are
homotopy equivalent, and $\mathsf F_i$ is homeomorphic to $\operatorname{Nil}^a_i$ by
\ref{es:Rig}.
\end{proof}

The bundles of the gluing theorem are realizing maps that are fibre bundles over all of $Y$,
so the theorem applies to them at every point.

\begin{corollary}[The Fibres Are Nilmanifolds]\label{es:cor:fibresnil}
In the situation of \cref{es:thm:gluingall}, for all large $i$, every fibre of $F^\natural_i$
over a point of $\pr^{-1}(B_{R_0/4}(a))$, for a centre $a$, is homeomorphic to
$\operatorname{Nil}^a_i$, and so are the manifolds $\mathsf F_i$ of
\cref{es:thm:gluingall}\textup{(iii)} for the points $x\in B_{R_0/4}(a)$ and those of
\cref{es:thm:gluingall}\textup{(iv)} for the tubes at these points.
\end{corollary}
\begin{proof}
The maps $F^\natural_i$ are realizing maps and fibre bundles over $Y$. The compact
manifold $Y$ has finitely many components, and the fibres of $F^\natural_i$ over one
component are homeomorphic. For every component and every centre $a$ such that the
component meets $\pr^{-1}(B_{R_0/4}(a))$ we apply \cref{es:thm:fibrenil} at one point $y$ of
the intersection, with $W_*$ equal to that open connected component of $Y$. By \cref{es:thm:gluingall}\textup{(iii)} and
\textup{(iv)}, the manifolds $\mathsf F_i$ there are fibres of $F^\natural_i$.
\end{proof}

\subsection{Isotropy Descent}\label{es:sec:fibrestages}

The fibres of $f^\natural_i$ are the quotients of the nilmanifold fibres of $F^\natural_i$ by free actions of the isotropy groups. The following lemma applies because \cref{es:thm:gluingall}\textup{(iv)} already makes these quotients manifolds.

\begin{lemma}[Compact Quotients of Aspherical Nilpotent Manifolds]\label{es:lem:quotientfibre}
Let $\mathsf F$ be a closed connected aspherical topological manifold with finitely generated
nilpotent fundamental group $\Lambda$. Suppose that a compact Lie group $I$ acts freely
on $\mathsf F$, that $I_0$ is a torus, and that $\mathsf M:=\mathsf F/I$ is a topological manifold. Then $\mathsf M$ is
closed and aspherical, its fundamental group is torsion-free and virtually nilpotent,
and it is homeomorphic to an infranilmanifold. The orbit map induces an injection
$\pi_1(I_0)\to\Lambda$ with central image $P$, and there is an exact sequence
\begin{equation}\label{es:eq:quotientextension}
1\longrightarrow\Lambda/P\longrightarrow\pi_1(\mathsf M)\longrightarrow\pi_0(I)\longrightarrow1.
\end{equation}
\end{lemma}
\begin{proof}
Write $I_0=\mathbb T^m$. A free action of a compact Lie group is a principal bundle \cite{Pal61},
\cite[Chapter~II, Theorem~5.8]{Bre72}. Moreover $\mathsf M_0:=\mathsf F/I_0\to \mathsf M$ is a finite covering, since
$I/I_0$ acts freely on $\mathsf M_0$. In particular $\mathsf M_0$ is a closed manifold.
The homotopy sequence of $\mathbb T^m\to\mathsf F\to\mathsf M_0$ gives
\begin{equation*}
\pi_j(\mathsf M_0)=0\quad(j\ge3),\qquad
\pi_2(\mathsf M_0)\cong\ker(\mathbb Z^m\to\Lambda)\cong\mathbb Z^r
\end{equation*}
for some $r\ge0$. The universal cover of $\mathsf M_0$ has CW type, so it is
homotopy equivalent to $K(\mathbb Z^r,2)\simeq(\mathbb{CP}^{\infty})^r$
\cite[Proposition~4.30 and Example~4.50]{Hat02}. If $r>0$, its rational homology is non-zero in
arbitrarily high degrees, contradicting the finite dimension of this covering
manifold. Thus $r=0$, the orbit map is injective on $\pi_1(I_0)$, and
$\mathsf M_0$ and $\mathsf M$ are aspherical. The image $P$ of the orbit map is central: for a loop
$a$ in $I_0$ based at the identity and a loop $c$ in $\mathsf F$ based at $x$,
the square $(s,t)\mapsto a(t)c(s)$ fills the commutator of $c$ and the orbit loop
$t\mapsto a(t)x$. The principal bundle and finite covering give \eqref{es:eq:quotientextension};
hence $\pi_1(\mathsf M)$ is finitely generated and virtually nilpotent. It is torsion-free, by the
standard argument for groups of finite cohomological dimension
\cite[Chapter~VIII, Corollary~2.5]{Bro82}: a subgroup of prime order $p$ would give a covering
manifold of $\mathsf M$ that is a $K(\mathbb Z/p,1)$, whose mod-$p$ homology is non-zero in arbitrarily
high degrees.

A finitely generated torsion-free virtually nilpotent group is the fundamental group
of an infranilmanifold \cite[Corollary~4.3]{Dek17}. This model and $\mathsf M$ are aspherical
spaces of CW type with the same fundamental group, so they are homotopy equivalent.
They are homeomorphic by \ref{es:Rig}; in dimension zero both are a point.
\end{proof}

We now pass from the fibres of the bundles $F^\natural_i$ to the fibres of the induced maps
$f^\natural_i$.

\begin{theorem}[The Fibres over the Limit]\label{es:thm:infranil}
In the situation of \cref{es:thm:gluingall}, for all large $i$ the following holds
at every point $x'\in X$. Let $I'$ be the isotropy group of a point of $Y$ over $x'$,
and let $a$ be a centre with $x'\in B_{R_0/4}(a)$. Then the fibre of $f^\natural_i$ over $x'$ is a closed manifold of dimension $\dim G^a-\dim I'$. It is the quotient of a fibre of $F^\natural_i$ by a free action of
$I'$, and it is homeomorphic to an infranilmanifold.
\end{theorem}
\begin{proof}
\emph{The local form at every point.} Let $x_\alpha$ be a tube of the system of
\cref{es:lem:system}, and let $\mathsf F_i$ be the manifolds of
\cref{es:thm:gluingall}\textup{(iv)} for this tube. By that theorem, for all large $i$, with one
index bound for all $v\in B_V(r_\ell)$, the fibre of $f^\natural_i$ over
$x_v:=\bar\xi^{-1}(\mathsf R(K)v)$ is the quotient $\mathsf F_i/K_v$ by the free action of
$K_v$, whose image $I_v=\{A_k:k\in K_v\}$ is the isotropy group of $\Xi([\one,v])$, and it is a closed manifold of dimension
$\dim\mathsf F_i-\dim K_v$. The manifold $\mathsf F_i$ is homeomorphic to the fibre of
$F^\natural_i$ over $\Xi([\one,0])$, hence to the fibre over $\Xi([\one,v])$, because $F^\natural_i$ is
a fibre bundle and $t\mapsto\Xi([\one,tv])$ joins the two points.

\emph{Classification.} Let $a$ be a centre with $x_v\in B_{R_0/4}(a)$. By
\cref{es:cor:fibresnil}, the fibre of $F^\natural_i$ over $\Xi([\one,v])$, and
hence $\mathsf F_i$, is homeomorphic to $\operatorname{Nil}^a_i$, a closed aspherical manifold of
dimension $\dim G^a$ with the fundamental group $\Lambda^a_i$ (\cref{fl:lem:nilmanifold}). The
connected component $(K_v)_0$ is a torus: it is a
connected compact subgroup of $G_0$, whose maximal compact subgroup is a central
torus by \cref{m:fr:prop:limitgroup}\textup{(ii),(iii)}. Apply
\cref{es:lem:quotientfibre} to $\mathsf F_i$ and $K_v$. This proves the assertion at $x_v$, for
every $v$ and every such centre, without a further subsequence or index bound. The finitely
many tubes cover $X$ by \cref{es:lem:system}\textup{(i)}, so the conclusions hold
at every point of $X$.
\end{proof}

By \cref{es:thm:gluingall}\textup{(iii)} and \cref{es:cor:fibresnil}, near a point $x'$ as in
the theorem $f^\natural_i$ is the map
\begin{equation*}
\operatorname{Nil}^a_i\times_{I'}B\longrightarrow B/I'
\end{equation*}
for a free action of $I'$ on the nilmanifold. As in \cite{Fuk89,CFG92}, the
fibre dimension drops by $\dim I'$. Here the local products are topological,
and the action is not known to be affine. For finite $I'$ it is affine in suitable
coordinates, by the following remark.

\begin{remark}[Finite Isotropy Groups]\label{es:rem:affine}
Let $I'$ be finite, let $\mathsf F_i$ be the manifold of \cref{es:thm:gluingall}\textup{(iii)}
at $x'$, so that $\mathsf F_i/I'$ is the fibre of $f^\natural_i$ over $x'$, and let $h$ be a
homeomorphism of $\mathsf F_i/I'$ onto an infranilmanifold $\operatorname{Inf}_i$
(\cref{es:thm:infranil}). The covering of $\operatorname{Inf}_i$ that corresponds under $h$
to the covering $\mathsf F_i\to\mathsf F_i/I'$ is an infranilmanifold with a nilpotent
fundamental group. This group is a normal subgroup of the fundamental group of $\operatorname{Inf}_i$, because
the covering is regular. The translations in the fundamental group of an infranilmanifold
form a normal subgroup that is maximal among its nilpotent subgroups
\cite[Proposition~3.6]{Dek17}, and the product of two nilpotent normal subgroups of a group
is nilpotent, by the theorem of Fitting. So the fundamental group of the covering consists of
translations, and the covering space is a nilmanifold
$\mathsf N_i$, on which the group of deck transformations acts by affine
diffeomorphisms, and $h$ lifts to a homeomorphism of $\mathsf F_i$ onto $\mathsf N_i$
that commutes with the actions of $I'$. Thus near a point with a finite isotropy group
the map $f^\natural_i$ is the map $\mathsf N_i\times_{I'}B\to B/I'$ for a nilmanifold
$\mathsf N_i$ with a free affine action of $I'$ and a linear action on the ball. For isotropy groups of positive
dimension we have no corresponding statement.
\end{remark}

The last remark concerns the structure group, which the construction does not control.

\begin{remark}[Affine Structures]\label{es:rem:affinebundles}
The construction identifies fibres, not an affine structure group for the given bundle. In bounded-curvature collapse, the affine structure comes from local nilpotent isometric actions for a nearby invariant metric; no such actions on $X_i$ are constructed here. Even smooth bundles need not admit affine reduction: a torus diffeomorphism acting trivially on $\pi_1$ need not be isotopic to the identity in dimensions at least $6$ \cite{HS76,Hat78}, whereas an affine map with that action is a translation. Its mapping torus supplies the obstruction for smooth bundles; a topological obstruction would require the corresponding statement for homeomorphisms.

There is nevertheless an affine model up to fibre homotopy. Fix a centre and suppress its superscript. For $\operatorname{Nil}_i=\Lambda_i\backslash\cN_i$, the affine group has identity component $\cN_i/Z(\Lambda_i)$ and component group $\operatorname{Out}(\Lambda_i)$. The space of self-homotopy equivalences has the same components, each a $K(Z(\Lambda_i),1)$ \cite[Theorem~III.2]{Got65}, \cite[Remark~7.19]{Luc25}; the inclusion induces these identifications and is a weak homotopy equivalence. The classification theorem of \cite{Sta63} therefore makes $F_i^\natural$, over each component of $Y$, fibre homotopy equivalent to an affine nilmanifold bundle. This does not reduce the structure group of the original topological bundle.
\end{remark}

\subsection{Proof of the Main Theorem}\label{es:sec:proofsintro}

The proof collects the results of the preceding sections: the atlas and framed limit of
\cref{ch:coveratlas}, the manifold structure of \cref{ls:thm:limits}, the realizing maps of \cref{fl:thm:anr}, the gluing theorem of \cref{ch:globalbundle} and the fibre identification
above.

\begin{proof}[Proof of \cref{thm:main}]
Choose $\delta_1=\min\{\delta_*,\eps_{\mathrm{glue}}(n,N)/3\}$ and $c_1=c_*$.
\Cref{m:fr:thm:geomain,fl:prop:realization,fl:lem:invariant} give the frame
spaces and their compact realization. The threshold \eqref{es:eq:thresholdG} gives the hypothesis $3\delta\le\eps_L(n,N+1)$ of
\cref{ls:thm:limits}, since $\eps_{\mathrm{glue}}\le\eps_T\le\eps_L(n,N+1)$ by
\eqref{th:eq:threshold}; that theorem makes $Y$ a manifold with a locally linear
$O(n)$-action (\cref{ls:def:locallylinear}). \Cref{fl:thm:anr} supplies equivariant realizing maps.

Fix such maps and $\eps>0$. Uniform continuity of $\pr\colon Y\to X$ gives
$\eps'\in(0,\eps]$ such that $\sd(u,u')<\eps'$ implies
$d(\pr(u),\pr(u'))<\eps$. If infinitely many indices failed
\textup{(ii)--(iii)}, apply \cref{es:thm:gluingall} to that infinite set.
On a further subsequence it gives equivariant bundles within $\eps'$ of the
chosen maps, hence quotient maps within $\eps$.
Its part \textup{(iii)}, whose index bound does not depend on the point, gives for
all large $i$ and every $x\in X$ the neighbourhood $O_x:=\pr(\Xi(O(n)\times_IB))$ with
the twisted product and the associated $O(n)$-bundle, and part \textup{(iv)} gives the
fibres $\mathsf F/I_v$ with their dimension. The dimension of $B$ is $n-\dim\mathsf F+\dim I$ by
\cref{ls:thm:slices}\textup{(ii)}, because $\dim\mathsf F$ is the dimension of the limit
group of the chart. This contradicts the failure and proves \textup{(ii)} and
\textup{(iii)} along the whole subsequence of \cref{m:fr:thm:geomain}.

For \textup{(iv)}, restrict to the common subsequence of
\cref{fl:sec:malcev}. Apply the same contradiction argument to indices
failing the combined geometric and fibre conclusions.
The bundles supplied by \cref{es:thm:gluingall} are realizing maps, so
\cref{es:cor:fibresnil,es:thm:infranil} identify their fibres and the local
manifolds $\mathsf F$ as claimed. This yields the required contradiction.
\end{proof}

\section{Limits of Aspherical Spaces}\label{as:sec:aspherical}

This section proves \cref{as:thm:aspherical}. Under bounded sectional curvature, the affine
structure group makes lifted fibres into nilpotent orbits on the contractible universal cover,
where fixed-point theory excludes isotropy \cite[Remark~1.14 and Theorem~0.5]{Ron26}. Our bundle
has no affine structure group (\cref{es:rem:affinebundles}). Instead, we lift its Borel
construction, represent the kernel of the homomorphism from the fibre group to $\pi_1(X_i)$ by a
torus bundle, and
use Poincar\'e duality to force this kernel to vanish. The three subsections carry out the lift,
the torus argument, and the conclusion, including the type-$A_s$ case.

\emph{Standing notation.} Fix $i$ as in \cref{as:thm:aspherical}, with $X_i$ aspherical.
For $y\in Y$, write $I_y\subset O(n)$ for its isotropy group and
$\mathsf F_y:=(F^\natural_i)^{-1}(y)$ for the closed nilmanifold fibre. By
\cref{thm:main}\,\textup{(iii)}, the local form is
\begin{equation}\label{as:eq:localform}
O(n)\times_{I_y}(\mathsf F_y\times B_y)\longrightarrow O(n)\times_{I_y}B_y,
\end{equation}
where $B_y$ is a ball about $0$ in an orthogonal representation and $I_y$ acts diagonally,
freely on $\mathsf F_y$. Part~\textup{(iii)} gives this form at one point over each point of $X$;
the points of $Y$ over a point of $X$ form one orbit, so conjugation gives it at every $y$.
The fibre $\mathsf F_y/I_y$ of $f^\natural_i$ is an aspherical closed manifold
by \cref{es:thm:infranil}, and $X_i$ is a closed $n$-manifold, locally
$(\mathsf F_y\times B_y)/I_y$. Since $X$ is connected, $O(n)$ permutes the components of $Y$
transitively, and all $\mathsf F_y$ are homeomorphic. Put
\begin{equation*}
\Lambda:=\pi_1(\mathsf F_y),\qquad
\Lambda_\partial:=\ker[\Lambda\to\pi_1(X_i)],\qquad
\bar\Lambda:=\operatorname{im}[\Lambda\to\pi_1(X_i)],\qquad \sQ:=\pi_1(X_i)/\bar\Lambda.
\end{equation*}
The homomorphism is induced by $\mathsf F_y\hookrightarrow Y_i\to X_i$; $\bar\Lambda$ is normal by
\cref{as:lem:borel}. Write $\gamma\mapsto\bar\gamma$ for the quotient map and
$\tilde X_i$ for the universal cover, with its left deck action of $\pi_1(X_i)$.

\subsection{Borel Lifting}\label{as:sec:lifted}

For a compact Lie group $K$, write $S_{hK}:=EK\times_KS$, using the standard numerable
universal principal bundle $EK\to BK$ and the same $EK$ for closed subgroups.
All manifolds here are metrizable and have CW homotopy
type. A free compact Lie action on such a space is a numerable principal bundle, so
$S_{hK}\to S/K$ is a homotopy equivalence \cite{Pal61,Bre72}. Locally linear actions have
equivariant CW homotopy type \cite[Theorem~2.5]{Kha18}. On spaces of this type, the Borel
construction classifies equivariant principal bundles with discrete structure group, and also
with torus structure group: it gives a natural bijection with ordinary principal bundles on
$S_{hK}$ \cite[Theorems~1 and~3]{May90}, \cite{LMS83}. We use both forms below.

\begin{lemma}[Freeness from the Borel Construction]\label{as:lem:borelfree}
Let a compact Lie group $K$ act on a Hausdorff space $S$. If $S_{hK}$ is weakly contractible,
then the action is free.
\end{lemma}
\begin{proof}
A non-trivial stabilizer contains a cyclic subgroup $C_p$ of prime order. The fibre bundle
$K/C_p\to S_{hC_p}\to S_{hK}$ makes $S_{hC_p}$ weakly equivalent to the compact manifold $K/C_p$,
so its mod-$p$ cohomology is bounded in degree. But a $C_p$-fixed point splits the projection
to $BC_p$, injecting the unbounded cohomology
$H^*(BC_p;\mathbb F_p)$ into it, a contradiction.
\end{proof}

The Borel construction of $F^\natural_i$ compares $\Lambda$ with $\pi_1(X_i)$ and realizes the
kernel $\Lambda_\partial$ as the second homotopy group of $Y_{hO(n)}$.

\begin{lemma}[The Borel Fibration and Its Kernel]\label{as:lem:borel}\label{as:lem:homotopy}
The map $(F^\natural_i)_{hO(n)}$ is a fibre bundle with fibre $\mathsf F_y$ and total space
homotopy equivalent to $X_i$. Under this equivalence its fibre inclusion induces
$\Lambda\to\pi_1(X_i)$. The group $\bar\Lambda$ is normal, $\Lambda_\partial$ is free abelian
of finite rank, and
\begin{equation*}
\pi_1(Y_{hO(n)})=\sQ,\qquad \pi_2(Y_{hO(n)})=\Lambda_\partial,\qquad
\pi_j(Y_{hO(n)})=0\quad(j\ge3).
\end{equation*}
In particular the universal cover of $Y_{hO(n)}$ has homotopy type $K(\Lambda_\partial,2)$.
\end{lemma}
\begin{proof}
The local form \eqref{as:eq:localform} induces the associated bundle
\begin{equation*}
EO(n)\times_{I_y}(\mathsf F_y\times B_y)\longrightarrow EO(n)\times_{I_y}B_y.
\end{equation*}
These neighbourhoods cover the base. Since $O(n)$ acts freely on $Y_i$, its Borel construction
is homotopy equivalent to $X_i$, and, for $u\in EO(n)$, the map $[u,z]\mapsto O(n)z$ identifies
the fibre inclusion with
the stated map. Asphericity of $X_i$ and of $\mathsf F_y$ gives the exact sequence
\begin{equation*}
0\longrightarrow\pi_2(Y_{hO(n)})\longrightarrow\Lambda\longrightarrow\pi_1(X_i)
\longrightarrow\pi_1(Y_{hO(n)})\longrightarrow0
\end{equation*}
and vanishing of higher homotopy groups. Thus $\bar\Lambda$ is normal and
$\Lambda_\partial$ is abelian; as a subgroup of the finitely generated torsion-free nilpotent
group $\Lambda$, it is free abelian of finite rank. The CW-type assertion identifies the
universal cover with $K(\Lambda_\partial,2)$ \cite[Proposition~4.30]{Hat02}.
\end{proof}

Next we unwrap $\pi_1(Y_{hO(n)})\cong\sQ$ by an equivariant covering of $Y$.

\begin{lemma}[The Lifted Frame Limit]\label{as:lem:lifted}
There is an $O(n)$-equivariant principal $\sQ$-covering $\hat Y\to Y$, with right $\sQ$-action
commuting with $O(n)$, whose Borel construction is the universal cover of $Y_{hO(n)}$.
Put $\hat X:=\hat Y/O(n)$.
\begin{enumerate}[label=\textup{(\roman*)}]
\item The space $\hat Y$ is a manifold, $O(n)$ permutes its components transitively, and
$\hat X$ is connected and finite-dimensional.
\item If $\hat y$ lies over $y$, the homomorphism $\chi_{\hat y}\colon I_y\to \sQ$ defined by
$h\hat y=\hat y\chi_{\hat y}(h)$ has kernel $I_{\hat y}$ and finite image. Thus
$(I_y)_0=(I_{\hat y})_0$, and $I_y/I_{\hat y}$ is the stabilizer of $O(n)\hat y$ in $\sQ$.
An invariant neighbourhood of $O(n)\hat y$ is $O(n)\times_{I_{\hat y}}B_y$.
\item The action of $\sQ$ on $\hat X$ is properly discontinuous and cocompact, has finite
stabilizers, and has quotient $X$.
\end{enumerate}
\end{lemma}
\begin{proof}
Apply the discrete-group bundle classification to the universal covering of $Y_{hO(n)}$, whose
deck group is $\sQ$ by \cref{as:lem:borel}.
Its connectedness implies transitivity on components in (i), and a covering of $Y$ is a
manifold. The action of $I_y$ on the $\sQ$-torsor over $y$ defines $\chi_{\hat y}$; its image is
finite because $I_y$ is compact and $\sQ$ discrete. The stabilizer assertion follows from
$\hat y\bar\gamma=h\hat y$.

Lift the slice $B_y\to Y$ to $s\colon B_y\to\hat Y$ through $\hat y$. Uniqueness of lifts gives
$hs(v)=s(hv)\chi_{\hat y}(h)$. Consequently $[g,v]\mapsto gs(v)$ identifies
$O(n)\times_{I_{\hat y}}B_y$ with the part of the lifted tube containing the orbit:
an equality $gs(v)=g's(v')$ first gives $g'=gh^{-1}$, $v'=hv$ for $h\in I_y$, and then
$\chi_{\hat y}(h)=1$. This is an isomorphism of coverings onto that open part, proving the
local form and local linearity. Finally the $\sQ$-action on $\hat Y$ is proper, and
$\hat Y\to\hat X$ is proper because $O(n)$ is compact. Properness descends to $\hat X$, with
finite stabilizers as in (ii) and compact quotient $X$.
\end{proof}

Put $\tilde Y_i:=Y_i\times_{X_i}\tilde X_i$. It is a principal $O(n)$-bundle over $\tilde X_i$
and a principal $\pi_1(X_i)$-covering of $Y_i$, with commuting actions.

\begin{lemma}[The Lifted Bundle]\label{as:lem:liftedbundle}
There is an $O(n)$-equivariant lift $\hat F^\natural_i\colon\tilde Y_i\to\hat Y$ of
$F^\natural_i$ such that
\begin{equation*}
\hat F^\natural_i(\gamma z)=\hat F^\natural_i(z)\bar\gamma^{-1},\qquad
\hat f^\natural_i(\gamma x)=\hat f^\natural_i(x)\bar\gamma^{-1},
\end{equation*}
where $\hat f^\natural_i\colon\tilde X_i\to\hat X$ is its quotient. Over the tube at $\hat y$,
the lifted bundle has the form
\begin{equation}\label{as:eq:liftedform}
O(n)\times_{I_{\hat y}}(\tilde{\mathsf F}_{\hat y}\times B_y)
\longrightarrow O(n)\times_{I_{\hat y}}B_y.
\end{equation}
Here $\tilde{\mathsf F}_{\hat y}$ is the connected covering of $\mathsf F_y$ with subgroup
$\Lambda_\partial$, hence is aspherical. The fibre of $\hat f^\natural_i$ is
$\tilde{\mathsf F}_{\hat y}/I_{\hat y}$, and the lifted Borel fibration has boundary isomorphism
\begin{equation}\label{as:eq:boundary}
\partial\colon\pi_2(\hat Y_{hO(n)})\xrightarrow{\ \cong\ }
\pi_1(\tilde{\mathsf F}_{\hat y})=\Lambda_\partial.
\end{equation}
\end{lemma}
\begin{proof}
The pullback of $\hat Y\to Y$ to $\tilde Y_i$ has trivial Borel covering, since
$(\tilde Y_i)_{hO(n)}\simeq\tilde X_i$ is contractible. The discrete-group classification
therefore gives an equivariant section, hence $\hat F^\natural_i$. The original Borel map
induces the quotient $\gamma\mapsto\bar\gamma$, and its lift has the displayed equivariance
after choosing compatible basepoints. The inverse occurs because the source deck action is
on the left and the target action on the right.

The components of the inverse image of $\mathsf F_y$ in $\tilde Y_i$ are indexed by $\sQ$,
and each is the covering classified by $\Lambda_\partial$. Equivariance identifies them
bijectively with the $\sQ$-torsor over $y$, proving the fibre assertion. Over the contractible
slice $B_y$, each lifted component is $\tilde{\mathsf F}_{\hat y}\times B_y$; uniqueness of
lifts, with agreement at $0$, identifies its $I_{\hat y}$-action with the product action.
This gives \eqref{as:eq:liftedform} and its quotient. Finally the total Borel space is
contractible, so its homotopy sequence gives \eqref{as:eq:boundary}; naturality identifies
its image in $\Lambda$ with the original kernel $\Lambda_\partial$.
\end{proof}

\subsection{The Torus of the Kernel}\label{as:sec:torus}

Set $\mathbb T_\partial:=(\Lambda_\partial\otimes\R)/\Lambda_\partial$. We realize this torus
on a contractible manifold whose orbit space is $\hat X$. Its fixed set detects maximal
isotropy, and Poincar\'e duality forces that locus to be all of $\hat X$.

\begin{proposition}[An Equivariant Torus Bundle]\label{as:lem:torusbundle}\label{as:lem:torusaction}
There is an $O(n)$-equivariant principal $\mathbb T_\partial$-bundle $\cB\to\hat Y$ with the
following properties.
\begin{enumerate}[label=\textup{(\roman*)}]
\item The space $\cB_{hO(n)}$ is weakly contractible, and $O(n)$ acts freely on $\cB$.
\item The action on $\cB_{\hat y}$ defines an injective homomorphism
$\psi_{\hat y}\colon I_{\hat y}\to\mathbb T_\partial$. If
$b_{\hat y}\colon BI_{\hat y}\to\hat Y_{hO(n)}$ is the orbit inclusion, then
\begin{equation}\label{as:eq:orbitmarking}
(\psi_{\hat y})_*=\partial\circ(b_{\hat y})_*=
\bigl[\pi_1(I_{\hat y})\longrightarrow\pi_1(\tilde{\mathsf F}_{\hat y})\bigr],
\end{equation}
where the last map is induced by the orbit map. In particular
$\dim I_y\le\operatorname{rank}\Lambda_\partial$.
\item The quotient $\cB/O(n)$ is a contractible manifold with a locally linear
$\mathbb T_\partial$-action, orbit space $\hat X$, and local model
$(B_y\times\mathbb T_\partial)/I_{\hat y}$, where $h$ acts on the second factor by
translation through $\psi_{\hat y}(h)$.
\item Put $\hat Y_{\max}:=\{\hat y:\dim I_{\hat y}=\operatorname{rank}\Lambda_\partial\}$
and $\hat X_{\max}:=\hat Y_{\max}/O(n)$. The orbit map restricts to a homeomorphism
$(\cB/O(n))^{\mathbb T_\partial}\cong\hat X_{\max}$.
\end{enumerate}
\end{proposition}
\begin{proof}
Since $\hat Y_{hO(n)}\simeq K(\Lambda_\partial,2)=B\mathbb T_\partial$ by \cref{as:lem:borel}, choose a principal
$\mathbb T_\partial$-bundle on it whose boundary map is \eqref{as:eq:boundary}. The torus
classification gives an equivariant bundle $\cB\to\hat Y$ with this Borel construction.
Its homotopy sequence makes $\cB_{hO(n)}$ weakly contractible; \cref{as:lem:borelfree}
then makes the $O(n)$-action free. On each torus fibre write
$hw=w\psi_{\hat y}(h)$. Thus $\ker\psi_{\hat y}$ is a stabilizer in $\cB$, so it is trivial.

For the marking, pull back both this Borel bundle and the Borel construction of the lifted
bundle along
$b_{\hat y}$. They become $EO(n)\times_{I_{\hat y}}\mathbb T_\partial$ and
$EO(n)\times_{I_{\hat y}}\tilde{\mathsf F}_{\hat y}$ over $BI_{\hat y}$.
Their boundary maps are the two orbit homomorphisms. With the same associated-bundle
convention, naturality and the choice of $\partial$ give \eqref{as:eq:orbitmarking}.
Injectivity of $\psi_{\hat y}$ and $(I_y)_0=(I_{\hat y})_0$ give the rank bound.

The equivariant contraction of $B_y$ to $0$ and the same classification trivialize the
restriction of $\cB$ to the slice as $B_y\times\mathbb T_\partial$. Hence the local quotient
is $(B_y\times\mathbb T_\partial)/I_{\hat y}$. Since $\psi_{\hat y}$ is injective, this is a
smooth associated-bundle model with a smooth torus action: identify $I_{\hat y}$ with its
closed image in $\mathbb T_\partial$. It supplies the manifold structure and local linearity.
The homotopy equivalence $\cB/O(n)\simeq\cB_{hO(n)}$ and CW type now give contractibility.
Commuting quotients identify its torus orbit space with $\hat X$. The torus isotropy over
$\hat y$ is exactly $\psi_{\hat y}(I_{\hat y})$, so a point is fixed precisely when
$\dim I_{\hat y}=\operatorname{rank}\Lambda_\partial$. The orbit map is closed and injective
on the fixed set, proving (iv).
\end{proof}

The proof of \cref{as:prop:rank} uses the following property of the points of maximal
isotropy dimension.

\begin{lemma}[Contractible Lifted Fibres]\label{as:lem:contractible}
If $\dim I_{\hat y}=\operatorname{rank}\Lambda_\partial$, then $I_{\hat y}$ is connected and
$\tilde{\mathsf F}_{\hat y}/I_{\hat y}$ is contractible.
\end{lemma}
\begin{proof}
The injection $\psi_{\hat y}$ is onto the torus, hence an isomorphism. Thus
$I_{\hat y}=(I_y)_0$, and \eqref{as:eq:orbitmarking} identifies its orbit map with an
isomorphism onto $\Lambda_\partial$. The principal-bundle sequence makes
$\tilde{\mathsf F}_{\hat y}/I_{\hat y}$ simply connected. By \eqref{as:eq:liftedform}, this
quotient is a component of the inverse image in $\tilde X_i$ of the fibre $\mathsf F_y/I_y$ of
$f^\natural_i$, hence is a covering manifold of that fibre. The latter is aspherical by
\cref{es:lem:quotientfibre}; therefore the quotient is contractible.
\end{proof}

If the kernel vanishes, the torus bundle is $\hat Y$ itself, and the preceding results give
the following.

\begin{corollary}[Trivial Kernel]\label{as:lem:finite}
If $\Lambda_\partial=0$, then $O(n)$ acts freely on $\hat Y$, all $I_y$ are finite, $\hat X$
is a contractible manifold, and $\hat f^\natural_i$ is a fibre bundle with fibre the
universal cover of $\mathsf F_y$.
\end{corollary}
\begin{proof}
The torus $\mathbb T_\partial$ is trivial, so $\cB=\hat Y$ and the assertions follow from
\cref{as:lem:torusbundle,as:lem:liftedbundle}.
\end{proof}

We use two fixed-point and duality facts. Cohomology below is singular cohomology; it agrees
with \v{C}ech cohomology on the locally contractible metrizable spaces in question.

\begin{lemma}[Fixed Points of Tori]\label{as:lem:borelfixed}
If a torus $H$ acts locally linearly on a finite-dimensional, metrizable,
$\mathbb Q$-acyclic manifold $S$, then $S^H$ is non-empty and $\mathbb Q$-acyclic.
\end{lemma}
\begin{proof}
Induct over circle factors of $H$. At each step the preceding fixed set is again a
locally linear manifold. Let a circle act on this $\mathbb Q$-acyclic manifold $S'$, with
fixed set $S''$. Outside $S''$ its stabilizers are finite, so relative rational Borel comparison
\cite{Bor60,Bre72,Hsi75} gives
\begin{equation*}
H^j((S')_{hS^1},S''_{hS^1};\mathbb Q)
\cong H^j(S'/S^1,S'';\mathbb Q)=0\quad\text{for }j>\dim(S'/S^1)+1.
\end{equation*}
The orbit space is locally contractible by linear slices. Since $S'$ is acyclic and
$S^1$ acts trivially on $S''$, the two absolute Borel cohomologies are
$\mathbb Q[u]$ and $H^*(S'';\mathbb Q)\otimes\mathbb Q[u]$, with $\deg u=2$.
The relative vanishing makes their dimensions agree in high degrees. As the cohomology of
$S''$ is bounded in degree, its even Betti numbers sum to one and its odd ones to zero.
Thus $S''$ is non-empty and acyclic, completing the induction without any global finiteness
assumption on orbit types.
\end{proof}

We also use the following consequence of Poincar\'e duality \cite[Chapter~VIII]{Bro82}.
If a Poincar\'e duality group $\Pi$ of dimension $n$ acts freely, properly discontinuously
and cocompactly on a $k$-acyclic manifold $S$, then
\begin{equation}\label{as:eq:PD}
H_c^n(S;k)\cong H^n(\Pi;k\Pi)\cong k,\qquad \dim S\ge n.
\end{equation}
Indeed the Cartan--Leray spectral sequence with the group-ring local system collapses by
acyclicity. This applies to $\pi_1(X_i)$ and its finite-index subgroups, since $X_i$ is a
closed aspherical $n$-manifold.

\begin{proposition}[The Isotropy Rank Is Constant]\label{as:prop:rank}
Every isotropy group $I_y$ has dimension $\operatorname{rank}\Lambda_\partial$.
\end{proposition}
\begin{proof}
The torus $\mathbb T_\partial$ acts locally linearly on the contractible manifold $\cB/O(n)$.
By \cref{as:lem:borelfixed}, its fixed set is non-empty and $\mathbb Q$-acyclic. Thus \cref{as:lem:torusbundle}\,\textup{(iv)} gives the same
conclusion for $\hat X_{\max}$.

For $\hat y\in\hat Y_{\max}$, put $K:=I_{\hat y}$. By \cref{as:lem:contractible}, $K$ is
connected and $\tilde{\mathsf F}_{\hat y}/K$ is contractible. In the tube at $\hat y$,
the maximal-rank locus is $O(n)\times_KB_y^K$: a subgroup of the connected torus $K$ has
full dimension only if it equals $K$. Thus $\hat X_{\max}$ is locally $B_y^K$, and
\eqref{as:eq:liftedform} identifies its preimage
$\tilde X_{\max}:=(\hat f^\natural_i)^{-1}(\hat X_{\max})$ locally with
$(\tilde{\mathsf F}_{\hat y}/K)\times B_y^K$. These are manifolds, and the restricted map
is a bundle with contractible fibres. Hence $\tilde X_{\max}$ is $\mathbb Q$-acyclic.

The rank bound and the slice descriptions make $\hat Y_{\max}$ closed; it is also
$O(n)$- and $\sQ$-invariant. Thus $\tilde X_{\max}$ is closed and deck-invariant, with compact
quotient in $X_i$. By \eqref{as:eq:PD}, its dimension is at least $n$, and therefore equals
$n$. Invariance of domain makes it open in $\tilde X_i$. Connectedness gives
$\tilde X_{\max}=\tilde X_i$, proving the assertion.
\end{proof}

\subsection{Conclusion and the Type-\texorpdfstring{$A_s$}{As} Case}\label{as:sec:proof}

The remaining geometric input is that a non-trivial connected isotropy group of a limit chart
cannot fix an open tube, because an isometry of an $\RCD$ space that fixes a non-empty open set
is the identity. Together with the preceding constant-rank result,
this forces the kernel to vanish.

\begin{lemma}[Generic Finite Isotropy]\label{as:lem:generic}
The points of $Y$ with finite $O(n)$-isotropy are dense.
\end{lemma}
\begin{proof}
Use a limit chart $(Z,z,G)$ and its region $\Omega'$ from \ref{fl:S:limit}, with principal
frame projection and open chart projection $\tau$ from \ref{fl:S:frames}, \ref{fl:S:charts}.
In any non-empty open $U\subset\Omega'$, choose $q$ minimizing $\dim G_q$. By
\cref{ls:thm:slices}, shrink its slice $S$ into $U$. For $s\in S$, the inclusion
$G_s\subset G_q$ and minimality give $(G_s)_0=(G_q)_0=:H$. By
\cref{m:fr:prop:limitgroup}\,\textup{(ii),~(iii)}, $H$ is central in $G_0$, so it fixes
$G_0S$ pointwise. This tube is open, since $G_0G_q$ is open in $G$. By \ref{m:fr:R:isometries}, an isometry
that fixes a non-empty open set pointwise is the identity, so $H=1$. Thus points with finite $G_q$ are
dense in $\Omega'$. The required threshold follows from
$3\delta\le\eps_{\mathrm{glue}}\le\eps_T\le\eps_L(n,N+1)$, as in the proof of \cref{thm:main}.

The inverse image under $\tau$ of a non-empty open set projects to a non-empty open subset
of $\Omega'$. Choose there a point with finite $G_q$ and a frame above it in that inverse
image. Its image in $Y$ has isotropy $A_e(G_q)$ by \ref{fl:S:charts}, hence finite.
The chart images cover $Y$, proving density.
\end{proof}

For fundamental groups of type $A_s$, Smith theory controls the elements of $\sQ$ that have
fixed points.

\begin{lemma}[Finite Subgroups of the Quotient Group]\label{as:lem:As}
Suppose $\Lambda_\partial=0$ and $\pi_1(X_i)$ is of type $A_s$. Every element of $\sQ$ with a
fixed point in $\hat X$ acts trivially. Consequently the action kernel $\sQ_0$ is finite,
$\sQ/\sQ_0$ acts freely, and $X$ is a closed aspherical manifold with universal cover $\hat X$.
\end{lemma}
\begin{proof}
By \cref{as:lem:finite}, $\hat X$ is a contractible manifold and $\hat f^\natural_i$ is a
bundle with contractible fibres. The slice at $\hat y$ is now $B_y$, and its stabilizer in
$\sQ$ is $I_y$, acting through its orthogonal representation by \cref{as:lem:lifted}\,\textup{(ii)}.

If an element of a stabilizer acts non-trivially, a subgroup of its finite cyclic group
has image of prime order $p$ in the homeomorphism group; its Sylow $p$-subgroup $P$ has the
same image. Smith theory \cite{Bre72} makes $\operatorname{Fix}(P)$ non-empty and
$\mathbb F_p$-acyclic. The local linear models make it a manifold. The locus near which
$P$ acts trivially is open and closed in the connected space $\hat X$, hence empty for
this non-trivial action. Thus $\dim\operatorname{Fix}(P)\le\dim\hat X-1$.

The group $\bar\Lambda$ is normal and nilpotent, so type $A_s$ makes $N_{\sQ}(P)$ finite-index
in $\sQ$. Its preimage $\Pi\subset\pi_1(X_i)$ is a Poincar\'e duality group of dimension $n$
and preserves $\tilde X_P:=(\hat f^\natural_i)^{-1}(\operatorname{Fix}(P))$.
The restricted bundle makes $\tilde X_P$ an $\mathbb F_p$-acyclic manifold of dimension at
most $n-1$. The deck action of $\Pi$ on it is free and proper, with quotient a closed
subset of the finite covering $\tilde X_i/\Pi$ of $X_i$, hence compact. This contradicts
\eqref{as:eq:PD}. Therefore every element with a fixed point acts trivially. The kernel
$\sQ_0$ lies in each finite stabilizer, and $\sQ/\sQ_0$ acts freely and cocompactly on the
contractible manifold $\hat X$, with quotient $X$.
\end{proof}

\begin{proof}[Proof of \cref{as:thm:aspherical}]
By \cref{as:prop:rank,as:lem:generic}, the common isotropy dimension
$\operatorname{rank}\Lambda_\partial$ is zero. Thus $\Lambda_\partial=0$, and
\cref{as:lem:finite} gives finite isotropy and a bundle
$\hat f^\natural_i\colon\tilde X_i\to\hat X$ over a contractible manifold, whose fibre is
the universal cover of a nilmanifold, hence a Euclidean space.

(i) The kernel definition gives injectivity of $\Lambda\to\pi_1(X_i)$. The inverse image in
$\tilde X_i$ of the fibre $\mathsf F_y/I_y$ of $f^\natural_i$ consists of the fibres of
$\hat f^\natural_i$ over a discrete $\sQ$-orbit, each contractible. Therefore the fundamental
group of this fibre also injects into $\pi_1(X_i)$.

(ii) Apply \cref{as:lem:lifted}\,\textup{(iii)} and the preceding bundle description.

(iii) Apply \cref{as:lem:As}.
\end{proof}

\appendix
\clearpage
\section{Dependency and Parameter Guide for the Main Text}
\label[appendix]{guide:sec:dependencies}

This guide concerns \crefrange{ch:intro}{as:sec:aspherical} only. Arrows indicate
principal uses of results, with their stated hypotheses retained; they do not
replace the proofs or list every auxiliary input. Grey boxes import a result
from another part of the argument, and blue boxes mark conclusions. All result
numbers and defining locations are linked.

\begingroup
\newcommand{\DGref}[2]{\hyperref[#1]{#2~\ref*{#1}}}
\newcommand{\DGsec}[1]{\hyperref[#1]{\S\ref*{#1}}}
\newcommand{\DGpage}[2]{\hyperlink{page.\getpagerefnumber{#1}}{#2~\ref*{#1}}}
\tikzset{
  dgbox/.style={draw=black!55,rounded corners=2pt,fill=white,
    align=center,text width=6.3cm,inner sep=5pt,font=\small},
  dgwide/.style={dgbox,text width=13.9cm},
  dginput/.style={dgbox,fill=black!6},
  dgout/.style={dgbox,fill=MidnightBlue!7,draw=MidnightBlue!65},
  dgarrow/.style={-{Latex[length=2mm]},draw=black!65,line width=.5pt}
}

\subsection{Frames, Covers and the Framed Limit}\label{guide:map:geometry}

\begin{center}
\begin{tikzpicture}
\node[dgwide] (hyp) at (0,0)
  {\DGref{hyp:sequence}{Hypothesis}, with the dimensional smallness
  choices in \DGsec{guide:sec:choices}};
\node[dgbox] (frames) at (-3.8,-1.65)
  {Canonical frames, gauges and convergence\\
   \cref{m:fr:thm:bundle,m:fr:thm:gauges,m:fr:thm:frameconv}};
\node[dgbox] (covers) at (3.8,-1.65)
  {Completion, loop contraction and cover agreement\\
   \cref{m:fr:prop:completion,m:fr:thm:localBC,m:fr:lem:agreement,m:fr:prop:coverscharts}};
\node[dgbox] (free) at (-3.8,-3.7)
  {Localization and free frame actions\\
   \DGref{m:fr:prop:localization}{Proposition} and
   \DGpage{m:fr:prop:freeframes}{Proposition}};
\node[dgbox] (atlas) at (3.8,-3.7)
  {Compatible geometric atlas and framed convergence\\
   \DGref{m:fr:thm:geomain}{Theorem}};
\node[dgbox] (slices) at (-3.8,-5.4)
  {Equivariant harmonic charts and linear slices\\
   \cref{ls:prop:equivariant,ls:thm:slices}};
\node[dgbox] (limit) at (3.8,-5.4)
  {The framed limit is a manifold\\
   \DGref{ls:thm:limits}{Theorem}};
\node[dgout] (main) at (-3.8,-7.1)
  {Frame resolution\\
   \DGref{thm:main}{Theorem}\textup{(i)}};
\node[dgbox] (realize) at (3.8,-7.1)
  {Equivariant realizing maps\\
   \DGref{fl:thm:anr}{Theorem}};
\draw[dgarrow] (hyp.south) -- ++(0,-.25) -| (covers.north);
\draw[dgarrow] (frames) -- (free);
\draw[dgarrow] (frames.west) -- (-7.45,-1.65) |- (slices.west);
\draw[dgarrow] (covers) -- (atlas);
\draw[dgarrow] (free) -- (atlas);
\draw[dgarrow] (atlas) -- (limit);
\draw[dgarrow] (slices) -- (limit);
\draw[dgarrow] (limit.south) -- ++(0,-.4) -| (main.north);
\draw[dgarrow] (limit) -- (realize);
\end{tikzpicture}
\end{center}

The frame theorems do not use \cref{hyp:sequence}: the frame bundle exists on every metric
measure space, and gauges and the convergence of frames need only an $\RCD(\kappa,N)$ space
with a flat region at the frame scale (\cref{m:book:ch:07,m:book:ch:08}). The hypothesis
enters through the local covers. The atlas also supplies the compact realization used in the
frame-resolution conclusion (\cref{m:fr:sec:coverbundle,ls:sec:limits}). Realizing maps are the next
input to the bundle construction, not an additional hypothesis for part~\textup{(i)}.
The atlas and slice results feed the next two diagrams without being expanded again.

\medskip
\noindent\textbf{An independent application.}
The local bounded-covering theorem has its own hypotheses, not those of
\cref{hyp:sequence}, and does not use \cref{thm:main}.

\begin{center}
\begin{tikzpicture}[dgbox/.append style={text width=4.25cm}]
\node[dgbox] (a) at (-4.95,0)
  {Local bounded covering geometry\\
   \DGref{hz:def:lbcg}{Definition}\\
   Completion: \DGref{m:fr:prop:completion}{Proposition}};
\node[dgbox] (b) at (0,0)
  {Uniform covering regularity\\
   \DGref{hz:lem:scale}{Corollary}\\
   Smoothing and frames: \cite[Theorems~5.3, 1.7]{HZ26}};
\node[dgout,text width=4.25cm] (c) at (4.95,0)
  {Fibration under local bounded covering geometry\\
   \DGref{hz:thm:fibration}{Theorem}};
\draw[dgarrow] (a) -- (b);
\draw[dgarrow] (b) -- (c);
\end{tikzpicture}
\end{center}

\clearpage
\subsection{From Deck Labels to the Actual Bundle}\label{guide:map:bundle}

The first correction is algebraic: approximate labels become an exact cocycle.
The second step realizes that cocycle by sections of the actual frame spaces.
Only after this realization do the local modification and gluing arguments
produce the bundle.

\begin{center}
\begin{tikzpicture}
\node[dginput,text width=13.9cm] (input) at (0,0)
  {Atlas \ref{m:fr:thm:geomain}; realizing maps \ref{fl:thm:anr};
   quantitative slices and recentring
   \ref{th:lem:almostequiv}, \ref{th:lem:meets}, \ref{th:lem:recentre}};
\node[dgbox] (labels) at (-3.8,-1.65)
  {Transplanted charts and labels\\
   \DGref{th:lem:transplant}{Lemma} and \DGpage{es:lem:labelchart}{Lemma}};
\node[dgbox] (weights) at (3.8,-1.65)
  {Normalized deck weights\\
   \DGpage{th:lem:weights}{Lemma}};
\node[dgbox] (cocycle) at (0,-3.35)
  {Exact cocycle\\ \DGref{es:prop:cocycle}{Proposition}};
\node[dgbox] (coord) at (0,-4.95)
  {Equivariant coordinate and nearby-point model\\
   \cref{es:prop:coordinate,es:thm:nearbyprecise}};
\node[dgbox] (sections) at (0,-6.65)
  {Admissible sections of the actual frame spaces\\
   \DGref{es:lem:admissible}{Lemma}};
\node[dgbox] (modify) at (-3.8,-8.45)
  {Modification and transfer between tubes\\
   \cref{es:lem:modall,es:lem:twochartstransition}};
\node[dgbox] (system) at (3.8,-8.45)
  {Uniform data and a finite tube system\\
   \cref{es:lem:uniform,es:lem:system}};
\node[dgbox] (glue) at (-3.8,-10.25)
  {Global gluing\\ \DGref{es:thm:gluingall}{Theorem}};
\node[dginput] (trivial) at (3.8,-10.25)
  {Equivariant triviality over a ball\\
   \DGref{es:prop:trivial}{Proposition}};
\node[dgout] (result) at (0,-12.05)
  {Bundle and local quotient form\\
   \DGref{thm:main}{Theorem}\textup{(ii)--(iii)}};
\draw[dgarrow] (input.south) -- ++(0,-.25) -| (labels.north);
\draw[dgarrow] (input.south) -- ++(0,-.25) -| (weights.north);
\draw[dgarrow] (labels.south) -- ++(0,-.25) -| (cocycle.north);
\draw[dgarrow] (weights.south) -- ++(0,-.25) -| (cocycle.north);
\draw[dgarrow] (cocycle) -- (coord);
\draw[dgarrow] (coord) -- (sections);
\draw[dgarrow] (sections.south) -- ++(0,-.25) -| (modify.north);
\draw[dgarrow] (coord.east) -- (7.45,-4.95) |- (system.east);
\draw[dgarrow] (modify) -- (system);
\draw[dgarrow] (modify) -- (glue);
\draw[dgarrow] (system.south) -- ++(0,-.25) -| (glue.north);
\draw[dgarrow] (glue.south) -- ++(0,-.25) -| (result.north);
\draw[dgarrow] (trivial) -- (glue);
\end{tikzpicture}
\end{center}

The diagram ends with closed manifold fibres, not yet their nilmanifold
identification. The global modification preserves the regularity of the sections;
the finite overlap bound controls its accumulated loss. The independent lattice
argument in \cref{ch:lattices} is needed to determine the fibre topology.

\clearpage
\subsection{Fibre Topology and the Aspherical Application}\label{guide:map:topology}

The lattice branch starts from the atlas, not from the global bundle.
It meets the gluing branch only when the topology of the actual fibres is identified.

\begin{center}
\begin{tikzpicture}
\node[dginput] (atlas) at (-3.8,0)
  {Compatible atlas\\ \DGref{m:fr:thm:geomain}{Theorem}};
\node[dgbox] (lattice) at (3.8,0)
  {Marked lattices and period kernels\\
   \cref{m:fr:ext:physical-lattice-realization,m:fr:ext:period-convergence}};
\node[dgbox] (cover) at (-3.8,-1.85)
  {Unwrapped coverings and monodromy\\
   \cref{fl:prop:covering,fl:prop:monodromy}};
\node[dgbox] (carriers) at (3.8,-1.85)
  {Towers, carriers and liftable spheres\\
   \cref{fl:lem:towers,fl:lem:carriers,fl:lem:liftablespheres}};
\node[dginput] (glue) at (-3.8,-3.75)
  {The actual equivariant bundle\\
   \DGref{es:thm:gluingall}{Theorem}};
\node[dgbox] (groups) at (3.8,-3.75)
  {Intrinsic $\pi_1$ images; higher homotopy dies on enlargement\\
   \cref{fl:thm:fibre,fl:prop:preimagegroups,fl:thm:higher}};
\node[dgbox] (fibres) at (0,-5.65)
  {Closed fibres: asphericity and rigidity\\
   \cref{es:thm:fibrenil,es:cor:fibresnil}};
\node[dgout,text width=13.9cm] (quotient) at (0,-7.35)
  {Isotropy quotients: \cref{es:lem:quotientfibre,es:thm:infranil}\\
   Nilmanifold frame fibres and infranilmanifold quotient fibres:
   \DGref{thm:main}{Theorem}\textup{(iv)}};
\draw[dgarrow] (atlas) -- (lattice);
\draw[dgarrow] (lattice.south) -- ++(0,-.25) -| (cover.north);
\draw[dgarrow] (cover) -- (carriers);
\draw[dgarrow] (carriers) -- (groups);
\draw[dgarrow] (glue.south) -- ++(0,-.25) -| (fibres.north);
\draw[dgarrow] (groups.south) -- ++(0,-.25) -| (fibres.north);
\draw[dgarrow] (fibres) -- (quotient);
\end{tikzpicture}
\end{center}

The topology of nested pieces is not itself a bundle theorem.
In \cref{es:thm:fibrenil}, product neighbourhoods from the already constructed
bundle turn the enlargement maps into homotopy equivalences, identifying
the homotopy type of a closed fibre.

\medskip
\noindent\textbf{The aspherical application.}
Here the additional hypothesis is that $X_i$ is aspherical.
\begin{center}
\begin{tikzpicture}
\node[dginput] (main) at (-3.8,0)
  {\DGref{thm:main}{Theorem} and aspherical $X_i$\\
   Borel lifting and the torus of the kernel\\
   \cref{as:sec:lifted,as:lem:torusbundle}};
\node[dgbox] (rank) at (3.8,0)
  {Torus fixed points and duality give constant rank;
   generic isotropy is finite\\
   \cref{as:lem:borelfixed,as:prop:rank,as:lem:generic}};
\node[dgout] (finite) at (-3.8,-2.05)
  {$\Lambda_\partial=0$; \DGref{as:lem:finite}{Corollary}\\
   \DGref{as:thm:aspherical}{Theorem}\textup{(i)--(ii)}};
\node[dgout] (As) at (3.8,-2.05)
  {With type $A_s$: \DGref{as:lem:As}{Lemma}\\
   \DGref{as:thm:aspherical}{Theorem}\textup{(iii)}};
\draw[dgarrow] (main) -- (rank);
\draw[dgarrow] (rank.south) -- ++(0,-.25) -| (finite.north);
\draw[dgarrow] (finite) -- (As);
\end{tikzpicture}
\end{center}
The generic-isotropy step also uses the earlier limit-group and linear-slice
results; it does not assume the conclusion that all isotropy is finite.

\clearpage
\subsection{Order of the Choices}\label{guide:sec:choices}

The universal constants and the local construction scales are chosen in
different orders. The arrows below record that order, not new hypotheses.

\begin{center}
\begin{tikzpicture}
\node[dgbox] (flat) at (-3.8,0)
  {Atlas accuracy $\eps_\Omega$\\
   Flatness, free frame actions and gauges\\
   \DGsec{m:fr:sec:coverbundle}};
\node[dgbox] (budget) at (3.8,0)
  {Dimensional overlap budget\\
   $M_n=539^n$,\quad
   $\eps_*=[240((n+1)M_n+1)]^{-1}$\\
   \eqref{es:eq:kappasharp}};
\node[dgbox] (atlas) at (-3.8,-2.05)
  {$\delta_*=\min\{\delta_C(n,N,\eps_\Omega),\delta_P(n)/3\}$\\
   $c_*=c_C(n,N,\eps_\Omega)$\\
   \DGsec{m:fr:sec:coverbundle}};
\node[dgbox] (comparison) at (3.8,-2.05)
  {Comparison and slice accuracies\\
   $\eps_\flat,\theta_\flat,c_\flat$; then $\eps'_S,c'_S$\\
   \cref{th:sec:constants,th:sec:uniform}};
\node[dgbox] (glue) at (3.8,-4.05)
  {Final gluing accuracy $\eps_{\mathrm{glue}}(n,N)$\\
   Completed charts use the dimension bound $N+1$\\
   \eqref{th:eq:threshold}, \eqref{es:eq:thresholdG}};
\node[dgout,text width=13.9cm] (main) at (0,-6.05)
  {$\delta_1=\min\{\delta_*,\eps_{\mathrm{glue}}/3\}$,\qquad $c_1=c_*$\\
   Thresholds for \DGref{thm:main}{Theorem}, depending only on $n,N$};
\draw[dgarrow] (flat) -- (atlas);
\draw[dgarrow] (budget) -- (comparison);
\draw[dgarrow] (comparison) -- (glue);
\draw[dgarrow] (atlas.south) -- (atlas.south |- main.north);
\draw[dgarrow] (glue.south) -- (glue.south |- main.north);
\end{tikzpicture}
\end{center}

The overlap estimate is proved later in \cref{es:lem:system}, but its explicit
dimensional value fixes the accuracy budget before the analytic construction:
\[
  60\eps_*(n+1)M_n
  =\frac{(n+1)M_n}{4((n+1)M_n+1)}<\frac14.
\]
Thus the forward reference does not require choosing $\eps_*$ after gluing.
The factor $3$ in $\delta_1$ provides the limit-chart error required by
\eqref{es:eq:thresholdG}; the use of $N+1$ comes from the completed charts,
not from a change in the dimension $n$ of the frames.

\medskip
\noindent\textbf{Local choices after fixing the geometric setting.}
\begin{center}
\begin{tikzpicture}
\node[dgbox] (global) at (-3.8,0)
  {Fix $r\le c_*\rho$, the finite atlas and its subsequence\\
   Choose a common harmonic scale $s_0$\\
   \DGsec{es:sec:finitesetting}};
\node[dgbox] (local) at (3.8,0)
  {Fix the local orbit and isotropy data\\
   Choose $\ell_\alpha\le\ell_*(z_\alpha)$\\
   \eqref{th:eq:ellsmall}, \eqref{es:eq:ellone}};
\node[dgwide] (index) at (0,-1.95)
  {Retain the chosen tubes, then take $i$ sufficiently large\\
   For finitely many tubes, use a common subsequence and index bound};
\draw[dgarrow] (global) -- (local);
\draw[dgarrow] (local.south) -- ++(0,-.25) -| (index.north);
\end{tikzpicture}
\end{center}

The tube radii need not be equal. \Cref{th:lem:recentre} controls recentring
near points with the same isotropy group $K$; \cref{es:lem:uniform} supplies
a common tube scale on each compact set on which $k_0$ is constant, where $k_0$ is the
dimension of the fixed space $V_0$ of a tube.
Neither statement makes every orbit-dependent radius a dimensional constant.
The finitely many subsequent choices are recorded in \cref{es:lem:system}.

\clearpage
\subsection{Persistent Parameters and Notation}\label{guide:sec:register}

Only notation passed between results is included. ``Dimensional'' means
depending only on $n,N$ after the prescribed choices. Generic constants,
dummy variables and proof-local cut-offs are omitted; this table does not
alter their dependencies in the proofs.

\medskip
\noindent\textbf{Hypotheses, accuracies and scales.}

\small
\setlength{\tabcolsep}{4pt}
\renewcommand{\arraystretch}{1.16}
\begin{longtable}{@{}>{\raggedright\arraybackslash}p{.17\linewidth}
  >{\raggedright\arraybackslash}p{.25\linewidth}
  >{\raggedright\arraybackslash}p{.39\linewidth}
  >{\raggedright\arraybackslash}p{.12\linewidth}@{}}
\toprule
Symbol & Role & Scope and dependencies & Location\\
\midrule
$n,N,\kappa$ &
Frame dimension, curvature-dimension bound, curvature lower bound &
Fixed sequence data; $n$ need not equal $N$. &
\DGref{hyp:sequence}{Hyp.}\\
\addlinespace
$\delta,\rho$ &
Reifenberg error and covering radius &
Fixed throughout the sequence. Comparison scales are at most $\rho/3$;
$|\kappa|\rho^2\le1$ is assumed. &
\DGref{m:fr:def:localcovers}{Def.},
\DGref{hyp:sequence}{Hyp.}\\
\addlinespace
$r$ &
Frame-operator scale &
One fixed positive scale, at most $c_*\rho$ in the main construction.
It is not the harmonic or tube scale. &
\DGref{m:fr:def:operator}{Def.},
\DGsec{m:fr:sec:coverbundle}\\
\addlinespace
$R_B,R_0,R_1$ &
Agreement, chart and transition radii &
Fixed multiples of $\rho$:
$R_B=L_B\rho/10$, $R_0=L_B\rho/200$, $R_1=R_0/32$.
The frame choice ensures $128r\le R_1$. &
\DGsec{m:fr:sec:localcovers}\\
\addlinespace
$\eps_\Omega$;\newline $\delta_*,c_*$ &
Atlas accuracy and covering constants &
Dimensional; $\eps_\Omega$ includes the flatness, freeness and gauge
thresholds. Then $\delta\le\delta_*$ and $r\le c_*\rho$. &
\DGsec{m:fr:sec:coverbundle}\\
\addlinespace
$M_n,\eps_*$ &
Overlap and Lipschitz-error budget &
Fixed before the comparison estimates. At each point, at most $(n+1)M_n$
modifications, each of cost $60\eps_*$, have total cost below $1/4$. &
\eqref{es:eq:kappasharp},
\DGref{es:lem:system}{Lem.}\\
\addlinespace
$\eps_\flat,\theta_\flat,c_\flat$;\newline
$\eps'_S,c'_S$ &
Comparison and recentred slice thresholds &
Dimensional once $\eps_*$ is fixed. Use the completed-chart dimension
bound $N+1$ when applying them to the atlas. &
\DGsec{th:sec:constants},
\DGsec{th:sec:uniform},
\eqref{th:eq:threshold}\\
\addlinespace
$\eps_{\mathrm{glue}}$;\newline $\delta_1,c_1$ &
Final gluing and main-theorem thresholds &
Dimensional; $3\delta\le\eps_{\mathrm{glue}}$,
$\delta_1=\min\{\delta_*,\eps_{\mathrm{glue}}/3\}$, $c_1=c_*$. &
\DGsec{es:sec:finitesetting},
\DGref{thm:main}{Thm.}\\
\addlinespace
$s_0$ &
Harmonic-chart scale &
First used locally; chosen common across the finite atlas for gluing.
Its admissibility includes the curvature and neighbourhood buffers. &
\DGsec{th:sec:standing},
\DGsec{es:sec:finitesetting}\\
\addlinespace
$\ell_*,\ell_\alpha$;\newline $r_\ell=\ell/800$ &
Upper bound and chosen value of the tube radius &
Depend on local orbit/isotropy data and the earlier choices.
Fix $\ell_\alpha$ before taking $i$ large; it may differ between tubes. &
\DGsec{th:sec:scale},
\DGsec{es:sec:finitesetting}\\
\addlinespace
$\ell_{\mathcal S}$ &
Uniform radius on a stratum piece &
Uniform on a specified compact set on which $k_0$ is constant, not on
arbitrary unions of strata. &
\DGref{es:lem:uniform}{Lem.}\\
\bottomrule
\end{longtable}

\clearpage
\normalsize
\noindent\textbf{Objects carried between the geometric and topological arguments.}

\smallskip
The scale $r$ is suppressed in the frame notation after its definition.
Frames over a subset still use the ambient metric and measure: restricting
the base does not mean rebuilding the frame operator on the subset.

\small
\begin{longtable}{@{}>{\raggedright\arraybackslash}p{.17\linewidth}
  >{\raggedright\arraybackslash}p{.25\linewidth}
  >{\raggedright\arraybackslash}p{.39\linewidth}
  >{\raggedright\arraybackslash}p{.12\linewidth}@{}}
\toprule
Symbol & Role & Scope and distinctions & Location\\
\midrule
$\mu_x,\cK_x$;\newline $E_x,u_e$ &
Localized measure, frame operator, spectral space and frame map &
All belong to the same fixed frame scale $r$; the scale is not
changed by later tube shrinkings. &
\DGref{m:fr:def:operator}{Def.},
\DGref{m:fr:def:frames}{Def.}\\
\addlinespace
$Y_i,Y$ &
Frame spaces and framed limit &
$Y_i=\Fr_\rho(X_i)$ at the fixed scale $r$; the subscript $\rho$ is the
covering radius, not the frame scale.
$Y$ is the subsequential compact framed limit, with $Y/O(n)=X$. &
\DGref{m:fr:thm:geomain}{Thm.},
\DGref{ls:thm:limits}{Thm.}\\
\addlinespace
$Z_i,\Gamma_i$;\newline $Z,G$ &
Completed charts with deck groups, and their equivariant limits &
Chart-dependent; a superscript $a$ specifies the atlas centre.
The geometric data are those of the common atlas. &
\DGref{m:fr:thm:geomain}{Thm.},
\DGsec{fl:sec:setting}\\
\addlinespace
$K$;\newline $I$ &
Isotropy in the chart and in the framed limit &
$K$ is the isotropy group of $z$ in $G$, and $I$ that of $y$ in $O(n)$;
$k\mapsto A_k$ maps $K$ isomorphically onto $I$. &
\DGsec{th:sec:standing}\\
\addlinespace
$h,v_i,h_i$ &
Limit harmonic chart, harmonic charts on $Z_i$, transplanted charts &
$h_i=\hat h\circ v_i$ converges to $h$; a transplanted chart is not
asserted to be harmonic. &
\DGref{th:lem:stagecharts}{Lem.},
\DGref{th:def:transplant}{Def.}\\
\addlinespace
$k_\gamma,\varkappa_i$;\newline $\xi,\tilde\xi_i$ &
Labels, exact cocycle, limit and equivariant coordinates &
The labels are corrected to an exact cocycle before the equivariant
coordinate is built; limit objects and objects on $Z_i$ stay distinct. &
\DGref{es:lem:labels}{Lem.},
\DGref{es:prop:cocycle}{Prop.},
\DGref{es:prop:coordinate}{Prop.},
\DGref{es:lem:limit}{Lem.}\\
\addlinespace
$F_i,F_i^\natural$;\newline $f_i,f_i^\natural$ &
Realizing maps, bundle maps and their quotients &
$F_i$ is a continuous equivariant approximation.
$F_i^\natural$ is the nearby actual bundle; $f_i,f_i^\natural$
are the induced maps on $X_i$. &
\DGref{fl:def:realizing}{Def.},
\DGref{es:thm:gluingall}{Thm.}\\
\addlinespace
$\Lambda_i,\vartheta_i$;\newline $\ker\vartheta_i$ &
Fibre lattice, its marking and the period kernel &
The marking need not be injective. Unwrapping removes the periods;
it does not identify the lattice with its deck image. &
\DGref{m:fr:ext:physical-lattice-realization}{Thm.},
\DGref{m:fr:ext:period-convergence}{Thm.}\\
\addlinespace
$A_i(\mathbf b)$;\newline $V_i^{[j]}$ &
Sets of small lattice elements and the sets of a tower &
The box bounds control displacement; the tower length is set by the
dimension of the complex being filled, not the number of sheets visited. &
\DGsec{fl:sec:latticecomplex},
\DGpage{fl:def:tower}{Def.},
\DGref{fl:lem:carriers}{Lem.}\\
\addlinespace
$\Lambda_\partial$;\newline $\mathbb T_\partial$ &
The kernel and its torus &
Only in the aspherical application:
$\Lambda_\partial=\ker(\Lambda\to\pi_1(X_i))$.
This is not the period kernel $\ker\vartheta_i$. &
\DGref{as:lem:borel}{Lem.},
\DGref{as:lem:torusbundle}{Prop.}\\
\bottomrule
\end{longtable}
\endgroup

\section{Ricci-Flow Regularization under Local Rewinding Control}\label[appendix]{rf:app:flow}

A localized comparison with a complete comparison flow transfers curvature control to the
possibly incomplete Ricci flow lifted to a local cover; it gives a uniform lifetime on closed
manifolds and, by exhaustion, a complete flow even when the initial curvature is unbounded.

\subsection{Statements and Conventions}\label{rf:sec:results}

For the Ricci-flow argument we use the following all-scale smooth version of the
rewinding condition.

\begin{definition}[Local Rewinding Conditions]\label{rf:def:lrrp}
Let $(M^n,g)$ be a Riemannian manifold, $x\in M$ and $\rho>0$, and let
\begin{equation*}
\pi_x:\bigl(\tB_g(x,\rho),\tilde x\bigr)\longrightarrow\bigl(B_g(x,\rho),x\bigr)
\end{equation*}
be the Riemannian universal cover of the open ball. The point $x$ is a
\emph{$(\delta,\rho)$-local rewinding Reifenberg point}, abbreviated $(\delta,\rho)$-LRRP,
if
\begin{equation}\label{rf:eq:lrrp}
d_{GH}\bigl(B_{\pi_x^*g}(\tilde x,s),B_{\R^n}(0,s)\bigr)<\delta s
\qquad(0<s<\rho),
\end{equation}
for the pointed Gromov--Hausdorff distance, the centres of the balls being the base
points. The manifold has \emph{uniform $(\delta,\rho)$-LRRP} if every point
is a $(\delta,\rho)$-LRRP. The point $x$ is a
\emph{$(\rho,v)$-local rewinding volume point}, abbreviated $(\rho,v)$-LRVP, if
\begin{equation*}
\vol_{\pi_x^*g}B_{\pi_x^*g}(\tilde x,\rho)\ge v.
\end{equation*}
\end{definition}

The terminology follows \cite{HRW20,Ron22,Ron26}. The numerical constants of these appendices,
in particular the Reifenberg bounds of \cref{rf:sec:example}, use this pointed convention
\textup{(}\cref{rf:sec:conv}\textup{)}, not the approximation-map convention of
\cref{m:fr:def:localcovers}.

\begin{theorem}[Uniform Lifetime on Closed Manifolds]\label{rf:thm:main}
Fix $n\ge2$, $\rho>0$ and $\alpha\in(0,1)$. There are
$\delta_0=\delta_0(n,\alpha)>0$ and $T_0=T_0(n,\rho,\alpha)>0$ with the following
property. If $(M^n,g_0)$ is closed, $\Ric_{g_0}\ge-(n-1)g_0$, and
every point of $M$ is a $(\delta_0,\rho)$-LRRP, then the maximal Ricci flow $g(t)$
with $g(0)=g_0$ satisfies
\begin{equation}\label{rf:eq:main-bound}
T_{\max}>T_0,\qquad
\sup_M|\Rm_{g(t)}|\le\frac{\alpha}{t}\quad(0<t\le T_0).
\end{equation}
\end{theorem}

The second result allows complete non-compact data with unbounded initial curvature.

\begin{theorem}[Complete Manifolds]\label{rf:thm:noncompact}
Fix $n\ge2$, $\rho>0$ and $\alpha\in(0,1)$. There are
constants $\delta_{\mathrm{nc}}\in(0,\delta_0]$ and
$T_{\mathrm{nc}}\in(0,T_0]$, depending only on $(n,\alpha)$ and
$(n,\rho,\alpha)$, respectively, with the following
property. If $(M^n,g_0)$ is complete, $\Ric_{g_0}\ge-(n-1)g_0$, and
every point of $M$ is a $(\delta_{\mathrm{nc}},\rho)$-LRRP, then there is a
complete Ricci flow $g(t)$, smooth on $M\times[0,T_{\mathrm{nc}}]$, with $g(0)=g_0$ and
\begin{equation*}
\sup_M|\Rm_{g(t)}|\le\frac{\alpha}{t}\quad(0<t\le T_{\mathrm{nc}}).
\end{equation*}
It is the only complete Ricci flow on $M\times[0,T_{\mathrm{nc}}]$ with initial metric
$g_0$ whose curvature satisfies $\sup_M|\Rm_{g(t)}|\le A/t$ for some constant $A$.
\end{theorem}

The uniqueness assertion recovers \cite[Theorem~1.1]{Lee26} by the same comparison
argument. The smooth affine fibration consequences are given in
\cref{rf:sec:rong}.

\subsubsection{Conventions}\label{rf:sec:conv}
Manifolds are smooth and connected, and ``closed'' means compact without boundary.
Write $B_g(x,r)$ and $d_g$ for Riemannian balls and distance, and $B_Z(z,r)$ in a metric
space. A regular domain is relatively compact, open, and smoothly bounded; $\omega_n$
is the volume of the Euclidean unit ball. Ricci flows satisfy
$\partial_tg=-2\Ric_g$. Constants display their dependencies when needed; unadorned
$C,c$ depend only on the parameters of the current statement and may change by line.

A correspondence between metric spaces $(Z,d)$ and $(Z',d')$ is a subset
$\mathcal R\subset Z\times Z'$ whose two projections are onto, and its distortion is
\begin{equation*}
\dis\mathcal R:=\sup\bigl\{|d(z_1,z_2)-d'(z_1',z_2')|:
(z_1,z_1'),(z_2,z_2')\in\mathcal R\bigr\}.
\end{equation*}
The Gromov--Hausdorff distance is $d_{GH}(Z,Z')=\tfrac12\inf_{\mathcal R}\dis\mathcal R$
\cite[Theorem~7.3.25]{BBI01}. For pointed spaces we use
\begin{equation*}
d_{GH}\bigl((Z,z),(Z',z')\bigr):=\tfrac12\inf\{\dis\mathcal R:(z,z')\in\mathcal R\}.
\end{equation*}
This is at least the unpointed distance. For balls with displayed centres, including
\eqref{rf:eq:lrrp} and \eqref{rf:eq:hyp-gh}, $d_{GH}$ is pointed. If every point of $Z$ lies
within $s$ of $z$ and a correspondence between $Z$ and $B_{\R^n}(0,s)$ has distortion
at most $\eps$, every point related to $z$ has norm at most $\eps$; adding $(z,0)$
therefore at most doubles the distortion. Thus Rong's unpointed form of
\eqref{rf:eq:lrrp} implies ours with $2\delta$ in place of $\delta$.

\subsection{Local Analytic Inputs}\label{rf:sec:prelim}

Local completeness on the local covers and the estimates below supply the analytic inputs.

\subsubsection{Local Universal Covers}\label{rf:sec:covers}
The local covers are incomplete. For $\tilde z\in\tB_g(x,\rho)$ call
$\rho_{\tilde z}:=\rho-d_{\pi_x^*g}(\tilde x,\tilde z)$ its depth. Positive depth supplies
the following local substitute for completeness.

\begin{lemma}[Local Completeness of the Cover]\label{rf:lem:cover}
Let $(M^n,g)$ be a complete Riemannian manifold, $x\in M$ and $\rho>0$, and let
$\tilde z\in\tB_g(x,\rho)$ have depth $\rho_{\tilde z}>0$.
\begin{enumerate}[label=\textup{(\roman*)}]
\item The exponential map of $\pi_x^*g$ at $\tilde z$ is defined on the open ball of
      radius $\rho_{\tilde z}$ in the tangent space.
\item For $0<r<\rho_{\tilde z}$, the closed ball of radius $r$ about $\tilde z$ is
      compact, and each of its points is joined to $\tilde z$ by a minimizing geodesic.
\item If $\Ric_g\ge0$ and $0<r_1<r_2<\rho_{\tilde z}$, then
      $\vol B_{\pi_x^*g}(\tilde z,r_2)\le(r_2/r_1)^n\vol B_{\pi_x^*g}(\tilde z,r_1)$.
\end{enumerate}
\end{lemma}
\begin{proof}
For $z=\pi_x(\tilde z)$,
$d_g(x,z)\le d_{\pi_x^*g}(\tilde x,\tilde z)$. Thus every $g$-geodesic from $z$ of
length below $\rho_{\tilde z}$ stays in $B_g(x,\rho)$ and lifts, proving (i). The local
Hopf--Rinow proof \cite[Chapter~7, Theorem~2.8]{dC92} gives (ii), and the standard
Jacobian comparison along these minimizing geodesics gives (iii).
\end{proof}

\subsubsection{Local Entropy and Pseudolocality}\label{rf:sec:entropy}
For a metric $h$ on an $n$-manifold, a relatively compact ball $B=B_h(x,r)$ and $s>0$,
the local entropy $\nu(B,h,s)$ of \cite{Wan18} is the infimum of
\begin{equation*}
-n-\frac n2\log(4\pi s')
+\int_B\Bigl(s'\bigl(\Scal_h\psi^2+4|\nabla\psi|^2\bigr)-2\psi^2\log\psi\Bigr)\,d\vol_h
\end{equation*}
over $s'\in(0,s]$ and $\psi\in W^{1,2}_0(B)$ with $\psi\ge0$ and
$\int_B\psi^2\,d\vol_h=1$. Wang's pseudolocality theorem turns a lower bound for this
entropy into a scale-invariant curvature bound; \cref{rf:sec:selected} uses it to give the
comparison flow its Type-I bound.

\begin{theorem}[Wang's Pseudolocality Theorem {\cite[Theorem~1.2]{Wan20}}]\label{rf:thm:wang}
For every positive $\beta$ smaller than $1/(100n)$ there is
$\delta_W=\delta_W(n,\beta)>0$ with the
following property. Let $h(t)$, $0\le t\le T$, be a complete Ricci flow with bounded
curvature on an $n$-manifold, and let $x$ be a point of the manifold. If
\begin{equation*}
\nu\bigl(B_{h(0)}(x,\delta_W^{-1}\sqrt T),h(0),T\bigr)\ge-\delta_W^2,
\end{equation*}
then $t\,|\Rm_{h(t)}|\le\beta$ on $B_{h(t)}(x,\beta^{-1}\sqrt t)$ for every $t\in(0,T]$.
\end{theorem}

\subsubsection{Ricci--Harmonic Comparison Inputs}\label{rf:sec:lee-inputs}

We use the Ricci--harmonic gauge method of \cite[Sections~2--3]{Lee26}.
The following proposition states the local versions needed here;
\cref{rf:app:comparison-inputs} supplies the existence and derivative proofs.
Derivatives of a solution
$\partial_tF=\Delta_{g(t),h(t)}F$ use $g(t)$ and $h(t)$; when $F_t$ is a diffeomorphism
onto its image, $\widehat g=(F_t^{-1})^*g$ is Ricci--DeTurck relative to $h$.

\begin{proposition}[Local Comparison Estimates]\label{rf:prop:lee-inputs}
All source and target manifolds have dimension $n$. Let
$\eps_{\mathrm{DT}}(n)\in(0,10^{-2}]$ be a threshold below which the discrepancy
inequality \eqref{rf:eq:discrepancy-decay-pde} below holds, and set
$\eps_{\mathrm L}(n):=\eps_{\mathrm{DT}}(n)/10$.
\begin{enumerate}[label=\textup{(\roman*)}]
\item \emph{Local Existence.}  For every $\Lambda_0>0$ there is
$T_{\mathrm{ex}}=T_{\mathrm{ex}}(n,\Lambda_0)>0$ with the following property.  Let
$U\Subset M$ be a
smooth domain, let a possibly incomplete Ricci flow $g(t)$ be smooth on a
neighbourhood of $\overline U\times[0,T]$, let $N$ carry a complete Ricci flow $h(t)$,
and suppose, on $[0,T]$, that for some $r>0$
\begin{equation*}
 \sup_U|\Rm_g|+\sup_N|\Rm_h|\le r^{-2},\qquad
 \sup_U|\nabla\Rm_g|\le r^{-3},
\end{equation*}
while all curvature derivatives of $h(t)$ are globally bounded on $[0,T]$.
If $f:U\to N$ is
smooth on a neighbourhood of $\overline U$ and $|df|^2\le\Lambda_0$, then the
Ricci--harmonic map heat flow starting from $f$ exists, with $|dF|^2\le10\Lambda_0$,
on
\begin{equation*}
 U_r:=\{x\in U:B_{g(0)}(x,r)\Subset U\}
 \quad\hbox{for}\quad 0\le t<\min\{T,T_{\mathrm{ex}}r^2\}.
\end{equation*}
The proof is \cref{rf:lem:local-harmonic-existence}.

\item \emph{Local Derivatives.}  Given $A,\Lambda_0>0$ and an integer $q\ge0$,
there are $S_q=S_q(n,A,\Lambda_0,q)>0$ and
$C_q^{\mathrm{map}}=C_q^{\mathrm{map}}(n,A,\Lambda_0,q)<\infty$ such that the following holds.  If
the Ricci flows $g(t),h(t)$ satisfy $B_{g(0)}(x,4r)\Subset M$, the target
flow is complete, and
\begin{equation*}
 \sup_{B_{g(0)}(x,4r)}|\Rm_g|+\sup_N|\Rm_h|\le A/t,
 \qquad |dF|^2\le\Lambda_0,
\end{equation*}
where $F$ is a smooth Ricci--harmonic map heat flow on the radius-$4r$
source cylinder, then on $B_{g(0)}(x,r)$
\begin{equation*}
 |\nabla^j dF|^2\le C_q^{\mathrm{map}}t^{-j}\qquad(0\le j\le q)
\end{equation*}
for $0<t\le\min\{T,S_qr^2\}$; see \cref{rf:lem:local-map-derivatives}.

\item \emph{Ricci--DeTurck Derivatives.} Let $g$, $h$, $x$, $r$, $\Lambda_0$ and
$F$ be as in \textup{(ii)}, and suppose in addition that every $F_t$ is a
diffeomorphism onto its image with
\begin{equation*}
 (1-\eps_{\mathrm L})g(t)\le F_t^*h(t)\le(1+\eps_{\mathrm L})g(t)
 \quad\hbox{on }B_{g(0)}(x,4r).
\end{equation*}
Then $\widehat g=(F_t^{-1})^*g$ satisfies, on $F_t(B_{g(0)}(x,r))$ and for
$0<t\le\min\{T,S_qr^2\}$,
\begin{equation*}
 |(\nabla^h)^j\widehat g|_h\le C_q^{\mathrm{DT}}t^{-j/2}\qquad(1\le j\le q),
\end{equation*}
with $C_q^{\mathrm{DT}}=C_q^{\mathrm{DT}}(n,A,\Lambda_0,q)$; see
\cref{rf:lem:local-dt-derivatives}.

\item \emph{Ricci--DeTurck Decay.} Given $A,C_{\mathrm{dec}}>0$ and an integer $L\ge1$ there is
a constant $T_L^{\mathrm{dec}}>0$, depending only on $n$, $A$, $C_{\mathrm{dec}}$ and $L$, with the
following property.
Let $h(t)$ be a smooth Ricci flow on $N$ with $B_{h(0)}(y,4r)\Subset N$ and
$|\Rm_h|\le A/t$ there, and let $\widehat g$ be a smooth Ricci--DeTurck flow
relative to $h$ on a neighbourhood of the closed reference cylinder
$\overline{B_{h(0)}(y,4r)}\times[0,T]$, with $\widehat g(0)=h(0)$, satisfying
there
\begin{equation*}
 (1-\eps_{\mathrm{DT}})h\le\widehat g\le(1+\eps_{\mathrm{DT}})h,
 \qquad
 |(\nabla^h)^2\widehat g|_h\le\frac{C_{\mathrm{dec}}}{t}.
\end{equation*}
Then, on $B_{h(0)}(y,r)$ and for $0<t\le\min\{T,T_L^{\mathrm{dec}}r^2\}$,
\begin{equation*}
 |\widehat g-h|_h\le(t/r^2)^L,
 \qquad |\nabla^h(\widehat g-h)|_h^2\le r^{-2}(t/r^2)^L.
\end{equation*}
\end{enumerate}
\end{proposition}

\begin{proof}[Proof of part \textup{(iv)}]
Rescale $r=1$ and write
\begin{equation*}
 H=\widehat g-h,\qquad v=|H|_h^2,
\end{equation*}
so that $|(\nabla^h)^2H|\le C_{\mathrm{dec}}/t$ on $B_{h(0)}(y,4)$ by hypothesis.
The covariant discrepancy equation gives, for small $\eps_{\mathrm{DT}}$,
\begin{equation}\label{rf:eq:discrepancy-decay-pde}
 (\partial_t-\Delta_h)v\le\mathcal L_0v,\qquad
 \mathcal L_0=C_n\bigl(|\Rm_h|+|(\nabla^h)^2H|\bigr)\le b(n,A,C_{\mathrm{dec}})/t;
\end{equation}
see \cite[Lemma~2.5]{Lee26}; there the term $-|\nabla H|^2$ absorbs the quadratic
gradient terms, and replacing the operator of $\widehat g$ by $\Delta_h$ costs
$C_n|H|^2|\nabla^2H|$. Both flows are smooth on the closed cylinder, so $\mathcal L_0$ is
bounded there; only the bound $b/t$ degenerates at $t=0$. After enlarging $b$, one also
has $\Ric_h\le b/t$ and $v\le b/t$ for $t\le1$. For each $z\in B_{h(0)}(y,3/2)$, the ball
of radius $2$ about $z$ has closure in $B_{h(0)}(y,7/2)\Subset B_{h(0)}(y,4)$.
Since $v(0)=0$, \cite[Theorem~1.1]{LT22}, with an integer exponent
$k>\max\{2J,b+1\}$, therefore gives uniformly
\begin{equation*}
 \sup_{B_{h(0)}(y,3/2)}|H|_h\le t^J
\end{equation*}
for any prescribed integer $J$, on a time depending only on $n,A,C_{\mathrm{dec}},J$.

For $x\in B_{h(0)}(y,1)$, the shrinking-balls lemma \cite[Lemma~2.2]{ST21}, applied to
$B_{h(0)}(x,3/8)\Subset B_{h(0)}(y,3/2)$, gives $B_{h(t)}(x,1/4)\Subset B_{h(0)}(y,3/2)$
for a common smaller time. Along a unit-speed $h(t)$-geodesic of length
$s=t^{(J+1)/2}\le1/4$, which exists because $B_{h(t)}(x,1/4)$ is relatively compact,
express $H$ in a parallel frame and apply Taylor's formula to each component. This
gives
\begin{equation}\label{rf:eq:gradient-decay-interpolation}
 |\nabla^hH|(x,t)\le
 C_n\left(\frac{\sup|H|}{s}+s\sup|(\nabla^h)^2H|\right)
 \le C(n,A,C_{\mathrm{dec}},J)t^{(J-1)/2}.
\end{equation}
Choose $J=L+2$ and decrease the time to absorb the spare factor of $t$ in
the square of \eqref{rf:eq:gradient-decay-interpolation}. Restoring $r$ gives the two
bounds.
\end{proof}

\Cref{rf:sec:comparison} applies these estimates on successive time slabs.

\subsection{The Local Improvement and the First-Time Argument}\label{rf:sec:first-time}

We state the local improvement for incomplete flows and derive \cref{rf:thm:main};
\crefrange{rf:sec:selected}{rf:sec:bootstrap} prove the local statement.

\begin{proposition}[Local Type-I Improvement]\label{rf:prop:bootstrap}
Fix $n\ge2$, $A\ge1$ and $0<\alpha<A$. There are
$\delta_{\mathrm{imp}}=\delta_{\mathrm{imp}}(n,A,\alpha)>0$, $\kappa_{\mathrm{imp}}=\kappa_{\mathrm{imp}}(n,A,\alpha)>0$, and
$\Theta_{\mathrm{imp}}=\Theta_{\mathrm{imp}}(n,A,\alpha)\in(0,1)$ with the following property. Let $g(t)$,
$0\le t<T$, be a smooth Ricci flow on an
$n$-manifold $U$, let $p\in U$ and $r>0$, and assume that
\begin{enumerate}[label=\textup{(\alph*)}]
\item the closed ball $\overline{B_{g(0)}(p,4r)}$ is compact;
\item on $B_{g(0)}(p,4r)\times(0,T)$,
\begin{equation}\label{rf:eq:hyp-typeI}
|\Rm_{g(t)}|\le\frac{A}{t};
\end{equation}
\item on $B_{g(0)}(p,2r)$,
\begin{equation}\label{rf:eq:hyp-ricci}
\Ric_{g(0)}\ge-(n-1)\kappa_{\mathrm{imp}}r^{-2}g(0);
\end{equation}
\item
\begin{equation}\label{rf:eq:hyp-gh}
d_{GH}\bigl(B_{g(0)}(p,2r),B_{\R^n}(0,2r)\bigr)\le\delta_{\mathrm{imp}}r.
\end{equation}
\end{enumerate}
Then
\begin{equation}\label{rf:eq:improved}
|\Rm_{g(t)}|(p)\le\frac{\alpha}{t}
\quad\bigl(0<t<\min\{T,\Theta_{\mathrm{imp}} r^2\}\bigr),
\end{equation}
and consequently
\begin{equation}\label{rf:eq:improved-all}
|\Rm_{g(t)}|(p)\le\frac{\alpha}{t}+\frac{A}{\Theta_{\mathrm{imp}}}\,r^{-2}
\quad(0<t<T).
\end{equation}
\end{proposition}

All hypotheses concern the initial ball, because no Ricci lower bound is assumed at
positive time. Lifting the global flow to each local cover and applying the gain
$\alpha<A$ while $|\Rm|\le1/t$ gives the following first-time argument.

\begin{proof}[Proof of \cref{rf:thm:main}]
Apply \cref{rf:prop:bootstrap} with $A=1$ and the prescribed $\alpha$.
Let $\delta_{\mathrm{imp}},\kappa_{\mathrm{imp}},\Theta_{\mathrm{imp}}$ be its constants and put
\begin{equation*}
 r_0:=\min\left\{\frac{\rho}{16},
                  \frac{\sqrt{\kappa_{\mathrm{imp}}}}2\right\},
 \qquad
 \delta_0:=\frac{\delta_{\mathrm{imp}}}{2}.
\end{equation*}
Hamilton's theorem \cite{Ham82} starts the maximal flow for a positive, possibly
nonuniform time.

Fix $x\in M$ and lift the flow to the fixed $g_0$-cover
$\tB_{g_0}(x,\rho)$. Since $4r_0\le\rho/4$, \cref{rf:lem:cover} makes the closed
$4r_0$-ball about $\tilde x$ compact, while \eqref{rf:eq:lrrp} gives
\begin{equation*}
 d_{GH}\bigl(B_{\pi_x^*g_0}(\tilde x,2r_0),
          B_{\R^n}(0,2r_0)\bigr)
 <2\delta_0r_0\le\delta_{\mathrm{imp}}r_0.
\end{equation*}
The choice $r_0^2\le\kappa_{\mathrm{imp}}/4$ and the initial Ricci lower bound give
\begin{equation*}
 \Ric_{\pi_x^*g_0}\ge-(n-1)\kappa_{\mathrm{imp}}r_0^{-2}\pi_x^*g_0.
\end{equation*}

Define
\begin{equation}\label{rf:eq:TA}
 T_A:=\sup\left\{s<T_{\max}:
 t\sup_M|\Rm_{g(t)}|\le A
 \text{ for every }0<t\le s\right\}.
\end{equation}
Smoothness gives $T_A>0$. Before $T_A$, every lifted flow satisfies
\cref{rf:prop:bootstrap} with $p=\tilde x$ and $r=r_0$, hence
\begin{equation}\label{rf:eq:ft-improved}
 t\sup_M|\Rm_{g(t)}|\le\alpha<A
 \quad\bigl(0<t<\min\{T_A,\Theta_{\mathrm{imp}} r_0^2\}\bigr).
\end{equation}
If $T_A\le\Theta_{\mathrm{imp}} r_0^2$ and $T_A<T_{\max}$, the strict margin and continuity extend
\eqref{rf:eq:TA}; if $T_A=T_{\max}$, \eqref{rf:eq:ft-improved} bounds curvature on the final
half-interval and continuation extends the flow. Thus
\begin{equation*}
 T_{\max}>\Theta_{\mathrm{imp}} r_0^2,\qquad T_A>\Theta_{\mathrm{imp}} r_0^2.
\end{equation*}
One final application at every lifted point gives \eqref{rf:eq:main-bound} with
$T_0:=\frac12\Theta_{\mathrm{imp}} r_0^2$.
\end{proof}

\subsection{A Complete Comparison Flow}\label{rf:sec:selected}

We construct the complete comparison flow using Hochard's conformal completion
\cite{Hoc16,HW22} and \cref{rf:thm:wang}. The rescaling-invariant lemmas are stated for
reuse in \cref{rf:sec:noncompact}.

\subsubsection{Constants}\label{rf:sec:constants}
Set $\beta:=\min\{\alpha/32,1/(200n)\}$. This is the Type-I constant of the comparison
flow: it lies below Wang's threshold $1/(100n)$, and its gap to $\alpha$ absorbs the
comparison error in \eqref{rf:eq:boot-main}. Let $\delta_W=\delta_W(n,\beta)$ be from
\cref{rf:thm:wang}, and put $\delta_P:=1-e^{-\delta_W^2/2}$ and $\eps_P:=\delta_W/2$; in
\cref{rf:lem:selected} the isoperimetric deficit $\delta_P$ gives the entropy bound of
\cref{rf:thm:wang}, and $\eps_P$ fixes the lifetime of the comparison flow relative to the
isoperimetric scale. Choose the parameter $\theta=\theta(n,\beta)>0$ of the conformal
profile of \cref{rf:lem:completion} and then the isoperimetric deficit $\delta_I>0$ of the
initial metric so that
\begin{equation}\label{rf:eq:deltaI}
 \frac{1-\delta_I}{(1+2\theta)^{n(n-1)}}\ge1-\frac12\delta_P.
\end{equation}
All these constants depend only on $n$ and $\alpha$.

\subsubsection{Initial Isoperimetry}\label{rf:sec:iso}
The first lemma converts the initial hypotheses into definite-scale isoperimetry.

\begin{lemma}[Initial Isoperimetry]\label{rf:lem:isoperimetry}
For every $\eta\in(0,1)$ there are constants $c(\eta)\in(0,1/10)$, $\delta(\eta)>0$ and
$\kappa(\eta)>0$, depending only on $n$ and $\eta$, with the following property. Let $g$
be a Riemannian metric on an $n$-manifold, let $p$ be a point and $r>0$, and suppose
that the closed ball of radius $r/2$ about $p$ is compact and that
\begin{equation}\label{rf:eq:iso-hyp}
 \Ric_g\ge-(n-1)\kappa(\eta)r^{-2}g\ \hbox{ on }B_g(p,2r),\qquad
 d_{GH}\bigl(B_g(p,2r),B_{\R^n}(0,2r)\bigr)\le\delta(\eta)r .
\end{equation}
Then, with $R:=c(\eta)r$,
\begin{equation}\label{rf:eq:scal-R}
 \Scal_g\ge-R^{-2}\quad\hbox{on }B_g(p,R),
\end{equation}
and every regular domain $E\Subset B_g(p,R)$ satisfies
\begin{equation}\label{rf:eq:iso}
 \bigl(\area_g\partial E\bigr)^n
 \ge(1-\eta)n^n\omega_n
       \bigl(\vol_gE\bigr)^{n-1}.
\end{equation}
\end{lemma}

\begin{proof}
Scale to $r=1$. Let
$\bar\varepsilon_{0,n},\bar\eta_{0,n},\bar\delta_{0,n},C_{0,n}$ be the constants of
\cite[Theorem~1.1]{CM20} for lower Ricci bound $0$ and dimension $n$. Choose
$\eta_V\le\bar\eta_{0,n}$, $r_V\le\min\{\bar\delta_{0,n},1/10\}$, and
$\varepsilon_V\le\bar\varepsilon_{0,n}$ so that
\begin{equation*}
 C_{0,n}(r_V+\eta_V+\varepsilon_V)
 \le 1-(1-\eta)^{1/n}.
\end{equation*}
We use Colding's local volume continuity
\cite[Main Lemma~2.1 and Theorem~0.1]{Col97} in the following form. There are
$\delta_C,\kappa_C\in(0,1)$ and $r_C>1$ such that
$\Ric\ge-(n-1)\kappa_C$ on $B(y,r_C)$ and
$d_{GH}(B(y,r_C),B_{\R^n}(0,r_C))<\delta_C$ imply
$\vol B(y,1)\ge(1-\eta_V)\omega_n$; this is the form used in
\cite[proof of Lemma~2.4]{HW22}.

Set $c_C=(2r_C)^{-1}$ and $g_C=c_C^{-2}g$. Restricting a correspondence in
\eqref{rf:eq:iso-hyp} to concentric balls and rescaling gives
\begin{equation*}
 d_{GH}\bigl(B_{g_C}(p,r_C),B_{\R^n}(0,r_C)\bigr)
 \le c_C^{-1}\Psi(\delta(\eta)\mid n),
 \qquad
 \Ric_{g_C}\ge-(n-1)\kappa(\eta)c_C^2g_C,
\end{equation*}
where $\Psi(\delta\mid n)\to0$ with $\delta$. These balls correspond to relatively
compact $g$-balls in $B_g(p,1/2)$. Choose $\delta(\eta),\kappa(\eta)$ so that this error
is below $\delta_C$, $\kappa(\eta)c_C^2\le\kappa_C$, and
$(n-1)\kappa(\eta)c_C^2\le\varepsilon_V$. Thus
$\vol_{g_C}B_{g_C}(p,1)\ge(1-\eta_V)\omega_n$, and \cite[Theorem~1.1]{CM20} gives,
for finite-perimeter $E\subset B_{g_C}(p,r_V)$,
\begin{equation*}
 \operatorname{Per}_{g_C}(E)\ge n\omega_n^{1/n}
 \bigl(1-C_{0,n}(r_V+\eta_V+(n-1)\kappa(\eta)c_C^2)\bigr)
 \vol_{g_C}(E)^{(n-1)/n}.
\end{equation*}
The parameter choice, the $n$th power, and rescaling prove \eqref{rf:eq:iso} with
$c(\eta)=c_Cr_V$. Finally \eqref{rf:eq:iso-hyp} gives
$\Scal_g\ge-n(n-1)\kappa(\eta)$; decreasing $\kappa(\eta)$ until
$n(n-1)\kappa(\eta)\le c(\eta)^{-2}$ proves \eqref{rf:eq:scal-R}.
\end{proof}

For $\eta=\delta_I$, fix
\begin{equation}\label{rf:eq:cI}
 c_I:=c(\delta_I),\qquad\delta_{\mathrm{imp}}:=\delta(\delta_I),\qquad\kappa_{\mathrm{imp}}:=\kappa(\delta_I),
\end{equation}
the constants of \cref{rf:prop:bootstrap}, all depending only on $n,\alpha$.

\subsubsection{Conformal Completion}\label{rf:sec:completion}
The next conformal change is supported outside a smaller ball and nearly preserves
isoperimetry; its domain is a superlevel-set component of the cut-off.

\begin{lemma}[Conformal Completion]\label{rf:lem:completion}
For every $C_\chi\ge1$ there are constants $c_H\in(0,1/8]$ and $C_H>0$, depending only
on $n$, $\beta$ and $C_\chi$, with the following property. Let $g$ be a Riemannian
metric on an $n$-manifold, let $p$ be a point and $R>0$, and suppose that the closed
ball of radius $3R/4$ about $p$ is compact, that \eqref{rf:eq:scal-R} holds, that
\eqref{rf:eq:iso} holds with $\eta=\delta_I$, and that there is a smooth function $\chi$
with values in $[0,1]$, supported in $B_g(p,3R/4)$ and equal to one on $B_g(p,R/3)$,
such that
\begin{equation}\label{rf:eq:cutoff}
 R|\nabla\chi|+R^2|\Delta\chi|\le C_\chi.
\end{equation}
Then there are an open set $\Omega$ with
\begin{equation*}
 B_g(p,R/4)\Subset\Omega\Subset B_g(p,3R/4)
\end{equation*}
and a smooth function $u\ge0$ on $\Omega$ such that $h_0:=e^{2u}g$ has the following
properties, where $r_H:=c_HR$:
\begin{enumerate}[label=\textup{(\roman*)}]
\item $(\Omega,h_0)$ is complete, $h_0\ge g$, and
      $h_0=g$ on $B_g(p,R/4)$;
\item $h_0$ has bounded curvature;
\item $\Scal_{h_0}\ge-r_H^{-2}$;
\item every regular domain $E$ compactly contained in an $h_0$-ball of
      radius $r_H$ satisfies
      \begin{equation}\label{rf:eq:iso-P}
       \bigl(\area_{h_0}\partial E\bigr)^n
       \ge(1-\delta_P)n^n\omega_n\bigl(\vol_{h_0}E\bigr)^{n-1}.
      \end{equation}
\end{enumerate}
\end{lemma}

\begin{proof}
Choose a regular value $a\in(1/4,1/2)$ of $\chi$ and let $\Omega$ be the component of
$\{\chi>a\}$ containing $B_g(p,R/3)$; it has smooth compact closure in
$B_g(p,3R/4)$. Fix a smooth decreasing $\upsilon_0:[0,1]\to[0,1]$ with
$\upsilon_0=0$ near $1$, $\upsilon_0(0)=1$, $\upsilon_0'(0)<0$, and
$\upsilon_0(s)<1$ for $s>0$, and set
\begin{equation*}
 \varphi:=\upsilon_0\circ\frac{\chi-a}{1-a}.
\end{equation*}
Then $\varphi=0$ on $B_g(p,R/3)$, $0\le\varphi<1$, $\varphi\to1$ at
$\partial\Omega$, and \eqref{rf:eq:cutoff} gives
\begin{equation}\label{rf:eq:phi-bounds}
 R|\nabla\varphi|+R^2|\Delta\varphi|\le C_\varphi,
\end{equation}
where $C_\varphi=C_\varphi(C_\chi)$. The smooth extensions of $g,\varphi$ past
$\overline\Omega$ have finite curvature and Hessian norms there.

Let $w_\theta:[0,1)\to[0,\infty)$ be Hochard's conformal profile
\cite[proof of Lemma~6.2]{Hoc16}, smoothed as in
\cite[equations~(2.3)--(2.5)]{HW22}: it is nondecreasing and convex, vanishes on
$[0,1-3\theta/2]$, tends to infinity at $1$, and satisfies
\begin{equation}\label{rf:eq:w-deriv}
 e^{-w_\theta}w_\theta'\le C\theta^{-1},\qquad
 e^{-2w_\theta}\bigl(w_\theta''+(w_\theta')^2\bigr)
       \le C\theta^{-2}.
\end{equation}
For $1-3\theta/2<s<1$ and $b=\theta(1-s)/4$, it also satisfies
\begin{equation}\label{rf:eq:w-collar}
 e^{w_\theta(s+b)-w_\theta(s-b)}
       \le1+2\theta,\qquad
 b e^{w_\theta(s-b)}\ge c\theta^2.
\end{equation}
Set $u=w_\theta\circ\varphi$ and $h_0=e^{2u}g$. The second estimate in
\eqref{rf:eq:w-collar} and \eqref{rf:eq:phi-bounds} make every escaping curve have infinite
$h_0$-length. Thus $h_0$ is complete, dominates $g$, and equals $g$ on
$B_g(p,R/4)$. The conformal scalar-curvature formula,
\eqref{rf:eq:scal-R}, \eqref{rf:eq:phi-bounds}, and \eqref{rf:eq:w-deriv} give
\begin{equation*}
 \Scal_{h_0}\ge-C_HR^{-2}.
\end{equation*}
The full-curvature formula and the finite bounds on $\overline\Omega$ give bounded
sectional curvature.

Call the set $\{1-3\theta/2<\varphi<1\}$, where $u>0$, the collar. For a point $y$ of the
collar put $s=\varphi(y)$ and $b=\theta(1-s)/4$. The calculation
of \cite[equations~(2.23)--(2.33)]{HW22}, or the same first-exit argument using
\eqref{rf:eq:phi-bounds}--\eqref{rf:eq:w-collar}, gives, for
\mbox{$c_{\mathrm{col}}\le c\theta^2/C_\varphi$},
\begin{equation*}
 B_{h_0}(y,c_{\mathrm{col}}R)\subset\varphi^{-1}(s-b,s+b),\qquad
 e^{2w_\theta(s-b)}g\le h_0
 \le e^{2w_\theta(s+b)}g,
\end{equation*}
where the metric inequality holds on the displayed ball, which lies in $B_g(p,R)$ in
the original metric. Comparing boundary area and volume with the two conformal factors,
\eqref{rf:eq:iso} and \eqref{rf:eq:w-collar} give, for every regular domain $E$ compactly
contained in that ball,
\begin{equation}\label{rf:eq:iso-h0}
 \frac{\bigl(\area_{h_0}\partial E\bigr)^n}{\bigl(\vol_{h_0}E\bigr)^{n-1}}\ge
 (1-\delta_I)(1+2\theta)^{-n(n-1)}n^n\omega_n
 \ge(1-\tfrac12\delta_P)n^n\omega_n
 >(1-\delta_P)n^n\omega_n.
\end{equation}
The middle inequality is \eqref{rf:eq:deltaI}. Set
$c_H:=\min\{c_{\mathrm{col}}/2,1/8,C_H^{-1/2}\}$. Then (iii) holds. An $r_H$-ball
contained in $\{u=0\}$ satisfies \eqref{rf:eq:iso}; otherwise it meets the collar and lies
in a $2r_H$-ball to which \eqref{rf:eq:iso-h0} applies. Domain monotonicity proves (iv).
\end{proof}

The curvature bound in (ii) is qualitative; the next lemma uses it only to start the
Ricci flow of $h_0$.

\subsubsection{The Comparison Flow}\label{rf:sec:selected-flow}
We now apply \cref{rf:thm:wang} at every point.

\begin{lemma}[Complete Comparison Flow]\label{rf:lem:selected}
Let $h_0$ be a complete Riemannian metric with bounded curvature on an $n$-manifold
$\Omega$, and let $r_H>0$. Suppose that $\Scal_{h_0}\ge-r_H^{-2}$ and that every regular
domain compactly contained in an $h_0$-ball of radius $r_H$ satisfies
\eqref{rf:eq:iso-P}. Then $h_0$ has a complete bounded-curvature Ricci flow $h(t)$ on
$0\le t\le T_H:=(\eps_Pr_H)^2$, and
\begin{equation}\label{rf:eq:selected-bound}
 |\Rm_{h(t)}|\le\frac{\beta}{t}
 \quad(0<t\le T_H).
\end{equation}
\end{lemma}

\begin{proof}
Shi \cite{Shi89} starts a complete bounded-curvature flow. Let $0<T\le T_H$ precede
its maximal time. At each $x$, the scalar and isoperimetric hypotheses and
\cite[Lemma~3.5]{Wan18} give
\begin{equation*}
 \nu\bigl(B_{h_0}(x,r_H),h_0,s\bigr)
 \ge\log(1-\delta_P)-sr_H^{-2}\qquad(s>0).
\end{equation*}
Since $\delta_W^{-1}\sqrt T\le r_H/2$ and $\nu$ does not decrease when the ball
shrinks \cite[Proposition~2.1]{Wan18}, the definitions of $\delta_P,\eps_P$ give
\begin{equation*}
 \nu\bigl(B_{h_0}(x,\delta_W^{-1}\sqrt T),h_0,T\bigr)
 \ge-\tfrac12\delta_W^2-\eps_P^2=-\tfrac34\delta_W^2>-\delta_W^2 .
\end{equation*}
Thus \cref{rf:thm:wang} applies at every point and gives
\begin{equation}\label{rf:eq:selected-T}
 |\Rm_{h(t)}|\le\frac{\beta}{t}
 \quad(0<t\le T).
\end{equation}
If the maximal time were at most $T_H$, letting $T$ approach it would bound curvature
on its final half and Shi continuation would extend the flow. Hence the flow reaches
$T_H$, and \eqref{rf:eq:selected-T} becomes \eqref{rf:eq:selected-bound}.
\end{proof}

\subsection{Localized Comparison}\label{rf:sec:comparison}

The comparison flow of \cref{rf:sec:selected} must be compared with the given flow of
\cref{rf:prop:bootstrap}.  The two agree initially only on a ball, and only the comparison
flow is complete.  The following lemma is the local form of Lee's Ricci--harmonic map
comparison argument \cite{Lee26} that this requires.  It is the main analytic step of
the paper: on a time interval proportional to $\ell^2$, where $\ell$ is the radius scale
of the lemma, the given flow is, after a change of gauge, a perturbation of the
comparison flow that vanishes to any prescribed order at the initial time.

\begin{lemma}[Localized Type-I Comparison]\label{rf:lem:comparison}
Let $N\subset U$ be open subsets of one smooth manifold, with equality allowed,
let $g(t)$ be a
Ricci flow on $U\times[0,T]$, let $p\in U$ and $\ell>0$, and put
$D_j:=B_{g(0)}(p,j\ell)$.  Suppose
$D_{64}\Subset N$, and let $h(t)$, $0\le t\le T$, be a complete smooth Ricci flow
on $N$.  Assume
\begin{equation}\label{rf:eq:initial-order}
 h(0)=g(0)\ \hbox{on }D_{64},\qquad h(0)\ge g(0)
 \ \hbox{on }N,
\end{equation}
where $h(0)\ge g(0)$ is an inequality of quadratic forms, that is,
$h(0)(\xi,\xi)\ge g(0)(\xi,\xi)$ for every tangent vector $\xi$ at a point of
$N$,
and
\begin{equation}\label{rf:eq:comp-curv}
 \sup_{D_{64}}|\Rm_{g(t)}|+\sup_N|\Rm_{h(t)}|
       \le \frac{A_*}{t}.
\end{equation}
For every integer $L\ge1$ there are constants $\tau_L>0$, $C_{\mathrm{disp}}$ and
$C_L<\infty$, depending only on $n$, $A_*$ and $L$, with the following property. There
is a
Ricci--harmonic map heat flow $F$ from $D_8\times[0,\min\{T,\tau_L\ell^2\})$ to $N$ with
$F_0=\mathrm{id}$ such that, for $0<t<\min\{T,\tau_L\ell^2\}$, the map $F_t$ is a
diffeomorphism of $D_8$ onto its image and
\begin{equation}\label{rf:eq:comp-disp}
 d_{g(0)}(F_t(x),x)\le C_{\mathrm{disp}}\sqrt t,\qquad
 D_4\subset F_t(D_8)\subset D_{16}.
\end{equation}
Moreover, for the dimensional constant $\eps_*:=\eps_{\mathrm L}(n)$ fixed in
\cref{rf:prop:lee-inputs},
\begin{equation}\label{rf:eq:comp-equiv}
 (1-\eps_*)g(t)\le F_t^*h(t)\le(1+\eps_*)g(t)
 \quad\hbox{on }D_8.
\end{equation}
If $\widehat g=(F_t^{-1})^*g$ and $H=\widehat g-h$, then on $D_2$
\begin{equation}\label{rf:eq:comp-H}
 |H|_h\le C_L(t/\ell^2)^L,\qquad
 |\nabla^hH|_h^2\le C_L\ell^{-2}(t/\ell^2)^L,
\end{equation}
and on $D_1$
\begin{equation}\label{rf:eq:comp-Rm}
 t|\Rm_{\widehat g}-\Rm_h|_h
       \le C_L(t/\ell^2)^{L/2}.
\end{equation}
\end{lemma}

The proof occupies the rest of the section. We rescale so that $\ell=1$ and keep
the notation of the lemma. After fixing the constants and the induction hypothesis
(\cref{rf:sec:slab-constants}), we prove the step from one time slab to the next as a
separate statement (\cref{rf:lem:slab}), and then run the induction and pass from the
gauge estimates to the curvature comparison (\cref{rf:sec:slab-induction}). The
construction follows \cite[Sections~2--3]{Lee26}, with the local estimates of
\cref{rf:prop:lee-inputs} in place of global ones. No global bound on the initial
curvature of $h$ is needed: the initial Dirichlet solution uses only a compact source
and target neighbourhood, and every later extension starts at a positive time, when
completeness and the Type-I bound give global Shi estimates for $h$.

\subsubsection{Constants and the Induction Hypothesis}\label{rf:sec:slab-constants}
The constants are fixed in the following order; each depends only on $n$, $A_*$
and $L$, and on those chosen before it. With $\eps_*=\eps_{\mathrm L}(n)$ as in
\cref{rf:lem:comparison}, set $\Lambda_0:=2n$ and $\Lambda_1:=10\Lambda_0$. Put $a_*:=c(n)A_*$, so that
integration of $|\partial_th|_h\le c(n)A_*/t$ gives
\begin{equation}\label{rf:eq:growth}
 |\xi|_{h(t)}\le(t/s)^{a_*}|\xi|_{h(s)},\qquad 0<s\le t,
\end{equation}
for tangent vectors $\xi$, and choose an integer $m>2a_*+2$; it is the decay exponent in
\eqref{rf:eq:lee-decay} below, and the bound on $m$ makes the first integral in
\eqref{rf:eq:disp-int} converge.

A time slab beginning at a time $t_k>0$ is handled with the two scales
\begin{equation*}
 \sigma_k:=c_\sigma\sqrt{t_k},\qquad
 \sigma_k':=\frac{\Lambda_{\mathrm{mg}}}{8}\sqrt{t_k},
\end{equation*}
whose coefficients are fixed as follows. Choose $c_\sigma=c_\sigma(n,A_*)\in(0,1]$ so that the
curvature and first-curvature-derivative bounds \eqref{rf:eq:shi} below are at most
$\sigma_k^{-2}$ and $\sigma_k^{-3}$ on the slab. The metric equivalence
\eqref{rf:eq:ind-equiv} imposed in the induction gives
$|dF|^2\le n(1+\eps_*)<\Lambda_0$. Choose $\mu=\mu(n,A_*)\in(0,1]$ with
\begin{equation}\label{rf:eq:mu-choice}
 \mu<T_{\mathrm{ex}}(n,\Lambda_0)c_\sigma^2,\qquad
 C_{\mathrm{bar}}\mu\le\tfrac23\eps_*,
\end{equation}
where $C_{\mathrm{bar}}=C_{\mathrm{bar}}(n,A_*)$ is the constant of the tensor
barrier \eqref{rf:eq:barrier} below. The first inequality is the time restriction of
\cref{rf:prop:lee-inputs}\textup{(i)} at the scale $\sigma_k$, which
also supplies the bound $|dF|^2\le\Lambda_1$ on the new slab. Then choose
$\Lambda_{\mathrm{mg}}=\Lambda_{\mathrm{mg}}(n,A_*,L)\ge1$ so large that
\begin{equation}\label{rf:eq:Lambda-choice}
 \left(\frac{\Lambda_{\mathrm{mg}}}{8}\right)^2
 \min\bigl\{S_2(n,A_*,\Lambda_1,2),\,
 T_m^{\mathrm{dec}}(n,A_*,C_2^{\mathrm{DT}},m)\bigr\}\ge2,
 \qquad C_2^{\mathrm{DT}}:=C_2^{\mathrm{DT}}(n,A_*,\Lambda_1,2),
\end{equation}
so that part \textup{(ii)} of \cref{rf:prop:lee-inputs} at the scale $\sigma_k'$
covers the full interval $[0,t_{k+1}]$ whenever $t_{k+1}\le2t_k$, and parts
\textup{(iii)} and \textup{(iv)}, with exponent $m$ and scale $\sigma_k'$, cover
$[0,t_k]$. The constant $\Lambda_{\mathrm{mg}}$ is increased at four further places in the proof
of \cref{rf:lem:slab}, each time by a quantity depending only on $n$, $A_*$ and $L$.
Finally, the times $\tau'$ and $\tau_L$ are fixed in \cref{rf:lem:slab} and in its
application.

\begin{figure}[t]
\centering
\includegraphics[width=\textwidth]{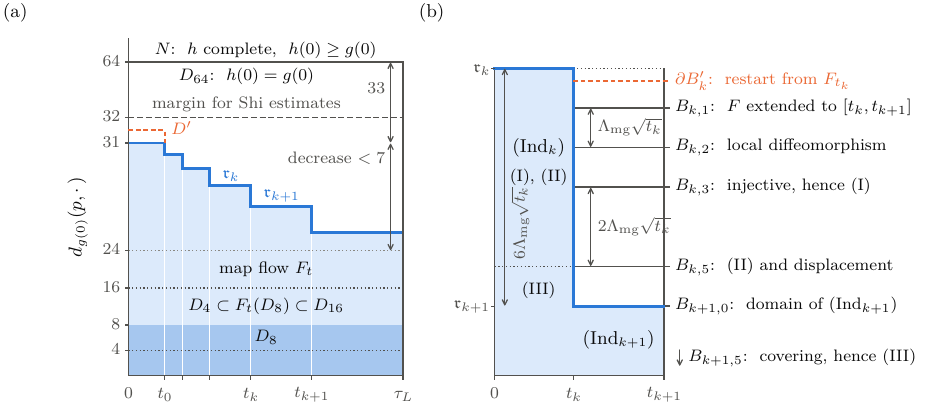}
\caption{The regions in the proof of \cref{rf:lem:comparison}, after the rescaling
$\ell=1$ (schematic, not to scale). (a)~Space-time. The map flow $F_t$ is constructed
on the shaded staircase: the initial Dirichlet flow lives on $D'$ up to the time $t_0$,
and at each time $t_k=(1+\mu)^kt_0$ the radius of the domain drops from $\mathfrak r_k$
to $\mathfrak r_{k+1}=\mathfrak r_k-6\Lambda_{\mathrm{mg}}\sqrt{t_k}$. The total decrease of the radius up to the time
$\tau_L$ is less than $7$, so the domain always contains the ball $D_8$ of the
statement (dark). The equality $h(0)=g(0)$ is used only inside $D_{32}$, and the $33$
units between the staircase and $\partial D_{64}$ contain the parabolic balls of the
local Shi estimates. (b)~The step from $(\mathrm{Ind}_k)$ to $(\mathrm{Ind}_{k+1})$ on
the slab $[t_k,t_{k+1}]$. The radius $6\Lambda_{\mathrm{mg}}\sqrt{t_k}$ lost in the step is spent in
margins of width $\Lambda_{\mathrm{mg}}\sqrt{t_k}$, and each property holds on the ball below the
corresponding line.}
\label{rf:fig:slabs}
\end{figure}

For $t_k>0$ and $\mathfrak r_k>0$ put
\begin{equation*}
 B_{k,j}:=B_{g(0)}\bigl(p,\mathfrak r_k-j\Lambda_{\mathrm{mg}}\sqrt{t_k}\bigr)
 \quad(j\ge0).
\end{equation*}
We say that a Ricci--harmonic map heat flow $F$ \emph{satisfies $(\mathrm{Ind}_k)$
for the data $(\mathfrak r_k,t_k)$} if $F$ is defined on $B_{k,0}\times[0,t_k]$
with $F_0=\mathrm{id}$ and has the following three properties there.
\begin{enumerate}[label=\textup{(\Roman*)}]
\item $F_t$ is a diffeomorphism onto its image;
\item
\begin{equation}\label{rf:eq:ind-equiv}
 (1-\eps_*)g(t)\le F_t^*h(t)\le(1+\eps_*)g(t);
\end{equation}
\item for $x\in B_{k,5}$ and $0\le t\le t_k$,
\begin{equation}\label{rf:eq:ind-disp}
 \begin{split}
 d_{g(0)}(x,F_t(x))&\le\frac14\Lambda_{\mathrm{mg}}\sqrt{t_k},\\
 B_{g(0)}(x,\Lambda_{\mathrm{mg}}\sqrt{t_k})
 &\subset F_t\bigl(B_{g(0)}(x,2\Lambda_{\mathrm{mg}}\sqrt{t_k})\bigr).
 \end{split}
\end{equation}
\end{enumerate}
\Cref{rf:fig:slabs} shows the regions that enter: panel~(a) the domains
$B_{k,0}\times[0,t_k]$ of the map flow along the induction, and panel~(b) the step
from $t_k$ to $t_{k+1}$.

\subsubsection{The Step from One Slab to the Next}\label{rf:sec:slab-step}
The next lemma extends the map flow across one time slab and recovers the induction
hypothesis on a domain shrunk by six margins. Its proof is the core of the argument:
the extension uses the local existence statement, a tensor barrier gives a rough
metric equivalence and injectivity on the new slab, and the decay of the
Ricci--DeTurck discrepancy on the earlier interval sharpens the equivalence back to
\eqref{rf:eq:ind-equiv}.

\begin{lemma}[One-Slab Extension]\label{rf:lem:slab}
Let $g$, $h$, $N$ and $p$ satisfy the hypotheses of \cref{rf:lem:comparison} with
$\ell=1$, and let the constants be as in \cref{rf:sec:slab-constants}. There is
$\tau'=\tau'(n,A_*,L)\in(0,1]$ with the following property. Let $0<t_k\le\tau'$ and
$24\le\mathfrak r_k\le31$, let $t_k<t_{k+1}\le(1+\mu)t_k$ with $t_{k+1}<T$, and
put $\mathfrak r_{k+1}:=\mathfrak r_k-6\Lambda_{\mathrm{mg}}\sqrt{t_k}$. If $F$ satisfies
$(\mathrm{Ind}_k)$ for the data $(\mathfrak r_k,t_k)$, then $F$ extends to a
Ricci--harmonic map heat flow on $B_{k,1}\times[0,t_{k+1}]$ that satisfies
$(\mathrm{Ind}_{k+1})$ for the data $(\mathfrak r_{k+1},t_{k+1})$.
\end{lemma}

\begin{proof}
\emph{Extension to the Next Slab.}
Since $\mathfrak r_k\le31$, one has
\begin{equation*}
 \operatorname{dist}_{g(0)}(B_{k,0},U\setminus D_{64})\ge33.
\end{equation*}
The shrinking-balls lemma \cite[Lemma~2.2]{ST21}, under the upper Ricci bound implied
by \eqref{rf:eq:comp-curv}, gives a dimensional constant $C_n$ such that, for $c>0$ and
$0\le t\le t_{k+1}$,
\begin{equation}\label{rf:eq:shrink}
 B_{g(t)}(x,c\sqrt{t_k})
 \subset B_{g(0)}
 \!\left(x,(c+C_n\sqrt{A_*})\sqrt{t_k}\right)
\end{equation}
whenever the ball on the right lies in $D_{64}$. Require
$\tau'\le33^2(1+C_n\sqrt{A_*})^{-2}$. Then, for $x\in B_{k,0}$ and
$t_k/2\le t\le t_{k+1}$, the ball $B_{g(t)}(x,\sqrt{t_k})$ lies in $D_{64}$, where the
Type-I bound holds. Local Shi estimates for $g$ on these parabolic neighbourhoods,
together with complete Shi estimates for $h$, give for $t_k\le t\le t_{k+1}$, on
$B_{k,0}$ for the source and globally for the target,
\begin{equation}\label{rf:eq:shi}
 |\nabla^j\Rm_g|+|\nabla^j\Rm_h|
       \le C_jt_k^{-1-j/2}.
\end{equation}
Choose a smooth domain $B_k'$ with
\begin{equation*}
 B_{g(0)}\bigl(p,\mathfrak r_k-\tfrac12\Lambda_{\mathrm{mg}}\sqrt{t_k}\bigr)
 \Subset B_k'\Subset B_{k,0}.
\end{equation*}
Apply \cref{rf:prop:lee-inputs}\textup{(i)} to the time-shifted flows $g(t_k+\cdot)$ and
$h(t_k+\cdot)$, at the scale $\sigma_k$, with initial map $F_{t_k}$ and source domain
$B_k'$.  Its curvature and first-curvature-derivative hypotheses hold on $B_k'$ by
\eqref{rf:eq:shi}, while \eqref{rf:eq:ind-equiv} bounds $|dF_{t_k}|$; the source metric and
$F_{t_k}$ are smooth on a neighbourhood of $\overline{B_k'}$, since
$\overline{B_k'}\subset B_{k,0}$; and \eqref{rf:eq:shi} gives the global bounds for the
curvature derivatives of $h$ on the new slab. Increase $\Lambda_{\mathrm{mg}}$, if necessary, so that
$\Lambda_{\mathrm{mg}}\ge2(c_\sigma+C_n\sqrt{A_*})$. Then \eqref{rf:eq:shrink} gives
\begin{equation*}
 B_{g(t_k)}(x,\sigma_k)\Subset B_k'
 \quad\hbox{for }x\in B_{k,1}.
\end{equation*}
Hence the corresponding interior set contains $B_{k,1}$. Since
$t_{k+1}-t_k\le\mu t_k<T_{\mathrm{ex}}\sigma_k^2$ by \eqref{rf:eq:mu-choice} and
$t_{k+1}<T$, the proposition extends $F$ to $B_{k,1}\times[0,t_{k+1}]$, with
$|dF|^2\le\Lambda_1$ on the new slab; the extended map is smooth, because the two
pieces agree at $t_k$ and the equation determines all their time derivatives there.

\emph{Lower Bound and Injectivity.}
It remains to bound $dF$ from below and to prove injectivity on a smaller domain.
Write
${\check h(s):=F_s^*h(s)}$.
For $x\in B_{k,2}$ one has
$\overline{B_{g(0)}(x,4\sigma_k')}\Subset B_{k,1}$.  Hence
\Cref{rf:prop:lee-inputs}\textup{(ii)}, applied on that ball, and
\eqref{rf:eq:Lambda-choice} give on $B_{k,2}\times[t_k,t_{k+1}]$
\begin{equation}\label{rf:eq:dF}
 |dF|\le C,\qquad |\nabla dF|\le Ct_k^{-1/2},\qquad
 |\nabla^2dF|\le Ct_k^{-1}.
\end{equation}
In particular, since $\partial_sF=\operatorname{tr}_{g(s)}\nabla dF$,
\begin{equation}\label{rf:eq:dtF}
 |\partial_sF|_{h(s)}\le Ct_k^{-1/2},\qquad
 |\nabla(\partial_sF)|_{h(s)}\le Ct_k^{-1}.
\end{equation}

These estimates supply a pointwise tensor barrier.  There
is $C_{\mathrm{bar}}=C_{\mathrm{bar}}(n,A_*)$ such that, if $0\le\eps<1$ and
on $B_{k,2}$ at time $t_k$
\begin{equation*}
 (1-\eps)g(t_k)\le\check h(t_k)\le(1+\eps)g(t_k),
\end{equation*}
then, for $t_k\le s\le t_{k+1}$ satisfying
$C_{\mathrm{bar}}(s-t_k)/t_k<1-\eps$,
\begin{equation}\label{rf:eq:barrier}
 \left(1-\eps-C_{\mathrm{bar}}\frac{s-t_k}{t_k}\right)g(s)
 \le\check h(s)\le
 \left(1+\eps+C_{\mathrm{bar}}\frac{s-t_k}{t_k}\right)g(s)
\end{equation}
on $B_{k,2}$. The time restriction keeps the lower scalar coefficient
$b_-(s):=1-\eps-C_{\mathrm{bar}}(s-t_k)/t_k$ positive; the upper coefficient
is automatically positive. To prove the lower bound, fix $x\in B_{k,2}$.
At a first null vector $\xi\in T_xU$, normalized by $g(s)(\xi,\xi)=1$,
of $\check h(s)-b_-(s)g(s)$, differentiating the tensor at fixed $x$ gives
\begin{equation}\label{rf:eq:contact}
 \begin{split}
 \partial_s\bigl[\check h(s)-b_-(s)g(s)\bigr](\xi,\xi)
 &=2h\bigl(dF(\xi),D_\xi\partial_sF\bigr)
   -2\Ric_h\bigl(dF(\xi),dF(\xi)\bigr)\\
 &\quad-b_-'(s)+2b_-(s)\Ric_g(\xi,\xi).
 \end{split}
\end{equation}
At contact $|dF(\xi)|_h^2=b_-(s)$; \eqref{rf:eq:comp-curv}, \eqref{rf:eq:dF}--\eqref{rf:eq:dtF}, and
$s\asymp t_k$ bound every term except $-b_-'(s)$ in absolute value by
$C(n,A_*)t_k^{-1}$.  Choosing $C_{\mathrm{bar}}>C(n,A_*)$ makes
\eqref{rf:eq:contact} positive, contrary to first contact.  Applying the same
calculation to $b_+(s)g(s)-\check h(s)$, where
\mbox{$b_+(s)=1+\eps+C_{\mathrm{bar}}(s-t_k)/t_k$}, proves the upper bound.
Replacing $\eps$ temporarily by $\eps+\eps'$ and letting $\eps'\downarrow0$
handles non-strict initial inequalities.

Apply \eqref{rf:eq:barrier} with $\eps=\eps_*$. Since
$C_{\mathrm{bar}}(s-t_k)/t_k\le C_{\mathrm{bar}}\mu\le\frac23\eps_*$ by
\eqref{rf:eq:mu-choice}, this gives
\begin{equation}\label{rf:eq:rough}
 (1-3\eps_*)g(s)\le\check h(s)\le(1+3\eps_*)g(s)
\end{equation}
on $B_{k,2}\times[t_k,t_{k+1}]$, so $F_s$ is a local diffeomorphism there.
It is injective on the smaller set $B_{k,3}$.  Indeed, suppose
$F_{s_0}(x_1)=F_{s_0}(x_2)=y_0$ for $x_1,x_2\in B_{k,3}$ and
$s_0\in(t_k,t_{k+1}]$. Local invertibility defines backward trajectories by
\begin{equation}\label{rf:eq:traj}
 \dot x_i(s)=-(dF_s|_{x_i(s)})^{-1}
                 (\partial_sF_s)(x_i(s)),\qquad
 x_i(s_0)=x_i.
\end{equation}
Here $dF_s$ is the spatial differential, and $\partial_sF_s$ is the time
derivative at a fixed source point. The chain rule gives
$\frac{d}{ds}[F_s(x_i(s))]=0$, so
$F_s(x_i(s))=F_{s_0}(x_i)=y_0$ along each trajectory.
As long as a trajectory remains in $B_{k,2}$, \eqref{rf:eq:dtF} and \eqref{rf:eq:rough} give
\begin{equation*}
 |\dot x_i(s)|_{g(s)}\le Ct_k^{-1/2}.
\end{equation*}
The Type-I bound compares $g(s)$ and $g(t_k)$ uniformly on this slab.
After integration, the local distance-distortion estimate between
$g(t_k)$ and $g(0)$ gives
\begin{equation*}
 d_{g(0)}(x_i(s),x_i(s_0))\le C_1(n,A_*,\mu)\sqrt{t_k}.
\end{equation*}
Choose $\Lambda_{\mathrm{mg}}>C_1$.  The margin of width $\Lambda_{\mathrm{mg}}\sqrt{t_k}$ between $B_{k,3}$
and the boundary of $B_{k,2}$ then prevents either trajectory from exiting
before $t_k$.  At $t_k$ the induction hypothesis gives
\mbox{$x_1(t_k)=x_2(t_k)$}, and uniqueness for the ODE \eqref{rf:eq:traj} gives $x_1=x_2$.
Thus $F_s|_{B_{k,3}}$ is a diffeomorphism onto its image.

\emph{Sharpening and Displacement.}
It remains to sharpen the rough estimate \eqref{rf:eq:rough} and recover
\eqref{rf:eq:ind-equiv}--\eqref{rf:eq:ind-disp}. Fix
$x\in B_{k,5}$.  By the induction hypothesis \eqref{rf:eq:ind-disp}, for every
$s\le t_k$,
\begin{equation*}
 B_{g(0)}(x,\Lambda_{\mathrm{mg}}\sqrt{t_k})
 \subset F_s\bigl(B_{g(0)}(x,2\Lambda_{\mathrm{mg}}\sqrt{t_k})\bigr).
\end{equation*}
Hence $\widehat g_s=(F_s^{-1})^*g(s)$ is defined on this fixed target ball.
It lies in $D_{32}\Subset D_{64}$, where $\widehat g_0=h(0)$.  For a fixed
$s\le t_k$ put $y=F_s(x)$.  The displacement estimate of $(\mathrm{Ind}_k)$ gives
$d_{g(0)}(x,y)\le\Lambda_{\mathrm{mg}}\sqrt{t_k}/4$.  Since $h(0)\ge g(0)$ and
$4\sigma_k'=\Lambda_{\mathrm{mg}}\sqrt{t_k}/2$,
\begin{equation*}
 \overline{B_{h(0)}(y,4\sigma_k')}
 \subset B_{g(0)}\left(x,\tfrac34\Lambda_{\mathrm{mg}}\sqrt{t_k}\right)
 \Subset B_{g(0)}(x,\Lambda_{\mathrm{mg}}\sqrt{t_k}).
\end{equation*}
Moreover, pushing forward \eqref{rf:eq:ind-equiv} and using
$\eps_*=\eps_{\mathrm{DT}}/10$ gives
\begin{equation*}
 (1-\eps_{\mathrm{DT}})h\le\widehat g_{s'}
 \le(1+\eps_{\mathrm{DT}})h\qquad(0\le s'\le t_k)
\end{equation*}
on this ball. Every point of $\overline{B_{h(0)}(y,4\sigma_k')}$ is of the form
$F_{s'}(z)$ with $z\in B_{g(0)}(x,2\Lambda_{\mathrm{mg}}\sqrt{t_k})$, for each $s'\le t_k$, and
$B_{g(0)}(z,4\sigma_k')\subset B_{g(0)}(x,\tfrac52\Lambda_{\mathrm{mg}}\sqrt{t_k})\subset B_{k,2}$.
So \cref{rf:prop:lee-inputs}\textup{(iii)}, applied at every such $z$ with the radius
$\sigma_k'$ and $q=2$, gives
$|(\nabla^h)^2\widehat g_{s'}|_h\le C_2^{\mathrm{DT}}/s'$ on
$\overline{B_{h(0)}(y,4\sigma_k')}\times(0,t_k]$, by \eqref{rf:eq:Lambda-choice}.
Thus the hypotheses of \cref{rf:prop:lee-inputs}\textup{(iv)} hold on the closed
cylinder over this ball, centred at $y$ with radius $\sigma_k'$, and applying it
with the exponent $m$ and evaluating at $y$ gives, for $0<s\le t_k$,
\begin{equation}\label{rf:eq:lee-decay}
 |\widehat g_s-h(s)|_h
 \le C\left(\frac{64s}{\Lambda_{\mathrm{mg}}^2t_k}\right)^m,
 \qquad
 |\nabla^h\widehat g_s|_h^2
 \le C\left(\frac{64}{\Lambda_{\mathrm{mg}}^2t_k}\right)^{m+1}s^m.
\end{equation}
Choose a dimensional $\delta_{\mathrm{met}}>0$ so small that
\begin{equation*}
 |\widehat g_s-h(s)|_h\le\delta_{\mathrm{met}}
 \quad\hbox{implies}\quad
 (1-\eps_*/3)\widehat g_s\le h(s)\le(1+\eps_*/3)\widehat g_s.
\end{equation*}
Choose $\Lambda_{\mathrm{mg}}$ so that the first
estimate in \eqref{rf:eq:lee-decay} is at most $\delta_{\mathrm{met}}$ at $s=t_k$.
Pulling back at $y=F_s(x)$ gives the $1\pm\eps_*/3$ estimate at $x$ for
$s\le t_k$.  Since $x\in B_{k,5}$ was arbitrary, \eqref{rf:eq:barrier}, now with
$\eps=\eps_*/3$ and $C_{\mathrm{bar}}\mu\le\frac23\eps_*$, gives exactly
\eqref{rf:eq:ind-equiv} on the new slab over $B_{k,5}$.  This is the required
estimate on the still smaller domain of $(\mathrm{Ind}_{k+1})$.

Moreover, the gauge identity $\partial_sF_s=-\mathcal V_s\circ F_s$ relating the two
flows \cite{Lee26}, in which the DeTurck vector field satisfies
$|\mathcal V_s|_h\le C|\nabla^h\widehat g_s|_h$, turns the second estimate in \eqref{rf:eq:lee-decay}
into
\begin{equation*}
 |\partial_sF_s|_{h(s)}
 \le C t_k^{-(m+1)/2}s^{m/2},\qquad0<s\le t_k,
\end{equation*}
while \cref{rf:prop:lee-inputs}\textup{(ii)} gives
$|\partial_sF_s|_{h(s)}\le Cs^{-1/2}$ on the new slab.  Measuring the
trajectory in $h(t)$, splitting its integral at $t_k$, and using \eqref{rf:eq:growth},
we obtain, for $t_k\le t\le t_{k+1}$,
\begin{equation}\label{rf:eq:disp-int}
 \begin{split}
 \int_0^{t_k}|\partial_sF_s|_{h(t)}\,ds
 &\le Ct^{a_*}t_k^{-(m+1)/2}\int_0^{t_k}s^{m/2-a_*}\,ds
 \le C\sqrt{t_k},\\
 \int_{t_k}^{t}|\partial_sF_s|_{h(t)}\,ds
 &\le C\sqrt{t_k}.
 \end{split}
\end{equation}
The first integral converges by the choice of $m$; on the second one the
factor $(t/s)^{a_*}\le(1+\mu)^{a_*}$ is absorbed into $C$.  By
\eqref{rf:eq:initial-order}, \eqref{rf:eq:disp-int} and the shrinking-balls lemma for $h$,
\begin{equation}\label{rf:eq:disp-conv}
 d_{g(0)}(x,F_t(x))
 \le d_{h(0)}(x,F_t(x))
 \le d_{h(t)}(x,F_t(x))+C(n,A_*)\sqrt t.
\end{equation}
Increasing $\Lambda_{\mathrm{mg}}$ once more makes the right-hand side of \eqref{rf:eq:disp-conv} at most
$\Lambda_{\mathrm{mg}}\sqrt{t_{k+1}}/4$; this proves the first line of \eqref{rf:eq:ind-disp} on the
new slab.

\emph{Covering.}
Fix
\begin{equation*}
 x\in B_{k+1,5},
 \qquad R_k:=\Lambda_{\mathrm{mg}}\sqrt{t_{k+1}},
\end{equation*}
and put $\mathcal B=B_{g(0)}(x,2R_k)$.  The definitions
$\mathfrak r_{k+1}=\mathfrak r_k-6\Lambda_{\mathrm{mg}}\sqrt{t_k}$ and
$t_{k+1}\le(1+\mu)t_k$ show that $\overline{\mathcal B}$ lies both in
$B_{k,3}$ and in the region $B_{k,5}$ where the new displacement estimate
holds.  Thus, for every $0\le t\le t_{k+1}$, $F_t$ is an embedding on a
neighbourhood of $\overline{\mathcal B}$ and
\begin{equation*}
 d_{g(0)}(x,F_t(z))
 \ge d_{g(0)}(x,z)-d_{g(0)}(z,F_t(z))
 \ge\frac74R_k>R_k
 \quad(z\in\partial\mathcal B),
\end{equation*}
whereas $F_t(x)\in B_{g(0)}(x,R_k)$.  Invariance of domain makes
$F_t(\mathcal B)$ open, and compactness of $\overline{\mathcal B}$ gives
$\partial F_t(\mathcal B)\subset F_t(\partial\mathcal B)$.  The connected
ball $B_{g(0)}(x,R_k)$ does not meet this boundary and contains $F_t(x)$;
hence it is contained in $F_t(\mathcal B)$.  This proves the second line of
\eqref{rf:eq:ind-disp} for the data $(\mathfrak r_{k+1},t_{k+1})$. Since
$B_{k+1,0}\subset B_{k,5}\subset B_{k,3}$, properties (I) and (II) hold on
$B_{k+1,0}\times[0,t_{k+1}]$ as well, so $F$ satisfies $(\mathrm{Ind}_{k+1})$.
\end{proof}

\subsubsection{Proof of the Localized Comparison}\label{rf:sec:slab-induction}
It remains to start the induction, to run it, and to pass from the gauge
estimates to the curvature comparison.

\begin{proof}[Proof of \cref{rf:lem:comparison}]
\emph{The Initial Solution.}
Let $D'$ be a smooth domain with $D_{31}\Subset D'\Subset D_{32}$. On $D_{64}$ the
metrics $g(0)$ and $h(0)$ coincide, and so do all their time derivatives
$\partial_t^kg(0)$ and $\partial_t^kh(0)$, since both are Ricci flows of the same
metric there and the time derivatives at $t=0$ are determined locally by the
initial metric. Hence the identity map solves the Ricci--harmonic map heat flow to
every order at $t=0$, and the Dirichlet problem on $D'$ with initial and boundary
value the identity has smooth compatible data. Short-time existence for this
problem is classical \cite{Ham75}; the argument there uses only the uniform parabolicity
of the system and the compactness of $\overline{D'}$ and of a target neighbourhood of it,
so it applies to time-dependent metrics. Let $F$ be the
solution on $\overline{D'}\times[0,t_0]$; it is smooth up to $t=0$ with
$F_0=\mathrm{id}$, and its time derivative is bounded on this compact cylinder, so
$d_{g(0)}(x,F_t(x))\le Ct$ and $dF_t$ is uniformly close to the identity for small
$t$. Decreasing $t_0$ therefore makes $F$ satisfy $(\mathrm{Ind}_0)$ for the data
$(31,t_0)$: (I) and (II) hold because a map $C^1$-close to the identity on the
compact set $\overline{D_{31}}$ is a diffeomorphism onto its image with the stated
metric equivalence, and (III) follows from $Ct_0\le\Lambda_{\mathrm{mg}}\sqrt{t_0}/4$ and the
covering argument at the end of the proof of \cref{rf:lem:slab}. The time $t_0$ may
depend on the initial geometry; the final time will not.

\emph{The Induction.}
Put $t_k:=(1+\mu)^kt_0$, $\mathfrak r_0:=31$ and
$\mathfrak r_{k+1}:=\mathfrak r_k-6\Lambda_{\mathrm{mg}}\sqrt{t_k}$. Since the times are geometric,
\begin{equation*}
 \sum_{j\le k}\sqrt{t_j}\le C(\mu)\sqrt{t_k}.
\end{equation*}
Choose $\tau_L=\tau_L(n,A_*,L)\le\tau'$ so small that, whenever
$t_k\le(1+\mu)\tau_L$, the total decrease $6\Lambda_{\mathrm{mg}}\sum_{j\le k}\sqrt{t_j}$ of the radius
is less than $7$ and $2\Lambda_{\mathrm{mg}}\sqrt{t_k}<4$; then $\mathfrak r_{k+1}>24$ and
$D_8\subset B_{k+1,5}$.
Decrease $t_0$ below $\tau_L$. Repeated application of \cref{rf:lem:slab} gives
$(\mathrm{Ind}_k)$ for every $k\ge1$ with $t_{k-1}<\min\{T,\tau_L\}$ and $t_k<T$:
at each step the lemma applies, since $t_{k-1}<\tau_L\le\tau'$ and
$\mathfrak r_k>24$. If instead $t_k\ge T$, the last step with $t_k$ replaced by any
$t'\in(t_{k-1},T)$ gives $(\mathrm{Ind}_k)$ for the data $(\mathfrak r_k,t')$. Hence,
for every $t<\min\{T,\tau_L\}$, the flow $F$ exists on $D_{24}\times[0,t]$ and
satisfies (I)--(III) for the data of the slab containing $t$: $(31,t_0)$ if
$t\le t_0$, and, if $t_{k-1}<t\le t_k$, either $(\mathfrak r_k,t_k)$ or
$(\mathfrak r_k,t')$ with $t\le t'<T$. In the second case the final time of the data
is at most $(1+\mu)t$; in every case it is at most $(1+\mu)\tau_L$, so the choice of
$\tau_L$ applies to it.

\emph{Conclusion.}
By the choice of $\tau_L$, the ball $D_8$ lies in the region $B_{k,5}$ where
(III) holds. Restricting to $D_8$ gives
$F_t(D_8)\subset D_{16}$ by the displacement estimate, and the
second line of \eqref{rf:eq:ind-disp} with centres in $D_4$ gives
$D_4\subset F_t(D_8)$. For $t>t_0$ the final time of the data is at most $(1+\mu)t$, so the first line
of \eqref{rf:eq:ind-disp} gives the stated $C_{\mathrm{disp}}\sqrt t$ bound.
For $t\le t_0$, the estimate \eqref{rf:eq:lee-decay}, which was derived
from $(\mathrm{Ind}_k)$ alone, applied to the data $(31,t_0)$, and
\eqref{rf:eq:growth} give, for $x\in D_8$,
\begin{equation*}
 \int_0^t|\partial_sF_s(x)|_{h(t)}\,ds
 \le C t_0^{-(m+1)/2}t^{(m+2)/2}\le C\sqrt t.
\end{equation*}
The shrinking-balls lemma for $h$ and \eqref{rf:eq:initial-order} then give the
same $g(0)$-displacement bound. This proves
\eqref{rf:eq:comp-disp}--\eqref{rf:eq:comp-equiv} for a time independent of
the initial local curvature.

\emph{The Ricci--DeTurck Estimates.}
Let $t<\min\{T,\tau_L\}$, with slab data $(\mathfrak r_k,t_k)$. Since
\mbox{$2\Lambda_{\mathrm{mg}}\sqrt{t_k}<4$}, the second line of \eqref{rf:eq:ind-disp} with centres in
$D_8$ gives $D_8\subset F_t(D_{16})$, and $B_{g(0)}(z,4)\subset D_{24}\subset B_{k,0}$
for $z\in D_{16}$. Hence \cref{rf:prop:lee-inputs}\textup{(iii)}, applied at every
$z\in D_{16}$ with $r=1$ and $q=4$, gives
\begin{equation*}
 |(\nabla^h)^j\widehat g|_h\le C_4^{\mathrm{DT}}t^{-j/2}
 \qquad(1\le j\le4)\quad\hbox{on }D_8,
\end{equation*}
after reducing $\tau_L$ below $S_4(n,A_*,\Lambda_1,4)$. On $D_8$ the metric
$\widehat g$ is a Ricci--DeTurck flow relative to $h$ with
$\widehat g(0)=h(0)$, and pushing forward \eqref{rf:eq:ind-equiv} from $D_{16}$ gives
the $\eps_{\mathrm{DT}}$-equivalence there. For $y\in D_2$ the closed ball
$\overline{B_{h(0)}(y,4)}$ lies in $D_8$ because $h(0)\ge g(0)$, so
\cref{rf:prop:lee-inputs}\textup{(iv)}, with $r=1$ and exponent $2L$, gives
\eqref{rf:eq:comp-H} on $D_2$ after a further reduction of $\tau_L$. By the
shrinking-balls lemma for $h$, an $h(t)$-ball of fixed radius about each point of $D_1$
stays inside $D_2$ for small $t$.  One-dimensional covariant interpolation along
geodesic segments in that ball, using $|(\nabla^h)^4H|\le Ct^{-2}$ and followed by
polarization, gives
\begin{equation*}
 t|(\nabla^h)^2H|\le C_Lt^{L/2}.
\end{equation*}
The antisymmetric part of $(\nabla^h)^2H$ is bounded by $|\Rm_h||H|$.
Finally,
\begin{equation*}
 |\Rm_{\widehat g}-\Rm_h|_h
 \le C\bigl(|(\nabla^h)^2H|_h+|\nabla^hH|_h^2
              +|\Rm_h|_h|H|_h\bigr),
\end{equation*}
which proves \eqref{rf:eq:comp-Rm}.  Restoring $\ell$ yields the displayed scale factors.
All constants depend only on $n,A_*,L$.
\end{proof}

\subsection{Completion of the Bootstrap}\label{rf:sec:bootstrap}

The comparison flow of \cref{rf:sec:selected} and the comparison lemma of
\cref{rf:sec:comparison} are now combined. At the centre of the ball where the two initial
metrics agree, $t|\Rm|$ for the given flow exceeds the bound $\beta$ of the comparison
flow only by terms that tend to zero with $t/\ell^2$, so it is at most $\alpha$ on a time
interval proportional to $r^2$.

\begin{proof}[Proof of \cref{rf:prop:bootstrap}]
Let $\beta$, $\eps_P$ and $\delta_I$ be the constants of \cref{rf:sec:constants}, and let
$c_I$, $\delta_{\mathrm{imp}}$ and $\kappa_{\mathrm{imp}}$ be as in \eqref{rf:eq:cI}.  By hypotheses (a), (c) and
(d), \cref{rf:lem:isoperimetry} applies to $g(0)$ with $\eta=\delta_I$.  Hence
\eqref{rf:eq:scal-R} and \eqref{rf:eq:iso} hold for $g(0)$, with $\eta=\delta_I$ and
\begin{equation*}
 R:=c_Ir.
\end{equation*}
The Ricci lower bound \eqref{rf:eq:hyp-ricci} and the local cut-off construction
\cite[Theorem~6.33]{CC96}, in the form used in \cite[proof of Lemma~2.2]{HW22}, give a
smooth function $\chi$ with values in $[0,1]$, supported in $B_{g(0)}(p,3R/4)$ and equal
to one on $B_{g(0)}(p,R/3)$, that satisfies \eqref{rf:eq:cutoff} with a constant
$C_\chi=C_\chi(n)$.  \Cref{rf:lem:completion} therefore gives a domain $\Omega$ and a
complete metric $h_0=e^{2u}g(0)$ on it, with constants $c_H$ and $C_H$ that depend only
on $n$ and $\beta$, and \cref{rf:lem:selected} gives the Ricci flow $h(t)$ of $h_0$ on
$[0,T_H]$, where $T_H=(\eps_Pr_H)^2$ and $r_H=c_HR$.

Put
\begin{equation*}
 \ell:=\frac{R}{512},\qquad N:=\Omega.
\end{equation*}
Then
\begin{equation*}
 D_{64}=B_{g(0)}(p,R/8)
       \Subset B_{g(0)}(p,R/4).
\end{equation*}
Since $u\ge0$ on $N$, the conformal factor is at least one and $h(0)\ge g(0)$ as
quadratic forms on all of $N$; and $u=0$ on $B_{g(0)}(p,R/4)$, so $h(0)=g(0)$ on
$D_{64}$.  This is \eqref{rf:eq:initial-order}.  Moreover
$D_{64}\Subset N\Subset B_{g(0)}(p,3R/4)\Subset U$.

By \eqref{rf:eq:hyp-typeI} and \eqref{rf:eq:selected-bound},
\begin{equation}\label{rf:eq:boot-curv}
 \sup_{D_{64}}|\Rm_{g(t)}|+\sup_N|\Rm_{h(t)}|\le\frac{A+\beta}{t}
 \qquad\bigl(0<t<\min\{T,T_H\}\bigr).
\end{equation}
Fix $T'<T$ and put $T_c=\min\{T',T_H\}$.  Apply
\cref{rf:lem:comparison} to $g$, $h$, and $N$ on $[0,T_c]$, with
$A_*=A+\beta$ and $L=4$; hypothesis \eqref{rf:eq:comp-curv} is \eqref{rf:eq:boot-curv}.
Put $\tau=t/\ell^2$.  When
$\tau<C_{\mathrm{disp}}^{-2}$, \eqref{rf:eq:comp-disp} places $F_t(p)$ in $D_1$.
Diffeomorphism invariance of curvature, \eqref{rf:eq:selected-bound}, and
\eqref{rf:eq:comp-H}--\eqref{rf:eq:comp-Rm} give
\begin{equation}\label{rf:eq:boot-main}
 t|\Rm_{g(t)}|(p)
 \le(1+C_{\mathrm{met}}\tau)\left(\beta+C_4\tau^2\right).
\end{equation}
Here $1+C_{\mathrm{met}}\tau$ is the norm-comparison factor between
$\widehat g$ and $h$, with
$C_{\mathrm{met}}=C_{\mathrm{met}}(n,A,\alpha)$, and $C_4$ is the constant $C_L$ of
\cref{rf:lem:comparison} for $L=4$.

Choose $\tau_0=\tau_0(n,A,\alpha)>0$ so that the right-hand side of
\eqref{rf:eq:boot-main} is at most $\alpha$ whenever $\tau\le\tau_0$, and also
\begin{equation}\label{rf:eq:tau0}
 \tau_0\le\min\left\{1,\tau_4,512^2(\eps_Pc_H)^2,
 \tfrac14C_{\mathrm{disp}}^{-2},(2C_{\mathrm{met}})^{-1}\right\},
\end{equation}
where $\tau_4$ is the constant $\tau_L$ of \cref{rf:lem:comparison} for $L=4$.
This choice depends only on $n,A,\alpha$.  Setting
\begin{equation*}
 \Theta_{\mathrm{imp}}:=\frac{\tau_0c_I^2}{512^2}
\end{equation*}
and using the $T_H$ restriction in \eqref{rf:eq:tau0}, we obtain the comparison estimate for
\begin{equation*}
 0<t<\min\{T',\Theta_{\mathrm{imp}}r^2\}.
\end{equation*}
Letting $T'\uparrow T$ proves \eqref{rf:eq:improved}.
For $t\ge\Theta_{\mathrm{imp}} r^2$, the assumed bound gives
\begin{equation*}
 |\Rm|\le\frac At\le\frac{A}{\Theta_{\mathrm{imp}}}\,r^{-2},
\end{equation*}
which proves \eqref{rf:eq:improved-all}.
\end{proof}

\subsection{Complete Non-Compact Initial Manifolds}\label{rf:sec:noncompact}

We now remove compactness of the initial manifold and prove \cref{rf:thm:noncompact}.
Throughout this section $(M^n,g_0)$ is complete and non-compact,
$\Ric_{g_0}\ge-(n-1)g_0$, and every point of $M$ is a
$(\delta_{\mathrm{nc}},\rho)$-LRRP, where $\delta_{\mathrm{nc}}$ is fixed in
\eqref{rf:eq:delta-nc}. We choose conformal completions whose initial metrics agree with
$g_0$ on an exhaustion of $M$. Uniform estimates for their Ricci flows on one common time
interval give a locally smooth diagonal limit, and the agreement of the initial metrics
with $g_0$ on the exhaustion identifies the initial value of that limit with $g_0$.

\subsubsection{Constants}\label{rf:sec:nc-constants}
The conformal change costs a factor close to one in the isoperimetric constant, so the
collar estimate needs slightly more room than was used in \cref{rf:sec:selected}, and the
constants are fixed accordingly. Let $\beta$, $\eps_P$, $\theta$ and $\delta_I$ be as in
\cref{rf:sec:constants}, and let $w_\theta$ be the conformal profile of
\cref{rf:lem:completion}. Put $\eta_0:=\delta_I/2$, and let
\begin{equation*}
 c_0:=c(\eta_0),\qquad\delta_\flat:=\delta(\eta_0),\qquad\kappa_\flat:=\kappa(\eta_0)
\end{equation*}
be the constants of \cref{rf:lem:isoperimetry} for $\eta=\eta_0$. The constant of
\cref{rf:thm:noncompact} is
\begin{equation}\label{rf:eq:delta-nc}
 \delta_{\mathrm{nc}}:=\min\{\delta_0,\delta_\flat/2\}.
\end{equation}
All these constants depend only on $n$ and $\alpha$.

\subsubsection{The Conformal Exhaustion}\label{rf:sec:nc-exhaustion}
The first lemma constructs the exhaustion and, for every $j$ and every point of
$\Omega_j$, a comparison flow on the local cover. The conformal change pushes the
artificial boundary to infinity, while isoperimetric control on the local cover survives
at the scale $e^{-u_j(x)}r_*$ of the original metric, which becomes a fixed scale for the
conformal metric. This scale tends to zero towards $\partial\Omega_j$, so the construction
uses the LRRP condition at all scales below $\rho$.

\begin{lemma}[Uniform Conformal Exhaustion]\label{rf:lem:nc-exhaustion}
There are relatively compact domains $\Omega_j\Subset M$ with smooth boundary, compact
sets $E_j\subset\Omega_j$ whose union is $M$, and smooth functions $u_j\ge0$ on
$\Omega_j$ such that $h_{j,0}:=e^{2u_j}g_0$ is a complete metric of bounded curvature
on $\Omega_j$ and
\begin{equation}\label{rf:eq:exhaustion}
 u_j=0\ \hbox{on }E_j,\qquad\overline{\Omega_j}\subset\Int E_{j+1}.
\end{equation}
Moreover, there are $r_E,T_E>0$, depending only on $n,\rho,\alpha$, with the following
property. For every $j$ and every $x\in\Omega_j$, let
$\widetilde h_{j,0}:=\pi_x^*h_{j,0}$ on the part of the local cover
$\tB_{g_0}(x,\rho)$ that lies over $\Omega_j$. There are an open set $N_{j,x}$ in this
part and a complete bounded-curvature Ricci flow $h^\sharp_{j,x}(t)$ on
$N_{j,x}\times[0,T_E]$ such that, with $\ell_E:=r_E/512$,
\begin{equation}\label{rf:eq:selected-E}
 \begin{split}
 &\overline{B_{\widetilde h_{j,0}}(\tilde x,64\ell_E)}
       \Subset N_{j,x},\\
 &h^\sharp_{j,x}(0)\ge\widetilde h_{j,0}\quad\hbox{on }N_{j,x},
 \qquad
 h^\sharp_{j,x}(0)=\widetilde h_{j,0}
 \quad\hbox{on }
 B_{\widetilde h_{j,0}}(\tilde x,128\ell_E),\\
 &|\Rm_{h^\sharp_{j,x}(t)}|\le\frac{\beta}{t}
 \quad(0<t\le T_E).
 \end{split}
\end{equation}
\end{lemma}

\begin{figure}[t]
\centering
\includegraphics[width=\textwidth]{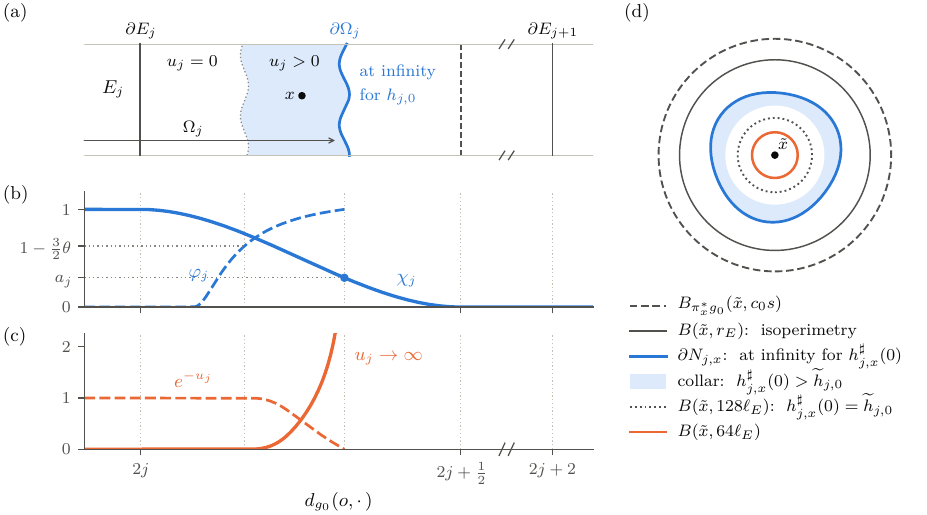}
\caption{The conformal exhaustion of \cref{rf:lem:nc-exhaustion} (schematic, not to
scale). Panels (a)--(c) show the construction along a ray from the base point $o$.
(a)~The regions. The metric $h_{j,0}=e^{2u_j}g_0$ equals $g_0$ where $u_j=0$, in
particular on $E_j$; the boundary of $\Omega_j$ is at infinite distance for it, and
$\overline{\Omega_j}$ lies in the interior of the next compact set $E_{j+1}$. (b)~The
cut-off $\chi_j$, the regular value $a_j$ that defines $\Omega_j$ as a superlevel set, and
the defining function $\varphi_j$. (c)~The conformal exponent
$u_j=w_\theta\circ\varphi_j$, which vanishes where $\varphi_j\le1-\frac32\theta$ and
tends to infinity at $\partial\Omega_j$, and the function $e^{-u_j}$, which is
Lipschitz and vanishes there. (d)~The local cover at a point $x\in\Omega_j$, with the
metric $\widetilde h_{j,0}$; balls without a subscript are taken for this metric. The
second completion takes place at the fixed scale $r_E=512\ell_E$, inside the ball
$B_{\pi_x^*g_0}(\tilde x,c_0s)$ of the lifted original metric, where
$s=e^{-u_j(x)}r_*$; the radius of this ball tends to zero as $x$ approaches
$\partial\Omega_j$. It produces the domain $N_{j,x}$ and the complete metric
$h^\sharp_{j,x}(0)$ that starts the comparison flow.}
\label{rf:fig:exhaustion}
\end{figure}

\begin{proof}
\emph{The Exhaustion.}
Fix $o\in M$ and $j\ge1$, and put $E_j:=\overline{B_{g_0}(o,2j)}$.  We first construct
a cut-off at unit scale, as in \cite[proof of Lemma~2.2]{HW22}.  Let $\{x_k\}$ be a
maximal subset of $E_j$ whose points are pairwise at distance at least $1/4$; it is
finite because $E_j$ is compact.  The balls $B_{g_0}(x_k,1/4)$ cover $E_j$, and the
balls $B_{g_0}(x_k,1/8)$ are disjoint.  By the Bishop--Gromov inequality, every
point of $M$ lies in at most $C(n)$ of the balls $B_{g_0}(x_k,1/2)$.  The
Ricci lower bound gives smooth cutoffs $b_k:M\to[0,1]$ with
\begin{equation*}
 b_k=1\ \hbox{on }B_{g_0}(x_k,1/4),\qquad
 \operatorname{supp}b_k\subset B_{g_0}(x_k,1/2),\qquad
 |\nabla b_k|+|\Delta b_k|\le C(n)
\end{equation*}
\cite[Theorem~6.33]{CC96}.  Fix a smooth nondecreasing
$\upsilon_1:[0,\infty)\to[0,1]$ with $\upsilon_1(0)=0$ and $\upsilon_1=1$ on
$[1,\infty)$, and put $\chi_j:=\upsilon_1\circ\sum_kb_k$.  Then
\begin{equation*}
 \chi_j=1\ \hbox{on }E_j,\qquad
 \operatorname{supp}\chi_j\subset\overline{B_{g_0}(o,2j+\tfrac12)},
 \qquad
 |\nabla\chi_j|+|\Delta\chi_j|\le C(n).
\end{equation*}
Choose a regular value $a_j\in(1/4,1/2)$ of $\chi_j$, let $\Omega_j$ be the component
of $\{\chi_j>a_j\}$ containing $E_j$, and put on $\Omega_j$, with the function
$\upsilon_0$ of the proof of \cref{rf:lem:completion},
\begin{equation*}
 \varphi_j:=\upsilon_0\circ\frac{\chi_j-a_j}{1-a_j},\qquad
 u_j:=w_\theta\circ\varphi_j,\qquad
 h_{j,0}:=e^{2u_j}g_0.
\end{equation*}
Thus $\varphi_j$ vanishes on $E_j$, tends to one at $\partial\Omega_j$, and satisfies
\begin{equation}\label{rf:eq:phi-j}
 |\nabla\varphi_j|+|\Delta\varphi_j|\le C(n).
\end{equation}
As in the proof of \cref{rf:lem:completion}, the last inequality in \eqref{rf:eq:w-collar}
and \eqref{rf:eq:phi-j} show that every curve leaving compact subsets of $\Omega_j$ has
infinite $h_{j,0}$-length, so $h_{j,0}$ is complete, and $h_{j,0}$ has bounded
curvature because $g_0$ and $\varphi_j$ extend smoothly past the compact set
$\overline{\Omega_j}$.  Since
\begin{equation*}
 \overline{\Omega_j}\subset\overline{B_{g_0}(o,2j+\tfrac12)}\subset\Int E_{j+1},
\end{equation*}
this gives \eqref{rf:eq:exhaustion}.  \Cref{rf:fig:exhaustion} shows these regions and
functions.

\emph{Uniform Bounds.}
By \eqref{rf:eq:w-deriv}, \eqref{rf:eq:phi-j}, the conformal scalar-curvature formula, and
the bound ${\Scal_{g_0}\ge-n(n-1)}$, there is a constant $C_E=C_E(n,\alpha)$ such that
\begin{equation}\label{rf:eq:CE}
 |\nabla^{h_{j,0}}u_j|_{h_{j,0}}\le C_E,\qquad
 \Scal_{h_{j,0}}\ge-C_E.
\end{equation}
Indeed the derivatives of $w_\theta$ occur only in the combinations of
\eqref{rf:eq:w-deriv}, and only the Laplacian of $\varphi_j$, not its full Hessian,
enters the scalar-curvature formula.

\emph{The Adaptive Scale.}
Fix $x\in\Omega_j$ and put
\begin{equation*}
 s:=e^{-u_j(x)}r_*,
\end{equation*}
where $r_*=r_*(n,\rho,\alpha)\in(0,1]$ will be fixed below.  From \eqref{rf:eq:CE},
\begin{equation*}
 |\nabla^{g_0}e^{-u_j}|_{g_0}=|\nabla^{h_{j,0}}u_j|_{h_{j,0}}
 \le C_E.
\end{equation*}
So $e^{-u_j}$ is $C_E$-Lipschitz for $d_{g_0}$ along curves in $\Omega_j$, and it
extends continuously by zero to $\partial\Omega_j$.  Hence
$d_{g_0}(x,\partial\Omega_j)\ge e^{-u_j(x)}/C_E$, and for every $c\in(0,1]$ with
$cC_Er_*\le\frac16$,
\begin{equation}\label{rf:eq:osc}
 \overline{B_{g_0}(x,cs)}\Subset\Omega_j,\qquad
 \bigl|e^{u_j(x)-u_j}-1\bigr|\le cC_Er_*\ \hbox{ on }B_{g_0}(x,cs).
\end{equation}

\emph{Isoperimetry on the Cover.}
Require
\begin{equation*}
 r_*\le\min\{\rho/20,\sqrt{\kappa_\flat}/4,1/20\}.
\end{equation*}
Since $s\le r_*$, the closed ball of radius $s/2$ about $\tilde x$ in the local cover
is compact, by \cref{rf:lem:cover}; the Ricci lower bound gives
$\Ric_{\pi_x^*g_0}\ge-(n-1)\kappa_\flat s^{-2}\pi_x^*g_0$; and the LRRP
condition at $x$ gives
\begin{equation*}
 d_{GH}\bigl(B_{\pi_x^*g_0}(\tilde x,2s),B_{\R^n}(0,2s)\bigr)
 <2\delta_{\mathrm{nc}}s\le\delta_\flat s .
\end{equation*}
So \cref{rf:lem:isoperimetry} applies to $\pi_x^*g_0$ at $\tilde x$, with $\eta=\eta_0$
and $r=s$: every regular domain $E\Subset B_{\pi_x^*g_0}(\tilde x,c_0s)$ satisfies
\begin{equation}\label{rf:eq:iso-lift}
 \bigl(\area_{\pi_x^*g_0}\partial E\bigr)^n
 \ge(1-\eta_0)n^n\omega_n
       \bigl(\vol_{\pi_x^*g_0}E\bigr)^{n-1}.
\end{equation}

\emph{The Conformal Metric on the Cover.}
Put $\eps_E:=c_0C_Er_*$ and require $\eps_E\le\frac16$.  By \eqref{rf:eq:osc} with
$c=c_0$, the ball $B_{\pi_x^*g_0}(\tilde x,c_0s)$ lies over $\Omega_j$, its closure is
compact, and on it the lift of $u_j$, still written $u_j$, satisfies
$|e^{u_j(x)-u_j}-1|\le\eps_E$.  Comparing boundary areas and volumes for
$\pi_x^*g_0$ and $\widetilde h_{j,0}=e^{2u_j}\pi_x^*g_0$, we get from
\eqref{rf:eq:iso-lift}, for every regular domain
$E\Subset B_{\pi_x^*g_0}(\tilde x,c_0s)$,
\begin{equation}\label{rf:eq:iso-conf}
 \bigl(\area_{\widetilde h_{j,0}}\partial E\bigr)^n
 \ge\Bigl(\frac{1-\eps_E}{1+\eps_E}\Bigr)^{n(n-1)}(1-\eta_0)n^n\omega_n
 \bigl(\vol_{\widetilde h_{j,0}}E\bigr)^{n-1}.
\end{equation}
Since $\eta_0=\delta_I/2$, we may decrease $r_*$ so that the constant on the right is
at least $(1-\delta_I)n^n\omega_n$.  Set
\begin{equation*}
 r_E:=\frac{c_0r_*}{2(1+\eps_E)},
\end{equation*}
and decrease $r_*$ once more so that $C_Er_E^2\le1$.  A curve from $\tilde x$ to the
boundary of $B_{\pi_x^*g_0}(\tilde x,c_0s)$ has length at least $c_0s$ for
$\pi_x^*g_0$, hence at least $e^{u_j(x)}(1+\eps_E)^{-1}c_0s=2r_E$ for
$\widetilde h_{j,0}$.  Therefore
$B_{\widetilde h_{j,0}}(\tilde x,2r_E)\subset B_{\pi_x^*g_0}(\tilde x,c_0s)$, the
closed $\widetilde h_{j,0}$-ball of radius $r_E$ about $\tilde x$ is compact, and by
\eqref{rf:eq:CE}
\begin{equation}\label{rf:eq:scal-E}
 \Scal_{\widetilde h_{j,0}}\ge-r_E^{-2}
 \quad\hbox{on }B_{\widetilde h_{j,0}}(\tilde x,r_E),
\end{equation}
while by \eqref{rf:eq:iso-conf} every regular domain
$E\Subset B_{\widetilde h_{j,0}}(\tilde x,r_E)$ satisfies \eqref{rf:eq:iso} for the metric
$\widetilde h_{j,0}$, with $\eta=\delta_I$.

\emph{The Cut-off on the Cover.}
The local cut-off construction \cite[Theorem~6.33]{CC96} for the metric $\pi_x^*g_0$,
whose Ricci curvature is at least $-(n-1)$, gives a smooth function $\tilde\chi$ with
values in $[0,1]$ such that
\begin{equation*}
 \tilde\chi=1\ \hbox{on }B_{\pi_x^*g_0}(\tilde x,c_0s/5),
 \qquad
 \operatorname{supp}\tilde\chi\subset
 B_{\pi_x^*g_0}(\tilde x,c_0s/4),
 \qquad
 s|\nabla\tilde\chi|+s^2|\Delta\tilde\chi|\le C,
\end{equation*}
where the derivatives are taken for $\pi_x^*g_0$ and $C=C(n,\alpha)$.
The conformal Laplacian identity
\begin{equation*}
 \Delta_{\widetilde h_{j,0}}\tilde\chi
 =e^{-2u_j}\bigl(\Delta_{\pi_x^*g_0}\tilde\chi
 +(n-2)\langle\nabla u_j,\nabla\tilde\chi\rangle_{\pi_x^*g_0}\bigr),
\end{equation*}
together with \eqref{rf:eq:CE}, \eqref{rf:eq:osc}, and $s=e^{-u_j(x)}r_*$, gives
\begin{equation}\label{rf:eq:cut-bounds}
 r_E|\nabla^{\widetilde h_{j,0}}\tilde\chi|_{\widetilde h_{j,0}}
 +r_E^2|\Delta_{\widetilde h_{j,0}}\tilde\chi|\le C_\chi',
\end{equation}
with $C_\chi'=C_\chi'(n,\alpha)$.  By the two-sided bound for $e^{u_j}$ on
$B_{\pi_x^*g_0}(\tilde x,c_0s)$, the $\widetilde h_{j,0}$-distance from $\tilde x$ to
the complement of $B_{\pi_x^*g_0}(\tilde x,c_0s/5)$ is at least $2r_E/5>r_E/3$, and
every point of $B_{\pi_x^*g_0}(\tilde x,c_0s/4)$ is at $\widetilde h_{j,0}$-distance
at most $c_0r_*/(4(1-\eps_E))<3r_E/4$ from $\tilde x$, because $\eps_E\le\frac16$.
Hence
\begin{equation}\label{rf:eq:cut-incl}
 \tilde\chi=1\ \hbox{on }B_{\widetilde h_{j,0}}(\tilde x,r_E/3),\qquad
 \operatorname{supp}\tilde\chi\subset
 B_{\widetilde h_{j,0}}(\tilde x,3r_E/4).
\end{equation}

\emph{The Comparison Flow.}
By \eqref{rf:eq:scal-E}, the isoperimetric inequality stated after it,
\eqref{rf:eq:cut-bounds} and \eqref{rf:eq:cut-incl}, \cref{rf:lem:completion} applies to the
metric $\widetilde h_{j,0}$, the point $\tilde x$ and the radius $R=r_E$, with
$C_\chi=C_\chi'$.  It gives a domain $N_{j,x}$ with
\begin{equation*}
 B_{\widetilde h_{j,0}}(\tilde x,r_E/4)\Subset N_{j,x}
 \Subset B_{\widetilde h_{j,0}}(\tilde x,3r_E/4)
\end{equation*}
and a complete metric of bounded curvature on $N_{j,x}$ that dominates
$\widetilde h_{j,0}$, agrees with it on $B_{\widetilde h_{j,0}}(\tilde x,r_E/4)$, and
satisfies the hypotheses of \cref{rf:lem:selected} with $r_H=c_H'r_E$, where
\mbox{$c_H'=c_H'(n,\alpha)$} is the constant $c_H$ of \cref{rf:lem:completion} for $C_\chi'$.
Panel~(d) of \cref{rf:fig:exhaustion} shows these regions.
Let $h^\sharp_{j,x}(t)$ be the flow of this metric given by \cref{rf:lem:selected}.  It
exists on $[0,T_E]$ with
\begin{equation*}
 T_E:=(\eps_Pc_H'r_E)^2,
\end{equation*}
and it satisfies the last line of \eqref{rf:eq:selected-E}.  Since
$128\ell_E=r_E/4$, the other two lines hold as well.
\end{proof}

\subsubsection{Uniform Lifetime}\label{rf:sec:nc-lifetime}
The second lemma is the first-time argument of \cref{rf:sec:first-time} for the flows of
the exhaustion metrics, with the flows $h^\sharp_{j,x}$ of \cref{rf:lem:nc-exhaustion} as
comparison flows.

\begin{lemma}[Uniform Lifetime of the Exhaustion Flows]\label{rf:lem:nc-lifetime}
There is $T_E'=T_E'(n,\rho,\alpha)>0$ such that, for every $j$, the complete
bounded-curvature Ricci flow $h_j(t)$ with $h_j(0)=h_{j,0}$ exists on $[0,T_E']$ and
satisfies
\begin{equation}\label{rf:eq:lifetime}
 \sup_{\Omega_j}|\Rm_{h_j(t)}|\le\frac{\alpha}{t}
 \quad(0<t\le T_E').
\end{equation}
\end{lemma}

\begin{proof}
Shi's theorem \cite{Shi89} starts a maximal complete bounded-curvature flow from each
$h_{j,0}$; its lifetime $T_{\max,j}$ may depend on $j$.  Define
\begin{equation*}
 T_j:=\sup\left\{s<T_{\max,j}:
 t\sup_{\Omega_j}|\Rm_{h_j(t)}|\le1
 \text{ for }0<t\le s\right\}.
\end{equation*}
This number is positive for each $j$.

Fix $x\in\Omega_j$ and lift $h_j(t)$ to the part of the local cover
$\tB_{g_0}(x,\rho)$ that lies over $\Omega_j$.  On every interval $[0,T]$ with
$T<\min\{T_j,T_E\}$, the lifted flow and the comparison flow $h^\sharp_{j,x}(t)$ satisfy
the hypotheses of \cref{rf:lem:comparison} with $N=N_{j,x}$, $p=\tilde x$,
$\ell=\ell_E$ and $A_*=1+\beta$: the first line of \eqref{rf:eq:selected-E} gives
$D_{64}\Subset N$, \eqref{rf:eq:initial-order} is its second line, and
\eqref{rf:eq:comp-curv} follows from the definition of $T_j$ and its third line.  As in the proof of \cref{rf:prop:bootstrap}, with
$\tau=t/\ell_E^2$ and $L=4$, this gives
\begin{equation*}
 t|\Rm_{h_j(t)}|(x)\le(1+C_{\mathrm{met}}\tau)\left(\beta+C_4\tau^2\right)
\end{equation*}
for $\tau<C_{\mathrm{disp}}^{-2}$.  Choose $\tau_E=\tau_E(n,\alpha)>0$ so that the
right-hand side is at most $\alpha$ for $\tau\le\tau_E$ and
$\tau_E\le\min\{1,\tau_4,\frac14C_{\mathrm{disp}}^{-2},(2C_{\mathrm{met}})^{-1}\}$, and put
\begin{equation*}
 T_E':=\min\{T_E,\tau_E\ell_E^2\}.
\end{equation*}
Then, for every $x\in\Omega_j$,
\begin{equation}\label{rf:eq:lifetime-improved}
 t|\Rm_{h_j(t)}|(x)\le\alpha
 \quad\bigl(0<t<\min\{T_j,T_E'\}\bigr).
\end{equation}
Suppose $T_j\le T_E'$.  If $T_j<T_{\max,j}$, the strict margin in
\eqref{rf:eq:lifetime-improved} extends the defining inequality of $T_j$ past $T_j$; the
required continuity of the curvature supremum follows from the global derivative
estimates of Shi on every compact subinterval of positive time.  If
$T_j=T_{\max,j}$, \eqref{rf:eq:lifetime-improved} bounds the curvature on the final half
of the maximal interval, and the continuation theorem for complete bounded-curvature
flows extends the flow.  Both alternatives are impossible.  Hence $T_j>T_E'$, and
\eqref{rf:eq:lifetime} follows.
\end{proof}

\subsubsection{The Limit Flow}\label{rf:sec:nc-limit}
The third lemma passes to the limit. It follows the compactness argument of Simon and
Topping \cite[Section~8, proof of Theorem~1.7]{ST21}. Their existence theorem is
three-dimensional and comes with a Ricci lower bound for positive time, but the two
local estimates that it invokes and its diagonal extraction need neither.

\begin{lemma}[Exhaustion Limit]\label{rf:lem:nc-limit}
After passing to a subsequence, the flows $h_j(t)$ converge locally smoothly on
$M\times[0,T_E']$ to a complete Ricci flow $g(t)$ with $g(0)=g_0$ and
\begin{equation*}
 \sup_M|\Rm_{g(t)}|\le\frac{\alpha}{t}
 \qquad(0<t\le T_E').
\end{equation*}
\end{lemma}

\begin{proof}
Fix $i\ge1$ and a radius $r\in(0,1]$ with
$\sup_{E_{i+2}}|\Rm_{g_0}|\le r^{-2}$, for the compact sets
$E_i=\overline{B_{g_0}(o,2i)}$ of the proof of \cref{rf:lem:nc-exhaustion}. Then
\begin{equation*}
 B_{g_0}(y,r)\subset E_{i+2}\qquad(y\in E_{i+1}).
\end{equation*}
For all sufficiently large $j$, $u_j=0$ on a neighbourhood of $E_{i+2}$.
Since $h_{j,0}\ge g_0$, the $h_{j,0}$-ball of radius $r$ about a point
$y\in E_{i+1}$ is then contained in $B_{g_0}(y,r)$, so the curvature of $h_{j,0}$
is at most $r^{-2}$ on it.  The flow $h_j(t)$ is complete, so its balls are
relatively compact in $\Omega_j$, and it satisfies \eqref{rf:eq:lifetime} with
$\alpha<1$.  These are the hypotheses of \cite[Theorem~3.1]{Chen09}, in the form
recorded as \cite[Lemma~8.1]{ST21}, at every centre $y\in E_{i+1}$.
Its conclusion at the centre is
\begin{equation}\label{rf:eq:chen}
 \sup_{E_{i+1}\times[0,T_E']}|\Rm_{h_j(t)}|\le e^{C(n)}r^{-2}
\end{equation}
for all large $j$.  The initial curvature derivatives are uniformly
bounded on $E_{i+2}$, because the initial metric there is exactly $g_0$.
Since $E_{i+1}$ is a compact neighbourhood of $E_i$, \eqref{rf:eq:chen} and the local
derivative estimate \cite[Lemma~8.2]{ST21} give uniform bounds for all curvature
derivatives on $E_i\times[0,T_E']$.  Working in fixed coordinate
charts, as in \cite[Section~8, paragraph following Lemma~8.2]{ST21},
Arzel\`a--Ascoli gives a smooth subsequential limit on
$E_i\times[0,T_E']$.  A diagonal argument over $i$ produces
a smooth Ricci flow $g(t)$ on all of $M\times[0,T_E']$.  Since the initial
metrics are equal to $g_0$ on every fixed compact set for all large $j$,
the convergence includes $t=0$ and $g(0)=g_0$.  Estimate \eqref{rf:eq:lifetime} passes to
the limit.

It remains to check completeness.  The shrinking-balls lemma
\cite[Lemma~2.2]{ST21}, applied to the complete flow $h_j(t)$ with
\eqref{rf:eq:lifetime} and with arbitrarily large radii, gives, for a constant
$C=C(n)$ and all $x,y\in\Omega_j$,
\begin{equation}\label{rf:eq:shrinking}
 d_{h_j(t)}(x,y)
 \ge d_{h_{j,0}}(x,y)-C\sqrt{\alpha t}
 \ge d_{g_0}(x,y)-C\sqrt{\alpha t}.
\end{equation}
For any curve $\gamma$ in $M$ from $x$ to $y$, local
smooth convergence and \eqref{rf:eq:shrinking} give
\begin{equation*}
 \operatorname{length}_{g(t)}(\gamma)
 =\lim_{j\to\infty}\operatorname{length}_{h_j(t)}(\gamma)
 \ge d_{g_0}(x,y)-C\sqrt{\alpha t}.
\end{equation*}
Taking the infimum over such curves gives the same lower bound for
$d_{g(t)}(x,y)$.  Thus every $g(t)$-ball about a fixed base point is
contained in a compact $g_0$-ball, and the Hopf--Rinow theorem implies that $(M,g(t))$
is complete for every $t\in[0,T_E']$.
\end{proof}

The comparison lemma also gives the uniqueness part of \cref{rf:thm:noncompact}.

\begin{corollary}[Complete Type-I Uniqueness]\label{rf:cor:typeI-uniqueness}
Let $g(t)$ and $h(t)$ be complete Ricci flows, smooth on
$M\times[0,T]$, with the same complete smooth initial metric $g_0$.
If $|\Rm_g|+|\Rm_h|\le A/t$ on $M\times(0,T]$ for some $A<\infty$,
then $g(t)=h(t)$ on $[0,T]$.
\end{corollary}

\begin{proof}
Fix a point $p$ in a connected component. For all sufficiently large $R$,
apply \cref{rf:lem:comparison} with $N=U=M$, $\ell=R$, $A_*=\max\{2A,1\}$
and $L=2$; each of the two curvature norms is at most $A/t$. Completeness of $g_0$ gives the compact balls required there,
and $T\le\tau_2R^2$, so the conclusions hold for $0<t<T$. The resulting maps
$F_{R,t}$ are defined on
$B_{g_0}(p,8R)$, satisfy $F_{R,0}=\mathrm{id}$, and obey
\begin{equation*}
 d_{g_0}(x,F_{R,t}(x))\le C\sqrt t,\qquad
 |(F_{R,t}^{-1})^*g(t)-h(t)|_{h(t)}
       \le C(t/R^2)^2
       \quad\hbox{on }B_{g_0}(p,2R).
\end{equation*}
The first estimate keeps the images of each compact source set in a fixed
compact target set, independently of $R$, while \eqref{rf:eq:comp-equiv}
uniformly bounds $|dF_R|$. On a positive-time compact
cylinder, \cref{rf:prop:lee-inputs}\textup{(ii)} and interior parabolic
regularity give uniform derivatives of $F_R$. Diagonal extraction therefore
gives a smooth Ricci--harmonic map heat flow $F$ on $M\times(0,T)$, continuous
at $t=0$ with $F_0=\mathrm{id}$. For any fixed compact set, its images lie
in $B_{g_0}(p,2R)$ for large $R$. Pulling back the second estimate, with
only a uniform factor from \eqref{rf:eq:comp-equiv}, and taking
the limit gives $F_t^*h(t)=g(t)$. Thus $F_t$ is a local isometry, so its
covariant Hessian and tension vanish. The Ricci--harmonic equation gives
$\partial_tF=0$; continuity at $t=0$ forces $F_t=\mathrm{id}$ and hence
$g(t)=h(t)$ for $t<T$, and at $t=T$ by continuity. Apply the argument on each
component.
\end{proof}

\begin{proof}[Proof of \cref{rf:thm:noncompact}]
Set $T_{\mathrm{nc}}:=\min\{T_0,T_E'\}$, and recall that
$\delta_{\mathrm{nc}}\le\delta_0$ by \eqref{rf:eq:delta-nc}.  If $M$ is closed,
existence is \cref{rf:thm:main} and uniqueness is Hamilton's uniqueness theorem for
closed manifolds \cite{Ham82}.  If $M$ is non-compact,
\cref{rf:lem:nc-exhaustion,rf:lem:nc-lifetime,rf:lem:nc-limit} give the flow $g(t)$ on
$[0,T_E']$.  If two complete flows with initial metric $g_0$ satisfy a bound
$|\Rm|\le A/t$, the sum of their curvature norms is at most $2A/t$, and the two flows
coincide by \cref{rf:cor:typeI-uniqueness}.
\end{proof}

\subsection{Local Ricci--Harmonic Comparison Estimates}\label{rf:app:comparison-inputs}
This appendix proves parts~\textup{(i)}--\textup{(iii)} of \cref{rf:prop:lee-inputs}.
\Cref{rf:lem:local-harmonic-existence} constructs the map heat flow on a conformal
completion of the source; its existence time depends only on $n$ and $\Lambda_0$, and the
bounds of the construction that depend on the geometry enter only through the
short-time existence lemma of \cite{BK22}.
\Cref{rf:lem:local-map-derivatives,rf:lem:local-dt-derivatives} give the interior derivative
bounds.

\subsubsection{Interior Existence and Uniform Continuation}
The flow is first constructed on a conformal completion of the source domain, and a
comparison for its energy density makes the existence time uniform.

\begin{lemma}[Interior Existence]\label{rf:lem:local-harmonic-existence}
Under the hypotheses of \cref{rf:prop:lee-inputs}\textup{(i)}, its conclusion
holds with $2\Lambda_0$ in place of $10\Lambda_0$.
\end{lemma}

\begin{proof}
Rescale $r=1$. We use a fixed conformal completion of the
source, as in \cite{Hoc16,Lee26}, and verify the bounds needed for a moving source. On a
component with empty boundary take $G=g$; compactness gives all the qualitative bounds
below. The following construction is for components with boundary.

Put $d(x)=\operatorname{dist}_{g(0)}(x,\partial U)$, with paths in
$\overline U$. For $0<s\le1/320$ and $d(x)>2s$, let $\mathsf k_s$ be a normalized
smooth radial kernel on $|v|<s$ in $T_xU$ and set
\begin{equation*}
 \mathcal M_sd(x)=\int_{T_xU}d(\exp_{x,g(0)}v)\,\mathsf k_s(v)\,dv,
\end{equation*}
the Riemannian convolution of Greene and Wu \cite{GW79}. These geodesics remain in $U$,
and the curvature bound excludes conjugate points along them; since the integral is
taken over the tangent space, no bound on the injectivity radius enters. We use the
estimates
\begin{equation}\label{rf:eq:map-smoothing-bounds}
 |\mathcal M_sd-d|\le s,\qquad |\nabla\mathcal M_sd|\le C_n,\qquad
 |\nabla^2\mathcal M_sd|\le C_n/s.
\end{equation}
The first is the Lipschitz bound on $d$. For the others, work in an orthonormal
trivialization $(\partial_a)$ that is parallel at $x$, and write
\begin{equation*}
 \mathcal E(x,v)=\exp_xv,\qquad \mathcal J_a=\partial_x\mathcal E(\partial_a),\qquad
 \widetilde{\mathcal J}_b=(\partial_v\mathcal E)^{-1}\mathcal J_b.
\end{equation*}
Differentiating a smooth Lipschitz approximation of $d$ twice in the centre
produces
\begin{equation*}
 \nabla^2d(\mathcal J_b,\mathcal J_a)+\langle\nabla d,\mathcal J_{ab}\rangle,
\end{equation*}
where $\mathcal J_{ab}$ is the covariant second centre variation. Since
$\mathcal J_b=\partial_v\mathcal E(\widetilde{\mathcal J}_b)$, the Hessian term is a vertical
derivative along $\widetilde{\mathcal J}_b$ minus
$\langle\nabla d,\nabla^{\mathcal E}_{\widetilde{\mathcal J}_b}\mathcal J_a\rangle$, and integration by parts in $v$
gives
\begin{equation*}
 (\nabla^2\mathcal M_sd)_{ab}
 =-\int\operatorname{div}_v(\mathsf k_s\widetilde{\mathcal J}_b)\,\langle\nabla d,\mathcal J_a\rangle
 +\int\mathsf k_s\,\langle\nabla d,\mathcal J_{ab}
 -\nabla^{\mathcal E}_{\widetilde{\mathcal J}_b}\mathcal J_a\rangle.
\end{equation*}
The Jacobi equation and its first variation, together with the normalization
\begin{equation*}
 |\Rm_{g(0)}|+|\nabla\Rm_{g(0)}|\le2,
\end{equation*}
bound $\mathcal J_a$, $\widetilde{\mathcal J}_b$, $\mathcal J_{ab}$ and
$\nabla^{\mathcal E}_{\widetilde{\mathcal J}_b}\mathcal J_a$, and
$\|\operatorname{div}_v(\mathsf k_s\widetilde{\mathcal J}_b)\|_{L^1}\le C_n/s$. Since $|\nabla d|\le1$,
this proves \eqref{rf:eq:map-smoothing-bounds}, also in the limit for the distance
function itself. Writing $\mathcal M_sd(x)=\int d(z)\,\mathsf K_s(x,z)\,dz$ with a kernel
that is smooth in $x$, because $\mathsf k_s$ vanishes to all orders at $|v|=s$, shows that
$\mathcal M_sd$ is smooth; repeated centre differentiation and vertical integration by
parts, with $\mathcal M_s1=1$, give $|\nabla^q\mathcal M_sd|\le C_qs^{1-q}$, where $C_q$
for $q\ge3$ also depends on the curvature derivatives of the compact smooth extension of
the metric.

Take $\ell_j=2^{-j}/16$ and a fixed smooth nondecreasing function $\upsilon_2$ equal to
zero below $2$ and one above $4$. Define
\begin{equation*}
 \hat d(x)=\sum_{j\ge0}\ell_j\,\upsilon_2\bigl(\mathcal M_{\ell_j/20}d(x)/\ell_j\bigr).
\end{equation*}
Each summand is defined where $d>\ell_j/10$ and extended by zero; it already vanishes
for $d<3\ell_j/2$, so $\hat d$ is smooth. Only a bounded number of scales have non-zero
derivatives at a point, and the smaller-scale tail is locally constant. Consequently
\begin{equation}\label{rf:eq:map-boundary-radius}
 \hat d\asymp\min\{d,1\},\qquad
 |\nabla\hat d|+\hat d|\nabla^2\hat d|\le C_n,\qquad
 |\nabla^q\hat d|\le C_q\hat d^{1-q}.
\end{equation}
Also $\hat d=1/8$ on $U_1$. The constants $C_q$ for $q\ge3$ may depend on the geometry.

The source assumptions give $|\partial_tg|+|\nabla^g\partial_tg|\le C_n$. For a
dimensional $s_0>0$, integration of the metric and connection equations gives
$g(t)\asymp g(0)$ and $|\Gamma(g(t))-\Gamma(g(0))|_{g(0)}\le C_n$ for $t\le s_0$. Thus
\eqref{rf:eq:map-boundary-radius} holds with $g(t)$-derivatives as well. Fix a smooth
non-negative nonincreasing function $\upsilon_3$ with $\upsilon_3(z)=-\log z$ for
$z\le1/2$ and $\upsilon_3(z)=0$ for $z\ge1$, and set
\begin{equation*}
 u=\upsilon_3(64\hat d),\qquad G(t)=e^{2u}g(t).
\end{equation*}
Then $G\ge g$, $G=g$ on $U_1$, and $G(t)$ is complete because $e^u\asymp d^{-1}$ near
$\partial U$. The conformal curvature formula gives
\begin{equation}\label{rf:eq:map-completed-zero-order}
 |\Rm_G|_G+|\partial_tG|_G\le C_n
 \qquad(0\le t\le\min\{T,s_0\}).
\end{equation}
Indeed $e^{-u}|du|+e^{-2u}|\nabla^2u|\le C_n$ and $|\partial_tG|_G=|\partial_tg|_g$.

All the following bounds are finite:
\begin{equation}\label{rf:eq:map-completed-higher}
 \sup_{U\times[0,\min\{T,s_0\}]}
 \bigl(|(\nabla^G)^j\Rm_G|_G+|(\nabla^G)^j\partial_tG|_G\bigr)<\infty
 \qquad(j\ge0).
\end{equation}
For the first term, use derivatives of $u$ through order $j+2$; the second also uses
derivatives of $\partial_tg$ through order $j$. All are available from the compact smooth
extension and \eqref{rf:eq:map-boundary-radius}. Their boundary weights are respectively
$e^{-(j+2)u}(1+\hat d^{-j-2})$ and $e^{-ju}(1+\hat d^{-j})$, which are bounded. Unlike
\eqref{rf:eq:map-completed-zero-order}, these bounds depend on the geometry of the
extension; they enter only through \cite{BK22} below. The target Ricci flow has bounds of
the same kind by hypothesis.

Now $f:(U,G(0))\to(N,h(0))$ is $\sqrt{\Lambda_0}$-Lipschitz, and \cite[Lemma~A.8]{BK22}
gives a solution for the two complete evolving metric families on a short time interval,
whose length may depend on \eqref{rf:eq:map-completed-higher}. For $e_0=|dF|_{G,h}^2$, take
a $G$-orthonormal frame $(v_i)$. The Bochner identity for moving metrics is
\begin{equation*}
\begin{split}
 (\partial_t-\Delta_G)e_0={}&-2|\nabla dF|^2
 -\sum_{i,j}(\partial_tG+2\Ric_G)(v_i,v_j)\,h\bigl(dF(v_i),dF(v_j)\bigr)\\
 &+\sum_i(\partial_th)\bigl(dF(v_i),dF(v_i)\bigr)
 +2\sum_{i,j}\Rm_h\bigl(dF(v_i),dF(v_j),dF(v_j),dF(v_i)\bigr),
\end{split}
\end{equation*}
where $\Rm_h(v_1,v_2,v_2,v_1)>0$ when $v_1,v_2$ span a positively curved plane. By
\eqref{rf:eq:map-completed-zero-order} and $|\Rm_h|\le1$,
$(\partial_t-\Delta_G)e_0\le C_n(e_0+e_0^2)$. The maximum principle for bounded solutions
on these complete families with bounded geometry, as in \cite[proofs of Lemmas~A.6 and~A.7]{BK22},
compares $e_0$ with the solution of $y'=C_n(y+y^2)$, $y(0)=\Lambda_0$. Thus
$e_0\le y<2\Lambda_0$ up to
\begin{equation*}
 T_{\mathrm{ex}}=\min\left\{s_0,\frac1{2C_n}
             \log\frac{2(1+\Lambda_0)}{1+2\Lambda_0}\right\}>0.
\end{equation*}

Let $b=\min\{T,T_{\mathrm{ex}}\}$. The bounds \eqref{rf:eq:map-completed-higher} and those of
the target are uniform on $[0,b]$, and $|dF_{t'}|_{G,h}^2\le2\Lambda_0$ for $t'<b$, so at
every time $t'<b$ \cite[Lemma~A.8]{BK22} supplies the same positive existence time for
the initial map $F_{t'}$. Restart within the preceding existence interval and
concatenate: the two solutions agree at the joining time, and the equation determines
all their time derivatives there from the spatial derivatives of the common map, so the
one-sided jets match and the concatenation is smooth. The energy comparison on every
finite concatenation continues to use the original initial time and gives $e_0\le y$.
Finitely many restarts cover $[0,b)$. Restrict to $U_1$, where $G=g$, and restore the
scale $r$.
\end{proof}

\subsubsection{Covariant Map Derivatives}
The interior derivative bounds follow from a Shi-type induction on the covariant
derivatives of the map.

\begin{lemma}[Local Map Derivatives]\label{rf:lem:local-map-derivatives}
Let $g(t)$ be a Ricci flow on an open $n$-manifold $M$ and $h(t)$ a complete
Ricci flow on an $n$-manifold $N$. Suppose $B_{g(0)}(x_0,4r)\Subset M$ and
\[
 \sup_{B_{g(0)}(x_0,4r)}|\Rm_g|+\sup_N|\Rm_h|\le A/t.
\]
If $F$ is a smooth Ricci--harmonic map heat flow on this source cylinder with
$|dF|^2\le\Lambda_0$, then, for every integer $q\ge0$,
\[
 |\nabla^jdF|^2\le C_q^{\mathrm{map}}t^{-j}\quad(0\le j\le q)
\]
on $B_{g(0)}(x_0,r)$ for $0<t\le\min\{T,S_qr^2\}$, where $C_q^{\mathrm{map}},S_q$
depend only on $n,A,\Lambda_0,q$.
\end{lemma}

\begin{proof}
Rescale $r=1$ and enlarge $A$ to $\max\{A,1\}$. Local Shi estimates give, for a common
smaller time and through order $q+1$,
\begin{equation*}
 |\nabla^a\Rm_g|\le C_at^{-1-a/2}\quad\hbox{on }B_{g(0)}(x_0,3);
\end{equation*} by the shrinking-balls lemma
\cite[Lemma~2.2]{ST21}, the parabolic neighbourhoods of radius $c\sqrt t$ used on
$[t/2,t]$ stay inside $B_{g(0)}(x_0,4)$. Complete Shi estimates give the same bounds for
all target curvature derivatives.

Put $e_k=|\nabla^kdF|^2$, and let $\mathcal D_t$ denote target covariant differentiation
along $F$ in time, leaving the source indices unchanged. For a tensor $\mathcal T$ with covariant
source slots and one target slot,
\begin{equation}\label{rf:eq:evolving-covariant-norm}
 \partial_t|\mathcal T|^2=2\langle\mathcal D_t\mathcal T,\mathcal T\rangle
 +2\sum_a\langle\Ric_g^{\sharp(a)}\mathcal T,\mathcal T\rangle
 -2\langle\Ric_h^\sharp\mathcal T,\mathcal T\rangle,
\end{equation}
the sum running over the source slots. The Ricci endomorphisms act on the indicated
source slot and on the target slot, respectively; the target tensor is pulled back along
$F$. Spatial metric compatibility therefore gives
\begin{equation}\label{rf:eq:ek-evolution}
 (\partial_t-\Delta)e_k
 =2\langle(\mathcal D_t-\Delta)\nabla^kdF,\nabla^kdF\rangle-2e_{k+1}
       +\mathcal Q_k,\qquad
 |\mathcal Q_k|\le C_kt^{-1}e_k,
\end{equation}
where $\mathcal Q_k$ collects the Ricci terms of \eqref{rf:eq:evolving-covariant-norm}.

Commuting $\nabla$ through the tension field, and using that $\partial_t\Gamma(g)$ is a
universal contraction of $g^{-1}$ and $\nabla\Ric_g$, together with the target analogue,
gives universal contractions of the forms
\begin{equation}\label{rf:eq:commutator-types}
 \nabla^a\Rm_g\otimes\nabla^bdF,\quad a+b=k;\qquad
 (\widetilde\nabla^a\Rm_h)\circ F\otimes
 \nabla^{j_1}dF\otimes\cdots\otimes\nabla^{j_m}dF,
 \quad a+\textstyle\sum_lj_l=k.
\end{equation}
The number of factors depends only on $k$, and the factors $dF$ are bounded.
These relations follow inductively: a derivative either raises a curvature
order or a jet index, with an extra undifferentiated $dF$ when it
differentiates target curvature. The moving-target connection terms have
the same relations. Consequently, if $|\nabla^jdF|\le Ct^{-j/2}$ for $j<k$,
\[
 |(\mathcal D_t-\Delta)\nabla^kdF|\le Ct^{-1}|\nabla^kdF|+Ct^{-(k+2)/2}.
\]
Terms containing $\nabla^kdF$ have curvature order zero and only $dF$ as other
jet factors; all other terms are controlled by the lower jets.
For $k=0$ the same bound follows directly using $|dF|^2\le\Lambda_0$.

Suppose the bounds hold through $e_k$. Equations
\eqref{rf:eq:ek-evolution}--\eqref{rf:eq:commutator-types} and Young's
inequality imply
\begin{align}
 (\partial_t-\Delta)e_k
   &\le-2e_{k+1}+Ct^{-k-1},\label{rf:eq:ek}\\
 (\partial_t-\Delta)e_{k+1}
   &\le-2e_{k+2}+Ct^{-1}e_{k+1}+Ct^{-k-2}.
       \label{rf:eq:ekplus}
\end{align}

For the localization take $r_k=3-k/(q+1)$.
At the step from $k$ to $k+1$, centre a cut-off at
$z\in B_{g(0)}(x_0,r_{k+1})$ and put
$s=(r_k-r_{k+1})/2$. A nonincreasing distance cut-off
\[
 \Upsilon=\upsilon_4\bigl((d_{g(t)}(\cdot,z)+c(n)\sqrt{At})/s\bigr)
\]
has compact support in $B_{g(0)}(z,s)\Subset B_{g(0)}(x_0,r_k)$,
equals one at $z$ for a common smaller time, and satisfies
\[
 |\nabla\Upsilon|^2\le C\Upsilon,\qquad
 (\partial_t-\Delta)\Upsilon\le C\Upsilon .
\]
Here $\upsilon_4$ is a fixed smooth nonincreasing function equal to one below $1/2$ and
zero above $1$, with $|\upsilon_4'|^2\le C\upsilon_4$ and $\upsilon_4''\ge-C\upsilon_4$.
The shifted distance is a parabolic supersolution on the transition
annulus, and by the shrinking-balls lemma the support stays in $B_{g(0)}(z,s)$. At
cut points the calculation is understood in the barrier sense.

Consider $\Xi=\mathsf B t^ke_k+t^{k+1}\Upsilon e_{k+1}$.
At a positive-time interior maximum,
\eqref{rf:eq:ek}--\eqref{rf:eq:ekplus} give
\[
 0\le(-2\mathsf B+C)t^ke_{k+1}
      -2t^{k+1}\Upsilon e_{k+2}+C(\mathsf B+1)t^{-1}
      -2t^{k+1}\langle\nabla\Upsilon,\nabla e_{k+1}\rangle .
\]
The last term is at most
$t^{k+1}\Upsilon e_{k+2}
 +Ct^{k+1}(|\nabla\Upsilon|^2/\Upsilon)e_{k+1}$.
For $t\le1$ and large $\mathsf B$, it follows that
$t^{k+1}e_{k+1}\le C$ at the maximum.
On the lateral boundary $\Upsilon=0$ and $t^ke_k$ is already bounded.
At $t=0$ the second term vanishes; the first is bounded by $\mathsf B\Lambda_0$
when $k=0$ and vanishes otherwise, since $F$ is smooth on the closed
cylinder. Thus the maximum principle bounds $\Xi$, and
$\Upsilon(z,t)=1$ gives $e_{k+1}(z,t)\le Ct^{-k-1}$.
The constants are uniform in $z$. Induction reaches
$B_{g(0)}(x_0,r_q)\supset B_{g(0)}(x_0,1)$.
Restoring $r$ proves the statement.
\end{proof}

\subsubsection{Ricci--DeTurck Derivatives}\label{rf:app:deturck-derivatives}
The derivative bounds for the Ricci--DeTurck flow are the map derivative bounds of
\cref{rf:lem:local-map-derivatives} in disguise: the gauge identity below expresses the
covariant derivatives of $\widehat g$ through the Hessian of $F$.

\begin{lemma}[DeTurck Derivatives from Map Derivatives]\label{rf:lem:local-dt-derivatives}
Let $g(t)$, $h(t)$, $x_0$, $r$, $\Lambda_0$ and $F$ be as in
\cref{rf:lem:local-map-derivatives}, and suppose that every $F_t$ is a
diffeomorphism onto its image with
$(1-\eps)g(t)\le F_t^*h(t)\le(1+\eps)g(t)$ on $B_{g(0)}(x_0,4r)$ for some
$\eps\le\frac12$. Then $\widehat g=(F_t^{-1})^*g$ satisfies, for every integer
$q\ge1$,
\[
 |(\nabla^h)^j\widehat g|_h\le C_q^{\mathrm{DT}}t^{-j/2}\qquad(1\le j\le q)
\]
on $F_t(B_{g(0)}(x_0,r))$ for $0<t\le\min\{T,S_qr^2\}$, where
$C_q^{\mathrm{DT}}$ depends only on $n,A,\Lambda_0,q$ and $S_q$ is the constant of
\cref{rf:lem:local-map-derivatives}.
\end{lemma}

\begin{proof}
Fix $t$ and write $\check h:=F_t^*h(t)$, a metric on $B_{g(0)}(x_0,4r)$ for which
$F_t$ is an isometry onto its image, and $g=F_t^*\widehat g$. Pull-back commutes
with covariant differentiation, so
\[
 F_t^*\bigl((\nabla^h)^j\widehat g\bigr)=(\nabla^{\check h})^jg,
 \qquad
 \bigl|(\nabla^h)^j\widehat g\bigr|_h\circ F_t
 =\bigl|(\nabla^{\check h})^jg\bigr|_{\check h}.
\]
The two connections on the source are related through the Hessian of $F_t$:
since $dF_t$ intertwines $\nabla^{\check h}$ with $\nabla^h$, for vector fields $v_1,v_2$ on
the source,
\[
 \nabla^{\check h}_{v_1}v_2-\nabla^g_{v_1}v_2=dF_t^{-1}\bigl(\nabla dF_t(v_1,v_2)\bigr),
 \qquad
 \nabla dF_t(v_1,v_2):=\nabla^h_{dF_tv_1}(dF_tv_2)-dF_t(\nabla^g_{v_1}v_2).
\]
Because $\nabla^gg=0$, this gives
\[
 (\nabla^{\check h}_{v_1}g)(v_2,v_3)
 =-g\bigl(dF_t^{-1}\nabla dF_t(v_1,v_2),v_3\bigr)-g\bigl(v_2,dF_t^{-1}\nabla dF_t(v_1,v_3)\bigr),
\]
so $\nabla^{\check h}g$ is a universal contraction of $g$, $dF_t^{-1}$ and
$\nabla dF_t$. Differentiating again, with the connection $\nabla^{\check h}$ written
as ${\nabla^g+dF_t^{-1}\nabla dF_t}$ and with
${\nabla(dF_t^{-1})=-dF_t^{-1}\circ\nabla dF_t\circ dF_t^{-1}}$, shows inductively that
$(\nabla^{\check h})^jg$ is a sum of universal contractions of $g$, $dF_t^{-1}$ and
tensors $\nabla^{i_1}dF_t,\dots,\nabla^{i_m}dF_t$ with $i_1+\dots+i_m=j$ and all
$i_l\ge1$: each differentiation raises the order of one factor by one or
introduces a factor $\nabla dF_t$. The metric equivalence bounds $dF_t^{-1}$ in operator
norm by $(1-\eps)^{-1/2}$, hence in tensor norm by $\sqrt n(1-\eps)^{-1/2}$, and makes
the norms taken with $g$ and with $\check h$ comparable, so \cref{rf:lem:local-map-derivatives} gives, on $B_{g(0)}(x_0,r)$ and
for $0<t\le\min\{T,S_qr^2\}$,
\[
 \bigl|(\nabla^{\check h})^jg\bigr|_{\check h}
 \le C\sum\prod_l|\nabla^{i_l}dF_t|
 \le C\sum\prod_l(C_q^{\mathrm{map}})^{1/2}t^{-i_l/2}
 \le C_q^{\mathrm{DT}}t^{-j/2},
\]
the sums running over the finitely many contractions. Composing with $F_t^{-1}$
proves the lemma.
\end{proof}

\section{Smooth Affine Fibrations after Ricci-Flow Regularization}\label[appendix]{rf:sec:rong}

At a fixed positive time, \cref{rf:thm:main} and the distortion estimate of
\cite[Appendix~A]{HW22} reduce the smooth collapsing problem to bounded-curvature
theory.

LRRP is as in \cref{rf:def:lrrp}, and the threshold $\delta_{\mathrm{aff}}$ below is distinct
from the threshold $\delta_1(n,N)$ of \cref{thm:main}. We write $\Fr^{\mathrm{orth}}(M,g)$ for the
Riemannian orthonormal frame bundle, or $\Fr^{\mathrm{orth}}(M)$ when the metric is understood, to
distinguish it from the frame spaces $\Fr_\rho$ of the main text.
Write $\dis f=\sup_{x,y}|d(f(x),f(y))-d(x,y)|$ for the distortion of a map;
Gromov--Hausdorff approximations use the convention of \cref{m:fr:sec:conv}.

\subsection{The Fibration Theorem for Bounded Curvature}\label{rf:sec:cfg}
We use the theorem of Fukaya \cite{Fuk87,Fuk89} and
Cheeger--Fukaya--Gromov \cite{CFG92} in the form of \cite[Theorem~1.13]{Ron26}.  For a metric $\bar g$, equip
$\Fr^{\mathrm{orth}}(M)$ with its canonical metric, so that $O(n)$ acts isometrically and
$\Fr^{\mathrm{orth}}(M)/O(n)=(M,d_{\bar g})$.

\begin{theorem}[Nilpotent Fibrations under Bounded Curvature]\label{rf:thm:cfg}
Let $(M_i,\bar g_i)$ be closed $n$-manifolds with $|\sec_{\bar g_i}|\le1$
converging in the Gromov--Hausdorff topology to a compact metric space
$\bar X$ with $\diam\bar X>0$ and $\dim\bar X<n$.  After passing to a subsequence,
$\bar X=\bar Y/O(n)$ for a closed $C^{1,\beta}$ manifold $\bar Y$, for some $\beta\in(0,1)$,
equipped with the invariant Riemannian metric supplied by the cited theorem and an
isometric $O(n)$-action. For $i$ large
there is an $O(n)$-equivariant fibre bundle $\hat f_i:\Fr^{\mathrm{orth}}(M_i)\to \bar Y$ whose fibre is a
nilmanifold and whose structure group is affine, inducing a continuous map $\bar f_i$
in the diagram
\begin{equation}\label{rf:eq:diagram}
 \begin{array}{ccc}
 \Fr^{\mathrm{orth}}(M_i) & \xrightarrow{\ \hat f_i\ } & \bar Y\\[3pt]
 \Big\downarrow & &
 \Big\downarrow\rlap{$\scriptstyle/O(n)$}\\[3pt]
 M_i & \xrightarrow{\ \bar f_i\ } & \bar X=\bar Y/O(n),
 \end{array}
\end{equation}
Every fibre of $\bar f_i$ is a positive-dimensional infranilmanifold.  Moreover, for
some $\eps_i\to0$,
\begin{enumerate}[label=\textup{(\roman*)}]
\item both $\hat f_i:\Fr^{\mathrm{orth}}(M_i)\to \bar Y$ and
      $\bar f_i:(M_i,d_{\bar g_i})\to\bar X$ are $\eps_i$-Gromov--Hausdorff
      approximations;
\item $e^{-\eps_i}|\xi|\le|d\hat f_i(\xi)|\le e^{\eps_i}|\xi|$ for $\xi$
      orthogonal to a fibre;
\item the fibres of $\hat f_i$ have second fundamental form at most $c(n)$.
\end{enumerate}
\end{theorem}

\subsection{The Precise Statement}\label{rf:sec:rong-statement}
The parameter $\vartheta$ below governs the H\"older exponents relating the two limits.

\begin{theorem}[Singular Nilpotent Fibrations, Precise Form]\label{rf:thm:rong}
Fix $n\ge2$, $\rho>0$ and $\vartheta\in(0,1)$, and put
\begin{equation*}
 \alpha:=\frac{\vartheta}{2\sqrt n\,(n-1)},\qquad
 \delta_{\mathrm{aff}}(n,\vartheta):=\delta_0(n,\alpha),
\end{equation*}
where $\delta_0$ is the constant of \cref{rf:thm:main}.  There is
$\bar t=\bar t(n,\rho,\vartheta)\in(0,1]$ with the following property.  Let
$(M_i,g_i)$ be closed $n$-manifolds with $\Ric_{g_i}\ge-(n-1)g_i$ and
uniform $(\delta_{\mathrm{aff}}(n,\vartheta),\rho)$-LRRP, converging in the Gromov--Hausdorff
topology to a compact metric space $X$ of positive diameter, with
$\vol_{g_i}(M_i)\to0$.  Let $g_i(t)$ be the Ricci flow with $g_i(0)=g_i$,
and put $\bar g_i:=(\alpha/\bar t)\,g_i(\bar t)$.  Then every subsequence
has a further subsequence along which the following hold.
\begin{enumerate}[label=\textup{(\alph*)}]
\item The sequence $(M_i,\bar g_i)$ satisfies the hypotheses of
      \cref{rf:thm:cfg} for some limit $\bar X$.  Let $\bar Y$, $\hat f_i$, and
      $\bar f_i:M_i\to\bar X=\bar Y/O(n)$ be as in that theorem.
\item There is a homeomorphism $\Phi:X\to\bar X$, the limit of the identity
      maps $(M_i,d_{g_i})\to(M_i,d_{\bar g_i})$, such that $\Phi$ and
      $\Phi^{-1}$ are H\"older continuous, with exponents $1-\vartheta$ and
      $1/(1+\vartheta)$ and with constants depending only on
      $n,\rho,\vartheta$, and $\diam X$.
\item The maps
      \begin{equation*}
       f_i:=\Phi^{-1}\circ\bar f_i:M_i\longrightarrow X
      \end{equation*}
      are continuous and surjective, every fibre of $f_i$ is an
      infranilmanifold of positive dimension, and $f_i$ is an
      $\eps_i'$-Gromov--Hausdorff approximation from $(M_i,d_{g_i})$ to
      $X$, with $\eps_i'\to0$.
\end{enumerate}
Consequently, for every $\eps>0$, all but finitely many $M_i$ of the original sequence
carry a map $f_i=\Phi^{-1}\circ\bar f_i$ of this kind that is an
$\eps$-Gromov--Hausdorff approximation from $(M_i,d_{g_i})$ to $X$, with $\hat f_i$
satisfying \textup{(i)--(iii)} of \cref{rf:thm:cfg} for $\eps$ in place of $\eps_i$.
Here $\bar Y$, $\bar X$, and $\Phi$ may depend on $i$, and the constants in
\textup{(b)} do not.
\end{theorem}

Thus $f_i$ fits into \eqref{rf:eq:diagram}, with $X$ identified with $\bar Y/O(n)$ through
$\Phi$.  The estimates of \cref{rf:thm:cfg} use the canonical metric of $\bar g_i$; the
$g_i$- and $\bar g_i$-frame bundles are $O(n)$-equivariantly diffeomorphic via
$b\mapsto b\,(b^*\bar g_i)^{-1/2}$. For $\vartheta=\frac12$ the theorem gives the following form
for the whole sequence.

\begin{corollary}[Singular Nilpotent Fibrations]\label{rf:thm:fibration}
For every $n\ge2$ there is $\delta_{\mathrm{aff}}=\delta_{\mathrm{aff}}(n)>0$ with the following property. Let
$\rho>0$, and let $(M_i,g_i)$ be closed $n$-manifolds with $\Ric_{g_i}\ge-(n-1)g_i$ and
uniform $(\delta_{\mathrm{aff}},\rho)$-LRRP. Suppose that they converge in the Gromov--Hausdorff
topology to a compact metric space $X$ of positive diameter, and that
$\vol_{g_i}(M_i)\to0$. Then there are $\eps_i\to0$ and, for $i$ large, continuous
surjective maps $f_i:M_i\to X$ such that
\begin{enumerate}[label=\textup{(\roman*)}]
\item $f_i$ is an $\eps_i$-Gromov--Hausdorff approximation from $(M_i,d_{g_i})$ to $X$;
\item every fibre of $f_i$ is an infranilmanifold of positive dimension;
\item $f_i$ is a singular nilpotent fibration: it is induced by an $O(n)$-equivariant
      fibre bundle, with nilmanifold fibre and affine structure group, from the orthonormal frame
      bundle $\Fr^{\mathrm{orth}}(M_i)$ onto a closed manifold $\bar Y_i$ with $\bar Y_i/O(n)$ homeomorphic to $X$.
\end{enumerate}
\end{corollary}

\emph{Comparison with \cref{thm:main}.} The pointed condition \eqref{rf:eq:lrrp} with error
$\delta$ implies \cref{m:fr:def:localcovers} with error $2\delta$ and the \emph{same}
covering radius $\rho$: for $0<t\le\rho/3$, a correspondence of distortion less than
$2\delta t$ gives a map with that distortion and $2\delta t$-dense image. For the closed
manifolds of \cref{rf:thm:rong,rf:thm:fibration}, normalized Riemannian volume gives
$\RCD(-(n-1),n)$ spaces; compactness upgrades their Gromov--Hausdorff convergence to
measured Gromov--Hausdorff convergence along a subsequence \cite[Lemma~3.32]{GMS15}. To
apply \cref{hyp:sequence}, one must still check its curvature-scale bound
$|\kappa|\rho^2\le1$ and the smallness threshold for the converted error. For the lower
bound $\kappa=-(n-1)$ this requires $(n-1)\rho^2\le1$. Rescaling does not remove that
requirement, and decreasing the covering radius changes the universal cover; its
Reifenberg property is not supplied by this comparison. The results of these appendices
retain their own hypotheses for arbitrary $\rho$ and are not restricted by this additional
bound.

\subsection{Uniform Bi-H\"older Distortion}\label{rf:sec:distortion}
The Ricci flow enters only through the following comparison of distances at one fixed time.

\begin{lemma}[Uniform Bi-H\"older Distortion]\label{rf:lem:distortion}
Let $n,\rho,\vartheta,\alpha$ and $\delta_{\mathrm{aff}}(n,\vartheta)$ be as in \cref{rf:thm:rong}. There are
$\bar t,\bar r\in(0,1]$ and $\bar C\ge1$, depending only on
$n,\rho,\vartheta$, with the following property.  Let $(M^n,g_0)$ be closed
with $\Ric_{g_0}\ge-(n-1)g_0$ and uniform $(\delta_{\mathrm{aff}}(n,\vartheta),\rho)$-LRRP, let $g(t)$
be its Ricci flow, and put $\bar g:=(\alpha/\bar t)\,g(\bar t)$,
$d:=d_{g_0}$, and $\bar d:=d_{\bar g}$.  Then $|\sec_{\bar g}|\le1$,
$\vol_{\bar g}(M)\le\bar C\vol_{g_0}(M)$, and for all $x,y\in M$
\begin{align}
 \bar d(x,y)&\le\bar C\,d(x,y)^{1-\vartheta}
   &&\text{if }d(x,y)<\bar r,                          \label{rf:eq:dist-upper}\\
 d(x,y)&\le\bar C\,\bar d(x,y)^{1/(1+\vartheta)}
   &&\text{if }\bar d(x,y)<\bar r,                     \label{rf:eq:dist-lower}\\
 \bar d(x,y)&\le\bar C\bigl(d(x,y)+1\bigr).           \label{rf:eq:dist-coarse}
\end{align}
\end{lemma}

\begin{proof}
Put $T:=\min\{T_0,1\}$, where $T_0=T_0(n,\rho,\alpha)$ is given by
\cref{rf:thm:main}.  Its curvature estimate implies
\begin{equation*}
 |\Ric_{g(t)}|_{g(t)}\le\frac{\sqrt n\,(n-1)\alpha}{t}
 =\frac{2(n-1)\alpha'}{t}\quad(0<t\le T),
 \qquad
 \alpha':=\frac{\vartheta}{4(n-1)}<\frac1{2(n-1)}.
\end{equation*}
The scalar lower bound $\Scal_{g_0}\ge-n(n-1)$ persists by the maximum principle, since
$M$ is closed, and Bishop--Gromov and
Buser's inequality \cite{Bus82} give dimensional doubling and Poincar\'e constants on
balls of radius at most $10$.  Thus \cite[Theorem~A.1]{HW22}, applied with arbitrary
base point and coefficient $\alpha'$, gives \mbox{$C_D=C_D(n)\ge1$} and
$T_D=T_D(n,T)\in(0,T]$ such that, for $r_D:=\sqrt{T_D}$ and $d_t:=d_{g(t)}$,
\begin{equation}\label{rf:eq:hw-distortion}
 C_D^{-1}d_0(x,y)^{1+\vartheta}\le d_t(x,y)
 \le C_Dd_0(x,y)^{1-\vartheta}
 \qquad\bigl(d_0(x,y)<r_D,\ 0\le t\le T_D\bigr).
\end{equation}
Here $4(n-1)\alpha'=\vartheta$; we use the estimate obtained at the end of the proof of
\cite[Theorem~A.1]{HW22}.

Set $\bar t:=T_D$, so that
$\bar d=(\alpha/\bar t)^{1/2}d_{\bar t}$ and $|\sec_{\bar g}|\le1$.
After taking $\bar r\le r_D$ and enlarging $\bar C$, the upper half of
\eqref{rf:eq:hw-distortion} gives \eqref{rf:eq:dist-upper}.  Subdividing a
minimal $g_0$-geodesic into $k\le2d(x,y)/r_D+1$ arcs of length at most $r_D/2$ gives
$d_{\bar t}(x,y)\le kC_Dr_D^{1-\vartheta}$ and hence \eqref{rf:eq:dist-coarse}, after
enlarging $\bar C$.

Put $c_D:=C_D^{-1}(r_D/2)^{1+\vartheta}$.  If $d_{\bar t}(x,y)<c_D$ but
$d(x,y)\ge r_D$, continuity of $d(x,\cdot)$ along a minimizing
$g(\bar t)$-geodesic gives a point $z$ with $d(x,z)=r_D/2$, and
\eqref{rf:eq:hw-distortion} gives the contradiction
$d_{\bar t}(x,y)\ge d_{\bar t}(x,z)\ge c_D$.  Thus the lower half of
\eqref{rf:eq:hw-distortion} applies to $x,y$ and yields \eqref{rf:eq:dist-lower} after
shrinking $\bar r\le(\alpha/\bar t)^{1/2}c_D$ and enlarging $\bar C$.

Finally, since $\Scal_{g(t)}\ge-n(n-1)$,
\begin{equation*}
 \frac{d}{dt}\vol_{g(t)}(M)=-\int_M\Scal_{g(t)}\,d\vol_{g(t)}
 \le n(n-1)\vol_{g(t)}(M).
\end{equation*}
Thus $\vol_{\bar g}(M)\le(\alpha/\bar t)^{n/2}e^{n(n-1)}\vol_{g_0}(M)$, and all
constants depend only on $n,\rho,\vartheta$.
\end{proof}

\subsection{Proofs of the Smooth Fibration Theorems}\label{rf:sec:rong-proof}
We now identify the two limits and transport the fibrations of the smoothed metrics.

\begin{proof}[Proof of \cref{rf:thm:rong}]
Let $\bar t,\bar r,\bar C$ be the constants of \cref{rf:lem:distortion}.  Write
$d_i:=d_{g_i}$ and $\bar d_i:=d_{\bar g_i}$, and fix $\eps_i$-Gromov--Hausdorff
approximations $\psi_i:(M_i,d_i)\to X$ with $\eps_i\to0$.  We start from a
given subsequence and pass to further subsequences without changing the
notation.

\emph{The Smoothed Limit.}  By \cref{rf:lem:distortion},
$|\sec_{\bar g_i}|\le1$, and \eqref{rf:eq:dist-coarse} bounds the diameter of
$(M_i,\bar g_i)$ by $\bar C(\diam X+2)$ for $i$ large.  By Gromov's
precompactness theorem, a subsequence of $(M_i,\bar d_i)$ converges to a
compact metric space $\bar X$.  Since
$\vol_{\bar g_i}(M_i)\le\bar C\vol_{g_i}(M_i)\to0$, the sequence collapses
and $\dim\bar X\le n-1$ \cite{CC97}.  Choose $x_i,y_i\in M_i$
with $d_i(x_i,y_i)\to\diam X$.  By \eqref{rf:eq:dist-lower},
\begin{equation*}
 \bar d_i(x_i,y_i)\ge
 \min\bigl\{\bar r,\,(d_i(x_i,y_i)/\bar C)^{1+\vartheta}\bigr\},
\end{equation*}
so $\diam\bar X>0$.  Thus \cref{rf:thm:cfg} applies and gives
(a).  In particular the maps $\bar f_i$ are continuous, surjective
$\eps_i$-Gromov--Hausdorff approximations from $(M_i,\bar d_i)$ to $\bar X$, after
enlarging $\eps_i$ to leave strict slack.

\emph{The Homeomorphism.}  For $z\in X$ choose $x_i(z)\in M_i$ with
$d_X(\psi_i(x_i(z)),z)<\eps_i$, so that
\begin{equation}\label{rf:eq:lift-points}
 \bigl|d_i(x_i(z),x_i(z'))-d_X(z,z')\bigr|<3\eps_i
 \qquad(z,z'\in X).
\end{equation}
On a countable dense set $X'\subset X$, compactness and a diagonal subsequence define
\begin{equation*}
 \Phi(z):=\lim_i\bar f_i(x_i(z)).
\end{equation*}
Passing
\eqref{rf:eq:dist-upper}--\eqref{rf:eq:dist-lower} and \eqref{rf:eq:lift-points} to the limit gives
\begin{align}
 d_{\bar X}(\Phi(z),\Phi(z'))&\le\bar C\,d_X(z,z')^{1-\vartheta},\label{rf:eq:Phi-upper}\\
 d_X(z,z')&\le\bar C\,d_{\bar X}(\Phi(z),\Phi(z'))^{1/(1+\vartheta)}.\label{rf:eq:Phi-lower}
\end{align}
The first holds when $d_X(z,z')<\bar r$, and the second when the corresponding
$\bar X$-distance is less than $\bar r$.  Thus $\Phi$ extends continuously to $X$,
both inequalities persist under these restrictions, and $\Phi$ is injective.

We claim that
\begin{equation}\label{rf:eq:Phi-close}
 \sup_{x\in M_i}d_{\bar X}\bigl(\Phi(\psi_i(x)),\bar f_i(x)\bigr)
 \longrightarrow0.
\end{equation}
Indeed, fix $0<\eps<\bar r/3$ and a finite $\eps$-dense set
$X'_\eps\subset X'$.  For large $i$, $\eps_i<\eps$ and the defining convergence is
within $\eps$ on
$X'_\eps$.  Given $x\in M_i$, choose $z\in X'_\eps$ with
$d_X(\psi_i(x),z)<\eps$.  Then $d_i(x,x_i(z))<3\eps$, and
\eqref{rf:eq:dist-upper}, \eqref{rf:eq:Phi-upper}, and the triangle inequality bound the distance in
\eqref{rf:eq:Phi-close} by $2\bar C(3\eps)^{1-\vartheta}+2\eps$.  Letting
$\eps\downarrow0$ proves the claim.

Since $\bar f_i$ is surjective, \eqref{rf:eq:Phi-close} makes the compact set $\Phi(X)$
dense in $\bar X$.  Thus $\Phi$ is a homeomorphism, with the global H\"older bounds in
(b) following from \eqref{rf:eq:Phi-upper}--\eqref{rf:eq:Phi-lower} and the diameter bounds.

\emph{The Fibrations.}  Put $f_i:=\Phi^{-1}\circ\bar f_i$.  It is continuous and
surjective with the same positive-dimensional infranilmanifold fibres, and
\eqref{rf:eq:Phi-close} and \eqref{rf:eq:Phi-lower} give
\begin{equation*}
 \sup_{x\in M_i}d_X\bigl(\psi_i(x),f_i(x)\bigr)\longrightarrow0.
\end{equation*}
Hence $f_i$ is an $\eps_i'$-Gromov--Hausdorff approximation, where $\eps_i'$ is $\eps_i$
plus twice this supremum, and its diagram is \eqref{rf:eq:diagram} composed with
$\Phi^{-1}$.

Finally, if the whole-sequence assertion failed for some $\eps>0$, its bad
subsequence would have a good further subsequence by the argument above, a
contradiction.
\end{proof}

\begin{proof}[Proof of \cref{rf:thm:fibration}]
Apply \cref{rf:thm:rong} with $\vartheta=\frac12$, so that
\mbox{$\delta_{\mathrm{aff}}(n):=\delta_{\mathrm{aff}}(n,\tfrac12)$}. By its last
assertion there are indices $i_1<i_2<\cdots$ such that every $M_i$ with
$i\ge i_k$ carries a map $f_i$ of the kind described there that is a
$\frac1k$-Gromov--Hausdorff approximation.  For $i_k\le i<i_{k+1}$ choose
such a map and put $\eps_i=\frac1k$.  Properties (i) and (ii) are part~(c)
of \cref{rf:thm:rong}, and (iii) is the diagram \eqref{rf:eq:diagram}, in which the
manifold $\bar Y$ may depend on $i$.
\end{proof}

\subsection{A Consequence for Aspherical Manifolds}\label{rf:sec:aspherical}
The homeomorphism $\Phi$ transfers the manifold conclusion of
\cite[Conjecture~0.7]{Ron26} under the hypotheses of case~(0.6.1) of
\cite[Conjecture~0.6]{Ron26}.

\begin{corollary}[Aspherical Limits]\label{rf:cor:aspherical}
Let $(M_i,g_i)$ and $X$ be as in \cref{rf:thm:fibration}. If every $M_i$ is aspherical and
its fundamental group is of type $A_s$ in the sense of \cite[Definition~0.4]{Ron26},
then $X$ is homeomorphic to an aspherical manifold.
\end{corollary}
\begin{proof}
Apply \cref{rf:thm:rong} with $\vartheta=\frac12$.  The bounded-curvature sequence
$(M_i,\bar g_i)$ collapses to $\bar X$ with positive diameter, so
\cite[Theorem~A]{Ron26} makes $\bar X$ an aspherical manifold; now use
$\Phi:X\to\bar X$.
\end{proof}

The conclusion makes $X$ abstractly smoothable by transport of a smooth structure
from the bounded-curvature limit.

\section{Obstructions under Rewinding Volume Control}\label[appendix]{rf:sec:example}

Products of Nifa's four-manifolds with shrinking circles keep the volumes of the balls in
the local covers bounded below, but admit no continuous Gromov--Hausdorff approximations
with small error onto their limit. The argument is conditional on the unrefereed
construction \cite{Nif26}, and \cref{rf:sec:nifa} lists precisely the properties used.

\subsection{Statements and Their Scope}\label{rf:sec:obstruction-statement}

The obstruction under LRVP concerns continuous maps $f_i:M_i\to X$ whose
approximation errors tend to zero.
Case~(0.6.2) of Conjecture~0.6 in \cite{Ron22} asks in addition for finite infranil-quotient
fibres capturing the short loops. Open Problem~4.4 of \cite{Ron26} fixes $\eps,\rho_0>0$,
writes $X_\eps$ for the $\eps$-regular region and puts
\begin{equation*}
 X_0:=\{x\in X_\eps:d(x,X\setminus X_\eps)\ge\rho_0\}.
\end{equation*}
A smooth infranil bundle is already available over $X_0$, and the problem asks whether it
is the restriction of such a map. Our examples show that under LRVP the continuous
approximation itself need not exist, so both statements fail regardless of their
requirements on the fibres or over $X_0$.

\begin{theorem}[No Continuous Approximation]\label{rf:thm:example}
Let $(\Sigma_j,g_j)$ be the closed four-manifolds of \cite[Theorem~1.3]{Nif26}, and let
$\Sigma_\infty$ be their Gromov--Hausdorff limit. Let $\lambda_j\downarrow0$ with
$\lambda_j<2/3$, and give $P_j:=\Sigma_j\times(\R/\lambda_j\Z)$ the product metric.
Then the following hold.
\begin{enumerate}[label=\textup{(\roman*)}]
\item $\Ric_{P_j}\ge0$, and there is $v_*>0$ such that every point of every $P_j$ is a
      $(1,v_*)$-LRVP.
\item $P_j\to\Sigma_\infty$ in the Gromov--Hausdorff topology,
      $\diam\Sigma_\infty>0$, and $\vol(P_j)\to0$.
\item There are $\eps_0>0$ and $j_0$ such that $\dis f\ge\eps_0$ for every continuous
      map $f:P_j\to\Sigma_\infty$ with $j\ge j_0$. In particular, for $\eps<\eps_0$ and
      $j\ge j_0$ there is no continuous $\eps$-Gromov--Hausdorff approximation from
      $P_j$ to $\Sigma_\infty$.
\item $\Sigma_\infty$ is not locally contractible, and it is not homeomorphic to the
      orbit space of a smooth action of a compact Lie group on a manifold.
\item Every point of every $P_j$ is a $(\delta,\rho)$-LRRP for all
      $\delta>\frac{11}{12}$ and all $\rho>0$, while for no $\delta\le0.44$ and no
      $\rho>0$ do all the $P_j$ have uniform $(\delta,\rho)$-LRRP.
\end{enumerate}
\end{theorem}

Part~\textup{(iii)} excludes any proposed fibration structure that requires continuous
Gromov--Hausdorff approximations with errors tending to zero.

\begin{corollary}[The Statements under LRVP]\label{rf:cor:conjectures}
The manifolds $P_j$ of \cref{rf:thm:example} are closed five-manifolds with $\Ric\ge-4$
and uniform $(1,v_*)$-LRVP that collapse to the compact space $\Sigma_\infty$. For this
sequence: \textup{(a)} there is no singular nilpotent fibration
$P_j\to\Sigma_\infty$ for $j$ large, so \cref{rf:thm:fibration} fails with LRRP replaced by
LRVP, or by $(\delta,\rho)$-LRRP with $\delta>\frac{11}{12}$;
\textup{(b)} case~(0.6.2) of \cite[Conjecture~0.6]{Ron22} fails;
\textup{(c)} \cite[Open Problem~4.4]{Ron26} has a negative answer; and
\textup{(d)} \cite[Conjecture~3.3]{HRW20} fails whenever a singular fibration is, in
particular, a continuous Gromov--Hausdorff approximation onto the limit with error
tending to zero.
\end{corollary}

The \emph{Reifenberg constant} of a sequence is the infimum of the $\delta$ for which
some $\rho>0$ gives uniform $(\delta,\rho)$-LRRP. For $P_j$ it lies between $0.44$ and
$\frac{11}{12}$; the upper bound holds for every product of a circle with a complete
four-manifold of non-negative Ricci curvature.

\subsection{Proofs and Quantitative Bounds}\label{rf:sec:obstruction-proofs}
We recall the properties of Nifa's sequence that the argument uses, verify the hypotheses
for the products with circles, prove that no continuous approximation exists and that the
limit is not an orbit space, and bound the Reifenberg constant from both sides; the last
subsection assembles the proofs.

\subsubsection{Nifa's Sequence}\label{rf:sec:nifa}
After the subsequences in \cite[Theorem~1.3]{Nif26}, write
$(\Sigma_j,g_j,p_j)\to(\Sigma_\infty,d_\infty,p_\infty)$ for Nifa's pointed sequence
and its limit.  In \cite{Nif26} they are denoted $M_j$ and $Y$.  Write $\Cone(Z)$ for
the cone over $Z$ and $S^1_{2\pi {\mathfrak q}}$ for the circle of length $2\pi {\mathfrak q}$; the parameter ${\mathfrak q}$ is denoted $q$ in \cite{Nif26}.  We use exactly
the following cited properties.

\begin{enumerate}[label=\textup{(N\arabic*)},leftmargin=3.2em]
\item \emph{The Sequence.}  Each $\Sigma_j$ is a closed four-manifold, and there are
  $c_{\mathrm{Ric}},v_0>0$ with
  \begin{equation*}
   \Ric_{g_j}\ge c_{\mathrm{Ric}}\,g_j,\qquad \vol_{g_j}B_{\Sigma_j}(p_j,1)\ge v_0 .
  \end{equation*}
  The space $\Sigma_\infty$ is compact, and the convergence is realized by isometric
  embeddings of all $\Sigma_j$ and of $\Sigma_\infty$ into one compact metric space,
  called the ambient space below, in which $\Sigma_j\to\Sigma_\infty$ in the Hausdorff
  distance and $p_j\to p_\infty$ \cite[Theorem~1.3, Lemma~9.3, proof of
  Proposition~2.3]{Nif26}.
\item \emph{The Singular Point.}  The tangent cone of $\Sigma_\infty$ at
  $p_\infty$ is unique and equal to $\R^2\times\Cone(S^1_{2\pi {\mathfrak q}})$, where
  $0<{\mathfrak q}<2\cdot10^{-6}$.  The inclusion of any neighbourhood of $p_\infty$ into
  $\Sigma_\infty$ induces a map on $H_2(\,\cdot\,;\Z)$ whose image has infinite rank; in
  particular $\Sigma_\infty$ is not locally contractible at $p_\infty$
  \cite[Theorem~1.3 and~(7.1)]{Nif26}.
\item \emph{The Cores.}  There are compact Riemannian four-manifolds with boundary
  $W_i$, $i\ge1$, each diffeomorphic to $\CP^2$ minus an open ball, such that for all
  $j\ge i$ the manifold $\Sigma_j$ contains an isometric copy of $W_i$ as a smooth
  domain, and the copies of different $W_i$ in $\Sigma_j$ are disjoint
  \cite[(9.3)]{Nif26}.  Let $Q_i:=W_i/\partial W_i$, a closed topological
  four-manifold homeomorphic to $\CP^2$, with base point $*$ the image of
  $\partial W_i$, and let $\varpi_{i,j}:\Sigma_j\to Q_i$ be the quotient map on $W_i$
  and the constant map to $*$ on $\Sigma_j\setminus\Int W_i$.
\item \emph{The Limit Maps.}  There are continuous maps
  $\iota_i:W_i\to\Sigma_\infty$ and $\varpi_i:\Sigma_\infty\to Q_i$ such that
  $\varpi_i\circ\iota_i$ is the quotient map, $\iota_i$ is injective on $\Int W_i$, the
  inclusions $W_i\to\Sigma_j$ converge to $\iota_i$ uniformly in the ambient space as
  $j\to\infty$, and $\varpi_{i,j}(y_j)\to\varpi_i(y)$ whenever $\Sigma_j\ni y_j\to y$
  in the ambient space \cite[Proposition~2.3 and its proof, Lemma~2.5]{Nif26}.
\item \emph{The Regular Part.}  The set
  $\Sigma_\infty^\circ:=\Sigma_\infty\setminus\{p_\infty\}$ is a connected topological
  four-manifold.  It carries a smooth structure and a smooth Riemannian metric
  $g_\infty$ such that every point of $\Sigma_\infty^\circ$ has a neighbourhood on which
  $d_\infty$ is the distance of $g_\infty$.  Moreover there are a compact domain
  $\mathcal K\subset\Sigma_\infty^\circ$ with non-empty interior and smooth embeddings
  $\phi_j:\mathcal K\to\Sigma_j$ such that $\phi_j^*g_j\to g_\infty$ in
  $C^\infty(\mathcal K)$ and $\phi_j$ converges to the inclusion of $\mathcal K$
  uniformly in the ambient space \cite[Lemma~14.1 and its proof]{Nif26}.
\end{enumerate}

Two consequences of (N3)--(N4) will be used.  For $i\ne k$ and
$y\in\Sigma_\infty$,
\begin{equation}\label{rf:eq:wedge}
 \varpi_i(y)=*\quad\hbox{or}\quad\varpi_k(y)=*.
\end{equation}
Indeed, for $y_j\to y$, one of the two disjoint cores is avoided along a subsequence,
and the corresponding quotient converges to $*$.  Also
\begin{equation}\label{rf:eq:Wi}
 \mathcal W_i:=\varpi_i^{-1}(Q_i\setminus\{*\})=\iota_i(\Int W_i).
\end{equation}
The reverse inclusion follows because $\varpi_i\circ\iota_i$ is the quotient map.  For
the forward inclusion, an approximating sequence for a point with $\varpi_i(y)\ne*$
eventually lies in a compact subset of $\Int W_i$; the inverse is $\iota_i$ composed
with the inverse of the quotient there. Thus
$\varpi_i|_{\mathcal W_i}:\mathcal W_i\to Q_i\setminus\{*\}$ is a homeomorphism, and
$\mathcal W_i$ is an open four-manifold,
disjoint from $p_\infty$ by (N2).

The lower bound for the Reifenberg constant also uses the shape of $\Sigma_j$ near the
base point. The next lemma extracts what is needed from the construction.

\begin{lemma}[Simply Connected Neighbourhoods of the Base Point]\label{rf:lem:nbhd}
There are numbers $s_m\downarrow0$ and a constant $C_\Sigma\ge1$ such that for all
$j\ge m$ there is an open, connected and simply connected set
$\mathcal S_{m,j}\subset\Sigma_j$ with
\begin{equation*}
 B_{\Sigma_j}(p_j,s_m/64)\subset \mathcal S_{m,j}\subset B_{\Sigma_j}(p_j,C_\Sigma s_m).
\end{equation*}
\end{lemma}
\begin{proof}
We use the construction of \cite[Section~14]{Nif26}, in which $p_j=p$ for all $j$ and
$\Sigma_j$ is obtained from the unit sphere $S^4$ by replacing disjoint round balls with
the cores $W_i$ of \textup{(N3)}. Let $d_{S^4}$ be the round distance, and let $s_m$ be
the numbers of the proof of \cite[Lemma~14.1]{Nif26}. They satisfy $s_{m+1}=s_m/2$, the
round ball that $W_i$ replaces is centred at $d_{S^4}$-distance $\frac43s_i$ from $p$, and
the open $d_{S^4}$-ball of radius $s_m$ about $p$ contains the round balls replaced by
$W_i$ for $i>m$ and is disjoint from those replaced by $W_i$ for $i\le m$
\cite[Section~14 and the proof of Lemma~14.1]{Nif26}. Let $\mathcal S_{m,j}\subset\Sigma_j$ be the
open set obtained from this ball by the replacements of index $i>m$ made in $\Sigma_j$.
It is connected and simply connected by van Kampen, since finitely many disjoint balls
are replaced by copies of $\CP^2$ minus a ball.

Let $\mathsf f_j$ be the function on $\Sigma_j$ given by \cite[Lemma~13.2]{Nif26}, where
it is denoted $\lambda_j$. It vanishes at $p$, it is $32$-Lipschitz, it equals
$\frac43s_i$ on $W_i$, it lies between one half and twice $d_{S^4}(p,\cdot)$ outside the
$W_i$, and $d_{\Sigma_j}(p,\cdot)\le C_*\mathsf f_j$ for a constant $C_*$ independent of
$j$ \cite[Lemma~13.2, (13.5) and~(13.6)]{Nif26}. Moreover $\mathsf f_j\le2s_m$ on
$\mathcal S_{m,j}$ \cite[(14.4)]{Nif26}. Outside $\mathcal S_{m,j}$, a point lies either on some $W_i$ with
$i\le m$, where $\mathsf f_j=\frac43s_i\ge\frac43s_m$, or outside the $W_i$ at
$d_{S^4}$-distance at least $s_m$ from $p$; in both cases $\mathsf f_j\ge s_m/2$. Since
$\mathsf f_j\le32\,d_{\Sigma_j}(p,\cdot)$, every point at distance less than $s_m/64$
from $p$ lies in $\mathcal S_{m,j}$, and every point of $\mathcal S_{m,j}$ is at distance at most
$2C_*s_m$ from $p$. This proves the inclusions with $C_\Sigma=2C_*+1$.
\end{proof}

\subsubsection{The Products}\label{rf:sec:products}
Let $\lambda_j\downarrow0$, with $\lambda_j<2/3$, and equip
\begin{equation*}
 P_j:=\Sigma_j\times(\R/\lambda_j\Z)
\end{equation*}
with the product metric; write $[z,t]$ for the image of $(z,t)\in\Sigma_j\times\R$. The
next proposition verifies for $P_j$ the hypotheses of the statements recalled in the
introduction.

\begin{proposition}[The Hypotheses]\label{rf:prop:products}
Put $\mathfrak d:=\pi\sqrt{3/c_{\mathrm{Ric}}}$ and
$v_*:=2v_0/(3^5(\mathfrak d+1)^4)$.
\begin{enumerate}[label=\textup{(\alph*)}]
\item $\Ric_{P_j}\ge0$.
\item Every point of every $P_j$ is a $(1,v_*)$-LRVP.
\item $d_{GH}(P_j,\Sigma_j)\le\lambda_j/2$ and
  $\vol(P_j)\le\lambda_j\,\omega_4\mathfrak d^4$.  Hence $P_j\to\Sigma_\infty$ and
  $\vol(P_j)\to0$.
\end{enumerate}
\end{proposition}

\begin{proof}
(a) The circle factor is flat and $\Ric_{g_j}>0$.

(b) By Myers' theorem and (N1), $\diam\Sigma_j\le\mathfrak d$.  For
$z\in\Sigma_j$ the ball $B_{\Sigma_j}(z,\mathfrak d+1)$ contains
$B_{\Sigma_j}(p_j,1)$, and the Bishop--Gromov inequality for $\Ric_{g_j}\ge0$ gives
\begin{equation}\label{rf:eq:allpoint}
 \vol B_{\Sigma_j}(z,\tfrac13)\ge
 \frac{\vol B_{\Sigma_j}(z,\mathfrak d+1)}{3^4(\mathfrak d+1)^4}
 \ge\frac{v_0}{3^4(\mathfrak d+1)^4}.
\end{equation}
For $x=[z,0]$, let $B'$ be the component over $B_{P_j}(x,1)$ in
$\Sigma_j\times\R$ containing $(z,0)$.  Since $\lambda_j<2/3$, the connected cylinder
$B_{\Sigma_j}(z,\frac13)\times\R$ lies in $B'$, and its rectangle with
$|t|<\frac13$ lies in the intrinsic unit ball about $(z,0)$.  The local-isometric
cover $\tB_{P_j}(x,1)\to B'$ maps the unit ball about $\tilde x$ onto that ball by path
lifting.
Writing $g$ for the product metric, \eqref{rf:eq:allpoint} gives
\begin{equation*}
 \vol B_{\pi_x^*g}(\tilde x,1)\ge
 \vol\bigl(B_{\Sigma_j}(z,\tfrac13)\times(-\tfrac13,\tfrac13)\bigr)
 =\frac23\vol B_{\Sigma_j}(z,\tfrac13)\ge v_* .
\end{equation*}

(c) The projection $P_j\to\Sigma_j$ is onto and changes distances by at most
the diameter $\lambda_j/2$ of the circle.  By the Bishop inequality,
$\vol(\Sigma_j)\le\omega_4\mathfrak d^4$.
\end{proof}

\subsubsection{No Continuous Approximation}\label{rf:sec:nogha}
The obstruction is topological and concerns only $(\Sigma_j,\Sigma_\infty)$.

\begin{theorem}[No Continuous Approximation, General Form]\label{rf:thm:nogha}
There are $\eps_0>0$ and $j_0$ such that for $j\ge j_0$
\begin{equation*}
 \dis f\ge\eps_0\qquad\hbox{for every continuous map }f:\Sigma_j\to\Sigma_\infty .
\end{equation*}
The same holds, with the same $\eps_0$ and $j_0$, for every continuous map
$f:\Sigma_j\times(\R/\lambda\Z)\to\Sigma_\infty$ and every $\lambda>0$.
\end{theorem}

All degrees below are modulo two.  For a proper map between four-manifolds, write
$\deg_2(f;y)$ for the image of the fundamental class around the compact preimage of
$y$ in the target's local homology $\Z_2$.  It is locally constant in $y$, and on
closed connected four-manifolds it is the multiplier on $H^4(\,\cdot\,;\Z_2)$
\cite[Chapter~VIII, \S4]{Dol95}. The degree over a point is detected by the following
form of the Borsuk--Ulam theorem.

\begin{lemma}[Borsuk]\label{rf:lem:borsuk}
Let a continuous map from the closed ball $\overline{B_{\R^4}(0,r)}$ to $\R^4$ take
distinct values at any two points whose distance is at least $r$.  Then the image of
the centre is not in the image of the boundary sphere, and the degree of the map over
the image of the centre is odd.
\end{lemma}

\begin{proof}
Denote the map by $b$.  For $|x|=r$ and $0\le t\le1$ we have $|x-(-tx)|\ge r$, so
$b(x)\ne b(-tx)$.  Thus the normalized map below is always defined; at $t=0$ it also
shows that $b(0)$ is missed by the boundary and gives the degree map over $b(0)$:
\[
x\longmapsto\frac{b(x)-b(-tx)}{|b(x)-b(-tx)|}.
\]
At $t=1$ this map is odd, so its degree, and hence the degree of $b$ over $b(0)$, is
odd.
\end{proof}

With it, a continuous map of small distortion is shown to have degree one over the
regular part of $\Sigma_\infty$.

\begin{lemma}[The Degree Is One]\label{rf:lem:deg}
There are $\eps_0>0$ and $j_0$ such that for $j\ge j_0$ and every
continuous $f:\Sigma_j\to\Sigma_\infty$ with $\dis f<\eps_0$, the proper map
\begin{equation*}
 f:f^{-1}(\Sigma_\infty^\circ)\longrightarrow\Sigma_\infty^\circ
\end{equation*}
has $\deg_2=1$.
\end{lemma}

\begin{proof}
The restricted map is proper because the preimage of a compact subset of
$\Sigma_\infty^\circ$ is closed in $\Sigma_j$.  Its target is a
connected four-manifold by (N5), so it suffices to find one point over
which the degree is one.

Fix two points $y_1\ne y_2$ in the interior of $\mathcal K$ and smooth coordinate
charts of $\Sigma_\infty^\circ$ centred at them, with domains in the
interior of $\mathcal K$.  By (N5), as in the
proof of \cite[Lemma~14.1]{Nif26}, there are $0<c_1<c_2$, $r_2>0$ and $j_1$ such that
for $j\ge j_1$ the composition $\phi_j^{(1)}$ of the inverse of the chart at $y_1$ with
$\phi_j$ is defined on $B_{\R^4}(0,r_2)$ and satisfies
\begin{equation*}
 \sqrt{c_1}\,|x-y|\le d_{\Sigma_j}(\phi_j^{(1)}(x),\phi_j^{(1)}(y))\le\sqrt{c_2}\,|x-y|,
 \qquad
 B_{\Sigma_j}\bigl(\phi_j^{(1)}(0),\tfrac12\sqrt{c_1}\,r_2\bigr)
 \subset\phi_j^{(1)}\bigl(B_{\R^4}(0,\tfrac12r_2)\bigr),
\end{equation*}
and the same holds for the chart at $y_2$.  Put $c_3:=d_\infty(y_1,y_2)>0$.  By
(N5) there is $j_2$ with $d_{\Sigma_j}(\phi_j(y_1),\phi_j(y_2))\ge c_3/2$ for
$j\ge j_2$.

The set of points of $\Sigma_\infty$ at distance at least $c_3/8$ from $p_\infty$ is a
compact subset of the manifold $\Sigma_\infty^\circ$.  By the Lebesgue
number lemma there is $r_1>0$ such that every ball $B_{\Sigma_\infty}(y,r_1)$
about a point $y$ of this set lies in an open subset of
$\Sigma_\infty^\circ$ homeomorphic to $\R^4$.  Choose
\begin{equation*}
 0<r<\frac{r_2}2\ \hbox{ with }\ 2\sqrt{c_2}\,r\le r_1,
 \qquad
 \eps_0:=\min\Bigl\{\frac{c_3}8,\ \frac{\sqrt{c_1}\,r}2,\ \frac{r_1}2\Bigr\},
 \qquad j_0:=\max\{j_1,j_2\}.
\end{equation*}

Let $j\ge j_0$ and $\dis f<\eps_0$.  The images of $\phi_j(y_1)$ and $\phi_j(y_2)$
under $f$ are at distance more than $c_3/2-\eps_0\ge3c_3/8$, so they are not both
within $c_3/8$ of $p_\infty$.  Say $y':=f(\phi_j(y_1))$ is at distance at least
$c_3/8$ from $p_\infty$, and let $V$ be an open set homeomorphic to $\R^4$ that
contains $B_{\Sigma_\infty}(y',r_1)$.  Consider
$b:=f\circ\phi_j^{(1)}$ on the closed ball of radius $r$.  Its image lies in
$B_{\Sigma_\infty}(y',\sqrt{c_2}\,r+\eps_0)\subset V$.  If $|x-y|\ge r$, then
$d_{\Sigma_j}(\phi_j^{(1)}(x),\phi_j^{(1)}(y))\ge\sqrt{c_1}\,r\ge2\eps_0$, so $b(x)\ne b(y)$.
By \cref{rf:lem:borsuk}, applied in $V\cong\R^4$, the degree of $b$ over
$y'$ is odd.  Finally, every point of $f^{-1}(y')$ is within $\eps_0$ of
$\phi_j(y_1)=\phi_j^{(1)}(0)$ in $\Sigma_j$.  Since $\eps_0<\frac12\sqrt{c_1}\,r_2$, it lies in
$\phi_j^{(1)}(B_{\R^4}(0,\frac12r_2))$, and its coordinate has norm less than
$\eps_0/\sqrt{c_1}<r$.  So $f^{-1}(y')\subset\phi_j^{(1)}(B_{\R^4}(0,r))$, and
$\deg_2(f;y')=1$.
\end{proof}

\begin{proof}[Proof of \cref{rf:thm:nogha}]
Let $\eps_0,j_0$ be as in \cref{rf:lem:deg}, let $j\ge j_0$, and suppose
that $f:\Sigma_j\to\Sigma_\infty$ is continuous with $\dis f<\eps_0$.  For $i\ge1$
put
\begin{equation*}
 \zeta_i:=(\varpi_i\circ f)^*\bar\zeta_i\in H^2(\Sigma_j;\Z_2),
\end{equation*}
where $\bar\zeta_i$ generates $H^2(Q_i;\Z_2)$.  Here $i$ runs over all positive
integers, not only over the cores present in $\Sigma_j$: the maps $\varpi_i$ are
defined on $\Sigma_\infty$.

\emph{Squares.}  Let $z\in Q_i\setminus\{*\}$, and let $y\in\mathcal W_i$ be its
preimage under $\varpi_i$.  By \eqref{rf:eq:Wi},
${(\varpi_i\circ f)^{-1}(z)=f^{-1}(y)}$, and
$\varpi_i$ is a homeomorphism near $y$.  So
$\deg_2(\varpi_i\circ f;z)=\deg_2(f;y)=1$ by \cref{rf:lem:deg}.  Since $\bar\zeta_i^2$
generates $H^4(Q_i;\Z_2)$,
\begin{equation*}
 \zeta_i^2=(\varpi_i\circ f)^*(\bar\zeta_i^2)\ne0
 \quad\hbox{in }H^4(\Sigma_j;\Z_2)=\Z_2 .
\end{equation*}

\emph{Products.}  For $i\ne k$ the map
$(\varpi_i,\varpi_k):\Sigma_\infty\to Q_i\times Q_k$ takes values in the wedge
$Q_i\vee Q_k$, by \eqref{rf:eq:wedge}.  The class $\bar\zeta_i\times\bar\zeta_k$ restricts
to zero on the wedge, because it vanishes on each of the two factors.  Hence
$\varpi_i^*\bar\zeta_i\cup\varpi_k^*\bar\zeta_k=0$ on $\Sigma_\infty$, and
$\zeta_i\zeta_k=0$.

\emph{Count.}  If $\sum_ic_i\zeta_i=0$ is a finite relation with
$c_i\in\Z_2$, multiplication by $\zeta_k$ gives $c_k\zeta_k^2=0$, so $c_k=0$.
The classes $\zeta_1,\zeta_2,\dots$ are therefore linearly independent in
$H^2(\Sigma_j;\Z_2)$, which is finite-dimensional.  This is a contradiction.

For the product, restrict a continuous map
$f:\Sigma_j\times(\R/\lambda\Z)\to\Sigma_\infty$ to a slice
$\Sigma_j\times\{\mathrm{pt}\}$.  The slice carries the distance of $\Sigma_j$, so the
restriction has distortion at most $\dis f$.
\end{proof}

\subsubsection{The Limit Is Not an Orbit Space}\label{rf:sec:orbit}
The limit cannot be the compact-group quotient occurring in a
Cheeger--Fukaya--Gromov fibration.

\begin{proposition}[Orbit Spaces Are Locally Contractible]\label{rf:prop:orbit}
Let a compact Lie group $G$ act smoothly on a smooth manifold $M$.  Then
$M/G$ is locally contractible.  Consequently $\Sigma_\infty$ is not homeomorphic
to $M/G$ for any such $M$ and $G$.
\end{proposition}

\begin{proof}
By the slice theorem \cite[Chapter~VI, Section~2]{Bre72}, a neighbourhood of
$[z]\in M/G$ is modelled on a neighbourhood of $0$ in $V/G_z$.  A
$G_z$-invariant ball in $V$ contracts equivariantly to $0$, so its quotient contracts;
such quotients form a neighbourhood basis.  Apply (N2).
\end{proof}

\subsubsection{The Reifenberg Constant: Upper Bound}\label{rf:sec:upper}
We quantify the facts that a ball in the local cover of a product contains a line
segment through its centre and obeys Bishop--Gromov comparison.

\begin{lemma}[A Crude Correspondence]\label{rf:lem:crude}
Let $(Z,d)$ be a metric space with a point $z_0$ such that $d(z_0,z)<s$ for
all $z\in Z$.  Let $0<\eps\le1$, let $m_1,m_2\ge1$ be integers, and suppose that
\begin{enumerate}[label=\textup{(\alph*)}]
\item $S^{n-1}$ is the union of $m_1$ sets of diameter at most $2-\eps$;
\item there are $z_1,\dots,z_{m_1}\in Z$ with $d(z_0,z_i)\le(1-\eps)s$ and
  $d(z_i,z_k)\ge\eps s$ for $i\ne k$;
\item $Z$ is the union of $m_2$ sets of diameter at most $(2-\eps)s$;
\item there are $y_1,\dots,y_{m_2}\in\R^n$ with $|y_k|\le(1-\eps)s$ and
  $|y_k-y_l|\ge\eps s$ for $k\ne l$.
\end{enumerate}
Then there is a correspondence between $Z$ and $B_{\R^n}(0,s)$ that contains
$(z_0,0)$ and has distortion at most $(2-\eps)s$.  Hence the pointed
Gromov--Hausdorff distance between $(Z,z_0)$ and $(B_{\R^n}(0,s),0)$ is at
most $(1-\frac\eps2)s$.
\end{lemma}

\begin{proof}
Let $\mathcal C_1,\dots,\mathcal C_{m_1}$ be the sets in (a) and
$Z_1,\dots,Z_{m_2}$ those in (c), made disjoint.  Relate
$y\in B_{\R^n}(0,s)\setminus\{0\}$ to $z_i$ if $y/|y|\in\mathcal C_i$, relate
$z\in Z$ to $y_k$ if $z\in Z_k$, and relate $z_0$ to $0$.  Both projections
are onto.  For two pairs from the first rule, equal indices give distortion at most
$(2-\eps)s$ by radial convexity and (a), while unequal indices compare numbers in
$[\eps s,2(1-\eps)s]$ and $[0,2s)$.  The second rule is symmetric.  For one pair of
each kind, both distances are less than $(2-\eps)s$; a pair with $(z_0,0)$ compares
numbers in $[0,s)$.  Thus the distortion is at most $(2-\eps)s$, and the pointed bound
follows from \cref{rf:sec:conv}.
\end{proof}

The next lemma bounds the number of sets needed in hypothesis~(c).

\begin{lemma}[Covering by Radial Retraction]\label{rf:lem:radial}
Let $(M^n,g)$ be a complete Riemannian manifold with $\Ric_g\ge0$, let $x\in M$ and
$0<s<\rho$, and let $0<\eps<2a$ with $3a+\frac12(a-\frac\eps2)\le1$.  Then the ball
$B_{\pi_x^*g}(\tilde x,s)$ in the local cover is the union of at most
\begin{equation*}
 \Bigl(\frac{4a}{a-\eps/2}+1\Bigr)^n
\end{equation*}
sets of diameter less than $(2-\eps)s$.
\end{lemma}

\begin{proof}
Put $r:=(a-\frac\eps2)s$ and write $d$ for the lifted distance.  Radially retract $y$
along a minimizing geodesic, available
by \cref{rf:lem:cover}, to $y^\flat$ at distance $a\,d(\tilde x,y)$ from $\tilde x$.
Then
\begin{equation*}
 d(y,y')<2(1-a)s+d\bigl(y^\flat,y'^\flat\bigr).
\end{equation*}
Let $z_1,\dots,z_k$ be a maximal $r$-separated subset of
$B_{\pi_x^*g}(\tilde x,as)$.  The balls of radius $r$ about the $z_i$ cover
$B_{\pi_x^*g}(\tilde x,as)$, and their radial preimages have diameter less than
$2(1-a)s+2r=(2-\eps)s$.  The $r/2$-balls are disjoint.  If $z_*$ has the least
such volume, their union lies in $B(z_*,2as+r/2)$; moreover $z_*$ has depth at least
$\rho-as\ge2as+r/2$.  Hence \cref{rf:lem:cover}\,(iii) gives
\begin{equation*}
 k\le\frac{\vol B_{\pi_x^*g}(z_*,2as+\frac r2)}{\vol B_{\pi_x^*g}(z_*,\frac r2)}
 \le\Bigl(\frac{4as}{r}+1\Bigr)^n .\qedhere
\end{equation*}
\end{proof}

The two lemmas give the upper bound.

\begin{proposition}[The Upper Bound]\label{rf:prop:upper}
Let $(M^4,g)$ be a complete Riemannian manifold with $\Ric_g\ge0$, let
$\lambda>0$, and let $P=M\times(\R/\lambda\Z)$ with the product metric.  For every
$x\in P$, every $\rho>0$ and every $0<s<\rho$, the pointed
Gromov--Hausdorff distance between the ball of radius $s$ about $\tilde x$ in the local
cover at $x$ and $B_{\R^5}(0,s)$ is at most $\frac{11}{12}s$.  In particular every point of
$P$ is a $(\delta,\rho)$-LRRP for every $\delta>\frac{11}{12}$ and
every $\rho>0$.
\end{proposition}

\begin{proof}
We verify the hypotheses of \cref{rf:lem:crude} with $n=5$, $\eps=\frac16$, $m_1=10$,
and $Z$ the ball of radius $s$ about $z_0=\tilde x$ in the local cover at $x$.
For (a), the ten sets $\{y\in S^4:\pm y_i=\max_k|y_k|\}$ cover $S^4$ and have
diameter at most $\sqrt{16/5}<11/6$: in the first set,
$\langle y,y'\rangle\ge2/5-1$.  For (b), the circle projection lifts to a
$1$-Lipschitz function $\varsigma$ on $\tB_P(x,\rho)$ with
$\varsigma(\gamma(t))=t$ along the arclength lift $\gamma$ from $\tilde x$, for
$|t|<\rho$.  Thus the lifted-circle points
$\gamma(ks/6)$, $-5\le k\le4$, are pairwise $s/6$-separated and lie within $5s/6$
of $\tilde x$.  For (c), \cref{rf:lem:radial} with $a=0.29$ gives
$m_2\le12\,646$ sets.  Finally, the $16\,875$ points $(s/6)k$ with
$k\in\Z^5$ and $|k|\le5$ verify (d), since they are $s/6$-separated and lie in the
ball of radius $5s/6$.  Thus \cref{rf:lem:crude} gives $(1-1/12)s$.
\end{proof}

\subsubsection{The Reifenberg Constant: Lower Bound}\label{rf:sec:lower}
A quantitative Borsuk--Ulam theorem separates the tangent model
$\R^3\times\Cone(S^1_{2\pi {\mathfrak q}})$ from $\R^5$.

\begin{lemma}[Local Covers of Products]\label{rf:lem:product-cover}
Let $\Sigma$ be a closed Riemannian manifold, $\lambda>0$, and
$x=[z,0]\in\Sigma\times(\R/\lambda\Z)$.  Let $0<s\le\rho/4$, and suppose that
$B_\Sigma(z,2s)$ is contained in a connected and simply connected set
$S\subset\Sigma$ of diameter less than $\rho/4$.  Then the ball of radius $s$ about
$\tilde x$ in the cover $\tB(x,\rho)$ is isometric to the ball of radius $s$ about
$(z,0)$ in $\Sigma\times\R$ with the product distance, by an isometry that takes
$\tilde x$ to $(z,0)$.
\end{lemma}

\begin{proof}
Let $B'$ be the component over $B(x,\rho)$ in $\Sigma\times\R$ containing $(z,0)$.
The local cover factors through a local isometry $\pi':\tB(x,\rho)\to B'$, and $B'$
contains the product $\rho$-ball.  Product geodesics between points of the $s$-ball
have length below $2s$ and remain in the $3s$-ball, so the product and intrinsic
distances agree there; the same argument from $(z,0)$ identifies the intrinsic and
product balls of radius at most $2s$.  Moreover every loop in the $2s$-ball lies in the
simply connected set $S\times(-2s,2s)\subset B'$ and hence lifts closed.  Path lifting therefore maps
the $s$-ball about $\tilde x$ onto the product $s$-ball.  Given two lifted points, join their
images by paths of length tending to their product distance; the resulting loops with
paths from the base point lie in the $2s$-ball and lift closed, so the lifted distance
is no larger.  The reverse inequality follows from local isometry.
\end{proof}

We use unpointed Gromov--Hausdorff distance below; it is no larger than the pointed
distance in \eqref{rf:eq:lrrp}.

\begin{proposition}[The Model]\label{rf:prop:bu}
Let $0<{\mathfrak q}\le\frac12$, and let $\overline B_{\mathfrak q}$ be the closed unit ball of
$\R^3\times\Cone(S^1_{2\pi {\mathfrak q}})$ about a point of the singular axis, with the restricted
distance.  Put $c_{\mathfrak q}:=1-\cos\pi {\mathfrak q}$, and let $\eps_{\mathfrak q}$ be the smallest
positive root of
\begin{equation*}
 \bigl(\tfrac32-\tfrac52c_{\mathfrak q}\bigr)\eps^2+10(1+c_{\mathfrak q})\eps-10=0 .
\end{equation*}
Then the Gromov--Hausdorff distance between $\overline B_{\mathfrak q}$ and the closed unit ball of
$\R^5$ is at least $\delta_{\mathrm{low}}({\mathfrak q}):=\eps_{\mathfrak q}/2$.  The function
$\delta_{\mathrm{low}}$ is decreasing, with
\begin{equation*}
 \delta_{\mathrm{low}}({\mathfrak q})>0.4413\ \ \hbox{for }{\mathfrak q}\le\tfrac1{100},\qquad
 \delta_{\mathrm{low}}(\tfrac14)>0.37,\qquad
 \delta_{\mathrm{low}}(\tfrac12)>0.25 .
\end{equation*}
\end{proposition}

\begin{proof}
Let $\mathcal R$ be a correspondence between the closed unit ball of $\R^5$
and $\overline B_{\mathfrak q}$ with distortion at most $\eps$; we may assume $\eps<2$.  For
$x\in S^4$ choose a point of $\overline B_{\mathfrak q}$ that is related to $x$, and write it as
$(y(x),c(x))$, with $y(x)\in\R^3$ and $c(x)$ in the cone at
distance $\varrho(x)$ from the vertex.  Points of $\overline B_{\mathfrak q}$ satisfy
$|y|^2+\varrho^2\le1$, and the distance $d$ of two of them is given by
\begin{equation*}
 d^2=|y-y'|^2+\varrho^2+\varrho'^2-2\varrho\varrho'\cos\angle(c,c'),
 \qquad 0\le\angle(c,c')\le\pi {\mathfrak q}\le\tfrac\pi2 ,
\end{equation*}
where $\angle(c,c')$ is the angle between the directions of $c$ and $c'$.
Put $\Pi(x):=(y(x),\varrho(x))\in\R^4$.  Then, for $x,x'\in S^4$ and the
distance $d$ of the two chosen points,
\begin{equation}\label{rf:eq:Pi}
 d^2-2c_{\mathfrak q}\varrho(x)\varrho(x')\le|\Pi(x)-\Pi(x')|^2\le d^2 .
\end{equation}

\emph{Antipodes.}  For $x$ and $-x$ the distance of the chosen points is at least
$2-\eps$.  Since the cosine above is non-negative and
$|y-y'|^2\le2|y|^2+2|y'|^2\le4-2\varrho^2-2\varrho'^2$,
we get $(2-\eps)^2\le4-\varrho^2-\varrho'^2$, that is,
$\varrho(x)^2+\varrho(-x)^2\le4\eps-\eps^2$.  By \eqref{rf:eq:Pi},
\begin{equation*}
 |\Pi(x)-\Pi(-x)|^2\ge(2-\eps)^2-c_{\mathfrak q}(4\eps-\eps^2)=:R^2 .
\end{equation*}

\emph{Jumps.}  If $R^2\le0$, the inequality $4\eps^2\ge\frac52R^2$
obtained below is trivial, so assume $R>0$.  The map
\begin{equation*}
 x\longmapsto\frac{\Pi(x)-\Pi(-x)}{|\Pi(x)-\Pi(-x)|}
\end{equation*}
from $S^4$ to $S^3$ sends antipodes to antipodes.  By the theorem of Dubins
and Schwarz \cite{DS81}, in the form of \cite[Corollary~5.4]{LMS23}, its
modulus of discontinuity is at least $\arccos(-\frac14)$: for every
$\eps'>0$ there are $x',x''\in S^4$, arbitrarily close to each other,
whose images are at distance at least $\sqrt{5/2}-\eps'$.  On the
other hand, by \eqref{rf:eq:Pi} and the distortion bound,
\begin{equation*}
 |\Pi(x')-\Pi(x'')|\le|x'-x''|+\eps,
\end{equation*}
and the same holds at the antipodes.  So the numerator changes by at
most $2|x'-x''|+2\eps$ between $x'$ and $x''$.  Its norm is at least $R$, and the
radial projection onto the sphere of radius $R$ does not increase distances
outside that sphere, so the images of $x'$ and $x''$ are at distance at most
$(2|x'-x''|+2\eps)/R$.  Hence
$2\eps/R\ge\sqrt{5/2}$, that is, $4\eps^2\ge\frac52R^2$.  Written out,
the quadratic polynomial in the statement is non-negative at $\eps$, and it
is negative at $0$.  So $\eps\ge\eps_{\mathfrak q}$, and the Gromov--Hausdorff distance
is half the infimum of the distortions.  For fixed $\eps\in(0,4)$ the
polynomial increases with $c_{\mathfrak q}$, so its smallest positive root
decreases as ${\mathfrak q}$ increases.  The numerical values follow from the formula
for the root.
\end{proof}

Together with the construction, \cref{rf:prop:bu} gives the lower bound in part~(v) of
\cref{rf:thm:example}.

\begin{proposition}[The Lower Bound]\label{rf:prop:lower}
For every $\rho>0$ and every $\delta<\delta_{\mathrm{low}}({\mathfrak q})$ there is $j_1$
such that for $j\ge j_1$ the point $[p_j,0]\in P_j$ is not a
$(\delta,\rho)$-LRRP.  In particular, for no $\delta\le0.44$ and no
$\rho>0$ do all the $P_j$ have uniform $(\delta,\rho)$-LRRP.
\end{proposition}

\begin{proof}
By \cref{rf:lem:nbhd} choose $m$ with $2C_\Sigma s_m<\rho/4$, and let
$0<s\le s_m/128$.  For $j\ge m$ the ball $B_{\Sigma_j}(p_j,2s)$ lies in the
connected and simply connected set $\mathcal S_{m,j}$, whose diameter is less than $\rho/4$.  By
\cref{rf:lem:product-cover} the ball of radius $s$ about $\tilde x_j$ in the local cover at
$x_j:=[p_j,0]$ is isometric to the ball of radius $s$ about $(p_j,0)$ in
$\Sigma_j\times\R$.  By (N1), as $j\to\infty$ these balls converge in the
Gromov--Hausdorff distance to the ball of radius $s$ about $(p_\infty,0)$ in
$\Sigma_\infty\times\R$.  By (N2), as $s\to0$ the latter balls, rescaled by
$1/s$, converge to $\overline B_{\mathfrak q}$.  Therefore, by \cref{rf:prop:bu},
\begin{equation*}
 \liminf_{s\to0}\ \liminf_{j\to\infty}\
 \frac1s\,d_{GH}\bigl(B_{\pi_{x_j}^*g}(\tilde x_j,s),B_{\R^5}(0,s)\bigr)
 \ge\delta_{\mathrm{low}}({\mathfrak q}),
\end{equation*}
where $g$ is the metric of $P_j$ and the distance may be taken unpointed.  This proves
the first assertion.  The second follows because ${\mathfrak q}<\frac1{100}$ by (N2).
\end{proof}

\subsubsection{Proofs of the Obstruction Statements}\label{rf:sec:example-proof}
It remains to collect the pieces.

\begin{proof}[Proof of \cref{rf:thm:example}]
Parts (i) and (ii) are \cref{rf:prop:products}; the diameter of
$\Sigma_\infty$ is positive because $\Sigma_\infty$ has a tangent cone of dimension
four.  Part~(iii) is \cref{rf:thm:nogha}.  Part~(iv) is (N2) and
\cref{rf:prop:orbit}.  Part~(v) is
\cref{rf:prop:upper}, applied to $M=\Sigma_j$, and
\cref{rf:prop:lower}.
\end{proof}

\begin{proof}[Proof of \cref{rf:cor:conjectures}]
By parts~(i) and~(ii) of \cref{rf:thm:example}, the $P_j$ are closed five-manifolds with
$\Ric\ge0$ and uniform $(1,v_*)$-LRVP that collapse with vanishing volume to the compact,
positive-diameter limit $\Sigma_\infty$, which has dimension four. Thus they satisfy
the common hypotheses of the LRVP statements recalled in
\cref{rf:sec:obstruction-statement}. A singular
nilpotent fibration
in the sense of \cite[Section~1.3]{Ron26} is, in particular, a continuous
$\eps_j$-Gromov--Hausdorff approximation with $\eps_j\to0$: the frame-bundle map in
\eqref{rf:eq:diagram} is one, and the vertical maps are quotients by isometric actions.
So each of the statements (a)--(d) would provide, for all large $j$, a continuous
$\eps_j$-Gromov--Hausdorff approximation $P_j\to\Sigma_\infty$ with $\eps_j\to0$.  This
contradicts part~(iii) of \cref{rf:thm:example}.  The last clause of~(a) follows from
part~(v).
\end{proof}

\begin{remark}[Other Dimensions]\label{rf:rem:dimensions}
Replacing the circle by $(\R/\lambda_j\Z)^k$, with
$\lambda_j<2/(3\sqrt k)$, and using the radius $1/(3\sqrt k)$ in place of $\frac13$ in the
proof of \cref{rf:prop:products} gives
uniform $(1,v)$-LRVP for some $v=v(k,c_{\mathrm{Ric}},v_0)>0$, and hence analogues of
parts~(i)--(iv) of \cref{rf:thm:example} and the LRVP conclusions of
\cref{rf:cor:conjectures}, excluding the numerical LRRP clause in part~(a), in every
dimension $n\ge5$. The obstruction follows by restriction to a slice.
\end{remark}

\subsection{Scope of the Obstruction}\label{rf:sec:disc}
We record what the obstruction leaves open, between the two thresholds and in neighbouring
cases.

\subsubsection{Between the Two Thresholds}\label{rf:sec:between}
In dimension five, \cref{rf:thm:fibration,rf:thm:example} leave open the range between
$\delta_{\mathrm{aff}}(5)$ and $0.44$.  The upper estimate $11/12$ of \cref{rf:prop:upper} is automatic
for the circle products, whereas \cref{rf:prop:bu} shows that the conical model cannot
produce a counterexample below $1/4$ when its cone angle is at most $\pi$.  A smaller
counterexample would need a larger angle or a different mechanism.  It is also open
whether a positive threshold can be chosen independently of dimension, or whether
$\Sigma_\infty$ can arise from another Ricci-lower-bounded sequence with uniform LRVP that does
admit continuous approximations.

\subsubsection{Neighbouring Cases}\label{rf:sec:neighbours}
The sectional-curvature case~(0.6.2) of
\cite[Conjecture~0.6]{Ron26}, with
$\sec\ge-1$, remains open: such limits are locally contractible Alexandrov spaces
\cite{Kap06}, unlike $\Sigma_\infty$, while sectional lower bounds and LRVP do not imply
small-constant LRRP, as shown by products of a shrinking circle with smooth convex
surfaces converging to a surface with a conical point.

\section*{Acknowledgments and Tool Disclosure}
The author used Anthropic's Claude and OpenAI's Codex for drafting and restructuring
the exposition, for checking arguments against cited sources, for the numerical
evaluation of the constants in \cref{rf:sec:example}, and for the preparation of the
figures. The author takes full responsibility for the mathematical claims and the exposition.
\printbibliography[heading=bibintoc]
\end{document}